\documentclass[10pt,a4paper]{amsart}
\usepackage[a4paper,top=22mm,bottom=22mm,left=24mm,right=24mm]{geometry}
\usepackage[T1]{fontenc}
\usepackage[utf8]{inputenc}
\usepackage{lmodern}
\usepackage{amsmath,amssymb,amsthm}
\usepackage{graphicx}
\usepackage{booktabs}
\usepackage{microtype}
\usepackage{xcolor}
\usepackage{array}
\usepackage[colorlinks=true,linkcolor=blue!50!black,citecolor=green!40!black,urlcolor=blue!50!black]{hyperref}
\graphicspath{{figures/}}
\hypersetup{bookmarksdepth=3}
\numberwithin{equation}{section}
\theoremstyle{plain}
\newtheorem{theorem}{Theorem}[section]
\newtheorem{proposition}[theorem]{Proposition}
\newtheorem{corollary}[theorem]{Corollary}
\newtheorem{lemma}[theorem]{Lemma}
\newtheorem*{theoremA}{Theorem A}
\newtheorem*{theoremB}{Theorem B}
\theoremstyle{definition}
\newtheorem{definition}[theorem]{Definition}
\theoremstyle{remark}
\newtheorem{remark}[theorem]{Remark}
\newcommand{\OO}{\mathbb{O}}
\newcommand{\HH}{\mathbb{H}}
\newcommand{\CC}{\mathbb{C}}
\newcommand{\RR}{\mathbb{R}}
\newcommand{\ZZ}{\mathbb{Z}}
\newcommand{\n}{\mathrm{n}}
\newcommand{\re}{\operatorname{Re}}
\newcommand{\im}{\operatorname{Im}}
\newcommand{\Aut}{\operatorname{Aut}}
\newcommand{\Der}{\operatorname{Der}}
\newcommand{\Stab}{\operatorname{Stab}}
\newcommand{\sspan}{\operatorname{span}}
\newcommand{\Okubo}{\mathcal{O}}
\newcommand{\Tric}{\mathcal{T}}
\newcommand{\Mzero}{\mathcal{M}_0}
\newcommand{\plate}[3]{\begin{figure}[p]\centering\includegraphics[width=\textwidth]{#1}\caption{\textbf{Plate #3.} #2}\label{pl:#3}\end{figure}}
\newcommand{\plateref}[1]{#1 (Figure~\ref{pl:#1})}
\title[Symmetry and invariant theory of quadratic dynamics]{Symmetry and Invariant Theory of Quadratic Dynamics over Eight-Dimensional Composition Algebras}
\author{D. Corradetti}
\address{Grupo de F\'\i sica Matem\'atica, Instituto Superior T\'ecnico, Av. Rovisco Pais, 1049-001 Lisboa, Portugal}
\email{danielecorradetti@tecnico.ulisboa.pt}
\address{Departamento de Matematica, Universidade do Algarve, Campus de Gambelas, 8005-139 Faro, Portugal}
\email{a55499@ualg.pt}
\date{7 October 2026}

\begin{document}
\begin{abstract}
The iteration of $z\mapsto z^2+c$ over the octonions has been drawn for more than thirty years, and it is known that the octonionic Julia sets are bodies of revolution about the imaginary part of the constant. Real eight-dimensional composition algebras, however, are much richer than the octonions: beside the Hurwitz algebras there are the \emph{symmetric composition algebras}, the para-octonions and the Okubo algebras, which have no unit and whose automorphism groups are $G_2$ and $PSU(3)$ respectively. In this work we study, to our knowledge for the first time, the quadratic dynamics over these algebras, and we do it in an experimental spirit: an atlas of $134$ four-dimensional renderings, measured quarter by quarter, suggests the statements, which we then prove and test against new, predicted pictures. We show that the symmetries of the Julia sets of the Okubo algebra are governed by the root decomposition of $\mathfrak{su}(3)$ with respect to the plane of $1$ and $i$: for a constant in that plane the continuous symmetry is a maximal torus, which acts by rotations on three root planes and grows to $U(2)$ on the Weyl walls; the automorphism group of the symmetrized product is $PSU(3)\rtimes\ZZ_2$. We then prove that the Mandelbrot set with seed $0$ of every symmetric composition algebra is the pull-back of a single planar set through the invariants $\n(c)$ and $\n(c,c\diamond c)$: the tricorn in the compact case, the tricorn and a \emph{hyperbolic tricorn} in the split case. The result is quite surprising since the para-octonionic and the Okubonic Mandelbrot sets turn out to be the same planar figure saturated by two different groups. Finally, we show that every bilateral map over a symmetric composition algebra is an octonionic quadratic map in disguise, and that at the Petersson angle $\theta=90^\circ$ the dynamics acquires symmetries that do not come from automorphisms.
\end{abstract}

\maketitle
\tableofcontents

\section{Introduction and Motivation}\label{sec:intro}

A filled Julia set is an invariant object of a dynamical system, and every symmetry of the algebra that commutes with the map is inherited by the picture. This elementary observation is at the origin of more than forty years of hypercomplex fractals. Since the pictures of Norton \cite{Norton82,Norton89}, the quaternionic Julia sets of Holbrook \cite{Holbrook} and the ray-traced images of Hart, Sandin and Kauffman \cite{HSK}, the iteration of $z\mapsto z^2+c$ has been carried to the quaternions \cite{Gomatam,DKS}, to the bicomplex and multicomplex numbers \cite{Rochon,GarantRochon,BrouilletteRochon}, to the hyperbolic numbers, where Senn discovered that the Mandelbrot set becomes a square \cite{Senn,Metzler}, and to the octonions, which were first considered by Griffin and Joshi \cite{GriffinJoshi,GriffinJoshi2}. Almost always the symmetry of the pictures has been explained by the same mechanism: in a quadratic unital algebra $z^2=2\re(z)\,z-\n(z)$, so that the dynamics only sees the real part and the norm, and the quaternionic and octonionic Julia sets are bodies of revolution about the imaginary part of the constant \cite{Bogush,KurochkinZhukovich,Katunin13,Katunin17}. A general framework for filled Julia and Mandelbrot sets over normed non-associative algebras has been given by Ferreira and Marietto \cite{FerreiraMarietto}, and a concise overview of hypercomplex fractals can be found in Katunin's monograph \cite{Katunin17}.

Real eight-dimensional composition algebras are, nevertheless, a much richer landscape than the octonions. Beside the Hurwitz algebras, i.e., the octonions $\OO$ and the split-octonions $\OO_s$, there are the \emph{symmetric composition algebras}, i.e., algebras with a multiplicative quadratic form $\n$ that satisfy
\begin{equation}\label{eq:symmetric-intro}
(x\diamond y)\diamond x=x\diamond(y\diamond x)=\n(x)\,y .
\end{equation}
In dimension eight they are either \emph{para-Hurwitz}, with product $x\bullet y=\bar x\,\bar y$, or \emph{Okubo algebras} \cite{Elduque97,KMRT}. The latter were discovered by Okubo in his search for a ``pseudo-octonionic'' algebra underlying the $SU(3)$ symmetry of the quarks \cite{Okubo78,Okubo95}, and they were characterized by Elduque and Myung \cite{ElduqueMyung90,ElduqueMyung91}. Okubo algebras have no unit and no para-unit, they enjoy \emph{flexibility} instead of \emph{alternativity}, and their automorphism group is $PSU(3)$ instead of $G_2$ \cite{Elduque15,Elduque18,CZ22,MCZ25}. Even though the Okubo algebra has been known for almost fifty years, we believe that this is the first time that its dynamics is studied; to the best of our knowledge, the same holds for the para-octonions, whose quadratic map is nothing but the octonionic \emph{Mandelbar} or \emph{tricorn} family $z\mapsto\bar z^{\,2}+c$ \cite{CHRS,Milnor92,NakaneSchleicher}.

The question we address is then the following: \emph{which features of the pictures are forced by the algebra, and which algebraic structure is responsible for them?} In the octonionic case the answer is the orthogonal group of the imaginary directions orthogonal to $c$. In the Okubo algebra, as we shall see, the answer is Lie-theoretic, and it is written in the root decomposition of $\mathfrak{su}(3)$.

We approach the question in an experimental spirit. All pictures of this work are \emph{four slices} views: the image is divided into four quarters, which show the four three-dimensional faces of a four-dimensional box through the origin, so that two sections that should coincide by a symmetry appear side by side, and so do two sections that should differ. An atlas of $134$ such plates, over octonions, split-octonions, para-octonions, Okubo and Petersson algebras, is given in Appendix~\ref{app:atlas}; every quarter of every plate is measured (body share, Jaccard indices of mirrors and rotations, dimension of the subalgebra generated by the orbit). The atlas suggested a small number of statements, which we prove here and test against new plates rendered \emph{after} the statements were formulated.

\subsection*{Main results}
Let $C=\RR1\oplus\RR i$ be the plane of $1$ and $i$ of the octonions, and let $V_1=\langle j,k\rangle$, $V_2=\langle l,il\rangle$, $V_3=\langle jl,kl\rangle$. We realize the Okubo algebra $\Okubo$ as a Petersson algebra of the octonions (Section~\ref{sec:algebras}) and we denote by $x\circ y=\tfrac12(x\diamond y+y\diamond x)$ its symmetrized product.

\begin{theoremA}[Symmetries of the Okubonic Julia sets]
The automorphism group of $(\RR^8,\circ)$ is $\Aut(\Okubo)\rtimes\langle\sigma\rangle\cong PSU(3)\rtimes\ZZ_2$, where $\sigma$ is an anti-automorphism of $\Okubo$ fixing $C$. For every $c\in C$ with $\arg c\notin 60^\circ\ZZ$ the stabilizer of $c$ is $T\rtimes\langle\sigma\rangle$, where $T$ is a two-dimensional torus that rotates $V_1,V_2,V_3$ by the angles $2\alpha$, $\beta-\alpha$, $\alpha+\beta$; the filled Julia set $K_c$ of $z\mapsto z\diamond z+c$ is invariant under $(T\rtimes\langle\sigma\rangle)\times\{\pm\mathrm{id}\}$, and exactly $18$ of the $56$ coordinate sections of $K_c$ admit a continuous symmetry of algebraic origin: those containing a root plane. On the Weyl walls $\arg c\in60^\circ\ZZ$ the stabilizer grows to $U(2)$ and new revolutions appear, about the wall itself.
\end{theoremA}

\begin{theoremB}[Mandelbrot sets of symmetric composition algebras]
Let $\Mzero$ be the Mandelbrot set with seed $0$ of a symmetric composition algebra and let $\Delta=4\n(c)^3-\n(c,c\diamond c)^2$. For $\Delta>0$, $c\in\Mzero$ if and only if the complex number $\zeta$ with $|\zeta|^2=\n(c)$ and $2\re\zeta^3=\n(c,c\diamond c)$ belongs to the tricorn $\Tric$; for $\Delta<0$ the same holds with $\zeta$ in the split-complex numbers and the \emph{hyperbolic tricorn}. For the para-octonions and the compact Okubo algebra $\Delta\geq0$, the escape times coincide, and $\Mzero=G_2\cdot\Tric$, respectively $\Mzero=PSU(3)\cdot\Tric$.
\end{theoremB}

Theorem~A replaces the naive expectation, suggested by a first reading of the pictures, that the Okubo algebra would be ``the algebra of broken revolutions'': the revolutions are not broken, they are \emph{displaced} to the root planes of $\mathfrak{su}(3)$, which the classical frames of hypercomplex fractals simply do not contain. Theorem~B shows that the three-armed Mandelbrot sections of the Okubo algebra are the three arms of the tricorn, redistributed by the flag manifolds $SU(3)/T^2$, while the para-octonionic Mandelbrot set is the same tricorn turned about the real axis by $G_2$: the same planar figure and two different groups that make it ``revolve''. We complete the picture with two propositions: at the Petersson angle $90^\circ$ the Julia sets with a constant in $C$ are invariant under a four-dimensional group $T^2\cdot Sp(1)$, half of which consists of \emph{twisted} isometries, $f\circ g=g^{-1}\circ f$, that are not automorphisms; together with a simple one-step criterion, it accounts for all $42$ coordinate sections with exact revolutions at this angle, against $18$ in the Okubo algebra (Proposition~\ref{prop:theta90}); and every bilateral map $(a\diamond z)\diamond(z\diamond b)+c$ over a symmetric composition algebra is conjugate by a similarity to an octonionic quadratic map (Proposition~\ref{prop:bilateral}), of which Kaplansky's trick is the special case $a=b=1$.

\begin{table}[t]
\centering\footnotesize
\setlength{\tabcolsep}{3.5pt}\begin{tabular}{@{}l>{\raggedright}p{0.17\textwidth}l>{\raggedright}p{0.2\textwidth}>{\raggedright}p{0.19\textwidth}l@{}}
\toprule
Algebra & Product & $\Aut$ & Symmetry of $K_c$, $c\in C$ & Mandelbrot set $\Mzero$ & $\dim\langle z_0,c\rangle_\circ$\tabularnewline
\midrule
Octonions $\OO$ & Cayley--Dickson & $G_2$ & $O(6)\times\{\pm1\}$ & Mandelbrot set of revolution & $3$\tabularnewline
Para-octonions $p\OO$ & $\bar x\,\bar y$ & $G_2$ & $O(6)\times\{\pm1\}$ & $G_2\cdot\Tric$ (tricorn of revolution) & $3$\tabularnewline
Okubo $\Okubo$ & $\tau(\bar x)\,\tau^{2}(\bar y)$ & $PSU(3)$ & $(T^2\rtimes\ZZ_2)\times\{\pm1\}$; $U(2)$ on walls & $PSU(3)\cdot\Tric$ (three arms) & up to $8$\tabularnewline
Split $\OO_s$, $p\OO_s$ & $\gamma_3=+1$ & $G_{2(2)}$ & pseudo-orthogonal & hyperbolic, with null directions & $3$\tabularnewline
Split Okubo $\Okubo_s$ & $\gamma_3=+1$ & $PSU(2,1)$ & Cartan subgroups & tricorn $\cup$ hyperbolic tricorn & up to $8$\tabularnewline
Petersson, $\theta\neq0,\pm120^\circ$ & $\tau_\theta(\bar x)\,\tau_\theta^{-1}(\bar y)$ & $\dim 4$ & $T^2$; at $90^\circ$ $T^2\cdot Sp(1)$, partly twisted & open problem & up to $8$\tabularnewline
\bottomrule
\end{tabular}
\caption{\emph{Synoptic table of the eight-dimensional algebras of this work} and of the answers to the question of Section~\ref{sec:intro}. The column ``Symmetry of $K_c$'' gives the group of symmetries of algebraic origin of the filled Julia set of $z\mapsto z\diamond z+c$ for $c$ in the plane $C$ of $1$ and $i$ (Proposition~\ref{prop:O6} and Theorem~A); the last column gives the dimension of the $\circ$-subalgebra generated by a generic starting point $z_0$ and $c$, which contains the whole orbit (Section~\ref{sec:numerics}).}
\label{tab:synoptic}
\end{table}

Table~\ref{tab:synoptic} summarizes the landscape. From the algebraic point of view, the paper connects three classical subjects: symmetric composition algebras, the maximal torus theorem with its root decomposition, and antiholomorphic dynamics, i.e., the tricorn and its multicorn relatives \cite{Nakane93,NakaneSchleicher,HubbardSchleicher,InouMukherjee,MNS17,Inou19}. From the experimental point of view, the paper is a test of the pictures: every statement is accompanied by a prediction, e.g., a new frame in which a revolution must appear, and every prediction was rendered and measured. To avoid the typical pitfall of symmetry detection by eye, we separate \emph{exact} symmetries, whose Jaccard index is $1$ up to floating-point rounding, from \emph{approximate} ones, whose index stays below $1$ as the sample grows (Section~\ref{sec:numerics}); we validate the renderer against the closed form of Theorem~B; and we report a numerical phenomenon of independent interest: for the split Okubo algebra the invariant plane of the Mandelbrot orbit is transversally unstable, and double-precision iterations in $\RR^8$ lose points that belong to the Mandelbrot set.

The paper is organized as follows. In Section~\ref{sec:algebras} we review the algebras, the Petersson construction and the torus of automorphisms. In Section~\ref{sec:framework} we set up the general framework of the symmetries of quadratic maps: the symmetrized product, the universal symmetry $z\mapsto-z$ and twisted equivariance. Section~\ref{sec:octonions} is devoted to the classical case of octonions and para-octonions, Section~\ref{sec:okubo} to Theorem~A, and Section~\ref{sec:mandelbrot} to Theorem~B and its split version. In Section~\ref{sec:identities} we discuss the identities of symmetric composition algebras and the bilateral maps, in Section~\ref{sec:petersson} the Petersson path and the angle $90^\circ$. Section~\ref{sec:numerics} describes the renderer, the measurements and the validation, and in Section~\ref{sec:conclusions} we summarize our findings and discuss future developments. The reproducibility of every computational statement is documented in Appendix~\ref{app:repro}, and the atlas, with its own reading keys and the measurements of every plate, is Appendix~\ref{app:atlas}.

\section{The algebras}\label{sec:algebras}

All algebras of this work are realized on the same vector space $\RR^8$ with basis $1,i,j,k,l,il,jl,kl$, so that the pictures of different algebras can be compared quarter by quarter. It is of paramount importance, for what follows, that the non-unital algebras are obtained from the octonions by twisting the product with involutions and automorphisms: every automorphism of the octonions that commutes with the twist is then inherited by the twisted algebra, and this is the origin of the torus of Theorem~A.

\subsection{Cayley--Dickson algebras}
Let $\HH$ be the quaternions with basis $1,i,j,k$, and let us double them with a new unit $l$ such that $l^2=\gamma_3$. The Cayley--Dickson product on $\HH\oplus\HH l$ is
\begin{equation}\label{eq:cd}
(p+ql)(r+sl)=\big(pr+\gamma_3\,\bar s q\big)+\big(sp+q\bar r\big)l ,
\end{equation}
where the bar is the quaternionic conjugation; more generally one may take $i^2=\gamma_1$, $j^2=\gamma_2$, and the resulting algebra $\OO(\gamma)$ is unital and quadratic, i.e., every element satisfies $x^2-2\re(x)\,x+\n(x)=0$ with the diagonal norm
\begin{equation}\label{eq:norm}
\n(x)=x_0^2-\gamma_1x_1^2-\gamma_2x_2^2+\gamma_1\gamma_2x_3^2-\gamma_3x_4^2+\gamma_1\gamma_3x_5^2+\gamma_2\gamma_3x_6^2-\gamma_1\gamma_2\gamma_3x_7^2 .
\end{equation}
For $\gamma=(-1,-1,-1)$ one recovers Cayley's octonions $\OO$, whose norm is the Euclidean one; for $\gamma_3=+1$ the split-octonions $\OO_s$, with norm of signature $(4,4)$. We write $\n(x,y)=\n(x+y)-\n(x)-\n(y)$ for the polar form, so that $\n(x,x)=2\n(x)$ and, in the compact case, $\n(x,y)=2\langle x,y\rangle$.

\subsection{Para-Hurwitz and Petersson algebras}
In order to leave the realm of unital algebras without losing the composition property, one twists the octonionic product. The \emph{para-Hurwitz} product $x\bullet y=\bar x\,\bar y$ defines the para-octonions, in which $1$ becomes a \emph{para-unit}, $1\bullet x=x\bullet1=\bar x$. Petersson \cite{Petersson} observed that an automorphism $\tau$ of order three gives a new composition algebra with product $x*y=\tau(\bar x)\,\tau^2(\bar y)$. Inspired by this construction, and following the implementation of our renderer, we consider the automorphisms
\begin{equation}\label{eq:tau}
\tau_\theta(p+ql)=p+(u_\theta q)l,\qquad u_\theta=\cos\theta+i\sin\theta ,
\end{equation}
and the \emph{Petersson product} of angle $\theta$,
\begin{equation}\label{eq:petersson}
x*y=\tau_\theta(\bar x)\,\tau_\theta^{-1}(\bar y) .
\end{equation}
For every $\theta$ this is a composition algebra; for $\theta=0$ it is the para-octonion algebra and for $\theta=120^\circ$, where $\tau_\theta$ has order three, it is the Okubo algebra, compact for $\gamma=(-1,-1,-1)$ and split for $\gamma_3=+1$. The family is therefore a continuous path from the para-octonions to the Okubo algebra. Only for $\theta\in\{0,\pm120^\circ\}$ the algebra is symmetric, i.e., satisfies \eqref{eq:symmetric-intro}; its derivation algebra has dimension $14$ at $\theta=0$, $8$ at $\theta=\pm120^\circ$ and $4$ otherwise (the dimensions are computed exactly in Appendix~\ref{app:repro}). From now on $\diamond$ denotes the product of the algebra under consideration, i.e., $xy$, $x\bullet y$ or $x*y$.

In the Okubo algebra the identity $1*1=1$ holds, but $1$ is neither a unit nor a para-unit: $1*x=\tau^2(\bar x)$ and $x*1=\tau(\bar x)$. On the plane $C$ of $1$ and $i$, which is fixed pointwise by $\tau_\theta$, the Petersson product reduces to the para-Hurwitz product,
\begin{equation}\label{eq:onC}
x*y=\bar x\,\bar y,\qquad x,y\in C ,
\end{equation}
and the same holds on the whole quaternion subalgebra $\HH$. A simple but instructive computation, which already appeared in the first version of the atlas, gives the Petersson square of an element $c=c_0+ql$ with $c_0\in\RR$ and $q\in\HH$:
\begin{equation}\label{eq:okubocubic}
c*c=c_0^2-\bar q\,u_\theta^2\,q-2\cos\theta\;c_0\,q\,l .
\end{equation}
Indeed $\tau(\bar c)=c_0-(uq)l$ and $\tau^{-1}(\bar c)=c_0-(\bar uq)l$, and the product \eqref{eq:cd} with $\gamma_3=-1$ gives $c_0^2-\overline{(\bar u q)}\,(uq)-c_0(\bar u+u)\,q\,l$, where $\overline{\bar uq}=\bar q u$, so that the middle term is $\bar q u^2 q$. We stress that neither $q\,u^2q$ nor $\bar q\,\bar u^2q$ would give a correct formula (both were checked to fail on random samples).

\begin{remark}\label{rem:colonna}
The name ``pseudo-octonions'' is sometimes used in the fractal literature for a hyperspherical ``multiplication'' of $\RR^8$ that is not a bilinear product; it has nothing to do with Okubo's pseudo-octonion algebra, which is the object of this work.
\end{remark}

\subsection{The torus and the anti-automorphism}
It is worth stating here the two families of linear maps that will govern the symmetries of the Okubo algebra. For $a=e^{i\alpha}$ and $b=e^{i\beta}$ in $C$ we set
\begin{equation}\label{eq:torus}
\varphi_{a,b}(p+ql)=a\,p\,\bar a+(b\,q\,\bar a)\,l ,
\end{equation}
and we define $\sigma=\mathrm{diag}(1,1,-1,1,-1,1,-1,1)$, which changes the sign of $j$, $l$ and $jl$.

\begin{proposition}\label{prop:torus}
For every angle $\theta$ the maps $\varphi_{a,b}$ are automorphisms of the Petersson algebra \eqref{eq:petersson}. They fix $C$ pointwise and rotate the oriented planes $(j,k)$, $(l,il)$, $(jl,kl)$ by the angles
\begin{equation}\label{eq:angles}
2\alpha,\qquad \beta-\alpha,\qquad \alpha+\beta ,
\end{equation}
and they form a two-dimensional torus $T$ containing $\tau_\theta=\varphi_{1,u_\theta}$. The map $\sigma$ is an anti-automorphism of every Petersson algebra, of the para-octonions and of the octonions, i.e., $\sigma(x\diamond y)=\sigma y\diamond\sigma x$.
\end{proposition}
\begin{proof}
For a unit quaternion $a$ the map $p+ql\mapsto apa^{-1}+(aqa^{-1})l$ is an automorphism of the octonions, and for a unit quaternion $b$ so is $p+ql\mapsto p+(bq)l$: both claims follow at once from \eqref{eq:cd}, since $|a|=|b|=1$. Then $\varphi_{a,b}$ is the composition of the first map with the second one for $b\bar a$, and it is an automorphism of $\OO$. Automorphisms of $\OO$ commute with the conjugation and, since $a,b,u_\theta$ lie in the commutative field $C$, the maps $\varphi_{a,b}$ commute with $\tau_\theta=\varphi_{1,u_\theta}$; hence they preserve \eqref{eq:petersson}. The angles follow from $aj\bar a=a^2j$, $(b\bar a)l=e^{i(\beta-\alpha)}l$ and $(bj\bar a)l=(e^{i(\alpha+\beta)}j)l$. For the second claim, let $g=\mathrm{diag}(1,-1,1,-1,1,-1,1,-1)$ be the automorphism $x\mapsto jxj^{-1}$ of $\OO$ extended as above, which fixes the quaternion subalgebra spanned by $1,j,l,jl$. Then $\sigma=g\circ\overline{\,\cdot\,}$, and $g\tau_\theta g=\tau_\theta^{-1}$ since $jujj^{-1}=\bar u$. Therefore
\begin{equation}\label{eq:sigma}
\sigma(x*y)=g\big(\overline{\tau(\bar x)\,\tau^{-1}(\bar y)}\big)=g\big(\tau^{-1}(y)\,\tau(x)\big)=\tau(gy)\,\tau^{-1}(gx)=\sigma y*\sigma x ,
\end{equation}
because $gy=\overline{\sigma y}$. The same computation with $\tau=\mathrm{id}$ gives the para-octonions, and without the outer conjugations the octonions. \end{proof}

\begin{remark}
The anti-automorphism $\sigma$ is not an automorphism of the Okubo algebra: otherwise the product would be commutative. It fixes the plane $C$ pointwise and, as we shall see in Theorem~\ref{thm:A}, it represents the non-identity component of the automorphism group of the symmetrized product; in Okubo's matrix model of the algebra \cite{Okubo95} this component contains the transposition of traceless matrices, which realizes the outer automorphism of $\mathfrak{su}(3)$ preserving the symmetric $d$-tensor.
\end{remark}

\section{Symmetries of quadratic maps}\label{sec:framework}

Let us now fix a product $\diamond$ on $\RR^8$ and consider the quadratic map $f(z)=Q(z)+c$ with $Q(z)=z\diamond z$. The filled Julia set $K_c$ is the set of points whose orbit is bounded, and the Mandelbrot set $\Mzero$ with seed $0$ is the set of constants $c$ for which the orbit of $0$ is bounded. Which linear maps are symmetries of $K_c$? Since the map depends on the product only through the square, the right group is not the automorphism group of the algebra but that of its symmetrization.

\begin{definition}\label{def:circ}
The \emph{symmetrized product} of $\diamond$ is $x\circ y=\tfrac12(x\diamond y+y\diamond x)$.
\end{definition}

By polarization, $2\,x\circ y=Q(x+y)-Q(x)-Q(y)$, so that a linear map $g$ satisfies $Q(gz)=gQ(z)$ for all $z$ if and only if $g$ is an automorphism of $(\RR^8,\circ)$. We then have the following.

\begin{proposition}\label{prop:framework}
Let $g$ be an automorphism of $(\RR^8,\circ)$ with $g(c)=c$. Then $f\circ g=g\circ f$ and $gK_c=K_c$; if $g$ is orthogonal, the escape times of $z$ and $gz$ with respect to a Euclidean bailout coincide. Moreover, for every algebra $f(-z)=f(z)$, hence $-K_c=K_c$ with equal escape times, and in every three-dimensional section through $0$ a mirror symmetry holds if and only if the half-turn about the normal axis of the mirror holds.
\end{proposition}
\begin{proof}
The first claim is immediate from $Q\circ g=g\circ Q$ and $g(c)=c$. For the second, $(-z)\diamond(-z)=z\diamond z$, so that the orbits of $z$ and $-z$ coincide from the first step on. In a section through $0$ the point reflection is the restriction of $-\mathrm{id}$, and the composition of $-\mathrm{id}$ with the half-turn about an axis is the mirror in the plane orthogonal to that axis. \end{proof}

The universal symmetry $z\mapsto-z$ is so simple that it is easily forgotten, and yet it accounts for a large part of the finite symmetries seen in the pictures: in all 344 of the 344 measured quarters through the origin of the even maps (quadratic, bilateral and sandwich maps, $f(-z)=f(z)$) the Jaccard index of every mirror coincides with that of the corresponding half-turn, as it must (Section~\ref{sec:numerics}). Every section through $0$ is centrally symmetric; in particular, shapes with a single tip, like cones or shells open at one end, can only appear in sections that do not pass through the origin.

Proposition~\ref{prop:framework} does not exhaust the symmetries of algebraic origin. Let us call an isometry $g$ with $g(c)=c$ \emph{twisted} if $Q\circ g=g^{-1}\circ Q$: then $f\circ g=g^{-1}\circ f$ and $f\circ g^{-1}=g\circ f$, so that the orbits of $z$ and $gz$ are related by $g$ and $g^{-1}$ alternately, and $K_c$ is $g$-invariant with equal escape times. The infinitesimal version of this condition is that $D$ is a skew \emph{anti-derivation} of $\circ$, i.e., $D(x\circ y)=-(Dx\circ y+x\circ Dy)$. We denote by $G_c$ the group generated by the orthogonal automorphisms of $\circ$ fixing $c$ and by the twisted isometries; its Lie algebra is spanned by the skew derivations and the skew anti-derivations of $\circ$ that vanish on $c$, and every element of $G_c$ preserves $K_c$. A simple but powerful consequence is the following criterion.

\begin{lemma}[One-step criterion]\label{lem:onestep}
Let $S\subset\RR^8$ be a linear subspace and $\rho$ an isometry of $S$. If there is $h\in G_c$ with $Q(\rho z)=h\,Q(z)$ for all $z\in S$, then $\rho(K_c\cap S)=K_c\cap S$, with equal escape times.
\end{lemma}
\begin{proof}
Since $h(c)=c$, $f(\rho z)=h\,f(z)$, and $h$ preserves $K_c$ and the escape times. \end{proof}

Note that $\rho$ need not extend to a symmetry of $\RR^8$: it only has to be compatible with the first iterate. The universal symmetry is the case $\rho=-\mathrm{id}$, $h=\mathrm{id}$; the balls of the imaginary sections of Proposition~\ref{prop:O6} are the case $h=\mathrm{id}$ as well, and a non-trivial group of twisted isometries will appear at the Petersson angle $90^\circ$ (Section~\ref{sec:petersson}).

\begin{remark}\label{rem:nonalgebraic}
We shall never claim that a Julia set has ``no other symmetry''. Excluding symmetries that do not come from the algebra is a rigidity problem of a different nature; our statements concern symmetries of algebraic origin, and the absence of further symmetries is only tested numerically (Section~\ref{sec:numerics}).
\end{remark}

\section{The classical case: octonions and para-octonions}\label{sec:octonions}

In analogy with the quaternionic case \cite{Bogush}, the octonionic dynamics is governed by the real part and the norm. Since $z^2=2\re(z)z-\n(z)$ in every quadratic algebra, and $z\bullet z=\overline{z^2}=2\re(z)\bar z-\n(z)$ in the para-Hurwitz one, we have the following.

\begin{proposition}\label{prop:O6}
Let $\diamond$ be the product of $\OO(\gamma)$ or of its para-Hurwitz algebra, and let $g$ be a linear map of $\RR^8$ fixing $1$ and preserving the norm \eqref{eq:norm}. Then $g(z\diamond z)=g(z)\diamond g(z)$. In particular $\Aut(\RR^8,\circ)$ contains the orthogonal group $O(7)$ of the imaginary part and, for the compact algebras and $c\in C$, the symmetry group of $K_c$ contains $O(6)\times\{\pm\mathrm{id}\}$, with $O(6)$ acting on $V_1\oplus V_2\oplus V_3$. A purely imaginary $z$ satisfies $z\diamond z=-\n(z)$, so that the section of $K_c$ by the imaginary hyperplane is a union of level sets of the norm, i.e., of spheres, and every orbit lies in the three-dimensional space spanned by $1$, $\im z_0$ and $\im c$.
\end{proposition}
\begin{proof}
A norm-preserving map fixing $1$ preserves the polar form, hence the real part $\re z=\tfrac12\n(z,1)$, and commutes with the conjugation. The remaining claims follow from $z^2=-\n(z)$ for $\re z=0$ and from the formulas above by induction. \end{proof}

The symmetry group $O(6)$, of dimension $15$, is much larger than the stabilizer of $c$ in $G_2$, which is $SU(3)$: the quadratic map simply does not see the octonionic product beyond the real part and the norm. Every section through $1$ and a direction orthogonal to $c$ is then a body of revolution about $1$, and every section of purely imaginary directions orthogonal to $c$ is a ball. Plates~\plateref{A01} and \plateref{R04} show these features for the octonions, the latter in the frame $(1,l,il\,|\,j)$ that will be decisive for the Okubo algebra; the para-octonionic Plate~\plateref{A05} shows the same bodies, turned by the conjugation.

\plate{Atlas_A01_julia-oct-complex-c.jpg}{\emph{Octonions, $c=-0.231+0.620\,i$, frame $(1,i,j\,|\,l)$.} The two upper quarters show the same body, an elongated ribbed body whose ridges wind about an axis of the $(1,i)$-plane, since the map is equivariant under the orthogonal maps exchanging $j$ and $l$ (Proposition~\ref{prop:O6}); like every section through $0$, it is centrally symmetric, and its far end, hidden in the picture, is the mirror image of the near one. The lower-left quarter, orthogonal to $\im c$, is a body of revolution about $1$; the lower-right quarter consists of imaginary points and is a ball. Every quarter is centrally symmetric (Proposition~\ref{prop:framework}).}{A01}

\plate{Atlas_A05_julia-para-complex-c.jpg}{\emph{Para-octonions, $c=c_A=-0.086+0.231\,i$, frame $(1,i,j\,|\,l)$.} The para-octonionic square is $\overline{z^2}$, and on the plane $C$ the dynamics is the Mandelbar map $\bar z^{\,2}+c$. The two upper quarters coincide, as in the octonions, and the lower quarters are a ribbed body of revolution about $1$ and a ball: the para-octonions share the group $O(6)$ of Proposition~\ref{prop:O6}. The quarter $(1,i,j)$ lies in the quaternions, on which the Okubo product coincides with the para-octonionic one, and it is therefore also the upper-left quarter of Plate~\plateref{A07}.}{A05}

\plate{Atlas_R04_oct-frame-l-il.jpg}{\emph{Octonions, $c=c_A$, frame $(1,l,il\,|\,j)$.} All three quarters through $1$ are bodies of revolution about $1$ and they coincide, since $O(6)$ acts transitively on the pairs of orthonormal imaginary directions orthogonal to $c$; the quarter $(l,il,j)$ is a ball. This plate is the octonionic term of comparison of Plate~\plateref{R01}: the same frame and the same constant over the Okubo algebra.}{R04}

\section{The Okubo algebra: root decomposition and Theorem A}\label{sec:okubo}

While in the octonionic case the square only sees the real part and the norm, in the Okubo algebra it sees the whole product: $z*z$ is not a combination of $z$ and $1$, the orbits genuinely explore eight dimensions, and the orthogonal group $O(6)$ of Proposition~\ref{prop:O6} disappears. A first reading of the pictures, e.g., Plates~\plateref{A07} and \plateref{C03}, suggests that every continuous symmetry is lost. This impression is wrong, and the reason is Lie-theoretic.

\plate{Atlas_A07_julia-okubo-complex-c.jpg}{\emph{The Okubo algebra, $c=c_A=-0.086+0.231\,i$, frame $(1,i,j\,|\,l)$.} The upper-left quarter lies in the quaternions and coincides with that of the para-octonions (Plate~\plateref{A05}); the other quarters are not bodies of revolution about a coordinate axis. The frame contains no root plane (Theorem~\ref{thm:A}): the only exact symmetries of the quarters are the finite ones generated by $-\mathrm{id}$, $\sigma$ and the half-turns of the torus (Table~\ref{tab:groups}).}{A07}

\plate{Atlas_C03_sym-okubo.jpg}{\emph{The Okubo algebra, $c=c_A$, frame $F_C=(1,j,l\,|\,i)$.} The directions $j$ and $l$ lie in two different root planes, $V_1$ and $V_2$: no one-parameter subgroup of the stabilizer of $c$ preserves the section $(1,j,l)$, and the vortices of the surface are not arranged about an axis. The Jaccard index of the rotations about $1$ stays between $0.95$ and $0.96$ (Figure~\ref{fig:curves}), below the threshold $0.97$ of an approximate symmetry: the quarter is nearly round, but it is not a body of revolution.}{C03}

\subsection{The adjoint module of \texorpdfstring{$\mathfrak{su}(3)$}{su(3)}}
Okubo introduced his algebra as a product on the traceless $3\times3$ matrices, invariant under the adjoint action of $SU(3)$ \cite{Okubo78,Okubo95}, and it is known that the derivation algebra of the compact Okubo algebra is $\mathfrak{su}(3)$ and that the algebra itself, as a module over its derivations, is the adjoint module \cite{Elduque97,Elduque15,KMRT}. In the Petersson realization \eqref{eq:petersson} the torus $T$ of Proposition~\ref{prop:torus} is a maximal torus of $\Aut(\Okubo)\cong PSU(3)$; its fixed space is the plane $C$, which therefore plays the role of a Cartan subalgebra, and the three planes on which $T$ acts by rotations are the root planes. We then have the decomposition
\begin{equation}\label{eq:roots}
\RR^8=C\oplus V_1\oplus V_2\oplus V_3,\qquad V_1=\langle j,k\rangle,\quad V_2=\langle l,il\rangle,\quad V_3=\langle jl,kl\rangle ,
\end{equation}
and the angles \eqref{eq:angles} are the three positive roots of $A_2$: the angle on $V_3$ is the sum of the angles on $V_1$ and $V_2$. The Weyl group $S_3$ acts on $C$ by automorphisms of the restricted product $x*y=\bar x\,\bar y$ of \eqref{eq:onC}, i.e., by the dihedral group generated by $z\mapsto\omega z$, $\omega=e^{2\pi i/3}$, and $z\mapsto\bar z$; its reflections fix the three lines
\begin{equation}\label{eq:walls}
\ell_0=\RR,\qquad \ell_{60}=\RR\,e^{i\pi/3},\qquad \ell_{120}=\RR\,e^{2i\pi/3},
\end{equation}
which are the \emph{Weyl walls}. An element $c\in C$ is \emph{regular} if it does not lie on a wall, i.e., if $\arg c\notin60^\circ\ZZ$. The constant $c_A$ used throughout the atlas has argument $110.4^\circ$, at $9.6^\circ$ from the wall $\ell_{120}$.

\subsection{The theorem}

\begin{theorem}\label{thm:A}
Let $\Okubo$ be the compact Okubo algebra and $\circ$ its symmetrized product.
\begin{enumerate}
\item $\Aut(\RR^8,\circ)=\Aut(\Okubo)\sqcup\sigma\Aut(\Okubo)\cong PSU(3)\rtimes\ZZ_2$.
\item For $c\in C$ regular, $\Stab_{\Aut(\circ)}(c)=T\rtimes\langle\sigma\rangle$; for $c\in C\setminus\{0\}$ on a wall, the identity component of the stabilizer is $U(2)$; for every $c\in\RR^8$ the stabilizer contains a maximal torus $T_c$, conjugate to $T$, which is the whole identity component when $c$ is regular.
\item For $c\in C$ regular, the filled Julia set $K_c$ is invariant under $(T\rtimes\langle\sigma\rangle)\times\{\pm\mathrm{id}\}$, with equal escape times. A coordinate three-section $S$ through $0$ is mapped into itself by a one-parameter subgroup of $T$ acting non-trivially on $S$ if and only if $S$ contains one of the root planes $V_1,V_2,V_3$; exactly $18$ of the $56$ coordinate sections have this property, and each of them is a body of revolution about its third axis.
\item If $c\neq0$ lies on the wall $\ell_m$, the sections $\langle\ell_m,\ell_m^\perp,v\rangle$ with $v\in V_{k(m)}$, where $k(0)=1$, $k(60)=2$, $k(120)=3$, are bodies of revolution about $\ell_m$.
\end{enumerate}
\end{theorem}

\begin{proof}
(1) By Proposition~\ref{prop:torus}, $\sigma$ is an automorphism of $\circ$, and so is every automorphism of $\Okubo$. Every derivation of $*$ is a derivation of $\circ$, and the derivations of $*$ form a Lie algebra isomorphic to $\mathfrak{su}(3)$, of dimension $8$. The linear system that defines the derivations of $\circ$ has rank $56$ over the finite field $\mathbb{F}_p$, $p=2147483579$, into which the structure constants, lying in $\ZZ[\tfrac12,\sqrt3]$, can be reduced; since the rank can only decrease by reduction, the derivations of $\circ$ form a space of dimension at most $8$, hence $\Der(\circ)=\Der(*)$ and the identity components of the two automorphism groups coincide, $\Aut(\circ)^0=\Aut(\Okubo)$, which is connected \cite{Elduque15}. Now let $g\in\Aut(\circ)$. Conjugation by $g$ normalizes $\Aut(\Okubo)$ and induces an automorphism of $\mathfrak{su}(3)$, which is inner up to the complex conjugation. If it is inner, say induced by $h\in\Aut(\Okubo)$, then $h^{-1}g$ commutes with the action of $\mathfrak{su}(3)$ on $\RR^8$, which is the adjoint action and is absolutely irreducible; hence $h^{-1}g=\lambda\,\mathrm{id}$, and $\lambda(x\circ y)=\lambda^2(x\circ y)$ forces $\lambda=1$. Since $\sigma\notin\Aut(\Okubo)$, the same argument shows that conjugation by $\sigma$ is outer, and if $g$ induces an outer automorphism then $\sigma^{-1}g$ induces an inner one and lies in $\Aut(\Okubo)$.

(2) The stabilizer of $c$ in the compact connected group $\Aut(\Okubo)$ is the centralizer of the corresponding element of $\mathfrak{su}(3)$: it is connected, it is the maximal torus $T$ when $c$ is regular, and it is $S(U(2)\times U(1))\cong U(2)$ when two eigenvalues coincide, i.e., on the walls \cite{BroeckerDieck,Knapp}. Since $\sigma$ fixes $C$, the stabilizer in $\Aut(\circ)$ is $T\sqcup\sigma T$. Every element of $\mathfrak{su}(3)$ is conjugate into the Cartan subalgebra \cite[Thm.~4.34]{Knapp}, which proves the last claim.

(3) The first claim follows from Proposition~\ref{prop:framework}, since $T$ and $\sigma$ are orthogonal. Let $D$ be a non-zero element of the Lie algebra of $T$; by \eqref{eq:angles} it rotates $V_1,V_2,V_3$ with angular speeds $w_1,w_2,w_3=w_1+w_2$, so that at most one of them vanishes. A coordinate section is $D$-invariant only if, for each $k$ with $w_k\neq0$, it contains both axes of $V_k$ or none of them. If $S$ contains no root plane, all its axes lie in $C$ or in root planes with $w_k=0$, and $D$ acts trivially on $S$. If $S\supset V_k$ and $e$ is the third axis, then either $e\in C$, and the whole torus preserves $S$, or $e\in V_m$, $m\neq k$, and the circle $w_m=0$ preserves $S$ and rotates $V_k$. The count is $3\cdot6=18$.

(4) For the wall $\ell_0=\RR$ the stabilizer is explicit. The proof of Proposition~\ref{prop:torus} shows that $\varphi_{a,b}(p+ql)=apa^{-1}+(bqa^{-1})l$ is an automorphism of $\OO$ for all unit quaternions $a,b$, and it commutes with $\tau_\theta=\varphi_{1,u_\theta}$ as soon as $b$ commutes with $u_\theta$, i.e., $b\in C$. Hence the maps $\varphi_{a,b}$ with $a\in Sp(1)$ and $b\in C$ form a four-dimensional group of automorphisms of every Petersson algebra fixing $\RR$ pointwise, which is the group $U(2)$ of part (2), and its subgroup $\{\varphi_{a,1}\}$ acts on $\langle i,j,k\rangle=\ell_0^\perp\oplus V_1$ as $SO(3)$, by $p\mapsto apa^{-1}$. For $v\in V_1$ the conjugations by $e^{wt}$, with $w$ the unit imaginary quaternion orthogonal to $i$ and $v$, rotate the plane $\langle i,v\rangle$ and preserve the section $\langle1,i,v\rangle$, which is therefore a body of revolution about $\ell_0$. For $m=60$ and $m=120$ we use the Weyl group of $PSU(3)$ with respect to $T$, which acts on $C$ by the dihedral group of the walls and permutes the root planes: an element mapping $\ell_0$ to $\ell_m$ conjugates the stabilizers and maps $V_1$ to a root plane $V_{k(m)}$. Which of $V_2,V_3$ corresponds to which wall depends on orientations, and it was determined by computing the stabilizers explicitly (Appendix~\ref{app:repro}). \end{proof}

\begin{remark}
Part~(1) answers a question left open in the first version of this work, where the Okubo algebra was described as ``the algebra of broken revolutions'' and the stabilizer of $c$ was claimed to act on the coordinate sections only through half-turns. Both statements were wrong: the stabilizer of a regular $c$ is a two-dimensional torus, which acts by rotations on the three root planes, and the revolutions are displaced, not broken. The comparison with the octonions is dimensional: $O(6)$, of dimension $15$, against $T^2$, of dimension $2$, and $U(2)$, of dimension $4$, on the walls.
\end{remark}

\subsection{The root planes in the pictures}
\begin{figure}[htbp]
\centering
\includegraphics[width=\textwidth]{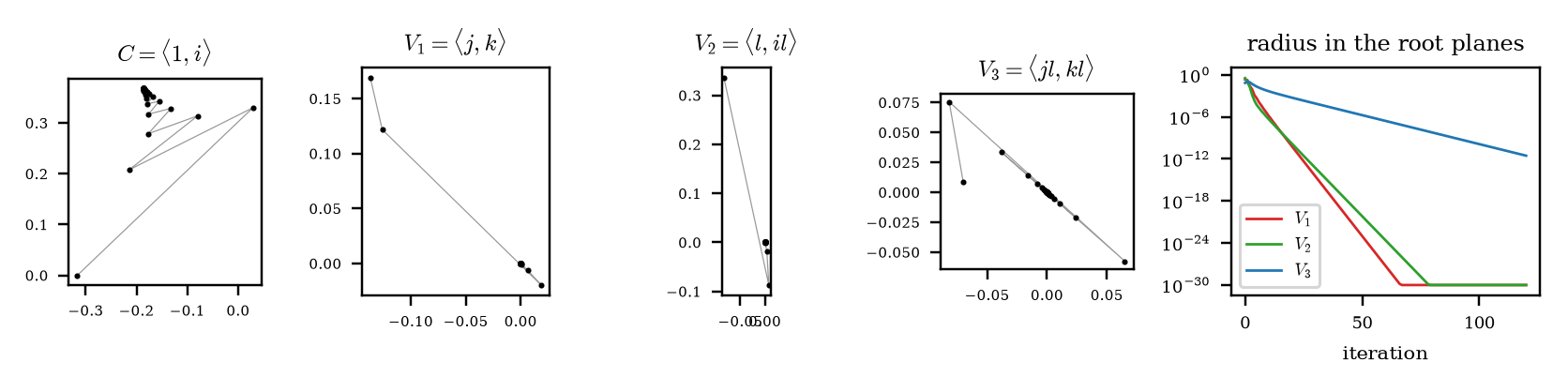}
\caption{\emph{An orbit of the Okubo algebra projected on the Cartan plane and on the three root planes} ($c=c_A$, the first $40$ iterations of a generic starting point), and the radii in the three root planes along $120$ iterations. The radii are invariant under the torus; they decay with three different rates, since the orbit converges to an attracting cycle in $C$ and the three roots see $c$ under different angles.}
\label{fig:orbits}
\end{figure}

The frames $F_A=(1,i,j\,|\,l)$ and $F_C=(1,j,l\,|\,i)$ used in the atlas contain no root plane, and this is why Plates~\plateref{A07} and \plateref{C03} show no revolution. Theorem~\ref{thm:A} predicts that in the frame $(1,l,il\,|\,j)$ the quarter $(1,l,il)$ is a body of revolution about $1$ and the quarter $(l,il,j)$ a body of revolution about $j$, while the quarters $(1,l,j)$ and $(1,il,j)$ are not. Plate~\plateref{R01} was rendered to test this prediction, and it confirms it: the Jaccard indices of the rotations by $37^\circ$ and $90^\circ$ are $1.000$ exactly where predicted. Plates~\plateref{R02} and \plateref{R03} show the other two root planes; the three profiles are different, since the three roots ``see'' the constant $c_A$ under different angles. This can be read off the dynamics as well: along a bounded orbit the radii in $V_1$, $V_2$, $V_3$ decay with three different rates (Figure~\ref{fig:orbits}).

\plate{Atlas_R01_okubo-root-plane-l-il.jpg}{\emph{The displaced revolution: the Okubo algebra, $c=c_A$, frame $(1,l,il\,|\,j)$.} The upper-left quarter contains the root plane $V_2=\langle l,il\rangle$ and is a body of revolution about $1$; the lower-right quarter contains $V_2$ and the axis $j$, fixed by the circle of $T$ that fixes $V_1$ pointwise, and it is a body of revolution about $j$. The two mixed quarters $(1,l,j)$ and $(1,il,j)$ are congruent, being exchanged by a torus element, but they are not bodies of revolution. Compare with the octonionic Plate~\plateref{R04} and with Plate~\plateref{C03}.}{R01}

\begin{figure}[htbp]
\centering
\includegraphics[width=0.49\textwidth]{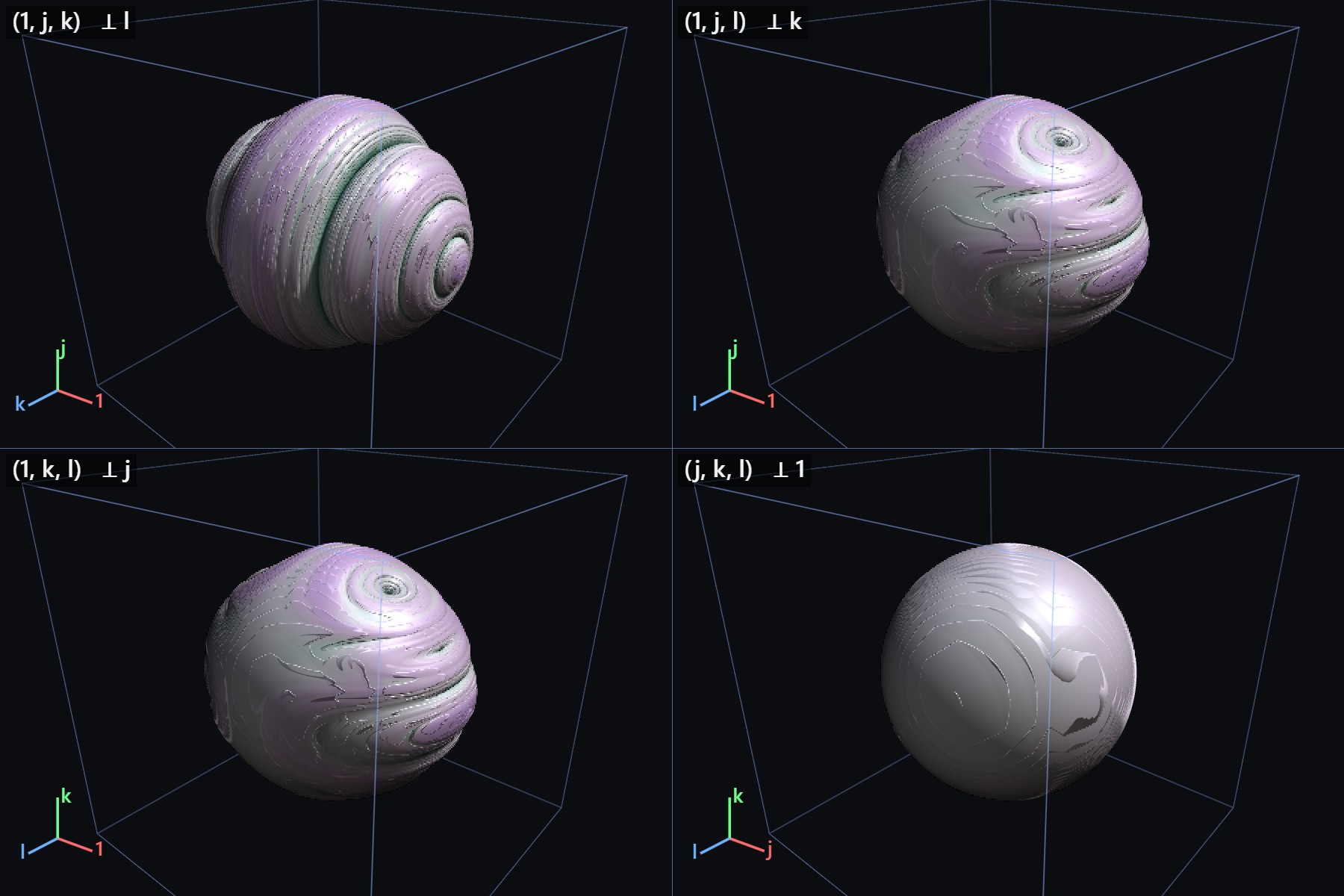}\hfill
\includegraphics[width=0.49\textwidth]{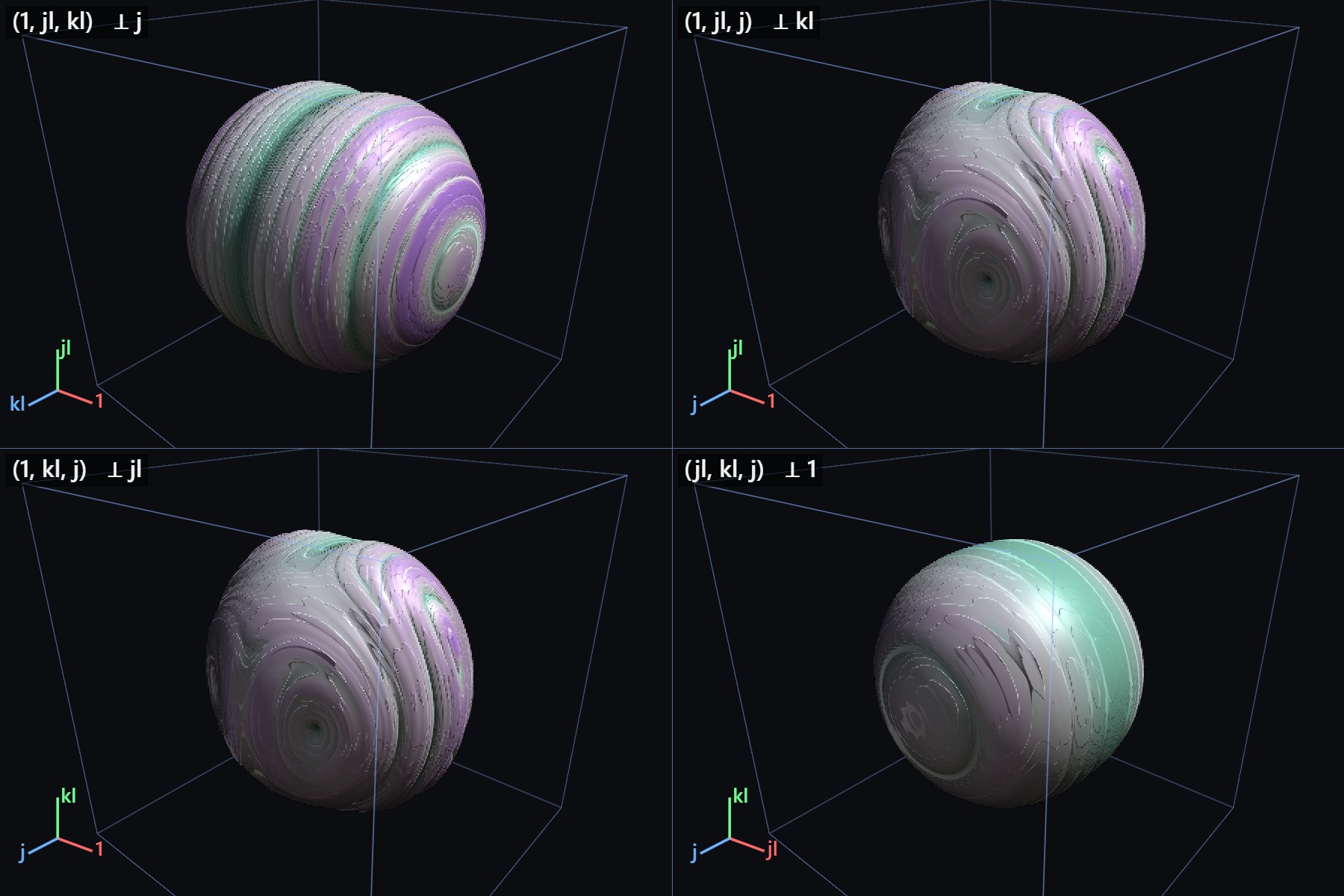}
\caption{\textbf{Plates R02 and R03.} \emph{The other two root planes: the Okubo algebra, $c=c_A$, frames $(1,j,k\,|\,l)$ (left) and $(1,jl,kl\,|\,j)$ (right).} In each plate the upper-left quarter contains a root plane and is a body of revolution about $1$, and the lower-right quarter, which contains the root plane and a fixed axis of a circle of $T$, is a body of revolution about its third axis. The profiles differ from each other and from that of Plate~\plateref{R01}.}
\label{pl:R02}\label{pl:R03}
\end{figure}

\subsection{Weyl walls}
\begin{figure}[htbp]
\centering
\includegraphics[width=\textwidth]{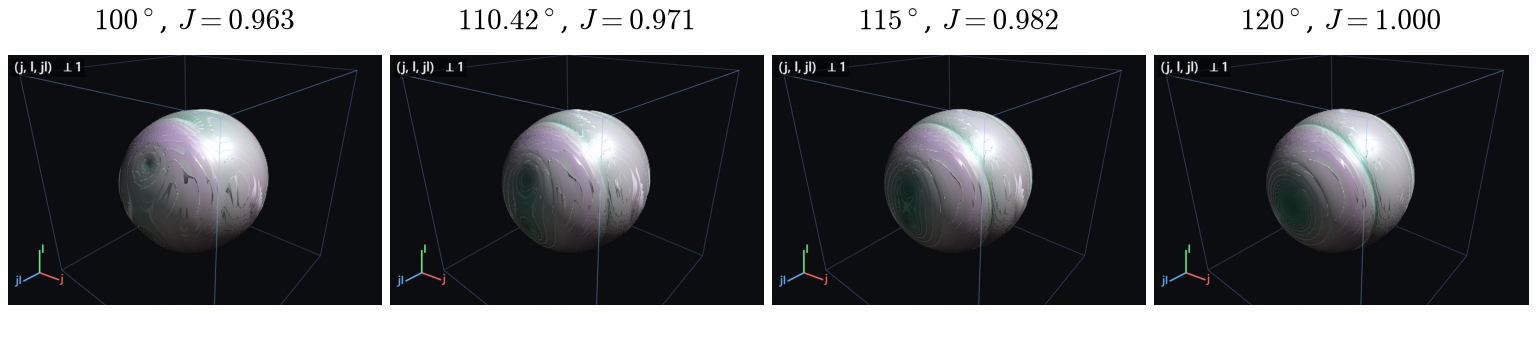}
\caption{\emph{The birth of the $U(2)$ symmetry.} The quarter $(j,l,jl)$ of the Julia set of the Okubo algebra for $|c|=|c_A|$ and $\arg c=100^\circ$, $110.4^\circ$ ($c=c_A$), $115^\circ$ and $120^\circ$ (plates R05--R08), with the Jaccard index $J$ of the rotation by $90^\circ$ about $jl$. The body becomes a body of revolution only on the Weyl wall.}
\label{fig:wall}
\end{figure}

When $c$ approaches a wall, the stabilizer grows from $T^2$ to $U(2)$ and new symmetries appear. For $c$ on $\ell_{120}$ the section $(j,l,jl)$, which mixes three root planes and contains none, becomes a body of revolution about $jl$: the corresponding circle of $U(2)$ rotates the plane $\langle j,l\rangle$ as a complex line of $V_1\oplus V_2\cong\CC^2$ and fixes $jl$. Plates R05--R08 of Appendix~\ref{app:atlas} show this section for $|c|=|c_A|$ and $\arg c=100^\circ$, $110.4^\circ$, $115^\circ$, $120^\circ$, and Figure~\ref{fig:wall} collects the relevant quarters with the Jaccard index of the rotation by $90^\circ$ about $jl$: it grows monotonically from $0.963$ to $0.971$, $0.982$ and finally $1.000$ on the wall. The near-symmetries observed in the first version of the atlas for $c_A$ are then explained by the proximity of $c_A$ to a wall. Part (4) of Theorem~\ref{thm:A} also predicts \emph{tilted} revolutions, about the line $\ell_m$ itself, which no coordinate frame can show; they were confirmed numerically for all three walls and all six choices of $v$, and only for the pairs $(\ell_0,V_1)$, $(\ell_{60},V_2)$, $(\ell_{120},V_3)$.

\subsection{A generic constant}
For a generic $c\in\RR^8$ the stabilizer is again a two-dimensional torus $T_c$, conjugate to $T$, and the statement ``the generic Okubonic Julia set has no continuous symmetry'', suggested by the frame of Plate~\plateref{B07}, is an effect of the frame. The adapted frame is computed from the stabilizer: its common fixed plane contains $c$, and a generic element of its Lie algebra has three invariant planes, the root planes of $c$, on which the torus acts with angular speeds $\omega_1,\omega_2,\omega_3=\omega_1+\omega_2$; for the constant of Plate~\plateref{B07} they are $0.675$, $0.398$ and $0.277$ (in arbitrary normalization). In the frame $(\hat c,v_1,v_2\,|\,c^\perp)$, where $\langle v_1,v_2\rangle$ is a root plane and $c^\perp$ spans the fixed plane together with $\hat c=c/|c|$, Theorem~\ref{thm:A} predicts two exact revolutions, and Plate~\plateref{R09} shows them.

\plate{Atlas_B07_julia-okubo-generic-c.jpg}{\emph{The Okubo algebra with a generic constant, $c\approx0.466\,d$ with $d$ a fixed generic direction, frame $F_B=(1,l,il\,|\,jl)$.} No quarter is a body of revolution: the root planes of the stabilizer of $c$ are in general position with respect to the coordinate axes. The imaginary quarter is nearly, but not exactly, round (Jaccard indices between $0.96$ and $0.99$); like every section through $0$ it is centrally symmetric.}{B07}

\plate{Atlas_R09_okubo-generic-adapted-frame.jpg}{\emph{The same constant in the frame adapted to its stabilizer torus $T_c$.} The upper-left quarter $(\hat c,v_1,v_2)$ contains a root plane of $T_c$ and is a body of revolution about $\hat c$; the lower-right quarter $(v_1,v_2,c^\perp)$ is a body of revolution about $c^\perp$. The mixed quarters are congruent and not of revolution. The revolutions that Plate~\plateref{B07} does not show are there.}{R09}

\subsection{The scan of the coordinate sections}
Theorem~\ref{thm:A} excludes further continuous symmetries of automorphic origin, not every symmetry (Remark~\ref{rem:nonalgebraic}). In order to test the statement against the pictures, we measured all $56$ coordinate sections of $K_{c_A}$ through the origin, for rotations by $37^\circ$ and by $90^\circ$ about each of the three axes: for $c=c_A$ exactly 18 sections have an exact rotation, and they are exactly the 18 sections that contain a root plane; all of the 18 pairs (section, axis) predicted by the torus are exact, and the only further pairs are $(i,j,k)$ about $j$, $(i,j,k)$ about $k$, which come from the identity $z*z=-\n(z)$ on the imaginary quaternions (Lemma~\ref{lem:onestep} with $h=\mathrm{id}$) (Table~\ref{tab:scan}). On the Weyl wall $\arg c=120^\circ$ the pairs of $c_A$ remain exact and 12 new pairs appear, in 8 new sections: $(j,il,kl)$ about $kl$, $(j,k,il)$ about $j$, $(j,k,il)$ about $k$, $(j,k,l)$ about $j$, $(j,k,l)$ about $k$, $(j,l,il)$ about $il$, $(j,l,il)$ about $l$, $(j,l,jl)$ about $jl$, $(k,il,jl)$ about $jl$, $(k,l,il)$ about $il$, $(k,l,il)$ about $l$, $(k,l,kl)$ about $kl$.

\begin{table}[htbp]\centering\footnotesize
\begin{tabular}{@{}lrrrr>{\raggedright\arraybackslash}p{0.3\textwidth}@{}}\toprule
Julia set & sections & pairs & torus & one-step & best non-exact pairs (min.~$J$)\\\midrule
Okubo, $c=c_A$ & 18 & 20 & 18 & 20 (0) & $(k,l,il)$ about $il$ (0.982); $(j,l,il)$ about $il$ (0.982)\\
Okubo, $\arg c=120^\circ$ & 22 & 32 & 18 & -- & $(k,l,jl)$ about $jl$ (0.959); $(k,il,kl)$ about $kl$ (0.959)\\
Petersson $90^\circ$, $c=c_A$ & 42 & 64 & 18 & 64 (0) & $(1,k,l)$ about $k$ (0.910); $(1,k,kl)$ about $k$ (0.910)\\
octonions, $c=c_A$ & 50 & 120 & 18 & -- & $(1,l,kl)$ about $kl$ (0.937); $(1,l,jl)$ about $jl$ (0.937)\\
Petersson $60^\circ$, $c=c_A$ & 18 & 20 & 18 & 20 (0) & --\\
Petersson $150^\circ$, $c=c_A$ & 18 & 20 & 18 & 20 (0) & --\\
\bottomrule\end{tabular}
\caption{\emph{The scan of the $56$ coordinate sections through $0$} ($4\cdot10^4$ points each, $20$ iterations). A pair (section, axis) is exact if the rotations by $37^\circ$ and $90^\circ$ about the axis both have $J\geq0.9995$; the columns give the number of sections with an exact pair, the number of exact pairs, and how many of them are among the 18 pairs predicted by the torus $T$ alone (Theorem~\ref{thm:A}). The column ``one-step'' gives the number of pairs predicted by Lemma~\ref{lem:onestep} with the group $G_c$ and, in brackets, the number of disagreements with the scan. The rows $60^\circ$ and $150^\circ$ were scanned after the prediction was made. The last column shows the two highest non-exact pairs.}
\label{tab:scan}
\end{table}

The same scan also shows that the finite symmetries of the quarters of the atlas are exactly those predicted by the group $(T\rtimes\langle\sigma\rangle)\times\{\pm\mathrm{id}\}$. For each quarter of a plate of the Okubo or Petersson algebras with $c\in C$ and a coordinate frame, we listed the transformations among the $22$ tested ones (three mirrors, the rotations by $180^\circ,120^\circ,90^\circ,72^\circ,60^\circ,37^\circ$ about the three axes, the point reflection) that are realized by an element of the group, and we compared the list with the transformations whose Jaccard index is at least $0.9995$. The two lists coincide in 71 of the 76 quarters (Table~\ref{tab:groups}). No predicted transformation is missing: the only disagreements are exact symmetries beyond the group, and they occur exactly where the theory predicts a larger group: at $\theta=90^\circ$ (Plates~C05, R10, R15; Proposition~\ref{prop:theta90}) and on the Weyl wall (Plate~R08; Theorem~\ref{thm:A}(4)).

\begin{table}[htbp]\centering\footnotesize
\begin{tabular}{@{}llll>{\raggedright}p{0.33\textwidth}>{\raggedright\arraybackslash}p{0.2\textwidth}@{}}\toprule
Plate & $\theta$ & $\arg c$ & agree & exact, not predicted & predicted, not exact\\\midrule
A07 & $120^\circ$ & $110.4^\circ$ & 4/4 & -- & --\\
C03 & $120^\circ$ & $110.4^\circ$ & 4/4 & -- & --\\
C04 & $60^\circ$ & $110.4^\circ$ & 4/4 & -- & --\\
C05 & $90^\circ$ & $110.4^\circ$ & 3/4 & (j,l,i): Ry120, Ry37, Ry60, Ry72, Ry90 & --\\
C06 & $138^\circ$ & $110.4^\circ$ & 4/4 & -- & --\\
R01 & $120^\circ$ & $110.4^\circ$ & 4/4 & -- & --\\
R02 & $120^\circ$ & $110.4^\circ$ & 4/4 & -- & --\\
R03 & $120^\circ$ & $110.4^\circ$ & 4/4 & -- & --\\
R05 & $120^\circ$ & $100.0^\circ$ & 4/4 & -- & --\\
R06 & $120^\circ$ & $110.4^\circ$ & 4/4 & -- & --\\
R07 & $120^\circ$ & $115.0^\circ$ & 4/4 & -- & --\\
R08 & $120^\circ$ & $120.0^\circ$ & 3/4 & (j,l,jl): Rz120, Rz37, Rz60, Rz72, Rz90 & --\\
R10 & $90^\circ$ & $110.4^\circ$ & 2/4 & (1,l,kl): Rx120, Rx37, Rx60, Rx72, Rx90, Ry120, Ry37, Ry60, Ry72, Ry90, Rz120, Rz37, Rz60, Rz72, Rz90; (l,kl,i): Rz120, Rz37, Rz60, Rz72, Rz90 & --\\
R11 & $60^\circ$ & $110.4^\circ$ & 4/4 & -- & --\\
R12 & $120^\circ$ & $110.4^\circ$ & 4/4 & -- & --\\
R13 & $30^\circ$ & $110.4^\circ$ & 4/4 & -- & --\\
R14 & $60^\circ$ & $110.4^\circ$ & 4/4 & -- & --\\
R15 & $90^\circ$ & $110.4^\circ$ & 3/4 & (j,l,i): Ry120, Ry37, Ry60, Ry72, Ry90 & --\\
R16 & $150^\circ$ & $110.4^\circ$ & 4/4 & -- & --\\
\bottomrule\end{tabular}
\caption{\emph{Finite symmetries of the quarters against the group $(T\rtimes\langle\sigma\rangle)\times\{\pm\mathrm{id}\}$} for the plates of the Petersson algebras (compact, $c\in C$, coordinate frame centred at $0$). A transformation (M: mirror, R: rotation about the local axis by the given angle, P: point reflection) is \emph{predicted} if it is the restriction of an element of the group, and \emph{exact} if its Jaccard index is at least $0.9995$ on $10^5$ points; in brackets the measured index of the predicted transformations that are not exact.}
\label{tab:groups}
\end{table}

\section{Mandelbrot sets and invariant theory: Theorem B}\label{sec:mandelbrot}

The Mandelbrot sets with seed $0$ of the para-octonions and of the Okubo algebra are, at first sight, very different: the para-octonionic one is a body of revolution about the real axis, whose sections through $1$ look like tricorns (Plate~\plateref{F02}), while the Okubonic one shows three arms in a section through $1$, $j$ and $l$, a disc-like body in the section through $1$, $i$ and $j$, and a tilted body in the section through $1$, $i$ and $l$ (Plate~\plateref{F03}). Quite surprisingly, they are the same planar figure.

\subsection{The invariant plane}
In a symmetric composition algebra the orbit of the seed $0$ never leaves a plane. Indeed, the linearization of \eqref{eq:symmetric-intro} in $x$ gives
\begin{equation}\label{eq:linearized}
x\diamond(y\diamond w)+w\diamond(y\diamond x)=\n(x,w)\,y ,
\end{equation}
and for $x=y=c$, $w=c\diamond c$ we obtain, together with \eqref{eq:symmetric-intro},
\begin{equation}\label{eq:plane}
c\diamond(c\diamond c)=(c\diamond c)\diamond c=\n(c)\,c,\qquad (c\diamond c)\diamond(c\diamond c)=K\,c-\n(c)\,c\diamond c,\qquad K=\n(c,c\diamond c).
\end{equation}
The plane $P_c=\sspan(c,c\diamond c)$ is therefore a two-dimensional subalgebra, whose structure constants only depend on the two invariants $\n(c)$ and $K$, and every orbit of the seed $0$ lies in it. If $c\diamond c=\lambda c$, then $\n(c)^2=\lambda^2\n(c)$ and $K=2\lambda\n(c)$, so that the \emph{discriminant}
\begin{equation}\label{eq:Delta}
\Delta=4\n(c)^3-K^2
\end{equation}
vanishes; for $\Delta\neq0$ the plane is really two-dimensional. It remains to identify the two-dimensional algebras that occur, and there are exactly two of them.

\subsection{The theorem}
Let $\Tric\subset\CC$ be the tricorn, i.e., the set of $\zeta$ for which the orbit of $0$ under $w\mapsto\bar w^{\,2}+\zeta$ is bounded \cite{CHRS,Milnor92}, and let $\Tric_{\mathrm{hyp}}$ be the set of split-complex numbers $\zeta=x+ye$, $e^2=+1$, for which the orbit of $0$ under the same formula, with the split-complex conjugation $\overline{x+ye}=x-ye$, is bounded. We call $\Tric_{\mathrm{hyp}}$ the \emph{hyperbolic tricorn}; to our knowledge it has not been studied separately in the literature, in spite of the attention received by the hyperbolic Mandelbrot set \cite{Senn,Metzler,Blankers}.

\begin{theorem}\label{thm:B}
Let $(S,\diamond,\n)$ be a real symmetric composition algebra and $c\in S$ with $\Delta\neq0$.
\begin{enumerate}
\item If $\Delta>0$, then $\n(c)>0$, there is $\zeta\in\CC$ with
\begin{equation}\label{eq:zeta}
|\zeta|^2=\n(c),\qquad 2\re\zeta^3=\n(c,c\diamond c),
\end{equation}
unique up to the dihedral group generated by $\zeta\mapsto e^{2\pi i/3}\zeta$ and $\zeta\mapsto\bar\zeta$, and the linear map $c\mapsto\zeta$, $c\diamond c\mapsto\bar\zeta^{\,2}$ is an isometric isomorphism of $P_c$ onto $(\CC,\bar x\bar y)$. In particular $c\in\Mzero$ if and only if $\zeta\in\Tric$.
\item If $\Delta<0$, the same holds with a split-complex $\zeta=x+ye$ satisfying $x^2-y^2=\n(c)$ and $2\re\zeta^3=\n(c,c\diamond c)$, and with $\Tric_{\mathrm{hyp}}$ in place of $\Tric$. In the light-cone coordinates $u=x+y$, $v=x-y$ the hyperbolic tricorn map is $(u,v)\mapsto(v^2+\zeta_u,\,u^2+\zeta_v)$, and $\zeta\in\Tric_{\mathrm{hyp}}$ if and only if $0$ has a bounded orbit under both real quartic maps $t\mapsto(t^2+\zeta_v)^2+\zeta_u$ and $t\mapsto(t^2+\zeta_u)^2+\zeta_v$.
\item If $\n$ is positive definite, then $\Delta\geq0$ for every $c$ and, with a bailout on the norm, the escape time of $c$ equals the Mandelbar escape time of $\zeta$. For the para-octonions and for the compact Okubo algebra,
\begin{equation}\label{eq:saturation}
\Mzero=G\cdot\Tric,\qquad G=G_2,\ \text{resp.}\ G=PSU(3),
\end{equation}
where $\Tric$ is placed in the plane $C$, and the image of the invariant map $c\mapsto(\n(c),K)$ is the region $K^2\leq4\n^3$, whose boundary corresponds to the Weyl walls.
\end{enumerate}
\end{theorem}

\begin{proof}
(1) If $\Delta>0$ then $\n(c)^3>K^2/4\geq0$. Put $r=\n(c)^{1/2}$ and $\zeta=re^{i\phi}$ with $\cos3\phi=K/2r^3\in(-1,1)$. In $(\CC,\bar x\bar y)$ we have $\zeta\diamond\zeta=\bar\zeta^{\,2}$, $\zeta\diamond\bar\zeta^{\,2}=\bar\zeta^{\,2}\diamond\zeta=|\zeta|^2\zeta$ and
\begin{equation}\label{eq:zeta4}
\bar\zeta^{\,2}\diamond\bar\zeta^{\,2}=\zeta^4=(\zeta^3+\bar\zeta^3)\,\zeta-|\zeta|^2\bar\zeta^{\,2}=K\zeta-\n(c)\,\bar\zeta^{\,2},
\end{equation}
which are the structure constants \eqref{eq:plane}; $\zeta$ and $\bar\zeta^{\,2}$ are independent because $\zeta^3\notin\RR$. The map is an isometry, since $\n(c\diamond c)=\n(c)^2=|\zeta|^4$ and $\n(c,c\diamond c)=K=2\re(\zeta\,\zeta^2)$. The admissible $\zeta$ form one orbit of the dihedral group, under which the Mandelbar dynamics is equivariant, since $\overline{\omega w}^{\,2}+\omega\zeta=\omega(\bar w^{\,2}+\zeta)$.
(2) With $y^2=x^2-\n(c)$ the second condition reads $8x^3-6\n(c)x=K$, which always has a root with $x^2\geq\n(c)$: any real root if $\n(c)\leq0$, and $x=\n(c)^{1/2}\cosh t$ with $2\n(c)^{3/2}\cosh3t=|K|$ otherwise. The computation \eqref{eq:zeta4} holds verbatim in the commutative and associative algebra of split-complex numbers. In light-cone coordinates the conjugation exchanges $u$ and $v$ and the square acts coordinatewise, which gives the formula; the even and odd terms of the orbit are orbits of the two quartic maps.
(3) If $\n$ is definite, its restriction to $P_c$ cannot be isometric to the split plane, so $\Delta\geq0$; the isometry gives the equality of the escape times. For the para-octonions, an automorphism in $G_2$ maps $\im c$ to $|\im c|\,i$, hence $c$ into $C$; for the Okubo algebra, every element is conjugate into the Cartan plane $C$ (Theorem~\ref{thm:A}). On $C$ the product is $\bar x\bar y$ by \eqref{eq:onC}, so $\Mzero\cap C=\Tric$, and on $C$ one has $\n=r^2$ and $K=2r^3\cos3\phi$, which gives the image and its boundary $\sin3\phi=0$. \end{proof}

\begin{figure}[htbp]
\centering
\includegraphics[width=0.95\textwidth]{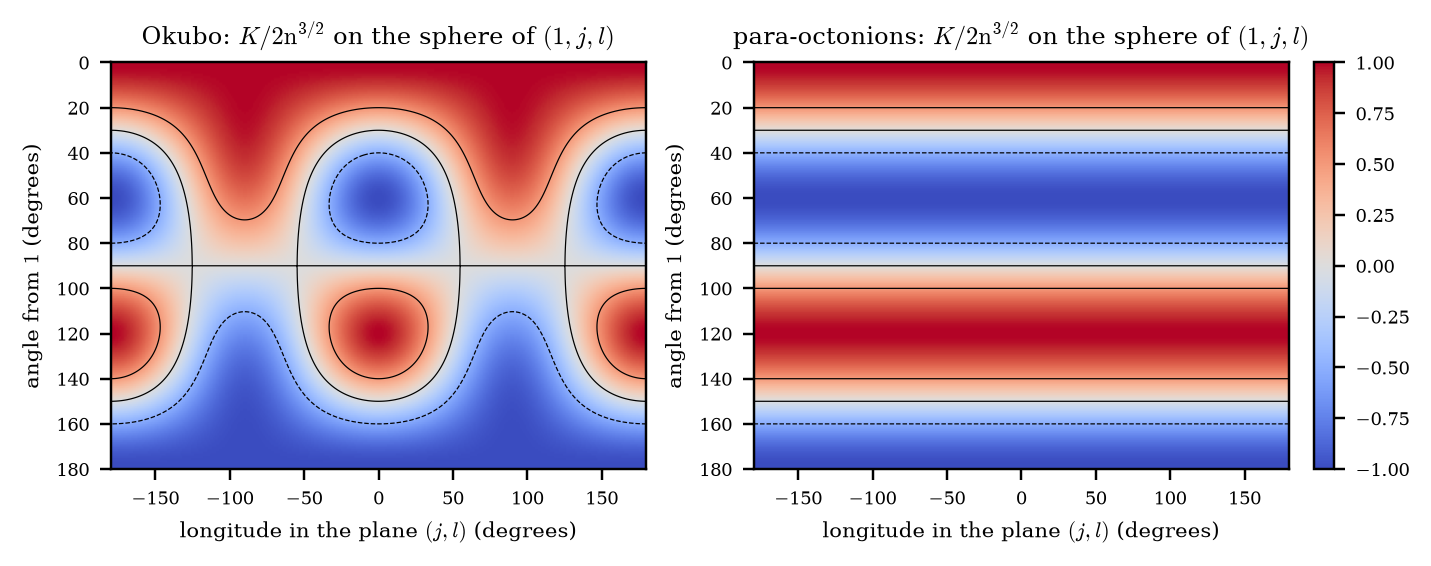}
\caption{\emph{Why three arms.} The normalized cubic invariant $K/2\n^{3/2}=\cos3\phi$ on the unit sphere of the section $(1,j,l)$, in spherical coordinates with the pole on the axis of $1$. In the para-octonions (right) it only depends on the angle from $1$, since $G_2$ acts transitively on the directions orthogonal to $1$; in the Okubo algebra (left) it has three maxima and three minima, which generate the three arms of the lower-left quarter of Plate~\plateref{F03}.}
\label{fig:cubic}
\end{figure}

As an implication of Theorem~\ref{thm:B} we recover, and simplify, the formulas of the first version of the atlas: in the para-octonions $K=2c_0(c_0^2-3|\im c|^2)=2\re\zeta^3$ with $\zeta=c_0+i|\im c|$, and the fact that the Okubonic cubic invariant is not a function of $c_0$ and $|\im c|$ alone is nothing but the fact that the orbits of $PSU(3)$ on $\RR^8$, of dimension at most $6$ and cut out by two invariants, are smaller than the orbits of $G_2$, which are cut out by $c_0$ and $|\im c|$ (Figure~\ref{fig:cubic}).

\begin{corollary}\label{cor:tilted}
For every wall line $\ell_m$ and every $v\in V_{k(m)}$ (Theorem~\ref{thm:A}(4)), the section of the Okubonic Mandelbrot set $\Mzero$ by $\langle1,i,v\rangle$ is a body of revolution about $\ell_m$. In particular the section by $\langle1,i,j\rangle$ is a body of revolution about $1$, and the section by $\langle1,i,l\rangle$ about the line of $e^{i\pi/3}$.
\end{corollary}
\begin{proof}
By \eqref{eq:saturation} $\Mzero$ is invariant under $PSU(3)$, which contains the $SO(3)$ of Theorem~\ref{thm:A}(4). \end{proof}

\plate{Atlas_F02_mandel-para-seed0.jpg}{\emph{The Mandelbrot set with seed $0$ of the para-octonions, centre $-0.35$, frame $(1,i,j\,|\,l)$.} Every quarter through $1$ is a tricorn turned about the real axis by $G_2$ (Theorem~\ref{thm:B}); the imaginary quarter, which only meets the set near the real axis, contains a tiny ball.}{F02}

\plate{Atlas_F03_mandel-okubo-seed0.jpg}{\emph{The Mandelbrot set with seed $0$ of the Okubo algebra, centre $-0.35$, frame $(1,i,j\,|\,l)$.} The upper-left quarter is a body of revolution about $1$, and the upper-right quarter $(1,i,l)$ a body of revolution about the tilted Weyl line of $e^{i\pi/3}$ (Corollary~\ref{cor:tilted}); the lower-left quarter $(1,j,l)$ shows the three arms of the tricorn, redistributed by the flag manifolds of $SU(3)$. By Theorem~\ref{thm:B} every point of every quarter can be decided without iterating in $\RR^8$, by evaluating $(\n(c),K)$ and the planar tricorn.}{F03}

\begin{figure}[htbp]
\centering
\includegraphics[width=\textwidth]{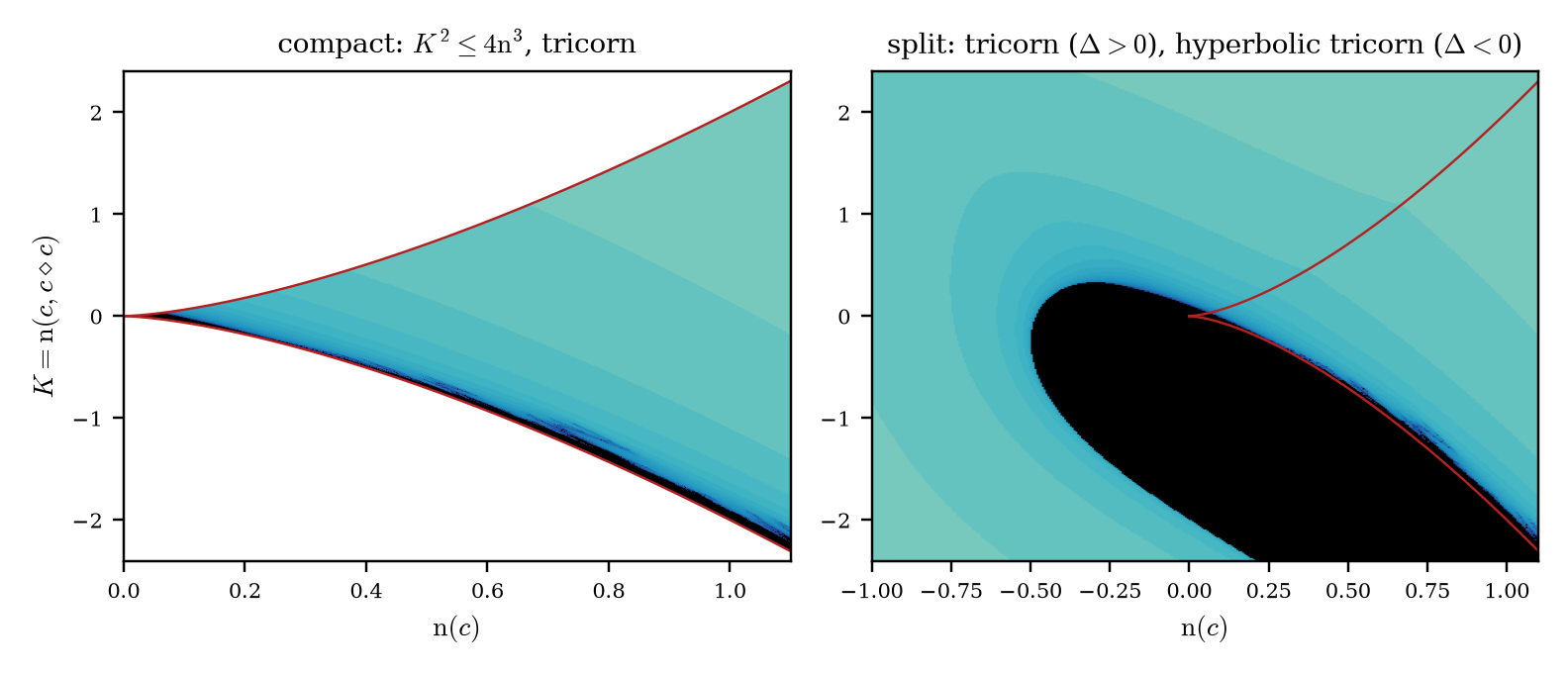}
\caption{\emph{The mother figure of all seed-$0$ Mandelbrot sets of symmetric composition algebras.} The plane of the invariants $(\n(c),K)$, coloured by the escape time of the planar model; bounded orbits in black, the curves $K=\pm2\n^{3/2}$ in red. Left: the compact case, where $(\n,K)$ ranges over $K^2\leq4\n^3$ and the planar model is the tricorn. Right: the split case, where the tricorn ($\Delta>0$) and the hyperbolic tricorn ($\Delta<0$) fill the whole plane. Every three-dimensional section of Plates~\plateref{F02}--\plateref{F06} is a pull-back of these pictures.}
\label{fig:orbitspace}
\end{figure}

\subsection{An external test of the renderer}
\begin{figure}[htbp]
\centering
\includegraphics[width=\textwidth]{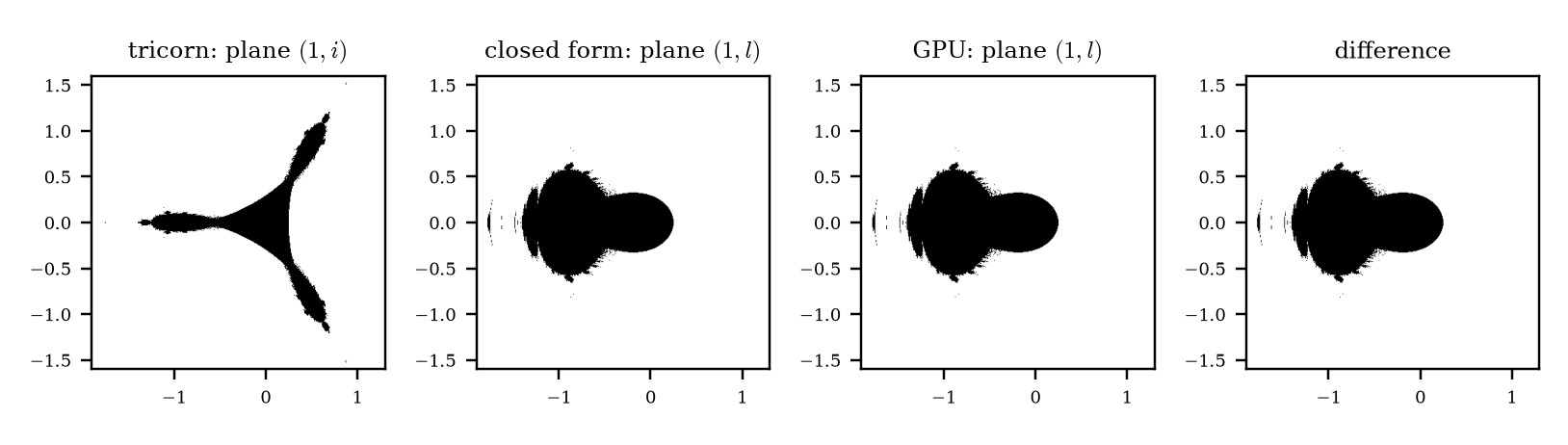}
\caption{\emph{Validation of the renderer against Theorem~\ref{thm:B}.} From left to right: the tricorn, i.e., the plane $(1,i)$ of the Okubonic Mandelbrot set; the plane $(1,l)$ from the closed form; the same plane rendered by the GPU; their difference (black: bounded in both; red or blue: bounded in one only). The planes $(1,l)$ and $(1,jl)$ coincide, as predicted by the $SU(3)$-invariance.}
\label{fig:validation}
\end{figure}

Theorem~\ref{thm:B} gives the Okubonic Mandelbrot set in closed form, and therefore an independent test of the renderer. We rendered nine plane sections of the seed-$0$ Mandelbrot sets of the Okubo algebra, of the para-octonions and of the split Okubo algebra with the GPU of the application (single precision), at $400\times400$ pixels and $100$ iterations, and we evaluated the same pixel centres with the application's double-precision kernel, with an independent \texttt{numpy} iteration in $\RR^8$ and with the closed form. The closed form, the \texttt{numpy} iteration and the application's CPU kernel agree on all $1.44\cdot10^6$ pixels; the GPU agrees with them on at least $99.994\%$ of the pixels of every section, the differences being a handful of boundary pixels per section, due to single-precision rounding (Figure~\ref{fig:validation}).

\subsection{Split algebras}
For the split Okubo algebra and the split para-octonions the norm is indefinite, the discriminant $\Delta$ takes both signs, and Theorem~\ref{thm:B} gives a description of $\Mzero$ as the union of two pull-backs: of the tricorn over $\Delta>0$ and of the hyperbolic tricorn over $\Delta<0$, up to the closed set $\Delta=0$, which has measure zero. This settles a conjecture formulated by the referee of the first version of this work. Two warnings are nevertheless in order, and they concern the pictures, not the theorem.

\begin{remark}\label{rem:bailout}
In a split algebra the Euclidean bailout of the renderer is not an escape certificate, since orbits may remain bounded while their Euclidean norm is large. Theorem~\ref{thm:B} is a statement on bounded orbits, which is independent of the choice of the norm, while the rendered pictures depend on the bailout.
\end{remark}

\begin{figure}[htbp]
\centering
\includegraphics[width=\textwidth]{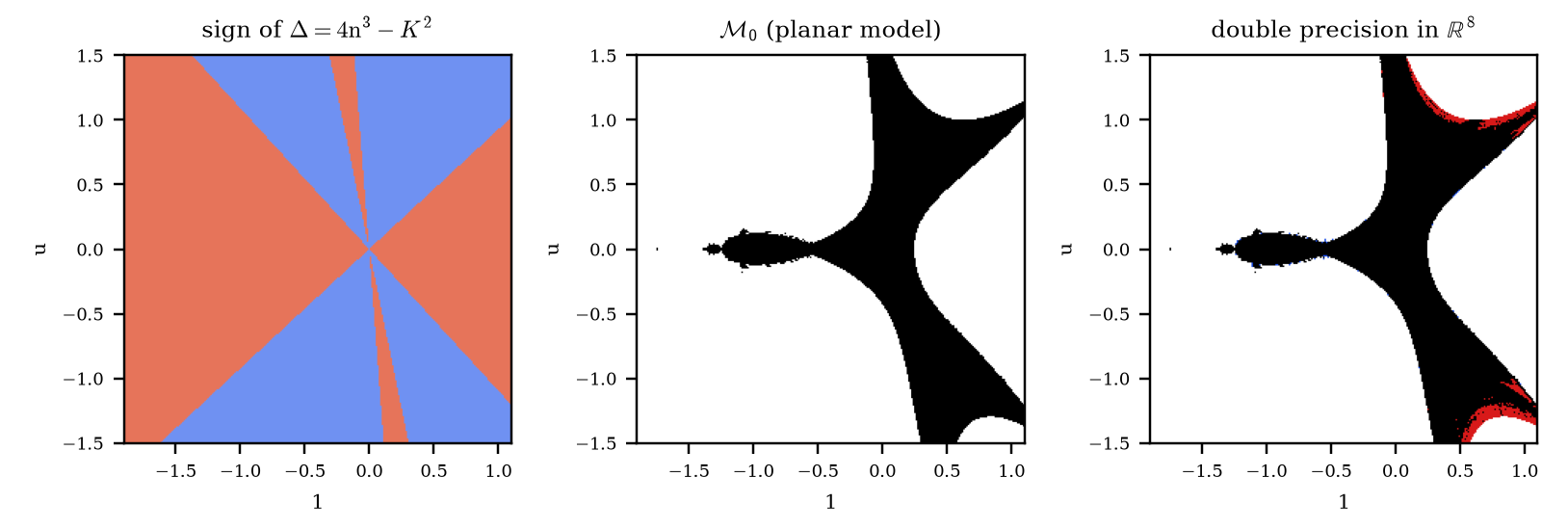}
\caption{\emph{The split Okubo algebra in a generic plane $(1,u)$.} Left: the sign of $\Delta$ (red: $\Delta>0$, tricorn; blue: $\Delta<0$, hyperbolic tricorn). Centre: $\Mzero$ from the planar model (Theorem~\ref{thm:B}). Right: the renderer's rule in double precision in $\RR^8$ ($200$ iterations, bailout $100$); in red the points of $\Mzero$ lost by transverse instability, in blue the points counted as bounded because their escape needs more than $200$ iterations.}
\label{fig:split}
\end{figure}

\begin{remark}\label{rem:transverse}
In the split Okubo algebra the invariant plane $P_c$ can be \emph{transversally unstable}: a double-precision iteration in $\RR^8$ leaves $P_c$ by rounding, the distance from the plane grows exponentially, from $10^{-16}$ to $10^{-1}$ in about $180$ iterations, and the orbit escapes although the exact orbit is bounded. Among $1500$ random constants, $18$ escape times computed in $\RR^8$ differ from those of the planar model, and all $18$ agree with the planar model when the iteration in $\RR^8$ is repeated with $80$ significant digits. In a generic plane $(1,u)$ the double-precision iteration loses 902 of the 14731 points of $\Mzero$ found by the planar model (Figure~\ref{fig:split}); no such loss occurs for the compact algebras. The rendered split Mandelbrot sets (Plate~\plateref{F06}) are therefore approximations of $\Mzero$ from inside.
\end{remark}

\plate{Atlas_F06_mandel-sokubo-seed0.jpg}{\emph{The Mandelbrot set with seed $0$ of the split Okubo algebra, centre $-0.35$, frame $(1,i,j\,|\,l)$.} The quarters inside the quaternions coincide with those of the split para-octonions; the quarters through $l$ cross the surfaces $\Delta=0$, where the planar model changes from the tricorn to the hyperbolic tricorn. The rendered body is an approximation from inside of $\Mzero$ (Remarks~\ref{rem:bailout} and \ref{rem:transverse}).}{F06}

\section{Identities of symmetric composition algebras}\label{sec:identities}

The identities \eqref{eq:symmetric-intro}, harmless in appearance, have drastic dynamical consequences, which are visible at a glance in the atlas. The first two were already noticed in the first version of this work.

\begin{proposition}\label{prop:sandwich}
In a symmetric composition algebra the sandwich map satisfies $(z\diamond a)\diamond z+c=\n(z)\,a+c$ and the cubic map $(z\diamond z)\diamond z+c=\n(z)\,z+c$. The norms $t_n=\n(z_n)$ of a sandwich orbit obey the real recursion
\begin{equation}\label{eq:sandwich}
t_{n+1}=\n(a)\,t_n^2+\n(a,c)\,t_n+\n(c),
\end{equation}
which is affinely conjugate to $s\mapsto s^2+c_{\mathrm{eff}}$ with $c_{\mathrm{eff}}=\n(a)\n(c)+\tfrac12\n(a,c)-\tfrac14\n(a,c)^2$; in the compact case the filled set is a ball, or a spherical shell.
\end{proposition}
\begin{proof}
The first two identities are \eqref{eq:symmetric-intro}. Taking norms of $z_{n+1}=t_na+c$ gives \eqref{eq:sandwich}, and the change of variable $s=\n(a)t+\tfrac12\n(a,c)$ gives the conjugation. \end{proof}

The sandwich balls of the para-octonions and of the Okubo algebra are therefore identical, and their radius is predicted by \eqref{eq:sandwich}: for the constants of Plate~I02 of Appendix~\ref{app:atlas} the prediction is $1.270$, and the body share measured on $10^5$ points, $58.4\%$ of the ball of radius $1.52$, gives the radius $1.270$. The bilateral map $(a\diamond z)\diamond(z\diamond b)+c$ looks more promising, since it does not degenerate into a function of the norm. Quite unexpectedly, however, it degenerates into the octonions.

\begin{proposition}\label{prop:bilateral}
Let $(S,\diamond,\n)$ be a symmetric composition algebra, $a,b\in S$ and $e=b\diamond a$. Then
\begin{equation}\label{eq:bilateral}
(a\diamond z)\diamond(z\diamond b)=\n(z,e)\,z-\n(z)\,e\qquad\text{for all }z\in S .
\end{equation}
The product $x\star y=(a\diamond x)\diamond(y\diamond b)$ satisfies $e\star y=y\star e=\n(a)\n(b)\,y$, so that for $\n(a)\n(b)=1$ it is a Hurwitz algebra with unit $e$. If $\n$ is positive definite and $Q$ is an orthogonal map with $Qe=|e|\,1$, then
\begin{equation}\label{eq:bilateral-conj}
|e|\,Q\big((a\diamond z)\diamond(z\diamond b)+c\big)=v^2+|e|\,Qc,\qquad v=|e|\,Qz,
\end{equation}
where $v^2$ is the octonionic square. In particular every bilateral filled Julia set over the para-octonions or the Okubo algebra is similar to an octonionic filled Julia set, every orbit lies in $\sspan(z_0,e,c)$, and Kaplansky's trick $(1\diamond x)\diamond(y\diamond1)=xy$ is the special case $a=b=1$.
\end{proposition}
\begin{proof}
The linearization of $(x\diamond y)\diamond x=\n(x)y$ reads $(x\diamond y)\diamond w+(w\diamond y)\diamond x=\n(x,w)\,y$. With $x=a$, $y=z$, $w=z\diamond b$, and since $(z\diamond b)\diamond z=\n(z)b$, we obtain $(a\diamond z)\diamond(z\diamond b)=\n(a,z\diamond b)\,z-\n(z)\,b\diamond a$; the associativity of the norm, $\n(x\diamond y,w)=\n(x,y\diamond w)$ \cite{KMRT,Elduque97}, gives $\n(a,z\diamond b)=\n(z,b\diamond a)$. Next, $a\diamond e=a\diamond(b\diamond a)=\n(a)b$, so $e\star y=\n(a)\,b\diamond(y\diamond b)=\n(a)\n(b)y$, and similarly on the right. Finally, with $v=|e|Qz$ one has $\n(z,e)=2|e|\re(Qz)$ and $\n(z)=\n(Qz)$, and \eqref{eq:bilateral-conj} follows from $v^2=2\re(v)\,v-\n(v)$. For $a=b=1$ one has $e=1$ both in the para-octonions and in the Okubo algebra, and \eqref{eq:bilateral} gives $2\re(z)z-\n(z)=z^2$. \end{proof}

It is mesmerizing how this simple identity accounts for the whole Part~J of the atlas: the orbit span of three dimensions measured for the bilateral maps of the para-octonions and of the Okubo algebra is exact, not an artefact of converging orbits, and the spiral ribs of Plate~\plateref{J06} are the ribs of an octonionic body of revolution, whose axis is tilted with respect to the coordinate frame. The escape times of the bilateral maps of Plates~J04 and J06, normalized to $\n(a)=\n(b)=1$, coincide with those of the octonionic map after the isometry $Q$ on $2\cdot10^4$ out of $2\cdot10^4$ points. The identity fails in the octonions and at the Petersson angle $90^\circ$, which are not symmetric.

\plate{Atlas_J06_bilateral-okubo.jpg}{\emph{The bilateral map $(a*z)*(z*b)+c$ over the Okubo algebra, generic $a,b,c$, frame $F_B=(1,l,il\,|\,jl)$.} By Proposition~\ref{prop:bilateral} the filled set is an octonionic filled Julia set moved by an isometry: every quarter is a section of a body of revolution about a tilted axis, which is why the ribs spiral about directions that are not coordinate axes.}{J06}

\begin{figure}[htbp]
\centering
\includegraphics[width=0.85\textwidth]{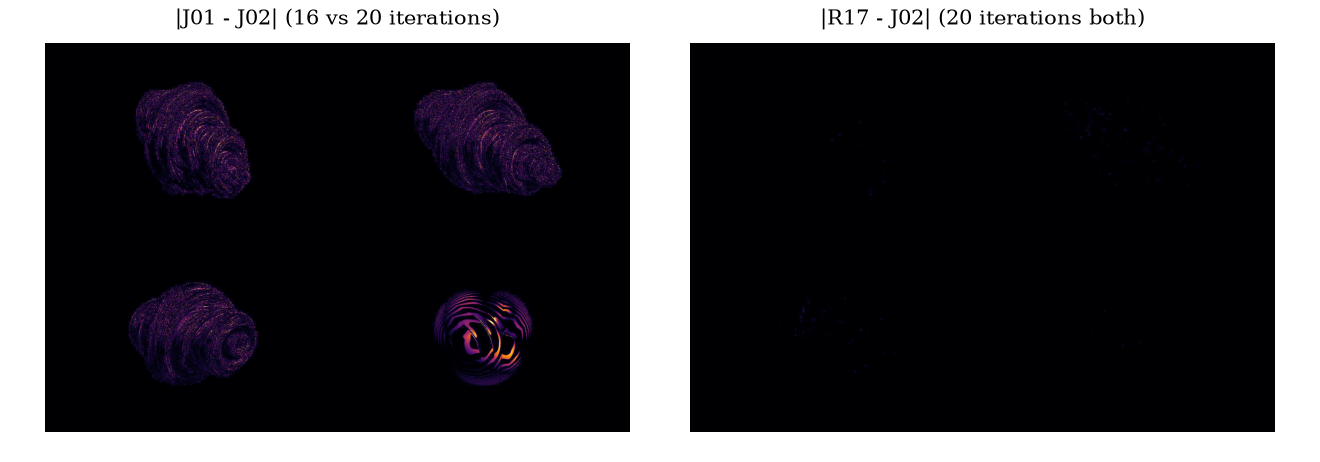}
\caption{\emph{Kaplansky's trick as a null test.} Left: the difference between the Okubonic bilateral map with $a=b=1$ at $16$ iterations (Plate~J01) and the octonionic Julia set at $20$ iterations (Plate~J02). Right: the difference between the bilateral map at $20$ iterations (Plate~R17) and Plate~J02, which is empty up to JPEG noise.}
\label{fig:kaplansky}
\end{figure}

Kaplansky's trick gives the cleanest null test of the whole atlas. In the first version the bilateral map with $a=b=1$ over the Okubo algebra (Plate~J01) and the octonionic Julia set with the same constant (Plate~J02) were rendered with $16$ and $20$ iterations, and their body shares differed by up to one percentage point on the $2400$ points of the first version ($0.3$--$0.7$ points on $10^5$ points); the difference was due to the number of iterations only. With $20$ iterations for both (Plate~R17), the two masks coincide on every one of the $4\cdot10^5$ sampled points, while J01 and J02 differ on about $1000$ of them, and the two images differ by a mean of $0.013$ grey levels, i.e., by JPEG noise (Figure~\ref{fig:kaplansky}).

\section{The Petersson path and the angle \texorpdfstring{$90^\circ$}{90\textdegree{}}}\label{sec:petersson}

Along the Petersson path \eqref{eq:petersson} from the para-octonions ($\theta=0$) to the Okubo algebra ($\theta=120^\circ$) the algebra is a composition algebra but not a symmetric one, and Theorem~\ref{thm:B} does not apply. The torus $T$ of Proposition~\ref{prop:torus} survives for every $\theta$, as does the anti-automorphism $\sigma$, so that the finite symmetry groups of the frames $F_A$ and $F_C$ are the same for every angle: this is why the plates of the path in these frames look alike. A single angle is exceptional. In the first version of this work we observed that for $\theta=90^\circ$ the square of $x_0+ql$, with $x_0\in\RR$ and $q\in\HH$, is $x_0^2+|q|^2$, so that every section contained in $\RR\oplus\HH l$ is the same ball (Plate~\plateref{E04}); and we observed an unexplained exact rotational symmetry in the quarter $(j,l,i)$ of Plate~\plateref{C05}. The scan of the $56$ coordinate sections (Table~\ref{tab:scan}) shows that these are only the tip of the iceberg: at $\theta=90^\circ$, with the same constant $c_A$, $42$ sections have exact rotations, against $18$ for the Okubo algebra. All of them are explained by the following.

\begin{proposition}\label{prop:theta90}
Let $\theta=90^\circ$. For $z=p+ql$ with $p,q\in\HH$,
\begin{equation}\label{eq:square90}
z*z=\bar p^{\,2}+|q|^2-2\,\big(i\,q\,\im p\big)\,l .
\end{equation}
Consequently, for every unit quaternion $u$ the isometry $g_u(p+ql)=p+(uq)l$ satisfies
\begin{equation}\label{eq:twisted}
Q(g_uz)=g_{iui^{-1}}\,Q(z),\qquad z\in\RR^8 ,
\end{equation}
and for $c\in C$ the map $f(z)=z*z+c$ satisfies $f\circ g_u=g_{iui^{-1}}\circ f$. For $u\in C$ the map $g_u=\varphi_{1,u}$ is an automorphism of the torus $T$; for $u=e^{j\varphi}$ or $u=e^{k\varphi}$ it is twisted, $f\circ g_u=g_u^{-1}\circ f$. Hence, for $c\in C$, the group $G_c$ of Section~\ref{sec:framework} contains the four-dimensional group generated by $T$ and $\{g_u\}$, and $K_c$ is invariant under it, with equal escape times. In particular:
\begin{enumerate}
\item every section of $K_c$ contained in $\RR\oplus\HH l$ is the same ball, since there $z*z=x_0^2+|q|^2$;
\item the subspace $W=\langle1,i\rangle\oplus\langle l,kl\rangle$ is closed under squaring, $z*z=\bar p^{\,2}+|q|^2+2\im(p)\,\bar q\,l$ for $p\in\langle1,i\rangle$, $q\in\langle1,k\rangle$, and $K_c\cap W$ is invariant under the rotations $R_\varphi=g_{e^{k\varphi}}$ of the plane $\langle l,kl\rangle$;
\item the section $\langle i,j,l\rangle$ of $K_c$ is a body of revolution about $l$, since for $z=x_1i+x_2j+x_4l$ one has $z*z=(x_4^2-x_1^2-x_2^2)+2x_4\,(x_1\,l-x_2\,kl)$, and a rotation of $(x_1,x_2)$ by $t$ rotates the square by $g_{e^{-kt}}$ (Lemma~\ref{lem:onestep}).
\end{enumerate}
\end{proposition}
\begin{proof}
With $u_\theta=i$ we have $\tau(\bar z)=\bar p-(iq)l$ and $\tau^{-1}(\bar z)=\bar p+(iq)l$, and \eqref{eq:cd} gives $\bar p^{\,2}-\overline{iq}\,(-iq)=\bar p^{\,2}+|q|^2$ for the first component and $iq\,\bar p-iq\,p=-2\,iq\im p$ for the second one, which is \eqref{eq:square90}. Replacing $q$ by $uq$ leaves $|q|^2$ unchanged and multiplies $-2\,iq\im p$ on the left by $iui^{-1}$, which is \eqref{eq:twisted}; the maps $g_u$ fix $C$ pointwise. For $u\in C$, $iui^{-1}=u$; for $u=e^{j\varphi}$ or $e^{k\varphi}$, $iui^{-1}=\bar u$ and $g_{\bar u}=g_u^{-1}$. The group generated by the $g_u$ and by $T=\{\varphi_{a,b}\}$ has dimension $4$: the $g_u$ form a copy of $Sp(1)$, which meets $T$ in the circle $u\in C$. For (1), $p=x_0$ is real and $\im p=0$. For (2), $iqi=-\bar q$ for $q\in\langle1,k\rangle$. For (3), apply Lemma~\ref{lem:onestep} with $h=g_{e^{-kt}}$. \end{proof}

The Lie algebra of this group consists of two skew derivations of $\circ$, the Lie algebra of $T$, and two skew \emph{anti-derivations}, the generators of $g_{e^{j\varphi}}$ and $g_{e^{k\varphi}}$; the stabilizer of $c$ in $\Aut(\circ)$ is still the two-dimensional torus, so that half of the continuous symmetry at $90^\circ$ does not come from automorphisms. The one-step criterion of Lemma~\ref{lem:onestep}, applied with the group $G_c$, predicts which pairs (section, axis) of the $56$ coordinate sections have exact rotations: 64 pairs at $90^\circ$, 20 pairs at $120^\circ$, 20 pairs at $60^\circ$, 20 pairs at $150^\circ$. The scans find exactly these pairs, with no exception, including the scans at $60^\circ$ and $150^\circ$, which were made after the prediction (Table~\ref{tab:scan}). Along the path, among the angles $0^\circ$, $30^\circ$, $60^\circ$, $90^\circ$, $120^\circ$, $150^\circ$, skew anti-derivations of $\circ$ exist only at $60^\circ$ and $90^\circ$, and only at $90^\circ$ do they vanish on $c_A$: for constants in $C$, the angle $90^\circ$ is the only exceptional angle among those tested.

\begin{figure}[htbp]
\centering
\includegraphics[width=\textwidth]{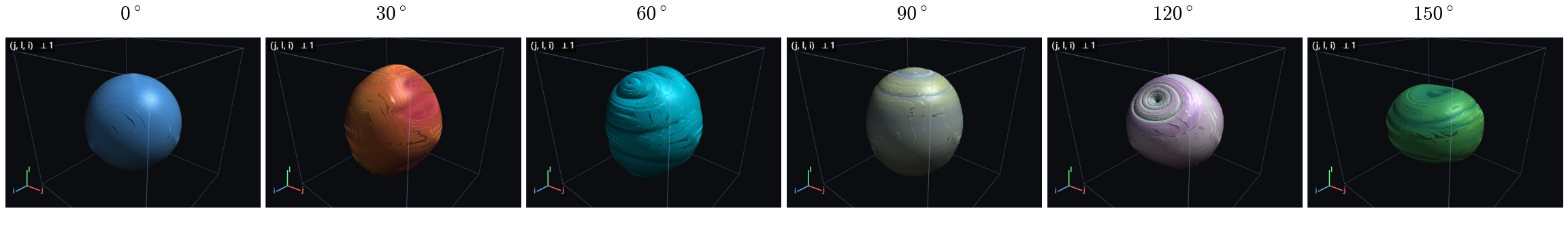}
\caption{\emph{The Petersson path at fixed $c=c_A$:} the quarter $(j,l,i)$ of the frame $F_C$ for $\theta=0^\circ$ (para-octonions, Plate~C02), $30^\circ$, $60^\circ$, $90^\circ$ (Plates R13--R15), $120^\circ$ (Okubo, Plate~C03) and $150^\circ$ (Plate R16). Only $\theta=0^\circ$ and $\theta=90^\circ$ give bodies of revolution about $l$.}
\label{fig:theta}
\end{figure}

Proposition~\ref{prop:theta90} predicts new exact revolutions, e.g., in the quarters $(1,l,kl)$, which is even a ball by the observation above, and $(l,kl,i)$ of the frame $(1,l,kl\,|\,i)$, and these were confirmed by Plate~\plateref{R10}; in the same frame, at $\theta=60^\circ$ and at $\theta=120^\circ$ (Plates R11 and R12 of Appendix~\ref{app:atlas}), the revolutions disappear. Figure~\ref{fig:theta} shows the section $(j,l,i)$ along the path at fixed $c=c_A$: bodies of revolution about $l$ at $0^\circ$, where the para-octonionic group $O(6)$ acts, and at $90^\circ$, and not elsewhere. The angle $90^\circ$ is therefore an exception, not a stage of a continuous transition.

\plate{Atlas_C05_sym-petersson-90.jpg}{\emph{The Petersson algebra at $\theta=90^\circ$, $c=c_A$, frame $F_C=(1,j,l\,|\,i)$.} The lower-right quarter $(j,l,i)$ is an exact body of revolution about $l$ (Proposition~\ref{prop:theta90}), which no automorphism explains; the other quarters have the finite symmetries of the torus.}{C05}

\plate{Atlas_R10_petersson-90-W.jpg}{\emph{The twisted revolutions of $W$: the Petersson algebra at $\theta=90^\circ$, $c=c_A$, frame $(1,l,kl\,|\,i)$.} The upper-left quarter $(1,l,kl)$ is a ball, since $z*z=x_0^2+|q|^2$ there; the lower-right quarter $(l,kl,i)$ lies in $W$ and is a body of revolution about $i$ by \eqref{eq:twisted}; the two mixed quarters $(1,l,i)$ and $(1,kl,i)$ are exchanged by $R_{90^\circ}$ and have the same shape.}{R10}

\plate{Atlas_E04_petersson-theta-90.jpg}{\emph{The Petersson algebra at $\theta=90^\circ$ with a generic constant, frame $F_B=(1,l,il\,|\,jl)$.} All four quarters lie in $\RR\oplus\HH l$, where $z*z=x_0^2+|q|^2$: the first iterate depends only on the norm and every quarter is the same ball.}{E04}

\section{Numerical methods and validation}\label{sec:numerics}

\subsection{The renderer}
All plates were rendered with HyperFractal~8D~1.0 \cite{HF8D}, an eight-dimensional explorer compatible with the Ultra Fractal formula language, in which the complex variables of a formula are read as elements of the chosen algebra. A frame is an orthonormal basis $E_0,\dots,E_7$ of $\RR^8$; the first three vectors are the visible axes $X,Y,Z$ and $E_3$ is the hidden axis $N$. In the \emph{four slices} mode the image is divided into quarters, which show the sections $(X,Y,Z)\perp N$, $(X,Y,N)\perp Z$, $(X,Z,N)\perp Y$ and $(Y,Z,N)\perp X$ through the centre of the box; every quarter carries its labels and a tripod of the local axes. The body is the set of points whose orbit does not leave the ball of squared radius $100$ within the number of iterations of the plate, and its surface is found by sphere tracing on the GPU in single precision. The frames used most often are $F_A=(1,i,j\,|\,l)$, $F_B=(1,l,il\,|\,jl)$ and $F_C=(1,j,l\,|\,i)$; the version used for the new plates accepts arbitrary orthonormal frames, such as the frame adapted to the stabilizer of a generic constant (Plate~\plateref{R09}).

\subsection{Measurements}
Every quarter of every plate without four-dimensional mappings was measured in double precision with an implementation of the algebras independent of the renderer, on $10^5$ uniform points of the ball of radius $0.95h$, where $h$ is the half size of the box. For each quarter we computed:
\begin{itemize}
\item the \emph{body share}, i.e., the fraction of the ball occupied by the body;
\item the \emph{Jaccard index} $J=|A\cap T(A)|/|A\cup T(A)|$ of the body $A$ and its image under $22$ transformations $T$ of the local coordinates (three mirrors, rotations by $180^\circ,120^\circ,90^\circ,72^\circ,60^\circ,37^\circ$ about the three axes, and the point reflection), with a $95\%$ Poisson-bootstrap interval;
\item the dimension of the subalgebra generated by a generic point of the section and the constants of the map, i.e., the $\circ$-subalgebra for the quadratic maps and the $\diamond$-subalgebra otherwise, which contains every orbit.
\end{itemize}
The key observation is that an \emph{exact} symmetry of algebraic origin, being an isometry that maps orbits to orbits, gives $J=1$ up to the handful of points whose escape is decided by rounding at the bailout, while an \emph{approximate} symmetry gives a value that stays below $1$ as the sample grows. We therefore classify a transformation as exact if $J\geq0.9995$ and as approximate if $0.97\leq J<0.9995$. On the 121 measured plates (484 quarters, $22$ transformations each) we find 2435 exact and 280 approximate symmetries. The exact indices cluster at $1$: all but 4 of them are at least $0.9999$. The largest approximate index, $0.99879$ (Plate~D08, the degenerate dual plane), has the $95\%$ interval $[0.99852,0.99906]$, which excludes the exact class; only 20 of the 280 approximate indices have an interval reaching $0.9995$; the median width of these intervals is $0.0026$. The complete tables, with the intervals, are in Appendix~\ref{app:atlas}.

\begin{figure}[htbp]
\centering
\includegraphics[width=0.85\textwidth]{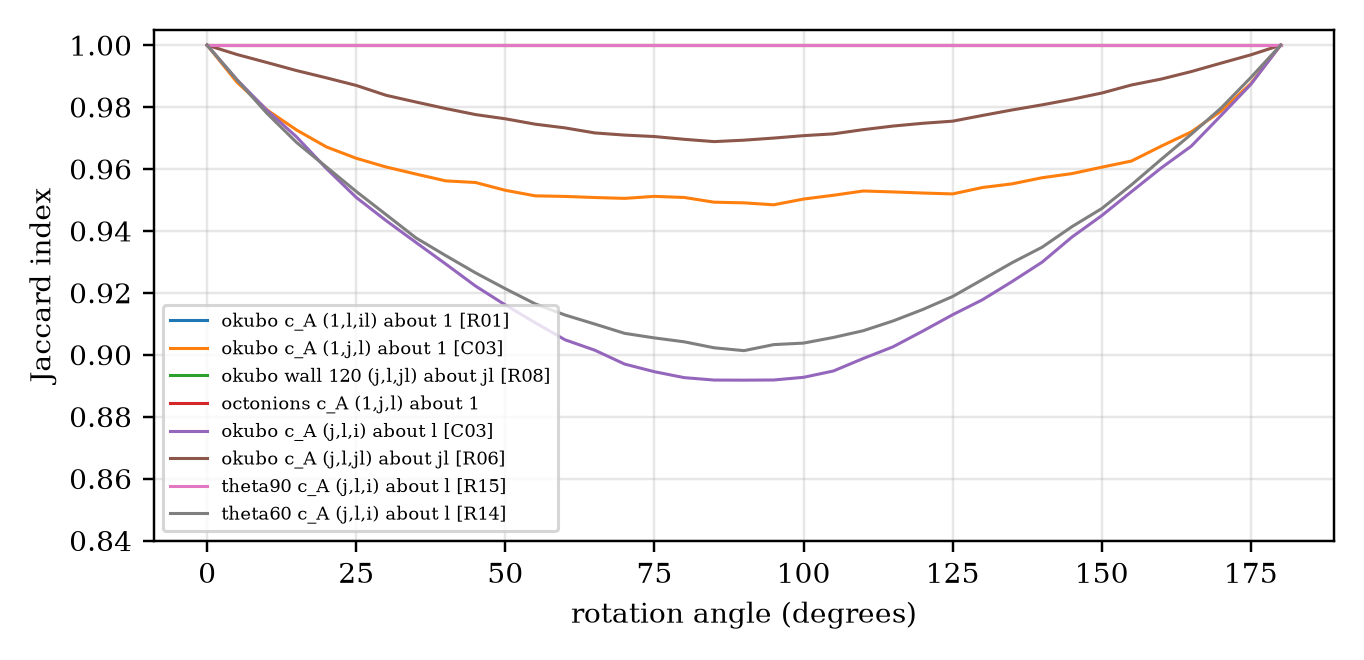}
\caption{\emph{Exact and approximate symmetries.} The Jaccard index of the rotation by an angle $0^\circ$--$180^\circ$ about an axis of a section, on $5\cdot10^4$ points. Exact symmetries (the root plane $(1,l,il)$ of the Okubo algebra, the octonionic $(1,j,l)$, the section $(j,l,i)$ at $\theta=90^\circ$, the section $(j,l,jl)$ on the Weyl wall) are flat at $1$; the others bend down.}
\label{fig:curves}
\end{figure}

The first version of the atlas used a threshold of $0.97$ on about $2400$ points, which cannot separate the two cases: for $c_A$ the rotation about $jl$ in the section $(j,l,jl)$ scores above $0.97$ but is not a symmetry (Figure~\ref{fig:wall}), and the rotations about $1$ in the quarter $(1,j,l)$ of Plate~\plateref{C03} scored above $0.97$ at $37^\circ$ only, while on $10^5$ points they all stay between $0.95$ and $0.96$. Figure~\ref{fig:curves} shows the Jaccard index as a function of the rotation angle for a few sections: an exact symmetry is a flat line at $1$, an approximate one a curve that bends down.

\begin{figure}[htbp]
\centering
\includegraphics[width=\textwidth]{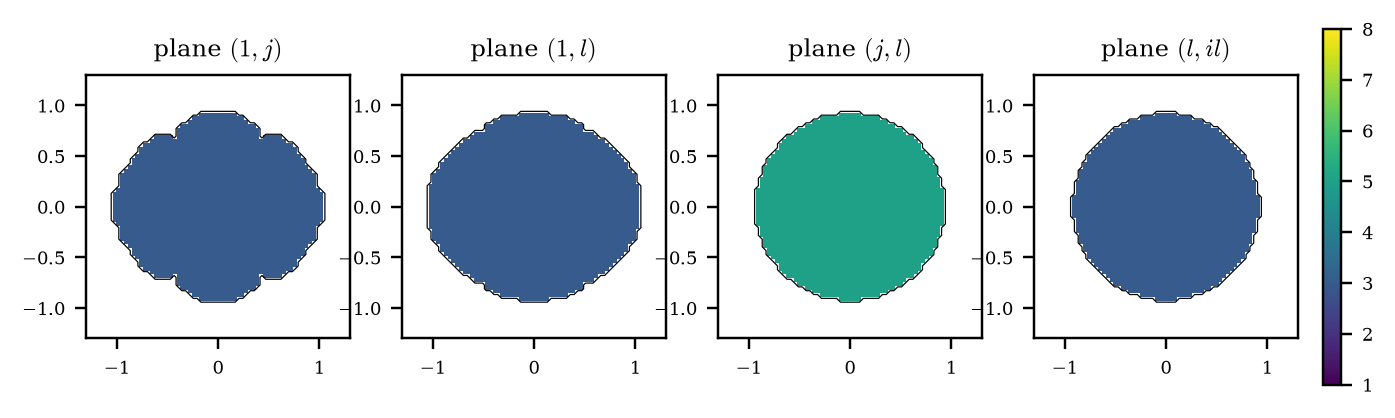}
\caption{\emph{The dimension of the $\circ$-subalgebra generated by $z_0$ and $c=c_A$} in the Okubo algebra, for $z_0$ in four coordinate planes, shown over the filled Julia set (outline in black). The dimension is constant on each plane.}
\label{fig:dimension}
\end{figure}

The \emph{orbit span} of the first version, i.e., the rank of the first $24$ points of six bounded orbits, was replaced by the dimension of the generated subalgebra, which is an exact invariant and does not collapse when an orbit converges to an attracting cycle. Figure~\ref{fig:dimension} shows it on four coordinate planes of the Okubo algebra with $c=c_A$: it is constant on each plane, equal to $3$ on the planes $(1,j)$, $(1,l)$ and $(l,il)$ and to $5$ on the plane $(j,l)$ at every sampled point of the filled Julia set, while a generic point of $\RR^8$ generates the whole algebra (the quarters of Plate~\plateref{B07} have dimension $8$). The orbits of the atlas thus explore subalgebras of dimension $3$, $5$ or $8$, depending on the position of the section with respect to the root decomposition.

\subsection{Dependence on the number of iterations}
With $16$--$20$ iterations the body of a Julia set near the boundary of the connectedness locus is not yet converged. Table~\ref{tab:iterations} gives the body share of six plates for $10$, $16$, $20$, $50$ and $100$ iterations. The share decreases monotonically with the number of iterations, and the largest change between $20$ and $100$ iterations is $1.4$ percentage points (F03). The symmetry statements do not depend on this: an exact symmetry maps the escape set of every finite iteration count to itself, while the body shares of the plates must be read as those of the finite-iteration bodies that are drawn.

\begin{table}[htbp]\centering\footnotesize
\begin{tabular}{@{}lrrrrr@{}}\toprule
Plate & 10 & 16 & 20 & 50 & 100\\\midrule
A01 (octonions) & 28.9 & 27.1 & 26.8 & 26.4 & 26.4\\
A07 (Okubo) & 35.6 & 35.1 & 35.0 & 35.0 & 35.0\\
C03 (Okubo) & 35.7 & 35.2 & 35.2 & 35.1 & 35.1\\
F03 (Okubo Mandelbrot) & 9.1 & 6.3 & 5.7 & 4.5 & 4.3\\
J01 (bilateral, Okubo) & 25.2 & 23.6 & 23.3 & 23.2 & 23.2\\
J02 (octonions) & 25.2 & 23.6 & 23.3 & 23.2 & 23.2\\
\bottomrule\end{tabular}
\caption{\emph{Body share (\%, mean over the four quarters, $5\cdot10^4$ points) against the number of iterations.}}
\label{tab:iterations}
\end{table}

\subsection{Validation}
The renderer was validated in three ways. First, the eight-dimensional shaders were compared with the double-precision CPU kernel of the application on $270$ scenes. Second, the independent \texttt{numpy} implementation used for the measurements agrees with the certified implementation of the algebras to machine precision ($0$ deviation on random samples). Third, the GPU was compared with the closed form of Theorem~\ref{thm:B} (Figure~\ref{fig:validation}). The identities used in the proofs, i.e., the torus angles, the anti-automorphism $\sigma$, formula \eqref{eq:okubocubic}, the invariant plane, the bilateral identity and the formulas of Proposition~\ref{prop:theta90}, were checked on random samples to $10^{-12}$ or better, and the derivation algebras and stabilizers were computed both in floating point and exactly, by reduction modulo a prime.

\section{Conclusions and Future Developments}\label{sec:conclusions}

In this work we have studied, for the first time to our knowledge, the quadratic dynamics over the para-octonions and the Okubo algebras, together with their split forms and the Petersson algebras that connect them. The pictures of a four-slices atlas suggested the statements, which we proved and then tested against new pictures, rendered on purpose. The octonions and the para-octonions appear as the algebras of the orthogonal group $O(6)$: the quadratic map only sees the real part and the norm. The Okubo algebra, on the contrary, appears as the algebra of the torus: its Julia sets are symmetric under a maximal torus of $PSU(3)$, which rotates the three root planes of $\mathfrak{su}(3)$ and grows to $U(2)$ on the Weyl walls, and the automorphism group of its symmetrized product is $PSU(3)\rtimes\ZZ_2$. The Mandelbrot sets with seed $0$ of all symmetric composition algebras are pulled back from a single planar picture, the tricorn, through the invariants $\n(c)$ and $\n(c,c\diamond c)$, and the split algebras add to it a hyperbolic tricorn, whose light-cone form is a pair of coupled real quartic maps. It is mesmerizing and surprising how faithfully the pictures report the algebra: the three arms of an Okubonic Mandelbrot section are the three arms of the tricorn, the displaced revolutions sit exactly in the root planes, and the symmetry that seemed unexplained at the Petersson angle $90^\circ$ belongs to a four-dimensional group $T^2\cdot Sp(1)$, half of which consists of twisted isometries that are not automorphisms, and which accounts for all the $42$ coordinate sections with exact revolutions at that angle.

The scope of the results must be stated with care. Theorem~\ref{thm:A} concerns the symmetries of algebraic origin; the absence of other symmetries is only tested numerically. Theorem~\ref{thm:B} concerns the seed $0$, and the Mandelbrot sets with a seed $a\neq0$ are not planar pull-backs. The rendered pictures are finite-iteration approximations and, for the split algebras, they are approximations from inside, because of the bailout and of the transverse instability of the invariant plane.

It would be then definitely interesting to analyze several directions.
\begin{itemize}
\item \emph{Seeded Mandelbrot sets.} With a seed $a\neq0$ the orbit depends on the invariants of the pair $(a,c)$ under the automorphism group. For the para-octonions the relevant group is $O(7)$ and five invariants suffice; for the Okubo algebra the joint invariants of two elements of $\mathfrak{su}(3)$ form a richer ring, and it would be interesting to know whether a finite-dimensional family of planar or low-dimensional models describes these sets.
\item \emph{The generic Petersson algebra.} For $\theta\notin\{0,\pm120^\circ\}$ the algebra is not symmetric and the orbit of the seed $0$ leaves the plane of $c$ and $c\diamond c$; the Mandelbrot set is still invariant under the torus, but it is not a planar pull-back. Among the angles we tested, only $90^\circ$ carries twisted isometries fixing the constants of $C$; at $60^\circ$ there are skew anti-derivations of $\circ$ which do not fix them. Which constants they fix, and whether a complete description of the twisted symmetries along the whole path exists, are natural questions.
\item \emph{Split algebras.} Theorem~\ref{thm:A} for the split Okubo algebra requires the classification of the Cartan subgroups of $PSU(2,1)$, compact and non-compact; the rendered pictures, on the other hand, require escape criteria adapted to indefinite norms and arithmetic that controls the transverse instability.
\item \emph{Rigidity.} Excluding symmetries of Julia sets that do not come from the algebra is a problem of a different nature, which we leave open.
\end{itemize}
The next steps of our investigation will clearly include the extension of this study to the exceptional Jordan algebras built over the Okubo algebra, whose projective geometry we have studied elsewhere \cite{CMZ24,CZ22}, and that we plan to cover in a forthcoming work.

\section*{Use of Generative AI}
The implementation of the rendering, measurement and verification scripts listed in Appendix~\ref{app:repro}, the production of the plates, the numerical checks of the statements and a first draft of the text were carried out with AI assistance (Claude, Anthropic) under the author's direction; the commentaries of the plates were checked against the images by independent AI agents and corrected. AI tools were also used for language editing during manuscript preparation. The author takes full responsibility for the content of this article.

\appendix
\section{Reproducibility}\label{app:repro}
\emergencystretch=1.5em
All scripts are ancillary files of the arXiv version of this article (folder \path{scripts}, whose \path{README.md} documents the commands); every run writes a new folder in \path{results} with a \path{manifest.json} that records the command, the environment, the input and output hashes, the conclusions and the limitations. The rendered images are the plates of Appendix~\ref{app:atlas}, and the renderer \cite{HF8D} is not among the ancillary files. The runs used in this work are the following.
\begin{itemize}
\item \path{run_verify_theory_2026-09-30T07-21-50} (\path{verify_theory_v2.py}): the torus, $\sigma$, $-\mathrm{id}$, the derivation algebras and stabilizers (in floating point and modulo $p=2147483579$), formula \eqref{eq:okubocubic}, Proposition~\ref{prop:theta90}, Theorem~\ref{thm:B} on random constants, including the split case with $80$-digit arithmetic, the invariant sections of Theorem~\ref{thm:A}, the adapted frame of Plate~\plateref{R09}; $14$ checks, all passed. The earlier run \path{...T07-19-23} of version~1 of the script gives the same results with a different generic test constant.
\item \path{run_verify_extra_2026-09-30T07-47-42} (\path{verify_extra_v1.py}): Proposition~\ref{prop:bilateral} and the tilted revolutions of Theorem~\ref{thm:A}(4) and Corollary~\ref{cor:tilted}; $5$ checks, all passed.
\item \path{run_atlas_render_2026-09-29T17-56-41} and \path{...T18-06-48} (\path{atlas_render_v1.mjs}): the plates of the first version; \path{run_atlas_render_2026-09-30T07-23-18} (\path{atlas_render_v2.mjs}, \path{plates_v3.mjs}): the new plates R01--R17.
\item \path{run_measure_v2_2026-09-30T13-09-48} (\path{measure_v2.py}): the measurements of Section~\ref{sec:numerics}.
\item \path{run_validate_gpu_2026-09-30T07-29-55} and \path{run_compare_gpu_2026-09-30T07-31-48}: the validation of the renderer against Theorem~\ref{thm:B}.
\item \path{run_figures_v2_2026-09-30T13-14-41} (\path{figures_v2.py}): the figures.
\item \path{run_verify_twisted_2026-09-30T13-15-21} (\path{verify_twisted_v1.py}): the skew derivations and anti-derivations of $\circ$ along the Petersson path, identity \eqref{eq:twisted}, the one-step criterion on the $56\times3$ pairs (section, axis) and the scans at $60^\circ$ and $150^\circ$.
\item \path{build_article_v2.py}: the assembly of this article; every number that depends on a run is inserted from that run, and the inserted values are recorded in \path{build_article_v2.json}.
\end{itemize}
The derivation algebras are computed as the null space of the linear system $D(e_a\circ e_b)=D(e_a)\circ e_b+e_a\circ D(e_b)$ in the $64$ entries of $D$; for the exact computation the structure constants, which lie in $\ZZ[\tfrac12,\sqrt3]$, are reduced modulo a prime $p\equiv11\pmod{12}$, in which $3$ is a square. The full atlas of $134$ plates, with the measurements of every quarter and the commentaries, is Appendix~\ref{app:atlas}.

\newgeometry{top=16mm,bottom=16mm,left=17mm,right=17mm,headsep=5mm,footskip=8mm}
\begingroup
\setlength{\parskip}{0.35em}
\makeatletter
\def\subsection{\@startsection{subsection}{2}\z@{.7\linespacing\@plus\linespacing}{.5\linespacing}{\normalfont\scshape\centering}}
\def\subsubsection{\@startsection{subsubsection}{3}\z@{.5\linespacing\@plus.7\linespacing}{-.5em}{\normalfont\bfseries}}
\makeatother
\section{The atlas: four slices of eight-dimensional fractals}\label{app:atlas}

This appendix is the atlas on which the work rests: $134$ plates, each one commented in the light of the algebra and accompanied by the measurements of its four quarters. It was conceived as a document that can be read on its own, and we kept it so: it has its own introduction, its reading keys (Section~\ref{atl:sec:algebras}) and its method (Section~\ref{atl:sec:method}); the plates follow in Parts~A--P and~R, and a synthesis closes it (Section~\ref{atl:sec:synthesis}). Beside the quadratic map of the main text, it covers the Mandelbrot family, the cubic, the ``sandwich'', the bilateral and the logistic map and a sine map, over octonions, split-octonions, the degenerate semi-octonions, para-octonions and the compact and split Okubo algebras, together with the Petersson algebras that interpolate between para-octonions and Okubo. The statements called Theorem~A and Theorem~B in this appendix are those of Section~\ref{sec:intro}.

\subsection{Introduction and Motivation}

Since the pictures of Norton \cite{Norton82} and the quaternionic Julia sets of Holbrook \cite{Holbrook} and of Hart, Sandin and Kauffman \cite{HSK}, the iteration of $z\mapsto z^2+c$ over algebras larger than the complex numbers has produced some of the most admired objects of mathematical visualization. The reason is simple and deep at the same time: a filled Julia set is an invariant object of the dynamics, and every symmetry of the algebra that commutes with the map is inherited by the picture. Griffin and Joshi were the first to iterate over the octonions \cite{GriffinJoshi}, Rochon and Garant-Pelletier developed a Fatou--Julia theory over bicomplex and multicomplex numbers \cite{Rochon,GarantRochon}, Senn discovered that over the hyperbolic numbers the Mandelbrot set becomes a square \cite{Senn}, and Crowe, Hasson, Rippon and Strain-Clark studied the ``Mandelbar'' set of $\bar z^2+c$ \cite{CHRS}. A concise overview of hypercomplex fractals is Katunin's monograph \cite{Katunin17}. In all these works, however, the algebra is either associative or alternative, and the images are usually a single three-dimensional slice.

Real eight-dimensional composition algebras are a much richer landscape. Besides the \emph{octonions} $\OO$ and the \emph{split-octonions} $\OO_s$, the unital Hurwitz algebras, there are the \emph{symmetric composition algebras}: the \emph{para-octonions} $p\OO$ and $p\OO_s$, with product $x\bullet y=\bar x\bar y$, and the \emph{Okubo algebras} $\mathcal O$ and $\mathcal O_s$, discovered by Okubo in the search for a ``pseudo-octonionic'' structure behind $SU(3)$ \cite{Okubo78} and characterized by Elduque and Myung \cite{ElduqueMyung90,Elduque15}. While octonions have a unit and para-octonions a para-unit, the Okubo algebra merely contains idempotent elements; it enjoys \emph{flexibility} instead of \emph{alternativity}, and its automorphism group is $PSU(3)$ instead of $G_2$ \cite{Elduque15,Elduque18,CZ22,MCZ25}. A comparison of these algebras is summarized in Table~\ref{atl:tab:algebras}. From a dynamical point of view the difference is dramatic: in a quadratic algebra such as $\OO$ the square $z^2=2\re(z)\,z-\n(z)$ never leaves the plane of $1$ and $z$, so that every orbit of $z^2+c$ lives in a three-dimensional subspace; in the Okubo algebra, on the contrary, $z_0$ and $c$ generically generate the whole algebra and the orbits genuinely explore eight dimensions.

\begin{table}[h]
\centering\footnotesize
\setlength{\tabcolsep}{3.5pt}\begin{tabular}{@{}llllll@{}}
\toprule
Algebra & Product & Unit & Identity & Automorphisms & Plates \\
\midrule
Octonions $\OO$ & Cayley--Dickson, $\gamma=(-1,-1,-1)$ & $1$ & alternative & $G_2$ (compact) & A01, B01, C01, \dots\\
Split-octonions $\OO_s$ & $\gamma=(-1,-1,1)$ or $(1,-1,-1)$ & $1$ & alternative & $G_{2(2)}$ & A02, A03, \dots\\
Semi-octonions & $\gamma=(-1,-1,0)$ (degenerate norm) & $1$ & alternative & non-reductive & A04, B04, D04\\
Para-octonions $p\OO$ & $x\bullet y=\bar x\,\bar y$ & para-unit & symmetric & $G_2$ & A05, B05, C02\\
Split para-octonions & $\bar x\,\bar y$, $\gamma_3=+1$ & para-unit & symmetric & $G_{2(2)}$ & A06, B06\\
Okubo $\mathcal O$ & $\tau(\bar x)\,\tau^2(\bar y)$, $\theta=120^\circ$ & none & symmetric & $PSU(3)$ & A07, B07, C03, R01--R09\\
Split Okubo $\mathcal O_s$ & $\gamma_3=+1$, $\theta=120^\circ$ & none & symmetric & $PSU(2,1)$ & A08, B08\\
Petersson, $\theta\ne0,\pm120^\circ$ & $\tau_\theta(\bar x)\,\tau_\theta^{-1}(\bar y)$ & none & composition only & $\dim 4$ & C04--C06, E, R10--R16\\
\bottomrule
\end{tabular}
\caption{\emph{Synoptic table of the eight-dimensional algebras of this atlas.} The Cayley--Dickson parameters $\gamma=(\gamma_1,\gamma_2,\gamma_3)$ are the squares $i^2,j^2,l^2$; $\tau_\theta$ is the Petersson twist of Section~\ref{atl:sec:algebras}. The last column lists representative plates.}
\label{atl:tab:algebras}
\end{table}

In this work we explore, systematically and in a comparative spirit, how these algebras shape the dynamics of \emph{very simple} maps. We deliberately avoid elaborate formulas: the quadratic map $z\mapsto z\diamond z+c$, where $\diamond$ is the product of the algebra, already tells most of the story. Around it we place its Mandelbrot family, the Mandelbar map, the cubic $(z\diamond z)\diamond z+c$, the ``sandwich'' $(z\diamond a)\diamond z+c$, the bilateral map $(a\diamond z)\diamond(z\diamond b)+c$, the logistic map $c\diamond(z-z\diamond z)$ and a sine map; we then vary the four-dimensional mappings applied to the sample point and the coloring algorithms. Every image is a \emph{four slices} view: the quarters show the four three-dimensional faces of the four-dimensional space spanned by three visible axes and one hidden axis. Two sections that should coincide by a symmetry appear side by side, and so do two sections that should differ; the images become an experiment.

The main output of this work is an atlas of $134$ plates, each one commented in the light of the algebra and accompanied by measurements: the $117$ plates of the first version, whose commentaries were checked again against the images and corrected, and the $17$ plates of Part~R, which test the predictions of the main text. A synoptic reading of the atlas is given in Section~\ref{atl:sec:synthesis}; the elementary propositions that explain the pictures are collected in Section~\ref{atl:sec:algebras}, with their proofs, and the two theorems that govern the Okubo algebra are proved in the main text. Among the findings we single out the following.
\begin{itemize}
\item In the octonions and para-octonions, with a complex constant $c$, every section through $1$ and a direction orthogonal to $c$ is a \emph{body of revolution} and every section through the purely imaginary directions is a \emph{ball}: the dynamics is equivariant under the full orthogonal group $O(6)$, much larger than the stabilizer of $c$ in $G_2$ (Proposition~\ref{atl:prop:equivariance}).
\item In the Okubo algebra, with the same constant, the frames $F_A$ and $F_C$ show no body of revolution outside the quaternions. The continuous symmetry, however, is not lost: a two-dimensional torus fixes $c$ and rotates the three root planes $V_1=\langle j,k\rangle$, $V_2=\langle l,il\rangle$, $V_3=\langle jl,kl\rangle$ of $\mathfrak{su}(3)$, and exactly the coordinate sections that contain a root plane are bodies of revolution (Proposition~\ref{atl:prop:torus}, Plates~R01--R03). On the Weyl walls the symmetry grows to $U(2)$ (Plates~R05--R08), and for a generic constant the revolutions reappear in a frame adapted to the stabilizer (Plate~R09).
\item The Mandelbrot sets with seed $0$ are pulled back from a planar set by the pair $(\n(c),\n(c,c\diamond c))$: in the compact symmetric composition algebras this set is the tricorn, and the three arms of the Okubonic sections are the three arms of the tricorn, redistributed by $PSU(3)$ (Proposition~\ref{atl:prop:okubocubic}).
\item The identities of symmetric composition algebras are visible at a glance: the sandwich degenerates into perfect balls and the cubic into bodies of revolution about the direction of $c$, identical in the para-octonions and in the Okubo algebra, and every bilateral map is an octonionic quadratic map in disguise, Kaplansky's trick $(1\diamond x)\diamond(y\diamond1)=xy$ being the special case $a=b=1$ (Propositions~\ref{atl:prop:symmetric} and \ref{atl:prop:kaplansky}).
\item The indefinite norm of the split algebras turns spheres into hyperboloids, bodies of revolution into horns and saddles; the degenerate semi-octonions produce exact cylinders; along the Petersson path a single angle, $\theta=90^\circ$, makes four different sections coincide in a perfect ball and produces exact revolutions that no automorphism explains (Proposition~\ref{atl:prop:theta90}).
\item Every filled Julia set of an even map is centrally symmetric, $z\mapsto-z$ being a symmetry of $z\diamond z+c$ in every algebra; in the sections through the origin every mirror is a half-turn composed with it (Remark~\ref{atl:rem:minus}).
\end{itemize}
In Sec.~\ref{atl:sec:algebras} we review the algebras and prove the reading keys. In Sec.~\ref{atl:sec:method} we describe the renderer, the four slices, the choice of the constants and the measurements, and we list the corrections with respect to the first version. Sections~\ref{atl:part:A}--\ref{atl:part:R} contain the plates, one per page. Finally, in Sec.~\ref{atl:sec:synthesis} we summarize the findings.

\subsection{The algebras and the reading keys}\label{atl:sec:algebras}

In this section we recall the eight-dimensional algebras used in the atlas, all realized on the same vector space $\RR^8$ with basis $1,i,j,k,l,il,jl,kl$, and we collect a small number of elementary facts that will be used, page after page, to read the plates. None of these facts is deep; their interest lies in the fact that each of them has a visible signature.

\subsubsection{Cayley--Dickson algebras with parameters}
Let $Q(\gamma_1,\gamma_2)$ be the quaternion-type algebra with basis $1,i,j,k=ij$, where $i^2=\gamma_1$, $j^2=\gamma_2$ and $ij=-ji$. Doubling it with a new unit $l$ such that $l^2=\gamma_3$, one obtains on $Q\oplus Ql$ the product
\begin{equation}\label{atl:eq:cd}
(p+ql)(r+sl)=\big(pr+\gamma_3\,\bar s q\big)+\big(sp+q\bar r\big)l ,
\end{equation}
where the bar is the quaternionic conjugation. The resulting algebra $\OO(\gamma)$, $\gamma=(\gamma_1,\gamma_2,\gamma_3)$, is unital and quadratic: every element $x$ satisfies $x^2-2\re(x)\,x+\n(x)=0$, with the diagonal norm
\begin{equation}\label{atl:eq:norm}
\n(x)=x\bar x=x_0^2-\gamma_1x_1^2-\gamma_2x_2^2+\gamma_1\gamma_2x_3^2-\gamma_3x_4^2+\gamma_1\gamma_3x_5^2+\gamma_2\gamma_3x_6^2-\gamma_1\gamma_2\gamma_3x_7^2 .
\end{equation}
For $\gamma=(-1,-1,-1)$ one recovers Cayley's octonions $\OO$; for $\gamma_3=+1$, or equivalently $\gamma_1=+1$, the split-octonions $\OO_s$, whose norm has signature $(4,4)$ and which contain zero divisors; for $\gamma_3=0$ the degenerate \emph{semi-octonions}. When all $\gamma_k$ are non-zero, $\OO(\gamma)$ is a composition algebra, i.e., $\n(xy)=\n(x)\n(y)$. The two realizations of the split-octonions used in the atlas are isomorphic as algebras, but they place the complex plane $\CC=\RR1\oplus\RR i$ differently: with $\gamma_3=+1$ the plane of $1$ and $i$ is still the complex plane, with $\gamma_1=+1$ it becomes the plane of the hyperbolic (split-complex) numbers.

\subsubsection{Para-Hurwitz and Petersson algebras}
In order to leave the realm of unital algebras without losing the composition property, one twists the product by involutions and automorphisms. The \emph{para-Hurwitz} product $x\bullet y=\bar x\,\bar y$ defines the para-octonions; the element $1$ becomes a \emph{para-unit}, $1\bullet x=x\bullet1=\bar x$. Petersson \cite{Petersson} observed that an automorphism $\tau$ of order three gives a new composition algebra with the product
\begin{equation}\label{atl:eq:petersson}
x*y=\tau(\bar x)\,\tau^2(\bar y).
\end{equation}
Following the implementation of the renderer, we use the automorphisms $\tau_\theta(p+ql)=p+(u_\theta q)l$ with $u_\theta=\cos\theta+i\sin\theta$, and the product $\tau_\theta(\bar x)\,\tau_\theta^{-1}(\bar y)$. For $\theta=120^\circ$ the automorphism has order three and \eqref{atl:eq:petersson} is the Okubo algebra, compact for $\gamma=(-1,-1,-1)$ and split for $\gamma_3=+1$; for $\theta=0$ it is the para-octonions. For a generic angle one still obtains a composition algebra, which is not symmetric; the family is a continuous path from para-octonions to Okubo.

Para-octonions and Okubo algebras are \emph{symmetric composition algebras}, i.e., they enjoy the identities
\begin{equation}\label{atl:eq:symmetric}
(x*y)*x=x*(y*x)=\n(x)\,y .
\end{equation}
It is of paramount importance to notice that these identities, harmless in appearance, have drastic dynamical consequences.

\subsubsection{Reading keys}
We now state the propositions that explain the plates. The proofs are short; all identities were also verified numerically on random samples by the script \path{verify_propositions_v3.py} (Section~\ref{atl:app:repro}).

\begin{proposition}\label{atl:prop:equivariance}
Let $\diamond$ be the product of $\OO(\gamma)$ or of its para-Hurwitz algebra, and let $g$ be a linear map of $\RR^8$ fixing $1$ and preserving the norm \eqref{atl:eq:norm}. Then $g(z\diamond z)=g(z)\diamond g(z)$. In particular, if $g(c)=c$, the map $z\mapsto z\diamond z+c$ commutes with $g$, its filled Julia set is $g$-invariant, and for the compact algebras and $c\in\CC$ the symmetry group contains $O(6)$ acting on the span of $j,k,l,il,jl,kl$.
\end{proposition}
\begin{proof}
Since $\OO(\gamma)$ is quadratic, $z^2=2\re(z)z-\n(z)$ and $\overline{z^2}=2\re(z)\bar z-\n(z)$. A norm-preserving map fixing $1$ preserves the polar form, hence the real part $\re z=\tfrac12\n(z,1)$, and commutes with the conjugation. The claim follows at once. \end{proof}

Proposition~\ref{atl:prop:equivariance} is the key to the octonionic plates: the symmetry of $z^2+c$ is not the stabilizer of $c$ in $G_2$, which is $SU(3)$, but the much larger group $O(6)$. It follows that every section through $1$, $i$ and a third direction orthogonal to $c$ is the same body, that sections through $1$ and two directions orthogonal to $c$ are bodies of revolution, and --- by the next proposition --- that sections of purely imaginary vectors are balls. For the split algebras the group is a pseudo-orthogonal group, and circles of revolution become hyperbolas.

\begin{proposition}\label{atl:prop:imaginary}
In $\OO(\gamma)$ and in its para-Hurwitz algebra, a purely imaginary $z$ satisfies $z\diamond z=-\n(z)$. Hence on the hyperplane $\re z_0=0$ the first iterate $z_1=c-\n(z_0)$ depends only on $\n(z_0)$, and the filled Julia set meets that hyperplane in a union of level sets of the norm: balls for the compact algebras, regions between hyperboloids for the split ones. In the Okubo algebra the identity fails.
\end{proposition}
\begin{proof}
If $\re z=0$ then $\bar z=-z$ and $z^2=-z\bar z=-\n(z)$; the para-Hurwitz square is $\overline{z^2}=-\n(z)$ as well. \end{proof}

\begin{proposition}\label{atl:prop:span}
In $\OO(\gamma)$ and in the para-Hurwitz algebras, every orbit of $z\mapsto z\diamond z+c$ lies in the three-dimensional space spanned by $1$, $\im z_0$ and $\im c$; every orbit of the Mandelbrot map with seed $0$ lies in the plane of $1$ and $c$. In the symmetric composition algebras the Mandelbrot orbit with seed $0$ lies in the plane spanned by $c$ and $c\diamond c$.
\end{proposition}
\begin{proof}
The first two claims follow from $z^2=2\re(z)z-\n(z)$ by induction. For the third, \eqref{atl:eq:symmetric} gives $c*(c*c)=(c*c)*c=\n(c)c$ and $(c*c)*(c*c)=\n(c,c*c)\,c-\n(c)\,c*c$, where $\n(x,y)=\n(x+y)-\n(x)-\n(y)$ is the polar form; so the span of $c$ and $c*c$ is invariant. \end{proof}

The escape of the seed-$0$ Mandelbrot orbit is therefore governed by the two invariants $\n(c)$ and $K=\n(c,c\diamond c)$, the coordinates of a planar map: all Mandelbrot sets of symmetric composition algebras are pulled back from one planar set. Theorem~B identifies it: when $\Delta=4\n(c)^3-K^2>0$, the plane of $c$ and $c\diamond c$ is isometric to $(\CC,\bar x\bar y)$ through $c\mapsto\zeta$ with $|\zeta|^2=\n(c)$ and $2\re\zeta^3=K$, so that $c$ belongs to the Mandelbrot set if and only if $\zeta$ belongs to the tricorn; when $\Delta<0$, which only happens in the split algebras, the tricorn is replaced by its analogue over the split-complex numbers, the \emph{hyperbolic tricorn}. What distinguishes the algebras is then the geometry of the cubic form $c\mapsto\n(c,c\diamond c)$ on the sections, i.e., the group that moves the tricorn: $G_2$ for the para-octonions, $PSU(3)$ for the Okubo algebra.

\begin{proposition}\label{atl:prop:okubocubic}
In the para-octonions $\n(c,c\bullet c)=2c_0(c_0^2-3|\im c|^2)$ depends only on $c_0$ and $|\im c|$. In the Petersson algebra of angle $\theta$, for $c=c_0+ql$ with $c_0\in\RR$ and $q\in\HH$,
\begin{equation}\label{atl:eq:okubocubic}
c*c=c_0^2-\bar q\,u_\theta^2\,q-2\cos\theta\;c_0\,q\,l ,\qquad \n(c,c*c)=2c_0\big(c_0^2-(\cos2\theta+2\cos\theta)\,|q|^2\big),
\end{equation}
which depends only on $c_0$ and $|q|$; but for general $c$ the Okubo cubic form is not a function of $c_0$ and $|\im c|$. Neither $q\,u_\theta^2\,q$ nor $\bar q\,\bar u_\theta^{\,2}\,q$ may replace the middle term of \eqref{atl:eq:okubocubic}.
\end{proposition}
\begin{proof}
For the para-octonions $c\bullet c=\overline{c^2}=2c_0\bar c-\n(c)$, and $\n(c,\bar c)=2(c_0^2-|\im c|^2)$ gives the formula. For $c=c_0+ql$ one has $\tau(\bar c)=c_0-(uq)l$ and $\tau^{-1}(\bar c)=c_0-(\bar uq)l$; the product \eqref{atl:eq:cd} with $\gamma_3=-1$ gives $c_0^2-\overline{(\bar uq)}(uq)-c_0(\bar u+u)\,q\,l$, where the first factor of the middle term is the \emph{conjugate} of $\bar uq$; since $\overline{\bar u q}=\bar q u$ this is \eqref{atl:eq:okubocubic}. The polar form with $c$ uses $\re(\bar q u^2 q)=\re(u^2)|q|^2$. The last claim follows from the fact that the orbits of $PSU(3)$, cut out by $\n(c)$ and $K$, are smaller than those of $G_2$, cut out by $c_0$ and $|\im c|$; it is visible in the plates F03 and F11. \end{proof}

\begin{proposition}\label{atl:prop:torus}
For every angle $\theta$ and every pair of unit complex numbers $a,b\in\CC$, the map $\varphi_{a,b}(p+ql)=ap\bar a+(bq\bar a)l$ is an automorphism of the Petersson algebra \eqref{atl:eq:petersson}. The quaternions $\HH=\RR1\oplus\RR i\oplus\RR j\oplus\RR k$ form a subalgebra of every Petersson algebra, on which the product is $\bar x\,\bar y$.
\end{proposition}
\begin{proof}
The maps $\varphi_{a,b}$ are automorphisms of the octonions, they commute with the conjugation and with $\tau_\theta$ because $a,b,u_\theta$ lie in the commutative field $\CC$; hence they preserve \eqref{atl:eq:petersson}. The automorphism $\tau_\theta$ fixes $\HH$ pointwise, so on $\HH$ the product is the para-quaternionic one. \end{proof}

This proposition sets the contrast between the octonionic and the Okubonic plates, and it must be read with care. With $c\in\CC$ the whole torus $\{\varphi_{a,b}\}$, $a=e^{i\alpha}$, $b=e^{i\beta}$, fixes $c$, and it acts by rotations on three coordinate planes,
\begin{equation}\label{atl:eq:rootangles}
V_1=\langle j,k\rangle:\ 2\alpha,\qquad V_2=\langle l,il\rangle:\ \beta-\alpha,\qquad V_3=\langle jl,kl\rangle:\ \alpha+\beta ,
\end{equation}
which are the three root planes of $\mathfrak{su}(3)$ with respect to the Cartan plane $\CC=\langle1,i\rangle$ (the angle on $V_3$ is the sum of the other two). A three-dimensional coordinate section is mapped into itself by a one-parameter subgroup of the torus if and only if it contains one of the root planes, and in that case it is a body of revolution about its third axis: $18$ of the $56$ coordinate sections have this property (Theorem~A). The frames $F_A=(1,i,j\,|\,l)$ and $F_C=(1,j,l\,|\,i)$ contain no root plane, and this is why the Okubonic plates in these frames show no revolution; the frame $(1,l,il\,|\,j)$ of Plate~R01 does. The continuous symmetry that the octonions owe to Proposition~\ref{atl:prop:equivariance}, of dimension $15$, is reduced to a torus of dimension $2$, but it is not lost: the revolutions are \emph{displaced} to the root planes. The torus is completed by the anti-automorphism $\sigma=\mathrm{diag}(1,1,-1,1,-1,1,-1,1)$, which reverses every product, $\sigma(x*y)=\sigma y*\sigma x$, and is therefore an automorphism of the symmetrized product and a symmetry of every quadratic map with $c\in\CC$. When $\arg c$ is a multiple of $60^\circ$, $c$ lies on a Weyl wall, the stabilizer of $c$ grows to $U(2)$, and new revolutions appear (Plates~R05--R08). For a generic $c\in\RR^8$ the stabilizer is a torus conjugate to $\{\varphi_{a,b}\}$, whose root planes are in general position with respect to the coordinate axes (Plates~B07 and R09). The sections contained in $\HH$ coincide with those of the para-octonions.

\begin{remark}\label{atl:rem:minus}
In every algebra $(-z)\diamond(-z)=z\diamond z$, so that every even map (the quadratic map, the bilateral and the sandwich maps, $z^p+c$ with $p$ even) has the symmetry $z\mapsto-z$, and its filled Julia set is centrally symmetric, with equal escape times. In a section through the origin the composition of $-\mathrm{id}$ with the half-turn about an axis is the mirror in the plane orthogonal to that axis: mirrors and half-turns come in pairs, and many of the finite symmetries of the tables are of this kind. Shapes with a single tip can appear only in sections that do not pass through the origin. The same conclusion holds for an odd map without additive constant, such as $c\sin z$, since $f(-z)=-f(z)$ maps the orbit of $z$ to the opposite of the orbit of $z$; and the logistic map $c\diamond(z-z\diamond z)$ is even about its critical point $\tfrac12$, since $z-z\diamond z=\tfrac14-w\diamond w$ for $z=\tfrac12+w$, so that its filled set is centrally symmetric about $\tfrac12$. The cubic map $(z\diamond z)\diamond z+c$ has none of these symmetries, and neither have the Mandelbrot sets, whose variable is the constant.
\end{remark}

\begin{proposition}\label{atl:prop:symmetric}
In a symmetric composition algebra the sandwich map satisfies $(z*a)*z+c=\n(z)\,a+c$ and the cubic map $(z*z)*z+c=\n(z)\,z+c$. The norms $t_n=\n(z_n)$ of a sandwich orbit obey the real quadratic recursion
\begin{equation}\label{atl:eq:sandwich}
t_{n+1}=\n(a)\,t_n^2+\n(a,c)\,t_n+\n(c),
\end{equation}
conjugate to $s\mapsto s^2+c_{\mathrm{eff}}$ with $c_{\mathrm{eff}}=\n(a)\n(c)+\tfrac12\n(a,c)-\tfrac14\n(a,c)^2$. The filled set of the sandwich map is the region $t_{\min}\le\n(z_0)\le t_{\max}$ of the real Julia interval, a ball in the compact case. The cubic map $\n(z)\,z+c$, on the contrary, is not a function of the norm alone, since $\n(z_1)$ also depends on $\n(z_0,c)$: it commutes with every isometry fixing $c$, and its filled set is a body of revolution about the direction of $c$.
\end{proposition}
\begin{proof}
The first two identities are \eqref{atl:eq:symmetric}. Taking norms of $z_{n+1}=t_na+c$ gives \eqref{atl:eq:sandwich}, and the affine change $s=\n(a)t+\tfrac12\n(a,c)$ conjugates it to $s^2+c_{\mathrm{eff}}$. \end{proof}

\begin{proposition}[Kaplansky's trick]\label{atl:prop:kaplansky}
In the para-octonions and in the Okubo algebra $(1*x)*(y*1)=xy$, the octonionic product. Hence the bilateral map with $a=b=1$ is not merely conjugate to the octonionic quadratic map: it coincides with it. More generally, in every symmetric composition algebra $(a*z)*(z*b)=\n(z,e)\,z-\n(z)\,e$ with $e=b*a$, so that every bilateral map is an octonionic quadratic map moved by a similarity (article, Proposition on bilateral maps).
\end{proposition}
\begin{proof}
For the para-octonions $(1\bullet x)\bullet(y\bullet1)=\bar{\bar x}\,\bar{\bar y}=xy$. For the Okubo algebra, $1*x=\tau^2(\bar x)$ and $y*1=\tau(\bar y)$, so $(1*x)*(y*1)=\tau(\tau^2(x))\,\tau^2(\tau(y))=xy$, since $\tau^3=1$. \end{proof}

\begin{proposition}\label{atl:prop:theta90}
For $\theta=90^\circ$ and $x=x_0+ql$ with $x_0\in\RR$, $q\in\HH$, the Petersson square is $x*x=x_0^2+|q|^2$. Consequently every section contained in $\RR\oplus\HH l$ meets the filled Julia set of $z*z+c$ in the same ball. Moreover the subspace $W=\langle1,i\rangle\oplus\langle l,kl\rangle$ is closed under squaring: for $z=p+ql$ with $p\in\langle1,i\rangle$, $q\in\langle1,k\rangle$,
\begin{equation}\label{atl:eq:W}
z*z=\bar p^{\,2}+|q|^2+2\im(p)\,\bar q\,l ,
\end{equation}
and the rotations $R_\varphi(p+ql)=p+(e^{k\varphi}q)l$ satisfy $f(R_\varphi z)=R_{-\varphi}f(z)$ on $W$ for $c\in\CC$. This \emph{twisted} equivariance makes $K_c\cap W$ a body of revolution, and since $z*z\in W$ for $z\in\langle i,j,l\rangle$, the section $\langle i,j,l\rangle$ is a body of revolution about $l$ (Plates~C05 and R10).
\end{proposition}
\begin{proof}
Put $\cos\theta=0$ and $u^2=-1$ in \eqref{atl:eq:okubocubic}. For \eqref{atl:eq:W} and the twisted equivariance see Proposition~\ref{prop:theta90}. \end{proof}

\begin{proposition}\label{atl:prop:rescaling}
For $\gamma_3<0$ the map $p+ql\mapsto p+\sqrt{-\gamma_3}\,ql$ is an isomorphism $\OO(\gamma_1,\gamma_2,\gamma_3)\to\OO(\gamma_1,\gamma_2,-1)$. For $\gamma_3=0$, $(p+ql)^2=p^2+2\re(p)\,q\,l$.
\end{proposition}
\begin{proof}
Both claims are direct substitutions in \eqref{atl:eq:cd}. \end{proof}

The first half of Proposition~\ref{atl:prop:rescaling} explains why, along the path $\gamma_3\in(-1,0)$, the octonionic bodies are simply stretched by $1/\sqrt{-\gamma_3}$ in the directions of $l$; the second half explains the semi-octonionic cylinders. The quaternionic part $p_n$ evolves autonomously as a quaternionic Julia set, while the $l$-part is driven linearly, $q_{n+1}=2\re(p_n)q_n+c_l$: on the hyperplane $\re p_0=0$ the first iterate forgets $q_0$ altogether, and the body is an exact cylinder along the $l$-directions.

\begin{remark}\label{atl:rem:dual}
In the dual plane ($\gamma_1=0$) the orbit $z_n=x_n+y_ni$ of $z^2+c$ satisfies $y_n=c_1\,\partial x_n/\partial c_0$: the dual part is the parameter derivative of the real orbit. Near the period-doubling point $c_0=-3/4$ this derivative is large, and along the ray of plate D08 the critical orbit escapes in a small window $-0.7605<\re c<-0.7468$ around it (verified numerically).
\end{remark}

\begin{remark}\label{atl:rem:mappings}
The four-dimensional mappings of the atlas (Part~M) are the transformations of HyperFractal Studio, which act on the first four coordinates $1,i,j,k$ of the eight-dimensional sample and leave $l,il,jl,kl$ untouched. A mapping that depends on the coordinate $k$ (``$W$'') is therefore the identity on sections with $x_3=0$; this explains why some mappings are invisible in some quarters. Since the mappings act on $1,i,j,k$ only, they also break, artificially, the symmetries of the root planes $V_2=\langle l,il\rangle$ and $V_3=\langle jl,kl\rangle$ of the Okubo algebra: the plates of Part~M do not respect the algebra, and they are not measured.
\end{remark}

\subsection{Method}\label{atl:sec:method}

\subsubsection{Four slices}
All images were rendered with HyperFractal~8D~1.0, an eight-dimensional extension of an Ultra Fractal-compatible explorer: every formula written for the complex plane is \emph{lifted} to $\RR^8$ by reading its complex variables as elements of the chosen algebra \cite{HF8D}. A frame is an orthonormal basis $E_0,\dots,E_7$ of $\RR^8$; the first three vectors are the visible axes $X,Y,Z$ and $E_3$ is the active hidden axis $N$. In the \emph{four slices} mode the image is divided into quarters, each showing a three-dimensional section through the center of the four-dimensional box spanned by $X,Y,Z,N$:
\begin{gather*}
\text{top left }(X,Y,Z)\perp N,\qquad \text{top right }(X,Y,N)\perp Z,\\
\text{bottom left }(X,Z,N)\perp Y,\qquad\text{bottom right }(Y,Z,N)\perp X .
\end{gather*}
In every quarter the first axis points to the right and forward, the second upwards and the third to the left and forward, as shown by the tripod in the corner; the label gives the three axes and the axis perpendicular to the section. The four frames used most often are $F_A=(1,i,j\,|\,l)$, $F_B=(1,l,il\,|\,jl)$, $F_C=(1,j,l\,|\,i)$ and $F_D=(i,j,l\,|\,il)$.

The body drawn in each quarter is the filled set of the map computed with the iteration count indicated on the plate: a point belongs to the body if its orbit does not escape the ball of squared radius $100$ (for the Mandelbrot maps the orbit of the seed with $c$ equal to the point). The surface is found by sphere tracing with a distance estimate derived from the escape potential, and it is colored by the \emph{outside} coloring evaluated near the surface, with soft shadows and ambient occlusion. The cube is the box of half size $h$ (the \emph{half-extent}); a body larger than the cube is cut by its faces, which then show flat cross-sections.

\subsubsection{Constants}
Most Julia sets are taken close to the boundary of their connectedness locus, where they are most articulate. The notation $c=0.97\,s^*$ means: on the ray of a fixed direction $d$, $s^*$ is the largest $s$ such that the orbit of the critical point $0$ of $z\mapsto z\diamond z+s d$ stays bounded for $400$ iterations, and $c$ is taken at $97\%$ of that distance. The two directions used are the complex direction $d_{\CC}\propto(-0.35,0.94)$ and a generic eight-dimensional direction $d_{\mathcal O}$ found by a parameter scan of the prototype HF8D~0.2. The same direction gives different constants in different algebras, since $s^*$ depends on the algebra. For $c\in\CC$ the critical orbit stays in $\CC$ in all these algebras; on $\CC$ the octonionic product is the complex one, while the para-octonionic and the Okubonic products coincide with $\bar z\,\bar w$ (Proposition~\ref{atl:prop:torus}). This is why the para-octonions and the Okubo algebra share their constants in Parts~A and~C, and the octonions have a different one.

\subsubsection{Measurements}
Every plate without four-dimensional mappings is accompanied by a table, computed in double precision by an implementation of the algebras independent of the renderer (\path{hf8np_v1.py}), on $10^5$ uniform points of the ball of radius $0.95h$ in each quarter, with the iteration count and bailout of the plate:
\begin{itemize}
\item \emph{Body}: the share of the ball occupied by the filled set, with its $95\%$ binomial interval.
\item \emph{Exact} and \emph{approximate} symmetries: for $22$ transformations $T$ of the local coordinates --- the mirrors $M_x,M_y,M_z$, the rotations $R_x,R_y,R_z$ by $180^\circ,120^\circ,90^\circ,72^\circ,60^\circ,37^\circ$ and the point reflection $P$ through the centre --- the Jaccard index $J=|A\cap T(A)|/|A\cup T(A)|$ of the body $A$, with a $95\%$ Poisson-bootstrap interval. An exact symmetry of algebraic origin maps orbits to orbits isometrically and gives $J=1$ up to the few points whose escape is decided by rounding; an approximate one gives a value that stays below $1$ as the sample grows. The table lists as \emph{exact} the transformations with $J\geq0.9995$ and as \emph{approximate} those with $0.97\leq J<0.9995$, with the smallest index of each group; ``all'' means that the rotations by every tested angle, including $37^\circ$, are exact.
\item \emph{Dim}: the dimension of the subalgebra generated by a generic point of the section and the constants of the map --- for the quadratic maps the subalgebra of the symmetrized product $x\circ y=\frac12(x\diamond y+y\diamond x)$, for the other maps that of $\diamond$ --- which contains every orbit. It is an exact invariant, computed by linear algebra on the structure constants, and it replaces the \emph{orbit span} of the first version, which collapsed when an orbit converged to an attracting cycle.
\end{itemize}
The threshold $0.97$ of the first version, on about $2400$ points, could not separate the two cases: for the constant of Plate~A07 the rotation about $jl$ of the section $(j,l,jl)$ scores about $0.97$, but it is not a symmetry, and it becomes one only when $\arg c$ reaches the Weyl wall at $120^\circ$ (Plates~R06 and R08). The measurements are, however, remarkably consistent with the propositions of Section~\ref{atl:sec:algebras}; for instance the sandwich balls of Plate~I02 have the radius predicted by Proposition~\ref{atl:prop:symmetric}.

\subsubsection{Corrections with respect to the first version}
The first version of this atlas was examined by a referee, whose report led to the main text and to the following corrections.
\begin{itemize}
\item The claim that the Okubo algebra ``loses every continuous symmetry'' and that ``only the half-turns'' of the torus leave a coordinate section invariant was wrong: the torus acts by rotations on the three root planes, and the sections containing a root plane are bodies of revolution (Proposition~\ref{atl:prop:torus}, Part~R).
\item The universal symmetry $z\mapsto-z$ was not stated; many mirrors and half-turns of the tables are this symmetry composed with an element of the torus (Remark~\ref{atl:rem:minus}).
\item The mirror/rotation pattern of the imaginary quarter of Plate~B07 came from scores oscillating about the threshold $0.97$; the quarter is nearly, but not exactly, round, and the explanation by an ``off-centre sphere'' was withdrawn.
\item Plate~J01 was rendered with $16$ iterations and Plate~J02 with $20$; the difference of their body shares came from the iteration count only. Plate~R17 repeats J01 with $20$ iterations: it coincides with J02.
\item The unexplained exact revolution of Plate~C05 is explained by the twisted equivariance of Proposition~\ref{atl:prop:theta90}.
\item All commentaries were checked again against the images, and the descriptions not supported by them were corrected.
\end{itemize}

\subsubsection{How to read a plate}
Each plate consists of the image, a line with the map, the algebra, the constants, the frame, the half-extent, the number of iterations, the coloring and the mapping, the table of measurements and a commentary. The commentary always follows the same order: what is seen in the four quarters, which quarters coincide or differ, and which property of the algebra accounts for it, in comparison with the octonionic case whenever this is meaningful.
\clearpage

\subsection{Part A. The quadratic map with a complex constant}\label{atl:part:A}
We begin with the simplest possible experiment: the quadratic map $z\mapsto z\diamond z+c$ with a complex constant $c$, drawn in the frame $F_A=(1,i,j\,|\,l)$ in all eight preset algebras. The constant is taken on the same complex direction $d_{\CC}$ in every algebra, at $97\%$ of the distance of the boundary of the connectedness locus. Since $c\in\CC$ lies in every one of these algebras, all the differences between the plates of this part come from the product alone. The top-left quarter, the quaternionic section $(1,i,j)$, is a built-in control: the quaternions are a subalgebra of every algebra of this part, with the product $pq$ in the octonions, split-octonions ($l^2=+1$) and semi-octonions, and the product $\bar p\,\bar q$ in the para-octonions and in the four Okubo-type algebras. The quarter is therefore shared by Plates~A01, A02, A04 on one side and by Plates~A05--A08 on the other. The other three quarters contain the direction $l$ and are where the algebras reveal themselves.

\subsubsection{A01. Octonions: the reference plate}\label{atl:pl:A01}\mbox{}\par\vspace{0.5mm}
\noindent\includegraphics[width=\textwidth]{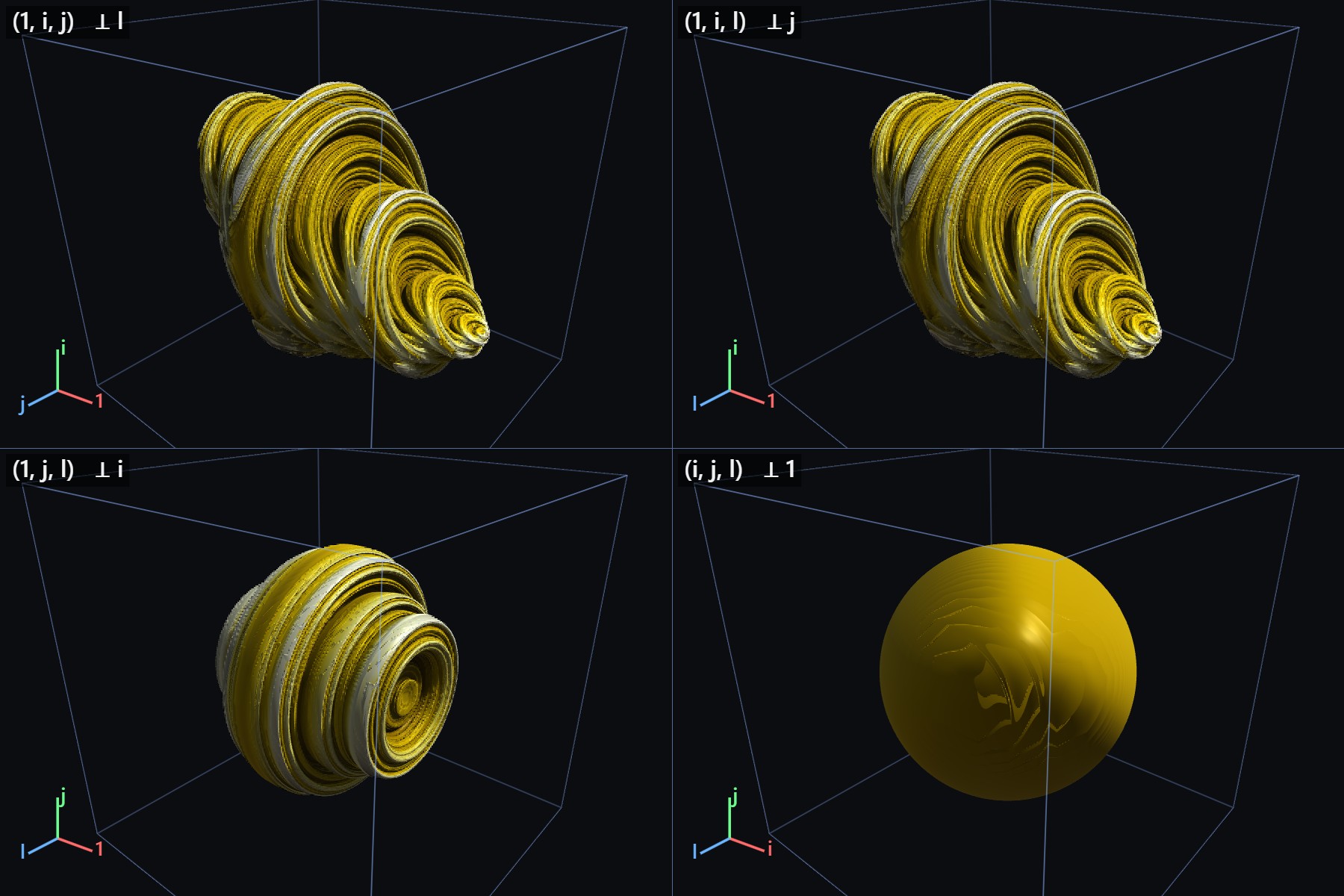}\par\vspace{1mm}
\noindent{\footnotesize \textbf{Map} $z\mapsto z\diamond z+c$ $\cdot$ \textbf{Algebra} Octonions (Cayley--Dickson, $\gamma=(-1,-1,-1)$) $\cdot$ $c=(-0.231,\,\allowbreak 0.62,\,\allowbreak 0,\,\allowbreak 0,\,\allowbreak 0,\,\allowbreak 0,\,\allowbreak 0,\,\allowbreak 0)$ ($=0.97\,s^*$, $s^*=0.6824$) $\cdot$ \textbf{Frame} $X,Y,Z,N=1,i,j,l$ $\cdot$ half-extent 1.4 $\cdot$ 20 iterations $\cdot$ \textbf{Coloring} 8D orbit trap (point), gradient Gold}\par\vspace{1mm}
\noindent{\footnotesize\begin{tabular}{@{}lllll@{}}\toprule Slice (X, Y, Z) & Body [95\%] & Exact ($J\geq0.9995$) & Approximate ($0.97\leq J<0.9995$) & Dim\\\midrule
($1, i, j$) $\perp l$ & 24.7\% {\scriptsize[24.5, 25.0]} & $M$: z; $R_z$: 180; $P$ & -- & 3 \\
($1, i, l$) $\perp j$ & 24.7\% {\scriptsize[24.5, 25.0]} & $M$: z; $R_z$: 180; $P$ & -- & 3 \\
($1, j, l$) $\perp i$ & 22.8\% {\scriptsize[22.5, 23.0]} & $M$: xyz; $R_x$: all; $R_y$: 180; $R_z$: 180; $P$ & -- & 3 \\
($i, j, l$) $\perp 1$ & 35.1\% {\scriptsize[34.8, 35.4]} & $M$: xyz; $R_x$: all; $R_y$: all; $R_z$: all; $P$ & -- & 3 \\
\bottomrule\end{tabular}}\par\vspace{2mm}
This is the reference plate of the atlas, and it is worth reading slowly. The constant is $c=-0.231+0.620\,i$, at $97\%$ of the boundary distance $s^*=0.682$ on the complex ray. The two upper quarters, $(1,i,j)$ and $(1,i,l)$, show the same body: an elongated, conch-like shell covered by fine spiral ridges, tapering to a small whorl at one end. The coincidence is not an accident of the picture: by Proposition~\ref{atl:prop:equivariance} the map $z^2+c$ commutes with every orthogonal transformation that fixes $1$ and $i$, in particular with the one exchanging $j$ and $l$. The measurements agree: the two quarters have the same body share ($25\%$) and the same exact symmetries, the mirror $M_z$, the half-turn $R_z(180^\circ)$ and the point reflection $P$. The last one is the universal symmetry $z\mapsto-z$ of Remark~\ref{atl:rem:minus}, and the mirror is its composition with the half-turn.

The lower-left quarter $(1,j,l)$ is orthogonal to $c$. Here the same proposition gives a whole circle of symmetries, the rotations of the plane of $j$ and $l$, and the body ($23\%$) is a body of revolution about the axis of $1$: the rotations about the local $x$ axis are exact for every tested angle, including $37^\circ$, and the ridges of the picture are circles around that axis, ending in a concentric whorl at the tip. The lower-right quarter $(i,j,l)$ consists of purely imaginary points, for which $z^2=-\n(z)$ (Proposition~\ref{atl:prop:imaginary}); the first iterate depends only on the norm, and the section is a ball ($35\%$), exact under every tested rotation.

All orbits lie in a subalgebra of dimension $3$, the span of $1$, $\im z_0$ and $\im c$ (Proposition~\ref{atl:prop:span}), as dictated by the quadratic identity $z^2=2\re(z)z-\n(z)$. This is the octonionic signature that the rest of the atlas will distort: identical quarters, a body of revolution and a ball.
\clearpage

\subsubsection{A02. Split-octonions with \texorpdfstring{$l^2=+1$}{l\texttwosuperior{}=+1}}\label{atl:pl:A02}\mbox{}\par\vspace{0.5mm}
\noindent\includegraphics[width=\textwidth]{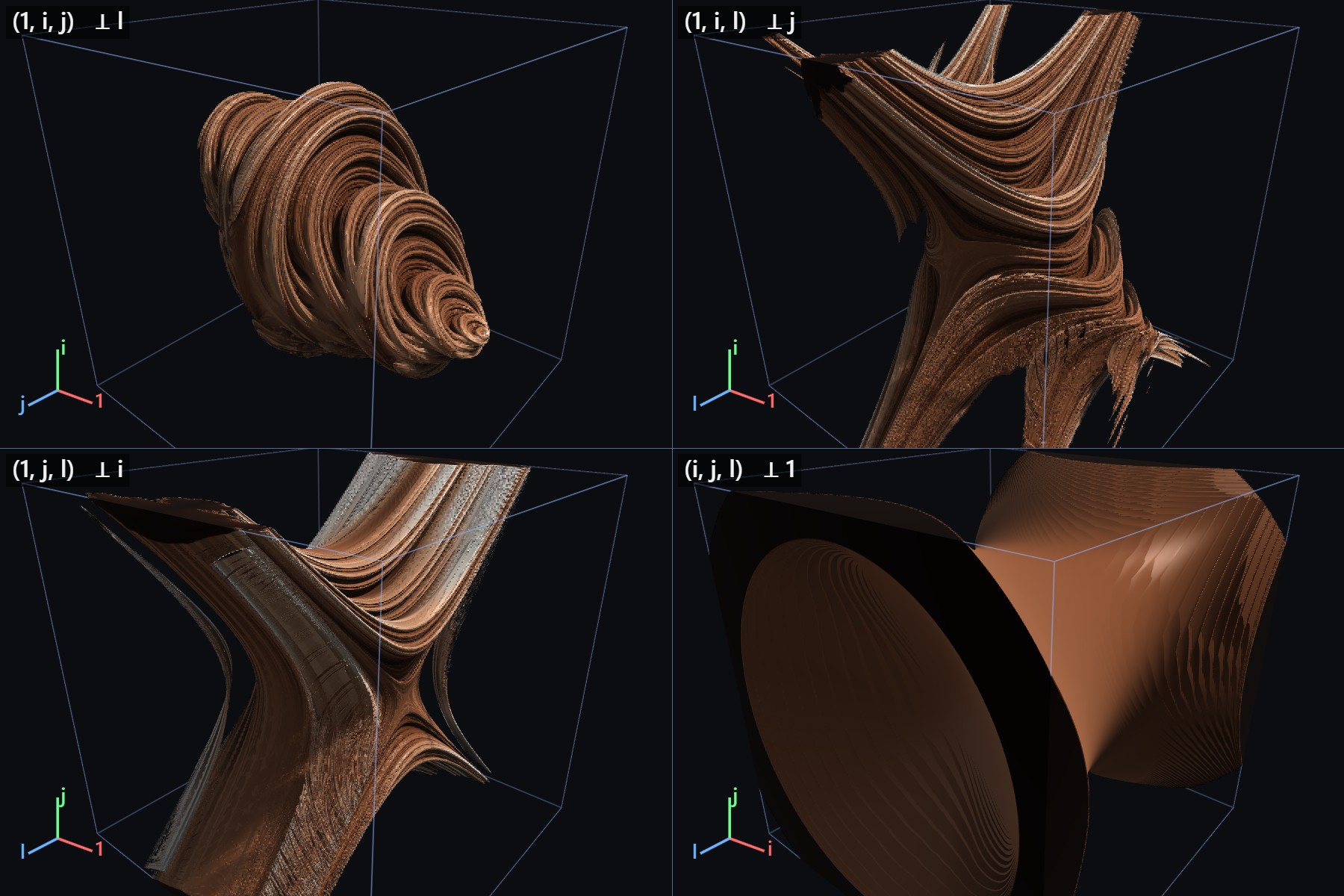}\par\vspace{1mm}
\noindent{\footnotesize \textbf{Map} $z\mapsto z\diamond z+c$ $\cdot$ \textbf{Algebra} Split-octonions ($l^2=+1$) (Cayley--Dickson, $\gamma=(-1,-1,1)$) $\cdot$ $c=(-0.231,\,\allowbreak 0.62,\,\allowbreak 0,\,\allowbreak 0,\,\allowbreak 0,\,\allowbreak 0,\,\allowbreak 0,\,\allowbreak 0)$ ($=0.97\,s^*$, $s^*=0.6824$) $\cdot$ \textbf{Frame} $X,Y,Z,N=1,i,j,l$ $\cdot$ half-extent 1.4 $\cdot$ 20 iterations $\cdot$ \textbf{Coloring} 8D orbit trap (point), gradient Copper}\par\vspace{1mm}
\noindent{\footnotesize\begin{tabular}{@{}lllll@{}}\toprule Slice (X, Y, Z) & Body [95\%] & Exact ($J\geq0.9995$) & Approximate ($0.97\leq J<0.9995$) & Dim\\\midrule
($1, i, j$) $\perp l$ & 24.8\% {\scriptsize[24.5, 25.0]} & $M$: z; $R_z$: 180; $P$ & -- & 3 \\
($1, i, l$) $\perp j$ & 27.4\% {\scriptsize[27.2, 27.7]} & $M$: z; $R_z$: 180; $P$ & -- & 3 \\
($1, j, l$) $\perp i$ & 30.2\% {\scriptsize[30.0, 30.5]} & $M$: xyz; $R_x$: 180; $R_y$: 180; $R_z$: 180; $P$ & -- & 3 \\
($i, j, l$) $\perp 1$ & 59.6\% {\scriptsize[59.3, 59.9]} & $M$: xyz; $R_x$: 180; $R_y$: 180; $R_z$: all; $P$ & -- & 3 \\
\bottomrule\end{tabular}}\par\vspace{2mm}
Changing only the sign of $l^2$ turns the octonions into the split-octonions, and the plate shows how differently the two halves of the algebra react. The upper-left quarter $(1,i,j)$ is the shell of Plate~A01, with the same body share ($25\%$): the quaternions are a common subalgebra, and the constant, which lies in $\CC$, is the same.

The three quarters containing $l$ are completely transformed. In $(1,i,l)$ and $(1,j,l)$ the shell is replaced by large bodies made of curved, finely striated sheets with a saddle-shaped profile, which leave the cube ($27\%$ and $30\%$ of the ball). The norm \eqref{atl:eq:norm} now contains $-x_4^2$ with the opposite sign, and its level sets, spheres in the compact case, become hyperboloids; the group of Proposition~\ref{atl:prop:equivariance} is now pseudo-orthogonal, and its circles become hyperbolas. Accordingly, the rotations about the axis of $1$ are no longer symmetries of $(1,j,l)$: only the three mirrors, the three half-turns and $P$ are exact. The lower-right quarter $(i,j,l)$ is a large spool bounded by a hyperboloid of one sheet whose axis is $l$, cut by the faces of the cube in a circular section ($60\%$): by Proposition~\ref{atl:prop:imaginary} the section is a region between level sets of $x_1^2+x_2^2-x_4^2$, and indeed the rotations about the local $z$ axis, the direction of $l$, are exact for every tested angle.

The contrast with Plate~A01 is a picture of the zero divisors: on the null cone $\n(z)=0$ the square $z^2=2\re(z)z$ loses the term $-\n(z)$ that keeps the compact bodies bounded, and the bodies extend beyond the cube along directions close to it. For the split algebras the Euclidean bailout of the renderer is not an escape criterion, so the bodies drawn here are approximations of the filled set.
\clearpage

\subsubsection{A03. Split-octonions with \texorpdfstring{$i^2=+1$}{i\texttwosuperior{}=+1}}\label{atl:pl:A03}\mbox{}\par\vspace{0.5mm}
\noindent\includegraphics[width=\textwidth]{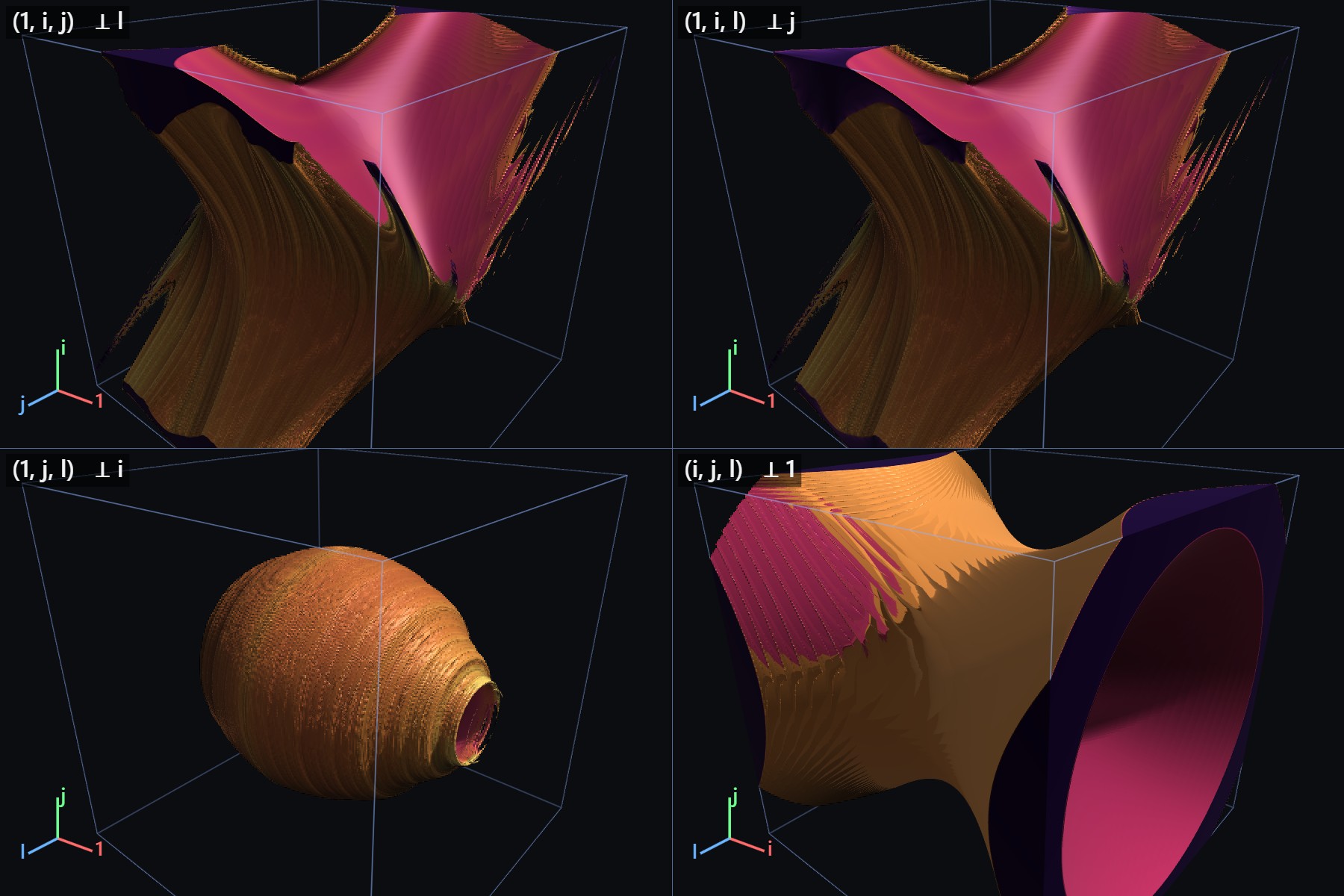}\par\vspace{1mm}
\noindent{\footnotesize \textbf{Map} $z\mapsto z\diamond z+c$ $\cdot$ \textbf{Algebra} Split-octonions ($i^2=+1$) (Cayley--Dickson, $\gamma=(1,-1,-1)$) $\cdot$ $c=(-0.144,\,\allowbreak 0.386,\,\allowbreak 0,\,\allowbreak 0,\,\allowbreak 0,\,\allowbreak 0,\,\allowbreak 0,\,\allowbreak 0)$ ($=0.97\,s^*$, $s^*=0.4251$) $\cdot$ \textbf{Frame} $X,Y,Z,N=1,i,j,l$ $\cdot$ half-extent 1.4 $\cdot$ 20 iterations $\cdot$ \textbf{Coloring} 8D orbit trap (point), gradient Sunset}\par\vspace{1mm}
\noindent{\footnotesize\begin{tabular}{@{}lllll@{}}\toprule Slice (X, Y, Z) & Body [95\%] & Exact ($J\geq0.9995$) & Approximate ($0.97\leq J<0.9995$) & Dim\\\midrule
($1, i, j$) $\perp l$ & 42.8\% {\scriptsize[42.5, 43.1]} & $M$: z; $R_z$: 180; $P$ & -- & 3 \\
($1, i, l$) $\perp j$ & 42.8\% {\scriptsize[42.5, 43.1]} & $M$: z; $R_z$: 180; $P$ & -- & 3 \\
($1, j, l$) $\perp i$ & 37.0\% {\scriptsize[36.7, 37.3]} & $M$: xyz; $R_x$: all; $R_y$: 180; $R_z$: 180; $P$ & -- & 3 \\
($i, j, l$) $\perp 1$ & 57.4\% {\scriptsize[57.1, 57.7]} & $M$: xyz; $R_x$: all; $R_y$: 180; $R_z$: 180; $P$ & -- & 3 \\
\bottomrule\end{tabular}}\par\vspace{2mm}
The split-octonions with $i^2=+1$ are isomorphic to those of Plate~A02, but the isomorphism moves the complex plane. Here the plane of $1$ and $i$ is the plane of the hyperbolic (split-complex) numbers, where the quadratic map splits, along the null lines $1\pm i$, into two real quadratic maps, and the Mandelbrot set becomes a square \cite{Senn}. It is worth noticing that the boundary distance $s^*=0.425$ is exactly the value at which one of the two real parameters, $c_0+c_1$, reaches $1/4$; the constant is $c=-0.144+0.386\,i$.

The upper two quarters, $(1,i,j)$ and $(1,i,l)$, are again identical ($43\%$ each, with the same exact symmetries $M_z$, $R_z(180^\circ)$, $P$), since $j$ and $l$ both have square $-1$ and are exchanged by a norm-preserving map fixing $1$ and $i$. Nothing of the shell remains: the quarters show a large body bounded by smooth saddle-like sheets that meet along sharp creases and leave the cube. The lower-left quarter $(1,j,l)$, on which the norm is positive definite, is the only compact section of the plate: an egg-shaped body of revolution about the axis of $1$ ($37\%$, exact rotations about the local $x$ axis for every tested angle), with a crater at its visible tip. The lower-right quarter $(i,j,l)$ is a spool bounded by a hyperboloid whose axis is $i$, the hyperbolic direction ($57\%$, exact rotations about the local $x$ axis): the norm there is $-x_1^2+x_2^2+x_4^2$.

The plate shows how the same abstract algebra looks different depending on where the complex plane is placed: the pictures record the signature of the norm on each section and on the plane of the constant.
\clearpage

\subsubsection{A04. Semi-octonions: the degenerate direction}\label{atl:pl:A04}\mbox{}\par\vspace{0.5mm}
\noindent\includegraphics[width=\textwidth]{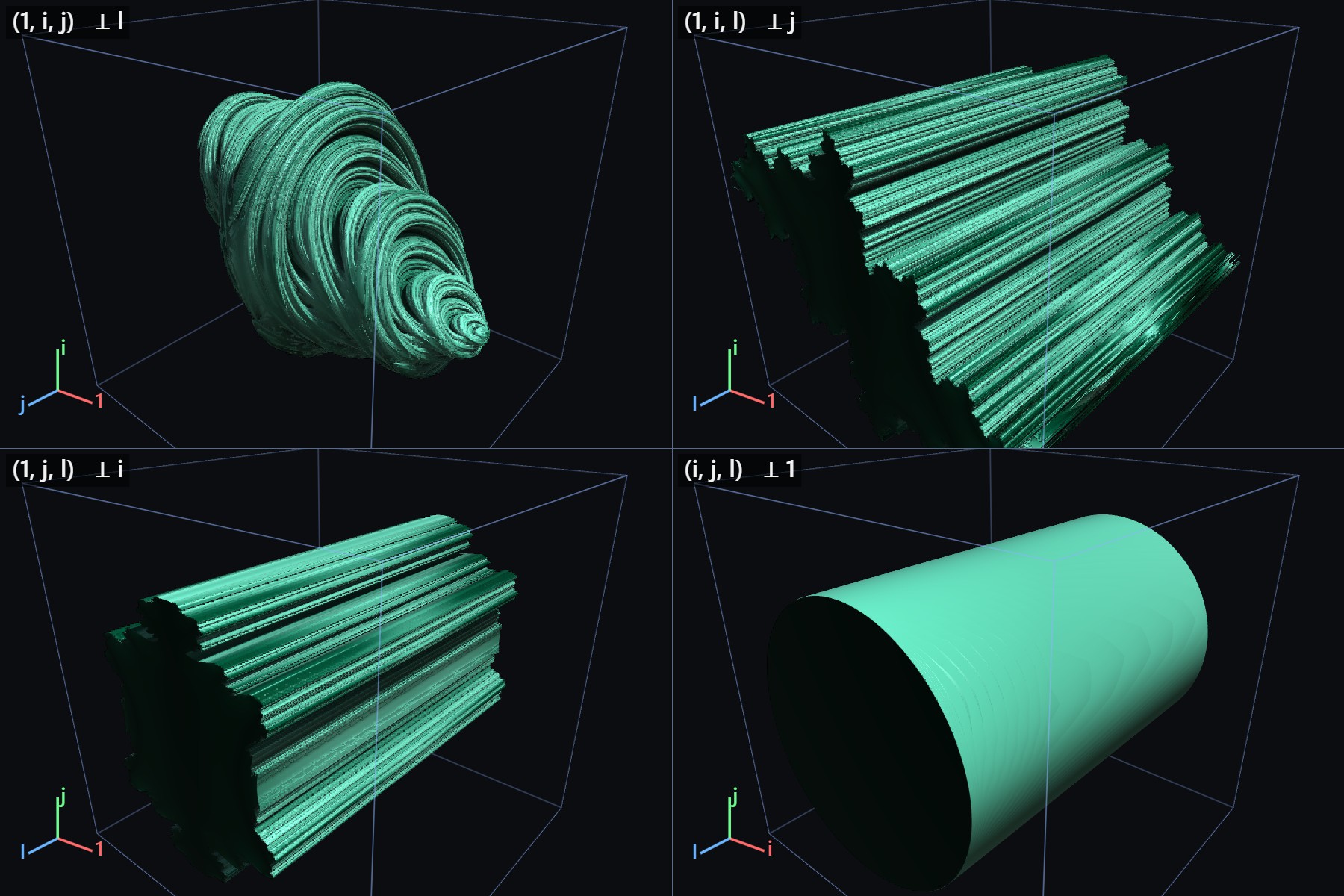}\par\vspace{1mm}
\noindent{\footnotesize \textbf{Map} $z\mapsto z\diamond z+c$ $\cdot$ \textbf{Algebra} Semi-octonions ($l^2=0$) (Cayley--Dickson, $\gamma=(-1,-1,0)$) $\cdot$ $c=(-0.231,\,\allowbreak 0.62,\,\allowbreak 0,\,\allowbreak 0,\,\allowbreak 0,\,\allowbreak 0,\,\allowbreak 0,\,\allowbreak 0)$ ($=0.97\,s^*$, $s^*=0.6824$) $\cdot$ \textbf{Frame} $X,Y,Z,N=1,i,j,l$ $\cdot$ half-extent 1.4 $\cdot$ 20 iterations $\cdot$ \textbf{Coloring} 8D orbit trap (point), gradient Emerald}\par\vspace{1mm}
\noindent{\footnotesize\begin{tabular}{@{}lllll@{}}\toprule Slice (X, Y, Z) & Body [95\%] & Exact ($J\geq0.9995$) & Approximate ($0.97\leq J<0.9995$) & Dim\\\midrule
($1, i, j$) $\perp l$ & 24.6\% {\scriptsize[24.3, 24.8]} & $M$: z; $R_z$: 180; $P$ & -- & 3 \\
($1, i, l$) $\perp j$ & 45.6\% {\scriptsize[45.3, 45.9]} & $M$: z; $R_z$: 180; $P$ & -- & 3 \\
($1, j, l$) $\perp i$ & 45.1\% {\scriptsize[44.8, 45.4]} & $M$: xyz; $R_x$: 180; $R_y$: 180; $R_z$: 180; $P$ & -- & 3 \\
($i, j, l$) $\perp 1$ & 64.3\% {\scriptsize[64.0, 64.6]} & $M$: xyz; $R_x$: 180; $R_y$: 180; $R_z$: all; $P$ & -- & 3 \\
\bottomrule\end{tabular}}\par\vspace{2mm}
The semi-octonions have $l^2=0$, a degenerate norm and a non-reductive automorphism group. The upper-left quarter $(1,i,j)$ is once again the quaternionic shell of Plate~A01 ($25\%$). The three quarters containing $l$, however, are cylinders along $l$ that cross the whole cube: in $(1,i,l)$ and $(1,j,l)$ they are fluted columns with jagged cut faces ($46\%$ and $45\%$), and in $(i,j,l)$ a smooth circular cylinder ($64\%$), for which the rotations about the local $z$ axis, the direction of $l$, are exact for every tested angle.

The explanation is Proposition~\ref{atl:prop:rescaling}: writing $z=p+ql$, one has $(p+ql)^2=p^2+2\re(p)\,q\,l$, so the quaternionic part evolves autonomously by $p\mapsto p^2+c$, and the $l$-part is transported linearly, $q_{n+1}=2\re(p_n)\,q_n$. On the hyperplane $\re p=0$ of the lower-right quarter the first iterate forgets $q$ entirely, and the body is an exact cylinder whose cross-section is the disc cut from the ball of Plate~A01 by the plane of $i$ and $j$. In the other two quarters $q_n$ is $q_0$ multiplied by the product of the factors $2\re p_k$ along the quaternionic orbit; where this product stays bounded every $q_0$ is allowed, and the body is a column along $l$ over the corresponding planar section of the quaternionic body.

This is perhaps the clearest picture of a degenerate direction: the dynamics does not see it. The orbits lie in a subalgebra of dimension $3$.
\clearpage

\subsubsection{A05. Para-octonions: the Mandelbar map in eight dimensions}\label{atl:pl:A05}\mbox{}\par\vspace{0.5mm}
\noindent\includegraphics[width=\textwidth]{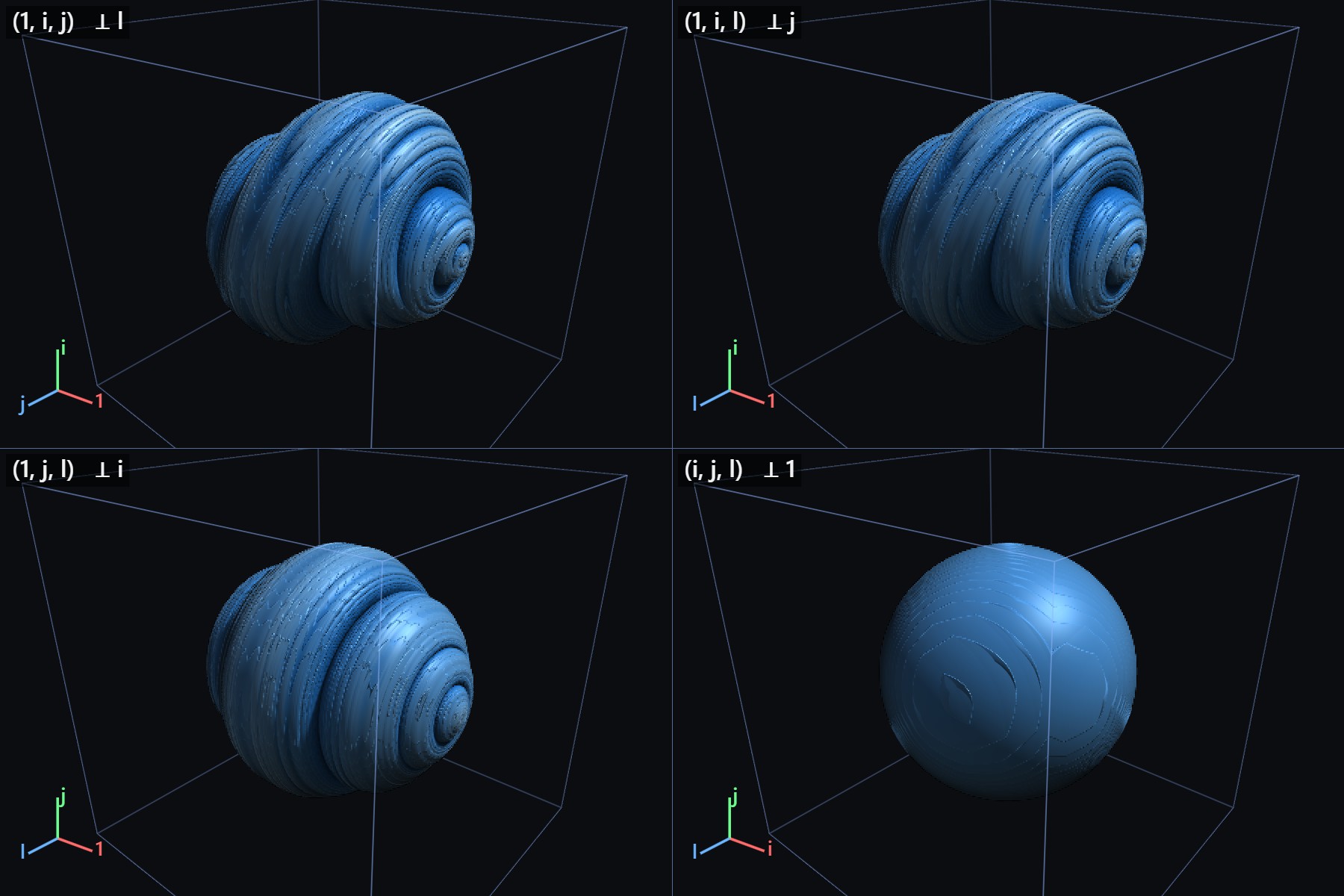}\par\vspace{1mm}
\noindent{\footnotesize \textbf{Map} $z\mapsto z\diamond z+c$ $\cdot$ \textbf{Algebra} Para-octonions (para-Hurwitz $\bar x\bar y$, $\gamma=(-1,-1,-1)$) $\cdot$ $c=(-0.086,\,\allowbreak 0.231,\,\allowbreak 0,\,\allowbreak 0,\,\allowbreak 0,\,\allowbreak 0,\,\allowbreak 0,\,\allowbreak 0)$ ($=0.97\,s^*$, $s^*=0.2540$) $\cdot$ \textbf{Frame} $X,Y,Z,N=1,i,j,l$ $\cdot$ half-extent 1.4 $\cdot$ 20 iterations $\cdot$ \textbf{Coloring} 8D orbit trap (point), gradient Deep Blue}\par\vspace{1mm}
\noindent{\footnotesize\begin{tabular}{@{}lllll@{}}\toprule Slice (X, Y, Z) & Body [95\%] & Exact ($J\geq0.9995$) & Approximate ($0.97\leq J<0.9995$) & Dim\\\midrule
($1, i, j$) $\perp l$ & 35.7\% {\scriptsize[35.4, 36.0]} & $M$: z; $R_z$: 180; $P$ & -- & 3 \\
($1, i, l$) $\perp j$ & 35.7\% {\scriptsize[35.4, 36.0]} & $M$: z; $R_z$: 180; $P$ & -- & 3 \\
($1, j, l$) $\perp i$ & 34.6\% {\scriptsize[34.3, 34.9]} & $M$: xyz; $R_x$: all; $R_y$: 180; $R_z$: 180; $P$ & -- & 3 \\
($i, j, l$) $\perp 1$ & 35.7\% {\scriptsize[35.4, 36.0]} & $M$: xyz; $R_x$: all; $R_y$: all; $R_z$: all; $P$ & -- & 3 \\
\bottomrule\end{tabular}}\par\vspace{2mm}
In the para-octonions the square is $z\bullet z=\overline{z^2}$, so that on the complex plane the map is the Mandelbar map $\bar z^2+c$ \cite{CHRS}. The constant is again complex, but on the same ray the boundary distance is much smaller, $s^*=0.254$ against $0.682$ in the octonions, and the constant is $c=-0.086+0.231\,i$: the connectedness locus of the Mandelbar family is a different set.

The symmetry pattern is exactly the octonionic one. The two upper quarters coincide ($36\%$ each, exact $M_z$, $R_z(180^\circ)$ and $P$); the lower-left quarter $(1,j,l)$ is a body of revolution about the axis of $1$ ($35\%$, exact rotations about the local $x$ axis for every tested angle); the lower-right quarter of purely imaginary points is a ball ($36\%$), because $\bar z=-z$ gives $z\bullet z=-\n(z)$ as in the octonions (Proposition~\ref{atl:prop:imaginary}). This is not surprising: the para-Hurwitz product is built from the octonionic one with the conjugation, which commutes with every orthogonal map fixing $1$, so Proposition~\ref{atl:prop:equivariance} holds verbatim. It is worth noticing that the three bodies through $1$ look almost alike, rounder and more compact than in Plate~A01, with spiral ridges ending in a whorl; only the measurements tell the body of revolution from the other two, whose rotations by $37^\circ$ score $0.87$.

The plate shows that the loss of the unit, by itself, does not break the symmetry: the orbits still lie in a subalgebra of dimension $3$, and the para-octonions are still an algebra of revolutions.
\clearpage

\subsubsection{A06. Split para-octonions}\label{atl:pl:A06}\mbox{}\par\vspace{0.5mm}
\noindent\includegraphics[width=\textwidth]{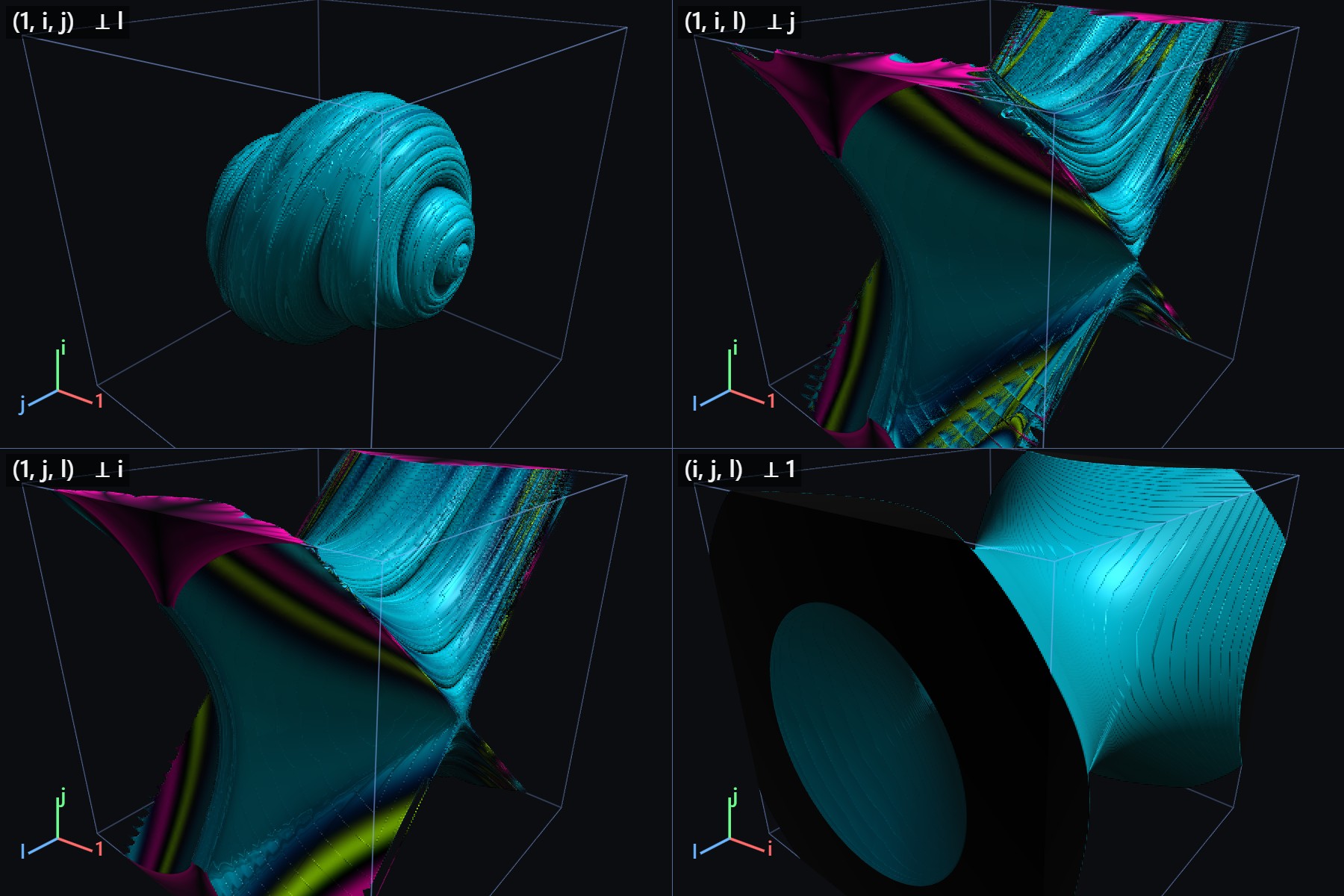}\par\vspace{1mm}
\noindent{\footnotesize \textbf{Map} $z\mapsto z\diamond z+c$ $\cdot$ \textbf{Algebra} Split para-octonions (para-Hurwitz $\bar x\bar y$, $\gamma=(-1,-1,1)$) $\cdot$ $c=(-0.086,\,\allowbreak 0.231,\,\allowbreak 0,\,\allowbreak 0,\,\allowbreak 0,\,\allowbreak 0,\,\allowbreak 0,\,\allowbreak 0)$ ($=0.97\,s^*$, $s^*=0.2540$) $\cdot$ \textbf{Frame} $X,Y,Z,N=1,i,j,l$ $\cdot$ half-extent 1.4 $\cdot$ 20 iterations $\cdot$ \textbf{Coloring} 8D orbit trap (point), gradient Neon}\par\vspace{1mm}
\noindent{\footnotesize\begin{tabular}{@{}lllll@{}}\toprule Slice (X, Y, Z) & Body [95\%] & Exact ($J\geq0.9995$) & Approximate ($0.97\leq J<0.9995$) & Dim\\\midrule
($1, i, j$) $\perp l$ & 35.8\% {\scriptsize[35.5, 36.1]} & $M$: z; $R_z$: 180; $P$ & -- & 3 \\
($1, i, l$) $\perp j$ & 52.5\% {\scriptsize[52.2, 52.8]} & $M$: z; $R_z$: 180; $P$ & -- & 3 \\
($1, j, l$) $\perp i$ & 52.8\% {\scriptsize[52.4, 53.1]} & $M$: xyz; $R_x$: 180; $R_y$: 180; $R_z$: 180; $P$ & -- & 3 \\
($i, j, l$) $\perp 1$ & 72.7\% {\scriptsize[72.4, 72.9]} & $M$: xyz; $R_x$: 180; $R_y$: 180; $R_z$: all; $P$ & -- & 3 \\
\bottomrule\end{tabular}}\par\vspace{2mm}
The split para-octonions combine the two previous deformations. The upper-left quarter repeats Plate~A05 ($36\%$), since the para-quaternions are a common subalgebra and the constant is the same. The quarters containing $l$ show hyperbolic sheets: saddle-shaped bodies in $(1,i,l)$ and $(1,j,l)$, bounded by smooth curved sheets and cut by the faces of the cube, and a spool in $(i,j,l)$, bounded by a hyperboloid of the split norm and cut by the cube in a circular section. The spool is an exact body of revolution about the $l$-axis (all tested rotations about the local $z$ axis), for the reason given in Plate~A02: by Proposition~\ref{atl:prop:imaginary} the imaginary section is a region between level sets of the norm.

The exact symmetries coincide quarter by quarter with those of Plate~A02, but the bodies are larger: $53\%$, $53\%$ and $73\%$ of the ball in the three quarters containing $l$, against $27\%$, $30\%$ and $60\%$. The constant is also much smaller in modulus ($c=-0.086+0.231\,i$ against $-0.231+0.620\,i$), which is presumably the main reason. The plate confirms that the para-Hurwitz construction transports the symmetries of the split-octonions, mirror by mirror; the orbits lie in a subalgebra of dimension $3$. As for all split algebras, the rendered bodies are approximations, the Euclidean bailout not being an escape criterion.
\clearpage

\subsubsection{A07. The Okubo algebra: no revolution in this frame}\label{atl:pl:A07}\mbox{}\par\vspace{0.5mm}
\noindent\includegraphics[width=\textwidth]{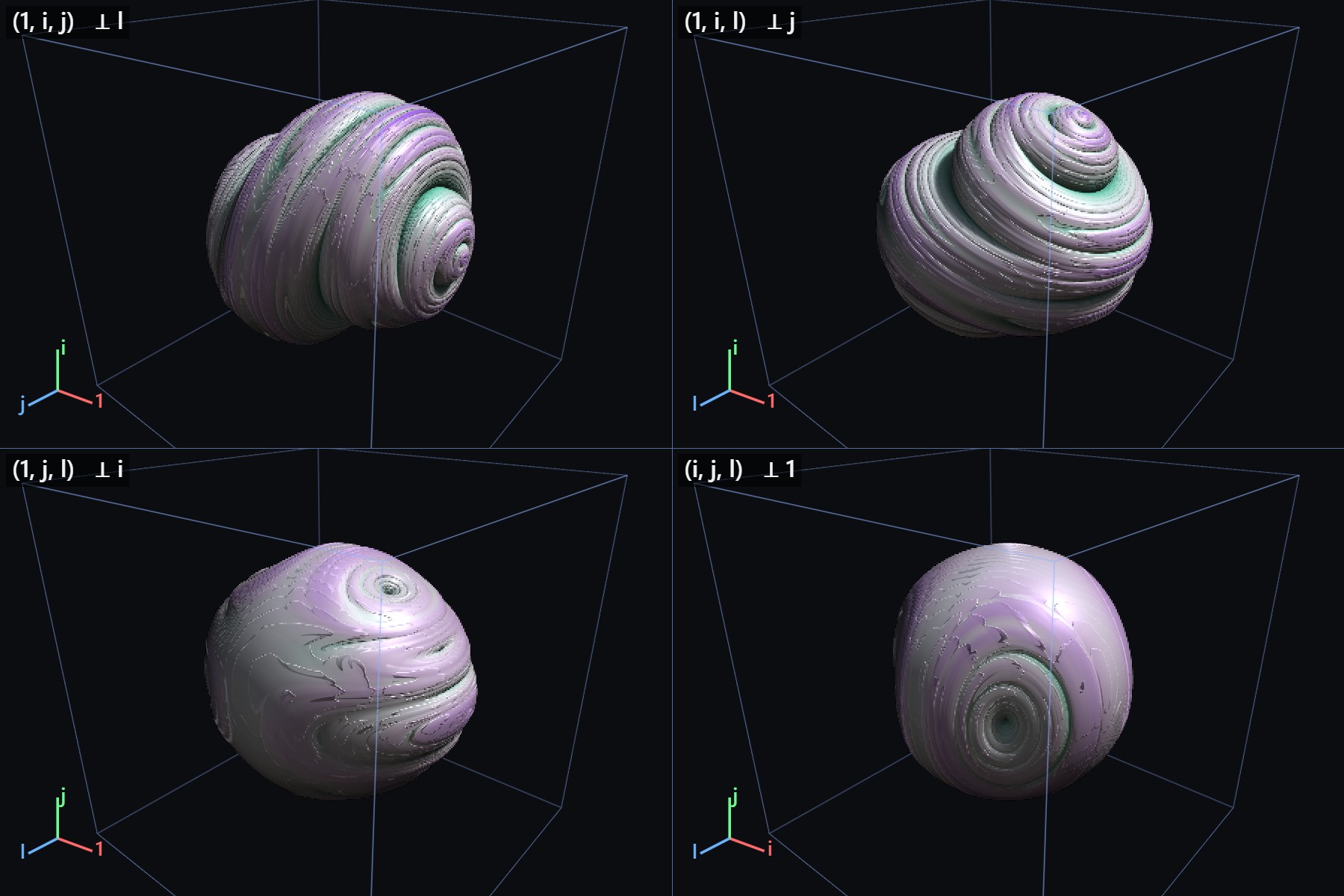}\par\vspace{1mm}
\noindent{\footnotesize \textbf{Map} $z\mapsto z\diamond z+c$ $\cdot$ \textbf{Algebra} Okubo algebra (Petersson, $\gamma=(-1,-1,-1)$, $\theta=120^\circ$) $\cdot$ $c=(-0.086,\,\allowbreak 0.231,\,\allowbreak 0,\,\allowbreak 0,\,\allowbreak 0,\,\allowbreak 0,\,\allowbreak 0,\,\allowbreak 0)$ ($=0.97\,s^*$, $s^*=0.2540$) $\cdot$ \textbf{Frame} $X,Y,Z,N=1,i,j,l$ $\cdot$ half-extent 1.4 $\cdot$ 20 iterations $\cdot$ \textbf{Coloring} 8D orbit trap (point), gradient Aurora}\par\vspace{1mm}
\noindent{\footnotesize\begin{tabular}{@{}lllll@{}}\toprule Slice (X, Y, Z) & Body [95\%] & Exact ($J\geq0.9995$) & Approximate ($0.97\leq J<0.9995$) & Dim\\\midrule
($1, i, j$) $\perp l$ & 35.9\% {\scriptsize[35.6, 36.2]} & $M$: z; $R_z$: 180; $P$ & -- & 3 \\
($1, i, l$) $\perp j$ & 34.6\% {\scriptsize[34.3, 34.9]} & $M$: z; $R_z$: 180; $P$ & -- & 3 \\
($1, j, l$) $\perp i$ & 37.3\% {\scriptsize[37.0, 37.6]} & $M$: xyz; $R_x$: 180; $R_y$: 180; $R_z$: 180; $P$ & -- & 5 \\
($i, j, l$) $\perp 1$ & 33.0\% {\scriptsize[32.7, 33.3]} & $M$: xyz; $R_x$: 180; $R_y$: 180; $R_z$: 180; $P$ & -- & 5 \\
\bottomrule\end{tabular}}\par\vspace{2mm}
This is the first Okubonic plate, and it should be compared quarter by quarter with Plates~A01 and A05. The constant is the same as in the para-octonions, $c=-0.086+0.231\,i$, because the critical orbit of a complex constant stays in $\CC$, where the Okubo product coincides with the para-octonionic one, $\bar z\,\bar w$. The upper-left quarter $(1,i,j)$ lies in $\HH$, on which every Petersson product is the para-quaternionic one (Proposition~\ref{atl:prop:torus}), and it is indeed the body of Plate~A05, with the same share ($36\%$) and the same exact symmetries.

The other quarters are new. The upper-right quarter $(1,i,l)$ is a similar ribbed shell ($35\%$), but turned: its spiral whorl faces upwards, and it is no longer a copy of the upper-left quarter, since the Okubo product does not enjoy the $O(6)$-equivariance that made the two quarters coincide in the octonions and para-octonions. The lower-left quarter $(1,j,l)$ is a rounded body with a small vortex on top ($37\%$), and the lower-right quarter of purely imaginary points is an egg-shaped body with a vortex facing the viewer ($33\%$): it is not a ball, because the identity $z*z=-\n(z)$ of Proposition~\ref{atl:prop:imaginary} fails.

In $(1,j,l)$ and $(i,j,l)$ the exact symmetries are the three mirrors, the three half-turns and $P$; no rotation by another angle is exact, and the best one, $R_x(37^\circ)$ in $(1,j,l)$, scores $0.957$: the body is nearly round but is not a body of revolution. This is exactly what Theorem~A predicts. The continuous symmetry of the octonions is reduced to the torus of Proposition~\ref{atl:prop:torus}, which fixes $c$ and rotates the root planes $V_1=\langle j,k\rangle$, $V_2=\langle l,il\rangle$, $V_3=\langle jl,kl\rangle$ by the angles \eqref{atl:eq:rootangles}. The frame $F_A$ contains no root plane (the section $(1,j,l)$ contains only one vector of $V_1$ and one of $V_2$), so no one-parameter subgroup maps its sections into themselves; what remains are the elements of the torus that turn $V_1$ or $V_2$ by $180^\circ$, acting as the mirrors $j\mapsto-j$ and $l\mapsto-l$, the anti-automorphism $\sigma$, and the universal symmetry $z\mapsto-z$ of Remark~\ref{atl:rem:minus}. The revolutions are not lost but displaced: the sections that contain a root plane are exact bodies of revolution (Plates~R01--R03).

The orbits started in $(1,j,l)$ and $(i,j,l)$ lie in a subalgebra of dimension $5$, those of the two upper quarters in a subalgebra of dimension $3$: the Okubo product mixes the quaternionic and the hidden half of the algebra. The plate is a small but complete picture of the difference between a para-unit and no unit at all.
\clearpage

\subsubsection{A08. The split Okubo algebra}\label{atl:pl:A08}\mbox{}\par\vspace{0.5mm}
\noindent\includegraphics[width=\textwidth]{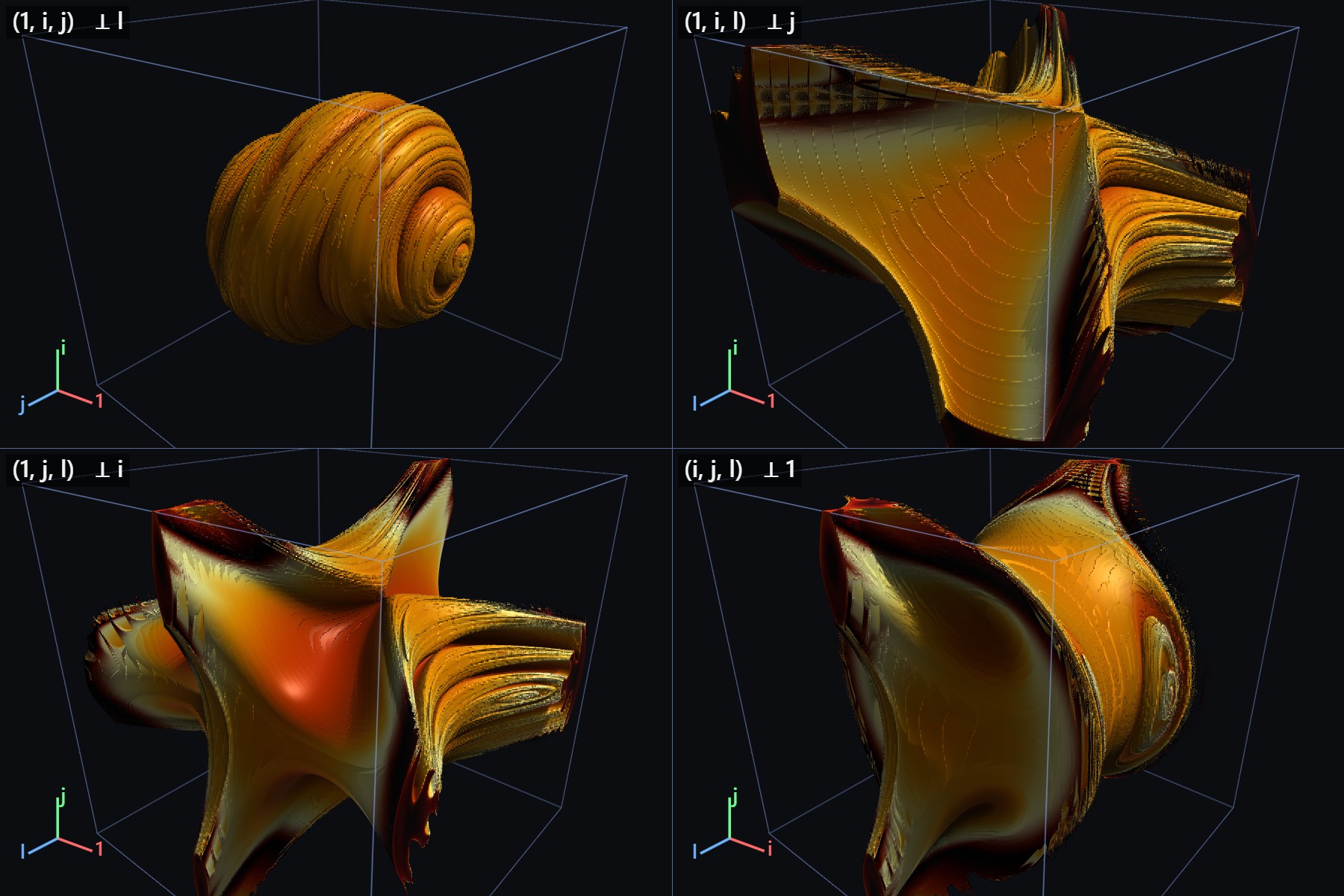}\par\vspace{1mm}
\noindent{\footnotesize \textbf{Map} $z\mapsto z\diamond z+c$ $\cdot$ \textbf{Algebra} Split Okubo algebra (Petersson, $\gamma=(-1,-1,1)$, $\theta=120^\circ$) $\cdot$ $c=(-0.086,\,\allowbreak 0.231,\,\allowbreak 0,\,\allowbreak 0,\,\allowbreak 0,\,\allowbreak 0,\,\allowbreak 0,\,\allowbreak 0)$ ($=0.97\,s^*$, $s^*=0.2540$) $\cdot$ \textbf{Frame} $X,Y,Z,N=1,i,j,l$ $\cdot$ half-extent 1.4 $\cdot$ 20 iterations $\cdot$ \textbf{Coloring} 8D orbit trap (point), gradient Fire}\par\vspace{1mm}
\noindent{\footnotesize\begin{tabular}{@{}lllll@{}}\toprule Slice (X, Y, Z) & Body [95\%] & Exact ($J\geq0.9995$) & Approximate ($0.97\leq J<0.9995$) & Dim\\\midrule
($1, i, j$) $\perp l$ & 35.9\% {\scriptsize[35.6, 36.2]} & $M$: z; $R_z$: 180; $P$ & -- & 3 \\
($1, i, l$) $\perp j$ & 54.5\% {\scriptsize[54.1, 54.8]} & $M$: z; $R_z$: 180; $P$ & -- & 3 \\
($1, j, l$) $\perp i$ & 45.6\% {\scriptsize[45.3, 45.9]} & $M$: xyz; $R_x$: 180; $R_y$: 180; $R_z$: 180; $P$ & -- & 5 \\
($i, j, l$) $\perp 1$ & 54.3\% {\scriptsize[54.0, 54.6]} & $M$: xyz; $R_x$: 180; $R_y$: 180; $R_z$: 180; $P$ & -- & 5 \\
\bottomrule\end{tabular}}\par\vspace{2mm}
The split Okubo algebra, with automorphism group $PSU(2,1)$, is to the Okubo algebra what the split para-octonions are to the para-octonions. The upper-left quarter repeats Plates~A05 and A07 ($36\%$, the same exact symmetries), as it should, since the constant is the same and the split and compact algebras share their quaternionic subalgebra.

The quarters containing $l$ show hyperbolic structures, as in Plate~A06, but of a different shape. The upper-right quarter $(1,i,l)$ is a large body ($55\%$) bounded by smooth curved sheets that meet along a sharp crease and are cut by the faces of the cube, with fibrous flanks on the right. The lower-left quarter $(1,j,l)$ ($46\%$) is a saddle-shaped body whose sheets reach out in several directions to the faces of the cube; the lower-right quarter $(i,j,l)$ ($54\%$) is made of curved sheets wrapped around a hollow, again cut by the cube. Unlike the spool of Plate~A06, this last quarter is not a body of revolution (no rotation other than the half-turns is exact; $R_z(37^\circ)$ scores $0.81$): as in the compact case, the identity $z*z=-\n(z)$ fails for imaginary $z$.

The exact symmetries are exactly those of Plate~A07: $M_z$, $R_z(180^\circ)$ and $P$ in the two upper quarters, the three mirrors, the three half-turns and $P$ in the two lower ones. The Petersson twist $\tau$, which multiplies the $l$-half of the algebra by $u=e^{2\pi i/3}$, and the split norm act together, and neither alone produces these shapes. The orbits lie in a subalgebra of dimension $3$ in the upper quarters and $5$ in the lower ones. As for all split algebras, the Euclidean bailout is not an escape criterion, and the bodies drawn are approximations of the filled set.
\clearpage

\subsection{Part B. The quadratic map with a generic eight-dimensional constant}\label{atl:part:B}
In this part the constant is a generic eight-dimensional vector, on the direction $d_{\mathcal O}$ found by the parameter scan of the prototype HF8D~0.2, at $95\%$ of the boundary distance. The frame is $F_B=(1,l,il\,|\,jl)$, entirely made of the ``hidden'' half of the algebra except for the unit. With a generic constant the octonionic symmetry of Proposition~\ref{atl:prop:equivariance} is reduced to the stabilizer of $\im c$, and the sections through $1$ are no longer bodies of revolution about a coordinate axis; nevertheless the octonionic orbits still lie in a three-dimensional subalgebra, while the Okubonic ones generate the whole algebra. In the Okubo algebra the stabilizer of a generic $c$ is again a two-dimensional torus, whose root planes are in general position with respect to the coordinate axes: the absence of revolutions in this part is an effect of the frame (Plate~R09).

\subsubsection{B01. Octonions with a generic constant}\label{atl:pl:B01}\mbox{}\par\vspace{0.5mm}
\noindent\includegraphics[width=\textwidth]{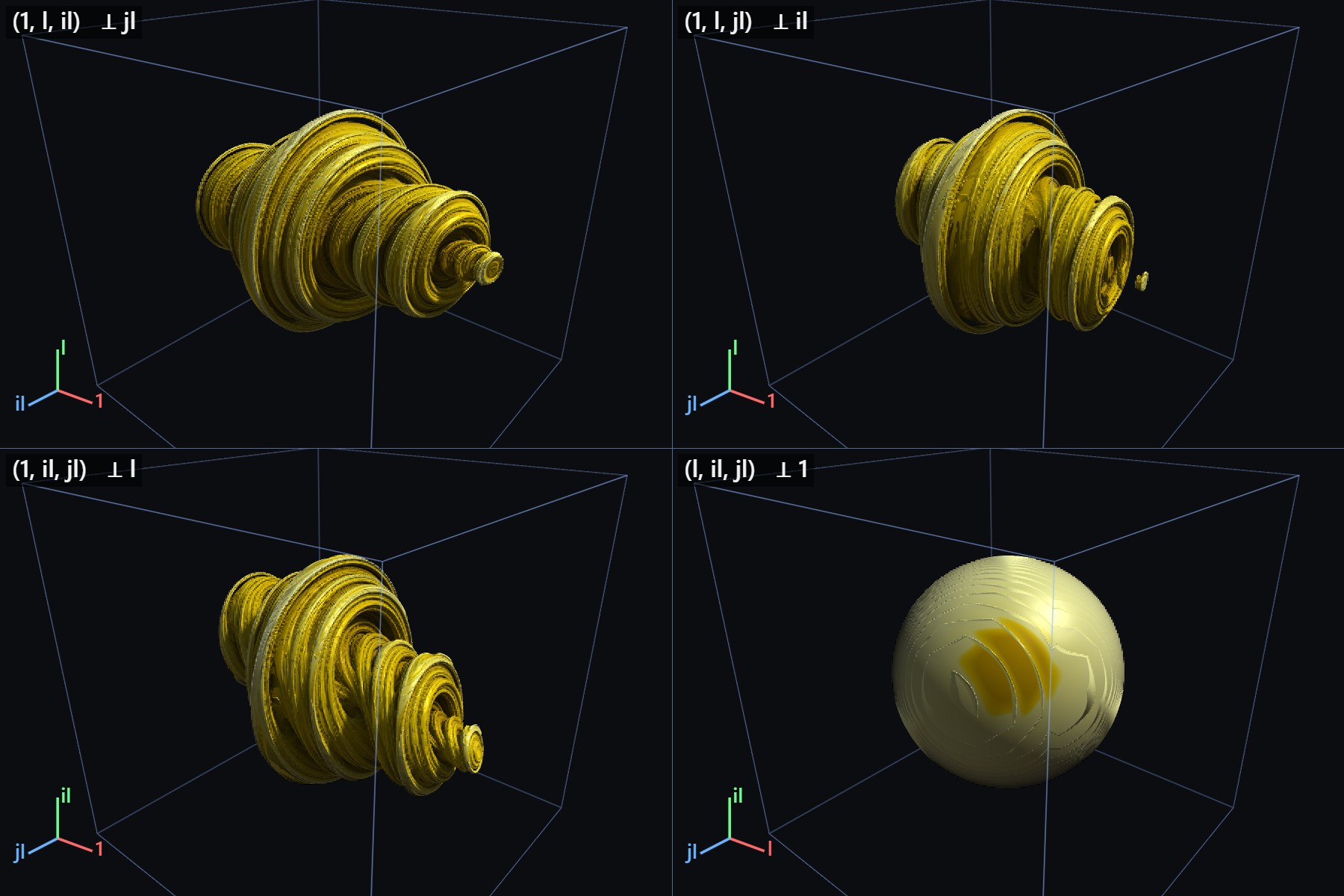}\par\vspace{1mm}
\noindent{\footnotesize \textbf{Map} $z\mapsto z\diamond z+c$ $\cdot$ \textbf{Algebra} Octonions (Cayley--Dickson, $\gamma=(-1,-1,-1)$) $\cdot$ $c=(-0.598,\,\allowbreak -0.032,\,\allowbreak 0.075,\,\allowbreak 0.025,\,\allowbreak 0.006,\,\allowbreak 0.241,\,\allowbreak -0.17,\,\allowbreak 0.2)$ ($=0.95\,s^*$, $s^*=0.7383$) $\cdot$ \textbf{Frame} $X,Y,Z,N=1,l,il,jl$ $\cdot$ half-extent 1.4 $\cdot$ 20 iterations $\cdot$ \textbf{Coloring} 8D orbit trap (point), gradient Gold}\par\vspace{1mm}
\noindent{\footnotesize\begin{tabular}{@{}lllll@{}}\toprule Slice (X, Y, Z) & Body [95\%] & Exact ($J\geq0.9995$) & Approximate ($0.97\leq J<0.9995$) & Dim\\\midrule
($1, l, il$) $\perp jl$ & 16.7\% {\scriptsize[16.5, 17.0]} & $P$ & $M$: y; $R_y$: 180 (min.\ 0.972) & 3 \\
($1, l, jl$) $\perp il$ & 16.3\% {\scriptsize[16.1, 16.5]} & $P$ & $M$: y; $R_y$: 180 (min.\ 0.974) & 3 \\
($1, il, jl$) $\perp l$ & 17.2\% {\scriptsize[16.9, 17.4]} & $P$ & -- & 3 \\
($l, il, jl$) $\perp 1$ & 21.9\% {\scriptsize[21.6, 22.1]} & $M$: xyz; $R_x$: all; $R_y$: all; $R_z$: all; $P$ & -- & 3 \\
\bottomrule\end{tabular}}\par\vspace{2mm}
With a generic constant, $c$ at $95\%$ of $s^*=0.738$ on the eight-dimensional ray $d_{\mathcal O}$, the octonionic map loses its large symmetry group in the coordinate sections, but it keeps its quadratic character. The three quarters through $1$, namely $(1,l,il)$, $(1,l,jl)$ and $(1,il,jl)$, show similar spiral shells with the characteristic conch profile ($17\%$, $16\%$ and $17\%$ of the ball). A small detached piece is visible near the tip of the second quarter. The lower-right quarter $(l,il,jl)$, orthogonal to $1$, is again a ball ($22\%$, exact under every tested rotation), because Proposition~\ref{atl:prop:imaginary} does not depend on the constant.

In the quarters through $1$ the only exact symmetry is the point reflection $P$, the universal symmetry $z\mapsto-z$ of Remark~\ref{atl:rem:minus}; the mirror $M_y$ and the half-turn $R_y(180^\circ)$ are approximate in the first two quarters ($J=0.972$ and $0.974$) and absent in the third. By Proposition~\ref{atl:prop:equivariance} the map still commutes with the group $O(6)$ of orthogonal maps fixing $1$ and $\im c$, but since $\im c$ is not a coordinate direction these maps move the coordinate sections out of themselves. Nevertheless every orbit lies in a subalgebra of dimension $3$, as predicted by Proposition~\ref{atl:prop:span}: the dynamics of the octonions is always three-dimensional, whatever the constant.
\clearpage

\subsubsection{B02. Split-octonions with a generic constant}\label{atl:pl:B02}\mbox{}\par\vspace{0.5mm}
\noindent\includegraphics[width=\textwidth]{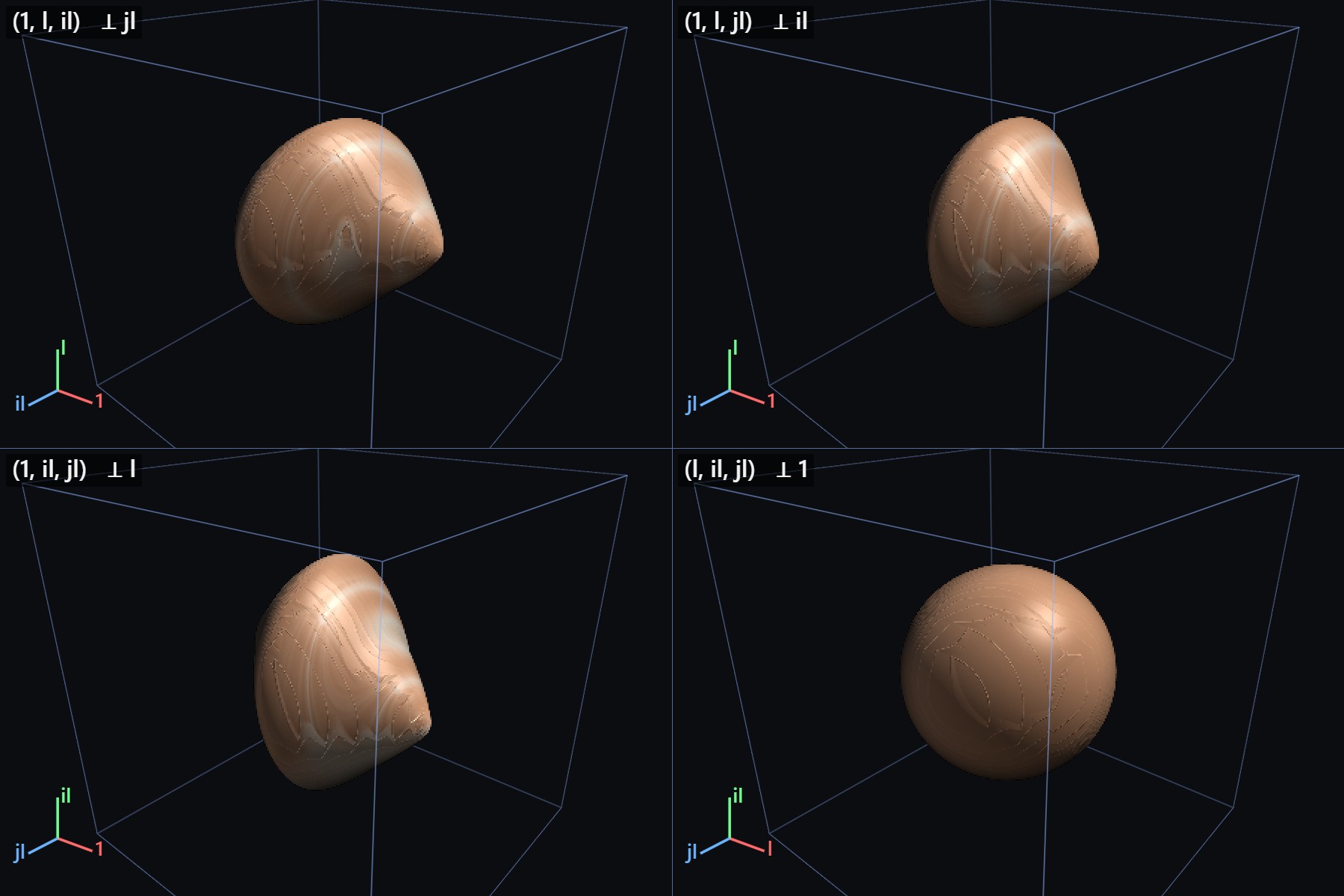}\par\vspace{1mm}
\noindent{\footnotesize \textbf{Map} $z\mapsto z\diamond z+c$ $\cdot$ \textbf{Algebra} Split-octonions ($l^2=+1$) (Cayley--Dickson, $\gamma=(-1,-1,1)$) $\cdot$ $c=(-1.203,\,\allowbreak -0.065,\,\allowbreak 0.15,\,\allowbreak 0.049,\,\allowbreak 0.012,\,\allowbreak 0.484,\,\allowbreak -0.342,\,\allowbreak 0.403)$ ($=0.95\,s^*$, $s^*=1.4856$) $\cdot$ \textbf{Frame} $X,Y,Z,N=1,l,il,jl$ $\cdot$ half-extent 2.4 $\cdot$ 20 iterations $\cdot$ \textbf{Coloring} 8D orbit trap (point), gradient Copper}\par\vspace{1mm}
\noindent{\footnotesize\begin{tabular}{@{}lllll@{}}\toprule Slice (X, Y, Z) & Body [95\%] & Exact ($J\geq0.9995$) & Approximate ($0.97\leq J<0.9995$) & Dim\\\midrule
($1, l, il$) $\perp jl$ & 15.9\% {\scriptsize[15.7, 16.1]} & $P$ & $M$: y; $R_y$: 180 (min.\ 0.991) & 3 \\
($1, l, jl$) $\perp il$ & 15.4\% {\scriptsize[15.2, 15.7]} & $P$ & $M$: y; $R_y$: 180 (min.\ 0.989) & 3 \\
($1, il, jl$) $\perp l$ & 16.4\% {\scriptsize[16.2, 16.7]} & $P$ & -- & 3 \\
($l, il, jl$) $\perp 1$ & 22.2\% {\scriptsize[21.9, 22.4]} & $M$: xyz; $R_x$: all; $R_y$: all; $R_z$: all; $P$ & -- & 3 \\
\bottomrule\end{tabular}}\par\vspace{2mm}
The split-octonions with a generic constant are drawn in a larger cube (half-extent $2.4$). The boundary distance on the ray is $s^*=1.486$, about twice the octonionic value, so that the constant is correspondingly larger. The three quarters through $1$ show smooth, rounded bodies without spiral ridges, narrowing towards a blunt tip ($16\%$, $15\%$ and $16\%$ of the ball). In this frame all three hidden directions $l,il,jl$ have square $+1$, and the norm restricted to these sections has signature $(1,2)$. The lower-right quarter, $(l,il,jl)$, is a ball ($22\%$, exact under every tested rotation): the section is purely imaginary, the norm is negative definite there, and Proposition~\ref{atl:prop:imaginary} applies with a change of sign.

The symmetry pattern is that of Plate~B01: $P$ is exact in every quarter, and the mirror $M_y$ with the half-turn $R_y(180^\circ)$ is approximate in the first two quarters ($J=0.991$ and $0.989$). The orbits lie in a subalgebra of dimension $3$, as in the compact octonions: the split-octonions are quadratic as well.
\clearpage

\subsubsection{B03. Split-octonions (\texorpdfstring{$i^2=+1$}{i\texttwosuperior{}=+1}) with a generic constant}\label{atl:pl:B03}\mbox{}\par\vspace{0.5mm}
\noindent\includegraphics[width=\textwidth]{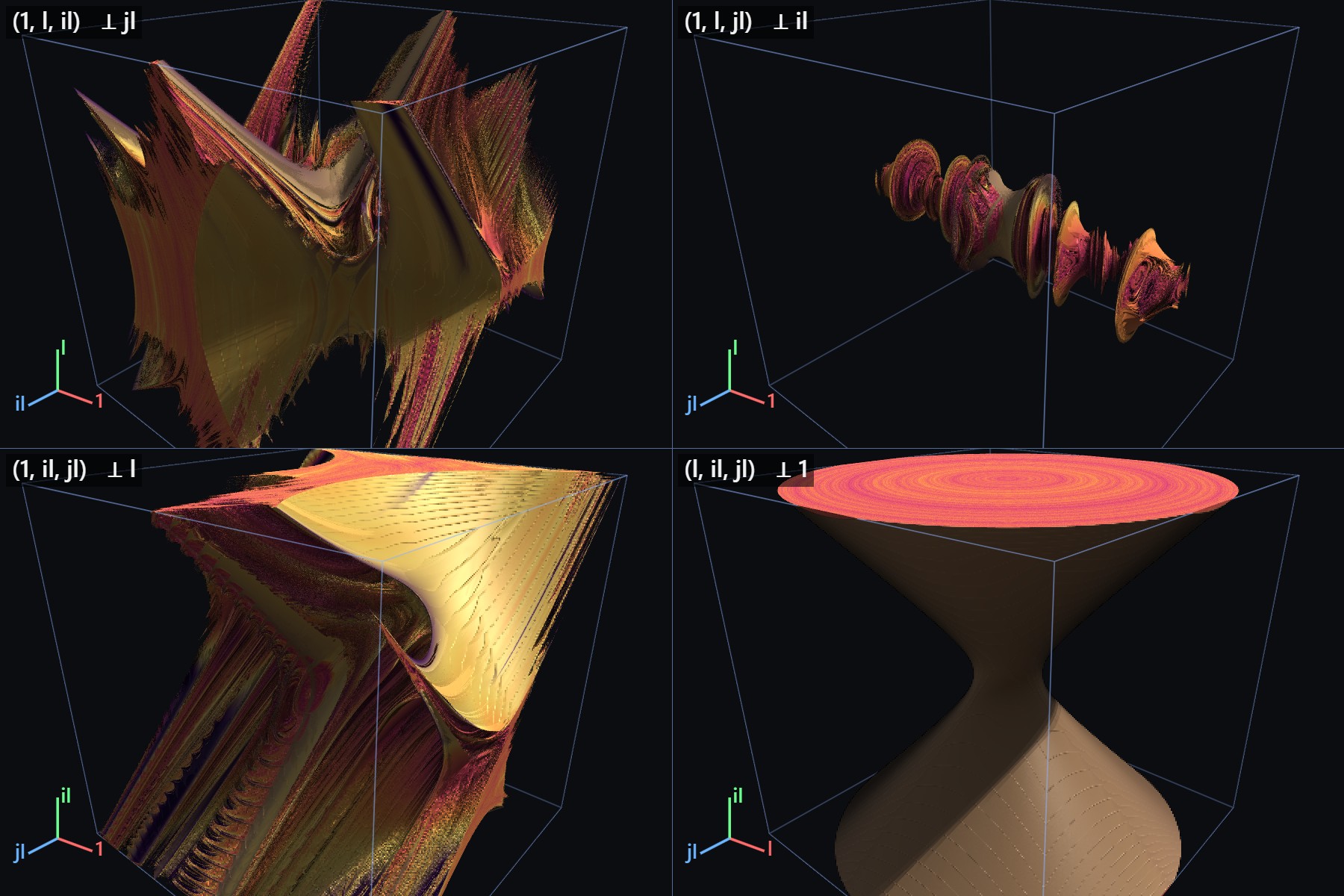}\par\vspace{1mm}
\noindent{\footnotesize \textbf{Map} $z\mapsto z\diamond z+c$ $\cdot$ \textbf{Algebra} Split-octonions ($i^2=+1$) (Cayley--Dickson, $\gamma=(1,-1,-1)$) $\cdot$ $c=(-1.332,\,\allowbreak -0.072,\,\allowbreak 0.166,\,\allowbreak 0.055,\,\allowbreak 0.013,\,\allowbreak 0.536,\,\allowbreak -0.379,\,\allowbreak 0.446)$ ($=0.95\,s^*$, $s^*=1.6438$) $\cdot$ \textbf{Frame} $X,Y,Z,N=1,l,il,jl$ $\cdot$ half-extent 1.4 $\cdot$ 20 iterations $\cdot$ \textbf{Coloring} 8D orbit trap (point), gradient Sunset}\par\vspace{1mm}
\noindent{\footnotesize\begin{tabular}{@{}lllll@{}}\toprule Slice (X, Y, Z) & Body [95\%] & Exact ($J\geq0.9995$) & Approximate ($0.97\leq J<0.9995$) & Dim\\\midrule
($1, l, il$) $\perp jl$ & 30.0\% {\scriptsize[29.7, 30.3]} & $P$ & -- & 3 \\
($1, l, jl$) $\perp il$ & 2.4\% {\scriptsize[2.3, 2.5]} & $P$ & -- & 3 \\
($1, il, jl$) $\perp l$ & 32.8\% {\scriptsize[32.5, 33.1]} & $P$ & -- & 3 \\
($l, il, jl$) $\perp 1$ & 33.4\% {\scriptsize[33.1, 33.7]} & $M$: xyz; $R_x$: 180; $R_y$: all; $R_z$: 180; $P$ & -- & 3 \\
\bottomrule\end{tabular}}\par\vspace{2mm}
In the split-octonions with $i^2=+1$ the frame $F_B$ mixes elliptic and hyperbolic directions: $l$ and $jl$ square to $-1$, $il$ to $+1$. The constant is large ($s^*=1.644$) and the picture reflects the signature unevenly. The quarter $(1,l,il)$ shows spiky hyperbolic sheets that leave the cube ($30\%$); the quarter $(1,l,jl)$, on which the norm is positive definite, is almost empty: a chain of small flared discs strung along a diagonal, $2.4\%$ of the ball; the quarter $(1,il,jl)$ is again a large body of curved sheets cut by the cube ($33\%$). The lower-right quarter $(l,il,jl)$ is an hourglass bounded by a hyperboloid of axis $il$ ($33\%$), exact under all tested rotations about the local $y$ axis: by Proposition~\ref{atl:prop:imaginary} it is a region between level sets of $x_4^2-x_5^2+x_6^2$.

In the three quarters through $1$ the only exact symmetry is $P$ (Remark~\ref{atl:rem:minus}), and no symmetry is even approximate. The orbits lie in a subalgebra of dimension $3$. The plate illustrates that, in a split algebra, the geometry of a section is decided by the signature of the norm on it rather than by the algebra as a whole: the same algebra of Plate~B02, in another presentation, gives different pictures. As for all split algebras, the rendered bodies are approximations.
\clearpage

\subsubsection{B04. Semi-octonions with a generic constant}\label{atl:pl:B04}\mbox{}\par\vspace{0.5mm}
\noindent\includegraphics[width=\textwidth]{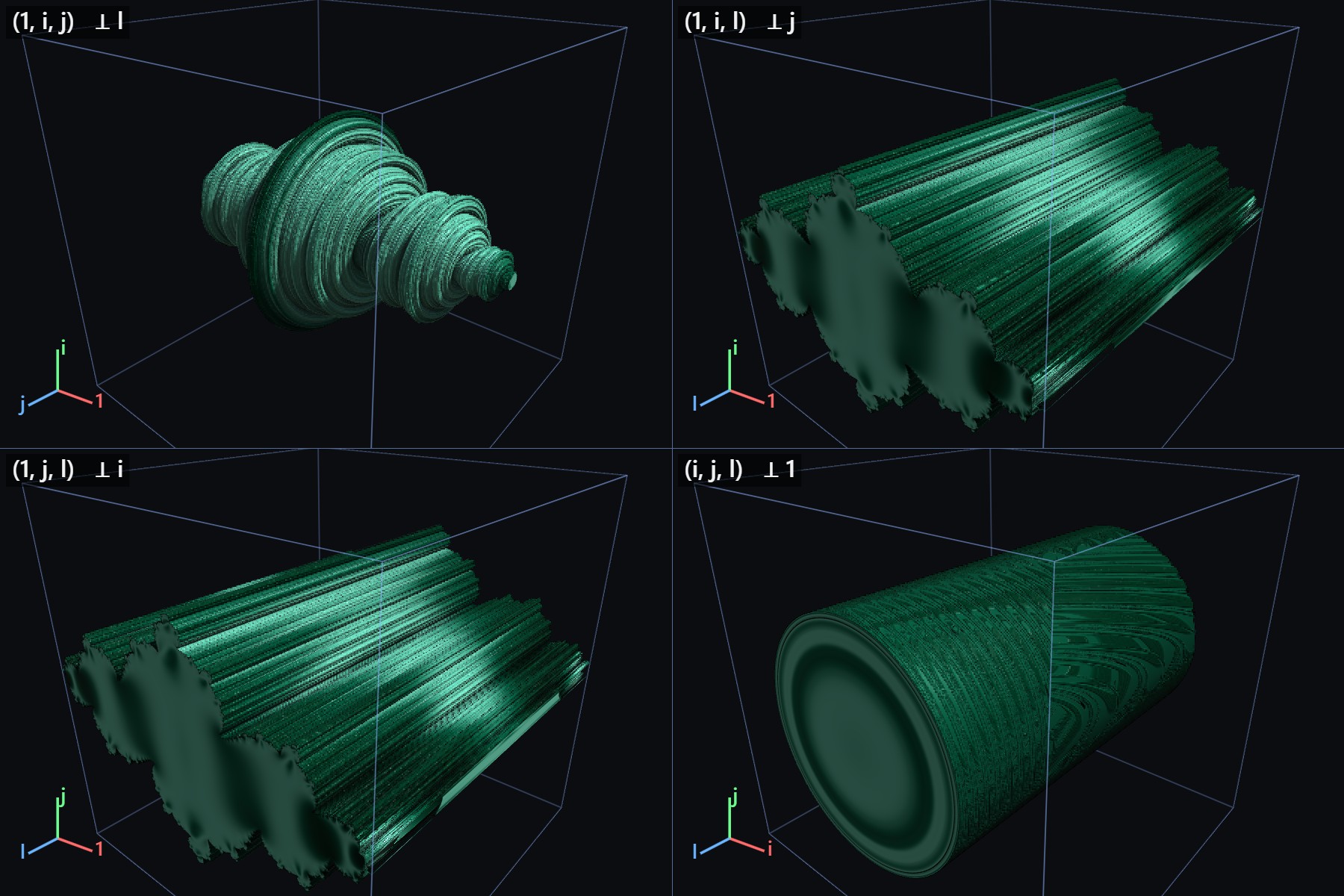}\par\vspace{1mm}
\noindent{\footnotesize \textbf{Map} $z\mapsto z\diamond z+c$ $\cdot$ \textbf{Algebra} Semi-octonions ($l^2=0$) (Cayley--Dickson, $\gamma=(-1,-1,0)$) $\cdot$ $c=(-0.705,\,\allowbreak -0.038,\,\allowbreak 0.088,\,\allowbreak 0.029,\,\allowbreak 0.007,\,\allowbreak 0.284,\,\allowbreak -0.201,\,\allowbreak 0.236)$ ($=0.95\,s^*$, $s^*=0.8708$) $\cdot$ \textbf{Frame} $X,Y,Z,N=1,i,j,l$ $\cdot$ half-extent 1.4 $\cdot$ 20 iterations $\cdot$ \textbf{Coloring} 8D orbit trap (point), gradient Emerald}\par\vspace{1mm}
\noindent{\footnotesize\begin{tabular}{@{}lllll@{}}\toprule Slice (X, Y, Z) & Body [95\%] & Exact ($J\geq0.9995$) & Approximate ($0.97\leq J<0.9995$) & Dim\\\midrule
($1, i, j$) $\perp l$ & 17.5\% {\scriptsize[17.2, 17.7]} & $P$ & -- & 3 \\
($1, i, l$) $\perp j$ & 49.1\% {\scriptsize[48.8, 49.4]} & $M$: z; $R_z$: 180; $P$ & -- & 3 \\
($1, j, l$) $\perp i$ & 49.0\% {\scriptsize[48.7, 49.3]} & $M$: z; $R_z$: 180; $P$ & -- & 3 \\
($i, j, l$) $\perp 1$ & 53.8\% {\scriptsize[53.5, 54.1]} & $M$: xyz; $R_x$: 180; $R_y$: 180; $R_z$: all; $P$ & -- & 3 \\
\bottomrule\end{tabular}}\par\vspace{2mm}
For the semi-octonions the frame $F_B$ is unusable, because every quarter contains an $l$-direction along which the body is a cylinder cut by the cube; we therefore use the frame $F_A$. The upper-left quarter $(1,i,j)$ is the quaternionic Julia set of the quaternionic part of the constant, a small ribbed shell ($17.5\%$) with $P$ as its only exact symmetry, since the imaginary part of $c_{\HH}$ is not a coordinate direction. The quarters $(1,i,l)$ and $(1,j,l)$ are fluted columns along $l$ that cross the whole cube, with jagged cut faces ($49\%$ each; exact $M_z$, $R_z(180^\circ)$ and $P$), and the lower-right quarter $(i,j,l)$ is a circular cylinder along $l$ with a ringed end face ($54\%$), exact under every tested rotation about the local $z$ axis.

By Proposition~\ref{atl:prop:rescaling} the quaternionic part evolves autonomously, $p_{n+1}=p_n^2+c_{\HH}$, while the $l$-part obeys the linear recursion $q_{n+1}=2\re(p_n)\,q_n+c_l$. On the hyperplane $\re p_0=0$ the first iterate is $q_1=c_l$, whatever $q_0$, so the body is an exact cylinder whose cross-section, by Proposition~\ref{atl:prop:imaginary}, is a disc. In the other quarters, wherever the product of the factors $2\re p_n$ decays along the quaternionic orbit, $q_n$ stays bounded for every $q_0$, and the body is a full column along $l$; this is also why the mirror $l\mapsto-l$ is exact there. The orbits lie in a subalgebra of dimension $3$, the span of $1$, $\im z_0$ and $\im c$.
\clearpage

\subsubsection{B05. Para-octonions with a generic constant}\label{atl:pl:B05}\mbox{}\par\vspace{0.5mm}
\noindent\includegraphics[width=\textwidth]{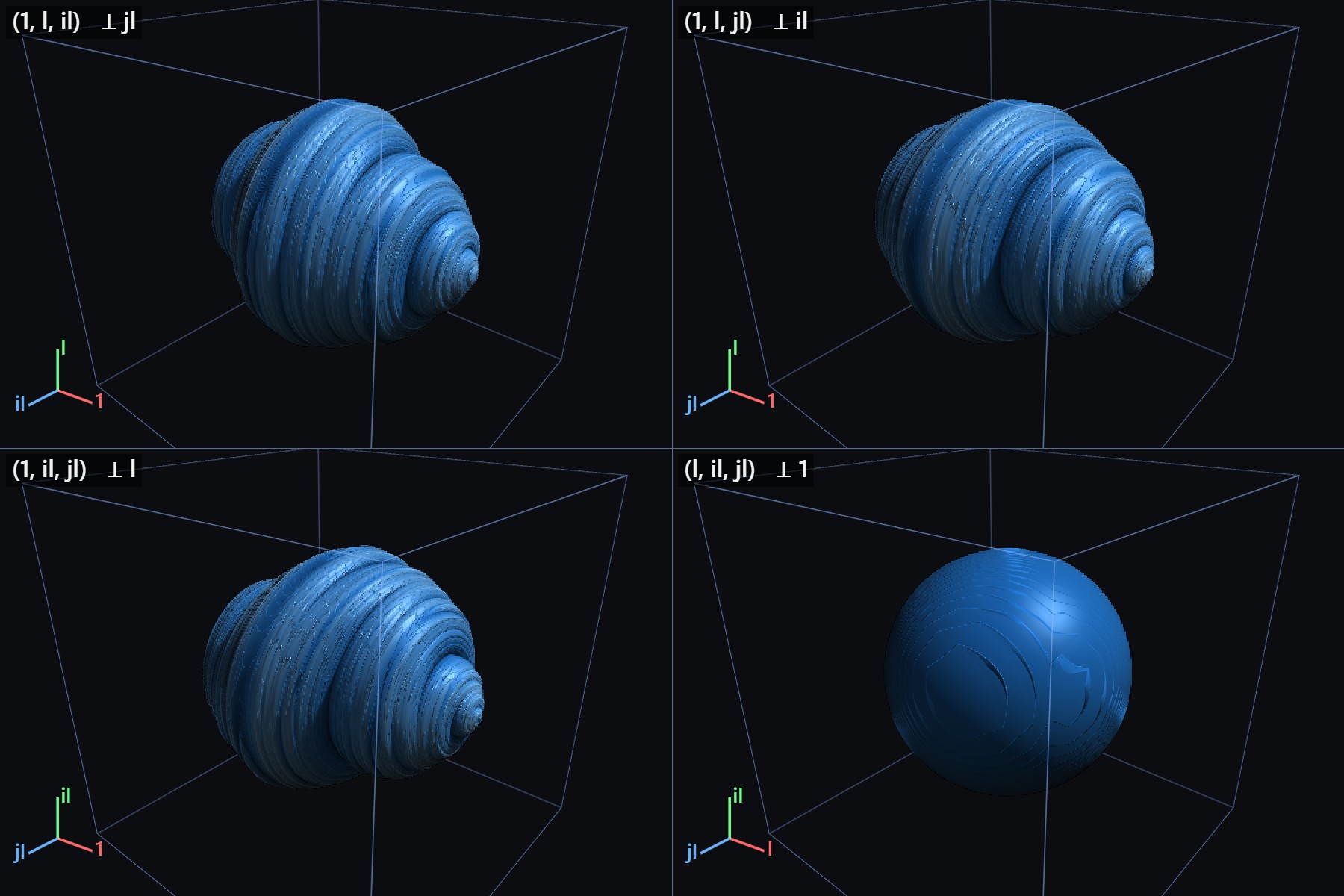}\par\vspace{1mm}
\noindent{\footnotesize \textbf{Map} $z\mapsto z\diamond z+c$ $\cdot$ \textbf{Algebra} Para-octonions (para-Hurwitz $\bar x\bar y$, $\gamma=(-1,-1,-1)$) $\cdot$ $c=(-0.235,\,\allowbreak -0.013,\,\allowbreak 0.029,\,\allowbreak 0.01,\,\allowbreak 0.002,\,\allowbreak 0.095,\,\allowbreak -0.067,\,\allowbreak 0.079)$ ($=0.95\,s^*$, $s^*=0.2901$) $\cdot$ \textbf{Frame} $X,Y,Z,N=1,l,il,jl$ $\cdot$ half-extent 1.4 $\cdot$ 20 iterations $\cdot$ \textbf{Coloring} 8D orbit trap (point), gradient Deep Blue}\par\vspace{1mm}
\noindent{\footnotesize\begin{tabular}{@{}lllll@{}}\toprule Slice (X, Y, Z) & Body [95\%] & Exact ($J\geq0.9995$) & Approximate ($0.97\leq J<0.9995$) & Dim\\\midrule
($1, l, il$) $\perp jl$ & 33.0\% {\scriptsize[32.7, 33.3]} & $P$ & $M$: y; $R_y$: 180 (min.\ 0.994) & 3 \\
($1, l, jl$) $\perp il$ & 32.9\% {\scriptsize[32.6, 33.2]} & $P$ & $M$: y; $R_y$: 180 (min.\ 0.995) & 3 \\
($1, il, jl$) $\perp l$ & 33.2\% {\scriptsize[32.9, 33.5]} & $P$ & -- & 3 \\
($l, il, jl$) $\perp 1$ & 32.2\% {\scriptsize[31.9, 32.5]} & $M$: xyz; $R_x$: all; $R_y$: all; $R_z$: all; $P$ & -- & 3 \\
\bottomrule\end{tabular}}\par\vspace{2mm}
The para-octonions with the generic constant give three similar spiral shells in the quarters through $1$ ($33\%$ each), as in Plate~B01, and a ball in the imaginary quarter ($32\%$, exact under every tested rotation). The shells are rounder and more compact than their octonionic counterparts, with a whorl at the tip; the map is not conjugate to the octonionic one, and on the complex plane it would be the Mandelbar map.

The boundary distance on the generic ray is $s^*=0.290$, much smaller than the octonionic $0.738$. The symmetry pattern is the octonionic one of Plate~B01: $P$ is exact in every quarter, and the mirror $M_y$ with the half-turn $R_y(180^\circ)$ is approximate in the first two quarters ($J=0.994$ and $0.995$). The orbits lie in a subalgebra of dimension $3$: Propositions~\ref{atl:prop:equivariance}, \ref{atl:prop:imaginary} and \ref{atl:prop:span} hold verbatim for the para-Hurwitz product.
\clearpage

\subsubsection{B06. Split para-octonions with a generic constant}\label{atl:pl:B06}\mbox{}\par\vspace{0.5mm}
\noindent\includegraphics[width=\textwidth]{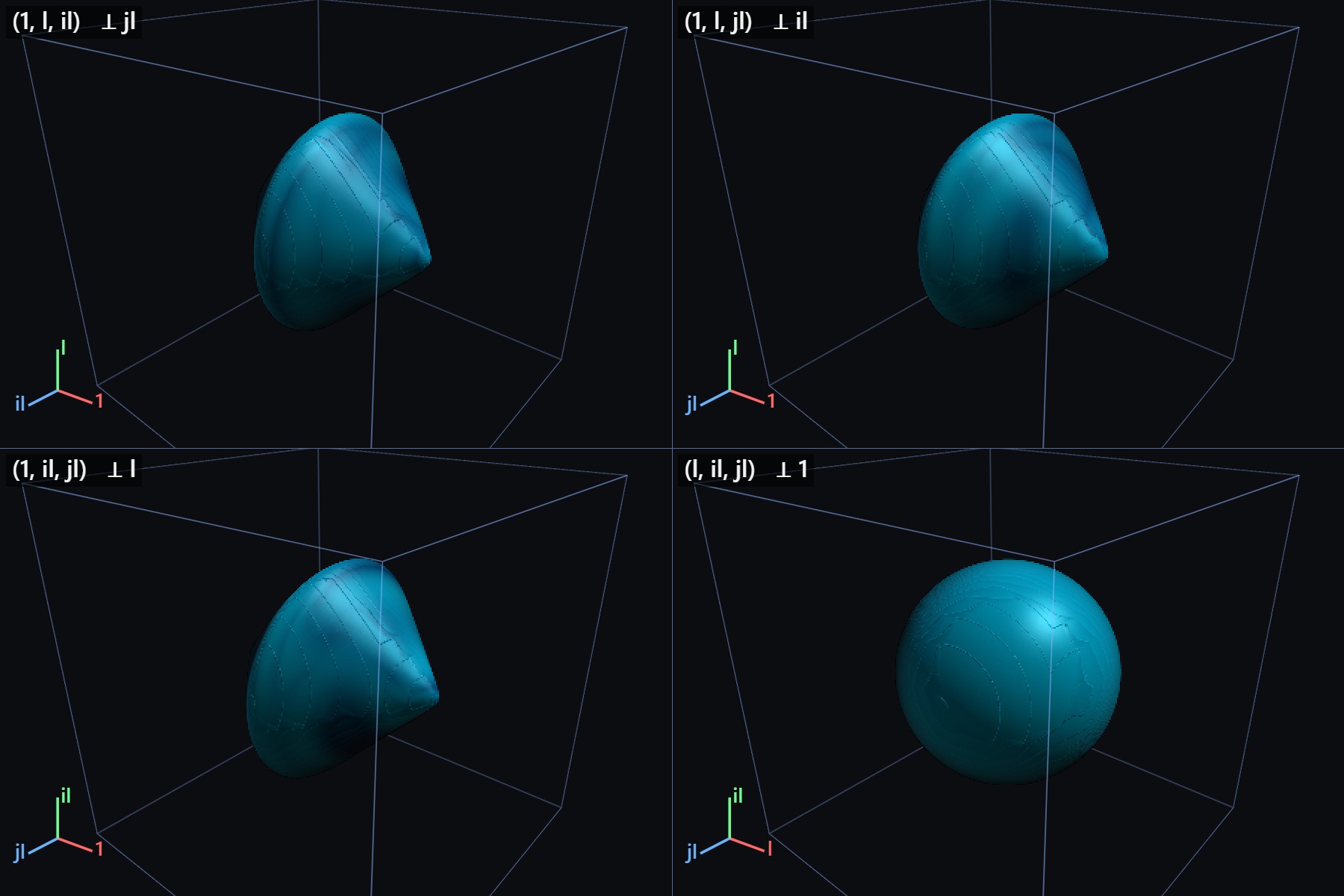}\par\vspace{1mm}
\noindent{\footnotesize \textbf{Map} $z\mapsto z\diamond z+c$ $\cdot$ \textbf{Algebra} Split para-octonions (para-Hurwitz $\bar x\bar y$, $\gamma=(-1,-1,1)$) $\cdot$ $c=(-0.869,\,\allowbreak -0.047,\,\allowbreak 0.108,\,\allowbreak 0.036,\,\allowbreak 0.009,\,\allowbreak 0.349,\,\allowbreak -0.247,\,\allowbreak 0.291)$ ($=0.95\,s^*$, $s^*=1.0725$) $\cdot$ \textbf{Frame} $X,Y,Z,N=1,l,il,jl$ $\cdot$ half-extent 2.4 $\cdot$ 20 iterations $\cdot$ \textbf{Coloring} 8D orbit trap (point), gradient Neon}\par\vspace{1mm}
\noindent{\footnotesize\begin{tabular}{@{}lllll@{}}\toprule Slice (X, Y, Z) & Body [95\%] & Exact ($J\geq0.9995$) & Approximate ($0.97\leq J<0.9995$) & Dim\\\midrule
($1, l, il$) $\perp jl$ & 14.7\% {\scriptsize[14.5, 15.0]} & $P$ & $M$: y; $R_y$: 180 (min.\ 0.995) & 3 \\
($1, l, jl$) $\perp il$ & 14.7\% {\scriptsize[14.5, 14.9]} & $P$ & $M$: y; $R_y$: 180 (min.\ 0.995) & 3 \\
($1, il, jl$) $\perp l$ & 14.9\% {\scriptsize[14.7, 15.1]} & $P$ & -- & 3 \\
($l, il, jl$) $\perp 1$ & 24.3\% {\scriptsize[24.1, 24.6]} & $M$: xyz; $R_x$: all; $R_y$: all; $R_z$: all; $P$ & -- & 3 \\
\bottomrule\end{tabular}}\par\vspace{2mm}
The split para-octonions with the generic constant ($s^*=1.072$), in the cube of half-extent $2.4$, show smooth, pointed bodies in the three quarters through $1$ ($15\%$ of the ball each), and a ball in the imaginary quarter ($24\%$, exact under every tested rotation). Compared with Plate~B02 the bodies are more pointed and slightly smaller. The imaginary quarter is a ball because the split norm is negative definite on the span of $l,il,jl$, and $z\bullet z=-\n(z)$ (Proposition~\ref{atl:prop:imaginary}).

The symmetry pattern is again that of Plates~B01 and B02: $P$ exact everywhere, and $M_y$ with $R_y(180^\circ)$ approximate in the first two quarters ($J=0.995$). The plate completes the comparison of the four quadratic algebras with a generic constant: octonions, split-octonions, para-octonions and split para-octonions all have orbits in a subalgebra of dimension $3$ and a ball in the imaginary quarter.
\clearpage

\subsubsection{B07. The Okubo algebra with a generic constant}\label{atl:pl:B07}\mbox{}\par\vspace{0.5mm}
\noindent\includegraphics[width=\textwidth]{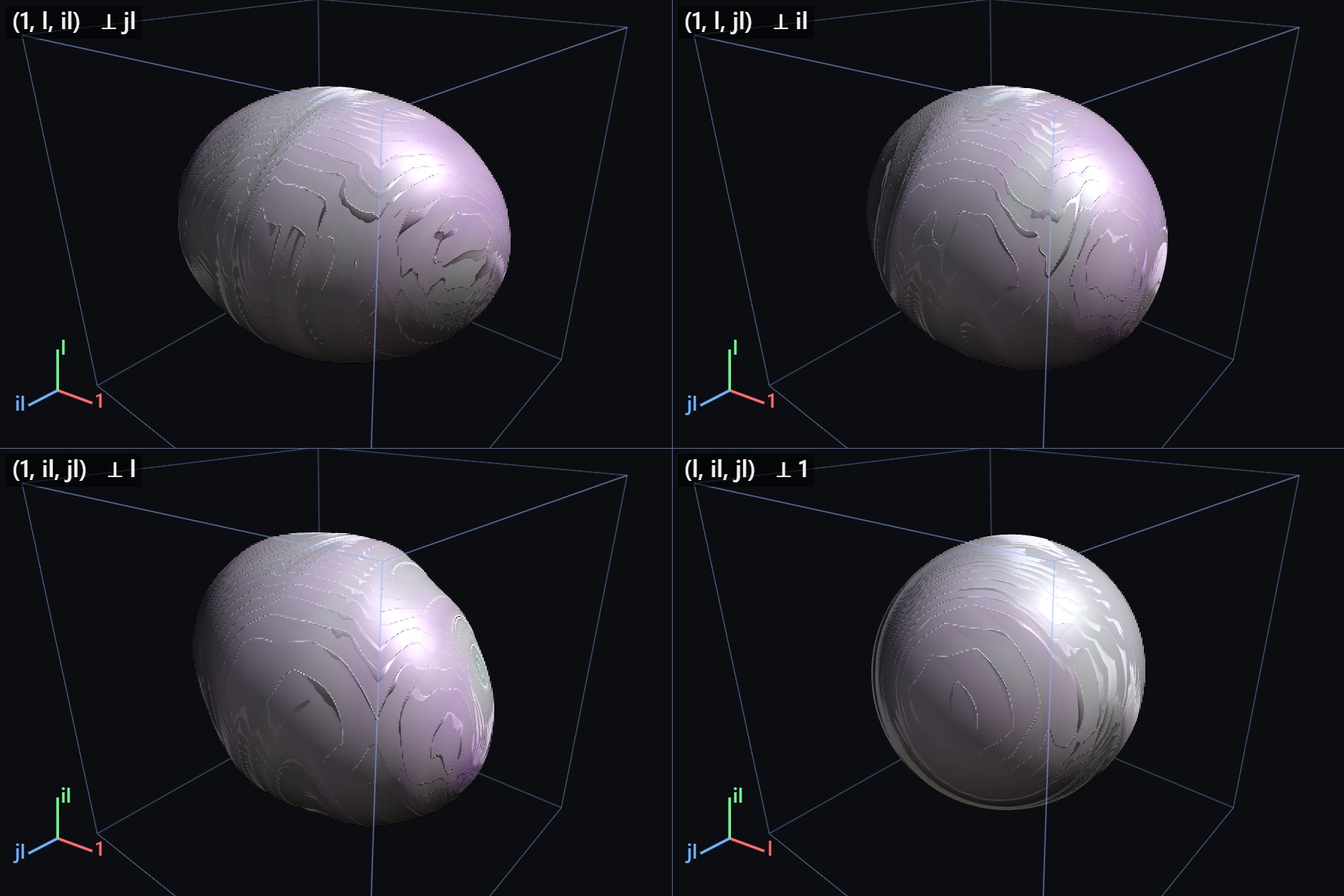}\par\vspace{1mm}
\noindent{\footnotesize \textbf{Map} $z\mapsto z\diamond z+c$ $\cdot$ \textbf{Algebra} Okubo algebra (Petersson, $\gamma=(-1,-1,-1)$, $\theta=120^\circ$) $\cdot$ $c=(-0.397,\,\allowbreak -0.021,\,\allowbreak 0.049,\,\allowbreak 0.016,\,\allowbreak 0.004,\,\allowbreak 0.16,\,\allowbreak -0.113,\,\allowbreak 0.133)$ ($=0.95\,s^*$, $s^*=0.4900$) $\cdot$ \textbf{Frame} $X,Y,Z,N=1,l,il,jl$ $\cdot$ half-extent 1.4 $\cdot$ 20 iterations $\cdot$ \textbf{Coloring} 8D orbit trap (point), gradient Aurora}\par\vspace{1mm}
\noindent{\footnotesize\begin{tabular}{@{}lllll@{}}\toprule Slice (X, Y, Z) & Body [95\%] & Exact ($J\geq0.9995$) & Approximate ($0.97\leq J<0.9995$) & Dim\\\midrule
($1, l, il$) $\perp jl$ & 52.6\% {\scriptsize[52.3, 52.9]} & $P$ & $M$: y; $R_y$: 180 (min.\ 0.986) & 8 \\
($1, l, jl$) $\perp il$ & 51.8\% {\scriptsize[51.5, 52.1]} & $P$ & -- & 8 \\
($1, il, jl$) $\perp l$ & 51.3\% {\scriptsize[51.0, 51.6]} & $P$ & -- & 8 \\
($l, il, jl$) $\perp 1$ & 42.2\% {\scriptsize[41.9, 42.5]} & $P$ & $M$: x; $R_x$: 180,37; $R_y$: 90,72,60,37; $R_z$: 90,72,60,37 (min.\ 0.972) & 8 \\
\bottomrule\end{tabular}}\par\vspace{2mm}
The Okubo algebra with the generic constant ($s^*=0.490$) produces fat, smooth, egg-like bodies in the three quarters through $1$ ($53\%$, $52\%$ and $51\%$ of the ball), covered by thin irregular ridges of the surface itself (they are the same in all the colorings of Part~N) and quite unlike the conch shells of the quadratic algebras. The measurements explain why the pictures look so different: in every quarter the orbits lie in a subalgebra of dimension $8$, i.e.\ $z_0$ and $c$ generate the whole algebra, so the orbit of $z_0$ under $z*z+c$ has no reason to stay in a subspace. The three-dimensional sections see a genuinely eight-dimensional dynamics.

The imaginary quarter $(l,il,jl)$ is a nearly round body ($42\%$). Its only exact symmetry is $P$, but many rotations are approximate, with $J$ between $0.97$ and $0.99$ (for instance $R_z(37^\circ)$ at $0.986$ and $R_y(37^\circ)$ at $0.982$): the body is close to a ball without being one, since the identity $z*z=-\n(z)$ of Proposition~\ref{atl:prop:imaginary} fails in the Okubo algebra. In the quarters through $1$ only $P$ is exact, with an approximate mirror $M_y$ in $(1,l,il)$ ($J=0.986$). The absence of revolutions is a frame effect: the stabilizer of a generic $c$ is a two-dimensional torus conjugate to that of Proposition~\ref{atl:prop:torus}, whose root planes are in general position with respect to the coordinate axes, and in a frame adapted to it the revolutions reappear (Plate~R09).
\clearpage

\subsubsection{B08. The split Okubo algebra with a generic constant}\label{atl:pl:B08}\mbox{}\par\vspace{0.5mm}
\noindent\includegraphics[width=\textwidth]{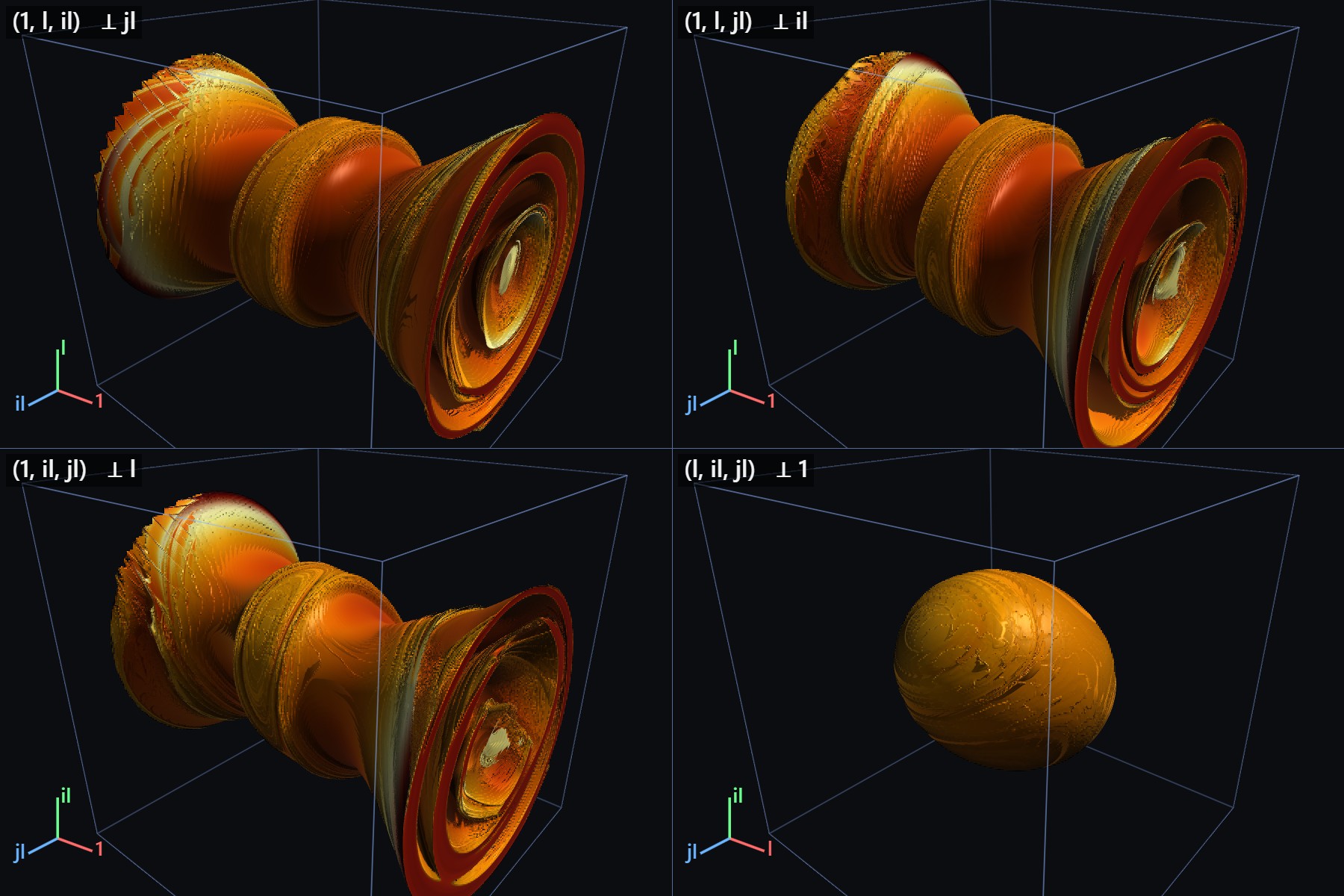}\par\vspace{1mm}
\noindent{\footnotesize \textbf{Map} $z\mapsto z\diamond z+c$ $\cdot$ \textbf{Algebra} Split Okubo algebra (Petersson, $\gamma=(-1,-1,1)$, $\theta=120^\circ$) $\cdot$ $c=(-0.433,\,\allowbreak -0.023,\,\allowbreak 0.054,\,\allowbreak 0.018,\,\allowbreak 0.004,\,\allowbreak 0.174,\,\allowbreak -0.123,\,\allowbreak 0.145)$ ($=0.95\,s^*$, $s^*=0.5339$) $\cdot$ \textbf{Frame} $X,Y,Z,N=1,l,il,jl$ $\cdot$ half-extent 1.4 $\cdot$ 20 iterations $\cdot$ \textbf{Coloring} 8D orbit trap (point), gradient Fire}\par\vspace{1mm}
\noindent{\footnotesize\begin{tabular}{@{}lllll@{}}\toprule Slice (X, Y, Z) & Body [95\%] & Exact ($J\geq0.9995$) & Approximate ($0.97\leq J<0.9995$) & Dim\\\midrule
($1, l, il$) $\perp jl$ & 35.5\% {\scriptsize[35.2, 35.8]} & $P$ & -- & 8 \\
($1, l, jl$) $\perp il$ & 37.4\% {\scriptsize[37.1, 37.7]} & $P$ & -- & 8 \\
($1, il, jl$) $\perp l$ & 38.2\% {\scriptsize[37.9, 38.5]} & $P$ & -- & 8 \\
($l, il, jl$) $\perp 1$ & 21.3\% {\scriptsize[21.1, 21.6]} & $P$ & $M$: x; $R_x$: 180 (min.\ 0.974) & 8 \\
\bottomrule\end{tabular}}\par\vspace{2mm}
The split Okubo algebra with the generic constant ($s^*=0.534$) gives three hourglass-shaped bodies in the quarters through $1$ ($36\%$, $37\%$ and $38\%$ of the ball): a narrow waist and two wide ends, cut by the faces of the cube, where the flat cross-sections show concentric rings. The imaginary quarter holds a small egg-shaped body ($21\%$). The flare of the hourglasses is the signature of the split norm, which has signature $(1,2)$ on the sections through $1$.

In every quarter the orbits lie in a subalgebra of dimension $8$, as in Plate~B07. The only exact symmetry is $P$ (Remark~\ref{atl:rem:minus}); in the imaginary quarter the mirror $M_x$ and the half-turn $R_x(180^\circ)$ are approximate ($J=0.974$). Compared with Plate~B07, the split norm opens the bodies along the null directions, while the constant, generic with respect to the coordinate frame, leaves no other symmetry visible in these sections. As for all split algebras, the rendered bodies are approximations.
\clearpage

\subsection{Part C. Symmetry tests in the plane of \texorpdfstring{$j$}{j} and \texorpdfstring{$l$}{l}}\label{atl:part:C}
The frame $F_C=(1,j,l\,|\,i)$ puts the complex direction $i$ of the constant in the hidden axis, so that the top-left quarter $(1,j,l)$ is \emph{orthogonal} to $\im c$. It is the ideal section to test rotational symmetry: in the octonions any rotation of the plane of $j$ and $l$ fixes $c$ and commutes with the map. The part compares the octonions, the para-octonions, the Okubo algebra and three Petersson algebras, and then rotates the frame by $30^\circ$ in the plane of $j$ and $l$: an experiment that the octonions pass and the Okubo algebra does not, because in the Okubo algebra $j$ and $l$ lie in two different root planes, $V_1$ and $V_2$, and no one-parameter subgroup of the stabilizer of $c$ rotates the plane of $j$ and $l$ (Proposition~\ref{atl:prop:torus}). The revolution is not broken but displaced: it appears in the frame $(1,l,il\,|\,j)$ of Plate~R01. In this frame the torus acts in the same way for every Petersson angle, which is why Plates~C04--C06 share the symmetry pattern of C03, with the single exception of the section $(j,l,i)$ at $\theta=90^\circ$.

\subsubsection{C01. Octonions: revolution in the plane of \texorpdfstring{$j$}{j} and \texorpdfstring{$l$}{l}}\label{atl:pl:C01}\mbox{}\par\vspace{0.5mm}
\noindent\includegraphics[width=\textwidth]{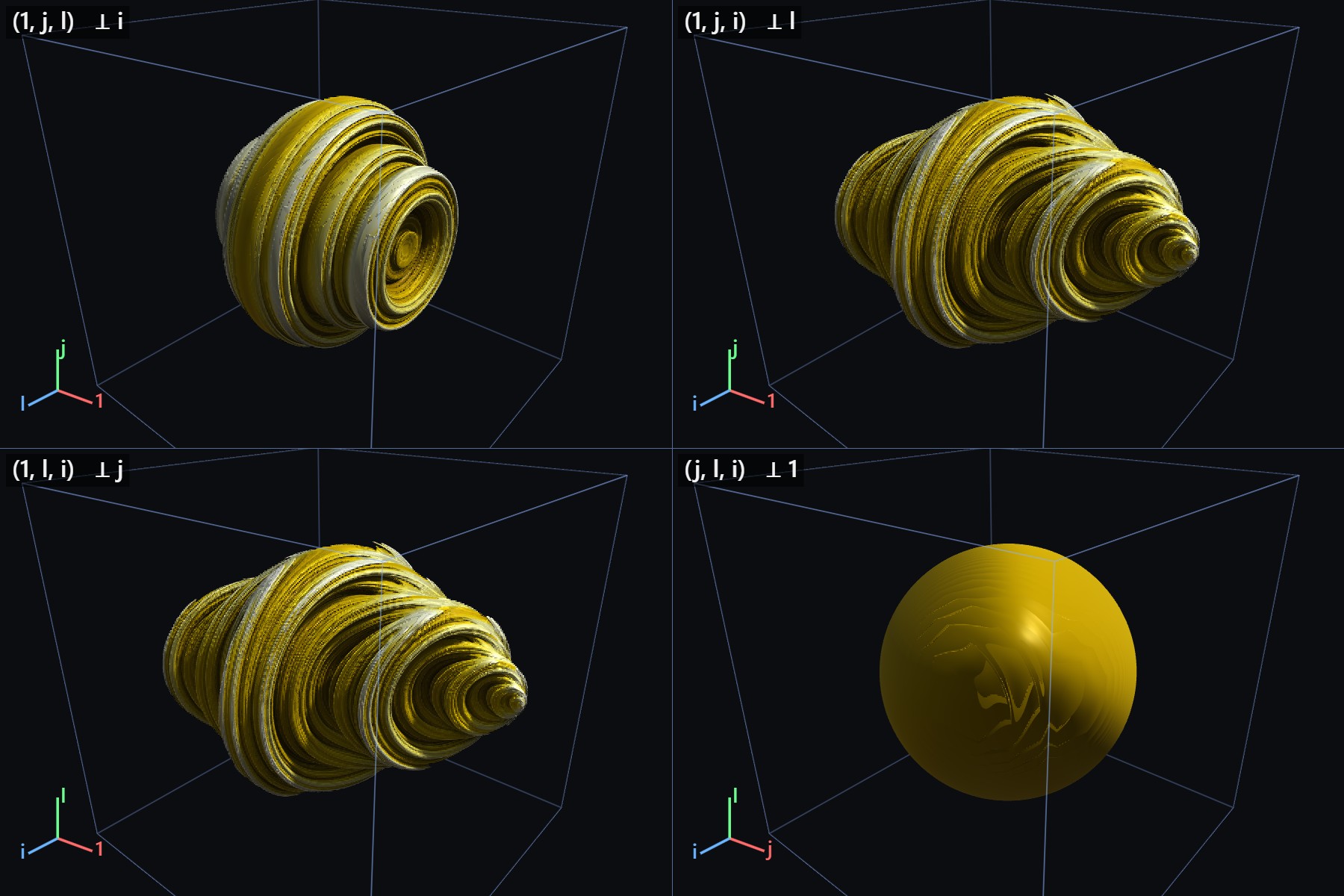}\par\vspace{1mm}
\noindent{\footnotesize \textbf{Map} $z\mapsto z\diamond z+c$ $\cdot$ \textbf{Algebra} Octonions (Cayley--Dickson, $\gamma=(-1,-1,-1)$) $\cdot$ $c=(-0.231,\,\allowbreak 0.62,\,\allowbreak 0,\,\allowbreak 0,\,\allowbreak 0,\,\allowbreak 0,\,\allowbreak 0,\,\allowbreak 0)$ ($=0.97\,s^*$, $s^*=0.6824$) $\cdot$ \textbf{Frame} $X,Y,Z,N=1,j,l,i$ $\cdot$ half-extent 1.4 $\cdot$ 20 iterations $\cdot$ \textbf{Coloring} 8D orbit trap (point), gradient Gold}\par\vspace{1mm}
\noindent{\footnotesize\begin{tabular}{@{}lllll@{}}\toprule Slice (X, Y, Z) & Body [95\%] & Exact ($J\geq0.9995$) & Approximate ($0.97\leq J<0.9995$) & Dim\\\midrule
($1, j, l$) $\perp i$ & 22.5\% {\scriptsize[22.3, 22.8]} & $M$: xyz; $R_x$: all; $R_y$: 180; $R_z$: 180; $P$ & -- & 3 \\
($1, j, i$) $\perp l$ & 24.4\% {\scriptsize[24.1, 24.7]} & $M$: y; $R_y$: 180; $P$ & -- & 3 \\
($1, l, i$) $\perp j$ & 24.4\% {\scriptsize[24.1, 24.7]} & $M$: y; $R_y$: 180; $P$ & -- & 3 \\
($j, l, i$) $\perp 1$ & 34.7\% {\scriptsize[34.4, 35.0]} & $M$: xyz; $R_x$: all; $R_y$: all; $R_z$: all; $P$ & -- & 3 \\
\bottomrule\end{tabular}}\par\vspace{2mm}
In the frame $F_C$ the hidden axis is $i$, the direction of $\im c$, so that the upper-left quarter $(1,j,l)$ is orthogonal to $\im c$. The octonionic body there is a rounded shell whose ridges are circles centred on the axis of $1$, and it is an exact body of revolution about that axis: the rotations about the local $x$ axis by all tested angles, $37^\circ$ included, are exact. This is Proposition~\ref{atl:prop:equivariance}: every rotation of the plane of $j$ and $l$ fixes $1$ and $c$ and commutes with the map.

The upper-right quarter $(1,j,i)$ and the lower-left quarter $(1,l,i)$ show the same shell, wound in concentric ridges, with identical body shares ($24.4\%$) and identical symmetries: the exchange of $j$ and $l$ is one of the rotations of Proposition~\ref{atl:prop:equivariance}. The quarter $(j,l,i)$ of purely imaginary points is a ball, exact under every tested rotation, because $z^2=-\n(z)$ there (Proposition~\ref{atl:prop:imaginary}). Every quarter is centrally symmetric (Remark~\ref{atl:rem:minus}), and the orbits lie in a subalgebra of dimension $3$.

The plate is the octonionic reference for the tests of this part: a body of revolution, two identical shells and a ball.
\clearpage

\subsubsection{C02. Para-octonions: the same test}\label{atl:pl:C02}\mbox{}\par\vspace{0.5mm}
\noindent\includegraphics[width=\textwidth]{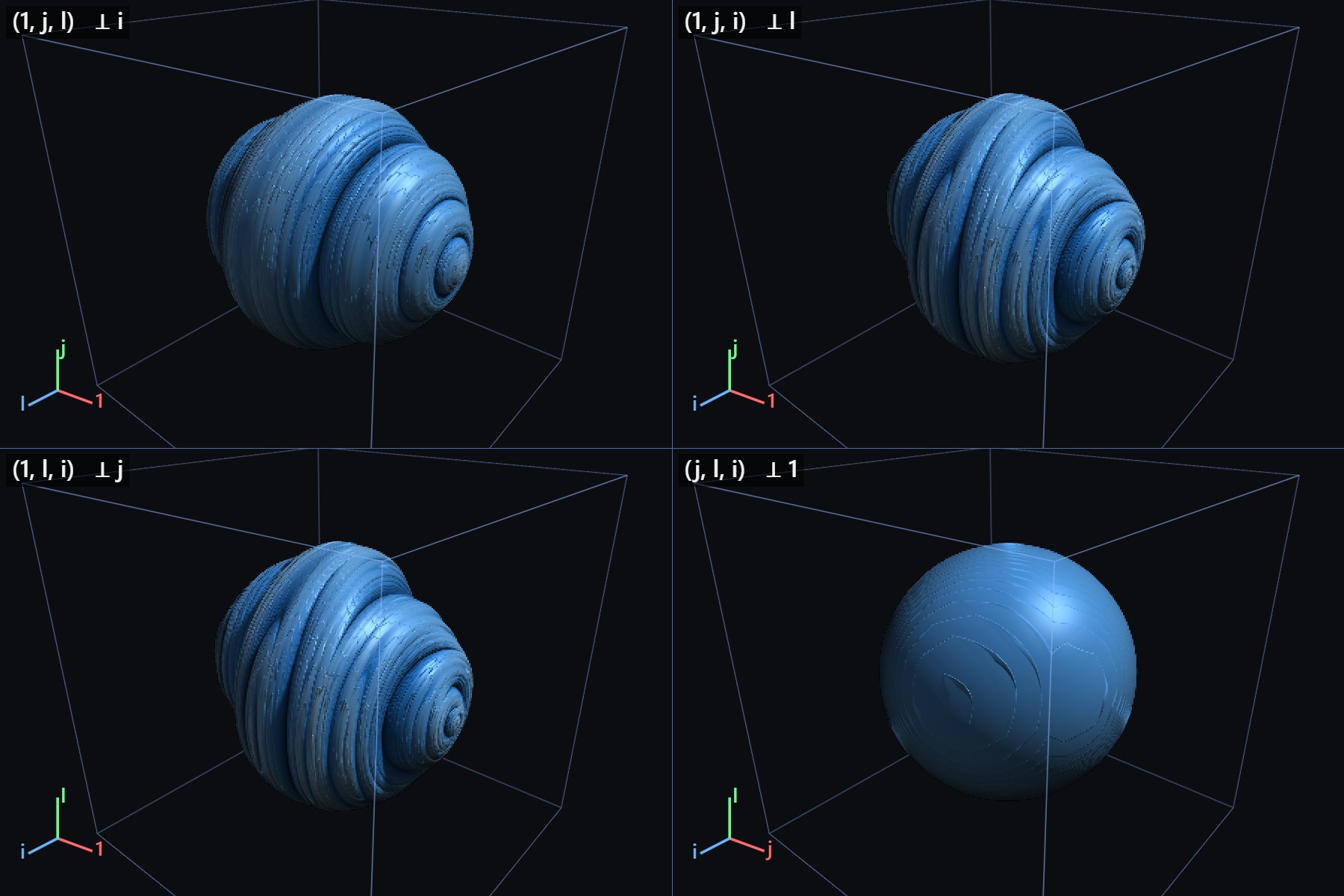}\par\vspace{1mm}
\noindent{\footnotesize \textbf{Map} $z\mapsto z\diamond z+c$ $\cdot$ \textbf{Algebra} Para-octonions (para-Hurwitz $\bar x\bar y$, $\gamma=(-1,-1,-1)$) $\cdot$ $c=(-0.086,\,\allowbreak 0.231,\,\allowbreak 0,\,\allowbreak 0,\,\allowbreak 0,\,\allowbreak 0,\,\allowbreak 0,\,\allowbreak 0)$ ($=0.97\,s^*$, $s^*=0.2540$) $\cdot$ \textbf{Frame} $X,Y,Z,N=1,j,l,i$ $\cdot$ half-extent 1.4 $\cdot$ 20 iterations $\cdot$ \textbf{Coloring} 8D orbit trap (point), gradient Deep Blue}\par\vspace{1mm}
\noindent{\footnotesize\begin{tabular}{@{}lllll@{}}\toprule Slice (X, Y, Z) & Body [95\%] & Exact ($J\geq0.9995$) & Approximate ($0.97\leq J<0.9995$) & Dim\\\midrule
($1, j, l$) $\perp i$ & 34.5\% {\scriptsize[34.2, 34.8]} & $M$: xyz; $R_x$: all; $R_y$: 180; $R_z$: 180; $P$ & -- & 3 \\
($1, j, i$) $\perp l$ & 35.6\% {\scriptsize[35.3, 35.8]} & $M$: y; $R_y$: 180; $P$ & -- & 3 \\
($1, l, i$) $\perp j$ & 35.6\% {\scriptsize[35.3, 35.8]} & $M$: y; $R_y$: 180; $P$ & -- & 3 \\
($j, l, i$) $\perp 1$ & 35.6\% {\scriptsize[35.3, 35.9]} & $M$: xyz; $R_x$: all; $R_y$: all; $R_z$: all; $P$ & -- & 3 \\
\bottomrule\end{tabular}}\par\vspace{2mm}
The para-octonions pass the same test. The quarter $(1,j,l)$ is an exact body of revolution about $1$, the quarters $(1,j,i)$ and $(1,l,i)$ coincide (body share $35.6\%$ and the same symmetries), and the imaginary quarter $(j,l,i)$ is a ball. The bodies are rounder and smoother than in Plate~C01, and larger ($34.5\%$ against $22.5\%$ in the quarter $(1,j,l)$), because the constant is smaller: on the same ray the para-octonionic boundary distance is much smaller than the octonionic one (Plate~A05). On $\CC$ the map is the Mandelbar map $\bar z^2+c$.

The explanation is the same as for the octonions. The para-Hurwitz product $\bar x\,\bar y$ is built from the octonionic product and the conjugation, which commutes with every orthogonal map fixing $1$, so that Proposition~\ref{atl:prop:equivariance} holds verbatim, and the orbits lie in a subalgebra of dimension $3$. This plate and the previous one show that the rotational symmetry in the plane of $j$ and $l$ is a feature of the whole para-Hurwitz family, not only of unital algebras.
\clearpage

\subsubsection{C03. Okubo: \texorpdfstring{$j$}{j} and \texorpdfstring{$l$}{l} in different root planes}\label{atl:pl:C03}\mbox{}\par\vspace{0.5mm}
\noindent\includegraphics[width=\textwidth]{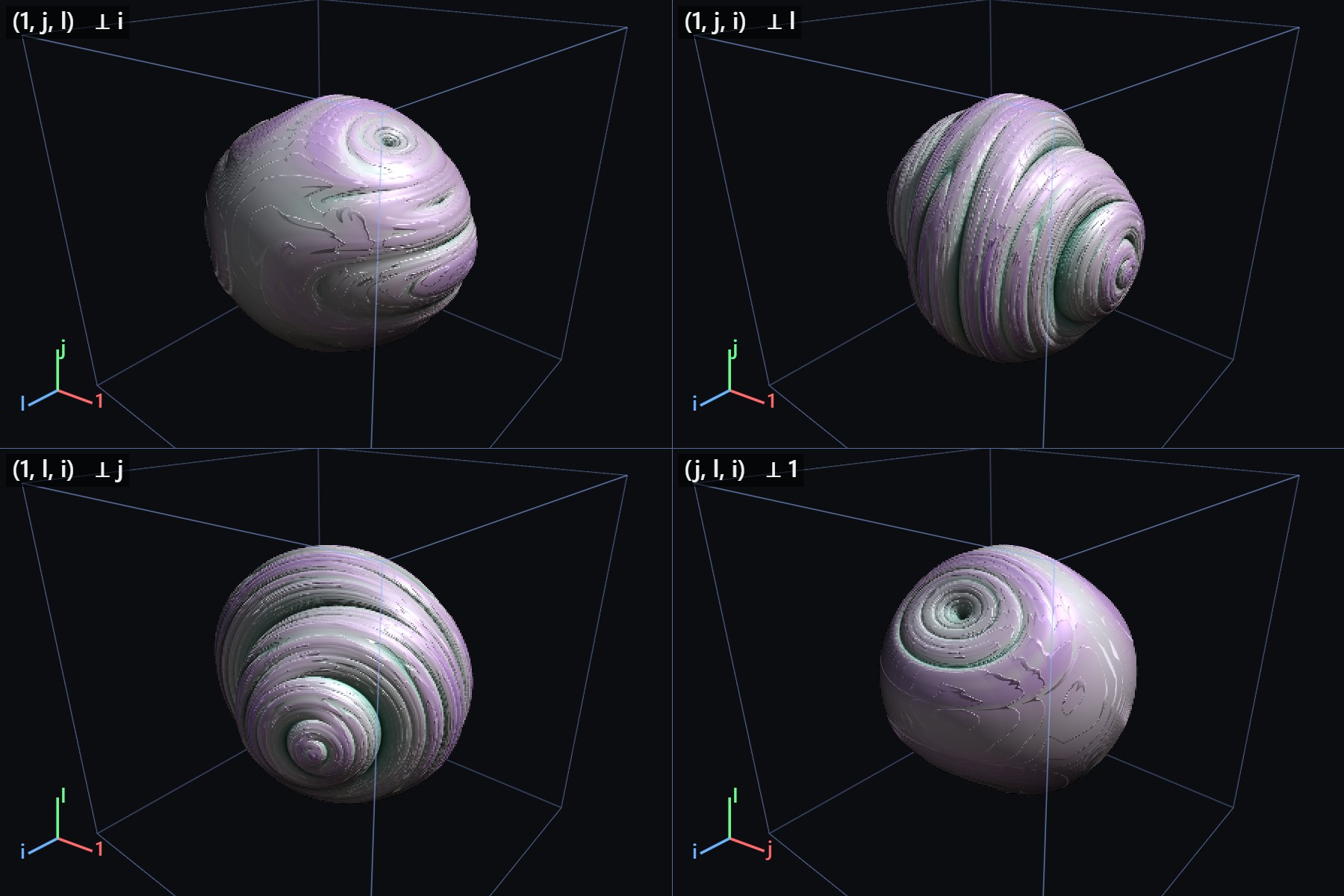}\par\vspace{1mm}
\noindent{\footnotesize \textbf{Map} $z\mapsto z\diamond z+c$ $\cdot$ \textbf{Algebra} Okubo algebra (Petersson, $\gamma=(-1,-1,-1)$, $\theta=120^\circ$) $\cdot$ $c=(-0.086,\,\allowbreak 0.231,\,\allowbreak 0,\,\allowbreak 0,\,\allowbreak 0,\,\allowbreak 0,\,\allowbreak 0,\,\allowbreak 0)$ ($=0.97\,s^*$, $s^*=0.2540$) $\cdot$ \textbf{Frame} $X,Y,Z,N=1,j,l,i$ $\cdot$ half-extent 1.4 $\cdot$ 20 iterations $\cdot$ \textbf{Coloring} 8D orbit trap (point), gradient Aurora}\par\vspace{1mm}
\noindent{\footnotesize\begin{tabular}{@{}lllll@{}}\toprule Slice (X, Y, Z) & Body [95\%] & Exact ($J\geq0.9995$) & Approximate ($0.97\leq J<0.9995$) & Dim\\\midrule
($1, j, l$) $\perp i$ & 37.4\% {\scriptsize[37.1, 37.7]} & $M$: xyz; $R_x$: 180; $R_y$: 180; $R_z$: 180; $P$ & -- & 5 \\
($1, j, i$) $\perp l$ & 36.0\% {\scriptsize[35.7, 36.3]} & $M$: y; $R_y$: 180; $P$ & -- & 3 \\
($1, l, i$) $\perp j$ & 34.8\% {\scriptsize[34.5, 35.1]} & $M$: y; $R_y$: 180; $P$ & -- & 3 \\
($j, l, i$) $\perp 1$ & 33.2\% {\scriptsize[32.9, 33.5]} & $M$: xyz; $R_x$: 180; $R_y$: 180; $R_z$: 180; $P$ & -- & 5 \\
\bottomrule\end{tabular}}\par\vspace{2mm}
The Okubo algebra, with the same constant as Plate~C02, does not pass the test. The quarter $(1,j,l)$ is a nearly round body with a vortex near the top and swirling bands, but it is not a body of revolution: the rotations about the axis of $1$ score only $J\approx0.95$--$0.96$ ($0.958$ at $37^\circ$), below the threshold of the approximate symmetries. What survives is exact but finite: the three mirrors, the three half-turns and the point reflection. The reason is Proposition~\ref{atl:prop:torus}. The torus $\varphi_{a,b}$ that fixes $c$ rotates $j$ inside the root plane $V_1=\langle j,k\rangle$ by $2\alpha$ and $l$ inside $V_2=\langle l,il\rangle$ by $\beta-\alpha$ \eqref{atl:eq:rootangles}: $j$ and $l$ lie in different root planes, and no element of the torus rotates the plane of $j$ and $l$ into itself, except the half-turns that change the sign of $j$ ($\alpha=\beta=90^\circ$) or of $l$ ($\alpha=0$, $\beta=180^\circ$). They give two of the mirrors, and the point reflection $z\mapsto-z$ of Remark~\ref{atl:rem:minus} gives the rest.

The quarter $(1,j,i)$ lies in $\HH$, where the product is the para-quaternionic one, and it is the shell of Plate~C02 ($36.0\%$ against $35.6\%$, with the same symmetries). The quarter $(1,l,i)$ is a similar banded shell, but it is no longer a copy of $(1,j,i)$: the centre of its concentric ridges is seen at the lower left, and its body share is $34.8\%$. The imaginary quarter $(j,l,i)$ is not a ball but an egg-shaped body with a vortex at the upper left, with the same finite symmetries as $(1,j,l)$: Proposition~\ref{atl:prop:imaginary} does not hold in the Okubo algebra. The two quarters that contain both $j$ and $l$ generate a subalgebra of dimension $5$, the other two of dimension $3$.

It is worth noticing what the plate does and does not show. The frame $F_C$ contains no root plane, so it cannot show a revolution; but the continuous symmetry is not lost, it is displaced to the sections that contain a root plane, such as $(1,l,il)$, which are exact bodies of revolution (Plates~R01--R03). Of the $15$-dimensional group of Proposition~\ref{atl:prop:equivariance} the Okubo algebra keeps a $2$-torus, and in this frame only its half-turns are visible.
\clearpage

\subsubsection{C04. Petersson algebra, \texorpdfstring{$\theta=60^\circ$}{theta=60\textdegree{}}}\label{atl:pl:C04}\mbox{}\par\vspace{0.5mm}
\noindent\includegraphics[width=\textwidth]{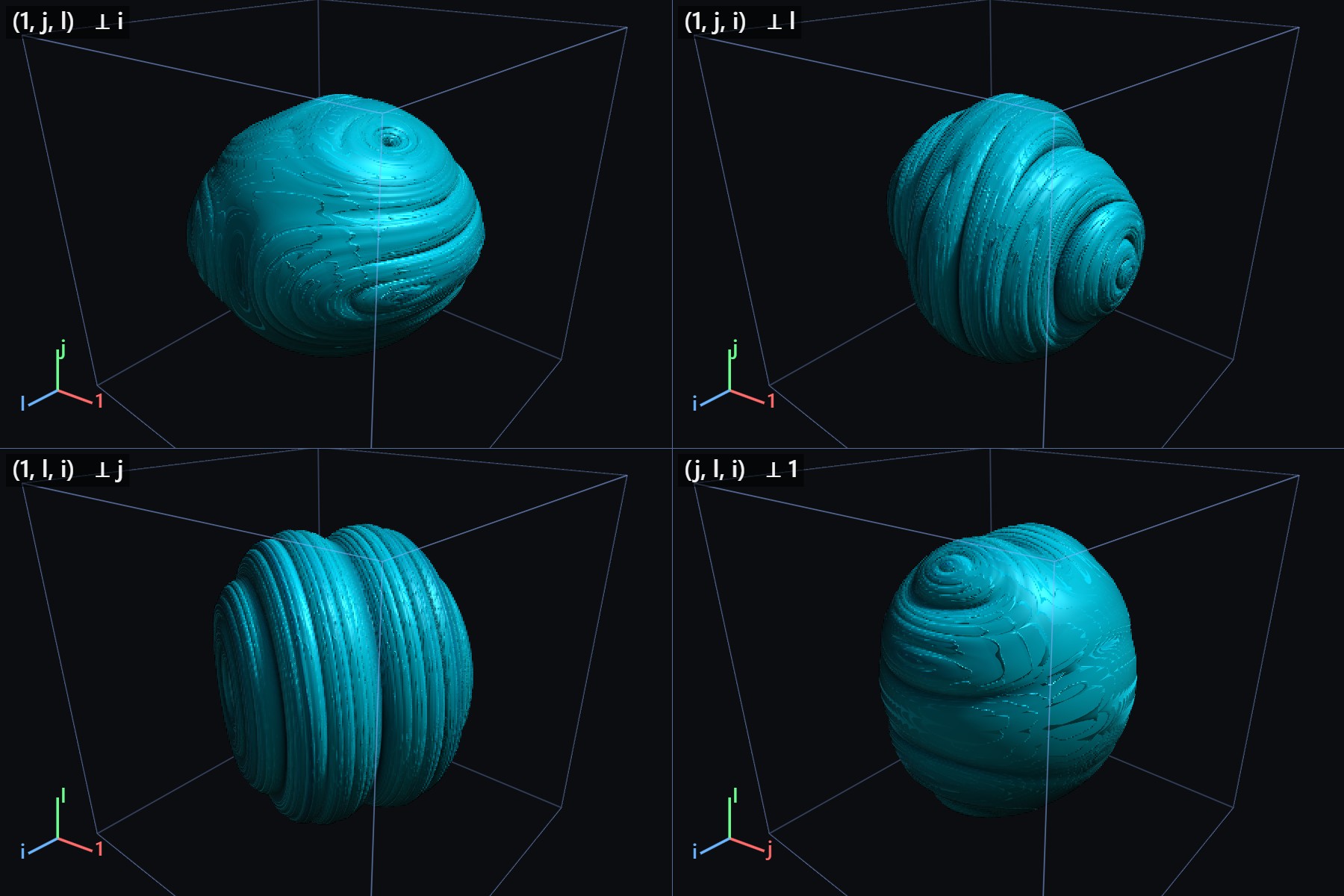}\par\vspace{1mm}
\noindent{\footnotesize \textbf{Map} $z\mapsto z\diamond z+c$ $\cdot$ \textbf{Algebra} custom: Petersson, $\gamma=(-1,-1,-1)$, $\theta=60^\circ$ $\cdot$ $c=(-0.086,\,\allowbreak 0.231,\,\allowbreak 0,\,\allowbreak 0,\,\allowbreak 0,\,\allowbreak 0,\,\allowbreak 0,\,\allowbreak 0)$ ($=0.97\,s^*$, $s^*=0.2540$) $\cdot$ \textbf{Frame} $X,Y,Z,N=1,j,l,i$ $\cdot$ half-extent 1.4 $\cdot$ 20 iterations $\cdot$ \textbf{Coloring} 8D orbit trap (point), gradient Neon}\par\vspace{1mm}
\noindent{\footnotesize\begin{tabular}{@{}lllll@{}}\toprule Slice (X, Y, Z) & Body [95\%] & Exact ($J\geq0.9995$) & Approximate ($0.97\leq J<0.9995$) & Dim\\\midrule
($1, j, l$) $\perp i$ & 42.0\% {\scriptsize[41.7, 42.3]} & $M$: xyz; $R_x$: 180; $R_y$: 180; $R_z$: 180; $P$ & -- & 5 \\
($1, j, i$) $\perp l$ & 35.7\% {\scriptsize[35.4, 36.0]} & $M$: y; $R_y$: 180; $P$ & -- & 3 \\
($1, l, i$) $\perp j$ & 45.9\% {\scriptsize[45.6, 46.3]} & $M$: y; $R_y$: 180; $P$ & -- & 3 \\
($j, l, i$) $\perp 1$ & 41.4\% {\scriptsize[41.1, 41.7]} & $M$: xyz; $R_x$: 180; $R_y$: 180; $R_z$: 180; $P$ & -- & 5 \\
\bottomrule\end{tabular}}\par\vspace{2mm}
The Petersson algebra with $\theta=60^\circ$ lies halfway between the para-octonions and the Okubo algebra. The quarter $(1,j,i)$ lies in $\HH$, which is fixed by every $\tau_\theta$, and it coincides with that of the para-octonions ($35.7\%$ against $35.6\%$). The quarter $(1,j,l)$ is a rounded, ball-like body with a vortex on top; the quarter $(1,l,i)$ shows two banded lobes side by side, separated by a groove, and it is the largest body of the plate ($45.9\%$); the imaginary quarter $(j,l,i)$ is an egg with a vortex at the upper left.

The table is the same as in Plate~C03: in $(1,j,l)$ and $(j,l,i)$ the three mirrors, the three half-turns and the point reflection, with no rotation scoring above $J=0.94$; in $(1,j,i)$ and $(1,l,i)$ a single mirror, its half-turn and the point reflection. This is expected: the torus $\varphi_{a,b}$ and the anti-automorphism $\sigma$ of Proposition~\ref{atl:prop:torus} exist for every $\theta$, and in the frame $F_C$ they act by the same half-turns. As in the Okubo algebra, the quarters containing both $j$ and $l$ generate a subalgebra of dimension $5$.

For $\theta\neq0,\pm120^\circ$ the Petersson algebra is a composition algebra but not a symmetric one: $(x*y)*x\neq\n(x)\,y$. The pictures therefore do not inherit the degeneracies of Parts~H and~I, while the symmetry pattern in this frame is that of the Okubo algebra.
\clearpage

\subsubsection{C05. Petersson algebra, \texorpdfstring{$\theta=90^\circ$}{theta=90\textdegree{}}}\label{atl:pl:C05}\mbox{}\par\vspace{0.5mm}
\noindent\includegraphics[width=\textwidth]{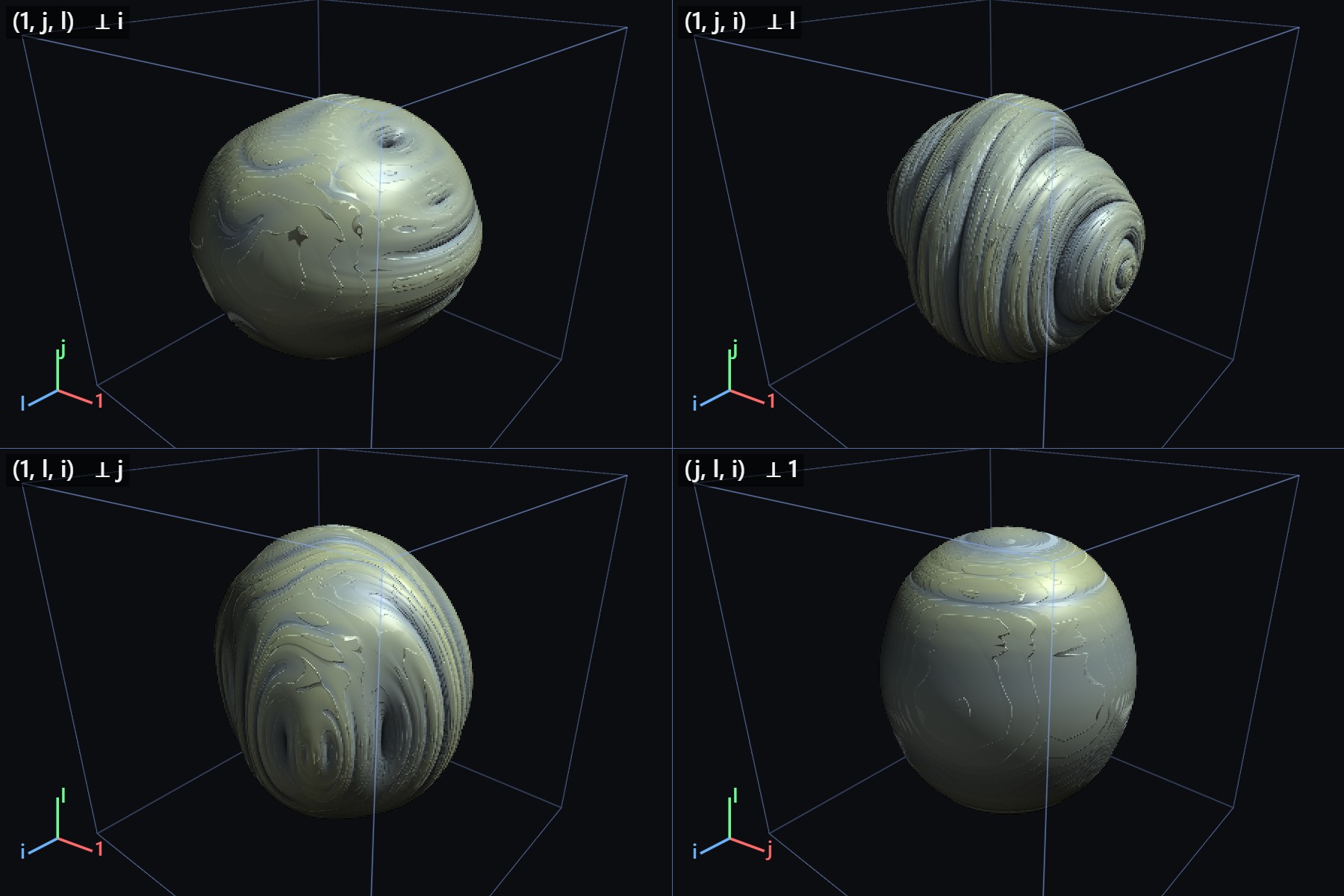}\par\vspace{1mm}
\noindent{\footnotesize \textbf{Map} $z\mapsto z\diamond z+c$ $\cdot$ \textbf{Algebra} custom: Petersson, $\gamma=(-1,-1,-1)$, $\theta=90^\circ$ $\cdot$ $c=(-0.086,\,\allowbreak 0.231,\,\allowbreak 0,\,\allowbreak 0,\,\allowbreak 0,\,\allowbreak 0,\,\allowbreak 0,\,\allowbreak 0)$ ($=0.97\,s^*$, $s^*=0.2540$) $\cdot$ \textbf{Frame} $X,Y,Z,N=1,j,l,i$ $\cdot$ half-extent 1.4 $\cdot$ 20 iterations $\cdot$ \textbf{Coloring} 8D orbit trap (point), gradient Pastel}\par\vspace{1mm}
\noindent{\footnotesize\begin{tabular}{@{}lllll@{}}\toprule Slice (X, Y, Z) & Body [95\%] & Exact ($J\geq0.9995$) & Approximate ($0.97\leq J<0.9995$) & Dim\\\midrule
($1, j, l$) $\perp i$ & 44.0\% {\scriptsize[43.6, 44.3]} & $M$: xyz; $R_x$: 180; $R_y$: 180; $R_z$: 180; $P$ & -- & 5 \\
($1, j, i$) $\perp l$ & 35.8\% {\scriptsize[35.5, 36.1]} & $M$: y; $R_y$: 180; $P$ & -- & 3 \\
($1, l, i$) $\perp j$ & 44.9\% {\scriptsize[44.6, 45.2]} & $M$: y; $R_y$: 180; $P$ & -- & 3 \\
($j, l, i$) $\perp 1$ & 41.1\% {\scriptsize[40.8, 41.4]} & $M$: xyz; $R_x$: 180; $R_y$: all; $R_z$: 180; $P$ & -- & 5 \\
\bottomrule\end{tabular}}\par\vspace{2mm}
At $\theta=90^\circ$ the twist $\tau_\theta$ has order four. The quarter $(1,j,i)$ is again the para-quaternionic shell ($35.8\%$). The quarter $(1,j,l)$ is a rounded body with a vortex dimple near the top and smaller dimples beside it; the quarter $(1,l,i)$ is an egg-shaped banded body whose lower half shows two vortices side by side; the imaginary quarter $(j,l,i)$ is an egg elongated along $l$, with circular bands around its top.

In the first three quarters the table shows the finite pattern of Plates~C03, C04 and C06. The imaginary quarter is different: it is an exact body of revolution about the local $y$ axis, the direction of $l$, with the rotations by all tested angles exact. Quite surprisingly, this revolution is not produced by any automorphism: the section $(i,j,l)$ contains no root plane, so the torus of Proposition~\ref{atl:prop:torus} cannot rotate it. It is explained by Proposition~\ref{atl:prop:theta90}. For $\theta=90^\circ$ the subspace $W=\langle1,i\rangle\oplus\langle l,kl\rangle$ is closed under squaring \eqref{atl:eq:W}, and on $W$ the map satisfies the twisted equivariance $f\circ R_\varphi=R_{-\varphi}\circ f$: the rotation is reversed at every step, so it is not a symmetry of the algebra, but it preserves the norms along the orbits and hence their escape, and $K_c\cap W$ is a body of revolution. Since the square of every point of $\langle i,j,l\rangle$ lies in $W$, the section $(i,j,l)$ inherits a revolution about $l$.

The same angle produces the four identical balls of Plate~E04, where every section contained in $\RR\oplus\HH l$ is the same ball because $x*x=x_0^2+|q|^2$. Here, as in the other Petersson plates, the quarters containing both $j$ and $l$ generate a subalgebra of dimension $5$, the other two of dimension $3$.
\clearpage

\subsubsection{C06. Petersson algebra at the golden angle}\label{atl:pl:C06}\mbox{}\par\vspace{0.5mm}
\noindent\includegraphics[width=\textwidth]{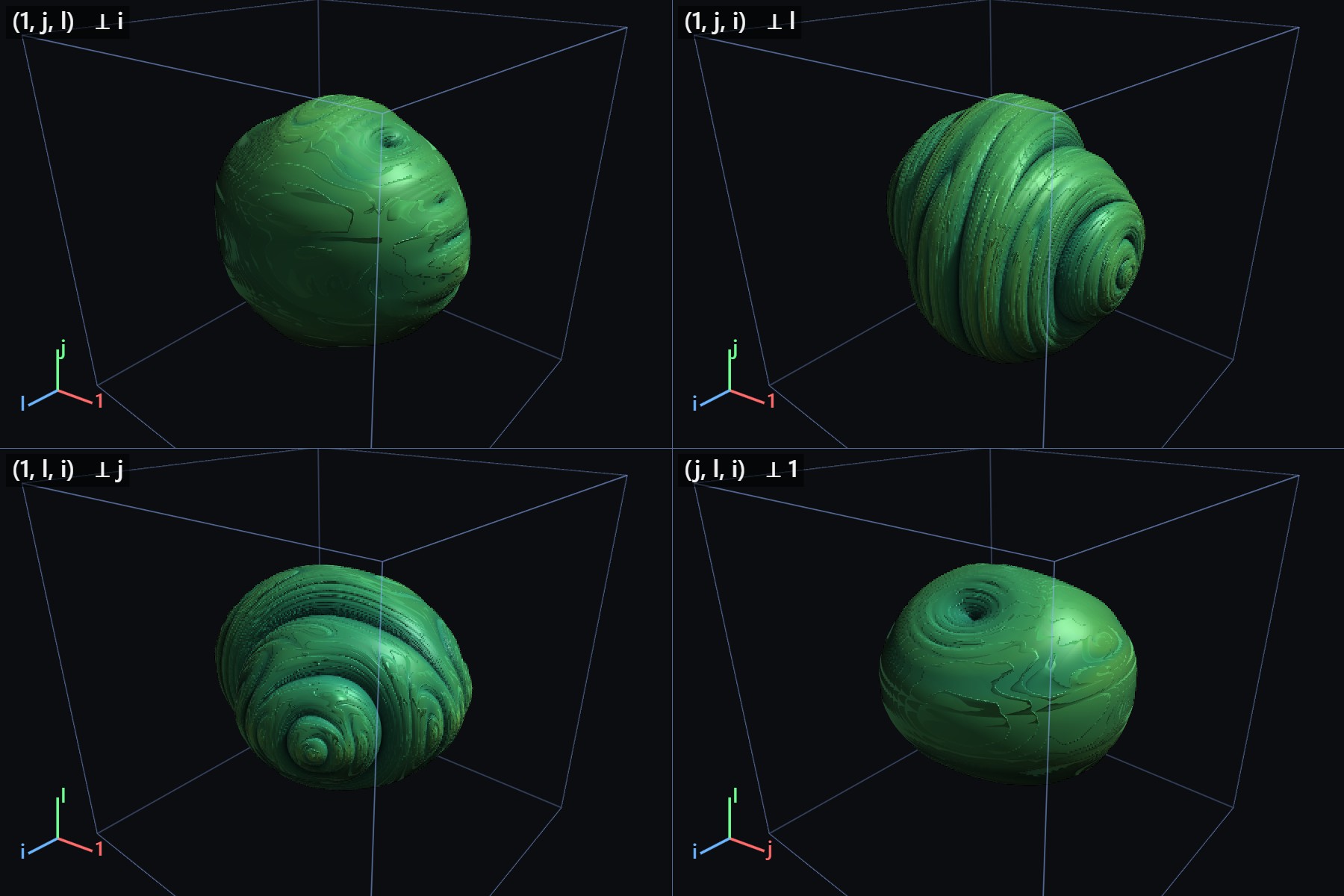}\par\vspace{1mm}
\noindent{\footnotesize \textbf{Map} $z\mapsto z\diamond z+c$ $\cdot$ \textbf{Algebra} custom: Petersson, $\gamma=(-1,-1,-1)$, $\theta=137.508^\circ$ $\cdot$ $c=(-0.086,\,\allowbreak 0.231,\,\allowbreak 0,\,\allowbreak 0,\,\allowbreak 0,\,\allowbreak 0,\,\allowbreak 0,\,\allowbreak 0)$ ($=0.97\,s^*$, $s^*=0.2540$) $\cdot$ \textbf{Frame} $X,Y,Z,N=1,j,l,i$ $\cdot$ half-extent 1.4 $\cdot$ 20 iterations $\cdot$ \textbf{Coloring} 8D orbit trap (point), gradient Viridis-like}\par\vspace{1mm}
\noindent{\footnotesize\begin{tabular}{@{}lllll@{}}\toprule Slice (X, Y, Z) & Body [95\%] & Exact ($J\geq0.9995$) & Approximate ($0.97\leq J<0.9995$) & Dim\\\midrule
($1, j, l$) $\perp i$ & 32.1\% {\scriptsize[31.8, 32.4]} & $M$: xyz; $R_x$: 180; $R_y$: 180; $R_z$: 180; $P$ & -- & 5 \\
($1, j, i$) $\perp l$ & 35.6\% {\scriptsize[35.3, 35.9]} & $M$: y; $R_y$: 180; $P$ & -- & 3 \\
($1, l, i$) $\perp j$ & 28.2\% {\scriptsize[27.9, 28.5]} & $M$: y; $R_y$: 180; $P$ & -- & 3 \\
($j, l, i$) $\perp 1$ & 28.6\% {\scriptsize[28.3, 28.9]} & $M$: xyz; $R_x$: 180; $R_y$: 180; $R_z$: 180; $P$ & -- & 5 \\
\bottomrule\end{tabular}}\par\vspace{2mm}
At the golden angle $\theta=137.5^\circ$ the twist has infinite order. The quarter $(1,j,i)$ is the usual para-quaternionic shell ($35.6\%$). The quarter $(1,j,l)$ is a rounded body with a vortex near the top and a few shallow dimples; the quarter $(1,l,i)$ is a banded shell whose concentric ridges are centred at the lower left, much as in Plate~C03; the imaginary quarter $(j,l,i)$ is a flattened rounded body with a vortex on top.

The table is, quarter by quarter, the same as in Plates~C03 and C04: the three mirrors, the three half-turns and the point reflection in $(1,j,l)$ and $(j,l,i)$, a single mirror with its half-turn in the other two, and a subalgebra of dimension $5$ wherever $j$ and $l$ appear together. These are the torus and the anti-automorphism $\sigma$ of Proposition~\ref{atl:prop:torus}, which exist for every $\theta$ and act in the frame $F_C$ by the same half-turns; only at $\theta=90^\circ$ (Plate~C05) does a further revolution appear, for a reason that is not an automorphism.

The comparison with Plate~C03 ($\theta=120^\circ$) shows that nothing special happens at the Okubo angle in this frame: the shapes and the body shares change with $\theta$ ($32.1\%$, $28.2\%$ and $28.6\%$ here against $37.4\%$, $34.8\%$ and $33.2\%$ in Plate~C03), but the symmetry pattern does not. The Okubo algebra is special for its identities, as Parts~H and~I will show, not for the symmetries of these sections.
\clearpage

\subsubsection{C07. Octonions, frame rotated by \texorpdfstring{$30^\circ$}{30\textdegree{}}}\label{atl:pl:C07}\mbox{}\par\vspace{0.5mm}
\noindent\includegraphics[width=\textwidth]{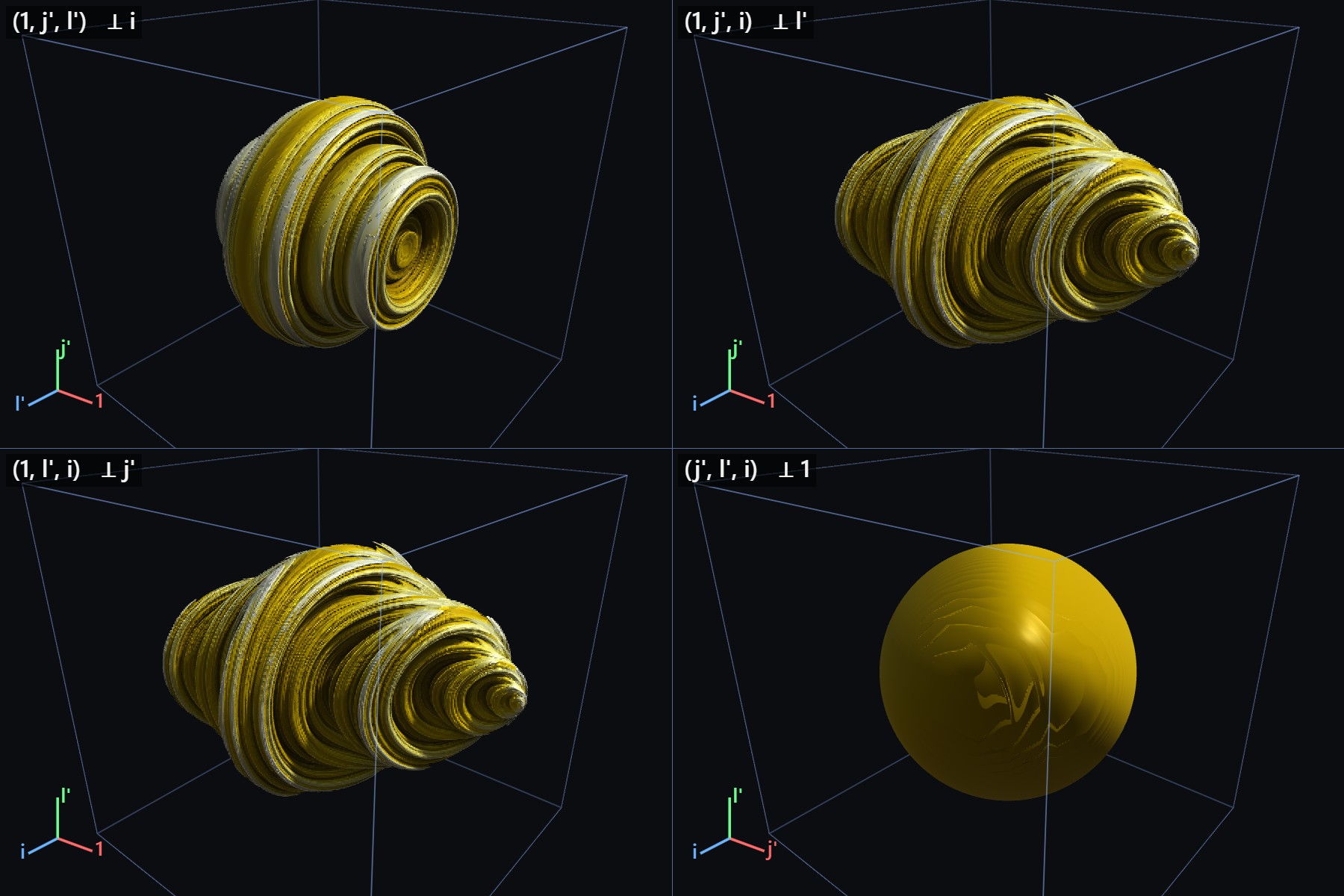}\par\vspace{1mm}
\noindent{\footnotesize \textbf{Map} $z\mapsto z\diamond z+c$ $\cdot$ \textbf{Algebra} Octonions (Cayley--Dickson, $\gamma=(-1,-1,-1)$) $\cdot$ $c=(-0.231,\,\allowbreak 0.62,\,\allowbreak 0,\,\allowbreak 0,\,\allowbreak 0,\,\allowbreak 0,\,\allowbreak 0,\,\allowbreak 0)$ ($=0.97\,s^*$, $s^*=0.6824$) $\cdot$ \textbf{Frame} $X,Y,Z,N=1,j',l',i$ (rotated by $30^\circ$ in the plane of $j,l$) $\cdot$ half-extent 1.4 $\cdot$ 20 iterations $\cdot$ \textbf{Coloring} 8D orbit trap (point), gradient Gold}\par\vspace{1mm}
\noindent{\footnotesize\begin{tabular}{@{}lllll@{}}\toprule Slice (X, Y, Z) & Body [95\%] & Exact ($J\geq0.9995$) & Approximate ($0.97\leq J<0.9995$) & Dim\\\midrule
($1, j', l'$) $\perp i$ & 22.8\% {\scriptsize[22.5, 23.0]} & $M$: xyz; $R_x$: all; $R_y$: 180; $R_z$: 180; $P$ & -- & 3 \\
($1, j', i$) $\perp l'$ & 24.8\% {\scriptsize[24.5, 25.0]} & $M$: y; $R_y$: 180; $P$ & -- & 3 \\
($1, l', i$) $\perp j'$ & 24.8\% {\scriptsize[24.5, 25.0]} & $M$: y; $R_y$: 180; $P$ & -- & 3 \\
($j', l', i$) $\perp 1$ & 34.9\% {\scriptsize[34.6, 35.2]} & $M$: xyz; $R_x$: all; $R_y$: all; $R_z$: all; $P$ & -- & 3 \\
\bottomrule\end{tabular}}\par\vspace{2mm}
Here the octonionic frame $F_C$ is rotated by $30^\circ$ in the plane of $j$ and $l$: the axes $j'$ and $l'$ are oblique combinations of $j$ and $l$. The quarters $(1,j',l')$ and $(j',l',i)$ are the same three-dimensional subspaces as $(1,j,l)$ and $(j,l,i)$, seen with axes rotated by $30^\circ$, while $(1,j',i)$ and $(1,l',i)$ are new sections. The plate is indistinguishable from Plate~C01, and so is the table: the same exact symmetries in every quarter, and body shares that agree within a few tenths of a percent, the size of the sampling error ($22.8\%$, $24.8\%$, $24.8\%$ and $34.9\%$ against $22.5\%$, $24.4\%$, $24.4\%$ and $34.7\%$).

This is a direct confirmation of Proposition~\ref{atl:prop:equivariance}: a rotation in the plane of $j$ and $l$ fixes $1$ and $c$ and commutes with the octonionic map, so the rotated frame gives the same sections, and the orbits still lie in a subalgebra of dimension $3$.
\clearpage

\subsubsection{C08. Okubo, frame rotated by \texorpdfstring{$30^\circ$}{30\textdegree{}}}\label{atl:pl:C08}\mbox{}\par\vspace{0.5mm}
\noindent\includegraphics[width=\textwidth]{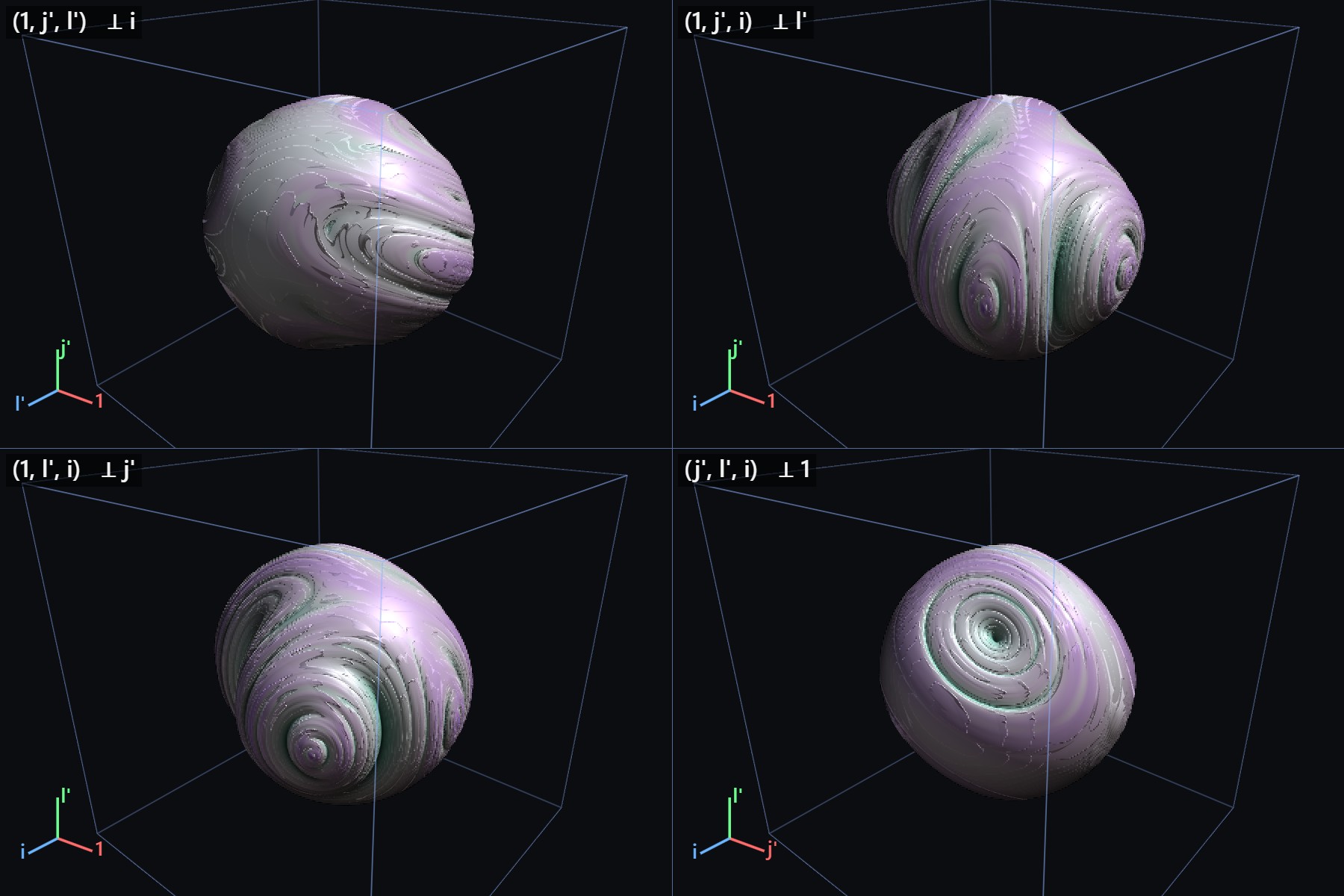}\par\vspace{1mm}
\noindent{\footnotesize \textbf{Map} $z\mapsto z\diamond z+c$ $\cdot$ \textbf{Algebra} Okubo algebra (Petersson, $\gamma=(-1,-1,-1)$, $\theta=120^\circ$) $\cdot$ $c=(-0.086,\,\allowbreak 0.231,\,\allowbreak 0,\,\allowbreak 0,\,\allowbreak 0,\,\allowbreak 0,\,\allowbreak 0,\,\allowbreak 0)$ ($=0.97\,s^*$, $s^*=0.2540$) $\cdot$ \textbf{Frame} $X,Y,Z,N=1,j',l',i$ (rotated by $30^\circ$ in the plane of $j,l$) $\cdot$ half-extent 1.4 $\cdot$ 20 iterations $\cdot$ \textbf{Coloring} 8D orbit trap (point), gradient Aurora}\par\vspace{1mm}
\noindent{\footnotesize\begin{tabular}{@{}lllll@{}}\toprule Slice (X, Y, Z) & Body [95\%] & Exact ($J\geq0.9995$) & Approximate ($0.97\leq J<0.9995$) & Dim\\\midrule
($1, j', l'$) $\perp i$ & 37.5\% {\scriptsize[37.2, 37.8]} & $M$: x; $R_x$: 180; $P$ & -- & 5 \\
($1, j', i$) $\perp l'$ & 36.4\% {\scriptsize[36.1, 36.7]} & $M$: y; $R_y$: 180; $P$ & -- & 5 \\
($1, l', i$) $\perp j'$ & 35.9\% {\scriptsize[35.6, 36.2]} & $M$: y; $R_y$: 180; $P$ & -- & 5 \\
($j', l', i$) $\perp 1$ & 33.2\% {\scriptsize[32.9, 33.5]} & $M$: z; $R_z$: 180; $P$ & -- & 5 \\
\bottomrule\end{tabular}}\par\vspace{2mm}
The same rotation applied to the Okubonic frame changes the picture. The quarters $(1,j',l')$ and $(j',l',i)$ are the sections $(1,j,l)$ and $(j,l,i)$ of Plate~C03, seen with the axes rotated by $30^\circ$ about $1$ and about $i$ respectively: their body shares are the same ($37.5\%$ and $33.2\%$ against $37.4\%$ and $33.2\%$), but these bodies are not bodies of revolution, and the images differ. Accordingly, the table keeps only the symmetries that commute with the rotation of the axes: in $(1,j',l')$ the half-turn about $1$, the mirror $M_x$ and the point reflection, in $(j',l',i)$ the half-turn about $i$, the mirror $M_z$ and the point reflection. The separate sign changes of $j$ and of $l$, which in Plate~C03 came from the half-turns of the torus in the root planes $V_1$ and $V_2$, are no longer mirrors in the new coordinate planes. The rotation about $1$ scores $J=0.960$ at $37^\circ$, as in Plate~C03.

The quarters $(1,j',i)$ and $(1,l',i)$ are genuinely new sections: $j'$ is no longer contained in $\HH$, and both quarters now generate a subalgebra of dimension $5$, like the other two. They are egg-shaped banded bodies; $(1,j',i)$ shows two vortices side by side in its lower half, $(1,l',i)$ a vortex at the lower left. Each keeps one exact mirror, the sign change of $j'$ (or $l'$), with its half-turn: it comes from the torus element that changes the signs of $j$ and $l$ together. The imaginary quarter shows a large vortex facing the viewer.

Comparing this plate with Plates~C03 and C07 makes the failure of the octonionic equivariance visible: in the Okubo algebra a rotation of the plane of $j$ and $l$ is not a symmetry, because $j$ and $l$ lie in different root planes (Proposition~\ref{atl:prop:torus}), and the choice of the section matters even inside a space orthogonal to the constant.
\clearpage

\subsection{Part D. The signature path: from octonions to split-octonions}\label{atl:part:D}
Along this path we deform the octonions by changing a single parameter, $\gamma_3=l^2$, from $-1$ (octonions) through $0$ (semi-octonions) to $+1$ (split-octonions), keeping the same constant $c\in\CC$ and the frame $F_C$. The quarter $(1,j,i)$ lies in $\HH$ and does not change at all: it is a built-in control. By Proposition~\ref{atl:prop:rescaling} the algebras with $\gamma_3<0$ are isomorphic to the octonions, and the path is visually a stretching; the transition through $\gamma_3=0$ turns ellipsoids into cylinders and cylinders into hyperboloids. Two further plates change $\gamma_1$ and $\gamma_2$ instead.

\subsubsection{D01. \texorpdfstring{$\gamma_3=-1$}{gamma3=-1}: the octonions}\label{atl:pl:D01}\mbox{}\par\vspace{0.5mm}
\noindent\includegraphics[width=\textwidth]{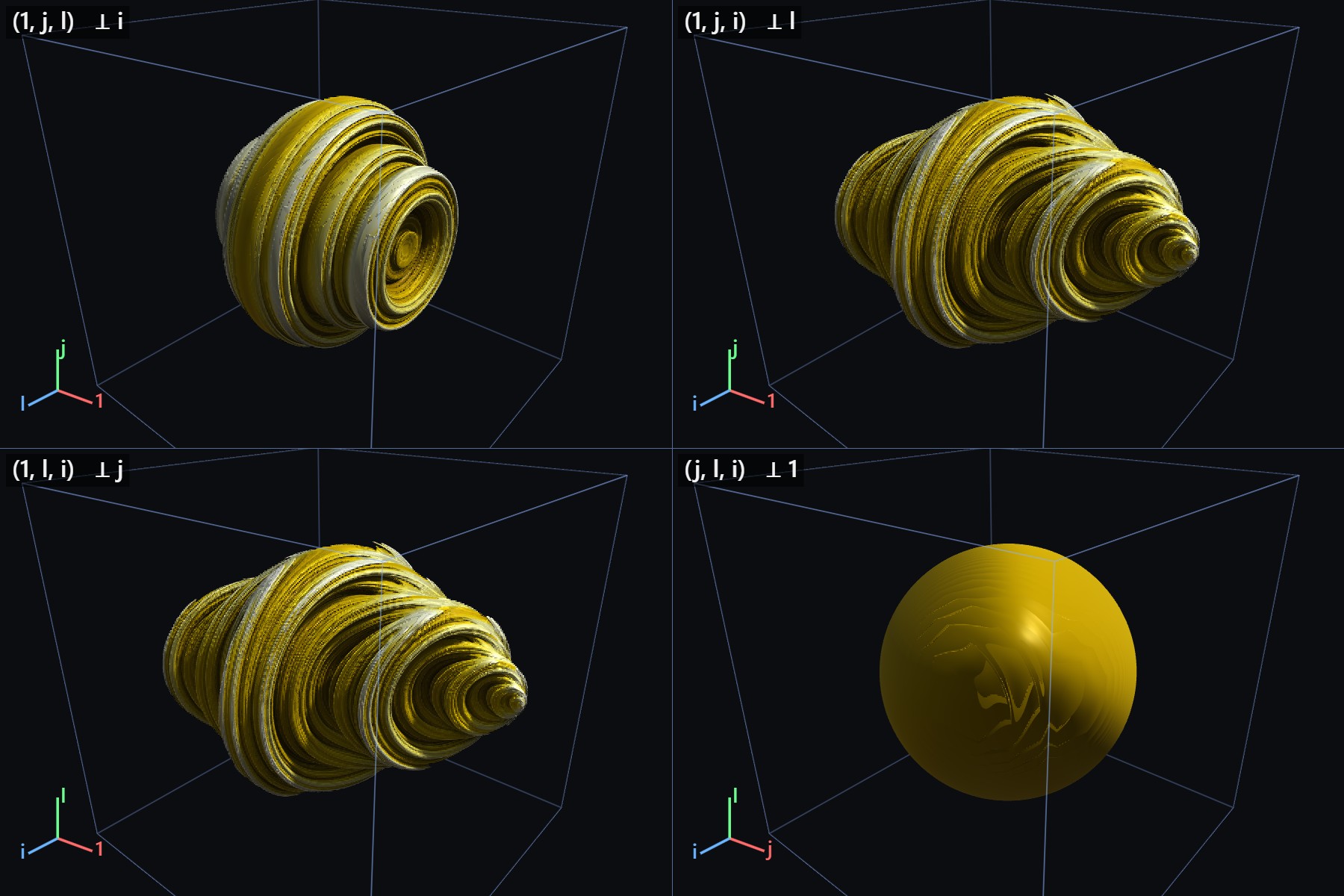}\par\vspace{1mm}
\noindent{\footnotesize \textbf{Map} $z\mapsto z\diamond z+c$ $\cdot$ \textbf{Algebra} custom: Cayley--Dickson, $\gamma=(-1,-1,-1)$ $\cdot$ $c=(-0.231,\,\allowbreak 0.62,\,\allowbreak 0,\,\allowbreak 0,\,\allowbreak 0,\,\allowbreak 0,\,\allowbreak 0,\,\allowbreak 0)$ ($=0.97\,s^*$, $s^*=0.6824$) $\cdot$ \textbf{Frame} $X,Y,Z,N=1,j,l,i$ $\cdot$ half-extent 1.4 $\cdot$ 20 iterations $\cdot$ \textbf{Coloring} 8D orbit trap (point), gradient Gold}\par\vspace{1mm}
\noindent{\footnotesize\begin{tabular}{@{}lllll@{}}\toprule Slice (X, Y, Z) & Body [95\%] & Exact ($J\geq0.9995$) & Approximate ($0.97\leq J<0.9995$) & Dim\\\midrule
($1, j, l$) $\perp i$ & 22.9\% {\scriptsize[22.6, 23.1]} & $M$: xyz; $R_x$: all; $R_y$: 180; $R_z$: 180; $P$ & -- & 3 \\
($1, j, i$) $\perp l$ & 24.6\% {\scriptsize[24.4, 24.9]} & $M$: y; $R_y$: 180; $P$ & -- & 3 \\
($1, l, i$) $\perp j$ & 24.6\% {\scriptsize[24.4, 24.9]} & $M$: y; $R_y$: 180; $P$ & -- & 3 \\
($j, l, i$) $\perp 1$ & 35.1\% {\scriptsize[34.8, 35.4]} & $M$: xyz; $R_x$: all; $R_y$: all; $R_z$: all; $P$ & -- & 3 \\
\bottomrule\end{tabular}}\par\vspace{2mm}
The first plate of the signature path is the octonionic reference in the frame $F_C$, the same as Plate~C01: an exact body of revolution about $1$ in $(1,j,l)$, two identical shells in $(1,j,i)$ and $(1,l,i)$ ($24.6\%$ each), and a ball in $(j,l,i)$. The constant is the same along the whole path, because it lies in $\CC$, where the product does not depend on $\gamma_3$. The quarter $(1,j,i)$, contained in $\HH$, is the same in all plates of the path (its body share stays between $24.5\%$ and $24.9\%$, with the same symmetries) and serves as a control.
\clearpage

\subsubsection{D02. \texorpdfstring{$\gamma_3=-0.5$}{gamma3=-0.5}: a stretched octonion algebra}\label{atl:pl:D02}\mbox{}\par\vspace{0.5mm}
\noindent\includegraphics[width=\textwidth]{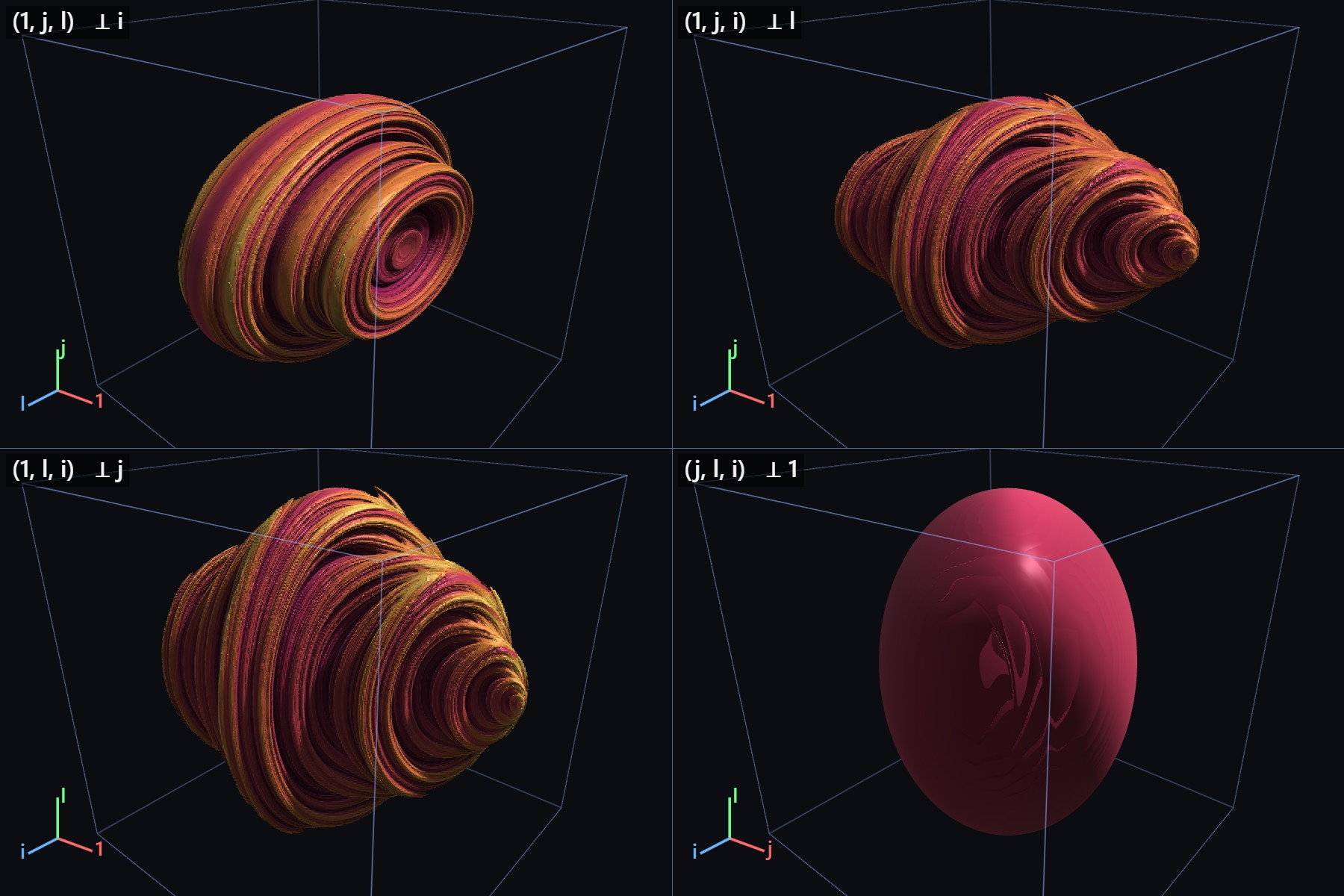}\par\vspace{1mm}
\noindent{\footnotesize \textbf{Map} $z\mapsto z\diamond z+c$ $\cdot$ \textbf{Algebra} custom: Cayley--Dickson, $\gamma=(-1,-1,-0.5)$ $\cdot$ $c=(-0.231,\,\allowbreak 0.62,\,\allowbreak 0,\,\allowbreak 0,\,\allowbreak 0,\,\allowbreak 0,\,\allowbreak 0,\,\allowbreak 0)$ ($=0.97\,s^*$, $s^*=0.6824$) $\cdot$ \textbf{Frame} $X,Y,Z,N=1,j,l,i$ $\cdot$ half-extent 1.4 $\cdot$ 20 iterations $\cdot$ \textbf{Coloring} 8D orbit trap (point), gradient Sunset}\par\vspace{1mm}
\noindent{\footnotesize\begin{tabular}{@{}lllll@{}}\toprule Slice (X, Y, Z) & Body [95\%] & Exact ($J\geq0.9995$) & Approximate ($0.97\leq J<0.9995$) & Dim\\\midrule
($1, j, l$) $\perp i$ & 32.3\% {\scriptsize[32.0, 32.5]} & $M$: xyz; $R_x$: 180; $R_y$: 180; $R_z$: 180; $P$ & -- & 3 \\
($1, j, i$) $\perp l$ & 24.8\% {\scriptsize[24.6, 25.1]} & $M$: y; $R_y$: 180; $P$ & -- & 3 \\
($1, l, i$) $\perp j$ & 35.0\% {\scriptsize[34.7, 35.3]} & $M$: y; $R_y$: 180; $P$ & -- & 3 \\
($j, l, i$) $\perp 1$ & 49.5\% {\scriptsize[49.2, 49.8]} & $M$: xyz; $R_x$: 180; $R_y$: all; $R_z$: 180; $P$ & -- & 3 \\
\bottomrule\end{tabular}}\par\vspace{2mm}
With $\gamma_3=-0.5$ the algebra is isomorphic to the octonions by the map $p+ql\mapsto p+\sqrt{0.5}\,ql$ (Proposition~\ref{atl:prop:rescaling}), and the plate shows exactly this: the quarters containing $l$ are the octonionic bodies of Plate~D01 stretched by a factor $\sqrt2\approx1.41$ along $l$. The imaginary quarter $(j,l,i)$ is a prolate spheroid elongated along $l$, still an exact body of revolution about $l$. The quarter $(1,j,l)$ keeps its three mirrors and half-turns but is no longer a body of revolution about $1$ ($J=0.70$ at $37^\circ$): a body of revolution stretched in a direction orthogonal to its axis is no longer one.

The body shares grow with the stretching ($32.3\%$, $35.0\%$ and $49.5\%$ in the quarters containing $l$, against $22.9\%$, $24.6\%$ and $35.1\%$ in Plate~D01), while the control quarter $(1,j,i)$ is unchanged ($24.8\%$).
\clearpage

\subsubsection{D03. \texorpdfstring{$\gamma_3=-0.15$}{gamma3=-0.15}: towards the cylinder}\label{atl:pl:D03}\mbox{}\par\vspace{0.5mm}
\noindent\includegraphics[width=\textwidth]{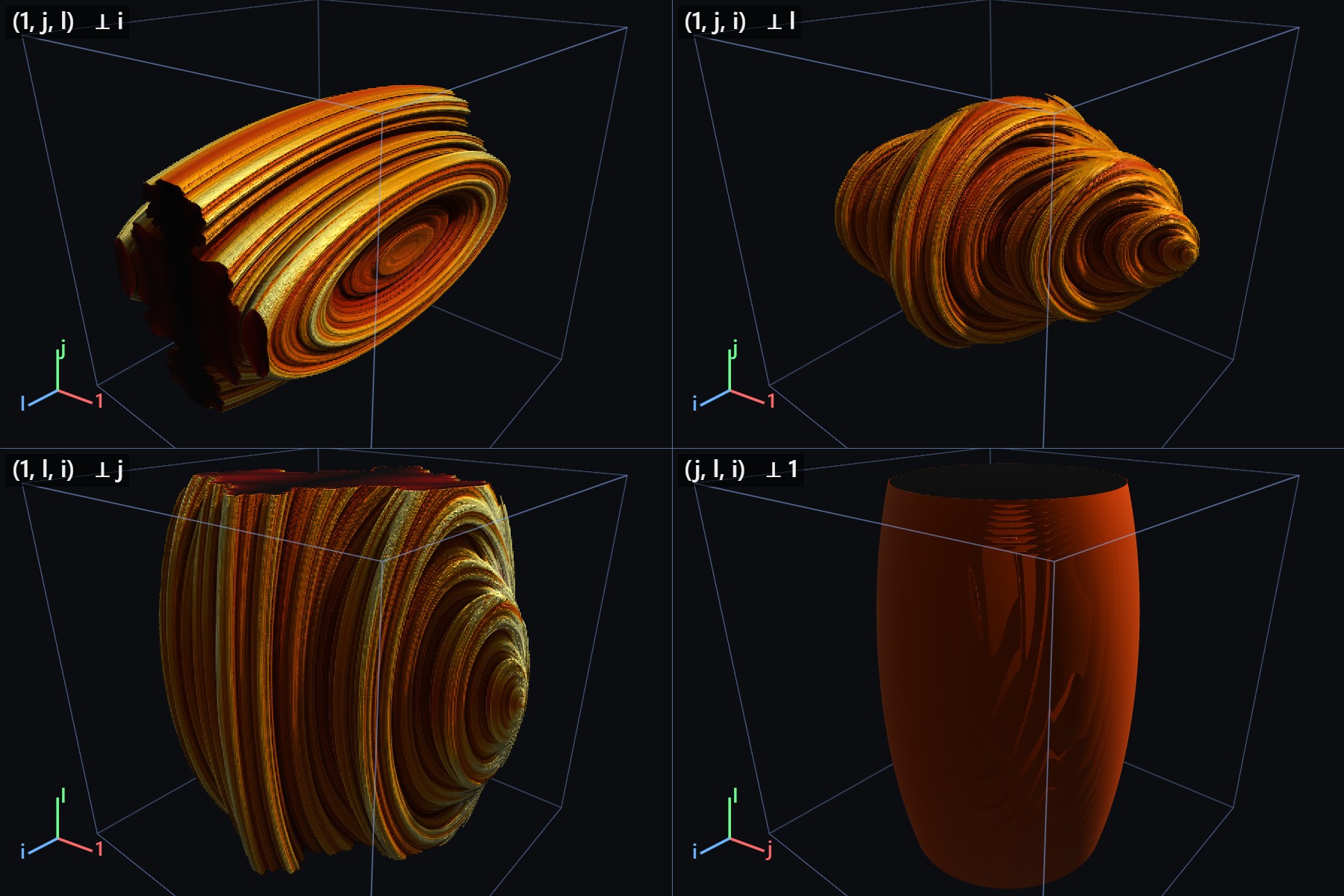}\par\vspace{1mm}
\noindent{\footnotesize \textbf{Map} $z\mapsto z\diamond z+c$ $\cdot$ \textbf{Algebra} custom: Cayley--Dickson, $\gamma=(-1,-1,-0.15)$ $\cdot$ $c=(-0.231,\,\allowbreak 0.62,\,\allowbreak 0,\,\allowbreak 0,\,\allowbreak 0,\,\allowbreak 0,\,\allowbreak 0,\,\allowbreak 0)$ ($=0.97\,s^*$, $s^*=0.6824$) $\cdot$ \textbf{Frame} $X,Y,Z,N=1,j,l,i$ $\cdot$ half-extent 1.4 $\cdot$ 20 iterations $\cdot$ \textbf{Coloring} 8D orbit trap (point), gradient Fire}\par\vspace{1mm}
\noindent{\footnotesize\begin{tabular}{@{}lllll@{}}\toprule Slice (X, Y, Z) & Body [95\%] & Exact ($J\geq0.9995$) & Approximate ($0.97\leq J<0.9995$) & Dim\\\midrule
($1, j, l$) $\perp i$ & 43.2\% {\scriptsize[42.9, 43.5]} & $M$: xyz; $R_x$: 180; $R_y$: 180; $R_z$: 180; $P$ & -- & 3 \\
($1, j, i$) $\perp l$ & 24.5\% {\scriptsize[24.2, 24.8]} & $M$: y; $R_y$: 180; $P$ & -- & 3 \\
($1, l, i$) $\perp j$ & 43.9\% {\scriptsize[43.6, 44.3]} & $M$: y; $R_y$: 180; $P$ & -- & 3 \\
($j, l, i$) $\perp 1$ & 61.0\% {\scriptsize[60.7, 61.3]} & $M$: xyz; $R_x$: 180; $R_y$: all; $R_z$: 180; $P$ & -- & 3 \\
\bottomrule\end{tabular}}\par\vspace{2mm}
With $\gamma_3=-0.15$ the stretching factor is $1/\sqrt{0.15}\approx2.58$, and the bodies no longer fit in the cube. The quarter $(1,j,l)$ is a long banded log along $l$, cut by the faces of the cube; the quarter $(1,l,i)$ is a banded column along $l$, cut by the top face; the imaginary quarter is a convex barrel with axis $l$, cut flat at the top, and an exact body of revolution about $l$. The body shares grow accordingly ($43.2\%$, $43.9\%$ and $61.0\%$). The algebra is still isomorphic to the octonions; the picture is approaching the degenerate cylinder of $\gamma_3=0$.
\clearpage

\subsubsection{D04. \texorpdfstring{$\gamma_3=0$}{gamma3=0}: the semi-octonions}\label{atl:pl:D04}\mbox{}\par\vspace{0.5mm}
\noindent\includegraphics[width=\textwidth]{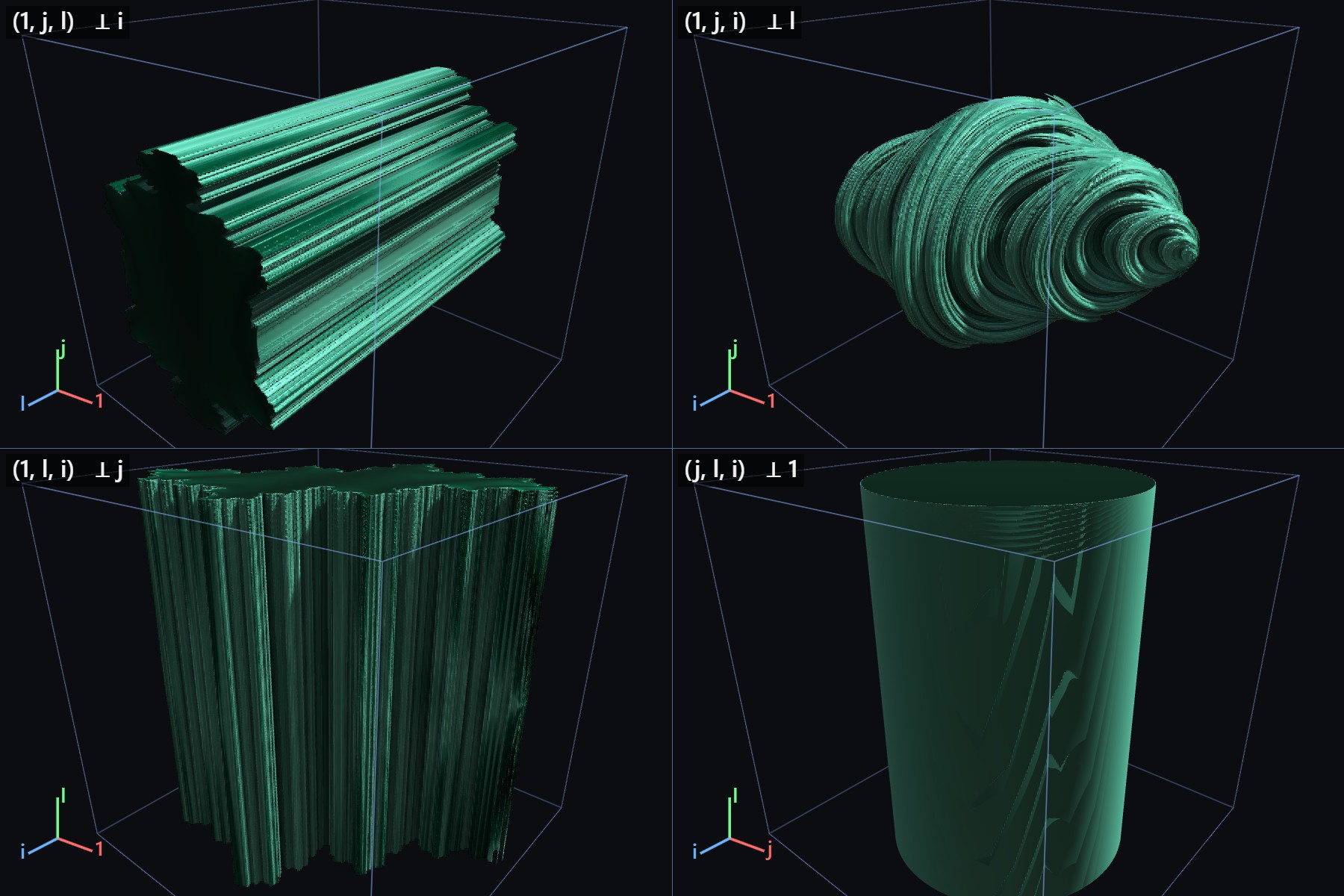}\par\vspace{1mm}
\noindent{\footnotesize \textbf{Map} $z\mapsto z\diamond z+c$ $\cdot$ \textbf{Algebra} custom: Cayley--Dickson, $\gamma=(-1,-1,0)$ $\cdot$ $c=(-0.231,\,\allowbreak 0.62,\,\allowbreak 0,\,\allowbreak 0,\,\allowbreak 0,\,\allowbreak 0,\,\allowbreak 0,\,\allowbreak 0)$ ($=0.97\,s^*$, $s^*=0.6824$) $\cdot$ \textbf{Frame} $X,Y,Z,N=1,j,l,i$ $\cdot$ half-extent 1.4 $\cdot$ 20 iterations $\cdot$ \textbf{Coloring} 8D orbit trap (point), gradient Emerald}\par\vspace{1mm}
\noindent{\footnotesize\begin{tabular}{@{}lllll@{}}\toprule Slice (X, Y, Z) & Body [95\%] & Exact ($J\geq0.9995$) & Approximate ($0.97\leq J<0.9995$) & Dim\\\midrule
($1, j, l$) $\perp i$ & 45.1\% {\scriptsize[44.7, 45.4]} & $M$: xyz; $R_x$: 180; $R_y$: 180; $R_z$: 180; $P$ & -- & 3 \\
($1, j, i$) $\perp l$ & 24.5\% {\scriptsize[24.2, 24.7]} & $M$: y; $R_y$: 180; $P$ & -- & 3 \\
($1, l, i$) $\perp j$ & 45.3\% {\scriptsize[45.0, 45.6]} & $M$: y; $R_y$: 180; $P$ & -- & 3 \\
($j, l, i$) $\perp 1$ & 64.1\% {\scriptsize[63.8, 64.4]} & $M$: xyz; $R_x$: 180; $R_y$: all; $R_z$: 180; $P$ & -- & 3 \\
\bottomrule\end{tabular}}\par\vspace{2mm}
At $\gamma_3=0$ the algebra is the semi-octonions, and the stretching becomes infinite. The quarters $(1,j,l)$ and $(1,l,i)$ are fluted columns along $l$, cut by the faces of the cube, and the imaginary quarter is a circular cylinder with axis $l$, an exact body of revolution about $l$. The cylinder is explained by Proposition~\ref{atl:prop:rescaling}: for $z=p+ql$ with $\re p=0$ one has $z^2=p^2$, so the first iterate forgets the $l$-component altogether. The body shares ($45.1\%$, $45.3\%$ and $64.1\%$) are close to those of Plate~D03. The plate is the midpoint of the path, where the ellipsoids of the octonions turn into the hyperboloids of the split-octonions.
\clearpage

\subsubsection{D05. \texorpdfstring{$\gamma_3=+0.15$}{gamma3=+0.15}: the signature flips}\label{atl:pl:D05}\mbox{}\par\vspace{0.5mm}
\noindent\includegraphics[width=\textwidth]{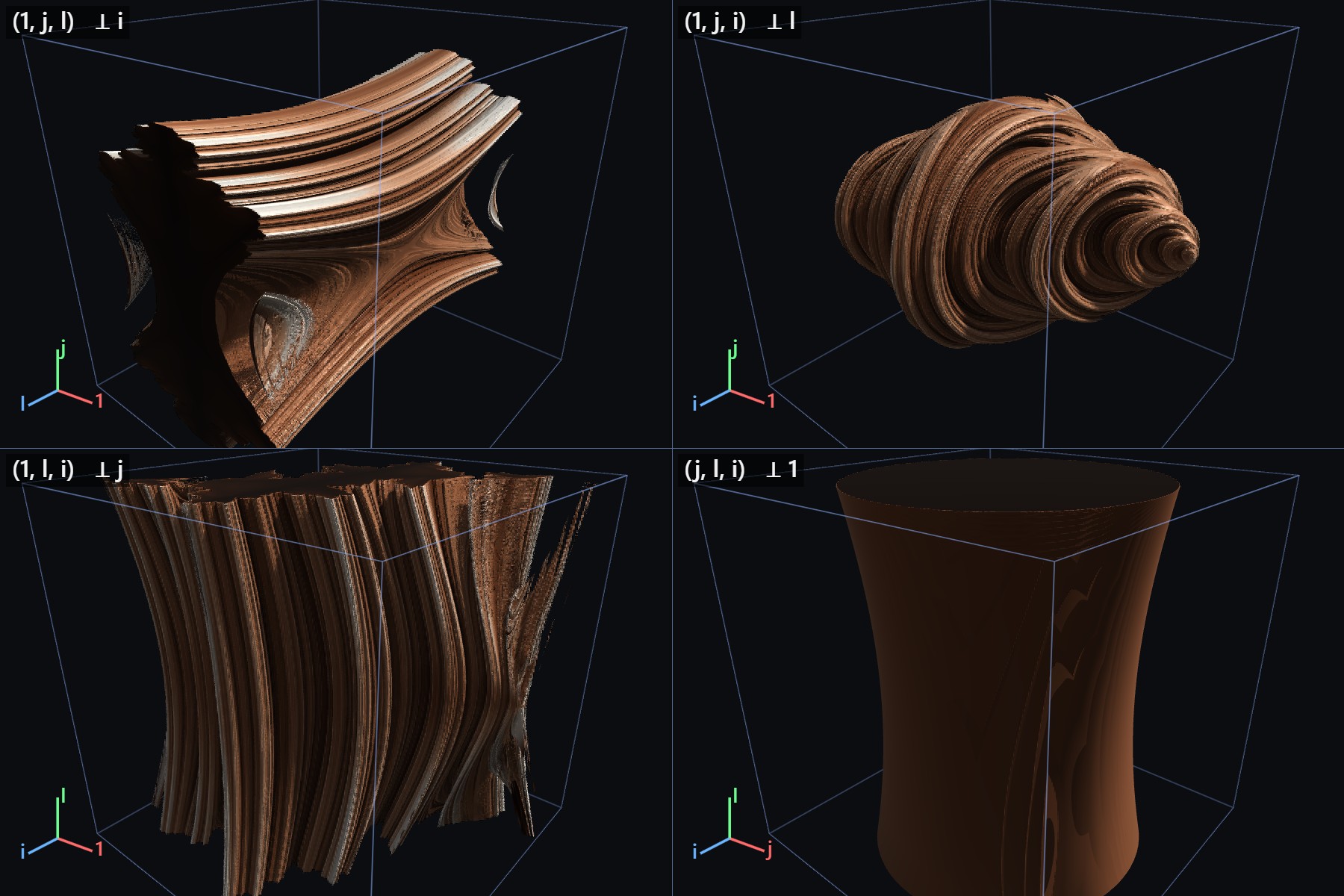}\par\vspace{1mm}
\noindent{\footnotesize \textbf{Map} $z\mapsto z\diamond z+c$ $\cdot$ \textbf{Algebra} custom: Cayley--Dickson, $\gamma=(-1,-1,0.15)$ $\cdot$ $c=(-0.231,\,\allowbreak 0.62,\,\allowbreak 0,\,\allowbreak 0,\,\allowbreak 0,\,\allowbreak 0,\,\allowbreak 0,\,\allowbreak 0)$ ($=0.97\,s^*$, $s^*=0.6824$) $\cdot$ \textbf{Frame} $X,Y,Z,N=1,j,l,i$ $\cdot$ half-extent 1.4 $\cdot$ 20 iterations $\cdot$ \textbf{Coloring} 8D orbit trap (point), gradient Copper}\par\vspace{1mm}
\noindent{\footnotesize\begin{tabular}{@{}lllll@{}}\toprule Slice (X, Y, Z) & Body [95\%] & Exact ($J\geq0.9995$) & Approximate ($0.97\leq J<0.9995$) & Dim\\\midrule
($1, j, l$) $\perp i$ & 45.0\% {\scriptsize[44.7, 45.3]} & $M$: xyz; $R_x$: 180; $R_y$: 180; $R_z$: 180; $P$ & -- & 3 \\
($1, j, i$) $\perp l$ & 24.5\% {\scriptsize[24.2, 24.8]} & $M$: y; $R_y$: 180; $P$ & -- & 3 \\
($1, l, i$) $\perp j$ & 45.1\% {\scriptsize[44.8, 45.4]} & $M$: y; $R_y$: 180; $P$ & -- & 3 \\
($j, l, i$) $\perp 1$ & 66.7\% {\scriptsize[66.4, 67.0]} & $M$: xyz; $R_x$: 180; $R_y$: all; $R_z$: 180; $P$ & -- & 3 \\
\bottomrule\end{tabular}}\par\vspace{2mm}
Just beyond the degenerate point, with $\gamma_3=+0.15$, the norm becomes indefinite, and the algebra is a split-octonion algebra stretched along $l$ (by the substitution of Proposition~\ref{atl:prop:rescaling} with $\sqrt{\gamma_3}$). The picture changes in a subtle but unmistakable way: the convex barrel of Plate~D03 has become a waisted spool with concave flanks, bounded by a hyperboloid of one sheet with axis $l$ and still an exact body of revolution about $l$, since by Proposition~\ref{atl:prop:imaginary} the section is a region between level sets of an indefinite norm. The quarter $(1,j,l)$ turns into a saddle whose banded sheets curve away from each other, with small separate pieces near the faces of the cube, and the quarter $(1,l,i)$ into a column of curved sheets with concave flanks. The signature flip is visible as the passage from convex to concave, while the body shares ($45.0\%$, $45.1\%$ and $66.7\%$) hardly change from Plate~D04.
\clearpage

\subsubsection{D06. \texorpdfstring{$\gamma_3=+0.5$}{gamma3=+0.5}: saddles and hyperboloids}\label{atl:pl:D06}\mbox{}\par\vspace{0.5mm}
\noindent\includegraphics[width=\textwidth]{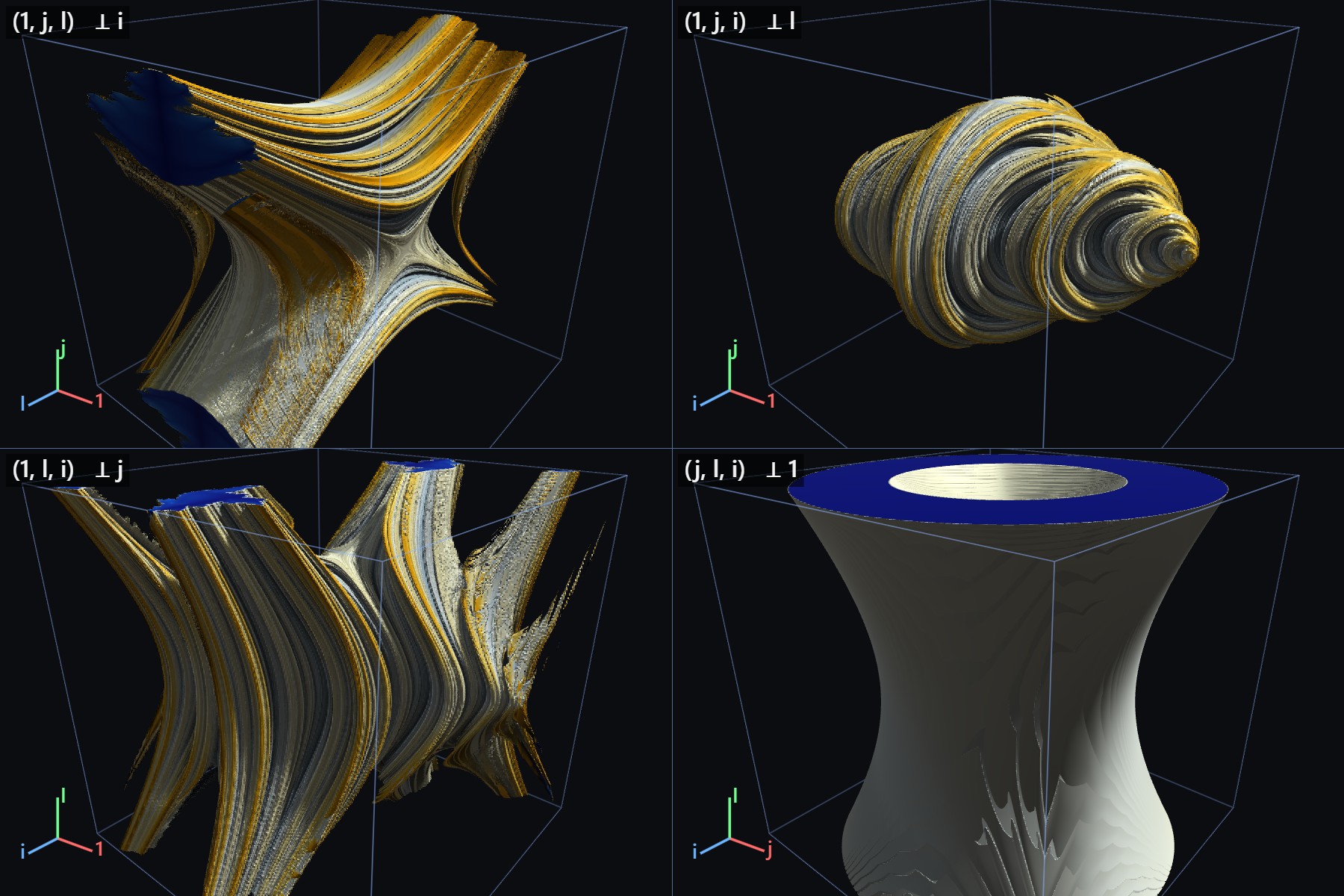}\par\vspace{1mm}
\noindent{\footnotesize \textbf{Map} $z\mapsto z\diamond z+c$ $\cdot$ \textbf{Algebra} custom: Cayley--Dickson, $\gamma=(-1,-1,0.5)$ $\cdot$ $c=(-0.231,\,\allowbreak 0.62,\,\allowbreak 0,\,\allowbreak 0,\,\allowbreak 0,\,\allowbreak 0,\,\allowbreak 0,\,\allowbreak 0)$ ($=0.97\,s^*$, $s^*=0.6824$) $\cdot$ \textbf{Frame} $X,Y,Z,N=1,j,l,i$ $\cdot$ half-extent 1.4 $\cdot$ 20 iterations $\cdot$ \textbf{Coloring} 8D orbit trap (point), gradient Classic}\par\vspace{1mm}
\noindent{\footnotesize\begin{tabular}{@{}lllll@{}}\toprule Slice (X, Y, Z) & Body [95\%] & Exact ($J\geq0.9995$) & Approximate ($0.97\leq J<0.9995$) & Dim\\\midrule
($1, j, l$) $\perp i$ & 37.3\% {\scriptsize[37.0, 37.6]} & $M$: xyz; $R_x$: 180; $R_y$: 180; $R_z$: 180; $P$ & -- & 3 \\
($1, j, i$) $\perp l$ & 24.6\% {\scriptsize[24.4, 24.9]} & $M$: y; $R_y$: 180; $P$ & -- & 3 \\
($1, l, i$) $\perp j$ & 35.3\% {\scriptsize[35.0, 35.6]} & $M$: y; $R_y$: 180; $P$ & -- & 3 \\
($j, l, i$) $\perp 1$ & 66.2\% {\scriptsize[65.9, 66.5]} & $M$: xyz; $R_x$: 180; $R_y$: all; $R_z$: 180; $P$ & -- & 3 \\
\bottomrule\end{tabular}}\par\vspace{2mm}
With $\gamma_3=+0.5$ the hyperbolic geometry is fully developed; the algebra is isomorphic to the split-octonions, stretched along $l$ by $\sqrt2$. The quarter $(1,j,l)$ is an X-shaped saddle whose banded sheets leave the cube; the quarter $(1,l,i)$ shows bundles of curved sheets that leave the cube through its top and bottom faces. The imaginary quarter is a waisted vase with axis $l$, an exact body of revolution about $l$, and its cut by the top face of the cube is a ring: near that face the body is hollowed out around its axis. This is what Proposition~\ref{atl:prop:imaginary} predicts for an indefinite norm, the section being the region between two level sets, a hyperboloid of one sheet outside and, near the axis, one of two sheets. The control quarter $(1,j,i)$ is the familiar quaternionic shell ($24.6\%$).
\clearpage

\subsubsection{D07. \texorpdfstring{$\gamma_3=+1$}{gamma3=+1}: the split-octonions}\label{atl:pl:D07}\mbox{}\par\vspace{0.5mm}
\noindent\includegraphics[width=\textwidth]{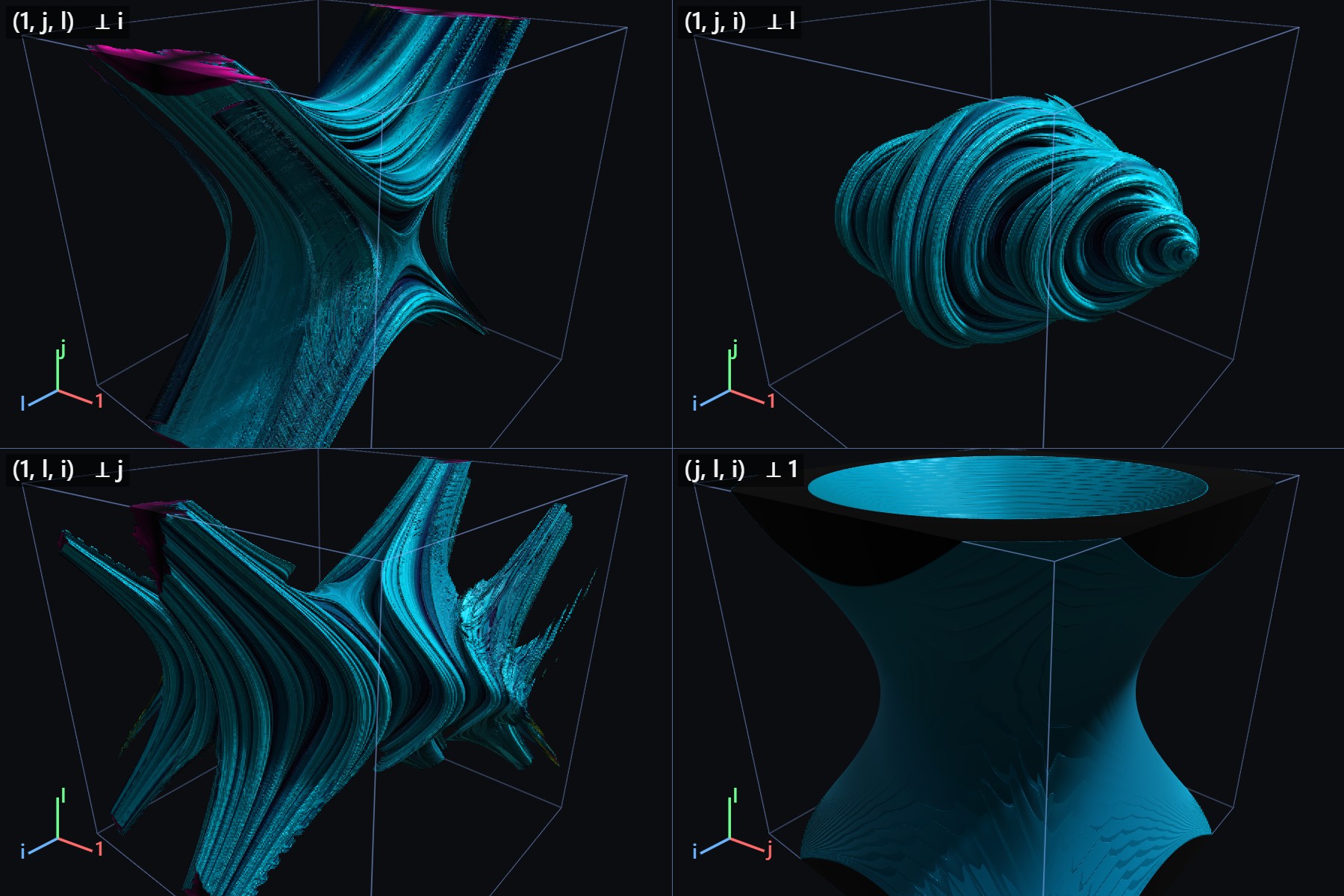}\par\vspace{1mm}
\noindent{\footnotesize \textbf{Map} $z\mapsto z\diamond z+c$ $\cdot$ \textbf{Algebra} custom: Cayley--Dickson, $\gamma=(-1,-1,1)$ $\cdot$ $c=(-0.231,\,\allowbreak 0.62,\,\allowbreak 0,\,\allowbreak 0,\,\allowbreak 0,\,\allowbreak 0,\,\allowbreak 0,\,\allowbreak 0)$ ($=0.97\,s^*$, $s^*=0.6824$) $\cdot$ \textbf{Frame} $X,Y,Z,N=1,j,l,i$ $\cdot$ half-extent 1.4 $\cdot$ 20 iterations $\cdot$ \textbf{Coloring} 8D orbit trap (point), gradient Neon}\par\vspace{1mm}
\noindent{\footnotesize\begin{tabular}{@{}lllll@{}}\toprule Slice (X, Y, Z) & Body [95\%] & Exact ($J\geq0.9995$) & Approximate ($0.97\leq J<0.9995$) & Dim\\\midrule
($1, j, l$) $\perp i$ & 30.3\% {\scriptsize[30.0, 30.6]} & $M$: xyz; $R_x$: 180; $R_y$: 180; $R_z$: 180; $P$ & -- & 3 \\
($1, j, i$) $\perp l$ & 24.9\% {\scriptsize[24.6, 25.2]} & $M$: y; $R_y$: 180; $P$ & -- & 3 \\
($1, l, i$) $\perp j$ & 27.4\% {\scriptsize[27.2, 27.7]} & $M$: y; $R_y$: 180; $P$ & -- & 3 \\
($j, l, i$) $\perp 1$ & 59.8\% {\scriptsize[59.5, 60.1]} & $M$: xyz; $R_x$: 180; $R_y$: all; $R_z$: 180; $P$ & -- & 3 \\
\bottomrule\end{tabular}}\par\vspace{2mm}
The path ends at the split-octonions, $\gamma_3=+1$. The four quarters are the sections of Plate~A02, drawn in the frame $F_C$: X-shaped saddles and bundles of curved sheets in $(1,j,l)$ and $(1,l,i)$, the untouched shell in $(1,j,i)$, and in the imaginary quarter a waisted spool with axis $l$, an exact body of revolution about $l$, whose top face shows a bowl-shaped cavity around the axis. The imaginary quarter coincides with that of Plate~D09, with the same body share ($59.8\%$) and the same symmetries, although there the hyperbolic direction is $j$ instead of $l$: by Proposition~\ref{atl:prop:imaginary} the section depends only on the norm restricted to it.

Looking back along Plates~D01--D07, the path is a continuous deformation of a single picture: stretching, degeneration into a cylinder and opening into hyperboloids, with the orbits always in a subalgebra of dimension $3$.
\clearpage

\subsubsection{D08. \texorpdfstring{$\gamma_1=0$}{gamma1=0}: a dual plane of \texorpdfstring{$1$}{1} and \texorpdfstring{$i$}{i}}\label{atl:pl:D08}\mbox{}\par\vspace{0.5mm}
\noindent\includegraphics[width=\textwidth]{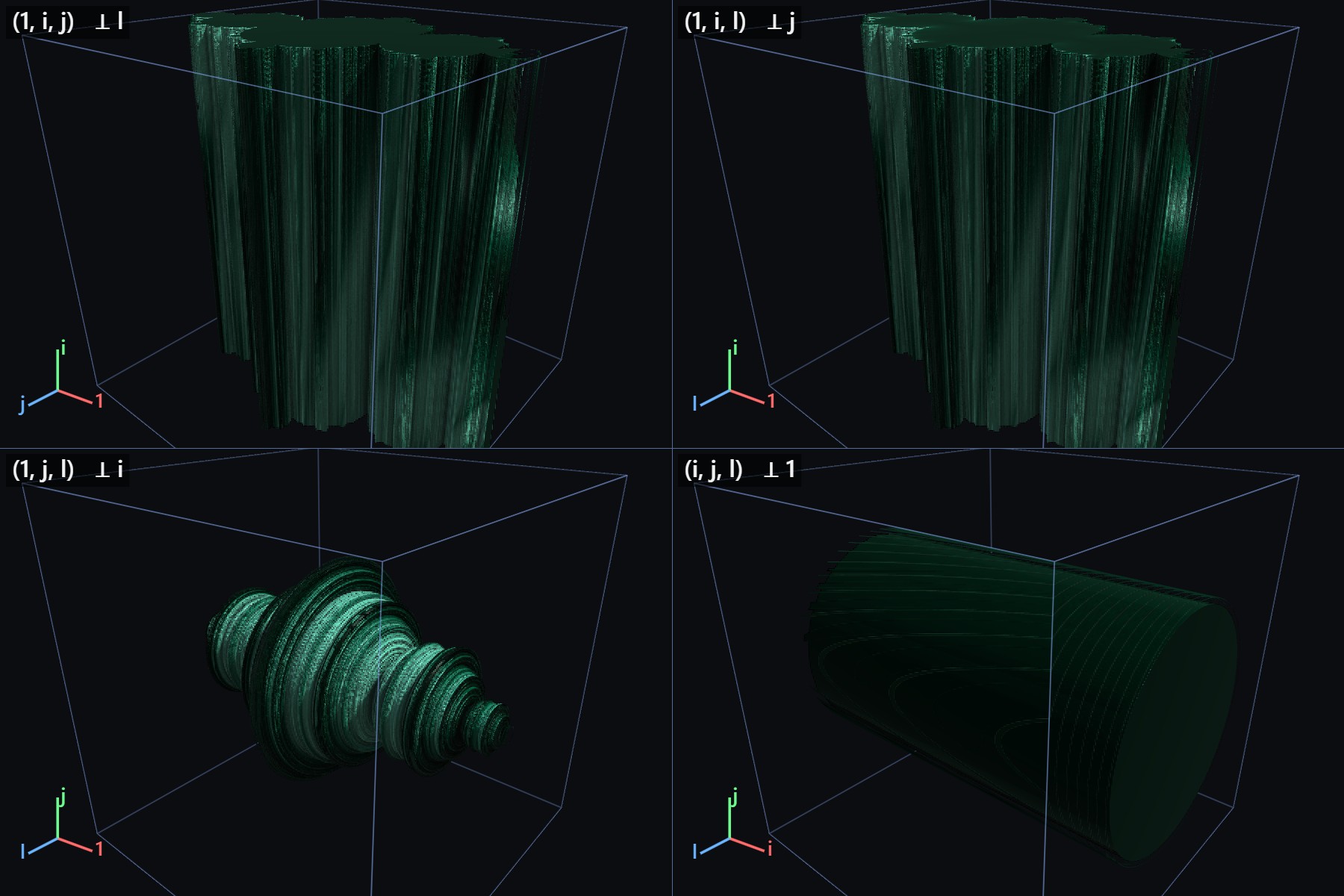}\par\vspace{1mm}
\noindent{\footnotesize \textbf{Map} $z\mapsto z\diamond z+c$ $\cdot$ \textbf{Algebra} custom: Cayley--Dickson, $\gamma=(0,-1,-1)$ $\cdot$ $c=(-0.724,\,\allowbreak 1.945,\,\allowbreak 0,\,\allowbreak 0,\,\allowbreak 0,\,\allowbreak 0,\,\allowbreak 0,\,\allowbreak 0)$ ($=0.97\,s^*$, $s^*=2.1401$) $\cdot$ \textbf{Frame} $X,Y,Z,N=1,i,j,l$ $\cdot$ half-extent 1.4 $\cdot$ 20 iterations $\cdot$ \textbf{Coloring} 8D orbit trap (point), gradient Emerald}\par\vspace{1mm}
\noindent{\footnotesize\begin{tabular}{@{}lllll@{}}\toprule Slice (X, Y, Z) & Body [95\%] & Exact ($J\geq0.9995$) & Approximate ($0.97\leq J<0.9995$) & Dim\\\midrule
($1, i, j$) $\perp l$ & 48.6\% {\scriptsize[48.3, 48.9]} & $M$: z; $R_z$: 180; $P$ & $M$: xy; $R_x$: 180; $R_y$: 180 (min.\ 0.999) & 3 \\
($1, i, l$) $\perp j$ & 48.6\% {\scriptsize[48.3, 48.9]} & $M$: z; $R_z$: 180; $P$ & $M$: xy; $R_x$: 180; $R_y$: 180 (min.\ 0.999) & 3 \\
($1, j, l$) $\perp i$ & 17.0\% {\scriptsize[16.7, 17.2]} & $M$: xyz; $R_x$: all; $R_y$: 180; $R_z$: 180; $P$ & -- & 3 \\
($i, j, l$) $\perp 1$ & 51.9\% {\scriptsize[51.6, 52.2]} & $M$: xyz; $R_x$: all; $R_y$: 180; $R_z$: 180; $P$ & -- & 3 \\
\bottomrule\end{tabular}}\par\vspace{2mm}
In this plate $\gamma_1=i^2=0$: the plane of $1$ and $i$ is the plane of the dual numbers. On the ray of the complex direction the boundary distance is $s^*=2.14$, very large, because in the dual plane the real part evolves as a real quadratic map and the dual part as its parameter derivative (Remark~\ref{atl:rem:dual}). The critical orbit escapes only in a small window around the period-doubling point $\re c=-3/4$; the chosen constant, $c=-0.724+1.945\,i$, lies just before that window.

The quarters $(1,i,j)$ and $(1,i,l)$ are identical columns along $i$, the dual direction, cut by the top face of the cube ($48.6\%$ each); besides the exact mirror in $j$ (or $l$) they are only approximately symmetric under the mirror in $i$ ($J=0.999$). The quarter $(1,j,l)$, orthogonal to $i$, is a small banded body ($17.0\%$), an exact body of revolution about $1$, since $j$ and $l$ are still elliptic and the rotations of their plane fix $1$ and $c$ (Proposition~\ref{atl:prop:equivariance}). The imaginary quarter $(i,j,l)$ is a circular cylinder along $i$, an exact body of revolution about $i$: on purely imaginary points the norm is $x_2^2+x_4^2$, which does not see the $i$-component, and by Proposition~\ref{atl:prop:imaginary} the section is a union of its level sets. The orbits lie in a subalgebra of dimension $3$.
\clearpage

\subsubsection{D09. \texorpdfstring{$\gamma_2=+1$}{gamma2=+1}: split-quaternions inside}\label{atl:pl:D09}\mbox{}\par\vspace{0.5mm}
\noindent\includegraphics[width=\textwidth]{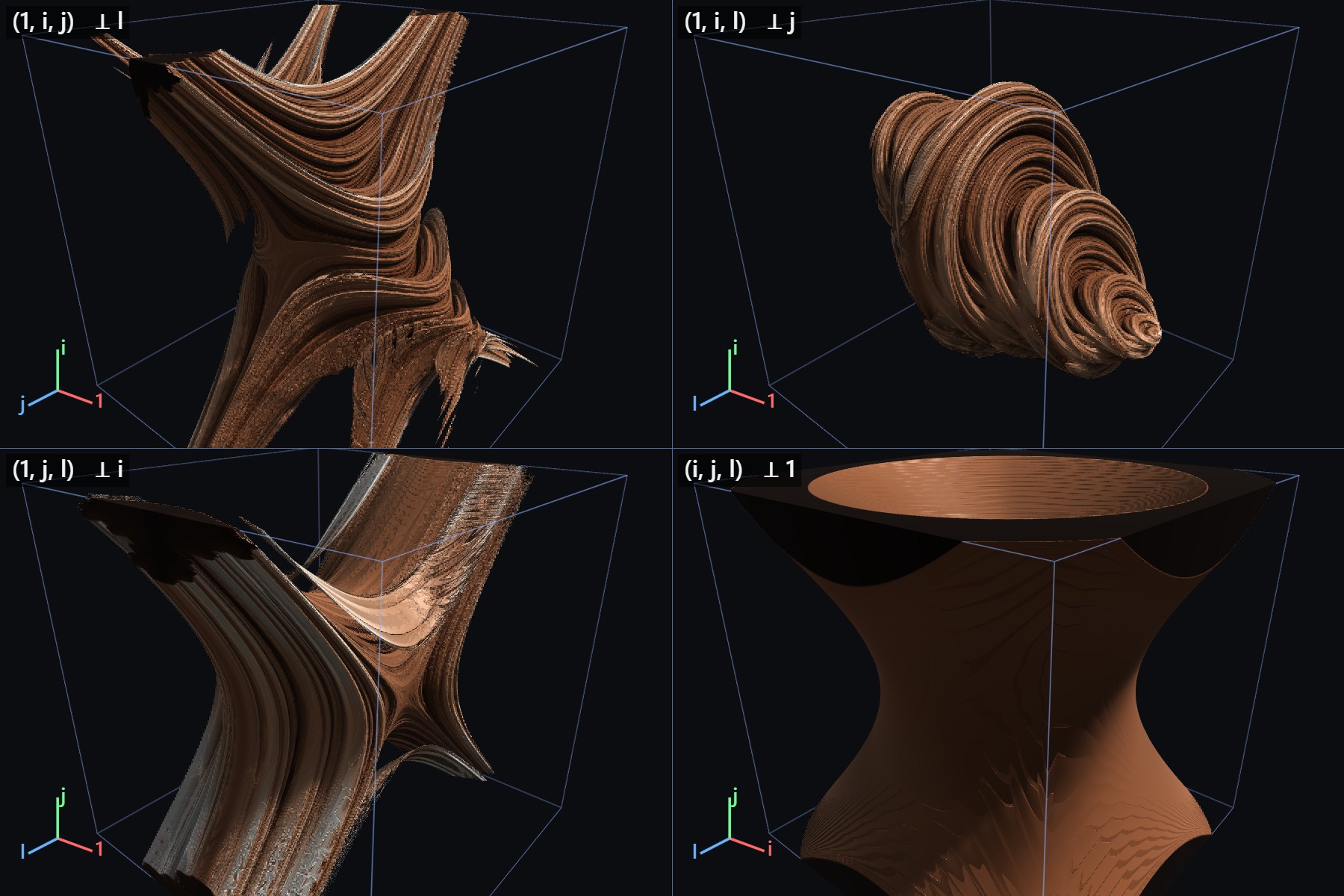}\par\vspace{1mm}
\noindent{\footnotesize \textbf{Map} $z\mapsto z\diamond z+c$ $\cdot$ \textbf{Algebra} custom: Cayley--Dickson, $\gamma=(-1,1,-1)$ $\cdot$ $c=(-0.231,\,\allowbreak 0.62,\,\allowbreak 0,\,\allowbreak 0,\,\allowbreak 0,\,\allowbreak 0,\,\allowbreak 0,\,\allowbreak 0)$ ($=0.97\,s^*$, $s^*=0.6824$) $\cdot$ \textbf{Frame} $X,Y,Z,N=1,i,j,l$ $\cdot$ half-extent 1.4 $\cdot$ 20 iterations $\cdot$ \textbf{Coloring} 8D orbit trap (point), gradient Copper}\par\vspace{1mm}
\noindent{\footnotesize\begin{tabular}{@{}lllll@{}}\toprule Slice (X, Y, Z) & Body [95\%] & Exact ($J\geq0.9995$) & Approximate ($0.97\leq J<0.9995$) & Dim\\\midrule
($1, i, j$) $\perp l$ & 27.3\% {\scriptsize[27.0, 27.6]} & $M$: z; $R_z$: 180; $P$ & -- & 3 \\
($1, i, l$) $\perp j$ & 24.8\% {\scriptsize[24.5, 25.0]} & $M$: z; $R_z$: 180; $P$ & -- & 3 \\
($1, j, l$) $\perp i$ & 30.2\% {\scriptsize[30.0, 30.5]} & $M$: xyz; $R_x$: 180; $R_y$: 180; $R_z$: 180; $P$ & -- & 3 \\
($i, j, l$) $\perp 1$ & 59.8\% {\scriptsize[59.5, 60.1]} & $M$: xyz; $R_x$: 180; $R_y$: all; $R_z$: 180; $P$ & -- & 3 \\
\bottomrule\end{tabular}}\par\vspace{2mm}
With $\gamma_2=j^2=+1$ the quaternion subalgebra $\langle1,i,j,k\rangle$ becomes the split-quaternions. The constant is complex and $i^2=-1$, so its orbit stays in $\CC$ and $s^*$ is the octonionic value. The quarter $(1,i,j)$ is now a hyperbolic saddle whose sheets leave the cube, with a single exact mirror (the sign change of $j$), its half-turn and the point reflection. The quarter $(1,i,l)$ lies in the subalgebra $\langle1,i,l,il\rangle$, which is still the ordinary quaternion algebra, and it is the usual octonionic shell ($24.8\%$). The quarter $(1,j,l)$ is a saddle with its three mirrors and half-turns but no revolution, and the imaginary quarter $(i,j,l)$ is a waisted spool with a cavity at the top, an exact body of revolution about the local $y$ axis, the hyperbolic direction $j$.

The plate is the counterpart of Plate~D07: there the hyperbolic direction was $l$, here it is $j$, and the imaginary quarters of the two plates coincide (body share $59.8\%$, same symmetries), as Proposition~\ref{atl:prop:imaginary} predicts. The quaternionic control quarter of the previous plates is deformed here, because the quaternion subalgebra has itself become split.
\clearpage

\subsection{Part E. The Petersson path: from para-octonions to Okubo and beyond}\label{atl:part:E}
The Petersson product $\tau_\theta(\bar x)\,\tau_\theta^{-1}(\bar y)$ interpolates continuously between the para-octonions ($\theta=0$) and the Okubo algebra ($\theta=120^\circ$), and continues to $\theta=180^\circ$. We follow it in seven steps, with the generic constant direction $d_{\mathcal O}$ and the frame $F_B$. The boundary distance $s^*$ and the dimension of the orbits change along the path, and one angle, $\theta=90^\circ$, is exceptional.

\subsubsection{E01. \texorpdfstring{$\theta=0^\circ$}{theta=0\textdegree{}}: the para-octonions}\label{atl:pl:E01}\mbox{}\par\vspace{0.5mm}
\noindent\includegraphics[width=\textwidth]{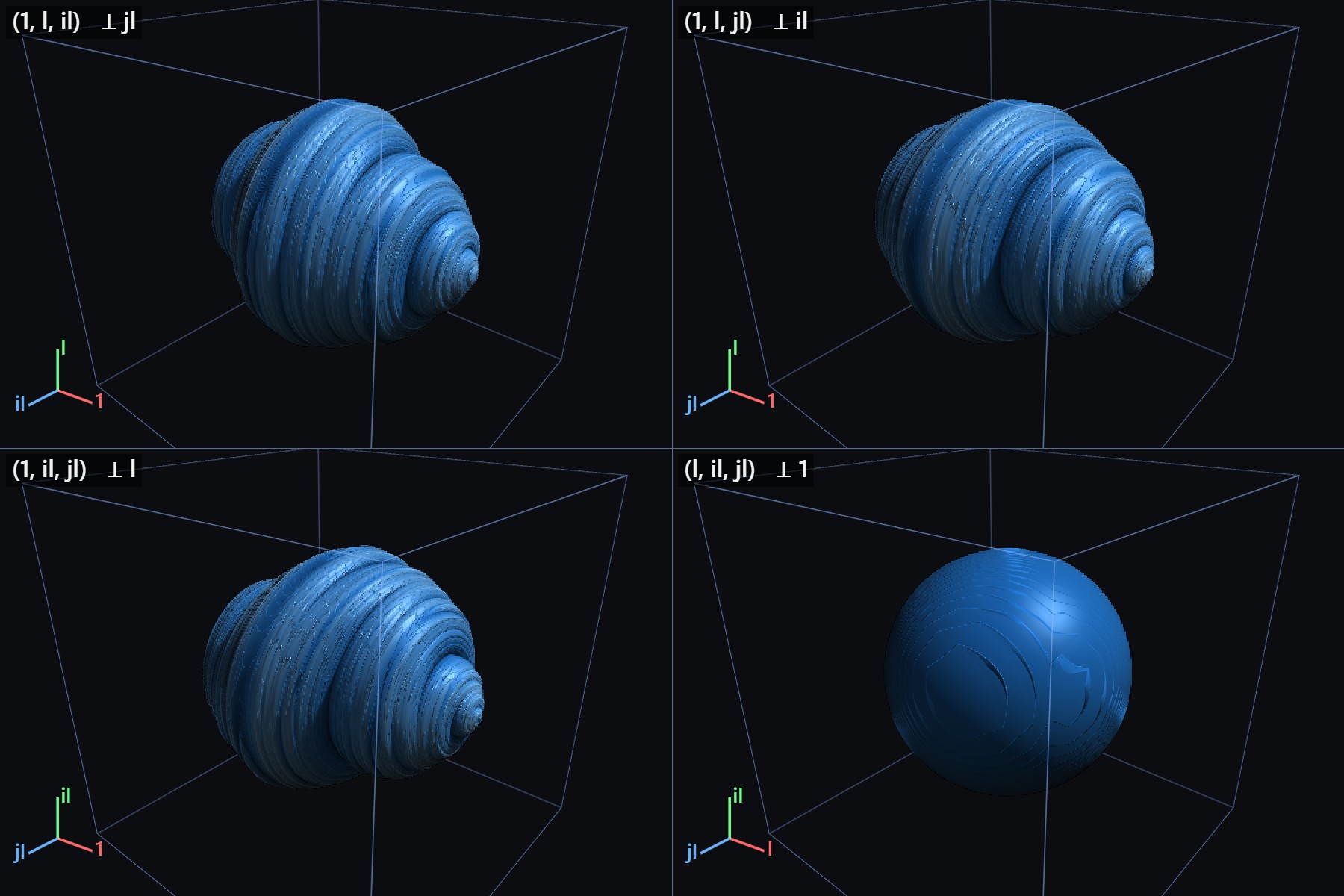}\par\vspace{1mm}
\noindent{\footnotesize \textbf{Map} $z\mapsto z\diamond z+c$ $\cdot$ \textbf{Algebra} custom: Petersson, $\gamma=(-1,-1,-1)$, $\theta=0^\circ$ $\cdot$ $c=(-0.235,\,\allowbreak -0.013,\,\allowbreak 0.029,\,\allowbreak 0.01,\,\allowbreak 0.002,\,\allowbreak 0.095,\,\allowbreak -0.067,\,\allowbreak 0.079)$ ($=0.95\,s^*$, $s^*=0.2901$) $\cdot$ \textbf{Frame} $X,Y,Z,N=1,l,il,jl$ $\cdot$ half-extent 1.4 $\cdot$ 20 iterations $\cdot$ \textbf{Coloring} 8D orbit trap (point), gradient Deep Blue}\par\vspace{1mm}
\noindent{\footnotesize\begin{tabular}{@{}lllll@{}}\toprule Slice (X, Y, Z) & Body [95\%] & Exact ($J\geq0.9995$) & Approximate ($0.97\leq J<0.9995$) & Dim\\\midrule
($1, l, il$) $\perp jl$ & 32.6\% {\scriptsize[32.3, 32.9]} & $P$ & $M$: y; $R_y$: 180 (min.\ 0.994) & 3 \\
($1, l, jl$) $\perp il$ & 32.6\% {\scriptsize[32.3, 32.9]} & $P$ & $M$: y; $R_y$: 180 (min.\ 0.995) & 3 \\
($1, il, jl$) $\perp l$ & 33.0\% {\scriptsize[32.7, 33.2]} & $P$ & -- & 3 \\
($l, il, jl$) $\perp 1$ & 31.9\% {\scriptsize[31.6, 32.1]} & $M$: xyz; $R_x$: all; $R_y$: all; $R_z$: all; $P$ & -- & 3 \\
\bottomrule\end{tabular}}\par\vspace{2mm}
At $\theta=0$ the twist is the identity and the Petersson algebra is the para-octonions: the plate reproduces Plate~B05, with the same constant ($c=0.95\,s^*$, $s^*=0.2901$) and, within the sampling error, the same measurements. The three quarters through $1$ are ribbed shells, wound in ridges around the axis of $1$ and tapering to an apex on the side of positive real part; each fills $33\%$ of the ball. The purely imaginary quarter $(l,il,jl)$ is a ball, exact under every tested rotation, because $z\bullet z=-\n(z)$ for purely imaginary $z$ (Proposition~\ref{atl:prop:imaginary}).

The shells are not bodies of revolution, since the constant is generic. The orthogonal maps that fix $1$ and $c$ commute with the map (Proposition~\ref{atl:prop:equivariance}), but in a section through $1$ and two hidden directions they only give the mirror in the plane of $1$ and of the projection of $\im c$, which is not a coordinate plane. In the quarters $(1,l,il)$ and $(1,l,jl)$ this projection is almost along $il$, respectively $jl$, because the $l$-component of $c$ is tiny, and the coordinate mirror $M_y$ is approximate ($J=0.994$ and $0.995$); in the quarter $(1,il,jl)$ it is oblique and no coordinate mirror comes close. The point reflection is exact in every quarter (Remark~\ref{atl:rem:minus}).

The orbits lie in a subalgebra of dimension $3$, spanned by $1$, $\im z_0$ and $\im c$ (Proposition~\ref{atl:prop:span}). This is the starting point of the Petersson path.
\clearpage

\subsubsection{E02. \texorpdfstring{$\theta=30^\circ$}{theta=30\textdegree{}}}\label{atl:pl:E02}\mbox{}\par\vspace{0.5mm}
\noindent\includegraphics[width=\textwidth]{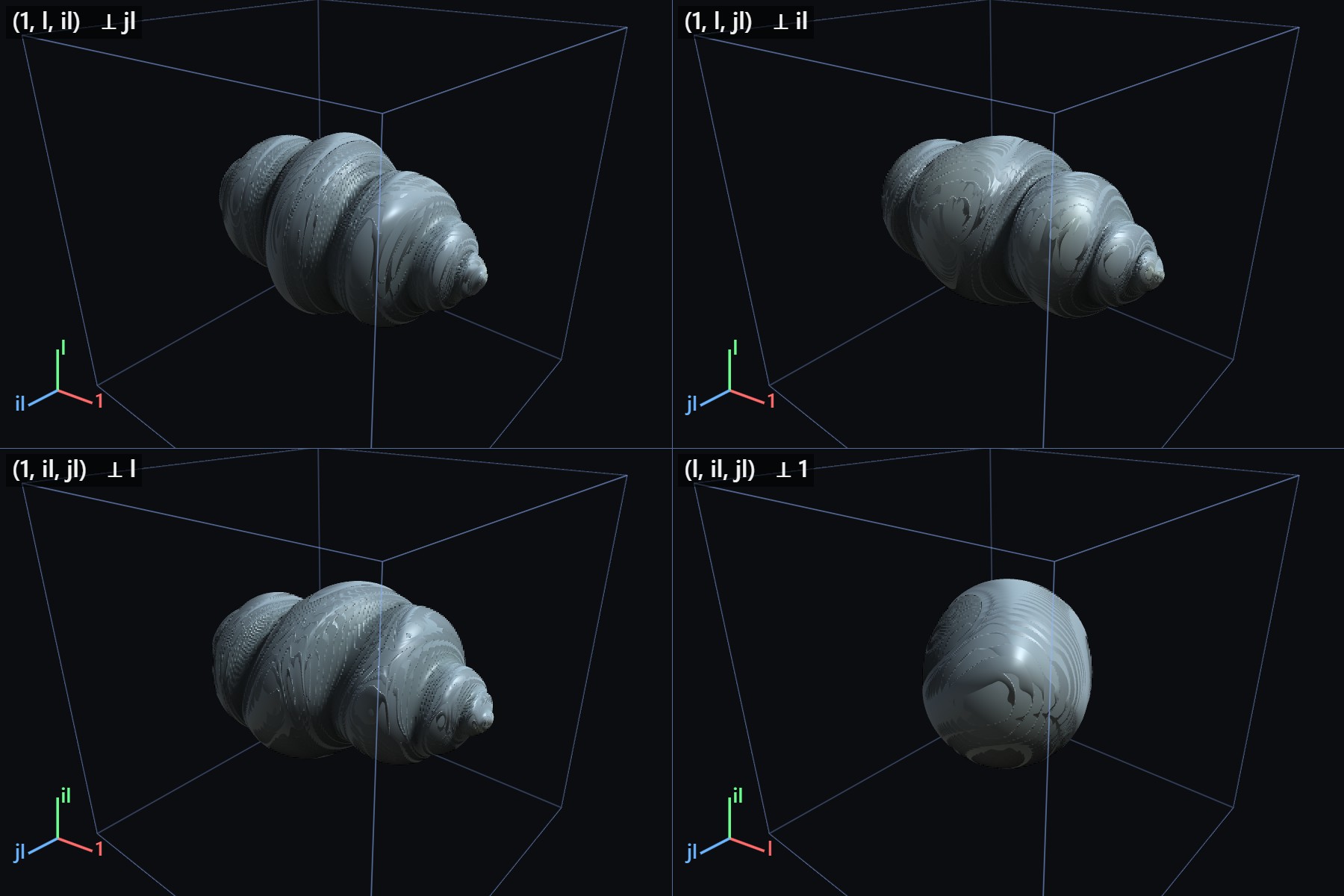}\par\vspace{1mm}
\noindent{\footnotesize \textbf{Map} $z\mapsto z\diamond z+c$ $\cdot$ \textbf{Algebra} custom: Petersson, $\gamma=(-1,-1,-1)$, $\theta=30^\circ$ $\cdot$ $c=(-0.353,\,\allowbreak -0.019,\,\allowbreak 0.044,\,\allowbreak 0.015,\,\allowbreak 0.003,\,\allowbreak 0.142,\,\allowbreak -0.1,\,\allowbreak 0.118)$ ($=0.95\,s^*$, $s^*=0.4358$) $\cdot$ \textbf{Frame} $X,Y,Z,N=1,l,il,jl$ $\cdot$ half-extent 1.4 $\cdot$ 20 iterations $\cdot$ \textbf{Coloring} 8D orbit trap (point), gradient Pastel}\par\vspace{1mm}
\noindent{\footnotesize\begin{tabular}{@{}lllll@{}}\toprule Slice (X, Y, Z) & Body [95\%] & Exact ($J\geq0.9995$) & Approximate ($0.97\leq J<0.9995$) & Dim\\\midrule
($1, l, il$) $\perp jl$ & 17.8\% {\scriptsize[17.6, 18.0]} & $P$ & -- & 8 \\
($1, l, jl$) $\perp il$ & 19.7\% {\scriptsize[19.5, 20.0]} & $P$ & -- & 8 \\
($1, il, jl$) $\perp l$ & 19.0\% {\scriptsize[18.7, 19.2]} & $P$ & -- & 8 \\
($l, il, jl$) $\perp 1$ & 14.5\% {\scriptsize[14.3, 14.7]} & $P$ & $M$: y; $R_y$: 180 (min.\ 0.971) & 8 \\
\bottomrule\end{tabular}}\par\vspace{2mm}
A twist of only $30^\circ$ ($s^*=0.436$) leaves the shapes recognizable but changes the algebra of the orbits at once. The three quarters through $1$ are still ribbed shells, like conch shells with a pointed apex on the side of positive real part, but they are smaller: $18$--$20\%$ of the ball, against $33\%$ at $\theta=0$. The imaginary quarter is no longer a ball but an egg, elongated along $il$; its best symmetry is the mirror $M_y$, only approximate ($J=0.971$).

The orbits now lie in a subalgebra of dimension $8$ in every quarter, against $3$ for the para-octonions: as soon as $\tau_\theta$ is not the identity, the product mixes the quaternionic half and the hidden half of the algebra. The Petersson algebra with $\theta=30^\circ$ is a composition algebra but not a symmetric one, and Proposition~\ref{atl:prop:imaginary} no longer applies. The only exact symmetry left in the table is the point reflection of Remark~\ref{atl:rem:minus}.

The jump of the dimension is immediate, while the change of the shapes, as the following plates show, is gradual.
\clearpage

\subsubsection{E03. \texorpdfstring{$\theta=60^\circ$}{theta=60\textdegree{}}}\label{atl:pl:E03}\mbox{}\par\vspace{0.5mm}
\noindent\includegraphics[width=\textwidth]{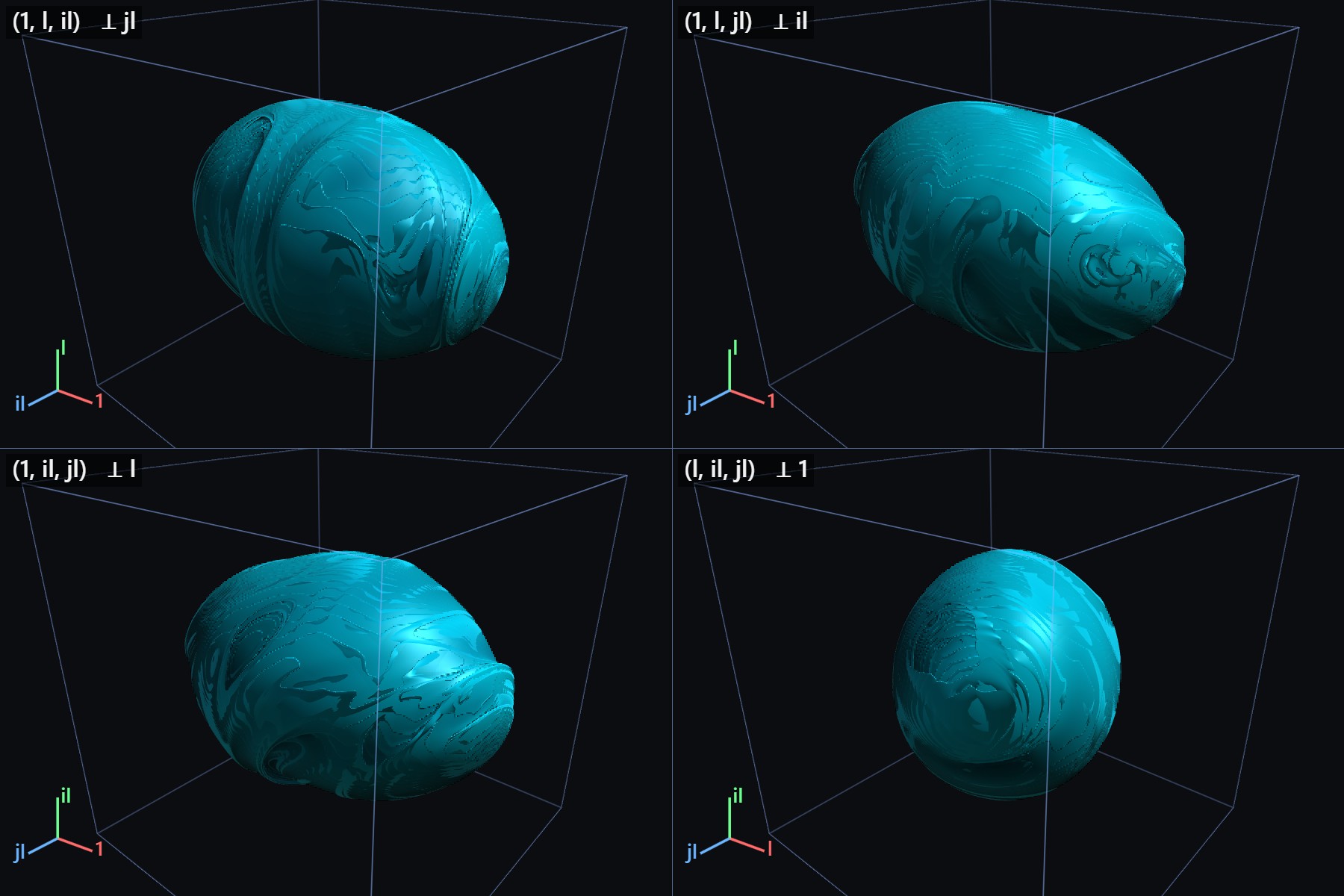}\par\vspace{1mm}
\noindent{\footnotesize \textbf{Map} $z\mapsto z\diamond z+c$ $\cdot$ \textbf{Algebra} custom: Petersson, $\gamma=(-1,-1,-1)$, $\theta=60^\circ$ $\cdot$ $c=(-0.503,\,\allowbreak -0.027,\,\allowbreak 0.063,\,\allowbreak 0.021,\,\allowbreak 0.005,\,\allowbreak 0.202,\,\allowbreak -0.143,\,\allowbreak 0.168)$ ($=0.95\,s^*$, $s^*=0.6206$) $\cdot$ \textbf{Frame} $X,Y,Z,N=1,l,il,jl$ $\cdot$ half-extent 1.4 $\cdot$ 20 iterations $\cdot$ \textbf{Coloring} 8D orbit trap (point), gradient Neon}\par\vspace{1mm}
\noindent{\footnotesize\begin{tabular}{@{}lllll@{}}\toprule Slice (X, Y, Z) & Body [95\%] & Exact ($J\geq0.9995$) & Approximate ($0.97\leq J<0.9995$) & Dim\\\midrule
($1, l, il$) $\perp jl$ & 42.9\% {\scriptsize[42.6, 43.2]} & $P$ & -- & 8 \\
($1, l, jl$) $\perp il$ & 47.7\% {\scriptsize[47.4, 48.0]} & $P$ & -- & 8 \\
($1, il, jl$) $\perp l$ & 45.4\% {\scriptsize[45.1, 45.7]} & $P$ & -- & 8 \\
($l, il, jl$) $\perp 1$ & 31.5\% {\scriptsize[31.2, 31.7]} & $P$ & -- & 8 \\
\bottomrule\end{tabular}}\par\vspace{2mm}
At $\theta=60^\circ$ ($s^*=0.621$) the bodies grow ($43$--$48\%$ of the ball in the quarters through $1$, $32\%$ in the imaginary quarter) and lose their ribs: they are smooth eggs, elongated along the axis of $1$, whose surface only carries the swirling patterns of the colouring. In appearance they are closer to the Okubonic bodies of Plate~E05 than to the para-octonionic shells of Plate~E01. The imaginary quarter is a rounded body, but not a ball.

No quarter has an exact or even an approximate mirror or rotation: the only exact symmetry is the point reflection of Remark~\ref{atl:rem:minus}, and the best scores ($J=0.958$ for $M_y$ in the quarter $(1,l,il)$, $0.957$ for $M_x$ in the quarter $(1,l,jl)$) stay below the threshold. The orbits lie in a subalgebra of dimension $8$.
\clearpage

\subsubsection{E04. \texorpdfstring{$\theta=90^\circ$}{theta=90\textdegree{}}: four identical balls}\label{atl:pl:E04}\mbox{}\par\vspace{0.5mm}
\noindent\includegraphics[width=\textwidth]{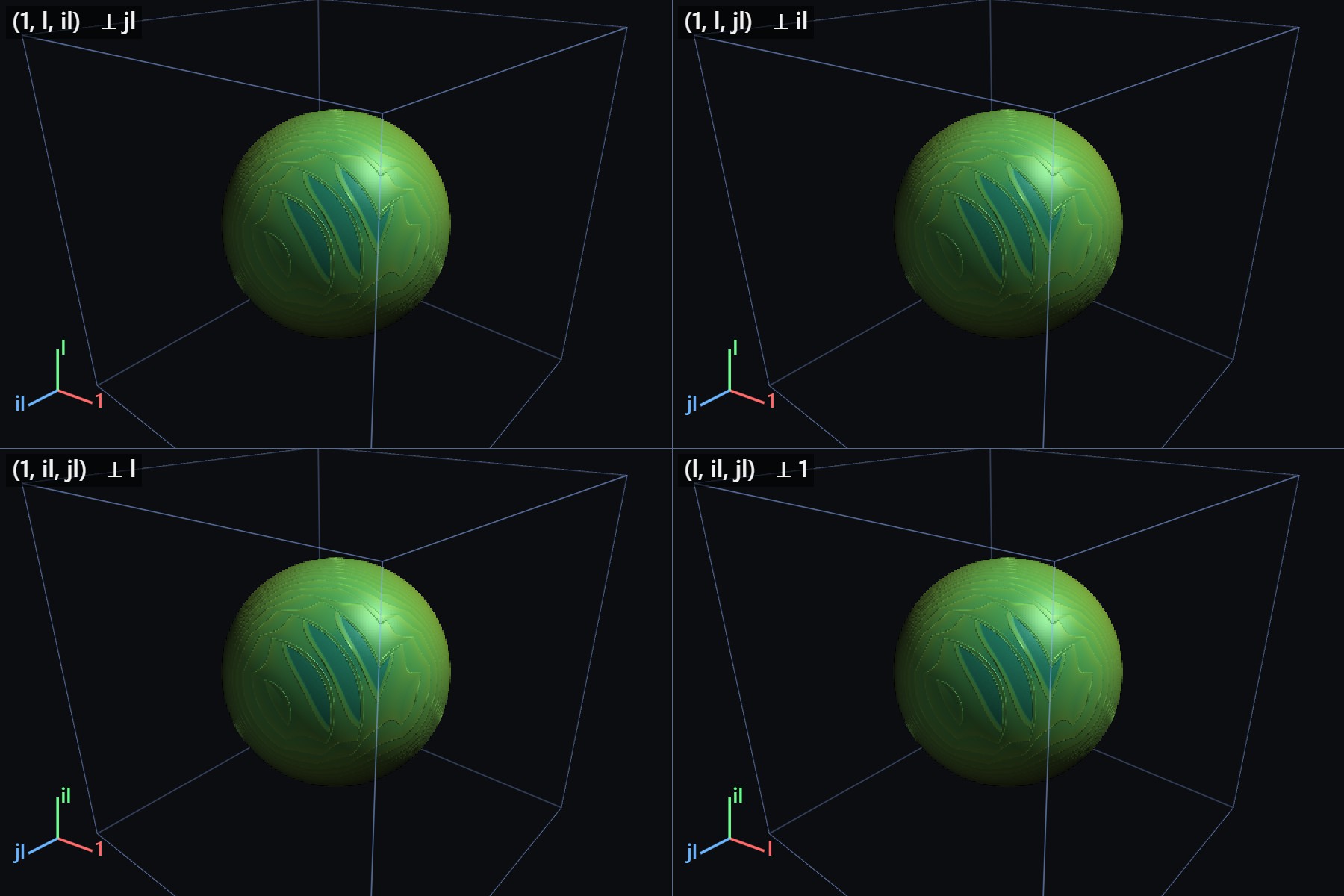}\par\vspace{1mm}
\noindent{\footnotesize \textbf{Map} $z\mapsto z\diamond z+c$ $\cdot$ \textbf{Algebra} custom: Petersson, $\gamma=(-1,-1,-1)$, $\theta=90^\circ$ $\cdot$ $c=(-0.756,\,\allowbreak -0.041,\,\allowbreak 0.094,\,\allowbreak 0.031,\,\allowbreak 0.007,\,\allowbreak 0.304,\,\allowbreak -0.215,\,\allowbreak 0.253)$ ($=0.95\,s^*$, $s^*=0.9326$) $\cdot$ \textbf{Frame} $X,Y,Z,N=1,l,il,jl$ $\cdot$ half-extent 2.4 $\cdot$ 20 iterations $\cdot$ \textbf{Coloring} 8D orbit trap (point), gradient Viridis-like}\par\vspace{1mm}
\noindent{\footnotesize\begin{tabular}{@{}lllll@{}}\toprule Slice (X, Y, Z) & Body [95\%] & Exact ($J\geq0.9995$) & Approximate ($0.97\leq J<0.9995$) & Dim\\\midrule
($1, l, il$) $\perp jl$ & 25.6\% {\scriptsize[25.3, 25.8]} & $M$: xyz; $R_x$: all; $R_y$: all; $R_z$: all; $P$ & -- & 6 \\
($1, l, jl$) $\perp il$ & 25.6\% {\scriptsize[25.3, 25.8]} & $M$: xyz; $R_x$: all; $R_y$: all; $R_z$: all; $P$ & -- & 6 \\
($1, il, jl$) $\perp l$ & 25.6\% {\scriptsize[25.3, 25.8]} & $M$: xyz; $R_x$: all; $R_y$: all; $R_z$: all; $P$ & -- & 6 \\
($l, il, jl$) $\perp 1$ & 25.6\% {\scriptsize[25.3, 25.8]} & $M$: xyz; $R_x$: all; $R_y$: all; $R_z$: all; $P$ & -- & 6 \\
\bottomrule\end{tabular}}\par\vspace{2mm}
The plate at $\theta=90^\circ$ is the most striking of the path: the four quarters show the same ball, and the four pictures coincide down to the faint pattern on its surface. The cube had to be enlarged to half-extent $2.4$ to contain it ($s^*=0.933$). In each quarter the ball fills $25.6\%$ of the ball of radius $0.95\cdot2.4=2.28$ in which the measurements are taken, so that its radius is about $1.45$, and each quarter is exactly invariant under all the $22$ tested transformations.

All the axes of the frame $F_B$ lie in $\RR\oplus\HH l$, and by Proposition~\ref{atl:prop:theta90} the Petersson square of $x=x_0+ql$ is the real number $x_0^2+|q|^2$: at $\theta=90^\circ$ the twist is $u_\theta=i$, and $u_\theta^2=-1$ turns the middle term of \eqref{atl:eq:okubocubic} into $+|q|^2$, while $\cos\theta=0$ removes the last one. The first iterate $z_1=|z_0|^2+c$ therefore depends only on the Euclidean norm of $z_0$, and the filled Julia set meets every section contained in $\RR\oplus\HH l$ in the same ball, whatever the constant. The orbits themselves are not simple: after the first step they leave $\RR\oplus\HH l$, and they lie in a subalgebra of dimension $6$.

It is worth noticing that these rotations are not automorphisms of the algebra: they come from an identity of the product on a subspace (Proposition~\ref{atl:prop:theta90}), not from its symmetry group. The Petersson algebra with $\theta=90^\circ$ is not symmetric, yet on $\RR\oplus\HH l$ its square is the Euclidean norm.
\clearpage

\subsubsection{E05. \texorpdfstring{$\theta=120^\circ$}{theta=120\textdegree{}}: the Okubo algebra}\label{atl:pl:E05}\mbox{}\par\vspace{0.5mm}
\noindent\includegraphics[width=\textwidth]{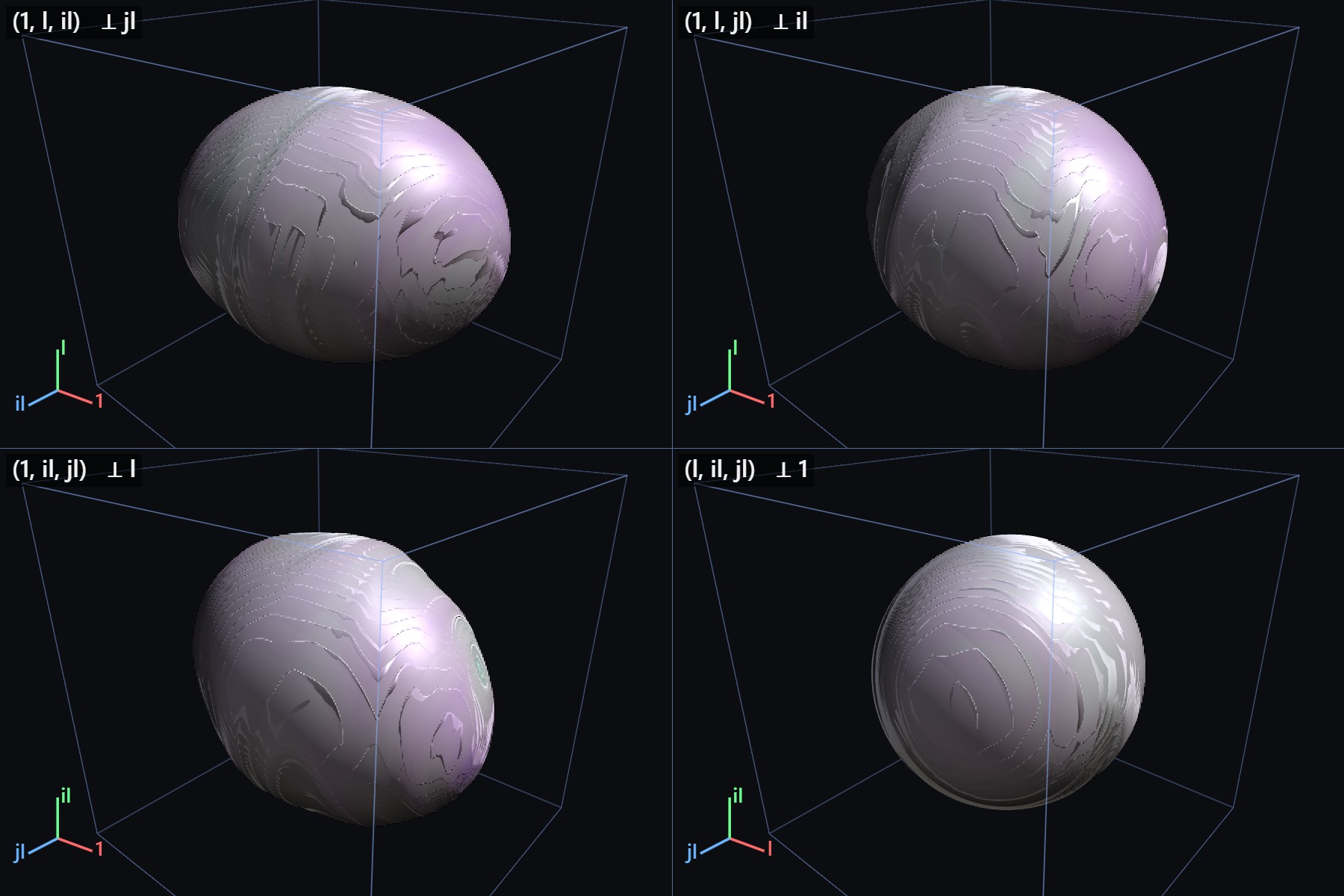}\par\vspace{1mm}
\noindent{\footnotesize \textbf{Map} $z\mapsto z\diamond z+c$ $\cdot$ \textbf{Algebra} custom: Petersson, $\gamma=(-1,-1,-1)$, $\theta=120^\circ$ $\cdot$ $c=(-0.397,\,\allowbreak -0.021,\,\allowbreak 0.049,\,\allowbreak 0.016,\,\allowbreak 0.004,\,\allowbreak 0.16,\,\allowbreak -0.113,\,\allowbreak 0.133)$ ($=0.95\,s^*$, $s^*=0.4900$) $\cdot$ \textbf{Frame} $X,Y,Z,N=1,l,il,jl$ $\cdot$ half-extent 1.4 $\cdot$ 20 iterations $\cdot$ \textbf{Coloring} 8D orbit trap (point), gradient Aurora}\par\vspace{1mm}
\noindent{\footnotesize\begin{tabular}{@{}lllll@{}}\toprule Slice (X, Y, Z) & Body [95\%] & Exact ($J\geq0.9995$) & Approximate ($0.97\leq J<0.9995$) & Dim\\\midrule
($1, l, il$) $\perp jl$ & 52.7\% {\scriptsize[52.4, 53.0]} & $P$ & $M$: y; $R_y$: 180 (min.\ 0.986) & 8 \\
($1, l, jl$) $\perp il$ & 51.9\% {\scriptsize[51.6, 52.2]} & $P$ & -- & 8 \\
($1, il, jl$) $\perp l$ & 51.4\% {\scriptsize[51.1, 51.7]} & $P$ & -- & 8 \\
($l, il, jl$) $\perp 1$ & 42.3\% {\scriptsize[42.0, 42.6]} & $P$ & $M$: x; $R_x$: 180,37; $R_y$: 90,72,60,37; $R_z$: 90,72,60,37 (min.\ 0.973) & 8 \\
\bottomrule\end{tabular}}\par\vspace{2mm}
At $\theta=120^\circ$ the path reaches the Okubo algebra, and the plate reproduces Plate~B07: same constant ($s^*=0.490$) and, within the sampling error, the same measurements. The quarters through $1$ are large smooth eggs ($51$--$53\%$ of the ball), similar to those at $60^\circ$; the imaginary quarter is a nearly round body ($42\%$). The orbits lie in a subalgebra of dimension $8$.

The point reflection is exact in every quarter (Remark~\ref{atl:rem:minus}); the mirror $M_y$ of the quarter $(1,l,il)$ is approximate ($J=0.986$), and the imaginary quarter is nearly, but not exactly, round: its mirror $M_x$, the half-turn about $l$ and several rotations score between $0.97$ and $0.99$, and none is exact. For a generic constant the automorphisms fixing $c$ form a torus conjugate to that of Proposition~\ref{atl:prop:torus}, whose root planes are in general position with respect to the coordinate axes: the absence of revolutions in the frame $F_B$ is a frame effect (Plate~R09).

Along the path, therefore, the Okubo algebra does not look special in this frame. What makes $\theta=120^\circ$ special, the symmetric identities \eqref{atl:eq:symmetric}, does not show in the Julia sets of the quadratic map with a generic constant; it becomes visible in the Mandelbrot sets with seed $0$ (Plates~F03 and F04) and in the sandwich and cubic maps (Parts~H and I).
\clearpage

\subsubsection{E06. \texorpdfstring{$\theta=150^\circ$}{theta=150\textdegree{}}}\label{atl:pl:E06}\mbox{}\par\vspace{0.5mm}
\noindent\includegraphics[width=\textwidth]{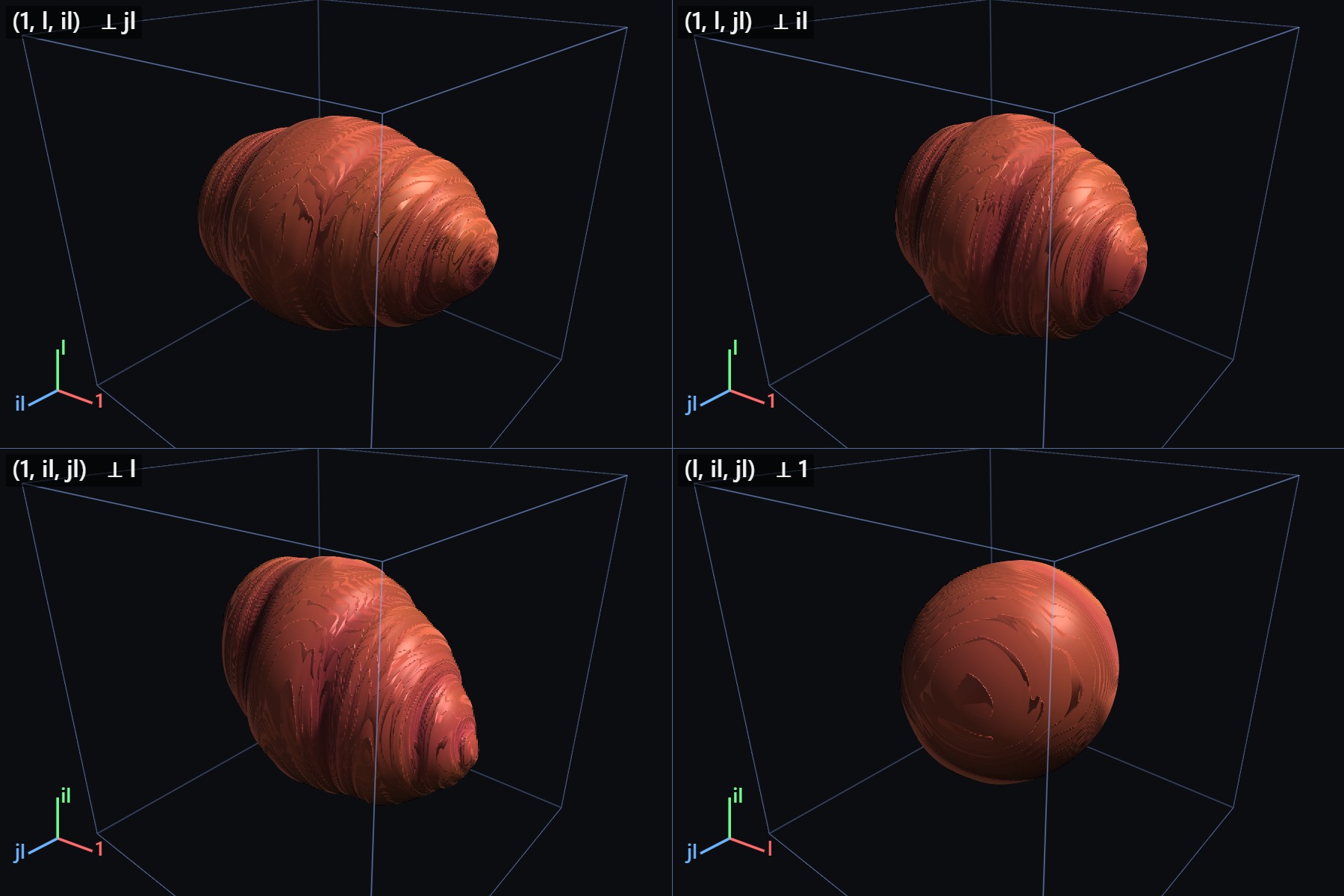}\par\vspace{1mm}
\noindent{\footnotesize \textbf{Map} $z\mapsto z\diamond z+c$ $\cdot$ \textbf{Algebra} custom: Petersson, $\gamma=(-1,-1,-1)$, $\theta=150^\circ$ $\cdot$ $c=(-0.415,\,\allowbreak -0.022,\,\allowbreak 0.052,\,\allowbreak 0.017,\,\allowbreak 0.004,\,\allowbreak 0.167,\,\allowbreak -0.118,\,\allowbreak 0.139)$ ($=0.95\,s^*$, $s^*=0.5127$) $\cdot$ \textbf{Frame} $X,Y,Z,N=1,l,il,jl$ $\cdot$ half-extent 1.4 $\cdot$ 20 iterations $\cdot$ \textbf{Coloring} 8D orbit trap (point), gradient Sunset}\par\vspace{1mm}
\noindent{\footnotesize\begin{tabular}{@{}lllll@{}}\toprule Slice (X, Y, Z) & Body [95\%] & Exact ($J\geq0.9995$) & Approximate ($0.97\leq J<0.9995$) & Dim\\\midrule
($1, l, il$) $\perp jl$ & 27.9\% {\scriptsize[27.7, 28.2]} & $P$ & -- & 8 \\
($1, l, jl$) $\perp il$ & 25.3\% {\scriptsize[25.0, 25.5]} & $P$ & -- & 8 \\
($1, il, jl$) $\perp l$ & 26.1\% {\scriptsize[25.8, 26.4]} & $P$ & -- & 8 \\
($l, il, jl$) $\perp 1$ & 20.5\% {\scriptsize[20.3, 20.8]} & $P$ & -- & 8 \\
\bottomrule\end{tabular}}\par\vspace{2mm}
At $\theta=150^\circ$ the twist $\tau_\theta$ has order twelve, and $s^*=0.513$. The quarters through $1$ are again ribbed bodies, elongated along the axis of $1$ with a rounded apex on the side of positive real part, fuller than the shells of Plate~E02 ($25$--$28\%$ of the ball). The imaginary quarter is a rounded body ($20.5\%$), close to a ball in appearance but without any measured symmetry: no rotation scores more than $J=0.94$.

The only exact symmetry of the table is the point reflection of Remark~\ref{atl:rem:minus}; no mirror or rotation reaches the threshold of the approximate symmetries (the best score is $J=0.947$, for $M_y$ in the quarter $(1,l,il)$). The orbits lie in a subalgebra of dimension $8$. Seen along the path, the ribs of the para-octonionic shells return as $\theta$ approaches $180^\circ$, while the smooth eggs belong to the angles $60^\circ$ and $120^\circ$.
\clearpage

\subsubsection{E07. \texorpdfstring{$\theta=180^\circ$}{theta=180\textdegree{}}: an involutive twist}\label{atl:pl:E07}\mbox{}\par\vspace{0.5mm}
\noindent\includegraphics[width=\textwidth]{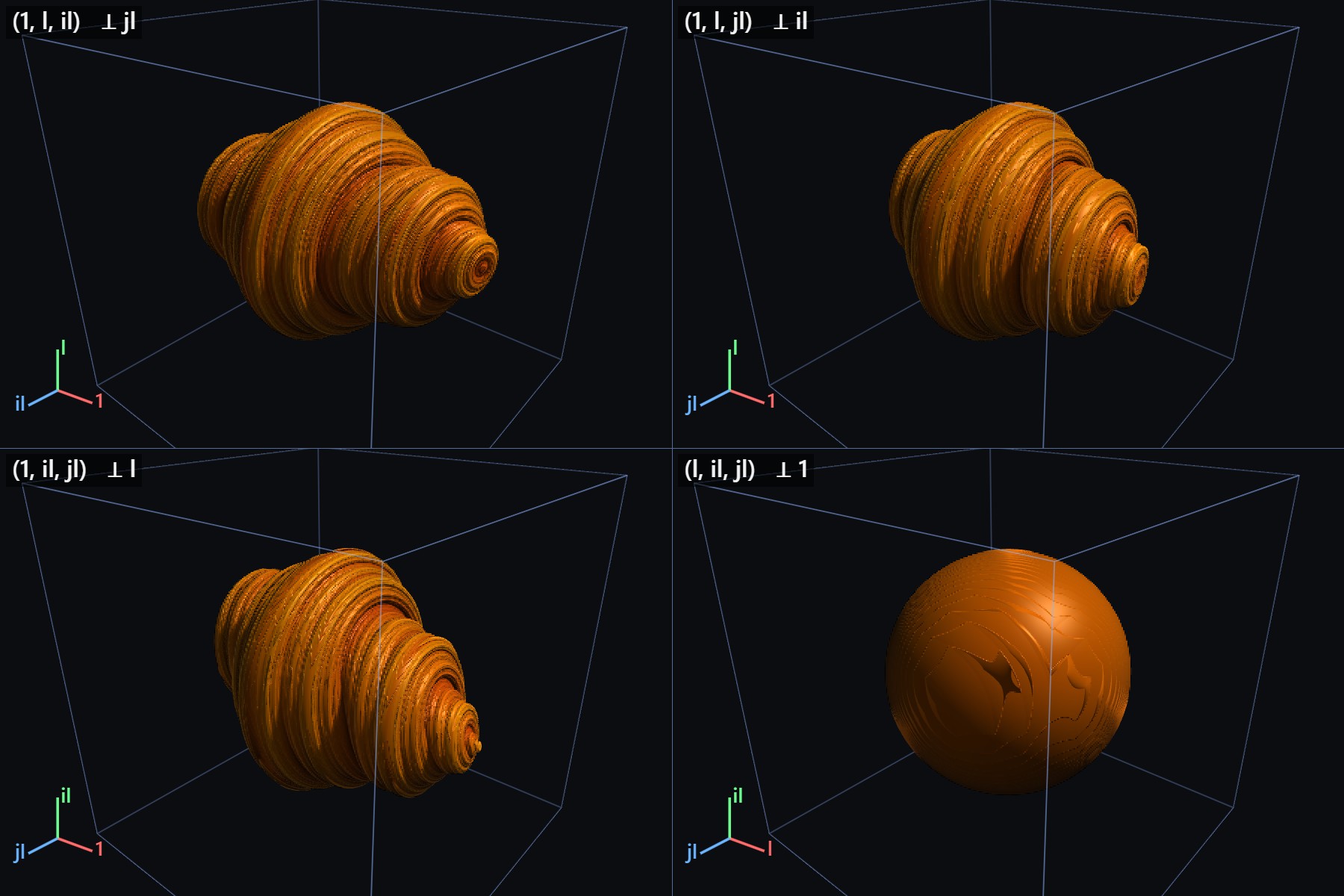}\par\vspace{1mm}
\noindent{\footnotesize \textbf{Map} $z\mapsto z\diamond z+c$ $\cdot$ \textbf{Algebra} custom: Petersson, $\gamma=(-1,-1,-1)$, $\theta=180^\circ$ $\cdot$ $c=(-0.419,\,\allowbreak -0.023,\,\allowbreak 0.052,\,\allowbreak 0.017,\,\allowbreak 0.004,\,\allowbreak 0.169,\,\allowbreak -0.119,\,\allowbreak 0.14)$ ($=0.95\,s^*$, $s^*=0.5172$) $\cdot$ \textbf{Frame} $X,Y,Z,N=1,l,il,jl$ $\cdot$ half-extent 1.4 $\cdot$ 20 iterations $\cdot$ \textbf{Coloring} 8D orbit trap (point), gradient Fire}\par\vspace{1mm}
\noindent{\footnotesize\begin{tabular}{@{}lllll@{}}\toprule Slice (X, Y, Z) & Body [95\%] & Exact ($J\geq0.9995$) & Approximate ($0.97\leq J<0.9995$) & Dim\\\midrule
($1, l, il$) $\perp jl$ & 26.9\% {\scriptsize[26.7, 27.2]} & $P$ & $M$: y; $R_y$: 180 (min.\ 0.992) & 4 \\
($1, l, jl$) $\perp il$ & 26.9\% {\scriptsize[26.7, 27.2]} & $P$ & $M$: y; $R_y$: 180 (min.\ 0.993) & 4 \\
($1, il, jl$) $\perp l$ & 27.2\% {\scriptsize[26.9, 27.4]} & $P$ & -- & 4 \\
($l, il, jl$) $\perp 1$ & 31.1\% {\scriptsize[30.8, 31.4]} & $M$: xyz; $R_x$: all; $R_y$: all; $R_z$: all; $P$ & -- & 4 \\
\bottomrule\end{tabular}}\par\vspace{2mm}
At $\theta=180^\circ$ the twist $\tau(p+ql)=p-ql$ is an involutive automorphism of the octonions, and the product $\tau(\bar x)\,\tau(\bar y)=\tau(\bar x\bar y)$ is the para-octonionic product followed by an automorphism: $z*z=\tau(z\bullet z)$. The plate is indeed close to Plate~E01: strongly ribbed conical shells in the quarters through $1$ ($27\%$ of the ball, against $33\%$), and an exact ball in the imaginary quarter ($31\%$), invariant under every tested rotation. The ball has the same origin as in the para-octonions: for purely imaginary $z$ the square $z\bullet z=-\n(z)$ is real, and $\tau$ fixes it.

The table repeats the pattern of Plate~E01: the point reflection is exact everywhere, and the mirror $M_y$ of the quarters $(1,l,il)$ and $(1,l,jl)$ is approximate ($J=0.992$ and $0.993$). The orbits lie in a subalgebra of dimension $4$, one more than for the para-octonions. The reason can be read off the product: $z\bullet z=2\re(z)\,\bar z-\n(z)$, hence $z*z=2\re(z)\,\tau(\bar z)-\n(z)$; for $z_0\in\RR\oplus\HH l$ the twist only changes the sign of $\im z_0$, but it separates the quaternionic part and the $\HH l$ part of $\im c$, which become two independent directions.

The Petersson path thus closes on itself: at $0^\circ$ and $180^\circ$ the dynamics is that of a para-Hurwitz product, with orbits in a subalgebra of dimension $3$ or $4$; in between the orbits fill all eight dimensions; and at $90^\circ$ the dynamics becomes radial on the five-dimensional subspace $\RR\oplus\HH l$.
\clearpage

\subsection{Part F. Mandelbrot sets}\label{atl:part:F}
The Mandelbrot map $z\mapsto z\diamond z+c$, with $c$ equal to the sample point and the seed $z_0=a$, is the parameter space of the quadratic family. With the seed $0$ its orbits lie in a plane (Proposition~\ref{atl:prop:span}), and the pictures are pulled back from a planar set, which is the tricorn in the compact symmetric composition algebras (Theorem~B); with a seed $a\ne0$ the dynamics becomes genuinely higher-dimensional. For the split algebras the Euclidean bailout of the renderer is not an escape criterion, and the rendered bodies are approximations. The center of the box is moved to $c_0=-0.35$ on the real axis.

\subsubsection{F01. Octonions, seed \texorpdfstring{$0$}{0}: the Mandelbrot set of revolution}\label{atl:pl:F01}\mbox{}\par\vspace{0.5mm}
\noindent\includegraphics[width=\textwidth]{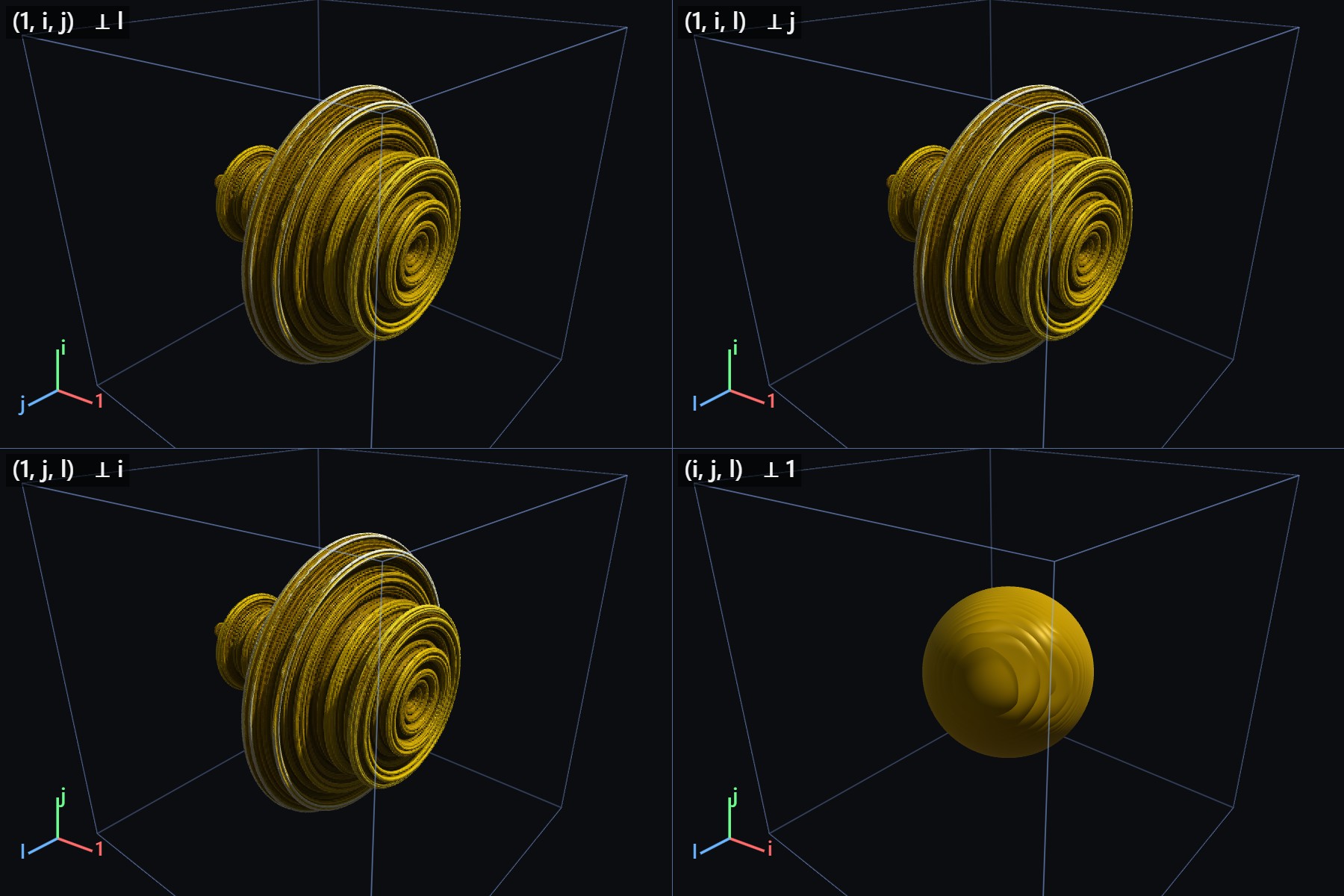}\par\vspace{1mm}
\noindent{\footnotesize \textbf{Map} $z\mapsto z\diamond z+c$, $c$ the sample point, $z_0=a$ $\cdot$ \textbf{Algebra} Octonions (Cayley--Dickson, $\gamma=(-1,-1,-1)$) $\cdot$ \textbf{Frame} $X,Y,Z,N=1,i,j,l$ $\cdot$ center $(-0.35,\,\allowbreak 0,\,\allowbreak 0,\,\allowbreak 0,\,\allowbreak 0,\,\allowbreak 0,\,\allowbreak 0,\,\allowbreak 0)$ $\cdot$ half-extent 1.5 $\cdot$ 16 iterations $\cdot$ \textbf{Coloring} 8D orbit trap (point), gradient Gold}\par\vspace{1mm}
\noindent{\footnotesize\begin{tabular}{@{}lllll@{}}\toprule Slice (X, Y, Z) & Body [95\%] & Exact ($J\geq0.9995$) & Approximate ($0.97\leq J<0.9995$) & Dim\\\midrule
($1, i, j$) $\perp l$ & 15.1\% {\scriptsize[14.9, 15.4]} & $M$: yz; $R_x$: all & -- & 2 \\
($1, i, l$) $\perp j$ & 15.1\% {\scriptsize[14.9, 15.4]} & $M$: yz; $R_x$: all & -- & 2 \\
($1, j, l$) $\perp i$ & 15.1\% {\scriptsize[14.9, 15.4]} & $M$: yz; $R_x$: all & -- & 2 \\
($i, j, l$) $\perp 1$ & 10.5\% {\scriptsize[10.3, 10.7]} & $M$: xyz; $R_x$: all; $R_y$: all; $R_z$: all; $P$ & -- & 2 \\
\bottomrule\end{tabular}}\par\vspace{2mm}
The octonionic Mandelbrot set with seed $0$ is the classical Mandelbrot set turned about the real axis. By Proposition~\ref{atl:prop:span} the orbit of $0$ lies in the plane of $1$ and $c$, a copy of $\CC$, and the escape depends only on $\re c$ and $|\im c|$. Hence the three quarters through $1$ are the same body ($15.1\%$ of the ball in each), an exact body of revolution about the axis of $1$, and the imaginary quarter at $\re c=-0.35$ is a ball, exact under every tested rotation.

The main cardioid appears as a large ribbed disc, whose front face is covered with the concentric rings swept by the bulbs attached to it; the period-two bulb is the smaller disc behind it, on the negative real side, followed by a small knob on the axis. The orbits lie in a subalgebra of dimension $2$, the plane of $1$ and $c$. The plate is the parameter-space counterpart of Plate~A01.
\clearpage

\subsubsection{F02. Para-octonions, seed \texorpdfstring{$0$}{0}: the tricorn of revolution}\label{atl:pl:F02}\mbox{}\par\vspace{0.5mm}
\noindent\includegraphics[width=\textwidth]{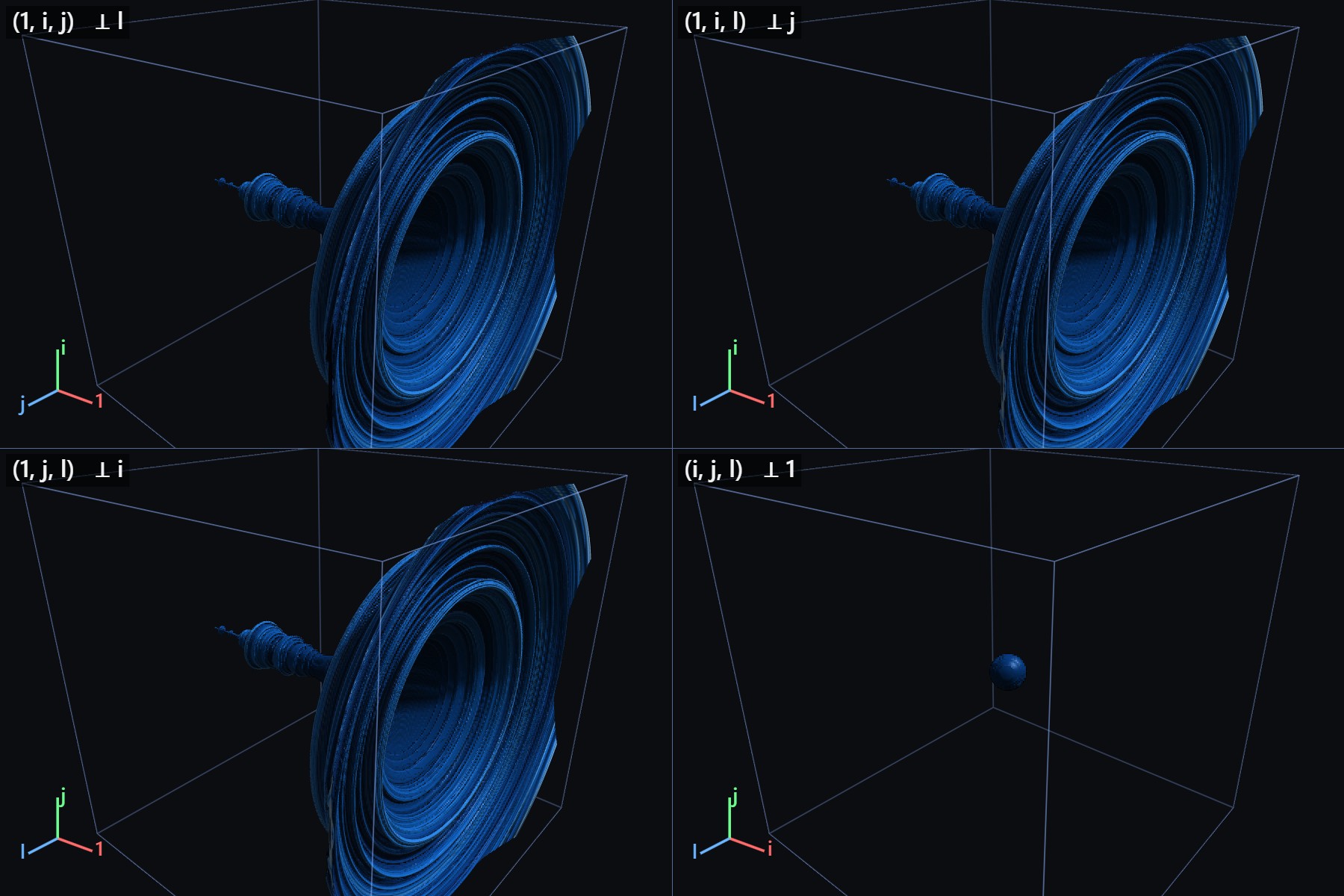}\par\vspace{1mm}
\noindent{\footnotesize \textbf{Map} $z\mapsto z\diamond z+c$, $c$ the sample point, $z_0=a$ $\cdot$ \textbf{Algebra} Para-octonions (para-Hurwitz $\bar x\bar y$, $\gamma=(-1,-1,-1)$) $\cdot$ \textbf{Frame} $X,Y,Z,N=1,i,j,l$ $\cdot$ center $(-0.35,\,\allowbreak 0,\,\allowbreak 0,\,\allowbreak 0,\,\allowbreak 0,\,\allowbreak 0,\,\allowbreak 0,\,\allowbreak 0)$ $\cdot$ half-extent 1.5 $\cdot$ 16 iterations $\cdot$ \textbf{Coloring} 8D orbit trap (point), gradient Deep Blue}\par\vspace{1mm}
\noindent{\footnotesize\begin{tabular}{@{}lllll@{}}\toprule Slice (X, Y, Z) & Body [95\%] & Exact ($J\geq0.9995$) & Approximate ($0.97\leq J<0.9995$) & Dim\\\midrule
($1, i, j$) $\perp l$ & 9.3\% {\scriptsize[9.1, 9.5]} & $M$: yz; $R_x$: all & -- & 2 \\
($1, i, l$) $\perp j$ & 9.3\% {\scriptsize[9.1, 9.5]} & $M$: yz; $R_x$: all & -- & 2 \\
($1, j, l$) $\perp i$ & 9.3\% {\scriptsize[9.1, 9.5]} & $M$: yz; $R_x$: all & -- & 2 \\
($i, j, l$) $\perp 1$ & 0.1\% {\scriptsize[0.1, 0.1]} & $M$: xyz; $R_x$: all; $R_y$: all; $R_z$: all; $P$ & -- & 2 \\
\bottomrule\end{tabular}}\par\vspace{2mm}
The para-octonionic Mandelbrot set with seed $0$ is the tricorn \cite{CHRS} turned about the real axis. The three quarters through $1$ are the same body ($9.3\%$ of the ball), an exact body of revolution about the axis of $1$: a wide funnel, open towards positive real part and cut by the faces of the cube, and a thin beaded stem along the negative real axis. The stem is swept by the arm of the tricorn that lies on the negative real axis, the funnel by the two arms that point towards $e^{\pm i\pi/3}$. The imaginary quarter at $\re c=-0.35$ contains only a tiny ball ($0.1\%$ of the section, exact under every tested rotation), because the vertical line $\re\zeta=-0.35$ meets the tricorn only in a short segment.

Theorem~B explains the plate at once. In a symmetric composition algebra the orbit of $0$ stays in the plane of $c$ and $c\bullet c$ (Proposition~\ref{atl:prop:span}), and when $4\n(c)^3>\n(c,c\bullet c)^2$ this plane is isometric to $(\CC,\bar x\,\bar y)$ through $c\mapsto\zeta$, with $|\zeta|^2=\n(c)$ and $2\re\zeta^3=\n(c,c\bullet c)$: the point $c$ belongs to the Mandelbrot set if and only if $\zeta$ belongs to the tricorn. In the para-octonions $\n(c,c\bullet c)=2c_0(c_0^2-3|\im c|^2)$ (Proposition~\ref{atl:prop:okubocubic}), so that one may take $\zeta=c_0+i|\im c|$: the Mandelbrot set is the tricorn moved by $G_2$, a tricorn of revolution. The orbits lie in a subalgebra of dimension $2$.

This plate and the previous one are the references for the Okubo algebra, whose Mandelbrot set is pulled back from the same tricorn, but moved by a different group.
\clearpage

\subsubsection{F03. Okubo, seed \texorpdfstring{$0$}{0}: the cubic form at work}\label{atl:pl:F03}\mbox{}\par\vspace{0.5mm}
\noindent\includegraphics[width=\textwidth]{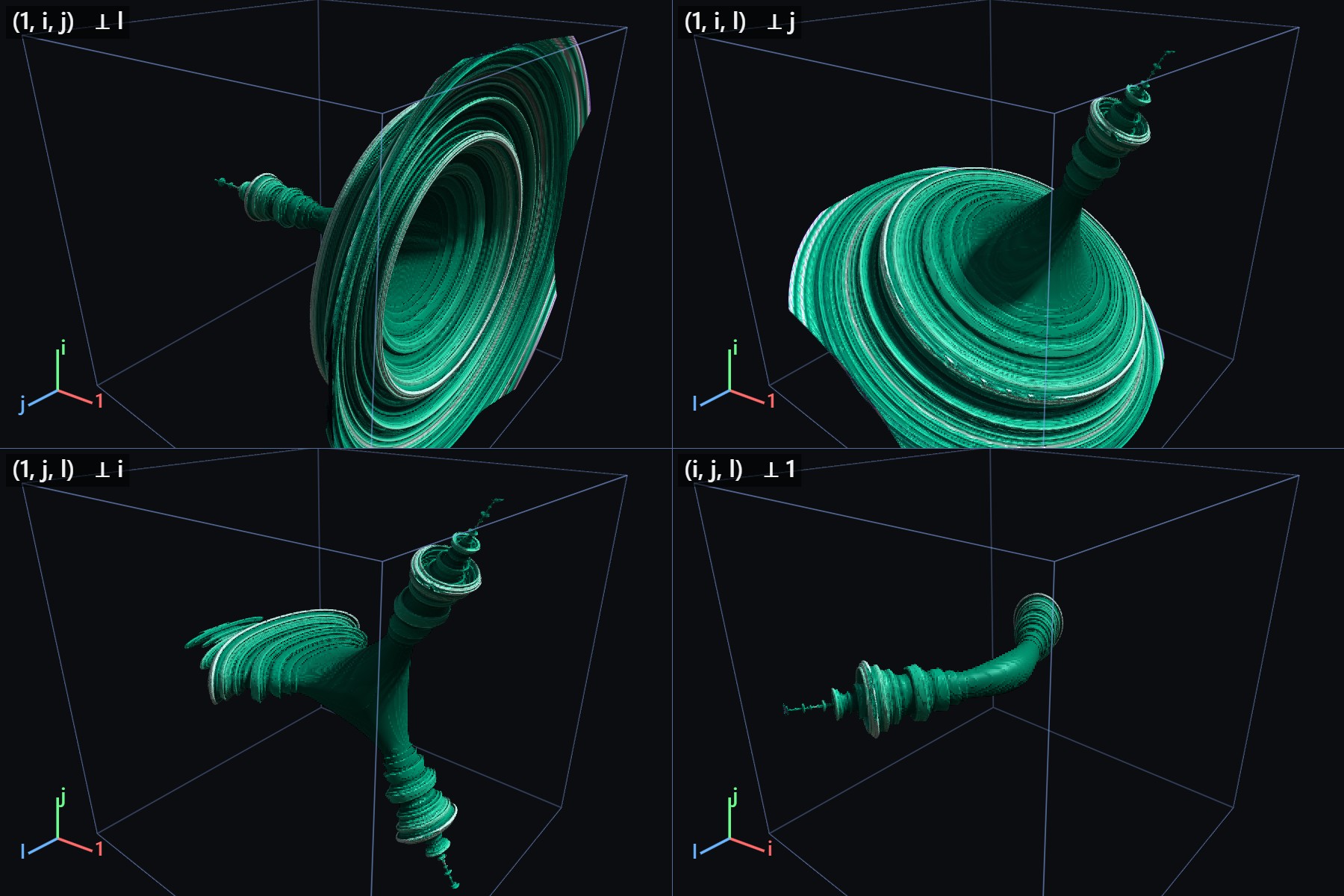}\par\vspace{1mm}
\noindent{\footnotesize \textbf{Map} $z\mapsto z\diamond z+c$, $c$ the sample point, $z_0=a$ $\cdot$ \textbf{Algebra} Okubo algebra (Petersson, $\gamma=(-1,-1,-1)$, $\theta=120^\circ$) $\cdot$ \textbf{Frame} $X,Y,Z,N=1,i,j,l$ $\cdot$ center $(-0.35,\,\allowbreak 0,\,\allowbreak 0,\,\allowbreak 0,\,\allowbreak 0,\,\allowbreak 0,\,\allowbreak 0,\,\allowbreak 0)$ $\cdot$ half-extent 1.5 $\cdot$ 16 iterations $\cdot$ \textbf{Coloring} 8D orbit trap (point), gradient Aurora}\par\vspace{1mm}
\noindent{\footnotesize\begin{tabular}{@{}lllll@{}}\toprule Slice (X, Y, Z) & Body [95\%] & Exact ($J\geq0.9995$) & Approximate ($0.97\leq J<0.9995$) & Dim\\\midrule
($1, i, j$) $\perp l$ & 9.3\% {\scriptsize[9.1, 9.5]} & $M$: yz; $R_x$: all & -- & 2 \\
($1, i, l$) $\perp j$ & 10.2\% {\scriptsize[10.1, 10.4]} & $M$: z & -- & 2 \\
($1, j, l$) $\perp i$ & 4.4\% {\scriptsize[4.3, 4.6]} & $M$: yz; $R_x$: 180 & -- & 2 \\
($i, j, l$) $\perp 1$ & 1.2\% {\scriptsize[1.2, 1.3]} & $M$: yz; $R_x$: 180 & -- & 2 \\
\bottomrule\end{tabular}}\par\vspace{2mm}
At first sight the Okubonic Mandelbrot set with seed $0$ is very different from the para-octonionic one. The upper-left quarter $(1,i,j)$ lies in $\HH$, where the Okubo product is $\bar x\,\bar y$ (Proposition~\ref{atl:prop:torus}), and it is the funnel of Plate~F02, with the same body share ($9.3\%$) and the same exact revolution about $1$. The upper-right quarter $(1,i,l)$ shows a funnel and a beaded stem again, but tilted: the stem points upwards and to the right, along a line of the plane of $1$ and $i$, and the funnel opens downwards to the left. The lower-left quarter $(1,j,l)$ is a three-armed body ($4.4\%$): a flared, ribbed fan on the left and two ribbed stems, upwards and downwards, joined by a thin web. The imaginary quarter at $\re c=-0.35$ is a bent, ribbed tube ($1.2\%$), ending in a small disc on one side and in a thin beaded spike on the other.

Quite surprisingly, all four quarters are pulled back from the same planar figure. By Theorem~B the Okubonic Mandelbrot set is the tricorn moved by $PSU(3)$: $c$ belongs to it if and only if the point $\zeta$ with $|\zeta|^2=\n(c)$ and $2\re\zeta^3=\n(c,c*c)$ belongs to the tricorn. What changes with respect to the para-octonions is the group: the orbits of $PSU(3)$, cut out by $\n(c)$ and $\n(c,c*c)$, are smaller than those of $G_2$, and the Okubo cubic form is not a function of $c_0$ and $|\im c|$ alone (Proposition~\ref{atl:prop:okubocubic}).

The revolutions are not lost but displaced. The section $\langle1,i,l\rangle$ is a body of revolution about the Weyl line of $e^{i\pi/3}$, a line of the plane of $1$ and $i$ at $60^\circ$ from the real axis: its stem is the arm of the tricorn that lies on this line, and its funnel is swept by the other two arms. This axis is neither a coordinate axis nor a line through the centre of the cube (it passes through the origin, $0.35$ to the right of the centre), so the table cannot detect the revolution, and it lists only the exact mirror $M_z$. In the section $\langle1,j,l\rangle$ the normalized cubic invariant $\n(c,c*c)/2\n(c)^{3/2}$ has three maxima and three minima on the unit sphere, and the three arms are the three arms of the tricorn, redistributed by $PSU(3)$; there, as in the imaginary quarter, the mirrors $M_y$, $M_z$ and the half-turn about the first axis are exact. The orbits lie in a subalgebra of dimension $2$, the plane of $c$ and $c*c$, as Proposition~\ref{atl:prop:span} predicts.
\clearpage

\subsubsection{F04. Okubo, seed \texorpdfstring{$0$}{0}, in the frame \texorpdfstring{$F_B$}{FB}}\label{atl:pl:F04}\mbox{}\par\vspace{0.5mm}
\noindent\includegraphics[width=\textwidth]{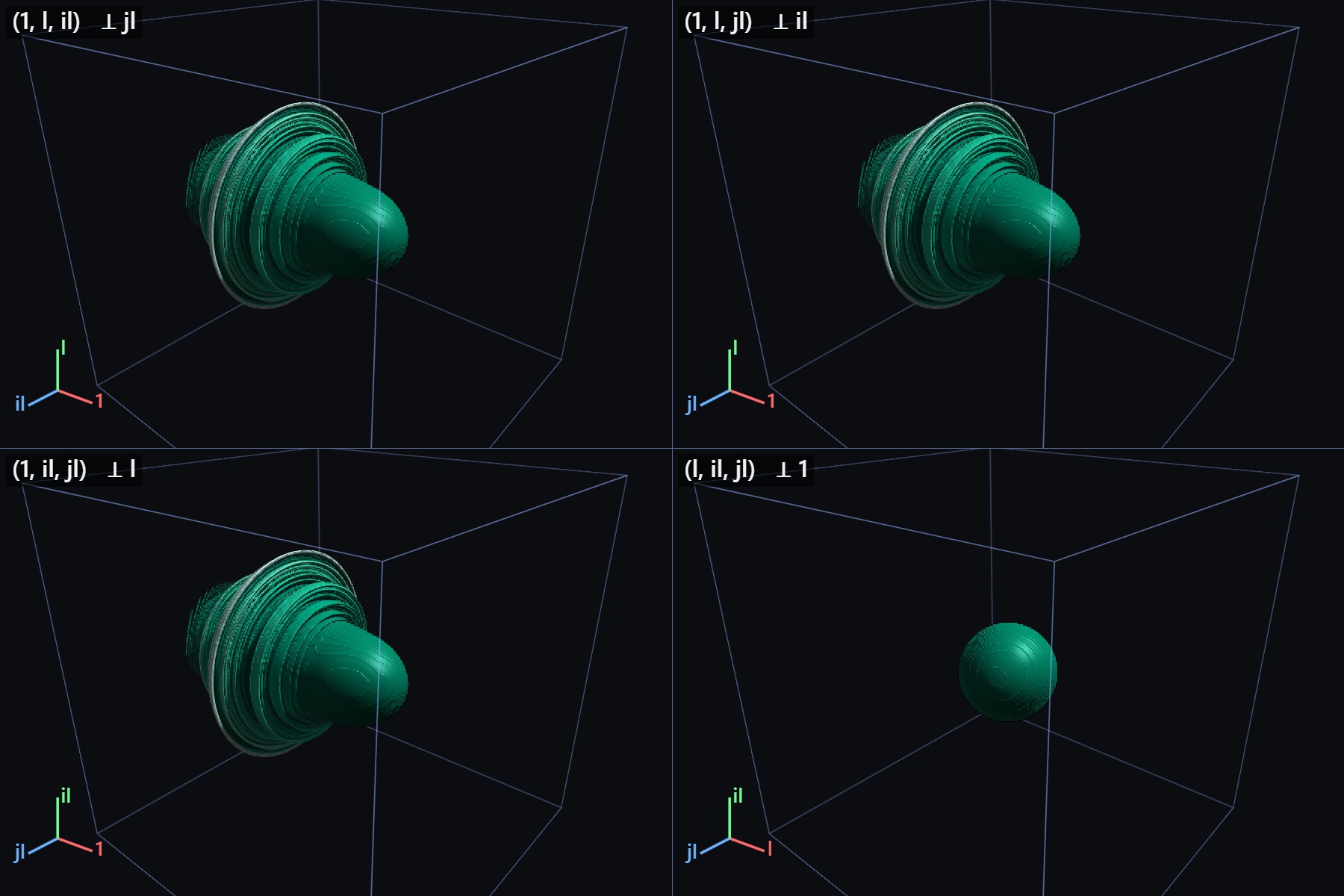}\par\vspace{1mm}
\noindent{\footnotesize \textbf{Map} $z\mapsto z\diamond z+c$, $c$ the sample point, $z_0=a$ $\cdot$ \textbf{Algebra} Okubo algebra (Petersson, $\gamma=(-1,-1,-1)$, $\theta=120^\circ$) $\cdot$ \textbf{Frame} $X,Y,Z,N=1,l,il,jl$ $\cdot$ center $(-0.35,\,\allowbreak 0,\,\allowbreak 0,\,\allowbreak 0,\,\allowbreak 0,\,\allowbreak 0,\,\allowbreak 0,\,\allowbreak 0)$ $\cdot$ half-extent 1.5 $\cdot$ 16 iterations $\cdot$ \textbf{Coloring} 8D orbit trap (point), gradient Aurora}\par\vspace{1mm}
\noindent{\footnotesize\begin{tabular}{@{}lllll@{}}\toprule Slice (X, Y, Z) & Body [95\%] & Exact ($J\geq0.9995$) & Approximate ($0.97\leq J<0.9995$) & Dim\\\midrule
($1, l, il$) $\perp jl$ & 10.8\% {\scriptsize[10.6, 11.0]} & $M$: yz; $R_x$: all & -- & 2 \\
($1, l, jl$) $\perp il$ & 10.8\% {\scriptsize[10.6, 11.0]} & $M$: yz; $R_x$: all & -- & 2 \\
($1, il, jl$) $\perp l$ & 10.8\% {\scriptsize[10.6, 11.0]} & $M$: yz; $R_x$: all & -- & 2 \\
($l, il, jl$) $\perp 1$ & 2.1\% {\scriptsize[2.0, 2.2]} & $M$: xyz; $R_x$: all; $R_y$: all; $R_z$: all; $P$ & -- & 2 \\
\bottomrule\end{tabular}}\par\vspace{2mm}
In the frame $F_B$ the same Okubonic Mandelbrot set shows three identical bodies of revolution about the axis of $1$ ($10.8\%$ of the ball each, with the rotations by all tested angles exact): a smooth rounded cap on the side of positive real part, followed on the negative side by a flared collar of ribbed rings. The imaginary quarter is a small ball ($2.1\%$), exact under every tested rotation. The orbits lie in a subalgebra of dimension $2$.

The explanation is Proposition~\ref{atl:prop:okubocubic}: on the subspace $\RR\oplus\HH l$, which contains all the quarters of $F_B$, the Okubo cubic form is $\n(c,c*c)=2c_0\big(c_0^2-(\cos2\theta+2\cos\theta)|q|^2\big)=2c_0\big(c_0^2+\tfrac32|q|^2\big)$, and, like $\n(c)=c_0^2+|q|^2$, it depends only on $c_0$ and $|q|$. By Theorem~B these two invariants decide whether $c$ belongs to the set, so every section through $1$ is a body of revolution about $1$, and the imaginary quarter, where $c_0=-0.35$ is fixed, is a ball. It is worth noticing that the quarters $(1,l,jl)$ and $(1,il,jl)$ contain no root plane, yet they are bodies of revolution: for a Mandelbrot set the variable is the constant itself, and the symmetry is not limited to the torus that fixes a given $c$.

The profile of these bodies is not the tricorn itself: the map $(c_0,|q|)\mapsto\zeta$ differs from the para-octonionic one $\zeta=c_0+i|\im c|$, and this is why the bodies of this plate have a cap where those of Plate~F02 have a stem. Plates~F03 and F04 together show that the Okubonic Mandelbrot set keeps its revolutions, about the axis of $1$ in the sections of $\RR\oplus\HH l$ and in $\langle1,i,j\rangle$, about a tilted Weyl line in $\langle1,i,l\rangle$, while in $\langle1,j,l\rangle$ they give way to the three arms of the tricorn.
\clearpage

\subsubsection{F05. Split-octonions, seed \texorpdfstring{$0$}{0}}\label{atl:pl:F05}\mbox{}\par\vspace{0.5mm}
\noindent\includegraphics[width=\textwidth]{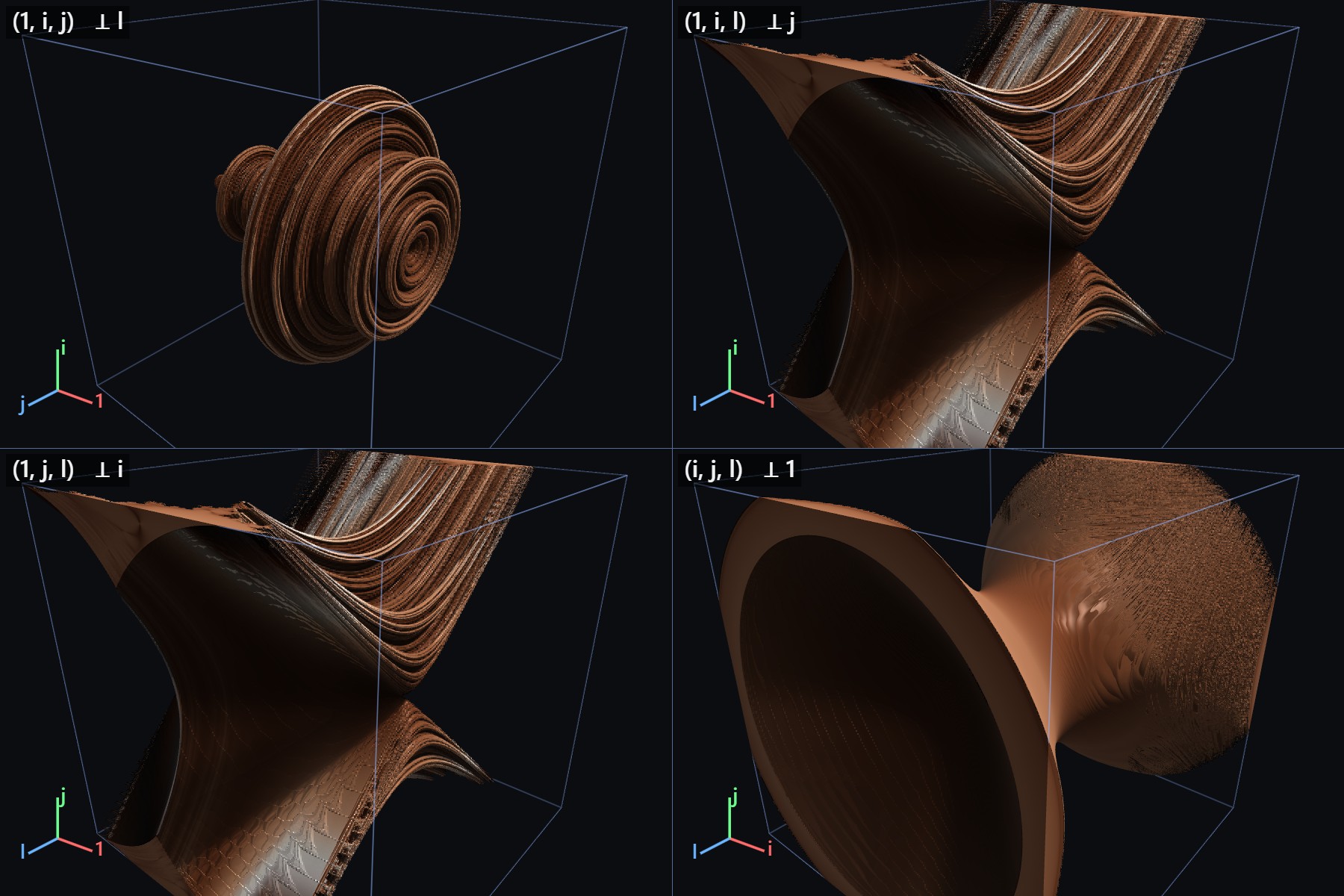}\par\vspace{1mm}
\noindent{\footnotesize \textbf{Map} $z\mapsto z\diamond z+c$, $c$ the sample point, $z_0=a$ $\cdot$ \textbf{Algebra} Split-octonions ($l^2=+1$) (Cayley--Dickson, $\gamma=(-1,-1,1)$) $\cdot$ \textbf{Frame} $X,Y,Z,N=1,i,j,l$ $\cdot$ center $(-0.35,\,\allowbreak 0,\,\allowbreak 0,\,\allowbreak 0,\,\allowbreak 0,\,\allowbreak 0,\,\allowbreak 0,\,\allowbreak 0)$ $\cdot$ half-extent 1.5 $\cdot$ 16 iterations $\cdot$ \textbf{Coloring} 8D orbit trap (point), gradient Copper}\par\vspace{1mm}
\noindent{\footnotesize\begin{tabular}{@{}lllll@{}}\toprule Slice (X, Y, Z) & Body [95\%] & Exact ($J\geq0.9995$) & Approximate ($0.97\leq J<0.9995$) & Dim\\\midrule
($1, i, j$) $\perp l$ & 15.2\% {\scriptsize[15.0, 15.4]} & $M$: yz; $R_x$: all & -- & 2 \\
($1, i, l$) $\perp j$ & 35.4\% {\scriptsize[35.1, 35.7]} & $M$: yz; $R_x$: 180 & -- & 2 \\
($1, j, l$) $\perp i$ & 35.4\% {\scriptsize[35.1, 35.7]} & $M$: yz; $R_x$: 180 & -- & 2 \\
($i, j, l$) $\perp 1$ & 35.3\% {\scriptsize[35.0, 35.6]} & $M$: xyz; $R_x$: 180; $R_y$: 180; $R_z$: all; $P$ & -- & 2 \\
\bottomrule\end{tabular}}\par\vspace{2mm}
The split-octonionic Mandelbrot set with seed $0$ ($l^2=+1$) keeps the classical Mandelbrot set of revolution in the quarter $(1,i,j)$, where the norm is positive definite: body share and symmetries coincide with those of Plate~F01 ($15.2\%$, all rotations about $1$ exact). The quarters $(1,i,l)$ and $(1,j,l)$ coincide ($35.4\%$ each): large saddle-shaped sheets, cut by the faces of the cube, with ribs along curved lines. The imaginary quarter is a trumpet with its waist in the middle and its axis along $l$, opening widely on both sides and hollow along the axis; it is exactly invariant under the rotations about $l$ by all tested angles.

By Proposition~\ref{atl:prop:span} the orbit of $0$ lies in the plane of $1$ and $c$, which is a copy of the complex plane when $\n(\im c)>0$ and of the split-complex plane when $\n(\im c)<0$; here $\n(\im c)=c_1^2+c_2^2+c_3^2-c_4^2-c_5^2-c_6^2-c_7^2$ is indefinite. The bounded set therefore depends only on $c_0$ and $\n(\im c)$: it is pulled back from the Mandelbrot set on one side of the null cone and from the hyperbolic Mandelbrot set, a square \cite{Senn}, on the other, and its level surfaces are hyperboloids. This is why, in the sections containing $l$, the rotations about $1$ are replaced by hyperbolic rotations, which are not Euclidean symmetries: only the mirrors $M_y$, $M_z$ and the half-turn about $1$ are exact there. In the imaginary quarter the body lies between two hyperboloids about the axis of $l$, which accounts for the trumpet and its hollow core.

The orbits lie in a subalgebra of dimension $2$. The rendered bodies depend on the Euclidean bailout, which is not an escape criterion in a split algebra: they are approximations of the bounded set.
\clearpage

\subsubsection{F06. Split Okubo, seed \texorpdfstring{$0$}{0}}\label{atl:pl:F06}\mbox{}\par\vspace{0.5mm}
\noindent\includegraphics[width=\textwidth]{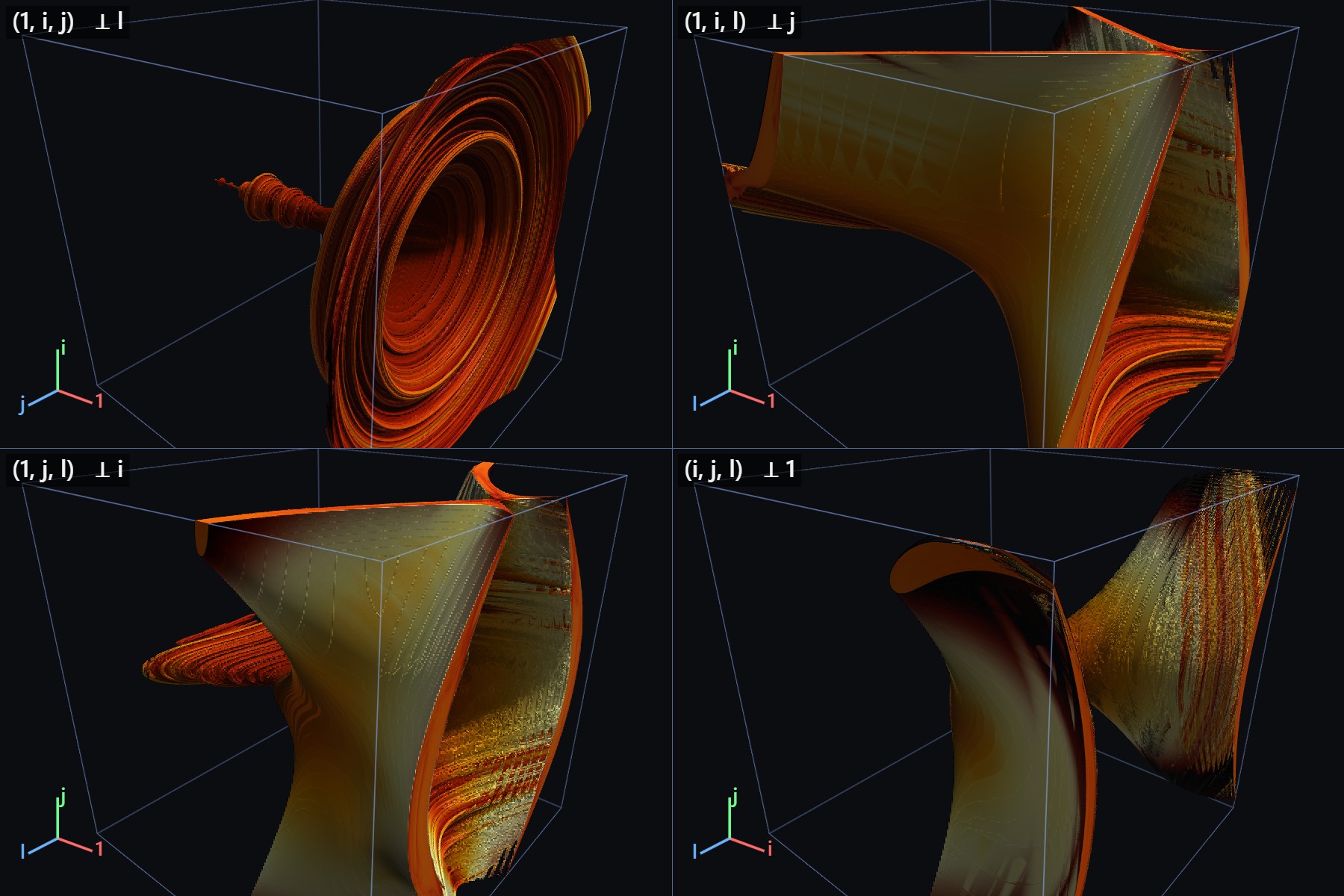}\par\vspace{1mm}
\noindent{\footnotesize \textbf{Map} $z\mapsto z\diamond z+c$, $c$ the sample point, $z_0=a$ $\cdot$ \textbf{Algebra} Split Okubo algebra (Petersson, $\gamma=(-1,-1,1)$, $\theta=120^\circ$) $\cdot$ \textbf{Frame} $X,Y,Z,N=1,i,j,l$ $\cdot$ center $(-0.35,\,\allowbreak 0,\,\allowbreak 0,\,\allowbreak 0,\,\allowbreak 0,\,\allowbreak 0,\,\allowbreak 0,\,\allowbreak 0)$ $\cdot$ half-extent 1.5 $\cdot$ 16 iterations $\cdot$ \textbf{Coloring} 8D orbit trap (point), gradient Fire}\par\vspace{1mm}
\noindent{\footnotesize\begin{tabular}{@{}lllll@{}}\toprule Slice (X, Y, Z) & Body [95\%] & Exact ($J\geq0.9995$) & Approximate ($0.97\leq J<0.9995$) & Dim\\\midrule
($1, i, j$) $\perp l$ & 9.2\% {\scriptsize[9.1, 9.4]} & $M$: yz; $R_x$: all & -- & 2 \\
($1, i, l$) $\perp j$ & 20.6\% {\scriptsize[20.3, 20.8]} & $M$: z & -- & 2 \\
($1, j, l$) $\perp i$ & 21.0\% {\scriptsize[20.8, 21.3]} & $M$: yz; $R_x$: 180 & -- & 2 \\
($i, j, l$) $\perp 1$ & 6.2\% {\scriptsize[6.1, 6.4]} & $M$: yz; $R_x$: 180 & -- & 2 \\
\bottomrule\end{tabular}}\par\vspace{2mm}
In the split Okubo algebra the quarter $(1,i,j)$ is the tricorn of revolution of Plate~F02 ($9.2\%$, all rotations about $1$ exact), since on $\HH$ the product is $\bar x\,\bar y$ as in the compact case. The three quarters through $l$ are dominated by large smooth sheets with sharp rims, cut by the faces of the cube: a curved, saddle-like sheet in $(1,i,l)$ ($20.6\%$), a similar sheet with a small ribbed fan on the left in $(1,j,l)$ ($21.0\%$), and two separate curved sheets in the imaginary quarter ($6.2\%$).

Theorem~B also covers this plate. The orbit of $0$ lies in the plane of $c$ and $c*c$ (dimension $2$), and the set is the pull-back of the tricorn where $\Delta=4\n(c)^3-\n(c,c*c)^2>0$ and of the \emph{hyperbolic tricorn}, its analogue over the split-complex numbers, where $\Delta<0$; the quarters through $l$ cross the surfaces $\Delta=0$, where the planar model changes. The finite symmetries are those of the compact Plate~F03: the mirror $M_z$ in the quarter $(1,i,l)$, the mirrors $M_y$, $M_z$ and the half-turn about the first axis in $(1,j,l)$ and in the imaginary quarter.

The picture must be read with two warnings. The Euclidean bailout of the renderer is not an escape criterion in a split algebra, and in double precision the invariant plane can be transversally unstable: rounding errors push the orbit out of the plane, the distance grows exponentially, and some points of the set escape although their exact orbit is bounded. The rendered bodies are therefore approximations of the Mandelbrot set, which Theorem~B describes exactly.
\clearpage

\subsubsection{F07. Semi-octonions, seed \texorpdfstring{$0$}{0}}\label{atl:pl:F07}\mbox{}\par\vspace{0.5mm}
\noindent\includegraphics[width=\textwidth]{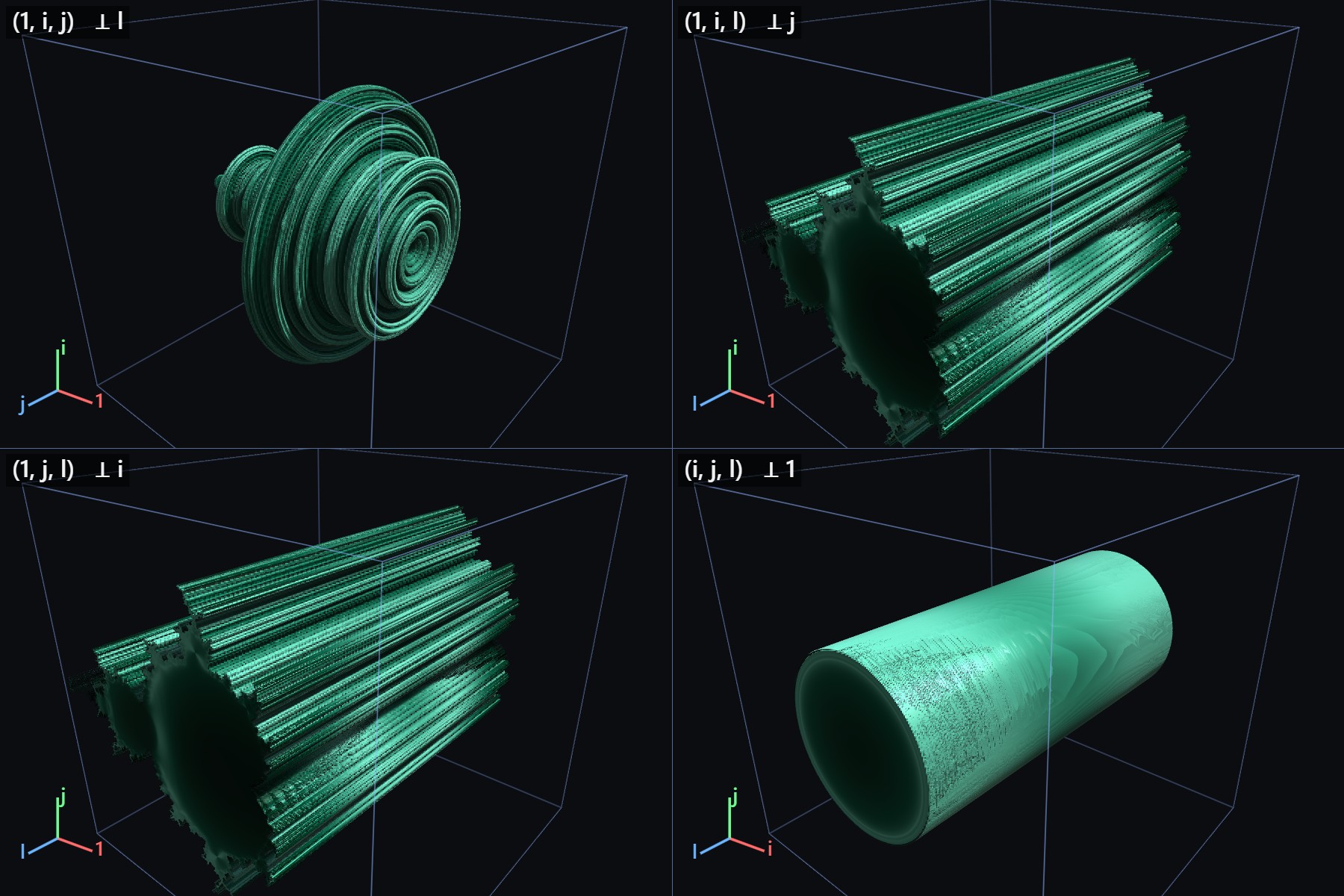}\par\vspace{1mm}
\noindent{\footnotesize \textbf{Map} $z\mapsto z\diamond z+c$, $c$ the sample point, $z_0=a$ $\cdot$ \textbf{Algebra} Semi-octonions ($l^2=0$) (Cayley--Dickson, $\gamma=(-1,-1,0)$) $\cdot$ \textbf{Frame} $X,Y,Z,N=1,i,j,l$ $\cdot$ center $(-0.35,\,\allowbreak 0,\,\allowbreak 0,\,\allowbreak 0,\,\allowbreak 0,\,\allowbreak 0,\,\allowbreak 0,\,\allowbreak 0)$ $\cdot$ half-extent 1.5 $\cdot$ 16 iterations $\cdot$ \textbf{Coloring} 8D orbit trap (point), gradient Emerald}\par\vspace{1mm}
\noindent{\footnotesize\begin{tabular}{@{}lllll@{}}\toprule Slice (X, Y, Z) & Body [95\%] & Exact ($J\geq0.9995$) & Approximate ($0.97\leq J<0.9995$) & Dim\\\midrule
($1, i, j$) $\perp l$ & 15.0\% {\scriptsize[14.8, 15.2]} & $M$: yz; $R_x$: all & -- & 2 \\
($1, i, l$) $\perp j$ & 39.0\% {\scriptsize[38.7, 39.3]} & $M$: yz; $R_x$: 180 & -- & 2 \\
($1, j, l$) $\perp i$ & 39.0\% {\scriptsize[38.7, 39.3]} & $M$: yz; $R_x$: 180 & -- & 2 \\
($i, j, l$) $\perp 1$ & 31.1\% {\scriptsize[30.8, 31.4]} & $M$: xyz; $R_x$: 180; $R_y$: 180; $R_z$: all; $P$ & -- & 2 \\
\bottomrule\end{tabular}}\par\vspace{2mm}
In the semi-octonions ($l^2=0$) the quarter $(1,i,j)$ is the classical Mandelbrot set of revolution ($15.0\%$, all rotations about $1$ exact, as in Plate~F01). The quarters $(1,i,l)$ and $(1,j,l)$ coincide ($39.0\%$): prisms along $l$, whose ridges run parallel to $l$ across the whole cube, and whose faces cut by the cube show the dark outline of the Mandelbrot set. The imaginary quarter is a circular cylinder along $l$, crossing the cube, exactly invariant under the rotations about $l$ by all tested angles ($31.1\%$).

The cylinders are explained by Proposition~\ref{atl:prop:rescaling}: for $\gamma_3=0$, $(p+ql)^2=p^2+2\re(p)\,q\,l$. The quaternionic part of the orbit is the orbit of the quaternionic Mandelbrot map and does not see the $l$-component of $c$; the $l$-part is driven linearly, $q_{n+1}=2\re(p_n)\,q_n+c_l$, and never feeds back. The body is therefore the quaternionic Mandelbrot set extended along $l$, as long as the linear growth of the $l$-part stays below the bailout. Only the mirrors $M_y$, $M_z$ and the half-turn about $1$ are exact in the quarters through $1$ and $l$, and the orbits lie in a subalgebra of dimension $2$.
\clearpage

\subsubsection{F08. Split-octonions (\texorpdfstring{$i^2=+1$}{i\texttwosuperior{}=+1}), seed \texorpdfstring{$0$}{0}}\label{atl:pl:F08}\mbox{}\par\vspace{0.5mm}
\noindent\includegraphics[width=\textwidth]{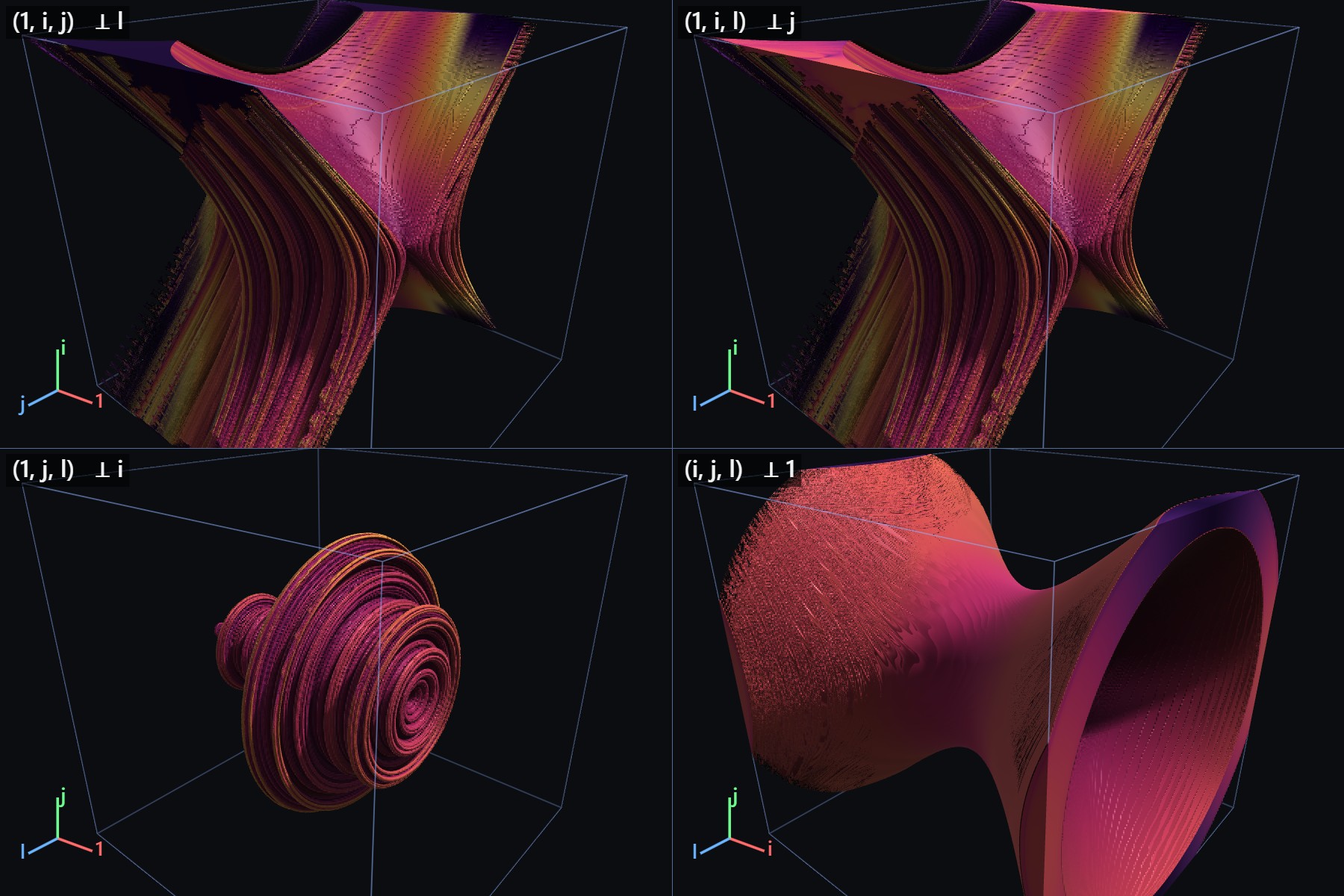}\par\vspace{1mm}
\noindent{\footnotesize \textbf{Map} $z\mapsto z\diamond z+c$, $c$ the sample point, $z_0=a$ $\cdot$ \textbf{Algebra} Split-octonions ($i^2=+1$) (Cayley--Dickson, $\gamma=(1,-1,-1)$) $\cdot$ \textbf{Frame} $X,Y,Z,N=1,i,j,l$ $\cdot$ center $(-0.35,\,\allowbreak 0,\,\allowbreak 0,\,\allowbreak 0,\,\allowbreak 0,\,\allowbreak 0,\,\allowbreak 0,\,\allowbreak 0)$ $\cdot$ half-extent 1.5 $\cdot$ 16 iterations $\cdot$ \textbf{Coloring} 8D orbit trap (point), gradient Sunset}\par\vspace{1mm}
\noindent{\footnotesize\begin{tabular}{@{}lllll@{}}\toprule Slice (X, Y, Z) & Body [95\%] & Exact ($J\geq0.9995$) & Approximate ($0.97\leq J<0.9995$) & Dim\\\midrule
($1, i, j$) $\perp l$ & 35.2\% {\scriptsize[34.9, 35.5]} & $M$: yz; $R_x$: 180 & -- & 2 \\
($1, i, l$) $\perp j$ & 35.2\% {\scriptsize[34.9, 35.5]} & $M$: yz; $R_x$: 180 & -- & 2 \\
($1, j, l$) $\perp i$ & 15.0\% {\scriptsize[14.8, 15.2]} & $M$: yz; $R_x$: all & -- & 2 \\
($i, j, l$) $\perp 1$ & 35.2\% {\scriptsize[34.9, 35.5]} & $M$: xyz; $R_x$: all; $R_y$: 180; $R_z$: 180; $P$ & -- & 2 \\
\bottomrule\end{tabular}}\par\vspace{2mm}
With $i^2=+1$ the plane of $1$ and $i$ is the split-complex plane, where the Mandelbrot set is a square \cite{Senn}. The quarters $(1,i,j)$ and $(1,i,l)$ coincide ($35.2\%$ each): X-shaped sheets, pinched in the middle like an hourglass and cut by the faces of the cube, with the mirrors $M_y$, $M_z$ and the half-turn about $1$ exact. The quarter $(1,j,l)$, where both hidden directions are elliptic, is the ordinary Mandelbrot set of revolution ($15.0\%$, all rotations about $1$ exact). The imaginary quarter is a trumpet with its axis along $i$, exactly invariant under the rotations about $i$ by all tested angles.

As in Plate~F05, the orbit of $0$ lies in the plane of $1$ and $c$, and the bounded set depends only on $c_0$ and $\n(\im c)=-c_1^2+c_2^2-c_3^2+c_4^2-c_5^2+c_6^2-c_7^2$: it is pulled back from the Mandelbrot set where $\n(\im c)>0$ and from the hyperbolic square where $\n(\im c)<0$. The plate places side by side, inside the same algebra, the two classical Mandelbrot sets of the plane, the complex one and the hyperbolic one. The orbits lie in a subalgebra of dimension $2$; as for every split algebra, the rendered bodies depend on the Euclidean bailout and are approximations.
\clearpage

\subsubsection{F09. Octonions with a seed}\label{atl:pl:F09}\mbox{}\par\vspace{0.5mm}
\noindent\includegraphics[width=\textwidth]{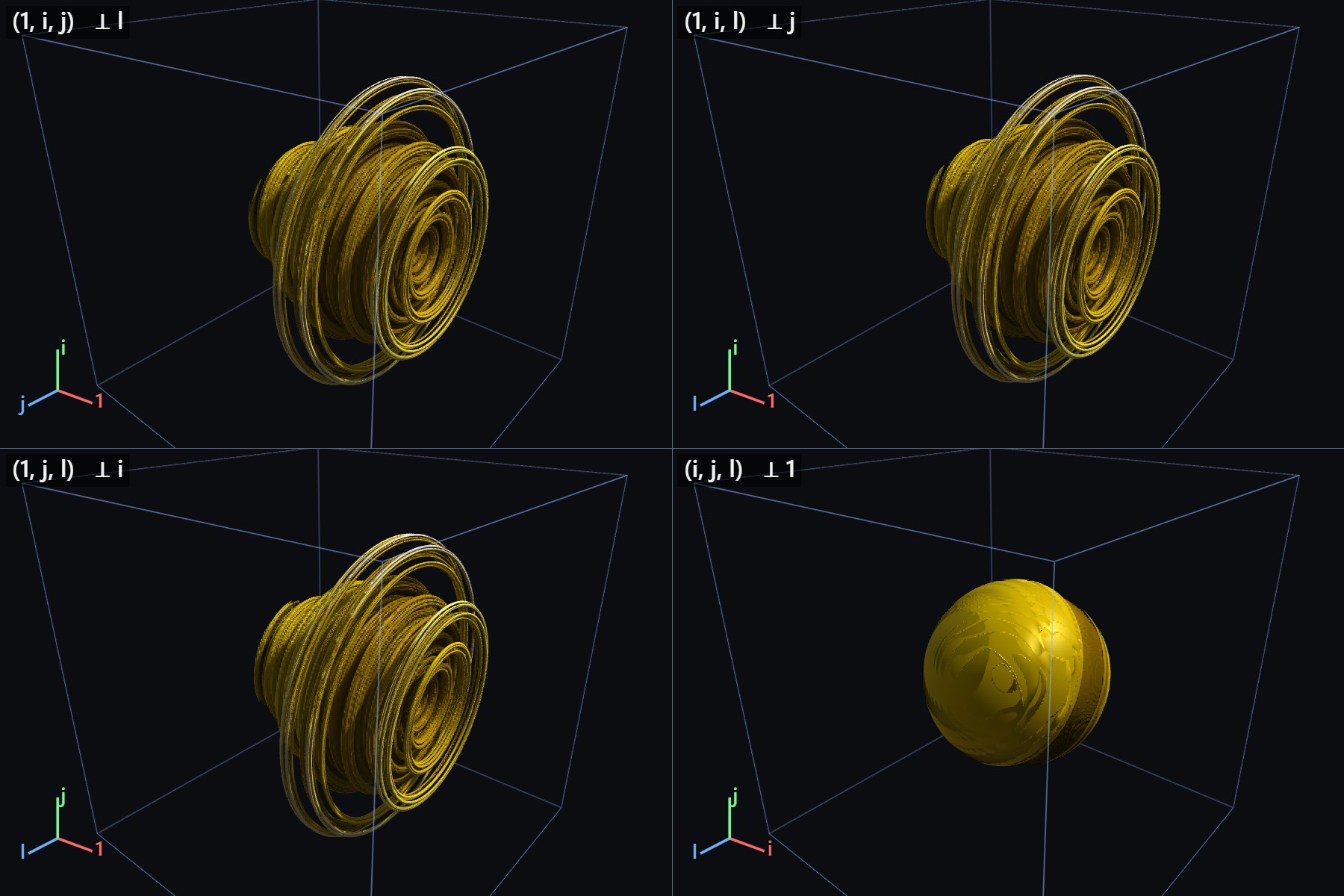}\par\vspace{1mm}
\noindent{\footnotesize \textbf{Map} $z\mapsto z\diamond z+c$, $c$ the sample point, $z_0=a$ $\cdot$ \textbf{Algebra} Octonions (Cayley--Dickson, $\gamma=(-1,-1,-1)$) $\cdot$ $a=(0.15,\,\allowbreak -0.064,\,\allowbreak 0.179,\,\allowbreak -0.243,\,\allowbreak 0.328,\,\allowbreak -0.138,\,\allowbreak 0.23,\,\allowbreak 0.066)$ $\cdot$ \textbf{Frame} $X,Y,Z,N=1,i,j,l$ $\cdot$ center $(-0.35,\,\allowbreak 0,\,\allowbreak 0,\,\allowbreak 0,\,\allowbreak 0,\,\allowbreak 0,\,\allowbreak 0,\,\allowbreak 0)$ $\cdot$ half-extent 1.5 $\cdot$ 16 iterations $\cdot$ \textbf{Coloring} 8D orbit trap (point), gradient Gold}\par\vspace{1mm}
\noindent{\footnotesize\begin{tabular}{@{}lllll@{}}\toprule Slice (X, Y, Z) & Body [95\%] & Exact ($J\geq0.9995$) & Approximate ($0.97\leq J<0.9995$) & Dim\\\midrule
($1, i, j$) $\perp l$ & 13.3\% {\scriptsize[13.1, 13.5]} & -- & -- & 3 \\
($1, i, l$) $\perp j$ & 13.2\% {\scriptsize[13.0, 13.4]} & -- & -- & 3 \\
($1, j, l$) $\perp i$ & 13.2\% {\scriptsize[13.0, 13.4]} & -- & -- & 3 \\
($i, j, l$) $\perp 1$ & 12.1\% {\scriptsize[11.9, 12.3]} & -- & -- & 3 \\
\bottomrule\end{tabular}}\par\vspace{2mm}
With a seed $a\ne0$ the first iterate is $z_1=a^2+c$: the seed shifts the first iterate by the fixed vector $a^2$, and the orbit is no longer confined to the plane of $1$ and $c$. It lies in the span of $1$, $\im a$ and $\im c$, a subalgebra of dimension $3$. The seed of this plate is the generic eight-dimensional seed of Plate~F11, of norm $\n(a)\approx0.30$. The three quarters through $1$ show the ribbed body of Plate~F01, less regular and surrounded by loose rings; they are similar, with nearly equal body shares ($13.2$--$13.3\%$), but not identical, because the seed projects differently onto each section. The imaginary quarter is no longer a ball: a rounded body ($12.1\%$) with a raised band on one side.

No symmetry of the table survives, not even approximately: the best scores stay below $J=0.91$. Proposition~\ref{atl:prop:equivariance} still acts, since every orthogonal map fixing $1$ and $a$ commutes with the map, but the symmetries it gives are not coordinate ones: in each section through $1$, the mirror in the plane of $1$ and of the projection of $\im a$; in the imaginary quarter, the rotations about the projection of $\im a$ on $\langle i,j,l\rangle$, a tilted axis, compatible with the raised band. None of them is tested by the table.
\clearpage

\subsubsection{F10. Para-octonions with a seed}\label{atl:pl:F10}\mbox{}\par\vspace{0.5mm}
\noindent\includegraphics[width=\textwidth]{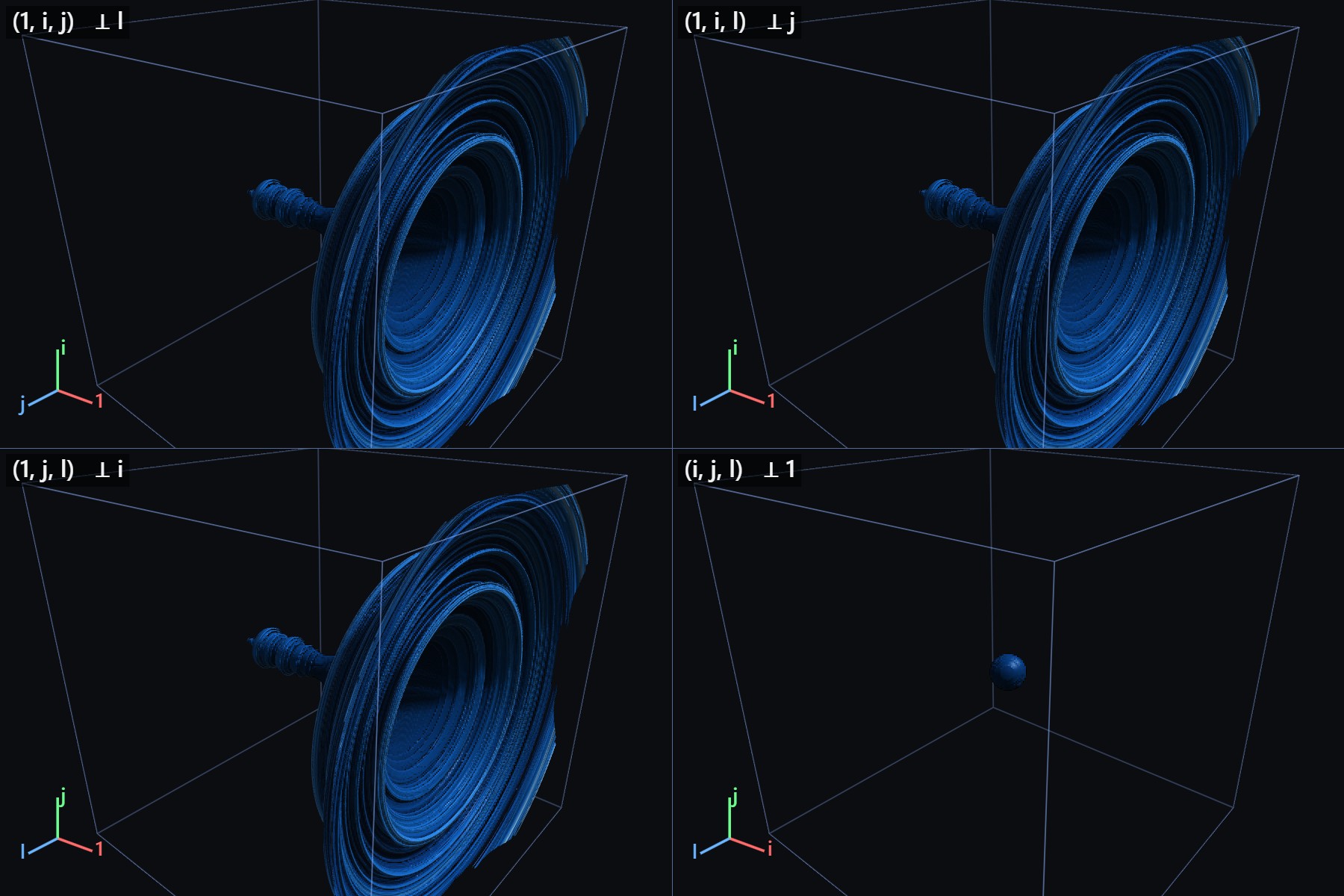}\par\vspace{1mm}
\noindent{\footnotesize \textbf{Map} $z\mapsto z\diamond z+c$, $c$ the sample point, $z_0=a$ $\cdot$ \textbf{Algebra} Para-octonions (para-Hurwitz $\bar x\bar y$, $\gamma=(-1,-1,-1)$) $\cdot$ $a=(0.209,\,\allowbreak -0.021,\,\allowbreak -0.035,\,\allowbreak -0.016,\,\allowbreak -0.04,\,\allowbreak 0.193,\,\allowbreak -0.022,\,\allowbreak 0.037)$ $\cdot$ \textbf{Frame} $X,Y,Z,N=1,i,j,l$ $\cdot$ center $(-0.35,\,\allowbreak 0,\,\allowbreak 0,\,\allowbreak 0,\,\allowbreak 0,\,\allowbreak 0,\,\allowbreak 0,\,\allowbreak 0)$ $\cdot$ half-extent 1.5 $\cdot$ 16 iterations $\cdot$ \textbf{Coloring} 8D orbit trap (point), gradient Deep Blue}\par\vspace{1mm}
\noindent{\footnotesize\begin{tabular}{@{}lllll@{}}\toprule Slice (X, Y, Z) & Body [95\%] & Exact ($J\geq0.9995$) & Approximate ($0.97\leq J<0.9995$) & Dim\\\midrule
($1, i, j$) $\perp l$ & 9.3\% {\scriptsize[9.1, 9.4]} & -- & -- & 3 \\
($1, i, l$) $\perp j$ & 9.3\% {\scriptsize[9.1, 9.4]} & -- & -- & 3 \\
($1, j, l$) $\perp i$ & 9.3\% {\scriptsize[9.1, 9.4]} & -- & -- & 3 \\
($i, j, l$) $\perp 1$ & 0.1\% {\scriptsize[0.1, 0.1]} & $R_x$: 37 & $M$: xyz; $R_x$: 180,120,90,72,60; $R_y$: all; $R_z$: all (min.\ 0.975) & 3 \\
\bottomrule\end{tabular}}\par\vspace{2mm}
The para-octonionic Mandelbrot set with the seed found by the parameter scan of HF8D~0.2 is very close to the tricorn of revolution of Plate~F02: the same funnels and beaded stems in the three quarters through $1$, with the same body share ($9.3\%$), and the same tiny ball in the imaginary quarter ($0.1\%$). Since $z_1=a\bullet a+c$, the seed only shifts the first iterate by a fixed vector of norm $\n(a)\approx0.09$, and its effect is a slight deformation.

The rotations about $1$ are no longer symmetries, but they come close: they score $J=0.95$--$0.96$, just below the threshold of the approximate symmetries. The reason is Proposition~\ref{atl:prop:equivariance}: the exact symmetries are the orthogonal maps fixing $1$ and $a$, and the projections of $\im a$ onto the three sections are short, since the seed is concentrated on $1$ and $il$. The orbits lie in a subalgebra of dimension $3$, spanned by $1$, $\im a$ and $\im c$. The imaginary quarter is too small for its symmetries to be measured reliably.
\clearpage

\subsubsection{F11. Okubo with a seed}\label{atl:pl:F11}\mbox{}\par\vspace{0.5mm}
\noindent\includegraphics[width=\textwidth]{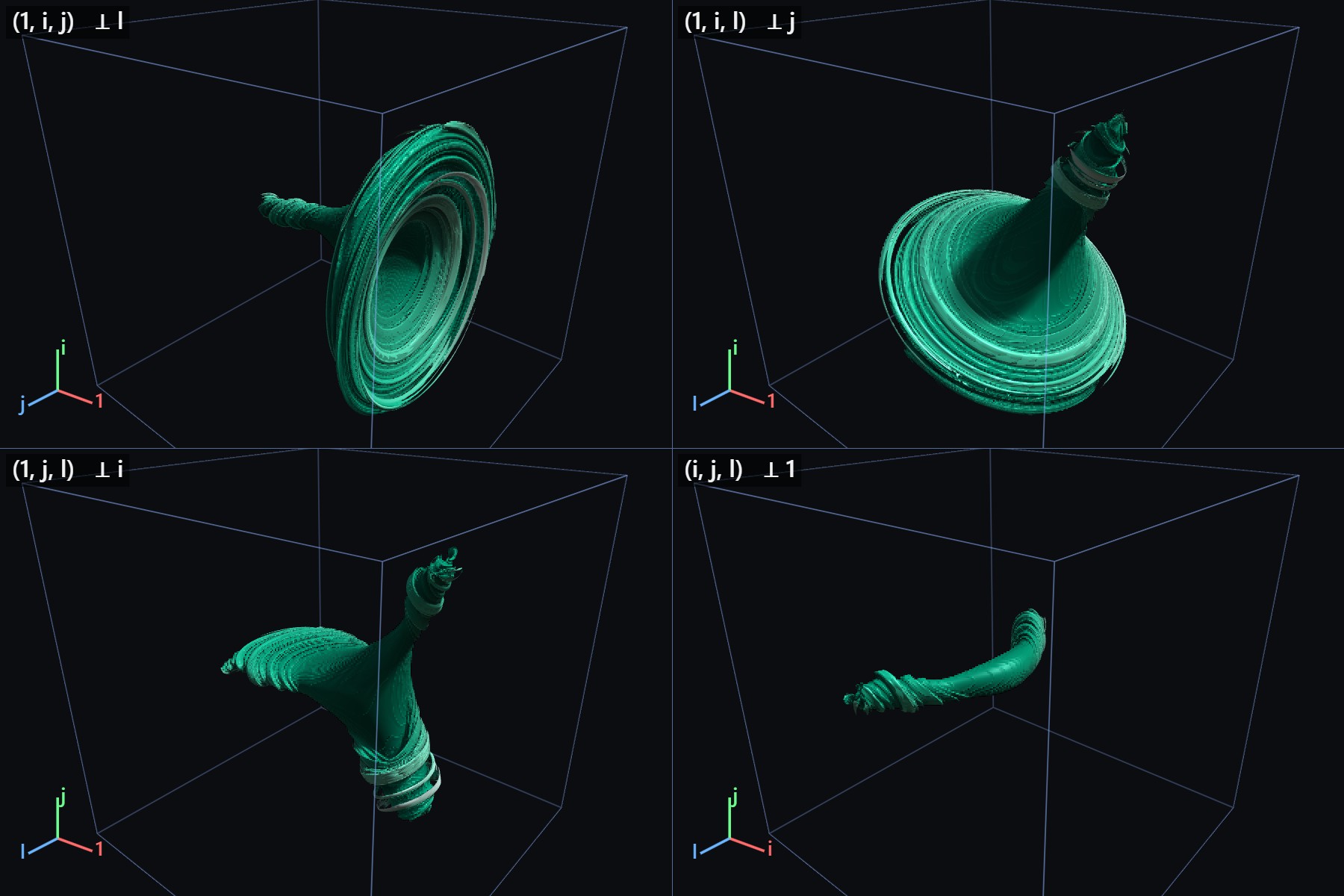}\par\vspace{1mm}
\noindent{\footnotesize \textbf{Map} $z\mapsto z\diamond z+c$, $c$ the sample point, $z_0=a$ $\cdot$ \textbf{Algebra} Okubo algebra (Petersson, $\gamma=(-1,-1,-1)$, $\theta=120^\circ$) $\cdot$ $a=(0.15,\,\allowbreak -0.064,\,\allowbreak 0.179,\,\allowbreak -0.243,\,\allowbreak 0.328,\,\allowbreak -0.138,\,\allowbreak 0.23,\,\allowbreak 0.066)$ $\cdot$ \textbf{Frame} $X,Y,Z,N=1,i,j,l$ $\cdot$ center $(-0.35,\,\allowbreak 0,\,\allowbreak 0,\,\allowbreak 0,\,\allowbreak 0,\,\allowbreak 0,\,\allowbreak 0,\,\allowbreak 0)$ $\cdot$ half-extent 1.5 $\cdot$ 16 iterations $\cdot$ \textbf{Coloring} 8D orbit trap (point), gradient Aurora}\par\vspace{1mm}
\noindent{\footnotesize\begin{tabular}{@{}lllll@{}}\toprule Slice (X, Y, Z) & Body [95\%] & Exact ($J\geq0.9995$) & Approximate ($0.97\leq J<0.9995$) & Dim\\\midrule
($1, i, j$) $\perp l$ & 7.9\% {\scriptsize[7.7, 8.0]} & -- & -- & 8 \\
($1, i, l$) $\perp j$ & 6.8\% {\scriptsize[6.6, 6.9]} & -- & -- & 8 \\
($1, j, l$) $\perp i$ & 3.5\% {\scriptsize[3.3, 3.6]} & -- & -- & 8 \\
($i, j, l$) $\perp 1$ & 0.8\% {\scriptsize[0.7, 0.8]} & -- & -- & 8 \\
\bottomrule\end{tabular}}\par\vspace{2mm}
The Okubonic Mandelbrot set with a generic seed keeps the shapes of Plate~F03, smaller and less regular: a funnel with a stem in the quarter $(1,i,j)$ ($7.9\%$), a tilted funnel with a stem pointing upwards and to the right in $(1,i,l)$ ($6.8\%$), a three-armed body in $(1,j,l)$ ($3.5\%$), and a bent tube in the imaginary quarter ($0.8\%$). Spiral ribbons wind around the stems.

The seed destroys every symmetry of the table, including the finite ones of Plate~F03: the best scores are $J=0.81$ for the rotations about $1$ in the first quarter and $0.72$ or less in the others. The Mandelbrot set with a seed is not a planar pull-back: the orbits lie in a subalgebra of dimension $8$, and the invariants of the pair $(a,c)$ that govern it are not known (an open problem, see the main text). The automorphisms that fix $a$ still commute with the map, but for a generic seed they form a torus whose root planes are in general position with respect to the coordinate axes, and they give no coordinate symmetry.
\clearpage

\subsubsection{F12. Okubo with a seed, in the frame \texorpdfstring{$F_B$}{FB}}\label{atl:pl:F12}\mbox{}\par\vspace{0.5mm}
\noindent\includegraphics[width=\textwidth]{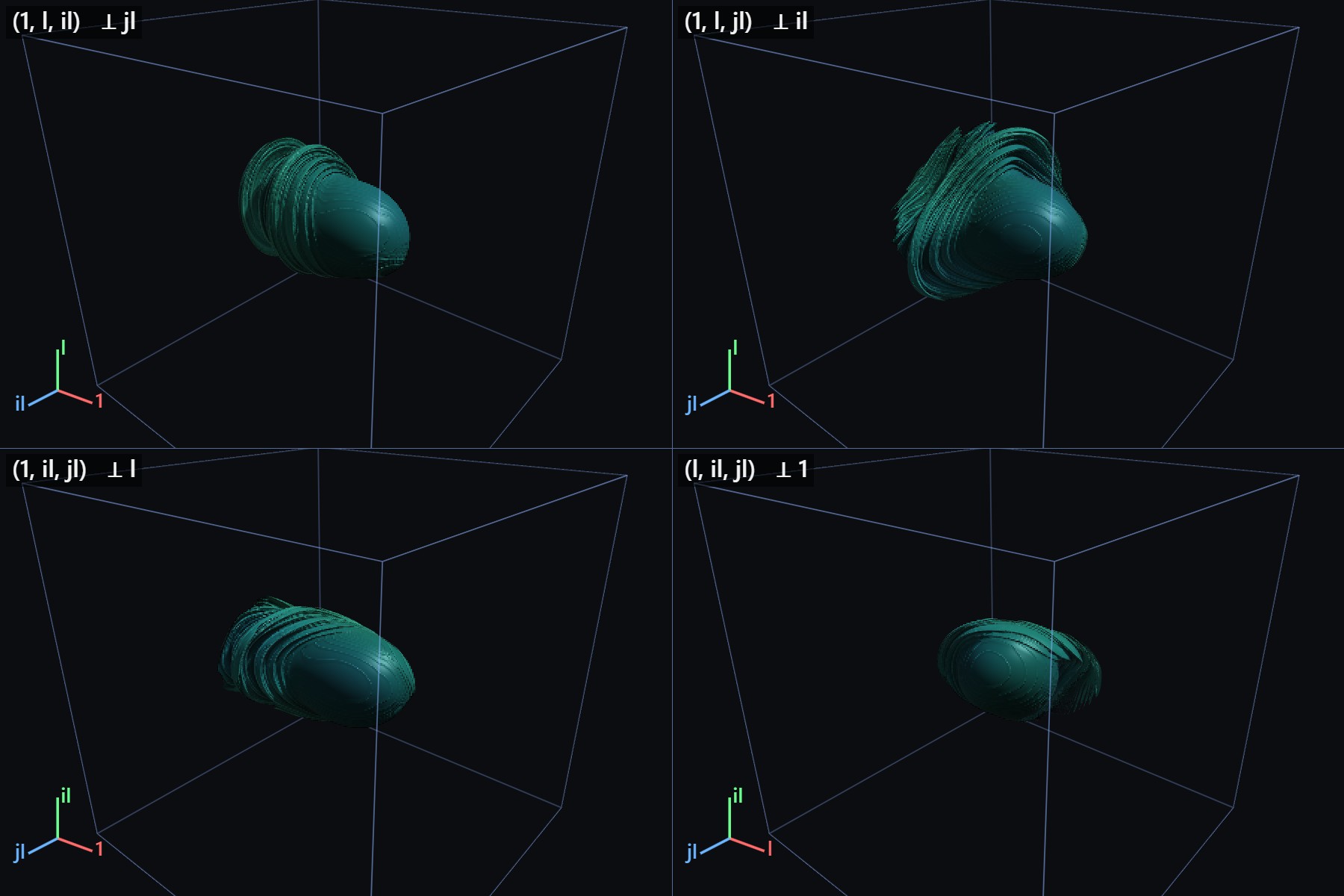}\par\vspace{1mm}
\noindent{\footnotesize \textbf{Map} $z\mapsto z\diamond z+c$, $c$ the sample point, $z_0=a$ $\cdot$ \textbf{Algebra} Okubo algebra (Petersson, $\gamma=(-1,-1,-1)$, $\theta=120^\circ$) $\cdot$ $a=(0.15,\,\allowbreak -0.064,\,\allowbreak 0.179,\,\allowbreak -0.243,\,\allowbreak 0.328,\,\allowbreak -0.138,\,\allowbreak 0.23,\,\allowbreak 0.066)$ $\cdot$ \textbf{Frame} $X,Y,Z,N=1,l,il,jl$ $\cdot$ center $(-0.35,\,\allowbreak 0,\,\allowbreak 0,\,\allowbreak 0,\,\allowbreak 0,\,\allowbreak 0,\,\allowbreak 0,\,\allowbreak 0)$ $\cdot$ half-extent 1.5 $\cdot$ 16 iterations $\cdot$ \textbf{Coloring} 8D orbit trap (point), gradient Viridis-like}\par\vspace{1mm}
\noindent{\footnotesize\begin{tabular}{@{}lllll@{}}\toprule Slice (X, Y, Z) & Body [95\%] & Exact ($J\geq0.9995$) & Approximate ($0.97\leq J<0.9995$) & Dim\\\midrule
($1, l, il$) $\perp jl$ & 5.3\% {\scriptsize[5.1, 5.4]} & -- & -- & 8 \\
($1, l, jl$) $\perp il$ & 6.3\% {\scriptsize[6.2, 6.5]} & -- & -- & 8 \\
($1, il, jl$) $\perp l$ & 6.5\% {\scriptsize[6.4, 6.7]} & -- & -- & 8 \\
($l, il, jl$) $\perp 1$ & 2.9\% {\scriptsize[2.8, 3.0]} & -- & -- & 8 \\
\bottomrule\end{tabular}}\par\vspace{2mm}
In the frame $F_B$ the seed breaks the revolutions of Plate~F04. The quarters through $1$ are small eggs ($5$--$7\%$ of the ball), with a smooth cap on the side of positive real part and ribbed, flaky layers behind; they are no longer identical, and the imaginary quarter is a small flattened egg ($2.9\%$) instead of a ball. No symmetry of the table survives, not even approximately (the best score is $J=0.87$).

The revolutions of Plate~F04 came from the fact that, on $\RR\oplus\HH l$, the two invariants $\n(c)$ and $\n(c,c*c)$ depend only on $c_0$ and $|q|$ (Proposition~\ref{atl:prop:okubocubic}). With a seed the orbit leaves the plane of $c$ and $c*c$ and lies in a subalgebra of dimension $8$, and these two invariants no longer decide the escape. The comparison with Plate~F04 isolates the effect of the seed.
\clearpage

\subsubsection{F13. Split para-octonions with a seed}\label{atl:pl:F13}\mbox{}\par\vspace{0.5mm}
\noindent\includegraphics[width=\textwidth]{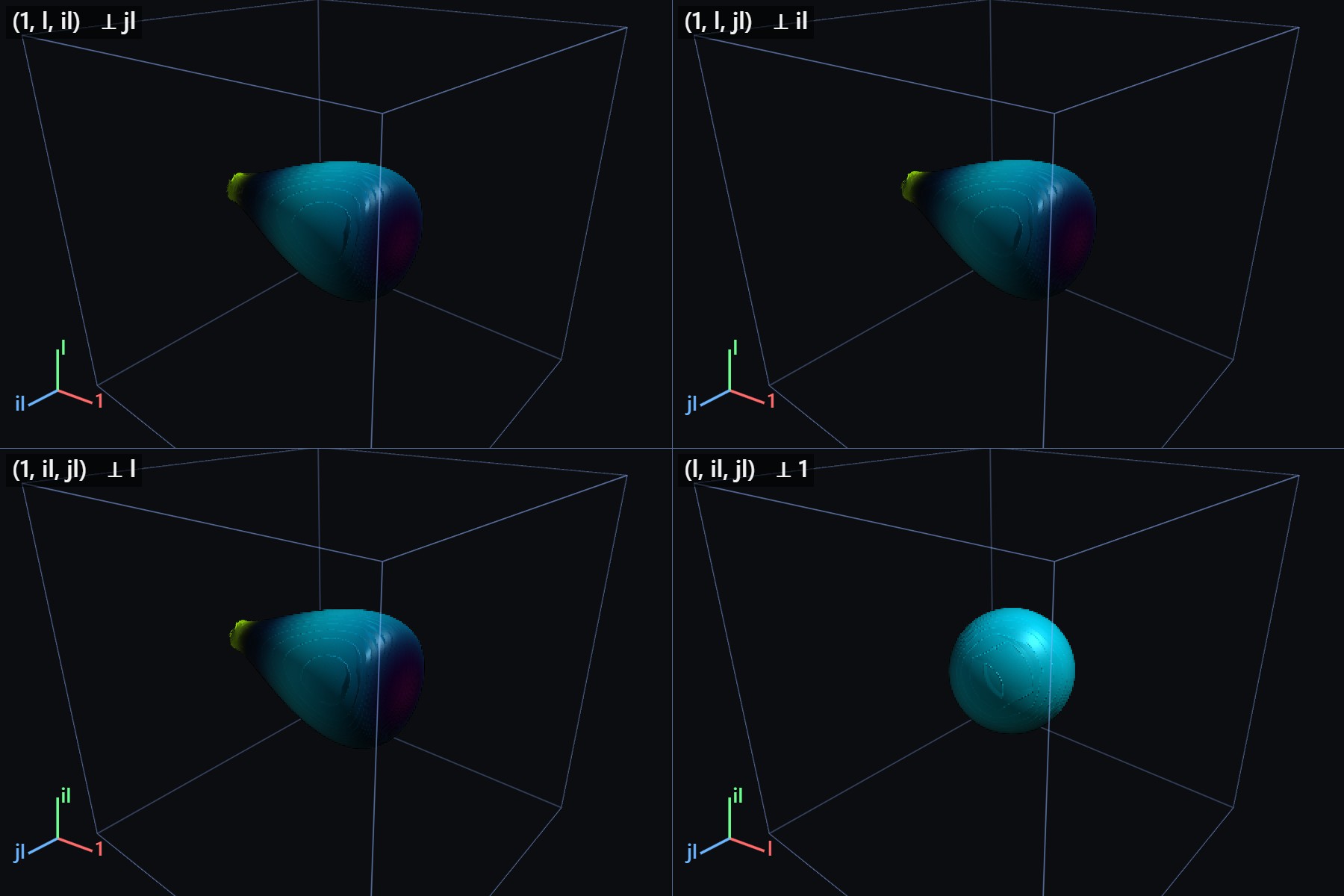}\par\vspace{1mm}
\noindent{\footnotesize \textbf{Map} $z\mapsto z\diamond z+c$, $c$ the sample point, $z_0=a$ $\cdot$ \textbf{Algebra} Split para-octonions (para-Hurwitz $\bar x\bar y$, $\gamma=(-1,-1,1)$) $\cdot$ $a=(-0.141,\,\allowbreak 0.012,\,\allowbreak -0.219,\,\allowbreak 0.535,\,\allowbreak -0.046,\,\allowbreak -0.005,\,\allowbreak 0.171,\,\allowbreak -0.003)$ $\cdot$ \textbf{Frame} $X,Y,Z,N=1,l,il,jl$ $\cdot$ center $(-0.35,\,\allowbreak 0,\,\allowbreak 0,\,\allowbreak 0,\,\allowbreak 0,\,\allowbreak 0,\,\allowbreak 0,\,\allowbreak 0)$ $\cdot$ half-extent 1.5 $\cdot$ 16 iterations $\cdot$ \textbf{Coloring} 8D orbit trap (point), gradient Neon}\par\vspace{1mm}
\noindent{\footnotesize\begin{tabular}{@{}lllll@{}}\toprule Slice (X, Y, Z) & Body [95\%] & Exact ($J\geq0.9995$) & Approximate ($0.97\leq J<0.9995$) & Dim\\\midrule
($1, l, il$) $\perp jl$ & 7.0\% {\scriptsize[6.9, 7.2]} & -- & $M$: z; $R_x$: 37 (min.\ 0.974) & 3 \\
($1, l, jl$) $\perp il$ & 7.1\% {\scriptsize[7.0, 7.3]} & -- & -- & 3 \\
($1, il, jl$) $\perp l$ & 7.1\% {\scriptsize[6.9, 7.3]} & -- & $M$: y (min.\ 0.990) & 3 \\
($l, il, jl$) $\perp 1$ & 4.5\% {\scriptsize[4.3, 4.6]} & -- & $M$: y; $R_z$: 60,37 (min.\ 0.972) & 3 \\
\bottomrule\end{tabular}}\par\vspace{2mm}
The split para-octonions with a seed, in the frame $F_B$, give three teardrop-shaped bodies in the quarters through $1$, with a small bright tip on the side of negative real part ($7.0$--$7.1\%$ of the ball), and a nearly round body in the imaginary quarter ($4.5\%$). In these sections $\im c$ lies in $\HH l$, where the norm is negative definite; without the seed, the sections would be pulled back from the hyperbolic tricorn by Theorem~B.

No symmetry is exact, but several are approximate, and Proposition~\ref{atl:prop:equivariance} explains which: the exact symmetries are the maps fixing $1$ and $a$ and preserving the norm, which include the rotations of $\HH l$ that fix the $\HH l$-part of $a$. This part is almost along $jl$. Accordingly, in the quarter $(1,l,il)$, where the seed has only a small projection, the mirror $M_z$ ($J=0.991$) and the rotation by $37^\circ$ about $1$ ($0.974$) are approximate; in $(1,il,jl)$ the mirror $M_y$ is approximate ($0.990$); in $(1,l,jl)$, where the projection is oblique, nothing comes close; and the imaginary quarter is nearly a body of revolution about $jl$ (mirror $M_y$ $0.995$, rotations about $jl$ by $37^\circ$ and $60^\circ$ $0.984$ and $0.972$). The orbits lie in a subalgebra of dimension $3$, spanned by $1$, $\im a$ and $\im c$.
\clearpage

\subsection{Part G. The Mandelbar map}\label{atl:part:G}
The Mandelbar map $z\mapsto\overline{z\diamond z}+c$ provides the cleanest check of the para-Hurwitz construction. Since $x\bullet x=\overline{x^2}$ in the para-octonions, the Mandelbar map of one algebra is the Mandelbrot map of the other.

\subsubsection{G01. The octonionic Mandelbar is the para-octonionic Mandelbrot}\label{atl:pl:G01}\mbox{}\par\vspace{0.5mm}
\noindent\includegraphics[width=\textwidth]{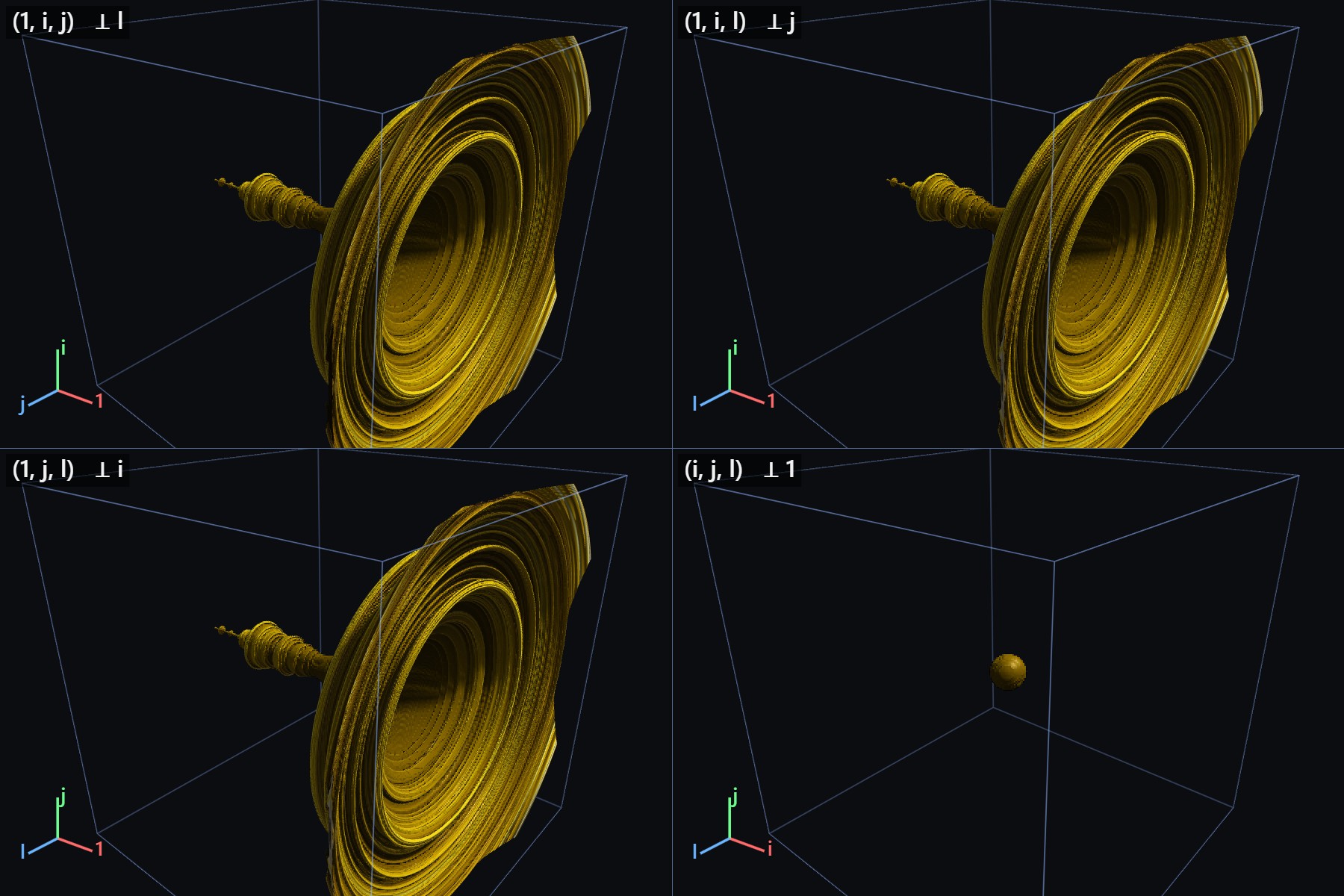}\par\vspace{1mm}
\noindent{\footnotesize \textbf{Map} $z\mapsto\overline{z\diamond z}+c$, $c$ the sample point $\cdot$ \textbf{Algebra} Octonions (Cayley--Dickson, $\gamma=(-1,-1,-1)$) $\cdot$ \textbf{Frame} $X,Y,Z,N=1,i,j,l$ $\cdot$ center $(-0.35,\,\allowbreak 0,\,\allowbreak 0,\,\allowbreak 0,\,\allowbreak 0,\,\allowbreak 0,\,\allowbreak 0,\,\allowbreak 0)$ $\cdot$ half-extent 1.5 $\cdot$ 16 iterations $\cdot$ \textbf{Coloring} 8D orbit trap (point), gradient Gold}\par\vspace{1mm}
\noindent{\footnotesize\begin{tabular}{@{}lllll@{}}\toprule Slice (X, Y, Z) & Body [95\%] & Exact ($J\geq0.9995$) & Approximate ($0.97\leq J<0.9995$) & Dim\\\midrule
($1, i, j$) $\perp l$ & 9.3\% {\scriptsize[9.1, 9.5]} & $M$: yz; $R_x$: all & -- & 2 \\
($1, i, l$) $\perp j$ & 9.3\% {\scriptsize[9.1, 9.5]} & $M$: yz; $R_x$: all & -- & 2 \\
($1, j, l$) $\perp i$ & 9.3\% {\scriptsize[9.1, 9.5]} & $M$: yz; $R_x$: all & -- & 2 \\
($i, j, l$) $\perp 1$ & 0.1\% {\scriptsize[0.1, 0.1]} & $M$: xyz; $R_x$: all; $R_y$: all; $R_z$: all; $P$ & -- & 2 \\
\bottomrule\end{tabular}}\par\vspace{2mm}
The octonionic Mandelbar map $\overline{z^2}+c$ is, by definition, the para-octonionic Mandelbrot map, since $z\bullet z=\bar z\,\bar z=\overline{z^2}$, and the plate reproduces Plate~F02. The three quarters through $1$ show the same body: a wide flared disc, cut by the faces of the cube, attached to a thin ribbed stalk that runs along the axis of $1$ towards negative real part. The imaginary quarter is the hyperplane $\re c=-0.35$ through the centre of the box, and it meets the body in a small ball.

By Theorem~B, the orbit of the seed $0$ stays in the plane of $1$ and $c$, and the para-octonionic Mandelbrot set is the tricorn turned about the real axis, the tricorn of revolution: every section through the axis of $1$ contains the same tricorn, and every section orthogonal to it is a ball. The tricorn has three long arms, at $180^\circ$ and at $\pm60^\circ$; the arm at $180^\circ$ becomes the stalk, and the two arms at $\pm60^\circ$ sweep the flared disc on the side of positive real part. At $\re c=-0.35$ the tricorn is only about $0.1$ high, whence the small ball. The measurements agree with Plate~F02: in each quarter through $1$ the rotations about $1$ by every tested angle are exact, the body share is $9.3\%$, and the imaginary quarter is a ball occupying $0.1\%$ of the measuring ball; the orbits lie in a subalgebra of dimension $2$. Since the variable of a Mandelbrot set is the constant, the point reflection $z\mapsto-z$ is not a symmetry here (Remark~\ref{atl:rem:minus}).
\clearpage

\subsubsection{G02. The para-octonionic Mandelbar is the octonionic Mandelbrot}\label{atl:pl:G02}\mbox{}\par\vspace{0.5mm}
\noindent\includegraphics[width=\textwidth]{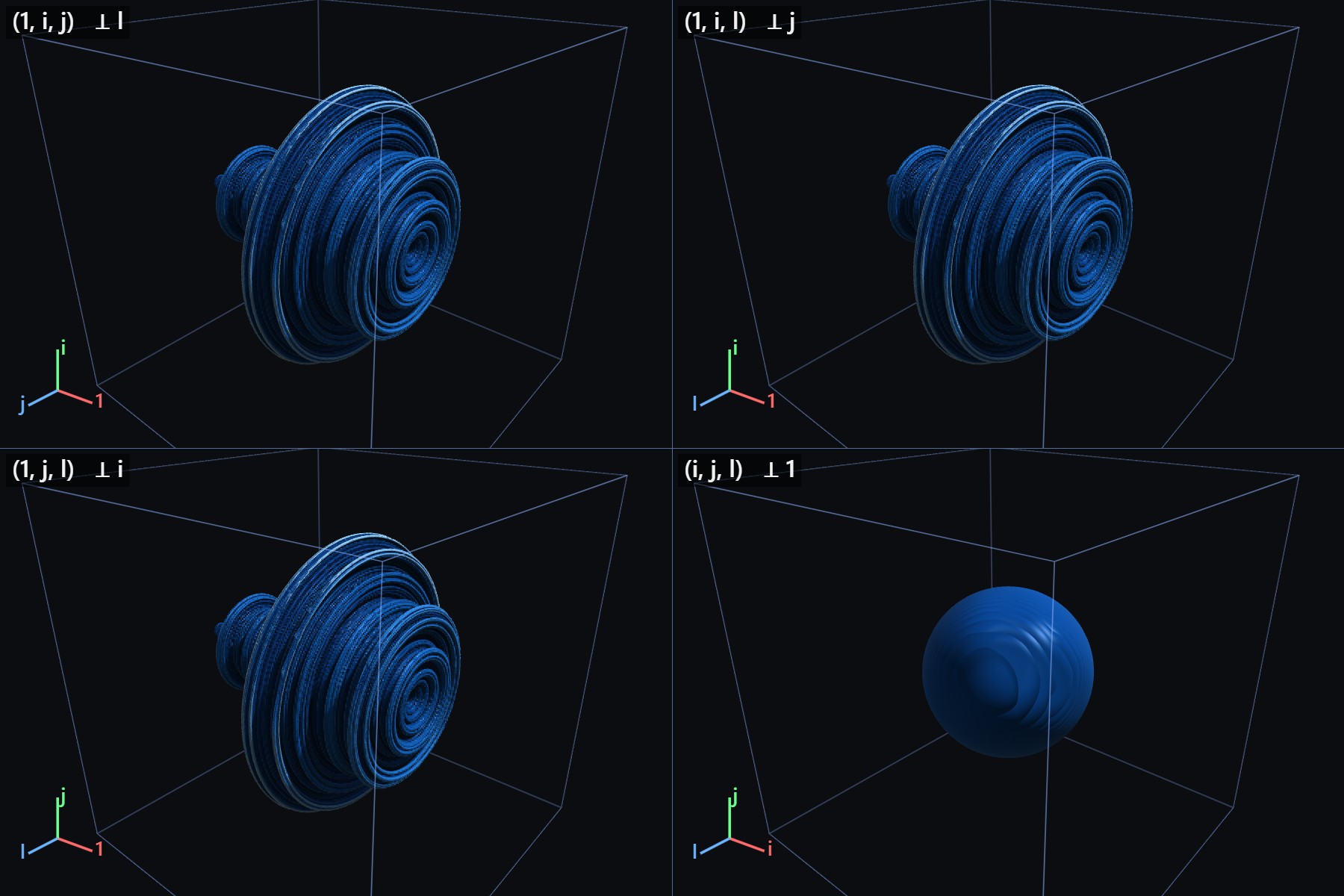}\par\vspace{1mm}
\noindent{\footnotesize \textbf{Map} $z\mapsto\overline{z\diamond z}+c$, $c$ the sample point $\cdot$ \textbf{Algebra} Para-octonions (para-Hurwitz $\bar x\bar y$, $\gamma=(-1,-1,-1)$) $\cdot$ \textbf{Frame} $X,Y,Z,N=1,i,j,l$ $\cdot$ center $(-0.35,\,\allowbreak 0,\,\allowbreak 0,\,\allowbreak 0,\,\allowbreak 0,\,\allowbreak 0,\,\allowbreak 0,\,\allowbreak 0)$ $\cdot$ half-extent 1.5 $\cdot$ 16 iterations $\cdot$ \textbf{Coloring} 8D orbit trap (point), gradient Deep Blue}\par\vspace{1mm}
\noindent{\footnotesize\begin{tabular}{@{}lllll@{}}\toprule Slice (X, Y, Z) & Body [95\%] & Exact ($J\geq0.9995$) & Approximate ($0.97\leq J<0.9995$) & Dim\\\midrule
($1, i, j$) $\perp l$ & 15.0\% {\scriptsize[14.8, 15.2]} & $M$: yz; $R_x$: all & -- & 2 \\
($1, i, l$) $\perp j$ & 15.0\% {\scriptsize[14.8, 15.2]} & $M$: yz; $R_x$: all & -- & 2 \\
($1, j, l$) $\perp i$ & 15.0\% {\scriptsize[14.8, 15.2]} & $M$: yz; $R_x$: all & -- & 2 \\
($i, j, l$) $\perp 1$ & 10.3\% {\scriptsize[10.1, 10.5]} & $M$: xyz; $R_x$: all; $R_y$: all; $R_z$: all; $P$ & -- & 2 \\
\bottomrule\end{tabular}}\par\vspace{2mm}
Conversely, the para-octonionic Mandelbar map $\overline{z\bullet z}+c=z^2+c$ is the octonionic Mandelbrot map, and the plate reproduces Plate~F01: the Mandelbrot set of revolution about the axis of $1$. The three quarters through $1$ show the same ribbed body, the main cardioid turned about the real axis with the bulbs stacked along it, and the imaginary quarter, which is the hyperplane $\re c=-0.35$, is a ball.

The orbit of the seed $0$ stays in the plane of $1$ and $c$ (Proposition~\ref{atl:prop:span}), so every section through the axis of $1$ contains the same complex Mandelbrot set, and the sections orthogonal to that axis are balls. The table agrees with Plate~F01 within the sampling error: the rotations about $1$ by every tested angle are exact in the quarters through $1$ ($15.0\%$), the imaginary quarter is an exact ball ($10.3\%$), and the orbits lie in a subalgebra of dimension $2$. Plates~G01 and G02 are two readings of the single identity $x\bullet y=\bar x\,\bar y$: conjugating the square exchanges the Mandelbrot and the Mandelbar maps of the two algebras.
\clearpage

\subsection{Part H. Cubic maps and higher powers}\label{atl:part:H}
Cubic maps are the first place where the symmetric composition algebras degenerate completely: $(z*z)*z=\n(z)\,z$. The comparison with the octonionic cubic, which does not degenerate, is immediate. The constants of this part are small; the fourth power computed by left bracketing, $((z*z)*z)*z$, is shown as well.

\subsubsection{H01. Octonions: a genuine cubic}\label{atl:pl:H01}\mbox{}\par\vspace{0.5mm}
\noindent\includegraphics[width=\textwidth]{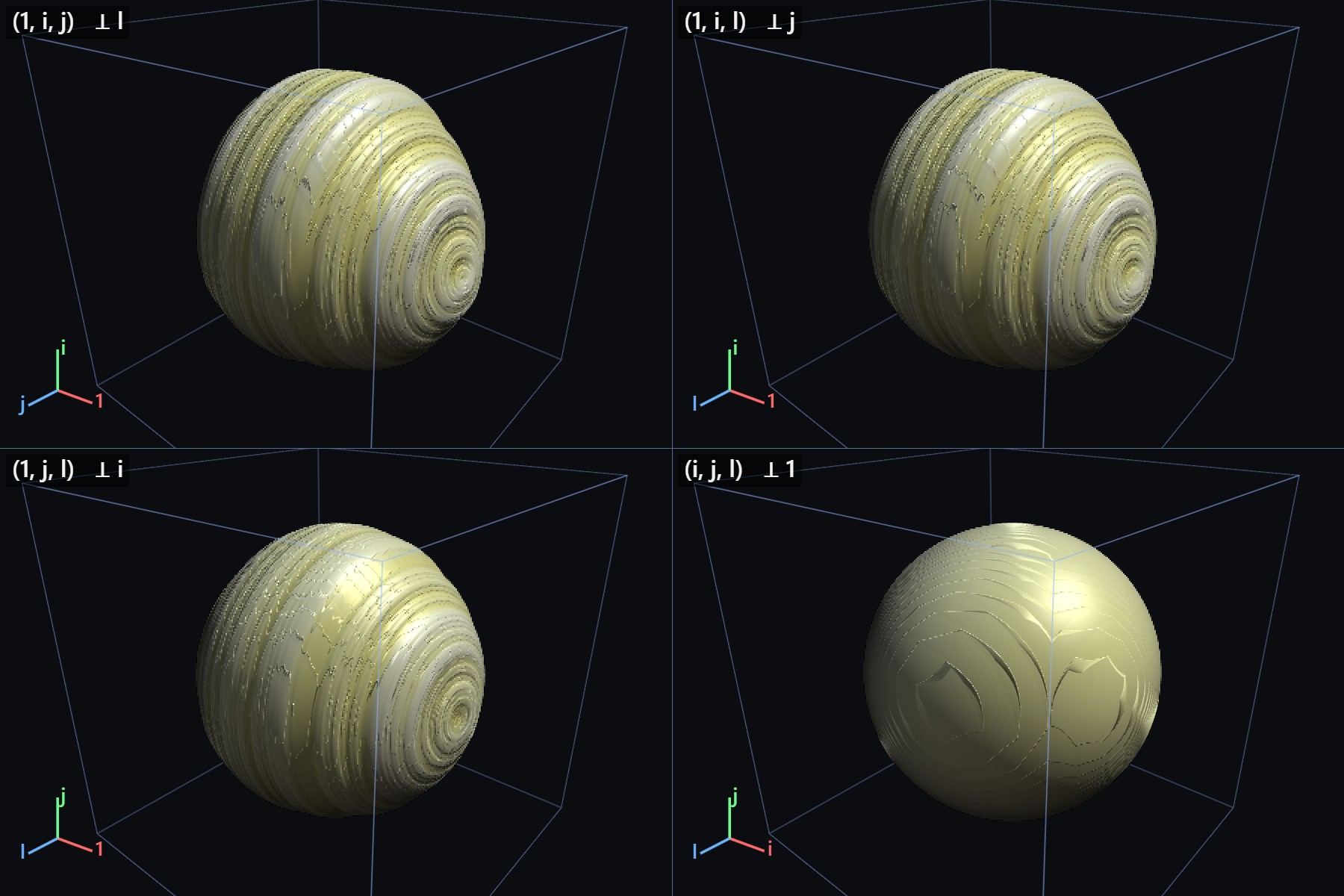}\par\vspace{1mm}
\noindent{\footnotesize \textbf{Map} $z\mapsto(z\diamond z)\diamond z+c$ $\cdot$ \textbf{Algebra} Octonions (Cayley--Dickson, $\gamma=(-1,-1,-1)$) $\cdot$ $c=(0.045,\,\allowbreak 0.203,\,\allowbreak 0.09,\,\allowbreak 0,\,\allowbreak 0.068,\,\allowbreak 0,\,\allowbreak 0,\,\allowbreak 0)$ $\cdot$ \textbf{Frame} $X,Y,Z,N=1,i,j,l$ $\cdot$ half-extent 1.3 $\cdot$ 16 iterations $\cdot$ \textbf{Coloring} 8D orbit trap (point), gradient Gold}\par\vspace{1mm}
\noindent{\footnotesize\begin{tabular}{@{}lllll@{}}\toprule Slice (X, Y, Z) & Body [95\%] & Exact ($J\geq0.9995$) & Approximate ($0.97\leq J<0.9995$) & Dim\\\midrule
($1, i, j$) $\perp l$ & 50.5\% {\scriptsize[50.2, 50.8]} & -- & -- & 4 \\
($1, i, l$) $\perp j$ & 50.5\% {\scriptsize[50.2, 50.8]} & -- & -- & 4 \\
($1, j, l$) $\perp i$ & 49.9\% {\scriptsize[49.5, 50.2]} & -- & -- & 4 \\
($i, j, l$) $\perp 1$ & 52.1\% {\scriptsize[51.8, 52.5]} & -- & $R_x$: 37 (min.\ 0.974) & 4 \\
\bottomrule\end{tabular}}\par\vspace{2mm}
The octonionic cubic $(z^2)z+c=z^3+c$ is a genuine cubic map, since the octonions are power-associative. The three quarters through $1$ show similar egg-shaped bodies, ribbed and with a spiral tip, with body shares of $50$--$51\%$; the imaginary quarter is a smoother, nearly round body ($52\%$).

The cube $z^3$ is a real polynomial in $\re z$, $\n(z)$ and $z$, so that, exactly as for the square (Proposition~\ref{atl:prop:equivariance}), the map commutes with every orthogonal map fixing $1$ and $c$. The constant $c=(0.045,0.2025,0.09,0,0.0675,0,0,0)$ has components along $i$, $j$ and $l$, and none of these symmetries is a coordinate mirror or rotation: the table finds no exact symmetry, and, the map being odd, not even the point reflection. In the imaginary quarter the symmetries fixing $1$ and $c$ include all the rotations about the projection of $\im c$, an axis tilted by about $29^\circ$ from $i$: the quarter is an exact body of revolution about this oblique axis, which the table sees only as an approximate rotation about $i$ ($J=0.974$). The subalgebra generated by $z_0$ and $c$ is a quaternion algebra, of dimension $4$.
\clearpage

\subsubsection{H02. Para-octonions: the cubic degenerates}\label{atl:pl:H02}\mbox{}\par\vspace{0.5mm}
\noindent\includegraphics[width=\textwidth]{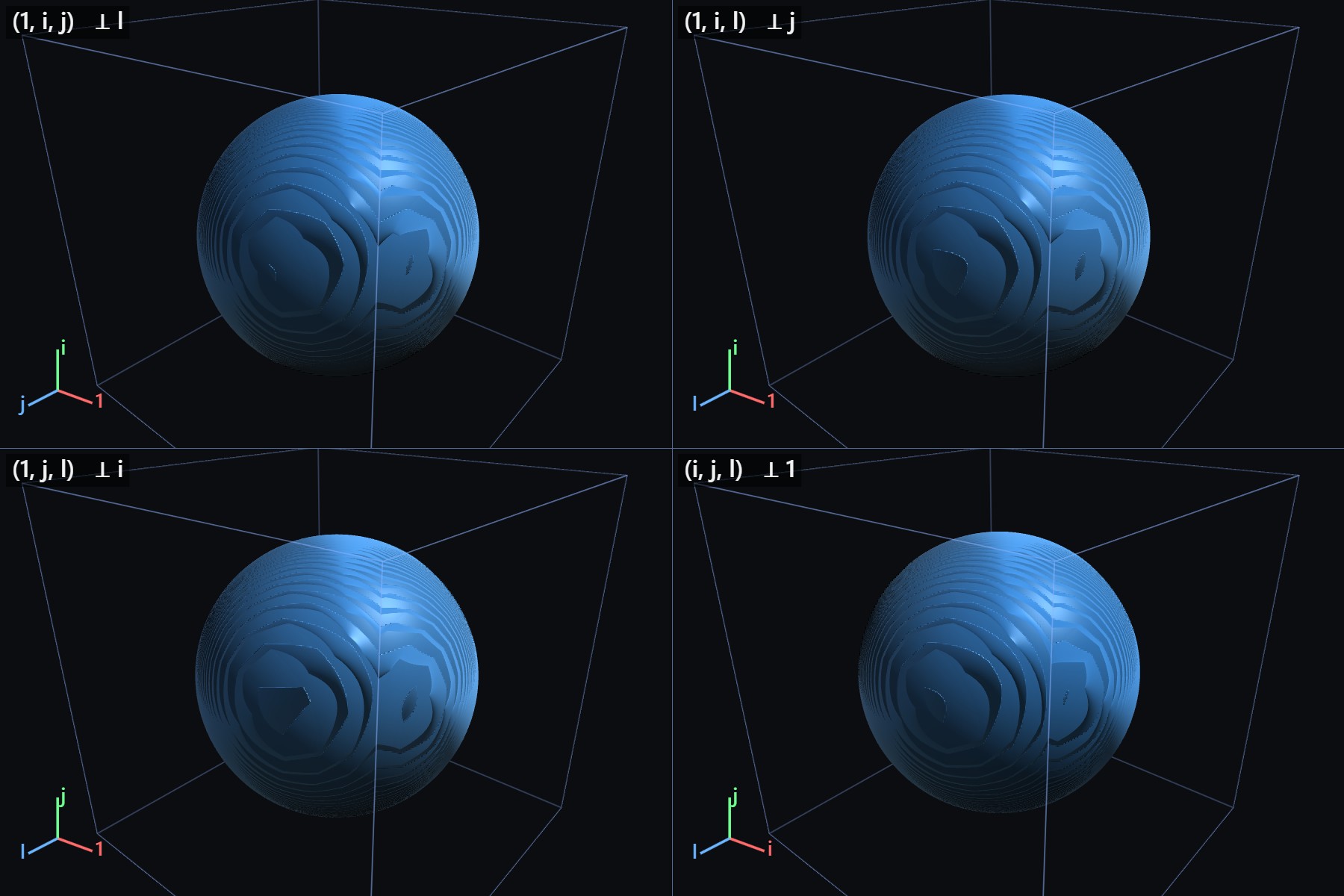}\par\vspace{1mm}
\noindent{\footnotesize \textbf{Map} $z\mapsto(z\diamond z)\diamond z+c$ $\cdot$ \textbf{Algebra} Para-octonions (para-Hurwitz $\bar x\bar y$, $\gamma=(-1,-1,-1)$) $\cdot$ $c=(0.045,\,\allowbreak 0.203,\,\allowbreak 0.09,\,\allowbreak 0,\,\allowbreak 0.068,\,\allowbreak 0,\,\allowbreak 0,\,\allowbreak 0)$ $\cdot$ \textbf{Frame} $X,Y,Z,N=1,i,j,l$ $\cdot$ half-extent 1.3 $\cdot$ 16 iterations $\cdot$ \textbf{Coloring} 8D orbit trap (point), gradient Deep Blue}\par\vspace{1mm}
\noindent{\footnotesize\begin{tabular}{@{}lllll@{}}\toprule Slice (X, Y, Z) & Body [95\%] & Exact ($J\geq0.9995$) & Approximate ($0.97\leq J<0.9995$) & Dim\\\midrule
($1, i, j$) $\perp l$ & 49.1\% {\scriptsize[48.8, 49.4]} & -- & -- & 4 \\
($1, i, l$) $\perp j$ & 49.0\% {\scriptsize[48.7, 49.3]} & -- & -- & 4 \\
($1, j, l$) $\perp i$ & 48.2\% {\scriptsize[47.9, 48.6]} & -- & -- & 4 \\
($i, j, l$) $\perp 1$ & 49.1\% {\scriptsize[48.7, 49.4]} & -- & -- & 4 \\
\bottomrule\end{tabular}}\par\vspace{2mm}
In the para-octonions the cubic degenerates: $(z\bullet z)\bullet z=\n(z)\,z$ by \eqref{atl:eq:symmetric}, and the map is $z\mapsto\n(z)\,z+c$. The four quarters show four nearly round bodies, with nearly identical trap contours and body shares of $48$--$49\%$.

The degenerate map commutes with every orthogonal transformation fixing $c$, a group isomorphic to $O(7)$, so the eight-dimensional body is a body of revolution about the axis of $c$. Each quarter is in turn a body of revolution about the projection of $c$ onto it (a numerical check with $37^\circ$ rotations about these axes gives $J=1$). Since $c$ is oblique to the coordinate axes, so are these axes, and the table, which tests only coordinate mirrors and rotations, reports no symmetry at all: the bodies are not balls, the largest score being $0.96$. The map is odd, and the point reflection is not a symmetry. The subalgebra generated by $z_0$ and $c$ has dimension $4$, but the orbits themselves stay in the plane of $z_0$ and $c$, since $\n(z)\,z+c$ lies in that plane.
\clearpage

\subsubsection{H03. Okubo: the same degenerate cubic}\label{atl:pl:H03}\mbox{}\par\vspace{0.5mm}
\noindent\includegraphics[width=\textwidth]{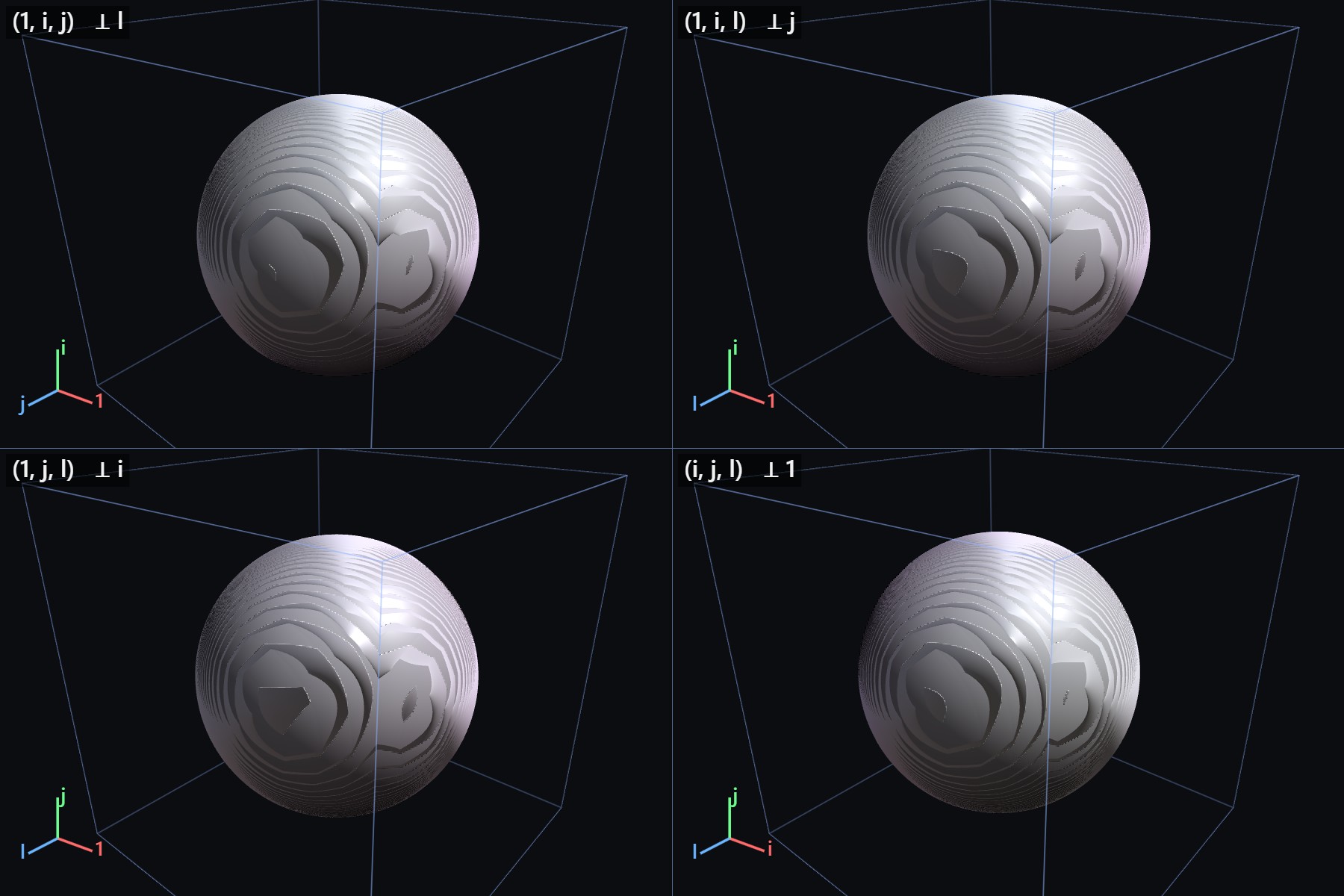}\par\vspace{1mm}
\noindent{\footnotesize \textbf{Map} $z\mapsto(z\diamond z)\diamond z+c$ $\cdot$ \textbf{Algebra} Okubo algebra (Petersson, $\gamma=(-1,-1,-1)$, $\theta=120^\circ$) $\cdot$ $c=(0.045,\,\allowbreak 0.203,\,\allowbreak 0.09,\,\allowbreak 0,\,\allowbreak 0.068,\,\allowbreak 0,\,\allowbreak 0,\,\allowbreak 0)$ $\cdot$ \textbf{Frame} $X,Y,Z,N=1,i,j,l$ $\cdot$ half-extent 1.3 $\cdot$ 16 iterations $\cdot$ \textbf{Coloring} 8D orbit trap (point), gradient Aurora}\par\vspace{1mm}
\noindent{\footnotesize\begin{tabular}{@{}lllll@{}}\toprule Slice (X, Y, Z) & Body [95\%] & Exact ($J\geq0.9995$) & Approximate ($0.97\leq J<0.9995$) & Dim\\\midrule
($1, i, j$) $\perp l$ & 49.5\% {\scriptsize[49.2, 49.8]} & -- & -- & 8 \\
($1, i, l$) $\perp j$ & 49.5\% {\scriptsize[49.2, 49.8]} & -- & -- & 8 \\
($1, j, l$) $\perp i$ & 48.7\% {\scriptsize[48.3, 49.0]} & -- & -- & 8 \\
($i, j, l$) $\perp 1$ & 49.3\% {\scriptsize[49.0, 49.6]} & -- & -- & 8 \\
\bottomrule\end{tabular}}\par\vspace{2mm}
The Okubonic cubic is the same degenerate map $\n(z)\,z+c$ (Proposition~\ref{atl:prop:symmetric}), and the plate is indistinguishable from Plate~H02; the body shares agree within the sampling error ($48.7$--$49.5\%$ against $48.2$--$49.1\%$), and on a common sample of points the two maps give identical bodies. Whatever the product, the cubic of a symmetric composition algebra forgets it, and only the norm survives.

The only visible trace of the Okubo product is in the table: the subalgebra generated by $z_0$ and $c$ has dimension $8$ instead of $4$. It contains the orbits, but it is much larger than them, since every orbit stays in the plane of $z_0$ and $c$. As in Plate~H02, each quarter is a body of revolution about the oblique projection of $c$, which the coordinate tests do not detect.
\clearpage

\subsubsection{H04. Split Okubo: a hyperbolic degenerate cubic}\label{atl:pl:H04}\mbox{}\par\vspace{0.5mm}
\noindent\includegraphics[width=\textwidth]{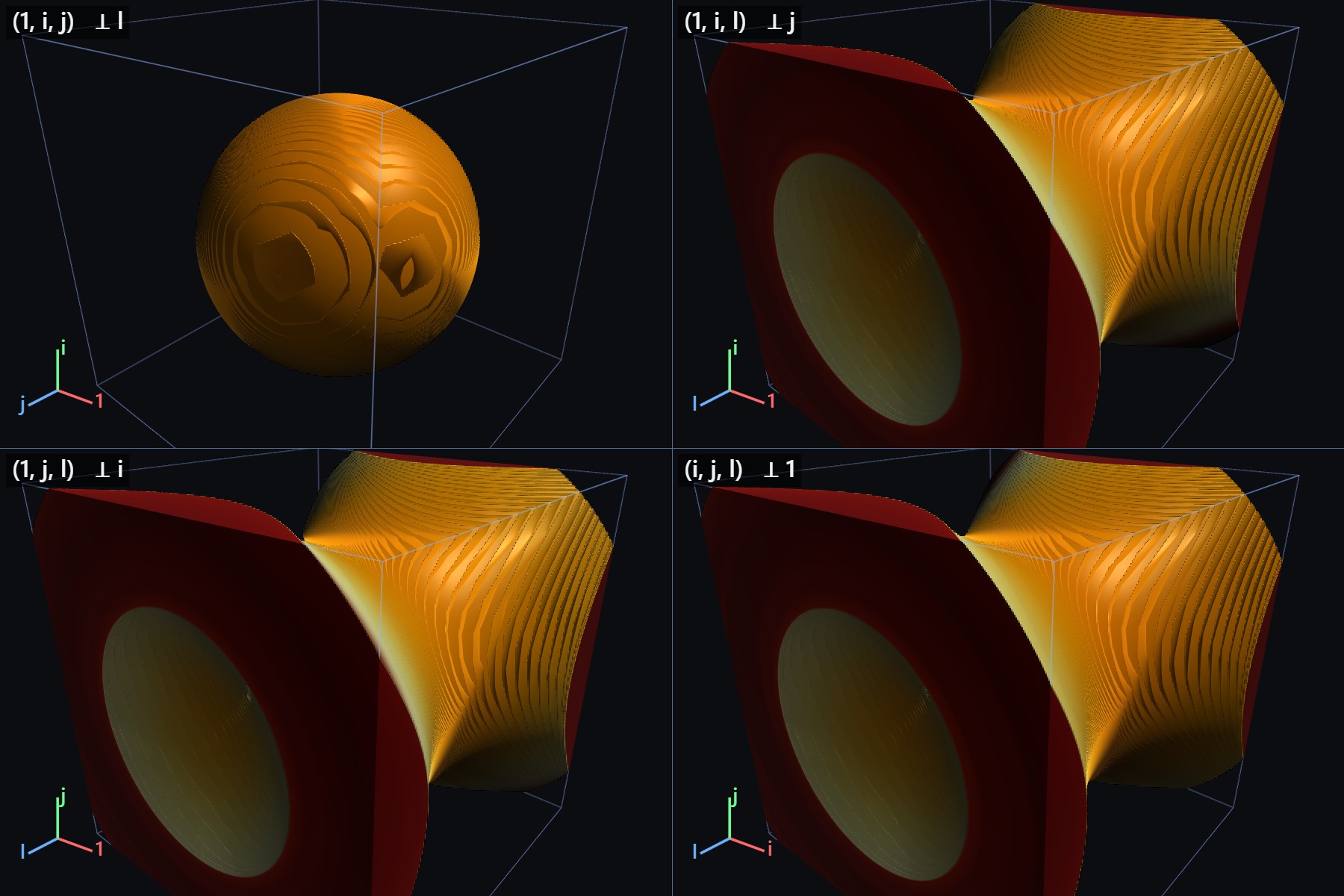}\par\vspace{1mm}
\noindent{\footnotesize \textbf{Map} $z\mapsto(z\diamond z)\diamond z+c$ $\cdot$ \textbf{Algebra} Split Okubo algebra (Petersson, $\gamma=(-1,-1,1)$, $\theta=120^\circ$) $\cdot$ $c=(0.045,\,\allowbreak 0.203,\,\allowbreak 0.09,\,\allowbreak 0,\,\allowbreak 0.068,\,\allowbreak 0,\,\allowbreak 0,\,\allowbreak 0)$ $\cdot$ \textbf{Frame} $X,Y,Z,N=1,i,j,l$ $\cdot$ half-extent 1.3 $\cdot$ 16 iterations $\cdot$ \textbf{Coloring} 8D orbit trap (point), gradient Fire}\par\vspace{1mm}
\noindent{\footnotesize\begin{tabular}{@{}lllll@{}}\toprule Slice (X, Y, Z) & Body [95\%] & Exact ($J\geq0.9995$) & Approximate ($0.97\leq J<0.9995$) & Dim\\\midrule
($1, i, j$) $\perp l$ & 50.3\% {\scriptsize[50.0, 50.6]} & -- & -- & 8 \\
($1, i, l$) $\perp j$ & 81.3\% {\scriptsize[81.0, 81.5]} & -- & $M$: z (min.\ 0.973) & 8 \\
($1, j, l$) $\perp i$ & 81.4\% {\scriptsize[81.1, 81.6]} & -- & $M$: z; $R_z$: 37 (min.\ 0.972) & 8 \\
($i, j, l$) $\perp 1$ & 81.3\% {\scriptsize[81.1, 81.6]} & -- & $M$: z (min.\ 0.974) & 8 \\
\bottomrule\end{tabular}}\par\vspace{2mm}
In the split Okubo algebra the degenerate cubic $\n(z)\,z+c$ involves an indefinite norm. The quaternionic quarter, where the norm is positive definite, is again a nearly round body ($50\%$), close to the one of Plate~H02. The three quarters containing $l$ look identical: bodies that fill most of the cube ($81\%$), cut by its faces into flat cross-sections, with a large concave bowl on the side of $l$ and saddle-shaped surfaces on the opposite side, the marks of the level sets of an indefinite norm.

The table reports no exact symmetry, only approximate mirrors in the $l$-direction ($J\approx0.97$) and, in the quarter $(1,j,l)$, an approximate rotation about $l$ ($J=0.977$). The map commutes with every linear map that preserves the split norm and fixes $c$, but in the quarters containing $l$ these are hyperbolic rotations, not Euclidean ones. The map depends on the product only through the norm, so the same picture would be obtained in the split para-octonions. The rendering uses a Euclidean bailout, which is not an escape criterion for an indefinite norm: the bodies of this plate are approximations.
\clearpage

\subsubsection{H05. Split-octonions: a genuine hyperbolic cubic}\label{atl:pl:H05}\mbox{}\par\vspace{0.5mm}
\noindent\includegraphics[width=\textwidth]{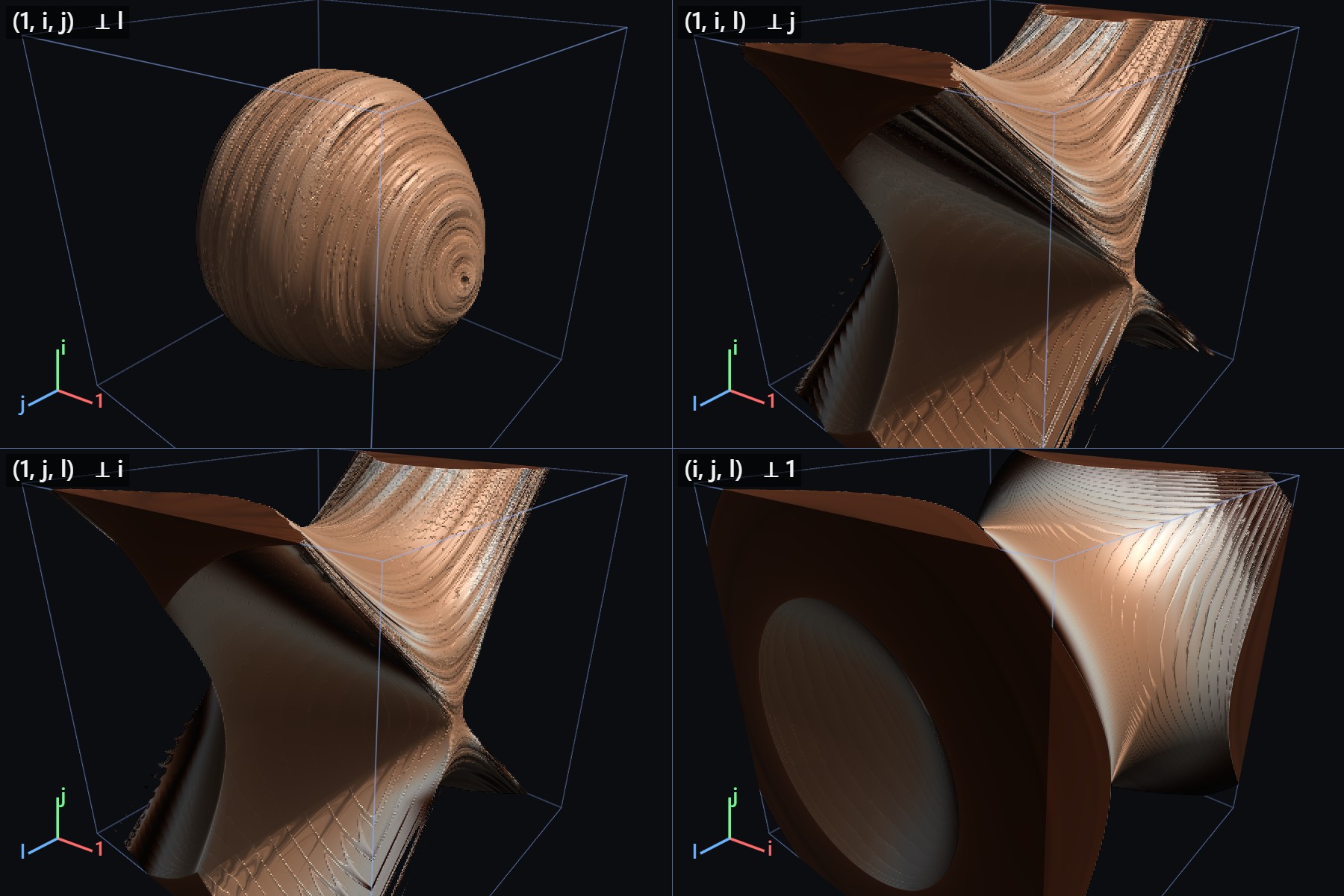}\par\vspace{1mm}
\noindent{\footnotesize \textbf{Map} $z\mapsto(z\diamond z)\diamond z+c$ $\cdot$ \textbf{Algebra} Split-octonions ($l^2=+1$) (Cayley--Dickson, $\gamma=(-1,-1,1)$) $\cdot$ $c=(0.045,\,\allowbreak 0.203,\,\allowbreak 0.09,\,\allowbreak 0,\,\allowbreak 0.068,\,\allowbreak 0,\,\allowbreak 0,\,\allowbreak 0)$ $\cdot$ \textbf{Frame} $X,Y,Z,N=1,i,j,l$ $\cdot$ half-extent 1.3 $\cdot$ 16 iterations $\cdot$ \textbf{Coloring} 8D orbit trap (point), gradient Copper}\par\vspace{1mm}
\noindent{\footnotesize\begin{tabular}{@{}lllll@{}}\toprule Slice (X, Y, Z) & Body [95\%] & Exact ($J\geq0.9995$) & Approximate ($0.97\leq J<0.9995$) & Dim\\\midrule
($1, i, j$) $\perp l$ & 51.1\% {\scriptsize[50.8, 51.4]} & -- & -- & 4 \\
($1, i, l$) $\perp j$ & 59.6\% {\scriptsize[59.3, 59.9]} & -- & -- & 4 \\
($1, j, l$) $\perp i$ & 60.2\% {\scriptsize[59.9, 60.5]} & -- & -- & 4 \\
($i, j, l$) $\perp 1$ & 82.8\% {\scriptsize[82.6, 83.1]} & -- & $M$: z; $R_z$: 37 (min.\ 0.974) & 4 \\
\bottomrule\end{tabular}}\par\vspace{2mm}
The split-octonionic cubic is a genuine cubic map, not degenerate. The quaternionic quarter is a ribbed egg with a spiral tip, as in Plate~H01, while the quarters $(1,i,l)$ and $(1,j,l)$ show large hyperbolic saddles and sheets cut by the faces of the cube; the imaginary quarter is a body that fills most of the cube ($83\%$), with a concave bowl on the side of $l$, very much like the corresponding quarter of Plate~H04.

The table reports no exact symmetry; in the imaginary quarter the mirror in $l$ and the rotation about $l$ are approximate ($J\approx0.975$). As in the octonions, the cube is a real polynomial in $\re z$, $\n(z)$ and $z$, but the maps that preserve the split norm and fix $1$ and $c$ are hyperbolic rotations in the directions of $l$, and they do not produce Euclidean revolutions. The comparison with Plate~H04 separates the effect of the split norm from the effect of the degeneration.
\clearpage

\subsubsection{H06. Okubo: the third power \texorpdfstring{$z^3$}{z\textthreesuperior{}}}\label{atl:pl:H06}\mbox{}\par\vspace{0.5mm}
\noindent\includegraphics[width=\textwidth]{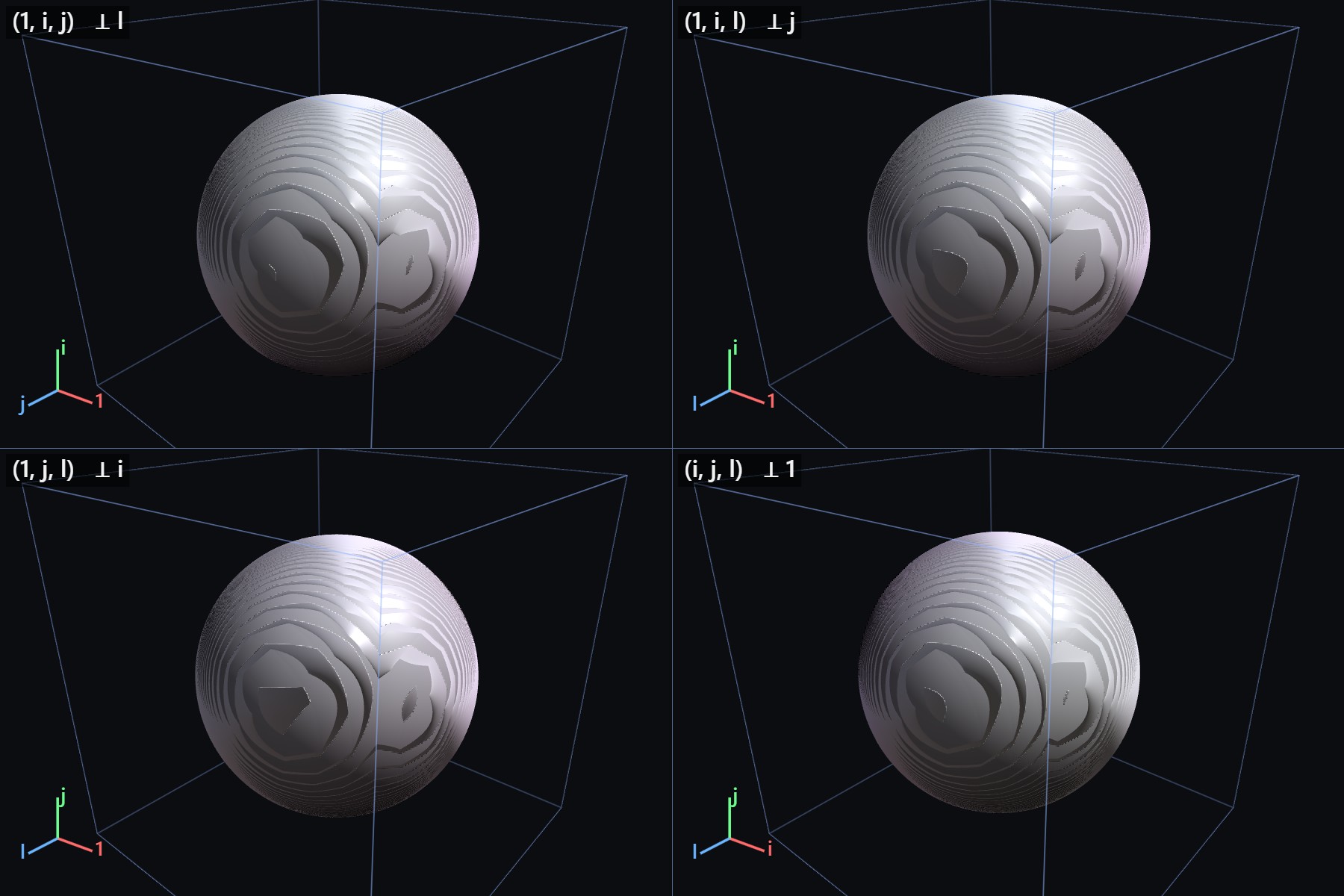}\par\vspace{1mm}
\noindent{\footnotesize \textbf{Map} $z\mapsto z^{p}+c$ (\texttt{Standard.ufm}) $\cdot$ \textbf{Algebra} Okubo algebra (Petersson, $\gamma=(-1,-1,-1)$, $\theta=120^\circ$) $\cdot$ $seed=(0.045,\,\allowbreak 0.203,\,\allowbreak 0.09,\,\allowbreak 0,\,\allowbreak 0.068,\,\allowbreak 0,\,\allowbreak 0,\,\allowbreak 0)$ $\cdot$ $p=3$ $\cdot$ \textbf{Frame} $X,Y,Z,N=1,i,j,l$ $\cdot$ half-extent 1.3 $\cdot$ 16 iterations $\cdot$ \textbf{Coloring} 8D orbit trap (point), gradient Aurora}\par\vspace{1mm}
\noindent{\footnotesize\begin{tabular}{@{}lllll@{}}\toprule Slice (X, Y, Z) & Body [95\%] & Exact ($J\geq0.9995$) & Approximate ($0.97\leq J<0.9995$) & Dim\\\midrule
($1, i, j$) $\perp l$ & 49.0\% {\scriptsize[48.7, 49.3]} & -- & -- & 8 \\
($1, i, l$) $\perp j$ & 48.9\% {\scriptsize[48.6, 49.2]} & -- & -- & 8 \\
($1, j, l$) $\perp i$ & 48.2\% {\scriptsize[47.9, 48.5]} & -- & -- & 8 \\
($i, j, l$) $\perp 1$ & 49.1\% {\scriptsize[48.8, 49.4]} & -- & -- & 8 \\
\bottomrule\end{tabular}}\par\vspace{2mm}
The third power computed by the standard formula $z^p+c$ with $p=3$ uses the left bracketing $(z*z)*z$ in the non-unital algebras, and gives exactly the degenerate cubic of Plate~H03: the plate is identical, ball-like body by ball-like body, and the body shares agree with those of Plates~H02 and H03 within the sampling error. In the Okubo algebra there is no other natural cube: by \eqref{atl:eq:symmetric}, $z*(z*z)=(z*z)*z=\n(z)\,z$.

As in Plates~H02 and H03, the table reports no symmetry: every quarter is a body of revolution about the projection of $c$, which is oblique to the coordinate axes, and the odd map has no point symmetry.
\clearpage

\subsubsection{H07. Octonions: the fourth power}\label{atl:pl:H07}\mbox{}\par\vspace{0.5mm}
\noindent\includegraphics[width=\textwidth]{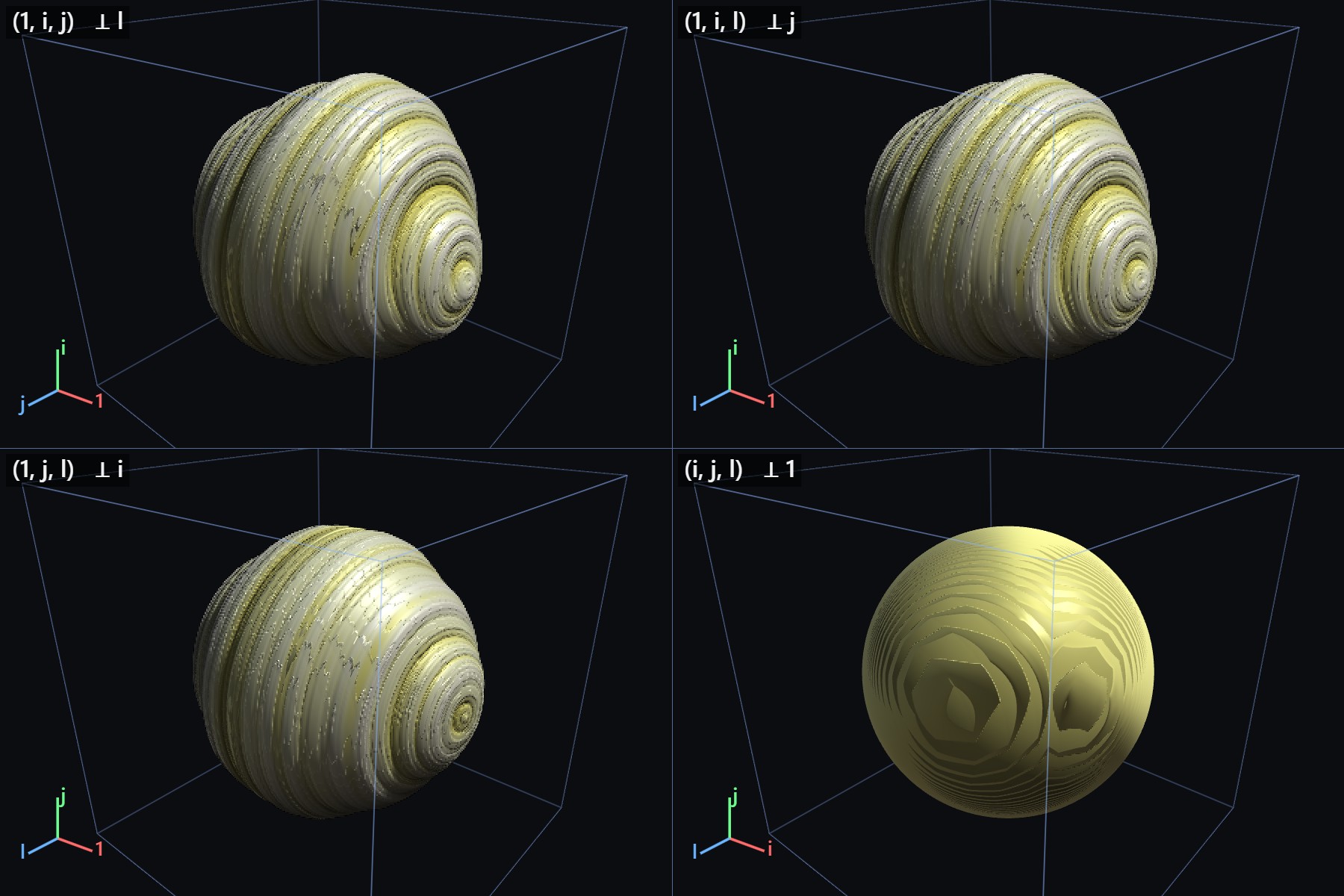}\par\vspace{1mm}
\noindent{\footnotesize \textbf{Map} $z\mapsto z^{p}+c$ (\texttt{Standard.ufm}) $\cdot$ \textbf{Algebra} Octonions (Cayley--Dickson, $\gamma=(-1,-1,-1)$) $\cdot$ $seed=(-0.1,\,\allowbreak 0.3,\,\allowbreak 0.05,\,\allowbreak 0,\,\allowbreak 0.1,\,\allowbreak 0,\,\allowbreak 0,\,\allowbreak 0)$ $\cdot$ $p=4$ $\cdot$ \textbf{Frame} $X,Y,Z,N=1,i,j,l$ $\cdot$ half-extent 1.3 $\cdot$ 16 iterations $\cdot$ \textbf{Coloring} 8D orbit trap (point), gradient Gold}\par\vspace{1mm}
\noindent{\footnotesize\begin{tabular}{@{}lllll@{}}\toprule Slice (X, Y, Z) & Body [95\%] & Exact ($J\geq0.9995$) & Approximate ($0.97\leq J<0.9995$) & Dim\\\midrule
($1, i, j$) $\perp l$ & 50.2\% {\scriptsize[49.9, 50.5]} & $P$ & -- & 4 \\
($1, i, l$) $\perp j$ & 50.2\% {\scriptsize[49.9, 50.5]} & $P$ & -- & 4 \\
($1, j, l$) $\perp i$ & 49.8\% {\scriptsize[49.5, 50.1]} & $P$ & $R_x$: 37 (min.\ 0.977) & 4 \\
($i, j, l$) $\perp 1$ & 50.6\% {\scriptsize[50.3, 51.0]} & $M$: xyz; $R_x$: all; $R_y$: all; $R_z$: all; $P$ & -- & 4 \\
\bottomrule\end{tabular}}\par\vspace{2mm}
The fourth power in the octonions is a genuine quartic map. The three quarters through $1$ are ribbed, rounded bodies whose ribs wind about a pole, very similar to each other ($50\%$); the imaginary quarter is a ball. For purely imaginary $z$ one has $z^2=-\n(z)$ and $z^4=\n(z)^2$, a real number, so that $z_1=\n(z_0)^2+c$ depends only on the norm, and the table confirms that the quarter is exactly invariant under every rotation.

The quartic is even, and the point reflection is exact in every quarter (Remark~\ref{atl:rem:minus}). As for the square, $z^4$ is a real polynomial in $\re z$, $\n(z)$ and $z$, and the map commutes with every orthogonal map fixing $1$ and $c=(-0.1,0.3,0.05,0,0.1,0,0,0)$. Since the components of $c$ along $j$ and $l$ are small, the rotation about $1$ in the quarter $(1,j,l)$ nearly fixes $c$, and the table reports it as approximate ($J=0.977$), not exact. The subalgebra generated by $z_0$ and $c$ has dimension $4$.
\clearpage

\subsubsection{H08. Okubo: the fourth power by left bracketing}\label{atl:pl:H08}\mbox{}\par\vspace{0.5mm}
\noindent\includegraphics[width=\textwidth]{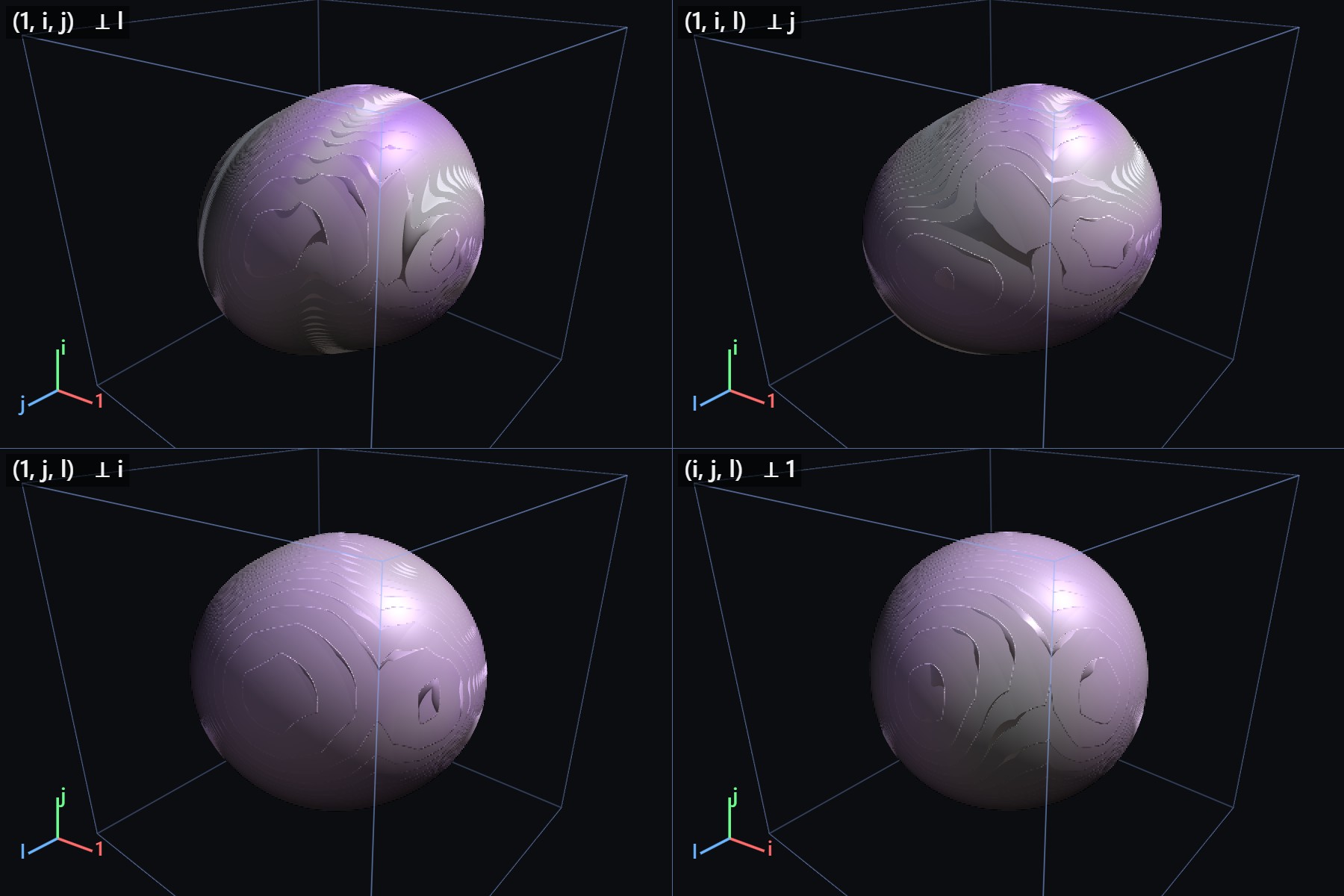}\par\vspace{1mm}
\noindent{\footnotesize \textbf{Map} $z\mapsto z^{p}+c$ (\texttt{Standard.ufm}) $\cdot$ \textbf{Algebra} Okubo algebra (Petersson, $\gamma=(-1,-1,-1)$, $\theta=120^\circ$) $\cdot$ $seed=(-0.1,\,\allowbreak 0.3,\,\allowbreak 0.05,\,\allowbreak 0,\,\allowbreak 0.1,\,\allowbreak 0,\,\allowbreak 0,\,\allowbreak 0)$ $\cdot$ $p=4$ $\cdot$ \textbf{Frame} $X,Y,Z,N=1,i,j,l$ $\cdot$ half-extent 1.3 $\cdot$ 16 iterations $\cdot$ \textbf{Coloring} 8D orbit trap (point), gradient Aurora}\par\vspace{1mm}
\noindent{\footnotesize\begin{tabular}{@{}lllll@{}}\toprule Slice (X, Y, Z) & Body [95\%] & Exact ($J\geq0.9995$) & Approximate ($0.97\leq J<0.9995$) & Dim\\\midrule
($1, i, j$) $\perp l$ & 47.9\% {\scriptsize[47.6, 48.2]} & $P$ & -- & 8 \\
($1, i, l$) $\perp j$ & 46.7\% {\scriptsize[46.4, 47.0]} & $P$ & -- & 8 \\
($1, j, l$) $\perp i$ & 47.7\% {\scriptsize[47.4, 48.0]} & $P$ & $R_x$: 37 (min.\ 0.976) & 8 \\
($i, j, l$) $\perp 1$ & 43.8\% {\scriptsize[43.5, 44.1]} & $P$ & $M$: y; $R_y$: 180 (min.\ 0.995) & 8 \\
\bottomrule\end{tabular}}\par\vspace{2mm}
In the Okubo algebra the fourth power computed by left bracketing is $((z*z)*z)*z=\n(z)\,(z*z)$, by \eqref{atl:eq:symmetric}: a quadratic map rescaled by the norm. The bodies are smooth eggs with trap contours, somewhat different from quarter to quarter, with body shares of $44$--$48\%$. This is not a degeneration into a function of the norm but a hybrid: the quartic of the Okubo algebra is the Okubonic square weighted by the norm.

The map is even, and the point reflection is exact in every quarter (Remark~\ref{atl:rem:minus}); the table finds no other exact symmetry, and only approximate ones: the rotation about $1$ in the quarter $(1,j,l)$ ($J=0.976$) and the mirror in $i$ of the imaginary quarter ($J=0.995$). The constant has components along $j$ and $l$, so that it does not lie in $\CC$ and the torus of Proposition~\ref{atl:prop:torus} does not fix it. The subalgebra generated by $z_0$ and $c$ has dimension $8$.
\clearpage

\subsection{Part I. The sandwich map}\label{atl:part:I}
The sandwich map $(z\diamond a)\diamond z+c$ is well defined in flexible algebras. In a symmetric composition algebra it collapses to $\n(z)\,a+c$, and its dynamics reduces to a real quadratic recursion on the norm (Proposition~\ref{atl:prop:symmetric}). The constants are chosen so that this recursion has a bounded interval of initial norms.

\subsubsection{I01. Octonions: the sandwich is genuine}\label{atl:pl:I01}\mbox{}\par\vspace{0.5mm}
\noindent\includegraphics[width=\textwidth]{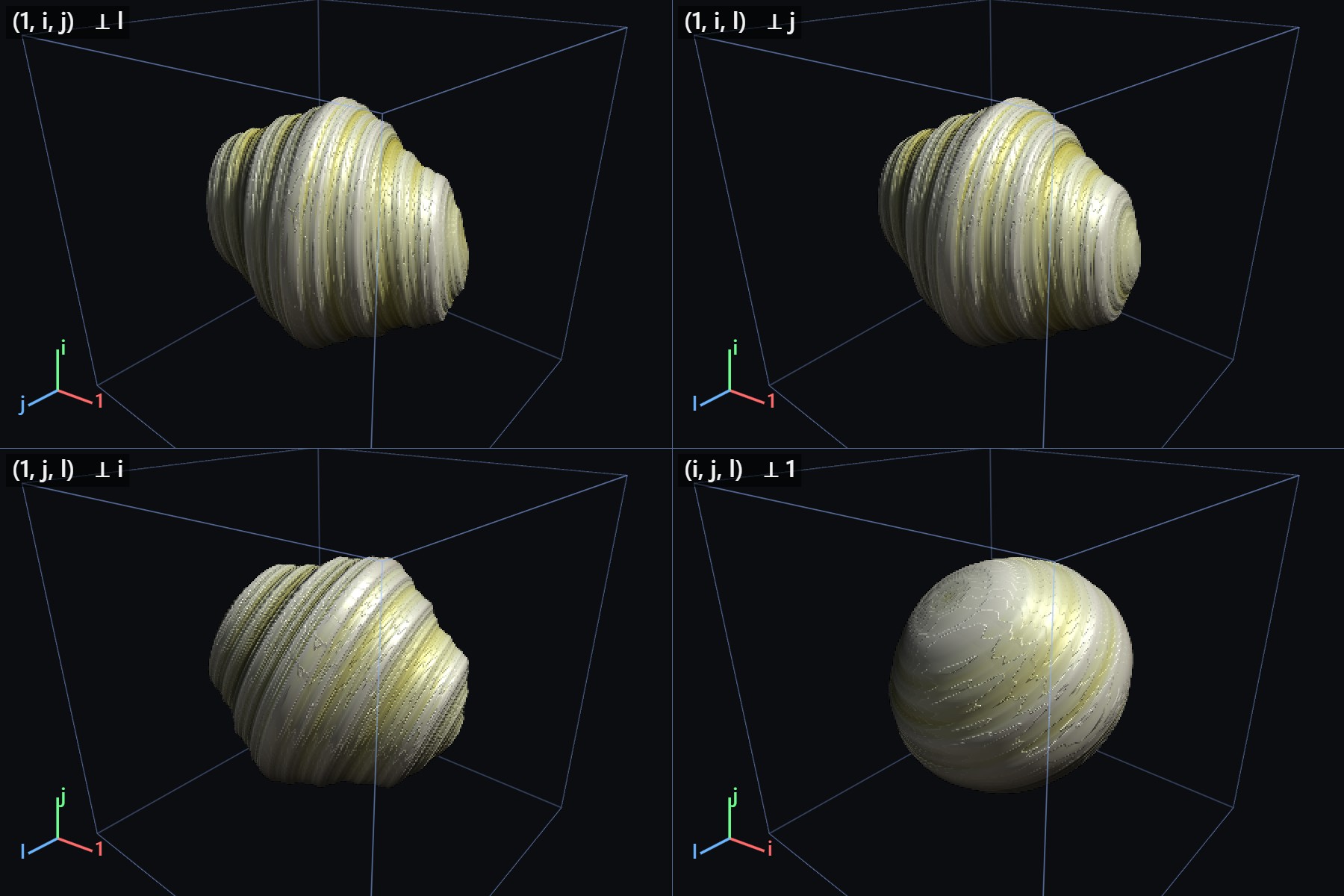}\par\vspace{1mm}
\noindent{\footnotesize \textbf{Map} $z\mapsto(z\diamond a)\diamond z+c$ $\cdot$ \textbf{Algebra} Octonions (Cayley--Dickson, $\gamma=(-1,-1,-1)$) $\cdot$ $c=(-0.3,\,\allowbreak 0.2,\,\allowbreak 0,\,\allowbreak 0.1,\,\allowbreak 0,\,\allowbreak 0.1,\,\allowbreak 0,\,\allowbreak 0)$ $\cdot$ $a=(0.8,\,\allowbreak 0,\,\allowbreak 0.3,\,\allowbreak 0,\,\allowbreak 0.2,\,\allowbreak 0,\,\allowbreak 0.3,\,\allowbreak 0)$ $\cdot$ \textbf{Frame} $X,Y,Z,N=1,i,j,l$ $\cdot$ half-extent 1.6 $\cdot$ 16 iterations $\cdot$ \textbf{Coloring} 8D orbit trap (point), gradient Gold}\par\vspace{1mm}
\noindent{\footnotesize\begin{tabular}{@{}lllll@{}}\toprule Slice (X, Y, Z) & Body [95\%] & Exact ($J\geq0.9995$) & Approximate ($0.97\leq J<0.9995$) & Dim\\\midrule
($1, i, j$) $\perp l$ & 28.9\% {\scriptsize[28.7, 29.2]} & $P$ & -- & 8 \\
($1, i, l$) $\perp j$ & 29.0\% {\scriptsize[28.8, 29.3]} & $P$ & -- & 8 \\
($1, j, l$) $\perp i$ & 29.0\% {\scriptsize[28.7, 29.2]} & $P$ & -- & 8 \\
($i, j, l$) $\perp 1$ & 28.1\% {\scriptsize[27.8, 28.4]} & $P$ & -- & 8 \\
\bottomrule\end{tabular}}\par\vspace{2mm}
In the octonions the sandwich map $(z\,a)\,z+c$ is a genuine quadratic-type map, well defined because the octonions are flexible. The three quarters through $1$ are ribbed, irregular bodies, and the imaginary quarter is a striped, nearly round body; the body shares are close to $29\%$ in every quarter.

The map is even, and the point reflection is exact in every quarter (Remark~\ref{atl:rem:minus}); since the constants $a$ and $c$ are generic, there is no other symmetry, and even the imaginary quarter is not a ball (the best score is $0.96$). The subalgebra generated by $z_0$, $a$ and $c$ is the whole algebra, of dimension $8$: unlike the quadratic map, the sandwich uses the non-associativity of the octonions.
\clearpage

\subsubsection{I02. Okubo: the sandwich is a ball}\label{atl:pl:I02}\mbox{}\par\vspace{0.5mm}
\noindent\includegraphics[width=\textwidth]{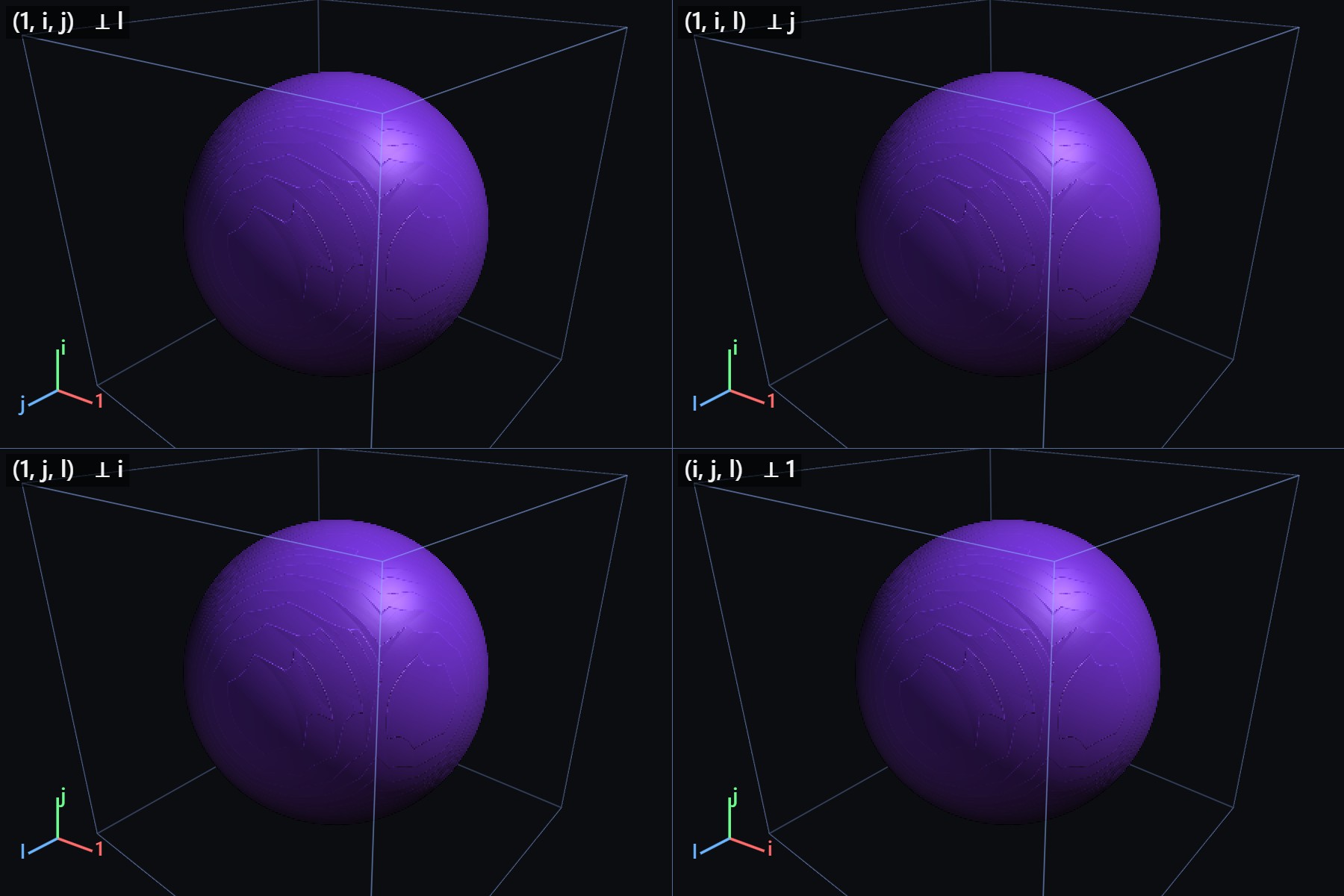}\par\vspace{1mm}
\noindent{\footnotesize \textbf{Map} $z\mapsto(z\diamond a)\diamond z+c$ $\cdot$ \textbf{Algebra} Okubo algebra (Petersson, $\gamma=(-1,-1,-1)$, $\theta=120^\circ$) $\cdot$ $c=(-0.3,\,\allowbreak 0.2,\,\allowbreak 0,\,\allowbreak 0.1,\,\allowbreak 0,\,\allowbreak 0.1,\,\allowbreak 0,\,\allowbreak 0)$ $\cdot$ $a=(0.8,\,\allowbreak 0,\,\allowbreak 0.3,\,\allowbreak 0,\,\allowbreak 0.2,\,\allowbreak 0,\,\allowbreak 0.3,\,\allowbreak 0)$ $\cdot$ \textbf{Frame} $X,Y,Z,N=1,i,j,l$ $\cdot$ half-extent 1.6 $\cdot$ 16 iterations $\cdot$ \textbf{Coloring} 8D orbit trap (point), gradient Aurora}\par\vspace{1mm}
\noindent{\footnotesize\begin{tabular}{@{}lllll@{}}\toprule Slice (X, Y, Z) & Body [95\%] & Exact ($J\geq0.9995$) & Approximate ($0.97\leq J<0.9995$) & Dim\\\midrule
($1, i, j$) $\perp l$ & 58.4\% {\scriptsize[58.1, 58.7]} & $M$: xyz; $R_x$: all; $R_y$: all; $R_z$: all; $P$ & -- & 8 \\
($1, i, l$) $\perp j$ & 58.4\% {\scriptsize[58.1, 58.7]} & $M$: xyz; $R_x$: all; $R_y$: all; $R_z$: all; $P$ & -- & 8 \\
($1, j, l$) $\perp i$ & 58.4\% {\scriptsize[58.1, 58.7]} & $M$: xyz; $R_x$: all; $R_y$: all; $R_z$: all; $P$ & -- & 8 \\
($i, j, l$) $\perp 1$ & 58.4\% {\scriptsize[58.1, 58.7]} & $M$: xyz; $R_x$: all; $R_y$: all; $R_z$: all; $P$ & -- & 8 \\
\bottomrule\end{tabular}}\par\vspace{2mm}
In the Okubo algebra the sandwich degenerates into $\n(z)\,a+c$ (Proposition~\ref{atl:prop:symmetric}), and its filled set is the region where the real recursion \eqref{atl:eq:sandwich} of the norms stays bounded. With the constants of this plate, $\n(a)=0.86$, $\n(c)=0.15$ and $\n(a,c)=-0.48$, the effective parameter is $c_{\mathrm{eff}}=-0.169$, and the bounded interval of initial norms is $0\le\n(z_0)\le1.613$: the filled set is the ball of radius $1.270$. The four quarters are four identical smooth balls, invariant under every tested rotation and mirror, with a body share of $58.4\%$ of the measuring ball of radius $0.95\cdot1.6$, which corresponds to a radius of $1.270$, in agreement with the prediction.

Of all plates of the atlas, this is perhaps the one where the algebra speaks most clearly: a map that in the octonions produces the complicated bodies of Plate~I01 is, in the Okubo algebra, a ball. The subalgebra generated by $z_0$, $a$ and $c$ has dimension $8$, but the orbit is trivial: from $z_1$ on it stays on the line $\{t\,a+c\}$.
\clearpage

\subsubsection{I03. Para-octonions: the same ball}\label{atl:pl:I03}\mbox{}\par\vspace{0.5mm}
\noindent\includegraphics[width=\textwidth]{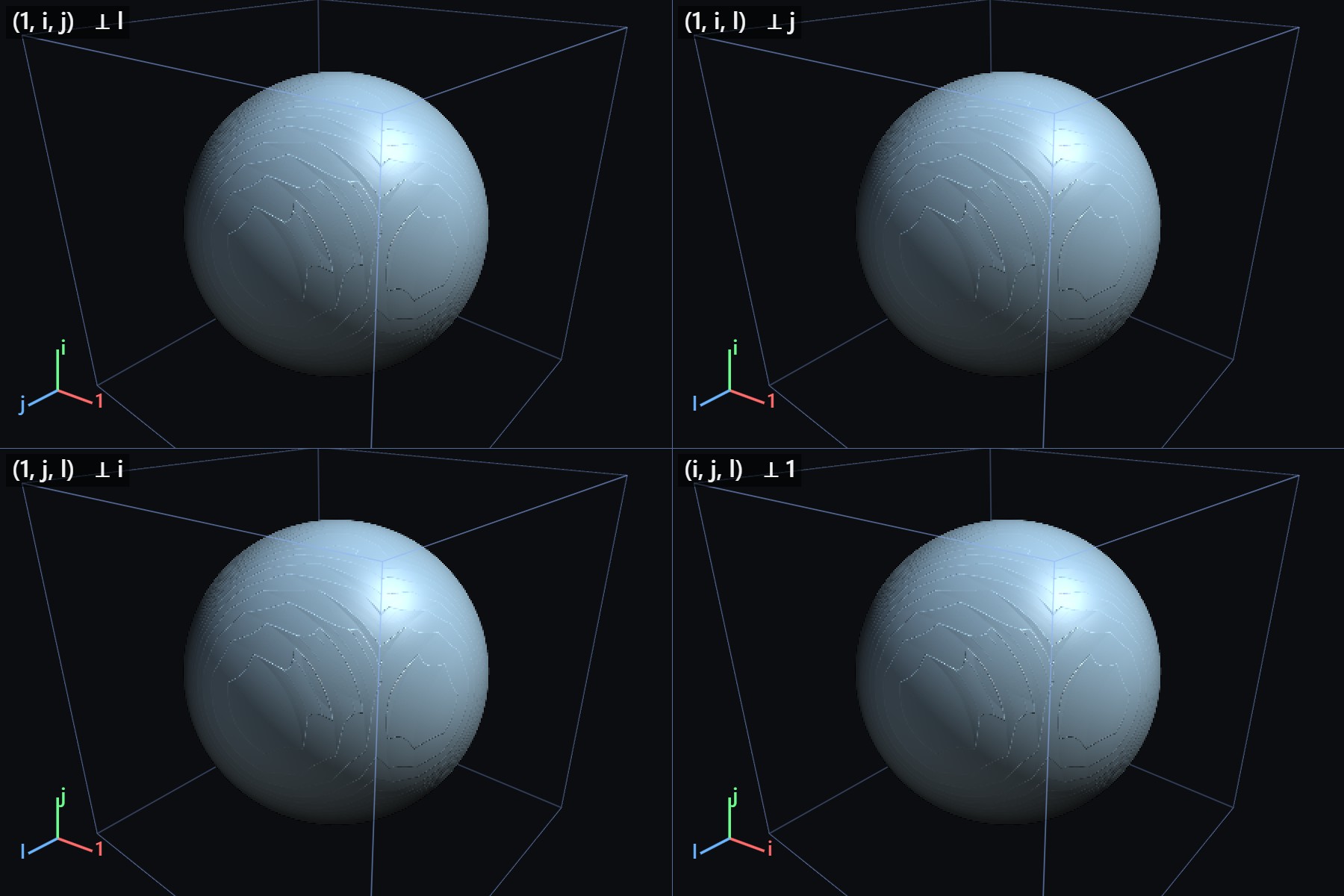}\par\vspace{1mm}
\noindent{\footnotesize \textbf{Map} $z\mapsto(z\diamond a)\diamond z+c$ $\cdot$ \textbf{Algebra} Para-octonions (para-Hurwitz $\bar x\bar y$, $\gamma=(-1,-1,-1)$) $\cdot$ $c=(-0.3,\,\allowbreak 0.2,\,\allowbreak 0,\,\allowbreak 0.1,\,\allowbreak 0,\,\allowbreak 0.1,\,\allowbreak 0,\,\allowbreak 0)$ $\cdot$ $a=(0.8,\,\allowbreak 0,\,\allowbreak 0.3,\,\allowbreak 0,\,\allowbreak 0.2,\,\allowbreak 0,\,\allowbreak 0.3,\,\allowbreak 0)$ $\cdot$ \textbf{Frame} $X,Y,Z,N=1,i,j,l$ $\cdot$ half-extent 1.6 $\cdot$ 16 iterations $\cdot$ \textbf{Coloring} 8D orbit trap (point), gradient Deep Blue}\par\vspace{1mm}
\noindent{\footnotesize\begin{tabular}{@{}lllll@{}}\toprule Slice (X, Y, Z) & Body [95\%] & Exact ($J\geq0.9995$) & Approximate ($0.97\leq J<0.9995$) & Dim\\\midrule
($1, i, j$) $\perp l$ & 58.0\% {\scriptsize[57.7, 58.3]} & $M$: xyz; $R_x$: all; $R_y$: all; $R_z$: all; $P$ & -- & 8 \\
($1, i, l$) $\perp j$ & 58.0\% {\scriptsize[57.7, 58.3]} & $M$: xyz; $R_x$: all; $R_y$: all; $R_z$: all; $P$ & -- & 8 \\
($1, j, l$) $\perp i$ & 58.0\% {\scriptsize[57.7, 58.3]} & $M$: xyz; $R_x$: all; $R_y$: all; $R_z$: all; $P$ & -- & 8 \\
($i, j, l$) $\perp 1$ & 58.0\% {\scriptsize[57.7, 58.3]} & $M$: xyz; $R_x$: all; $R_y$: all; $R_z$: all; $P$ & -- & 8 \\
\bottomrule\end{tabular}}\par\vspace{2mm}
The para-octonions give the same balls as the Okubo algebra, as Proposition~\ref{atl:prop:symmetric} predicts: the degenerate sandwich $\n(z)\,a+c$ depends only on the norm, which is the same in both algebras. The four quarters are identical balls, exactly invariant under every tested rotation and mirror; the body share, $58.0\%$, agrees with that of Plate~I02 ($58.4\%$) within the sampling error, and corresponds to the predicted radius $1.27$.
\clearpage

\subsubsection{I04. Split-octonions: a genuine hyperbolic sandwich}\label{atl:pl:I04}\mbox{}\par\vspace{0.5mm}
\noindent\includegraphics[width=\textwidth]{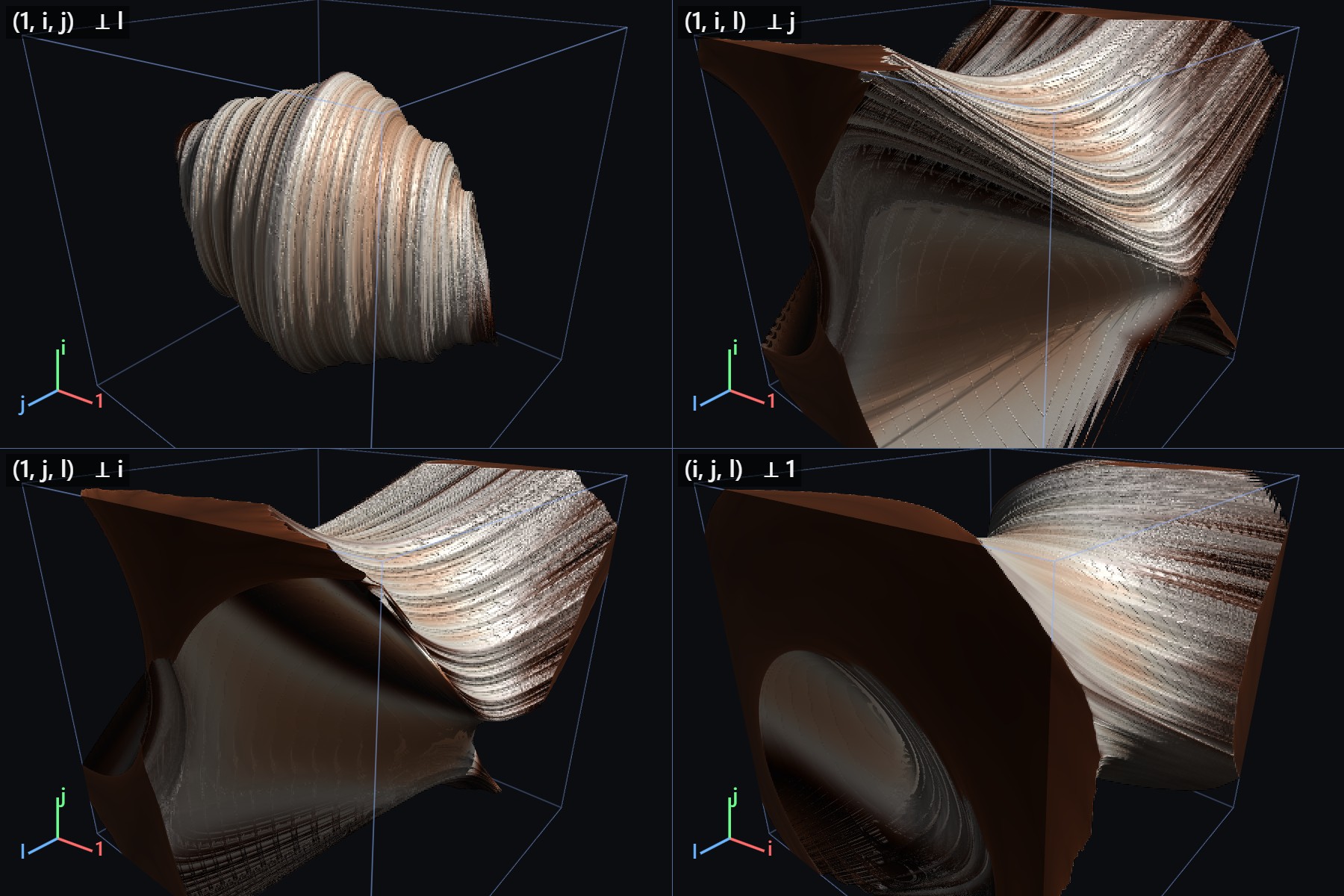}\par\vspace{1mm}
\noindent{\footnotesize \textbf{Map} $z\mapsto(z\diamond a)\diamond z+c$ $\cdot$ \textbf{Algebra} Split-octonions ($l^2=+1$) (Cayley--Dickson, $\gamma=(-1,-1,1)$) $\cdot$ $c=(-0.3,\,\allowbreak 0.2,\,\allowbreak 0,\,\allowbreak 0.1,\,\allowbreak 0,\,\allowbreak 0.1,\,\allowbreak 0,\,\allowbreak 0)$ $\cdot$ $a=(0.8,\,\allowbreak 0,\,\allowbreak 0.3,\,\allowbreak 0,\,\allowbreak 0.2,\,\allowbreak 0,\,\allowbreak 0.3,\,\allowbreak 0)$ $\cdot$ \textbf{Frame} $X,Y,Z,N=1,i,j,l$ $\cdot$ half-extent 1.6 $\cdot$ 16 iterations $\cdot$ \textbf{Coloring} 8D orbit trap (point), gradient Copper}\par\vspace{1mm}
\noindent{\footnotesize\begin{tabular}{@{}lllll@{}}\toprule Slice (X, Y, Z) & Body [95\%] & Exact ($J\geq0.9995$) & Approximate ($0.97\leq J<0.9995$) & Dim\\\midrule
($1, i, j$) $\perp l$ & 48.0\% {\scriptsize[47.7, 48.3]} & $P$ & -- & 8 \\
($1, i, l$) $\perp j$ & 69.8\% {\scriptsize[69.5, 70.1]} & $P$ & -- & 8 \\
($1, j, l$) $\perp i$ & 66.2\% {\scriptsize[65.9, 66.5]} & $P$ & -- & 8 \\
($i, j, l$) $\perp 1$ & 77.8\% {\scriptsize[77.5, 78.0]} & $P$ & -- & 8 \\
\bottomrule\end{tabular}}\par\vspace{2mm}
The split-octonionic sandwich is a genuine map, not degenerate, since the split-octonions are unital. The quaternionic quarter is a ribbed irregular body, as in Plate~I01 ($48\%$), while the quarters containing $l$ show large hyperbolic sheets and saddles, cut by the faces of the cube, that fill most of the box ($66$--$78\%$). The map is even, and the point reflection is the only exact symmetry (Remark~\ref{atl:rem:minus}); the orbits generate the whole algebra (dimension $8$). The comparison with Plate~I05 is the split analogue of the comparison between Plates~I01 and I02. For an indefinite norm the Euclidean bailout of the renderer is not an escape criterion, and the bodies are approximations.
\clearpage

\subsubsection{I05. Split Okubo: a ball and three carved cubes}\label{atl:pl:I05}\mbox{}\par\vspace{0.5mm}
\noindent\includegraphics[width=\textwidth]{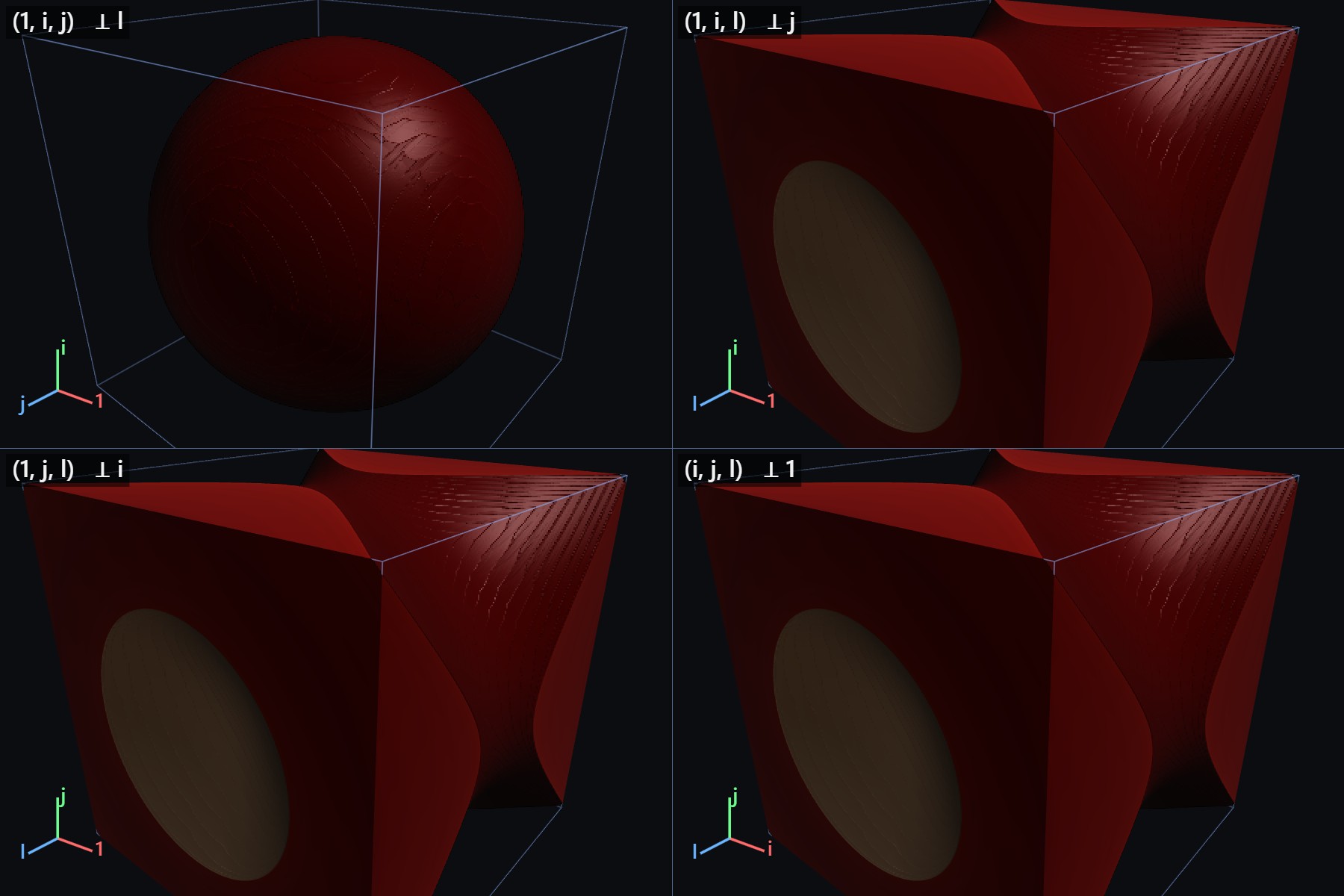}\par\vspace{1mm}
\noindent{\footnotesize \textbf{Map} $z\mapsto(z\diamond a)\diamond z+c$ $\cdot$ \textbf{Algebra} Split Okubo algebra (Petersson, $\gamma=(-1,-1,1)$, $\theta=120^\circ$) $\cdot$ $c=(-0.3,\,\allowbreak 0.2,\,\allowbreak 0,\,\allowbreak 0.1,\,\allowbreak 0,\,\allowbreak 0.1,\,\allowbreak 0,\,\allowbreak 0)$ $\cdot$ $a=(0.8,\,\allowbreak 0,\,\allowbreak 0.3,\,\allowbreak 0,\,\allowbreak 0.2,\,\allowbreak 0,\,\allowbreak 0.3,\,\allowbreak 0)$ $\cdot$ \textbf{Frame} $X,Y,Z,N=1,i,j,l$ $\cdot$ half-extent 1.6 $\cdot$ 16 iterations $\cdot$ \textbf{Coloring} 8D orbit trap (point), gradient Fire}\par\vspace{1mm}
\noindent{\footnotesize\begin{tabular}{@{}lllll@{}}\toprule Slice (X, Y, Z) & Body [95\%] & Exact ($J\geq0.9995$) & Approximate ($0.97\leq J<0.9995$) & Dim\\\midrule
($1, i, j$) $\perp l$ & 100.0\% {\scriptsize[100.0, 100.0]} & $M$: xyz; $R_x$: all; $R_y$: all; $R_z$: all; $P$ & -- & 8 \\
($1, i, l$) $\perp j$ & 98.0\% {\scriptsize[97.9, 98.1]} & $M$: xyz; $R_x$: 180; $R_y$: 180; $R_z$: all; $P$ & -- & 8 \\
($1, j, l$) $\perp i$ & 98.0\% {\scriptsize[97.9, 98.1]} & $M$: xyz; $R_x$: 180; $R_y$: 180; $R_z$: all; $P$ & -- & 8 \\
($i, j, l$) $\perp 1$ & 98.0\% {\scriptsize[97.9, 98.1]} & $M$: xyz; $R_x$: 180; $R_y$: 180; $R_z$: all; $P$ & -- & 8 \\
\bottomrule\end{tabular}}\par\vspace{2mm}
The split Okubo sandwich degenerates into $\n(z)\,a+c$ with the indefinite split norm (Proposition~\ref{atl:prop:symmetric}). Here \eqref{atl:eq:sandwich} reads $t\mapsto0.6\,t^2-0.48\,t+0.13$, and the filled set is the region $t_{\min}\le\n(z_0)\le t_{\max}$ with $t_{\min}=-1.575$ and $t_{\max}=2.375$. In the quaternionic quarter the norm is positive definite and the filled set is the ball of radius $\sqrt{2.375}=1.541$, nearly inscribed in the cube and larger than the measuring ball, which it fills entirely.

In the three quarters containing $l$ the norm has signature $(2,1)$, and the filled set is bounded by two hyperboloids: the two-sheeted $\n(z)=t_{\min}$, which carves the concave elliptic bowl seen on the side of $l$, and the one-sheeted $\n(z)=t_{\max}$, which shapes the saddle surfaces on the opposite side. The three quarters are therefore identical ``carved cubes'' ($98.0\%$ of the measuring ball), and each is an exact body of revolution about $l$, since the norm depends only on the distance from the $l$-axis and on the $l$-coordinate; the table confirms the rotations about $l$ by every tested angle. The Euclidean bailout of the renderer is harmless here, since the orbit after the first step lies on the line $\{t\,a+c\}$ and escapes together with $t$.
\clearpage

\subsubsection{I06. Split para-octonions: the same carved cubes}\label{atl:pl:I06}\mbox{}\par\vspace{0.5mm}
\noindent\includegraphics[width=\textwidth]{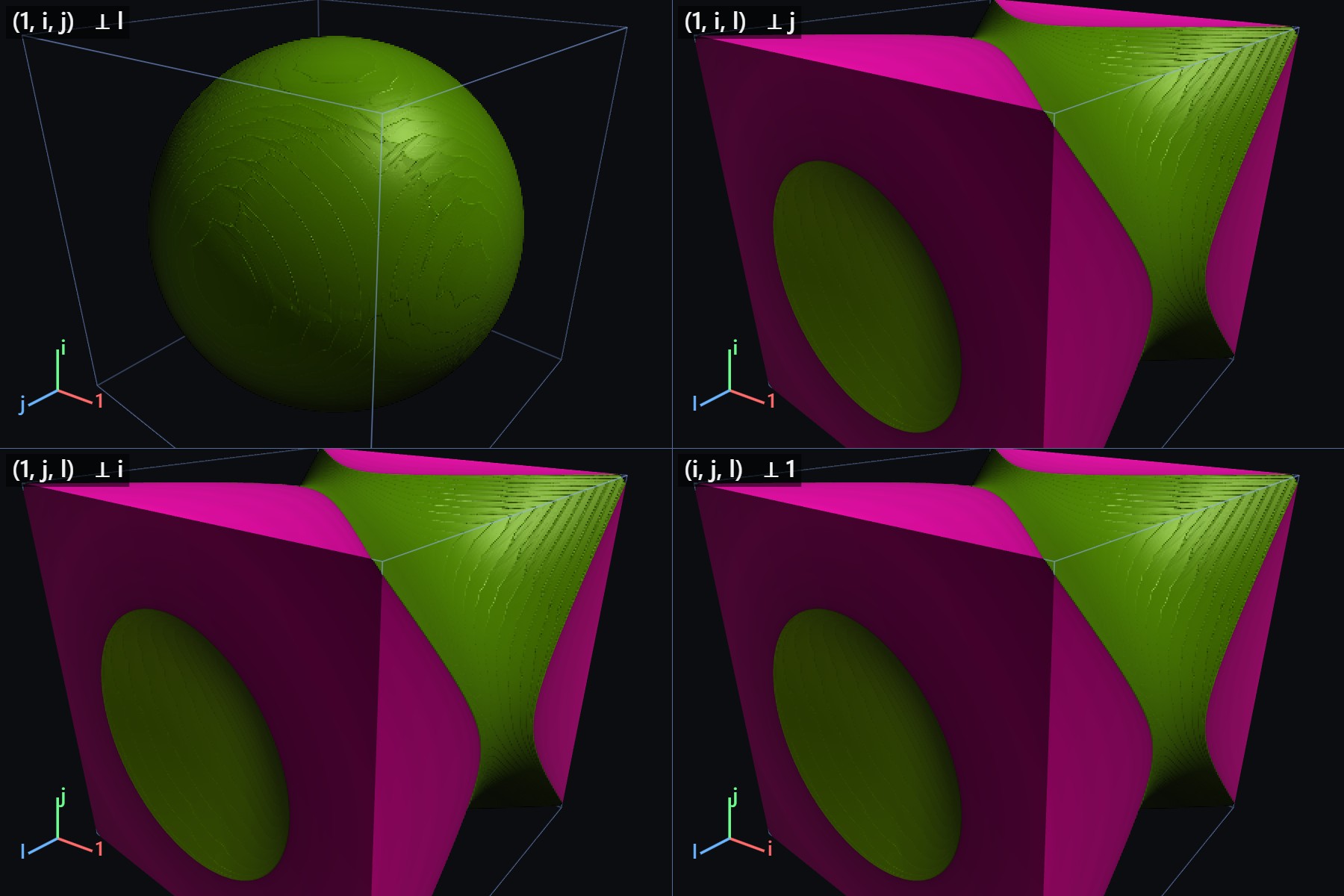}\par\vspace{1mm}
\noindent{\footnotesize \textbf{Map} $z\mapsto(z\diamond a)\diamond z+c$ $\cdot$ \textbf{Algebra} Split para-octonions (para-Hurwitz $\bar x\bar y$, $\gamma=(-1,-1,1)$) $\cdot$ $c=(-0.3,\,\allowbreak 0.2,\,\allowbreak 0,\,\allowbreak 0.1,\,\allowbreak 0,\,\allowbreak 0.1,\,\allowbreak 0,\,\allowbreak 0)$ $\cdot$ $a=(0.8,\,\allowbreak 0,\,\allowbreak 0.3,\,\allowbreak 0,\,\allowbreak 0.2,\,\allowbreak 0,\,\allowbreak 0.3,\,\allowbreak 0)$ $\cdot$ \textbf{Frame} $X,Y,Z,N=1,i,j,l$ $\cdot$ half-extent 1.6 $\cdot$ 16 iterations $\cdot$ \textbf{Coloring} 8D orbit trap (point), gradient Neon}\par\vspace{1mm}
\noindent{\footnotesize\begin{tabular}{@{}lllll@{}}\toprule Slice (X, Y, Z) & Body [95\%] & Exact ($J\geq0.9995$) & Approximate ($0.97\leq J<0.9995$) & Dim\\\midrule
($1, i, j$) $\perp l$ & 100.0\% {\scriptsize[100.0, 100.0]} & $M$: xyz; $R_x$: all; $R_y$: all; $R_z$: all; $P$ & -- & 8 \\
($1, i, l$) $\perp j$ & 97.9\% {\scriptsize[97.8, 97.9]} & $M$: xyz; $R_x$: 180; $R_y$: 180; $R_z$: all; $P$ & -- & 8 \\
($1, j, l$) $\perp i$ & 97.9\% {\scriptsize[97.8, 97.9]} & $M$: xyz; $R_x$: 180; $R_y$: 180; $R_z$: all; $P$ & -- & 8 \\
($i, j, l$) $\perp 1$ & 97.9\% {\scriptsize[97.8, 97.9]} & $M$: xyz; $R_x$: 180; $R_y$: 180; $R_z$: all; $P$ & -- & 8 \\
\bottomrule\end{tabular}}\par\vspace{2mm}
The split para-octonions give the same carved cubes as the split Okubo algebra, as predicted by Proposition~\ref{atl:prop:symmetric}: the degenerate sandwich depends only on the norm, and the norms of the two algebras coincide. The quaternionic quarter is the ball of radius $1.541$, which fills the measuring ball ($100\%$); the three quarters containing $l$ are identical ($97.9\%$) and exactly invariant under every rotation about $l$, the two-sheeted hyperboloid $\n(z)=t_{\min}$ removing only two small caps about the $l$-axis from the measuring ball.
\clearpage

\subsection{Part J. Bilateral maps and Kaplansky's trick}\label{atl:part:J}
The bilateral map $(a\diamond z)\diamond(z\diamond b)+c$ uses two different factors, so that non-associativity cannot be hidden. With $a=b=1$ Kaplansky's trick turns the Okubo and the para-octonionic maps into the octonionic quadratic map. Quite unexpectedly, the same holds for generic factors, found by the parameter scan of HF8D~0.2: in a symmetric composition algebra $(a\diamond z)\diamond(z\diamond b)=\n(z,e)\,z-\n(z)\,e$ with $e=b\diamond a$, so that every bilateral map of the para-octonions and of the Okubo algebra is an octonionic quadratic map moved by a similarity, whose axis of revolution is the tilted direction of $e$ (Proposition~\ref{atl:prop:kaplansky}). The octonionic bilateral map, on the contrary, is not of this form.

\subsubsection{J01. Kaplansky: the Okubo bilateral map with \texorpdfstring{$a=b=1$}{a=b=1}}\label{atl:pl:J01}\mbox{}\par\vspace{0.5mm}
\noindent\includegraphics[width=\textwidth]{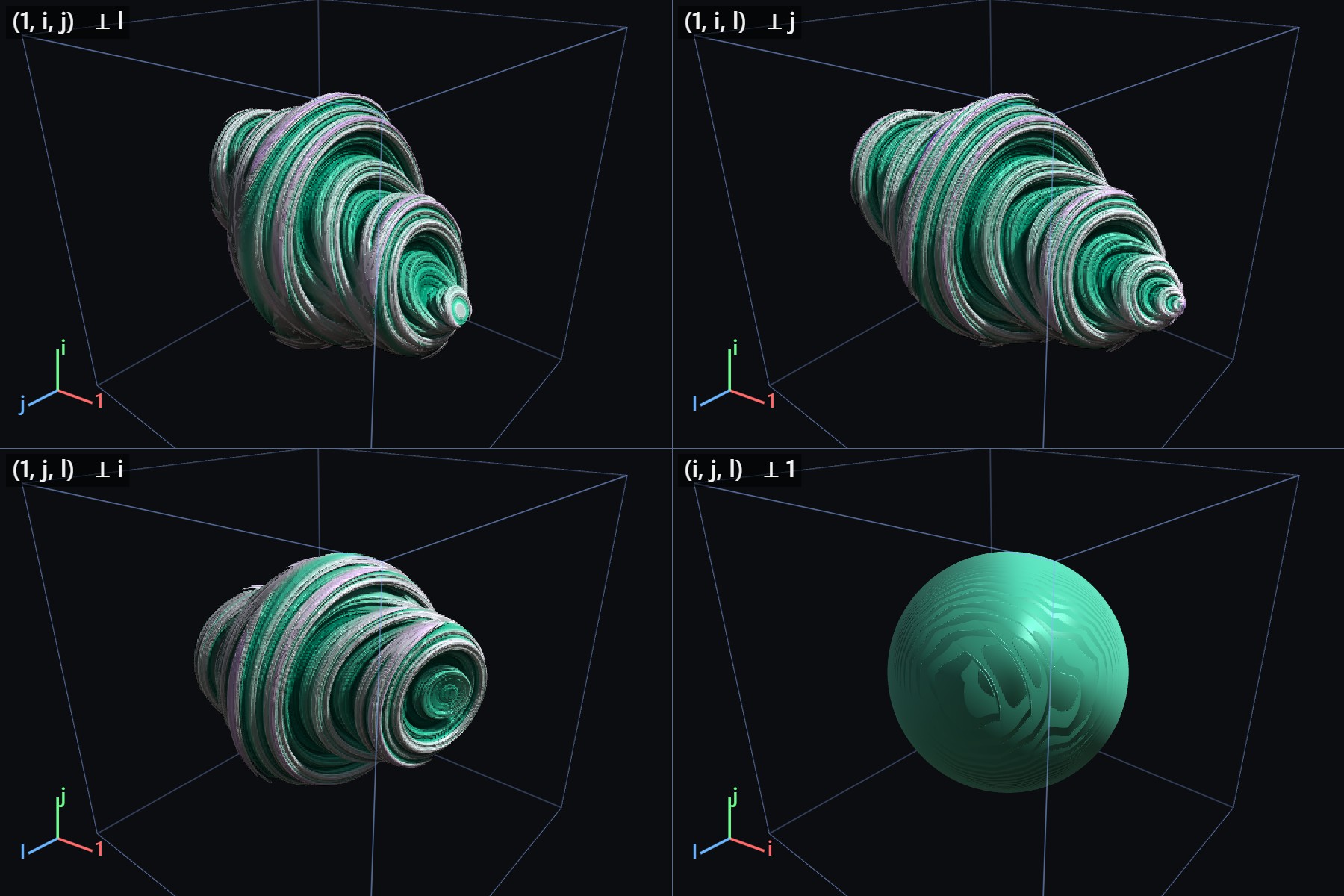}\par\vspace{1mm}
\noindent{\footnotesize \textbf{Map} $z\mapsto(a\diamond z)\diamond(z\diamond b)+c$ $\cdot$ \textbf{Algebra} Okubo algebra (Petersson, $\gamma=(-1,-1,-1)$, $\theta=120^\circ$) $\cdot$ $c=(-0.35,\,\allowbreak 0.45,\,\allowbreak 0,\,\allowbreak 0,\,\allowbreak 0.35,\,\allowbreak 0,\,\allowbreak 0,\,\allowbreak 0)$ $\cdot$ $a=(1,\,\allowbreak 0,\,\allowbreak 0,\,\allowbreak 0,\,\allowbreak 0,\,\allowbreak 0,\,\allowbreak 0,\,\allowbreak 0)$ $\cdot$ $b=(1,\,\allowbreak 0,\,\allowbreak 0,\,\allowbreak 0,\,\allowbreak 0,\,\allowbreak 0,\,\allowbreak 0,\,\allowbreak 0)$ $\cdot$ \textbf{Frame} $X,Y,Z,N=1,i,j,l$ $\cdot$ half-extent 1.4 $\cdot$ 16 iterations $\cdot$ \textbf{Coloring} 8D orbit trap (point), gradient Aurora}\par\vspace{1mm}
\noindent{\footnotesize\begin{tabular}{@{}lllll@{}}\toprule Slice (X, Y, Z) & Body [95\%] & Exact ($J\geq0.9995$) & Approximate ($0.97\leq J<0.9995$) & Dim\\\midrule
($1, i, j$) $\perp l$ & 22.1\% {\scriptsize[21.8, 22.3]} & $M$: z; $R_z$: 180; $P$ & -- & 8 \\
($1, i, l$) $\perp j$ & 22.7\% {\scriptsize[22.4, 23.0]} & $P$ & -- & 4 \\
($1, j, l$) $\perp i$ & 21.3\% {\scriptsize[21.0, 21.5]} & $M$: y; $R_y$: 180; $P$ & -- & 8 \\
($i, j, l$) $\perp 1$ & 29.4\% {\scriptsize[29.1, 29.6]} & $M$: xyz; $R_x$: all; $R_y$: all; $R_z$: all; $P$ & -- & 8 \\
\bottomrule\end{tabular}}\par\vspace{2mm}
With $a=b=1$ the Okubonic bilateral map $(1*z)*(z*1)+c$ \emph{coincides} with the octonionic map $z^2+c$, by Kaplansky's trick (Proposition~\ref{atl:prop:kaplansky}). The plate shows ribbed, spindle-shaped bodies in the three quarters through $1$ and a ball in the imaginary quarter, and it should be compared with the next one.

The symmetries are those of the octonionic Julia set with $c=(-0.35,0.45,0,0,0.35,0,0,0)$. The imaginary quarter is a ball (Proposition~\ref{atl:prop:imaginary}); the mirror in $j$, which fixes $1$ and $c$, is exact in the quarters $(1,i,j)$ and $(1,j,l)$ (Proposition~\ref{atl:prop:equivariance}), together with the half-turn that the point reflection pairs with it (Remark~\ref{atl:rem:minus}); in the quarter $(1,i,l)$, which contains $c$, only the point reflection is detected. The subalgebra generated by $z_0$ and the constants under the Okubo product has dimension $8$ ($4$ in the quarter $(1,i,l)$), but every orbit stays in the three-dimensional space spanned by $z_0$, $1$ and $c$. This plate was rendered with $16$ iterations and Plate~J02 with $20$; the body shares ($21$--$23\%$ through $1$, $29.4\%$ in the imaginary quarter) differ slightly from those of Plate~J02 for this reason only, and Plate~R17, which repeats this plate with $20$ iterations, coincides with Plate~J02.
\clearpage

\subsubsection{J02. Kaplansky: the octonionic Julia set}\label{atl:pl:J02}\mbox{}\par\vspace{0.5mm}
\noindent\includegraphics[width=\textwidth]{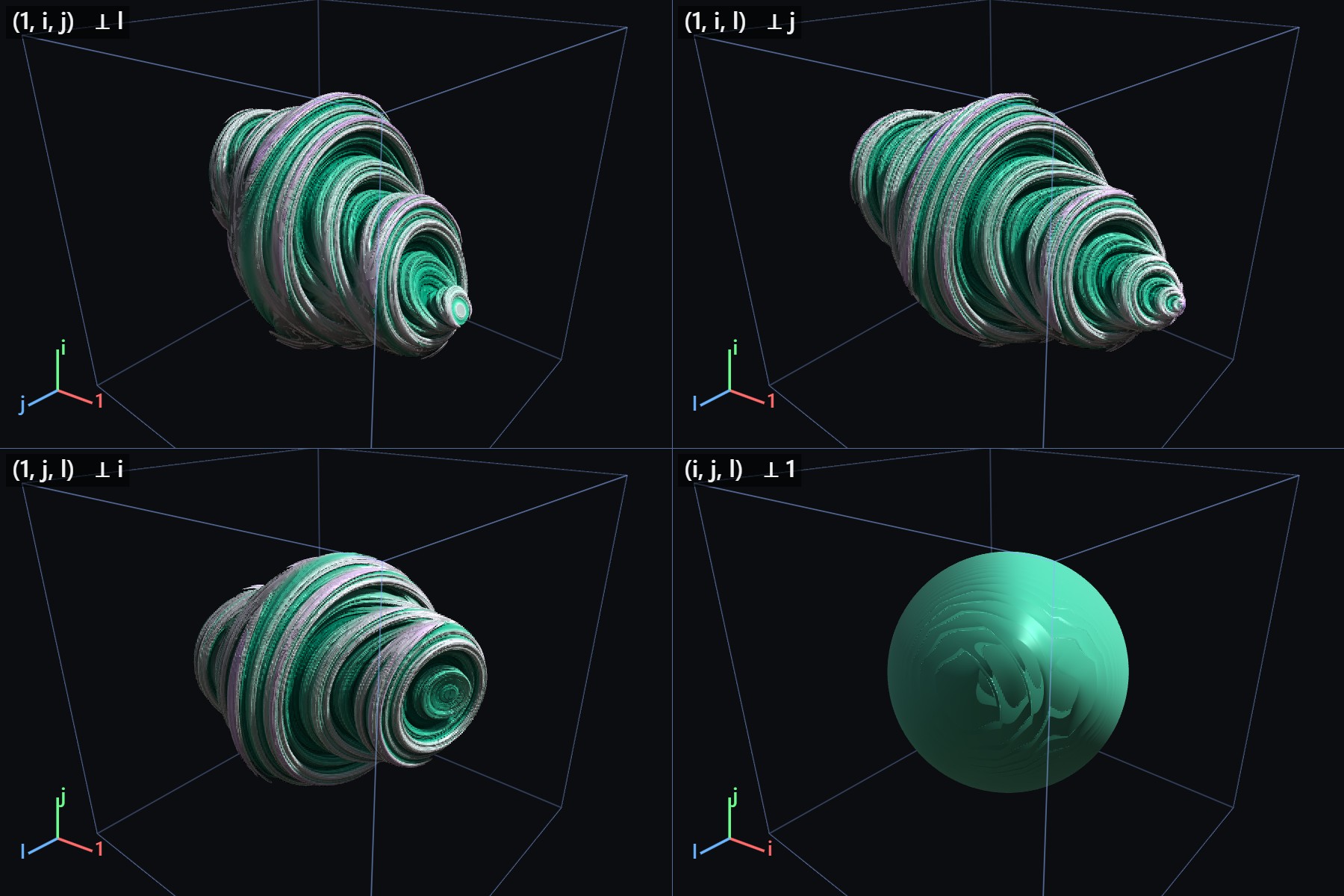}\par\vspace{1mm}
\noindent{\footnotesize \textbf{Map} $z\mapsto z\diamond z+c$ $\cdot$ \textbf{Algebra} Octonions (Cayley--Dickson, $\gamma=(-1,-1,-1)$) $\cdot$ $c=(-0.35,\,\allowbreak 0.45,\,\allowbreak 0,\,\allowbreak 0,\,\allowbreak 0.35,\,\allowbreak 0,\,\allowbreak 0,\,\allowbreak 0)$ $\cdot$ \textbf{Frame} $X,Y,Z,N=1,i,j,l$ $\cdot$ half-extent 1.4 $\cdot$ 20 iterations $\cdot$ \textbf{Coloring} 8D orbit trap (point), gradient Aurora}\par\vspace{1mm}
\noindent{\footnotesize\begin{tabular}{@{}lllll@{}}\toprule Slice (X, Y, Z) & Body [95\%] & Exact ($J\geq0.9995$) & Approximate ($0.97\leq J<0.9995$) & Dim\\\midrule
($1, i, j$) $\perp l$ & 21.4\% {\scriptsize[21.2, 21.7]} & $M$: z; $R_z$: 180; $P$ & -- & 3 \\
($1, i, l$) $\perp j$ & 22.3\% {\scriptsize[22.0, 22.5]} & $P$ & -- & 3 \\
($1, j, l$) $\perp i$ & 20.9\% {\scriptsize[20.7, 21.2]} & $M$: y; $R_y$: 180; $P$ & -- & 3 \\
($i, j, l$) $\perp 1$ & 29.1\% {\scriptsize[28.8, 29.4]} & $M$: xyz; $R_x$: all; $R_y$: all; $R_z$: all; $P$ & -- & 3 \\
\bottomrule\end{tabular}}\par\vspace{2mm}
This is the octonionic Julia set of $z^2+c$ with the same constant, rendered with $20$ iterations. The geometry is the same as in Plate~J01, with the same exact symmetries: a ball in the imaginary quarter, the mirror in $j$ and its paired half-turn in the quarters $(1,i,j)$ and $(1,j,l)$, and the point reflection everywhere. The body shares ($20.9$--$22.3\%$ through $1$, $29.1\%$ in the imaginary quarter) are slightly smaller than those of Plate~J01, which was rendered with only $16$ iterations: the difference comes from the iteration count, not from the algebra. With $20$ iterations the bilateral map gives Plate~R17, whose escape masks coincide with those of this plate on every sampled point. The orbits lie in a subalgebra of dimension $3$ (Proposition~\ref{atl:prop:span}).

The Okubo algebra has no unit, yet its product hides the octonionic one, and Kaplansky's trick exhibits it: the two plates are not merely similar, they are the same set.
\clearpage

\subsubsection{J03. Octonions: a bilateral map}\label{atl:pl:J03}\mbox{}\par\vspace{0.5mm}
\noindent\includegraphics[width=\textwidth]{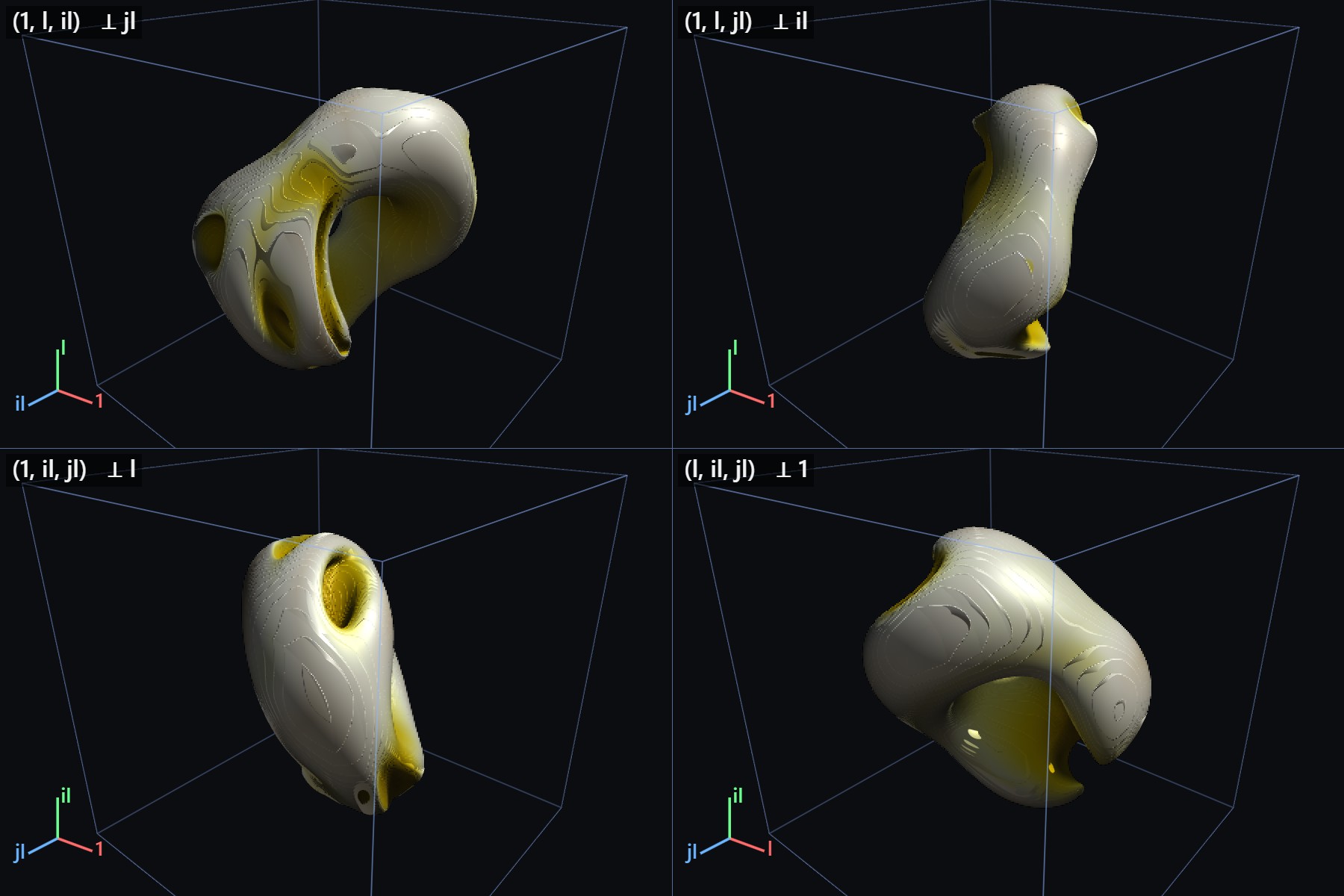}\par\vspace{1mm}
\noindent{\footnotesize \textbf{Map} $z\mapsto(a\diamond z)\diamond(z\diamond b)+c$ $\cdot$ \textbf{Algebra} Octonions (Cayley--Dickson, $\gamma=(-1,-1,-1)$) $\cdot$ $c=(-0.187,\,\allowbreak -0.071,\,\allowbreak -0.22,\,\allowbreak 0.581,\,\allowbreak 0.024,\,\allowbreak 0.037,\,\allowbreak -0.555,\,\allowbreak -0.128)$ $\cdot$ $a=(0.323,\,\allowbreak 0.114,\,\allowbreak -0.581,\,\allowbreak 0.499,\,\allowbreak 0.268,\,\allowbreak 0.164,\,\allowbreak 0.201,\,\allowbreak 0.396)$ $\cdot$ $b=(0.622,\,\allowbreak -0.284,\,\allowbreak -0.427,\,\allowbreak 0.433,\,\allowbreak -0.238,\,\allowbreak -0.154,\,\allowbreak 0.124,\,\allowbreak -0.258)$ $\cdot$ \textbf{Frame} $X,Y,Z,N=1,l,il,jl$ $\cdot$ half-extent 1.6 $\cdot$ 16 iterations $\cdot$ \textbf{Coloring} 8D orbit trap (point), gradient Gold}\par\vspace{1mm}
\noindent{\footnotesize\begin{tabular}{@{}lllll@{}}\toprule Slice (X, Y, Z) & Body [95\%] & Exact ($J\geq0.9995$) & Approximate ($0.97\leq J<0.9995$) & Dim\\\midrule
($1, l, il$) $\perp jl$ & 17.8\% {\scriptsize[17.6, 18.1]} & $P$ & -- & 8 \\
($1, l, jl$) $\perp il$ & 16.4\% {\scriptsize[16.2, 16.7]} & $P$ & -- & 8 \\
($1, il, jl$) $\perp l$ & 20.7\% {\scriptsize[20.5, 21.0]} & $P$ & -- & 8 \\
($l, il, jl$) $\perp 1$ & 21.4\% {\scriptsize[21.2, 21.7]} & $P$ & -- & 8 \\
\bottomrule\end{tabular}}\par\vspace{2mm}
A bilateral map with generic factors $a$ and $b$ in the octonions gives smooth, bent, bone-like bodies, with a tunnel in the quarter $(1,l,il)$, a hollow in $(1,il,jl)$ and an arch in the imaginary quarter; the body shares are $16$--$21\%$. The map is even, and the point reflection is exact in every quarter (Remark~\ref{atl:rem:minus}); no other symmetry is detected.

Unlike the para-octonionic and the Okubonic bilateral maps of Plates~J04--J07, the octonionic bilateral map is not an octonionic quadratic map in disguise: the identity of Proposition~\ref{atl:prop:kaplansky} requires a symmetric composition algebra, and in the octonions $(az)(zb)$ cannot be rebracketed. The orbits generate the whole algebra (dimension $8$), and the bodies have none of the ribbed structure of the quadratic Julia sets.
\clearpage

\subsubsection{J04. Para-octonions: a bilateral map}\label{atl:pl:J04}\mbox{}\par\vspace{0.5mm}
\noindent\includegraphics[width=\textwidth]{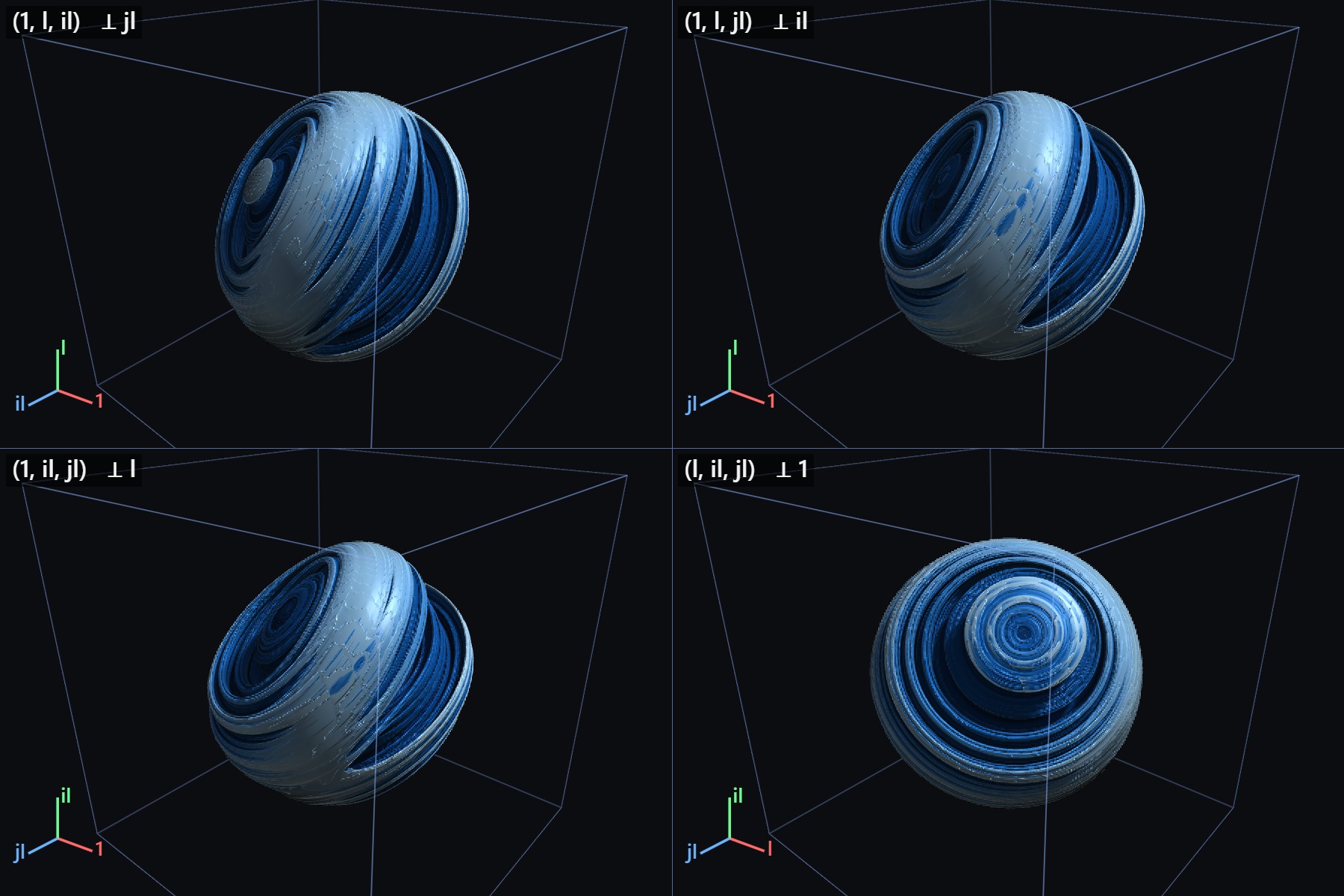}\par\vspace{1mm}
\noindent{\footnotesize \textbf{Map} $z\mapsto(a\diamond z)\diamond(z\diamond b)+c$ $\cdot$ \textbf{Algebra} Para-octonions (para-Hurwitz $\bar x\bar y$, $\gamma=(-1,-1,-1)$) $\cdot$ $c=(-0.482,\,\allowbreak -0.024,\,\allowbreak -0.141,\,\allowbreak -0.023,\,\allowbreak 0.102,\,\allowbreak 0.123,\,\allowbreak 0.314,\,\allowbreak 0.088)$ $\cdot$ $a=(0.334,\,\allowbreak 0.379,\,\allowbreak -0.115,\,\allowbreak -0.28,\,\allowbreak -0.414,\,\allowbreak 0.118,\,\allowbreak 0.006,\,\allowbreak -0.684)$ $\cdot$ $b=(0.438,\,\allowbreak -0.186,\,\allowbreak 0.498,\,\allowbreak 0.07,\,\allowbreak -0.222,\,\allowbreak -0.305,\,\allowbreak -0.452,\,\allowbreak -0.417)$ $\cdot$ \textbf{Frame} $X,Y,Z,N=1,l,il,jl$ $\cdot$ half-extent 1.6 $\cdot$ 16 iterations $\cdot$ \textbf{Coloring} 8D orbit trap (point), gradient Deep Blue}\par\vspace{1mm}
\noindent{\footnotesize\begin{tabular}{@{}lllll@{}}\toprule Slice (X, Y, Z) & Body [95\%] & Exact ($J\geq0.9995$) & Approximate ($0.97\leq J<0.9995$) & Dim\\\midrule
($1, l, il$) $\perp jl$ & 27.5\% {\scriptsize[27.2, 27.8]} & $P$ & -- & 8 \\
($1, l, jl$) $\perp il$ & 27.1\% {\scriptsize[26.8, 27.4]} & $P$ & -- & 8 \\
($1, il, jl$) $\perp l$ & 26.1\% {\scriptsize[25.8, 26.4]} & $P$ & -- & 8 \\
($l, il, jl$) $\perp 1$ & 28.5\% {\scriptsize[28.2, 28.8]} & $P$ & -- & 8 \\
\bottomrule\end{tabular}}\par\vspace{2mm}
The para-octonionic bilateral map produces balls with carved spiral grooves in the three quarters through $1$, and a ball with concentric grooves around a central cap in the imaginary quarter; the body shares are $26$--$29\%$. The point reflection is exact in every quarter (Remark~\ref{atl:rem:minus}), and no coordinate symmetry is detected.

By Proposition~\ref{atl:prop:kaplansky}, $(a\bullet z)\bullet(z\bullet b)=\n(z,e)\,z-\n(z)\,e$ with $e=b\bullet a$, and here $\n(e)=\n(a)\n(b)=1$: the filled set is an octonionic filled Julia set moved by an isometry that sends $1$ to $e$. It is therefore invariant under every rotation fixing $e$ and $c$, a body of revolution about the plane of $e$ and $c$, which is tilted with respect to the coordinate frame; the spiral grooves are the ribs of that body, and in each quarter the symmetry that survives is a mirror in an oblique plane, not tested by the table (a direct check gives $J=1$). The equivalent complex constant is $0.253+0.567\,i$, in the Mandelbrot set. The subalgebra generated by $z_0$, $a$, $b$ and $c$ has dimension $8$, but every orbit stays in the three-dimensional space spanned by $z_0$, $e$ and $c$: the ``orbit span of three'' of the first version is exact, and not an effect of converging orbits.
\clearpage

\subsubsection{J05. Split para-octonions: a bilateral map}\label{atl:pl:J05}\mbox{}\par\vspace{0.5mm}
\noindent\includegraphics[width=\textwidth]{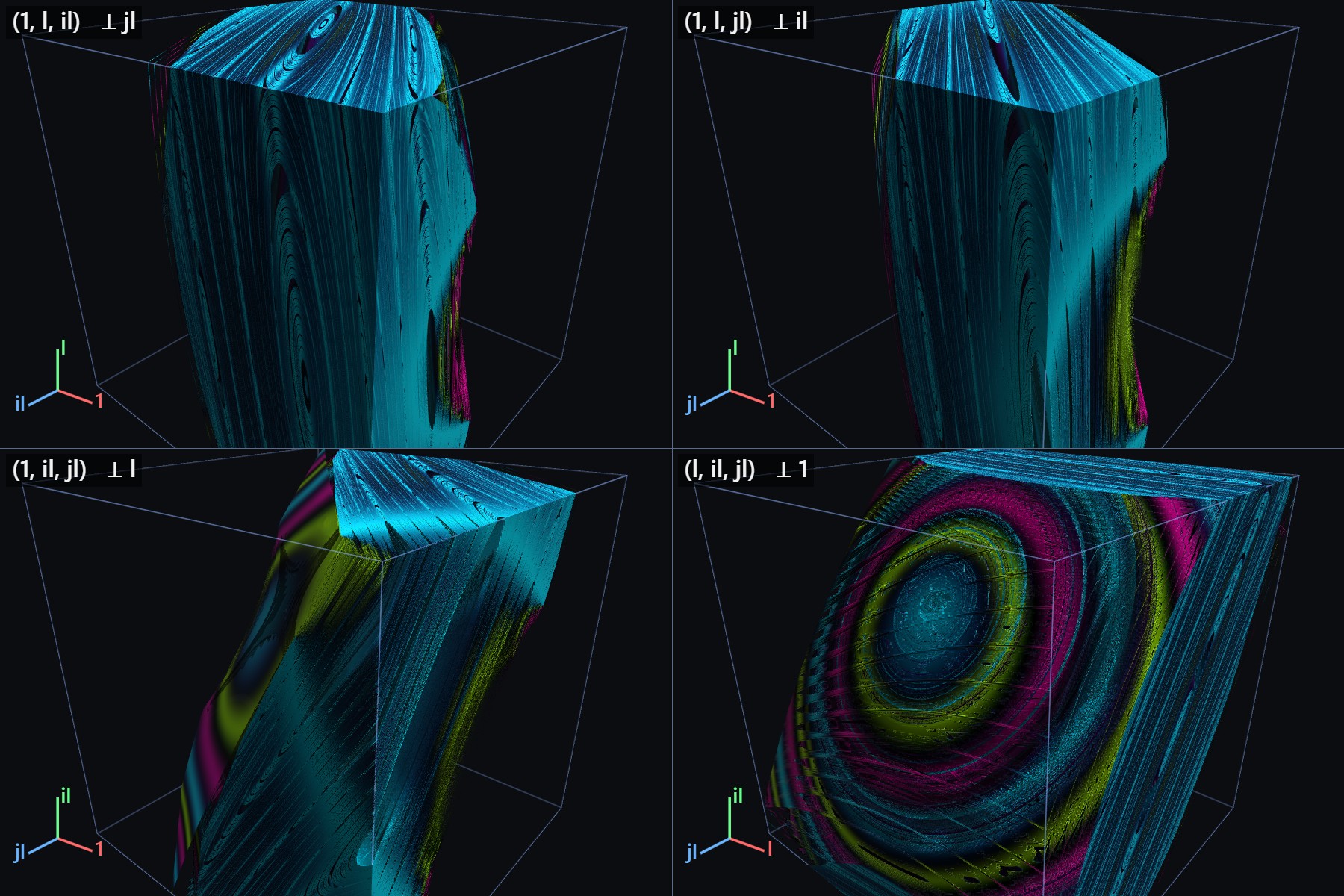}\par\vspace{1mm}
\noindent{\footnotesize \textbf{Map} $z\mapsto(a\diamond z)\diamond(z\diamond b)+c$ $\cdot$ \textbf{Algebra} Split para-octonions (para-Hurwitz $\bar x\bar y$, $\gamma=(-1,-1,1)$) $\cdot$ $c=(-0.465,\,\allowbreak 0.398,\,\allowbreak -0.461,\,\allowbreak -0.02,\,\allowbreak -0.058,\,\allowbreak 0.376,\,\allowbreak 0.148,\,\allowbreak -0.241)$ $\cdot$ $a=(0.539,\,\allowbreak -0.104,\,\allowbreak 0.222,\,\allowbreak -0.15,\,\allowbreak -0.129,\,\allowbreak -0.064,\,\allowbreak -0.401,\,\allowbreak 0.668)$ $\cdot$ $b=(0.361,\,\allowbreak 0.286,\,\allowbreak -0.383,\,\allowbreak -0.157,\,\allowbreak -0.077,\,\allowbreak -0.413,\,\allowbreak -0.507,\,\allowbreak 0.427)$ $\cdot$ \textbf{Frame} $X,Y,Z,N=1,l,il,jl$ $\cdot$ half-extent 1.6 $\cdot$ 16 iterations $\cdot$ \textbf{Coloring} 8D orbit trap (point), gradient Neon}\par\vspace{1mm}
\noindent{\footnotesize\begin{tabular}{@{}lllll@{}}\toprule Slice (X, Y, Z) & Body [95\%] & Exact ($J\geq0.9995$) & Approximate ($0.97\leq J<0.9995$) & Dim\\\midrule
($1, l, il$) $\perp jl$ & 65.6\% {\scriptsize[65.3, 65.9]} & $P$ & -- & 8 \\
($1, l, jl$) $\perp il$ & 56.1\% {\scriptsize[55.8, 56.4]} & $P$ & -- & 8 \\
($1, il, jl$) $\perp l$ & 47.9\% {\scriptsize[47.6, 48.3]} & $P$ & -- & 8 \\
($l, il, jl$) $\perp 1$ & 75.7\% {\scriptsize[75.5, 76.0]} & $P$ & -- & 8 \\
\bottomrule\end{tabular}}\par\vspace{2mm}
The split para-octonionic bilateral map, the best result of the parameter scan of HF8D~0.2, produces bodies that fill most of the box ($48$--$76\%$) and are cut by its faces: striated columns along $l$ and $il$ in the quarters through $1$, and a block in the imaginary quarter. On the flat face where the cube cuts the imaginary quarter, the coloring draws concentric bands, magenta and yellow-green on a cyan ground: these are bands of the orbit-trap coloring on a cross-section, not rings of the geometry. The point reflection is exact in every quarter (Remark~\ref{atl:rem:minus}), and no coordinate symmetry is detected.

The identity of Proposition~\ref{atl:prop:kaplansky} holds in every symmetric composition algebra, the split ones included. Here $e=b\bullet a$ has split norm $\n(e)=0.059>0$, so that the filled set is a split-octonionic filled Julia set moved by a map that preserves the split norm, followed by a dilation by $1/\sqrt{\n(e)}\approx4.1$: this large factor is why the bodies overflow the cube. The map preserving the split norm is not a Euclidean isometry, and the shapes are distorted accordingly; moreover the Euclidean bailout of the renderer is not an escape criterion for an indefinite norm, and the bodies are approximations.
\clearpage

\subsubsection{J06. Okubo: a bilateral map}\label{atl:pl:J06}\mbox{}\par\vspace{0.5mm}
\noindent\includegraphics[width=\textwidth]{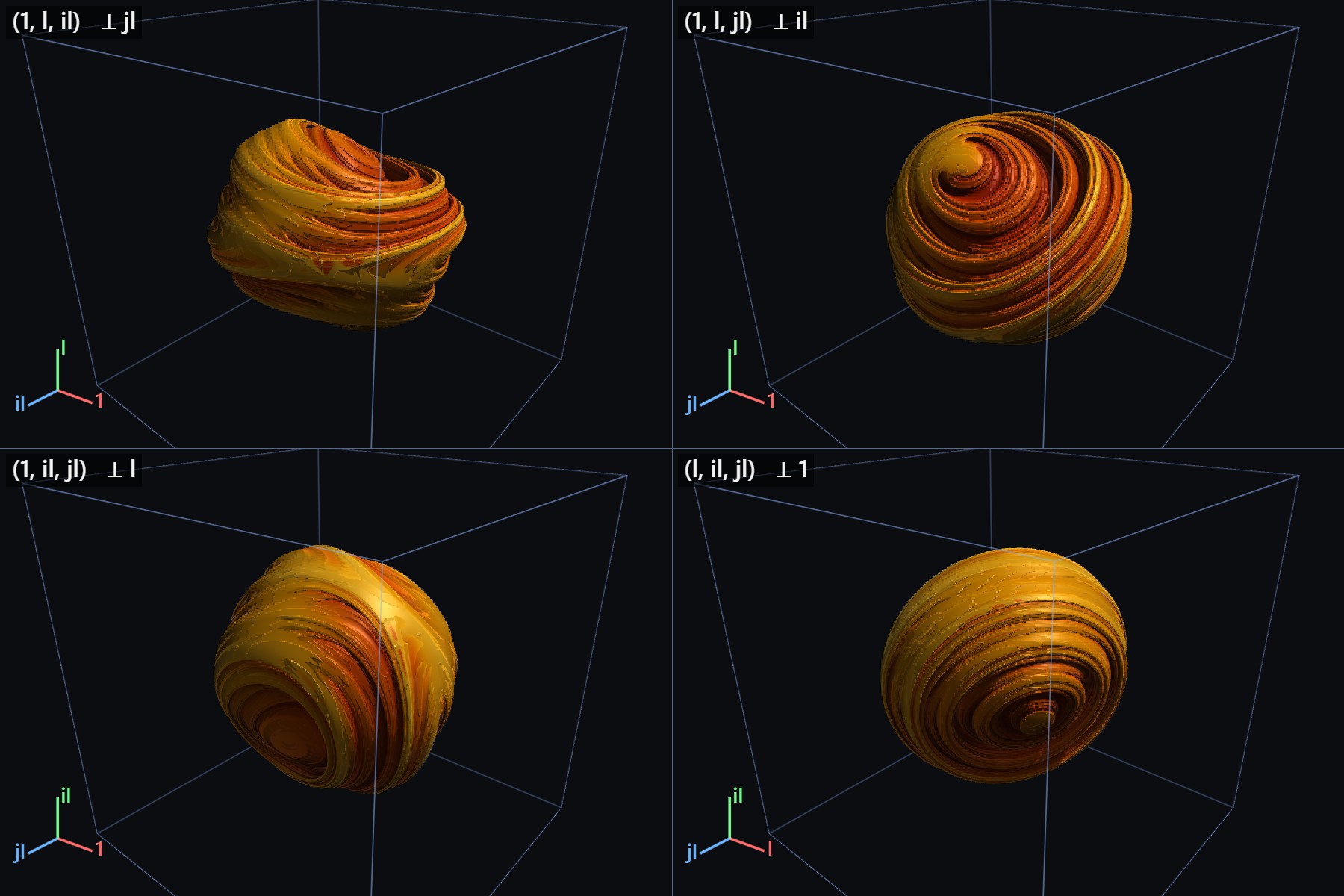}\par\vspace{1mm}
\noindent{\footnotesize \textbf{Map} $z\mapsto(a\diamond z)\diamond(z\diamond b)+c$ $\cdot$ \textbf{Algebra} Okubo algebra (Petersson, $\gamma=(-1,-1,-1)$, $\theta=120^\circ$) $\cdot$ $c=(-0.439,\,\allowbreak 0.022,\,\allowbreak -0.1,\,\allowbreak -0.34,\,\allowbreak -0.096,\,\allowbreak 0.312,\,\allowbreak -0.048,\,\allowbreak 0.057)$ $\cdot$ $a=(0.348,\,\allowbreak -0.5,\,\allowbreak 0.55,\,\allowbreak -0.335,\,\allowbreak -0.254,\,\allowbreak 0.189,\,\allowbreak -0.127,\,\allowbreak 0.312)$ $\cdot$ $b=(0.872,\,\allowbreak 0.094,\,\allowbreak 0.267,\,\allowbreak 0.082,\,\allowbreak 0.233,\,\allowbreak -0.26,\,\allowbreak 0.077,\,\allowbreak -0.159)$ $\cdot$ \textbf{Frame} $X,Y,Z,N=1,l,il,jl$ $\cdot$ half-extent 1.6 $\cdot$ 16 iterations $\cdot$ \textbf{Coloring} 8D orbit trap (point), gradient Fire}\par\vspace{1mm}
\noindent{\footnotesize\begin{tabular}{@{}lllll@{}}\toprule Slice (X, Y, Z) & Body [95\%] & Exact ($J\geq0.9995$) & Approximate ($0.97\leq J<0.9995$) & Dim\\\midrule
($1, l, il$) $\perp jl$ & 22.0\% {\scriptsize[21.8, 22.3]} & $P$ & -- & 8 \\
($1, l, jl$) $\perp il$ & 21.6\% {\scriptsize[21.4, 21.9]} & $P$ & -- & 8 \\
($1, il, jl$) $\perp l$ & 26.3\% {\scriptsize[26.0, 26.5]} & $P$ & -- & 8 \\
($l, il, jl$) $\perp 1$ & 20.4\% {\scriptsize[20.1, 20.6]} & $P$ & -- & 8 \\
\bottomrule\end{tabular}}\par\vspace{2mm}
The Okubonic bilateral map gives swirled balls with spiral ribs, with a clear spiral vortex in the quarter $(1,l,jl)$ and body shares of $20$--$26\%$. The point reflection is exact in every quarter (Remark~\ref{atl:rem:minus}), and no coordinate symmetry is detected.

By Proposition~\ref{atl:prop:kaplansky} the bilateral map is $\n(z,e)\,z-\n(z)\,e+c$ with $e=b*a$, and here $\n(e)=1$: the filled set is an octonionic filled Julia set moved by an isometry that sends $1$ to $e$, in spite of the Okubo product. The spiral ribs are the ribs of an octonionic body of revolution about the plane of $e$ and $c$, tilted with respect to the coordinate frame, and in each quarter the surviving symmetry is a mirror in an oblique plane, which the table does not test (a direct check gives $J=1$). The equivalent complex constant, $-0.031+0.656\,i$, lies just outside the Mandelbrot set, but its critical orbit needs $30$ iterations to escape, more than the $16$ of the rendering. The subalgebra generated by $z_0$ and the constants has dimension $8$, while every orbit stays in the three-dimensional space spanned by $z_0$, $e$ and $c$.
\clearpage

\subsubsection{J07. Split Okubo: a bilateral map}\label{atl:pl:J07}\mbox{}\par\vspace{0.5mm}
\noindent\includegraphics[width=\textwidth]{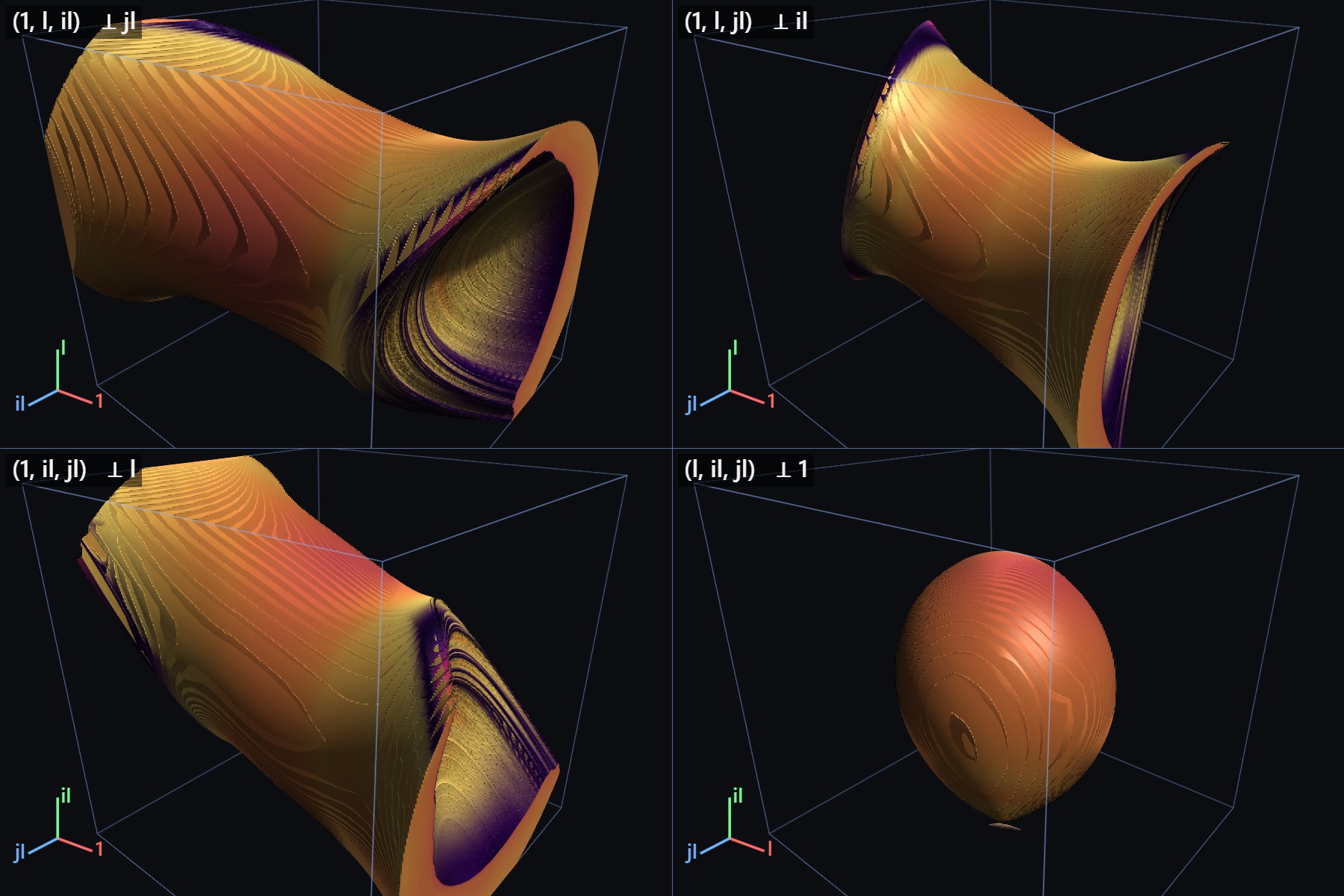}\par\vspace{1mm}
\noindent{\footnotesize \textbf{Map} $z\mapsto(a\diamond z)\diamond(z\diamond b)+c$ $\cdot$ \textbf{Algebra} Split Okubo algebra (Petersson, $\gamma=(-1,-1,1)$, $\theta=120^\circ$) $\cdot$ $c=(0.005,\,\allowbreak 0.153,\,\allowbreak 0.101,\,\allowbreak -0.494,\,\allowbreak 0.308,\,\allowbreak -0.167,\,\allowbreak 0.327,\,\allowbreak 0.488)$ $\cdot$ $a=(0.907,\,\allowbreak -0.072,\,\allowbreak -0.09,\,\allowbreak -0.142,\,\allowbreak -0.14,\,\allowbreak -0.004,\,\allowbreak -0.33,\,\allowbreak 0.126)$ $\cdot$ $b=(0.632,\,\allowbreak 0.564,\,\allowbreak 0.2,\,\allowbreak -0.44,\,\allowbreak 0.142,\,\allowbreak -0.075,\,\allowbreak 0.138,\,\allowbreak 0.068)$ $\cdot$ \textbf{Frame} $X,Y,Z,N=1,l,il,jl$ $\cdot$ half-extent 1.6 $\cdot$ 16 iterations $\cdot$ \textbf{Coloring} 8D orbit trap (point), gradient Sunset}\par\vspace{1mm}
\noindent{\footnotesize\begin{tabular}{@{}lllll@{}}\toprule Slice (X, Y, Z) & Body [95\%] & Exact ($J\geq0.9995$) & Approximate ($0.97\leq J<0.9995$) & Dim\\\midrule
($1, l, il$) $\perp jl$ & 44.7\% {\scriptsize[44.4, 45.0]} & $P$ & -- & 8 \\
($1, l, jl$) $\perp il$ & 48.5\% {\scriptsize[48.2, 48.9]} & $P$ & -- & 8 \\
($1, il, jl$) $\perp l$ & 46.5\% {\scriptsize[46.2, 46.8]} & $P$ & -- & 8 \\
($l, il, jl$) $\perp 1$ & 24.5\% {\scriptsize[24.2, 24.7]} & $P$ & -- & 8 \\
\bottomrule\end{tabular}}\par\vspace{2mm}
The split Okubo bilateral map produces flared, trumpet-like tubes in the quarters through $1$, bounded by hyperbolic surfaces and cut by the faces of the cube ($45$--$49\%$), and an egg in the imaginary quarter ($25\%$). The point reflection is exact in every quarter (Remark~\ref{atl:rem:minus}), and no coordinate symmetry is detected.

As in Plate~J05, the identity of Proposition~\ref{atl:prop:kaplansky} applies: $e=b*a$ has split norm $\n(e)=0.64$, and the filled set is a split-octonionic filled Julia set moved by a map preserving the split norm and dilated by $1/\sqrt{\n(e)}\approx1.25$. The trumpets are the signature of the split norm, the lack of coordinate symmetries the signature of the tilted element $e$. The Euclidean bailout of the renderer is not an escape criterion for the split norm, and the bodies are approximations.
\clearpage

\subsection{Part K. The logistic map}\label{atl:part:K}
The logistic map $c\diamond(z-z\diamond z)$ is conjugate to the quadratic map in the complex plane, but the conjugacy uses the unit and the associativity. In the algebras of this atlas it is an independent family.

\subsubsection{K01. Octonions: small spindles}\label{atl:pl:K01}\mbox{}\par\vspace{0.5mm}
\noindent\includegraphics[width=\textwidth]{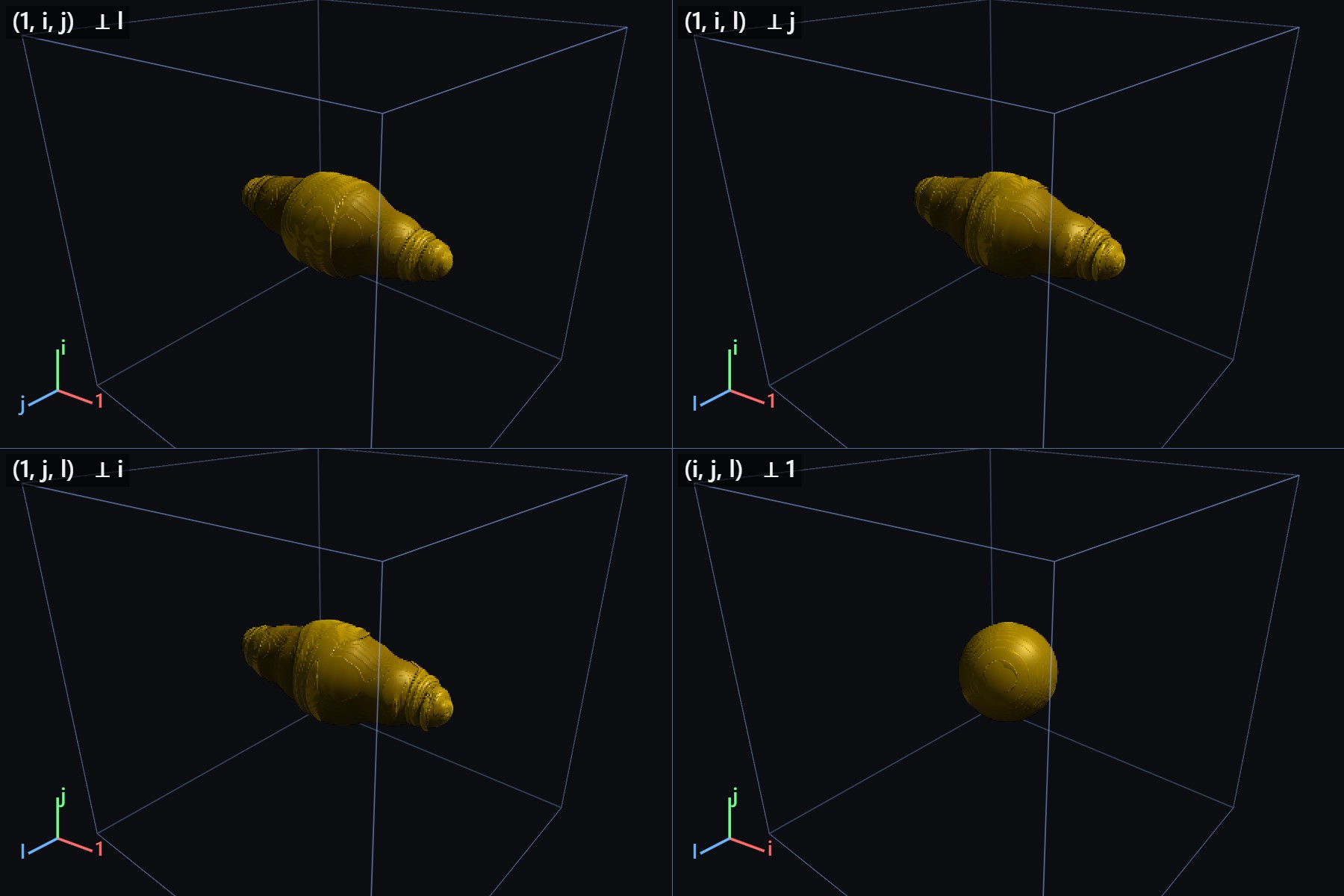}\par\vspace{1mm}
\noindent{\footnotesize \textbf{Map} $z\mapsto c\diamond(z-z\diamond z)$ $\cdot$ \textbf{Algebra} Octonions (Cayley--Dickson, $\gamma=(-1,-1,-1)$) $\cdot$ $c=(-0.824,\,\allowbreak -0.036,\,\allowbreak 0.006,\,\allowbreak 0.287,\,\allowbreak -0.125,\,\allowbreak -0.119,\,\allowbreak 0.031,\,\allowbreak 0.28)$ $\cdot$ \textbf{Frame} $X,Y,Z,N=1,i,j,l$ $\cdot$ center $(0.5,\,\allowbreak 0,\,\allowbreak 0,\,\allowbreak 0,\,\allowbreak 0,\,\allowbreak 0,\,\allowbreak 0,\,\allowbreak 0)$ $\cdot$ half-extent 1.8 $\cdot$ 16 iterations $\cdot$ \textbf{Coloring} 8D orbit trap (point), gradient Gold}\par\vspace{1mm}
\noindent{\footnotesize\begin{tabular}{@{}lllll@{}}\toprule Slice (X, Y, Z) & Body [95\%] & Exact ($J\geq0.9995$) & Approximate ($0.97\leq J<0.9995$) & Dim\\\midrule
($1, i, j$) $\perp l$ & 4.1\% {\scriptsize[4.0, 4.2]} & $P$ & $M$: z; $R_x$: 90,72,60,37; $R_z$: 180 (min.\ 0.972) & 4 \\
($1, i, l$) $\perp j$ & 4.2\% {\scriptsize[4.1, 4.4]} & $P$ & -- & 4 \\
($1, j, l$) $\perp i$ & 4.2\% {\scriptsize[4.1, 4.3]} & $P$ & $M$: y; $R_y$: 180 (min.\ 0.991) & 4 \\
($i, j, l$) $\perp 1$ & 2.1\% {\scriptsize[2.0, 2.2]} & $M$: xyz; $R_x$: all; $R_y$: all; $R_z$: all; $P$ & -- & 4 \\
\bottomrule\end{tabular}}\par\vspace{2mm}
The octonionic logistic map $c(z-z^2)$ produces small spindle-shaped bodies, like beads with fine ridges at both ends, elongated along the $1$-axis in the three quarters through $1$, and a small ball in the imaginary quarter. The spindles occupy only $4\%$ of the ball, the imaginary ball $2\%$.

The centre of the box is the point $\tfrac12$, the critical point of the logistic map, and this explains the two exact symmetries of the table. Writing $z=\tfrac12+w$ one has $z-z^2=\tfrac14-w^2$, an even function of $w$: the point reflection $P$ through the centre is exact in all four quarters, by the argument of Remark~\ref{atl:rem:minus} applied to the variable $w$. In the imaginary quarter $w$ is purely imaginary, $w^2=-\n(w)$, and the first iterate $c\big(\tfrac14+\n(w)\big)$ depends only on $\n(w)$, as in Proposition~\ref{atl:prop:imaginary}: the quarter is an exact ball, and all the tested mirrors and rotations are exact.

In the quaternionic quarter the rotations about the $1$-axis are approximate ($J=0.97$--$0.985$), the mirror in $j$ and the half-turn about $j$ nearly exact ($J=0.995$); in $(1,j,l)$ the mirror in $j$ is approximate ($0.991$), and $(1,i,l)$ has no symmetry beyond $P$. The near-revolution of the quaternionic quarter has a simple origin: the automorphisms of the octonions that preserve the quaternions and rotate the plane $\langle i,j\rangle$ can be chosen to fix the components of $c$ outside that plane, and the components of $c$ along $i$ and $j$ are very small ($-0.036$ and $0.006$). The orbits lie in a subalgebra of dimension $4$: $z-z^2$ lies in the plane of $1$ and $z$, and by Artin's theorem $c$ and $z_0$ generate an associative, quaternionic subalgebra.
\clearpage

\subsubsection{K02. Para-octonions: conch shells}\label{atl:pl:K02}\mbox{}\par\vspace{0.5mm}
\noindent\includegraphics[width=\textwidth]{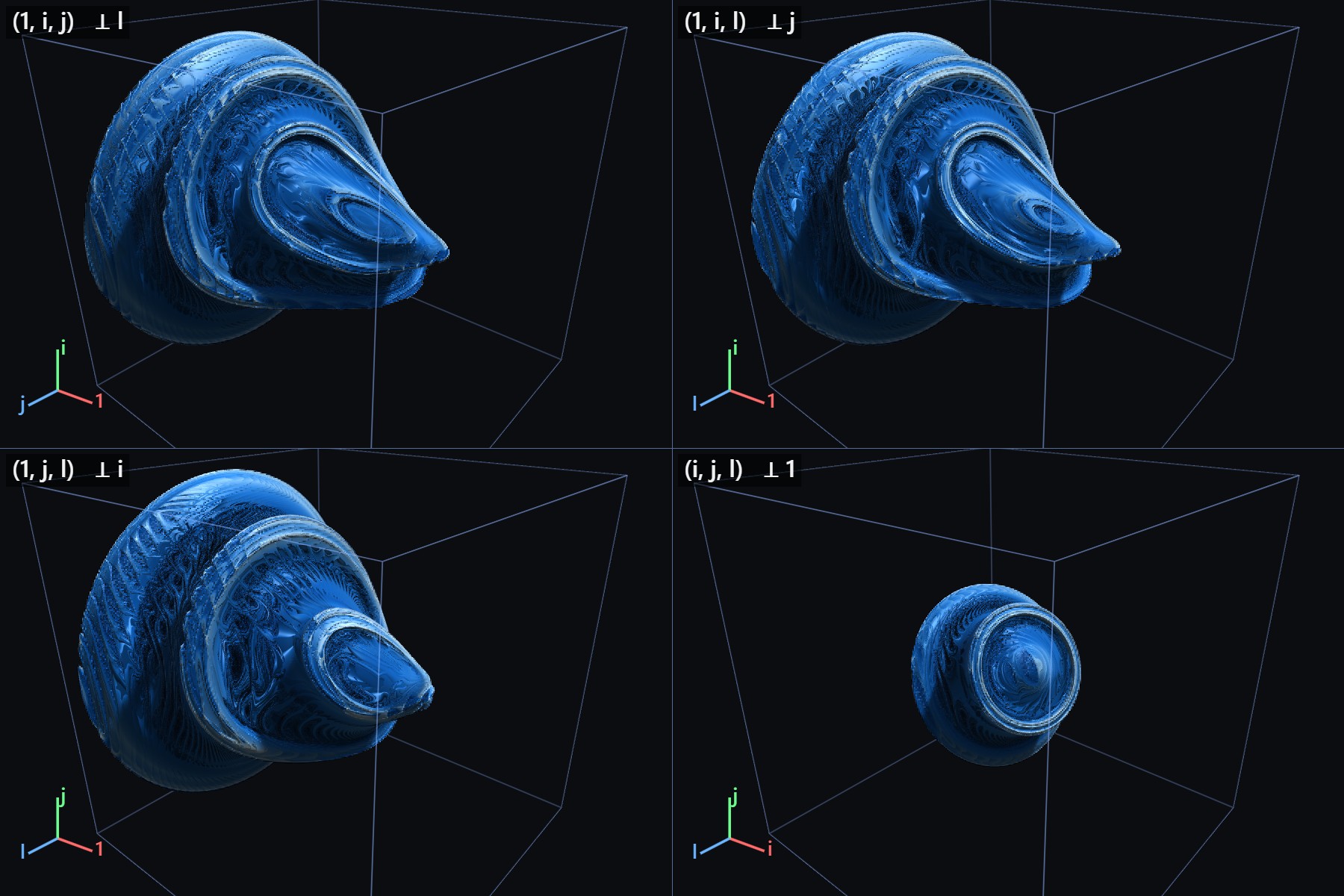}\par\vspace{1mm}
\noindent{\footnotesize \textbf{Map} $z\mapsto c\diamond(z-z\diamond z)$ $\cdot$ \textbf{Algebra} Para-octonions (para-Hurwitz $\bar x\bar y$, $\gamma=(-1,-1,-1)$) $\cdot$ $c=(-0.754,\,\allowbreak -0.37,\,\allowbreak -0.194,\,\allowbreak 0.239,\,\allowbreak -0.098,\,\allowbreak 0.195,\,\allowbreak -0.152,\,\allowbreak -0.115)$ $\cdot$ \textbf{Frame} $X,Y,Z,N=1,i,j,l$ $\cdot$ center $(0.5,\,\allowbreak 0,\,\allowbreak 0,\,\allowbreak 0,\,\allowbreak 0,\,\allowbreak 0,\,\allowbreak 0,\,\allowbreak 0)$ $\cdot$ half-extent 1.8 $\cdot$ 16 iterations $\cdot$ \textbf{Coloring} 8D orbit trap (point), gradient Deep Blue}\par\vspace{1mm}
\noindent{\footnotesize\begin{tabular}{@{}lllll@{}}\toprule Slice (X, Y, Z) & Body [95\%] & Exact ($J\geq0.9995$) & Approximate ($0.97\leq J<0.9995$) & Dim\\\midrule
($1, i, j$) $\perp l$ & 15.0\% {\scriptsize[14.7, 15.2]} & -- & -- & 4 \\
($1, i, l$) $\perp j$ & 15.3\% {\scriptsize[15.1, 15.6]} & -- & -- & 4 \\
($1, j, l$) $\perp i$ & 16.8\% {\scriptsize[16.6, 17.0]} & -- & -- & 4 \\
($i, j, l$) $\perp 1$ & 4.0\% {\scriptsize[3.9, 4.1]} & -- & -- & 4 \\
\bottomrule\end{tabular}}\par\vspace{2mm}
The para-octonionic logistic map gives conch-like shells with thick, rolled lips, similar in the three quarters through $1$ (body $15\%$ to $17\%$), and in the imaginary quarter a small rounded body with a circular rim ($4\%$). The shells are much larger than the octonionic spindles of Plate~K01, but the two plates use different constants, so that the comparison concerns the families rather than the algebras.

No symmetry is found, not even the point reflection: in the para-octonions $z\bullet z=\bar z^{\,2}$, the substitution $z=\tfrac12+w$ leaves a term linear in $w$, and the argument of Plate~K01 fails. The orbits lie in a subalgebra of dimension $4$, as in the octonions: $z-z\bullet z$ lies in the span of $1$ and $z$, and the product with $c$ stays in the quaternionic subalgebra generated by $c$ and $z_0$.
\clearpage

\subsubsection{K03. Split para-octonions: funnels and saddles}\label{atl:pl:K03}\mbox{}\par\vspace{0.5mm}
\noindent\includegraphics[width=\textwidth]{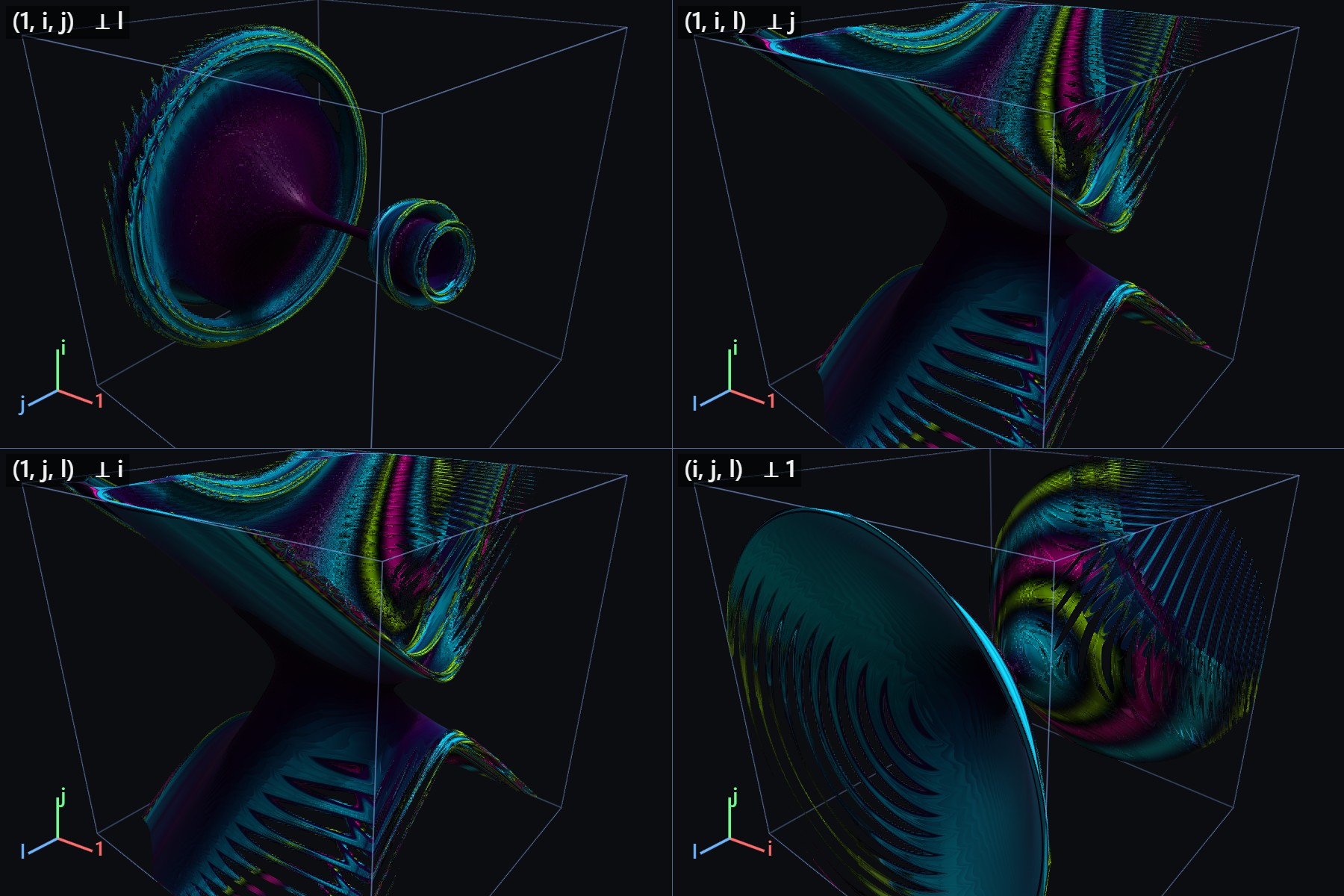}\par\vspace{1mm}
\noindent{\footnotesize \textbf{Map} $z\mapsto c\diamond(z-z\diamond z)$ $\cdot$ \textbf{Algebra} Split para-octonions (para-Hurwitz $\bar x\bar y$, $\gamma=(-1,-1,1)$) $\cdot$ $c=(0.891,\,\allowbreak 0.098,\,\allowbreak 0.086,\,\allowbreak -0.089,\,\allowbreak 0.007,\,\allowbreak 0.001,\,\allowbreak 0.107,\,\allowbreak 0.088)$ $\cdot$ \textbf{Frame} $X,Y,Z,N=1,i,j,l$ $\cdot$ center $(0.5,\,\allowbreak 0,\,\allowbreak 0,\,\allowbreak 0,\,\allowbreak 0,\,\allowbreak 0,\,\allowbreak 0,\,\allowbreak 0)$ $\cdot$ half-extent 1.8 $\cdot$ 16 iterations $\cdot$ \textbf{Coloring} 8D orbit trap (point), gradient Neon}\par\vspace{1mm}
\noindent{\footnotesize\begin{tabular}{@{}lllll@{}}\toprule Slice (X, Y, Z) & Body [95\%] & Exact ($J\geq0.9995$) & Approximate ($0.97\leq J<0.9995$) & Dim\\\midrule
($1, i, j$) $\perp l$ & 3.9\% {\scriptsize[3.7, 4.0]} & -- & -- & 4 \\
($1, i, l$) $\perp j$ & 4.4\% {\scriptsize[4.3, 4.5]} & -- & -- & 4 \\
($1, j, l$) $\perp i$ & 4.2\% {\scriptsize[4.1, 4.3]} & -- & -- & 4 \\
($i, j, l$) $\perp 1$ & 1.5\% {\scriptsize[1.4, 1.6]} & -- & -- & 4 \\
\bottomrule\end{tabular}}\par\vspace{2mm}
In the split para-octonions ($l^2=+1$) the quaternionic quarter shows a wide funnel narrowing into a thin stem, which ends in a small ringed cup. The quarters $(1,i,l)$ and $(1,j,l)$, nearly identical, show two large curved sheets meeting at a narrow waist, like a saddle or an hourglass cut by the faces of the cube; the imaginary quarter shows a wide bell facing a rounded body. The bodies are thin: they occupy only $1.5\%$ to $4.4\%$ of the ball.

No symmetry is found; the best candidate, the mirror in $l$ in the quarters containing $l$, stays at $J=0.91$. The open, hyperbolic shapes of the quarters containing $l$ are consistent with the indefinite norm of the algebra. The orbits lie in a subalgebra of dimension $4$, as in Plates~K01 and~K02.
\clearpage

\subsubsection{K04. Okubo: a trumpet, a dish and blades}\label{atl:pl:K04}\mbox{}\par\vspace{0.5mm}
\noindent\includegraphics[width=\textwidth]{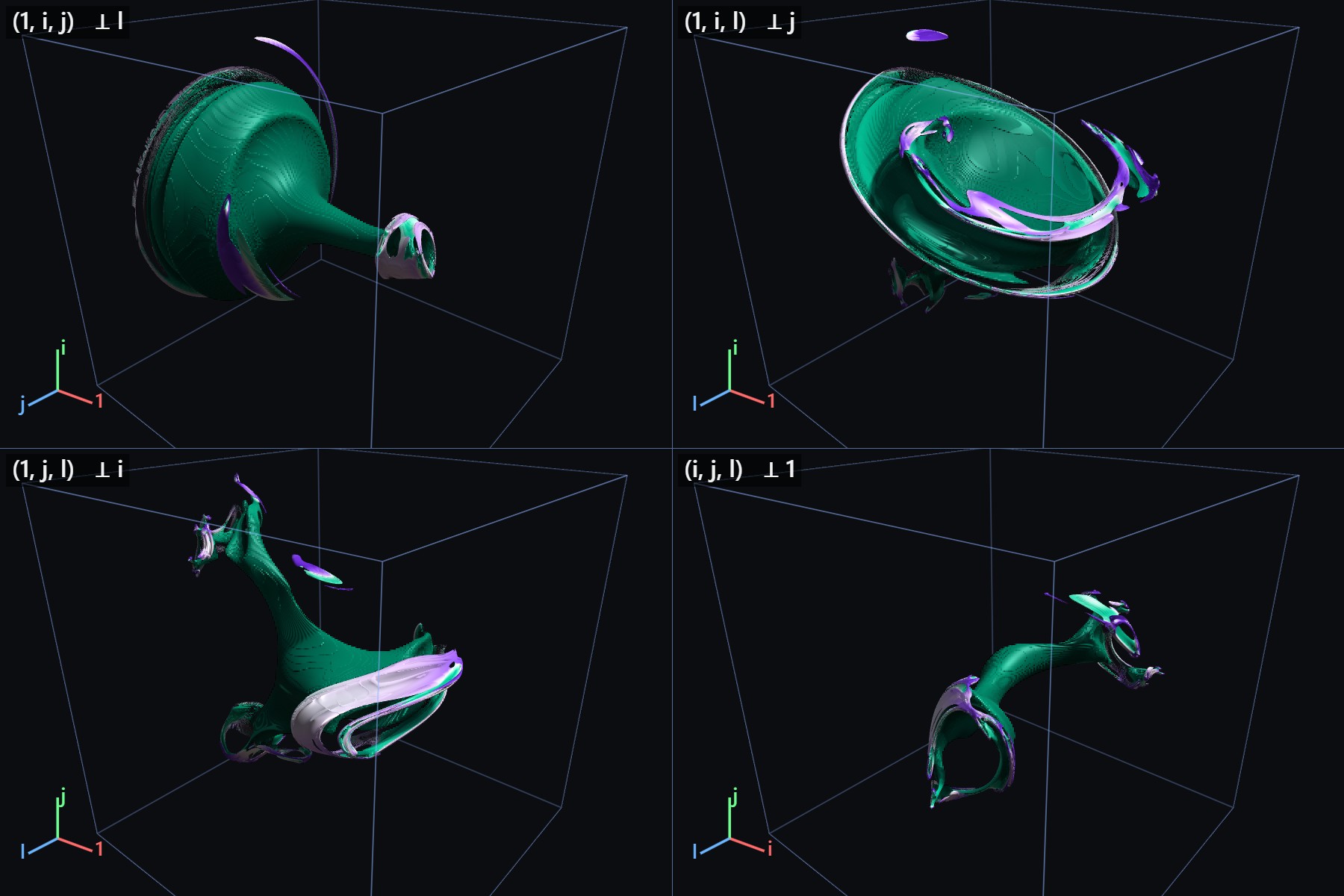}\par\vspace{1mm}
\noindent{\footnotesize \textbf{Map} $z\mapsto c\diamond(z-z\diamond z)$ $\cdot$ \textbf{Algebra} Okubo algebra (Petersson, $\gamma=(-1,-1,-1)$, $\theta=120^\circ$) $\cdot$ $c=(-0.824,\,\allowbreak -0.036,\,\allowbreak 0.006,\,\allowbreak 0.287,\,\allowbreak -0.125,\,\allowbreak -0.119,\,\allowbreak 0.031,\,\allowbreak 0.28)$ $\cdot$ \textbf{Frame} $X,Y,Z,N=1,i,j,l$ $\cdot$ center $(0.5,\,\allowbreak 0,\,\allowbreak 0,\,\allowbreak 0,\,\allowbreak 0,\,\allowbreak 0,\,\allowbreak 0,\,\allowbreak 0)$ $\cdot$ half-extent 1.8 $\cdot$ 16 iterations $\cdot$ \textbf{Coloring} 8D orbit trap (point), gradient Aurora}\par\vspace{1mm}
\noindent{\footnotesize\begin{tabular}{@{}lllll@{}}\toprule Slice (X, Y, Z) & Body [95\%] & Exact ($J\geq0.9995$) & Approximate ($0.97\leq J<0.9995$) & Dim\\\midrule
($1, i, j$) $\perp l$ & 6.6\% {\scriptsize[6.5, 6.8]} & -- & -- & 8 \\
($1, i, l$) $\perp j$ & 5.3\% {\scriptsize[5.2, 5.5]} & -- & -- & 8 \\
($1, j, l$) $\perp i$ & 2.8\% {\scriptsize[2.7, 2.9]} & -- & -- & 8 \\
($i, j, l$) $\perp 1$ & 0.9\% {\scriptsize[0.8, 1.0]} & -- & -- & 8 \\
\bottomrule\end{tabular}}\par\vspace{2mm}
The Okubonic logistic map produces the most intricate shapes of this part. The quaternionic quarter shows a trumpet-like body: a wide bell narrowing into a stem that ends in a small open ring, with a few thin curved blades nearby. The quarter $(1,i,l)$ is a large shallow dish, seen obliquely, with curled fragments along its rim; the quarter $(1,j,l)$ an irregular stalk that flares at its base into layered, curled blades; the imaginary quarter a bent tube, with a loop at one end. The bodies are small ($0.9\%$ to $6.6\%$ of the ball).

No symmetry is found, and the orbits generate the whole algebra (dimension $8$). The comparison with Plate~K01, which uses the same constant, is instructive: there the orbits stay in a quaternionic subalgebra and the logistic map is even about its critical point $\tfrac12$; here the element $1$ is not a unit, neither property survives, and no identity of Proposition~\ref{atl:prop:symmetric} applies to $c*(z-z*z)$, which, unlike the sandwich and the cubic map, does not reduce to the norm. The constant is generic, so that no rotation from the torus of Proposition~\ref{atl:prop:torus} is to be expected in this frame either.
\clearpage

\subsubsection{K05. Split Okubo: flames}\label{atl:pl:K05}\mbox{}\par\vspace{0.5mm}
\noindent\includegraphics[width=\textwidth]{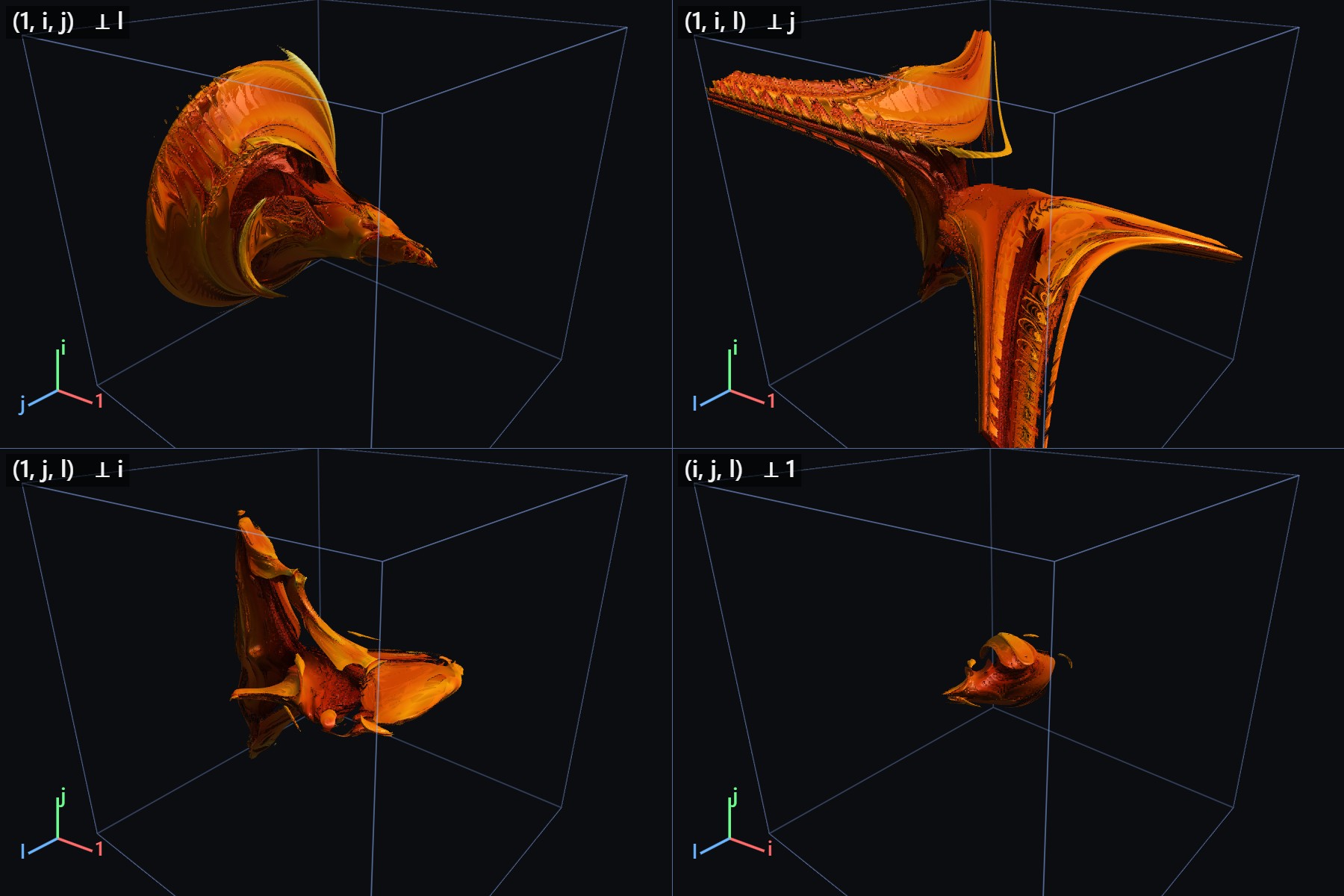}\par\vspace{1mm}
\noindent{\footnotesize \textbf{Map} $z\mapsto c\diamond(z-z\diamond z)$ $\cdot$ \textbf{Algebra} Split Okubo algebra (Petersson, $\gamma=(-1,-1,1)$, $\theta=120^\circ$) $\cdot$ $c=(-1.19,\,\allowbreak 0.1,\,\allowbreak -0.185,\,\allowbreak -0.11,\,\allowbreak 0.013,\,\allowbreak 0.194,\,\allowbreak 0.128,\,\allowbreak -0.116)$ $\cdot$ \textbf{Frame} $X,Y,Z,N=1,i,j,l$ $\cdot$ center $(0.5,\,\allowbreak 0,\,\allowbreak 0,\,\allowbreak 0,\,\allowbreak 0,\,\allowbreak 0,\,\allowbreak 0,\,\allowbreak 0)$ $\cdot$ half-extent 1.8 $\cdot$ 16 iterations $\cdot$ \textbf{Coloring} 8D orbit trap (point), gradient Fire}\par\vspace{1mm}
\noindent{\footnotesize\begin{tabular}{@{}lllll@{}}\toprule Slice (X, Y, Z) & Body [95\%] & Exact ($J\geq0.9995$) & Approximate ($0.97\leq J<0.9995$) & Dim\\\midrule
($1, i, j$) $\perp l$ & 6.1\% {\scriptsize[6.0, 6.3]} & -- & -- & 8 \\
($1, i, l$) $\perp j$ & 3.9\% {\scriptsize[3.8, 4.0]} & -- & -- & 8 \\
($1, j, l$) $\perp i$ & 2.4\% {\scriptsize[2.3, 2.5]} & -- & -- & 8 \\
($i, j, l$) $\perp 1$ & 0.5\% {\scriptsize[0.5, 0.6]} & -- & -- & 8 \\
\bottomrule\end{tabular}}\par\vspace{2mm}
The split Okubo logistic map gives thin, curved, flame-like blades: a curled flame ending in a pointed tip in the quaternionic quarter, a long crossed structure of ribbed blades in $(1,i,l)$, a cluster of twisted blades in $(1,j,l)$ and a small curled fragment in the imaginary quarter. The bodies are small and irregular ($0.5\%$ to $6.1\%$ of the ball).

No symmetry is found, and the orbits generate the whole algebra (dimension $8$) in all four quarters. The constant differs from that of Plate~K04, so the two plates show the two families rather than the effect of the split norm alone.
\clearpage

\subsection{Part L. A transcendental map in the unital algebras}\label{atl:part:L}
Transcendental functions are defined only in the unital algebras, through the complex function applied in the plane of $1$ and $\im z$. We draw the map $z\mapsto c\sin z$ in the four Cayley--Dickson algebras. The quarters where the body fills the cube show large regions of bounded orbits typical of the sine family.

\subsubsection{L01. Octonions: the sine family, nearly of revolution}\label{atl:pl:L01}\mbox{}\par\vspace{0.5mm}
\noindent\includegraphics[width=\textwidth]{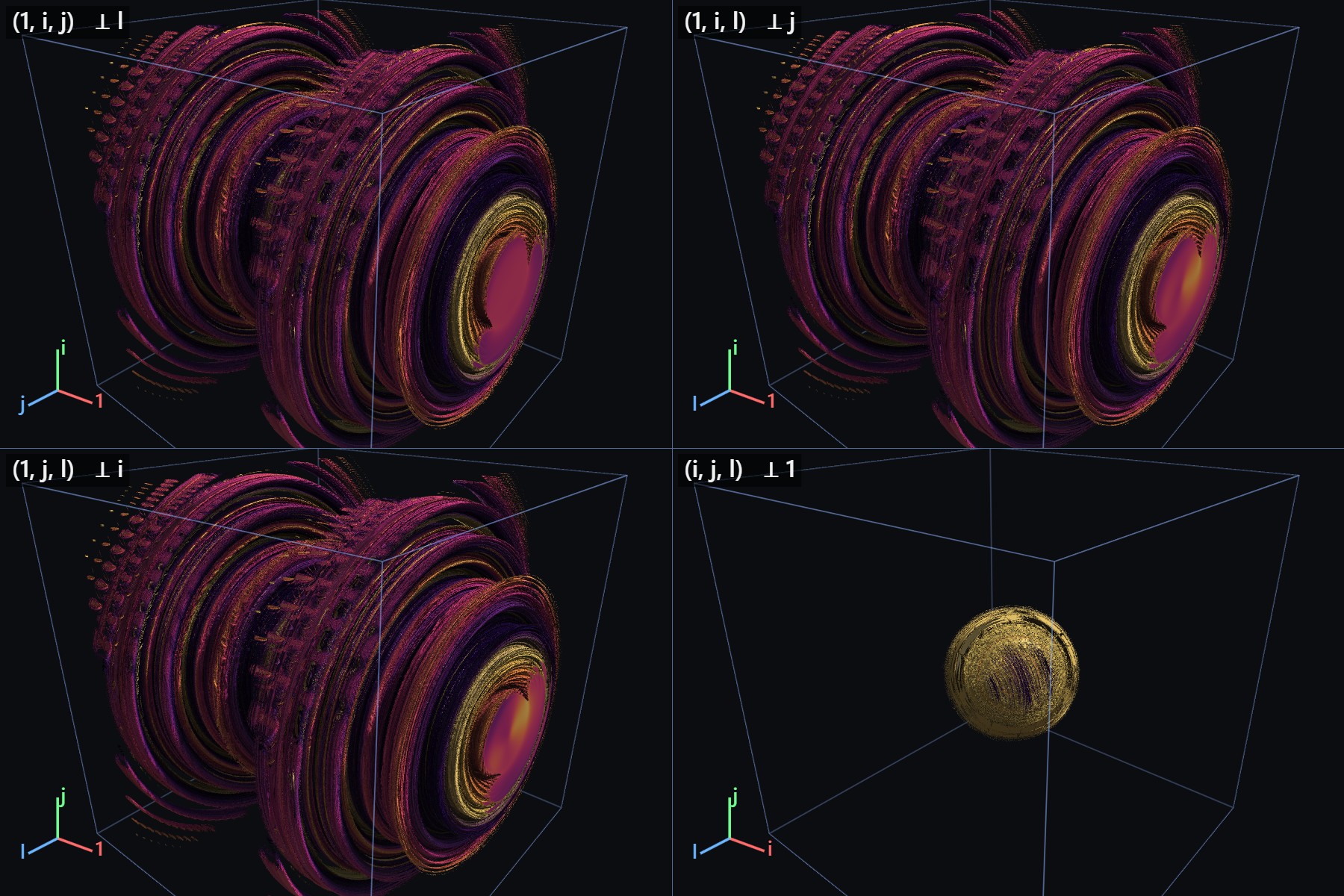}\par\vspace{1mm}
\noindent{\footnotesize \textbf{Map} $z\mapsto c\,\sin z$ $\cdot$ \textbf{Algebra} Octonions (Cayley--Dickson, $\gamma=(-1,-1,-1)$) $\cdot$ $c=(0.976,\,\allowbreak -0.091,\,\allowbreak -0.052,\,\allowbreak -0.054,\,\allowbreak -0.081,\,\allowbreak -0.157,\,\allowbreak -0.129,\,\allowbreak 0.021)$ $\cdot$ \textbf{Frame} $X,Y,Z,N=1,i,j,l$ $\cdot$ half-extent 3 $\cdot$ 12 iterations $\cdot$ \textbf{Coloring} 8D orbit trap (point), gradient Sunset}\par\vspace{1mm}
\noindent{\footnotesize\begin{tabular}{@{}lllll@{}}\toprule Slice (X, Y, Z) & Body [95\%] & Exact ($J\geq0.9995$) & Approximate ($0.97\leq J<0.9995$) & Dim\\\midrule
($1, i, j$) $\perp l$ & 53.7\% {\scriptsize[53.4, 54.0]} & $P$ & -- & 4 \\
($1, i, l$) $\perp j$ & 53.8\% {\scriptsize[53.4, 54.1]} & $P$ & -- & 4 \\
($1, j, l$) $\perp i$ & 53.7\% {\scriptsize[53.4, 54.0]} & $P$ & -- & 4 \\
($i, j, l$) $\perp 1$ & 3.0\% {\scriptsize[2.9, 3.1]} & $P$ & -- & 4 \\
\bottomrule\end{tabular}}\par\vspace{2mm}
The octonionic map $c\sin z$ is computed by applying the complex sine in the plane of $1$ and $\im z$, so that $\sin z$ lies in the plane of $1$ and $z$ and commutes with every rotation fixing $1$; the symmetry is broken only by the left multiplication by $c$, through the imaginary part of $c$, which for the constant of this plate is about a quarter of its real part. The three quarters through $1$ are indeed nearly identical (body $54\%$ in each) and show the periodic structure of the sine: nested shells and rings stacked along the real axis, with a smooth cap at one end. They look like bodies of revolution about the $1$-axis, but they are not: the best rotation about that axis, by $37^\circ$, reaches only $J=0.92$--$0.93$.

The imaginary quarter is a small, roundish, grainy body ($3\%$ of the ball). For purely imaginary $z$ the sine is purely imaginary, $\sin(tv)=v\sinh t$ for a unit imaginary $v$, but the first iterate $c\,v\sinh t$ still depends on the direction $v$; on a body this small the measurements do not confirm any rotation ($J\le0.64$).

The only exact symmetry, in all four quarters, is the point reflection $P$. It is not the symmetry of the even maps of Remark~\ref{atl:rem:minus}: here the map is odd, $c\sin(-z)=-c\sin z$, so the orbit of $-z_0$ is the opposite of the orbit of $z_0$ and escapes at the same time. The orbits lie in a subalgebra of dimension $4$: $\sin z$ stays in the plane of $1$ and $z$, and by Artin's theorem $c$ and $z_0$ generate an associative, quaternionic subalgebra of the octonions.
\clearpage

\subsubsection{L02. Split-octonions (\texorpdfstring{$l^2=+1$}{l\texttwosuperior{}=+1})}\label{atl:pl:L02}\mbox{}\par\vspace{0.5mm}
\noindent\includegraphics[width=\textwidth]{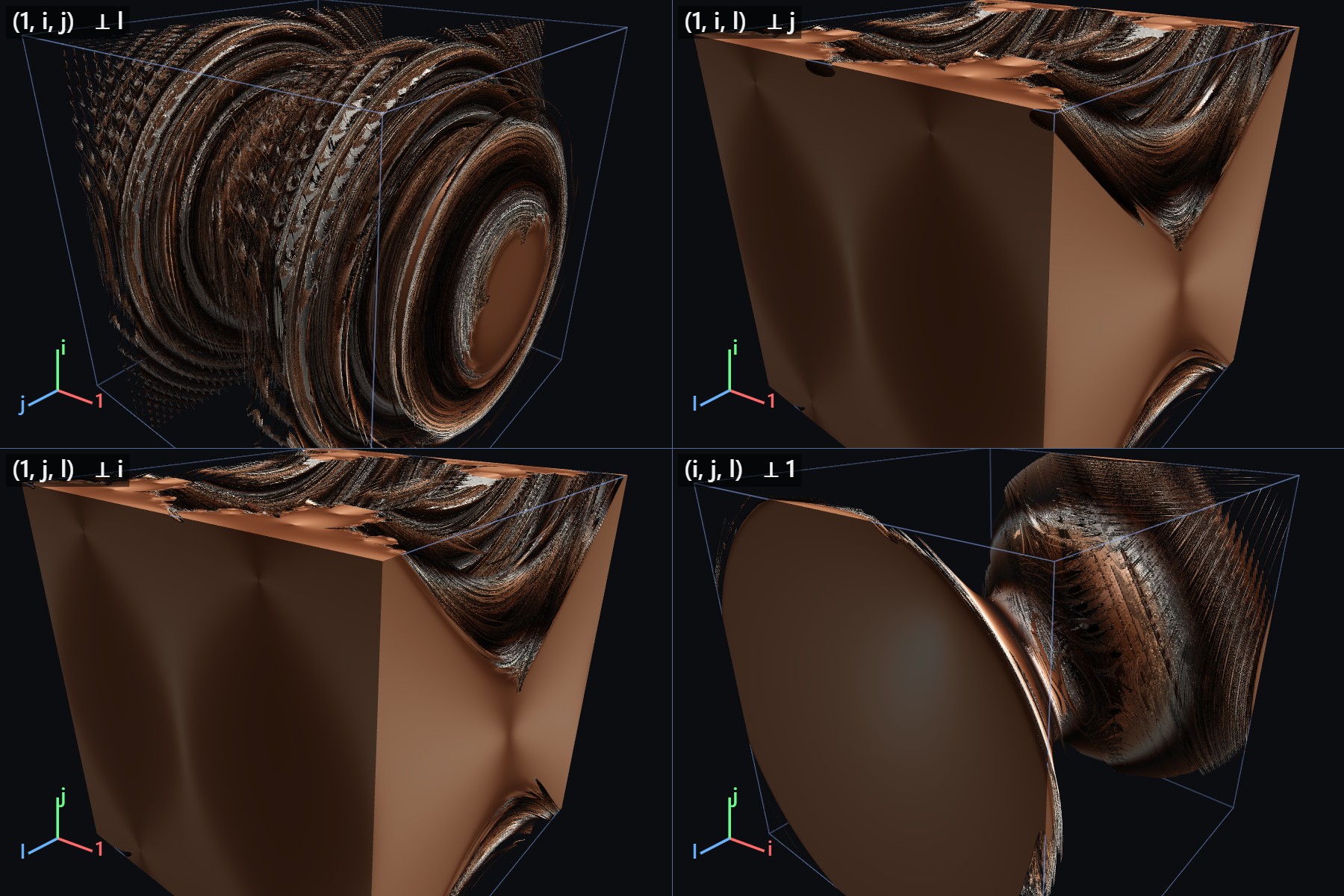}\par\vspace{1mm}
\noindent{\footnotesize \textbf{Map} $z\mapsto c\,\sin z$ $\cdot$ \textbf{Algebra} Split-octonions ($l^2=+1$) (Cayley--Dickson, $\gamma=(-1,-1,1)$) $\cdot$ $c=(0.976,\,\allowbreak -0.091,\,\allowbreak -0.052,\,\allowbreak -0.054,\,\allowbreak -0.081,\,\allowbreak -0.157,\,\allowbreak -0.129,\,\allowbreak 0.021)$ $\cdot$ \textbf{Frame} $X,Y,Z,N=1,i,j,l$ $\cdot$ half-extent 3 $\cdot$ 12 iterations $\cdot$ \textbf{Coloring} 8D orbit trap (point), gradient Copper}\par\vspace{1mm}
\noindent{\footnotesize\begin{tabular}{@{}lllll@{}}\toprule Slice (X, Y, Z) & Body [95\%] & Exact ($J\geq0.9995$) & Approximate ($0.97\leq J<0.9995$) & Dim\\\midrule
($1, i, j$) $\perp l$ & 63.5\% {\scriptsize[63.2, 63.8]} & $P$ & -- & 4 \\
($1, i, l$) $\perp j$ & 88.3\% {\scriptsize[88.1, 88.5]} & $P$ & $M$: z; $R_z$: 180 (min.\ 0.970) & 4 \\
($1, j, l$) $\perp i$ & 88.3\% {\scriptsize[88.1, 88.5]} & $P$ & -- & 4 \\
($i, j, l$) $\perp 1$ & 44.7\% {\scriptsize[44.4, 45.0]} & $P$ & -- & 4 \\
\bottomrule\end{tabular}}\par\vspace{2mm}
In the split-octonions the sine is applied in the plane of $1$ and $\im z$, which is a complex plane when $\n(\im z)>0$ and a split-complex plane when $\n(\im z)<0$; in the latter case $f(r+v)=\tfrac12\big(f(r+m)+f(r-m)\big)+\tfrac{v}{2m}\big(f(r+m)-f(r-m)\big)$ with $m^2=-\n(v)$, and near the null cone the implementation uses the first-order expansion $f(r)+v\,f'(r)$.

The quaternionic quarter shows again the nested shells and rings of the sine family (body $64\%$). The quarters $(1,i,l)$ and $(1,j,l)$ are filled cubes ($88\%$), whose top faces reveal fibrous structures; the imaginary quarter shows a large flat-faced disc joined through a narrow neck to a barrel-shaped body ($45\%$). In the hyperbolic directions the sine is evaluated at the real points $r\pm m$, where it is bounded, which may explain the large body shares.

As in Plate~L01 the map is odd and the point reflection $P$ is exact in every quarter; in the quarter $(1,i,l)$ the half-turn about $l$ and the mirror in $l$ are only approximate ($J=0.970$), and no other symmetry is found. The orbits lie in a subalgebra of dimension $4$.
\clearpage

\subsubsection{L03. Semi-octonions: columns along \texorpdfstring{$l$}{l}}\label{atl:pl:L03}\mbox{}\par\vspace{0.5mm}
\noindent\includegraphics[width=\textwidth]{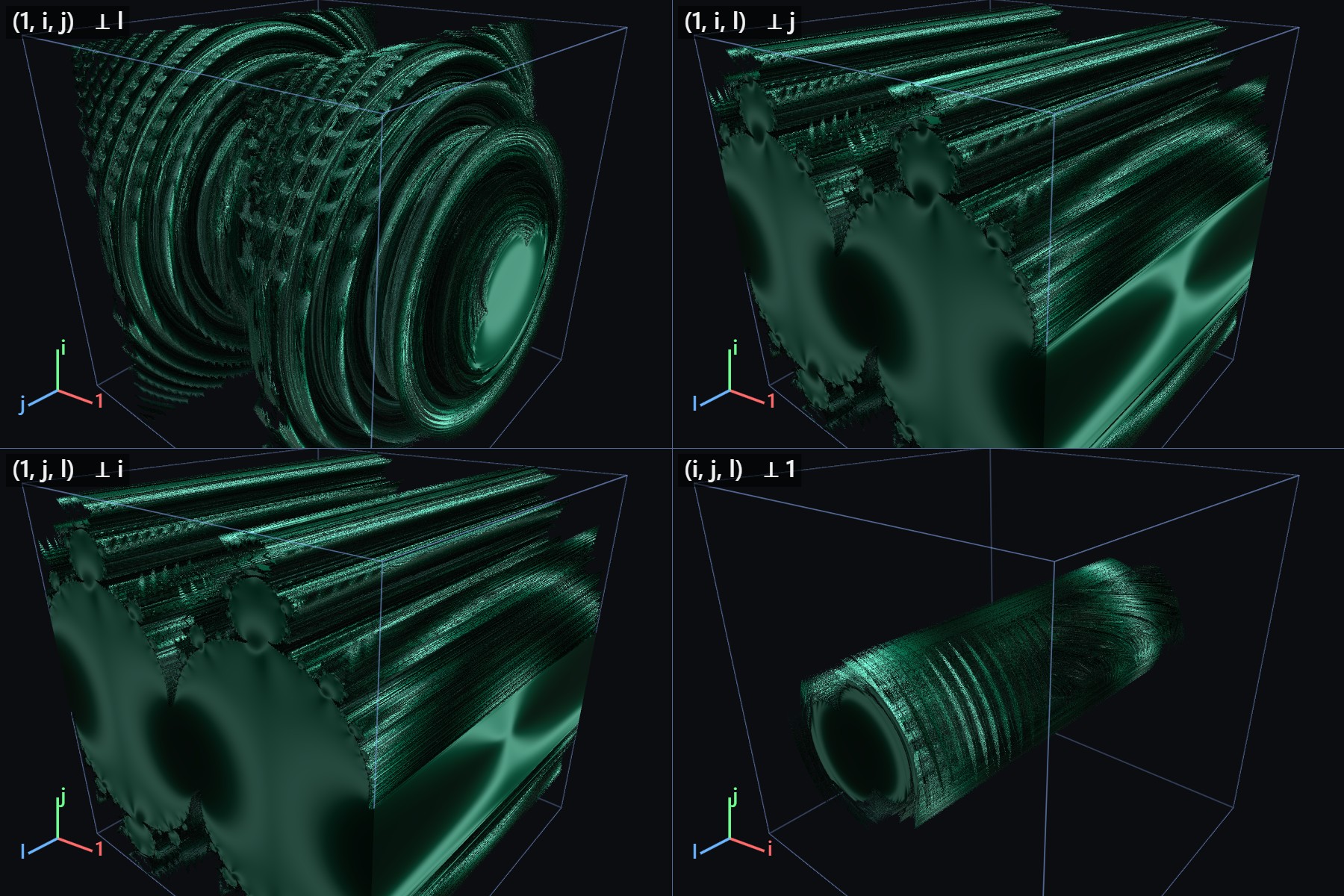}\par\vspace{1mm}
\noindent{\footnotesize \textbf{Map} $z\mapsto c\,\sin z$ $\cdot$ \textbf{Algebra} Semi-octonions ($l^2=0$) (Cayley--Dickson, $\gamma=(-1,-1,0)$) $\cdot$ $c=(0.976,\,\allowbreak -0.091,\,\allowbreak -0.052,\,\allowbreak -0.054,\,\allowbreak -0.081,\,\allowbreak -0.157,\,\allowbreak -0.129,\,\allowbreak 0.021)$ $\cdot$ \textbf{Frame} $X,Y,Z,N=1,i,j,l$ $\cdot$ half-extent 3 $\cdot$ 12 iterations $\cdot$ \textbf{Coloring} 8D orbit trap (point), gradient Emerald}\par\vspace{1mm}
\noindent{\footnotesize\begin{tabular}{@{}lllll@{}}\toprule Slice (X, Y, Z) & Body [95\%] & Exact ($J\geq0.9995$) & Approximate ($0.97\leq J<0.9995$) & Dim\\\midrule
($1, i, j$) $\perp l$ & 59.9\% {\scriptsize[59.6, 60.2]} & $P$ & -- & 4 \\
($1, i, l$) $\perp j$ & 80.5\% {\scriptsize[80.3, 80.7]} & $P$ & $M$: z; $R_z$: 180 (min.\ 0.998) & 4 \\
($1, j, l$) $\perp i$ & 80.3\% {\scriptsize[80.1, 80.6]} & $P$ & $M$: z; $R_z$: 180 (min.\ 0.998) & 4 \\
($i, j, l$) $\perp 1$ & 17.4\% {\scriptsize[17.2, 17.6]} & $P$ & -- & 4 \\
\bottomrule\end{tabular}}\par\vspace{2mm}
In the semi-octonions ($l^2=0$) the sine family is extruded along $l$. The quarters $(1,i,l)$ and $(1,j,l)$ are filled with striated columns parallel to the $l$-axis (body $80\%$); where the columns are cut by the faces of the cube, the cross-sections show the round bulbs of the filled Julia set of the sine. The imaginary quarter is a striated tube along $l$ ($17\%$), and the quaternionic quarter shows the shells and rings of the sine family ($60\%$).

The mechanism is the one described after Proposition~\ref{atl:prop:rescaling} for the quadratic map: the norm does not see the $l$-part, the quaternionic part of the orbit evolves autonomously, and the $l$-part is transported linearly by the derivative of the function. This also explains why the quaternionic quarter differs from the one of Plate~L01 ($54\%$): there the $l$-components of the constant feed back into the quaternionic part, here they do not.

The columns are not exact cylinders: the reflection $l\mapsto-l$ is only approximate ($J=0.998$ in both quarters), since the $l$-part is driven by the $l$-components of the constant and may grow where the derivative expands. The point reflection $P$ is exact, as for every odd map (Plate~L01), and the orbits lie in a subalgebra of dimension $4$.
\clearpage

\subsubsection{L04. Split-octonions (\texorpdfstring{$i^2=+1$}{i\texttwosuperior{}=+1})}\label{atl:pl:L04}\mbox{}\par\vspace{0.5mm}
\noindent\includegraphics[width=\textwidth]{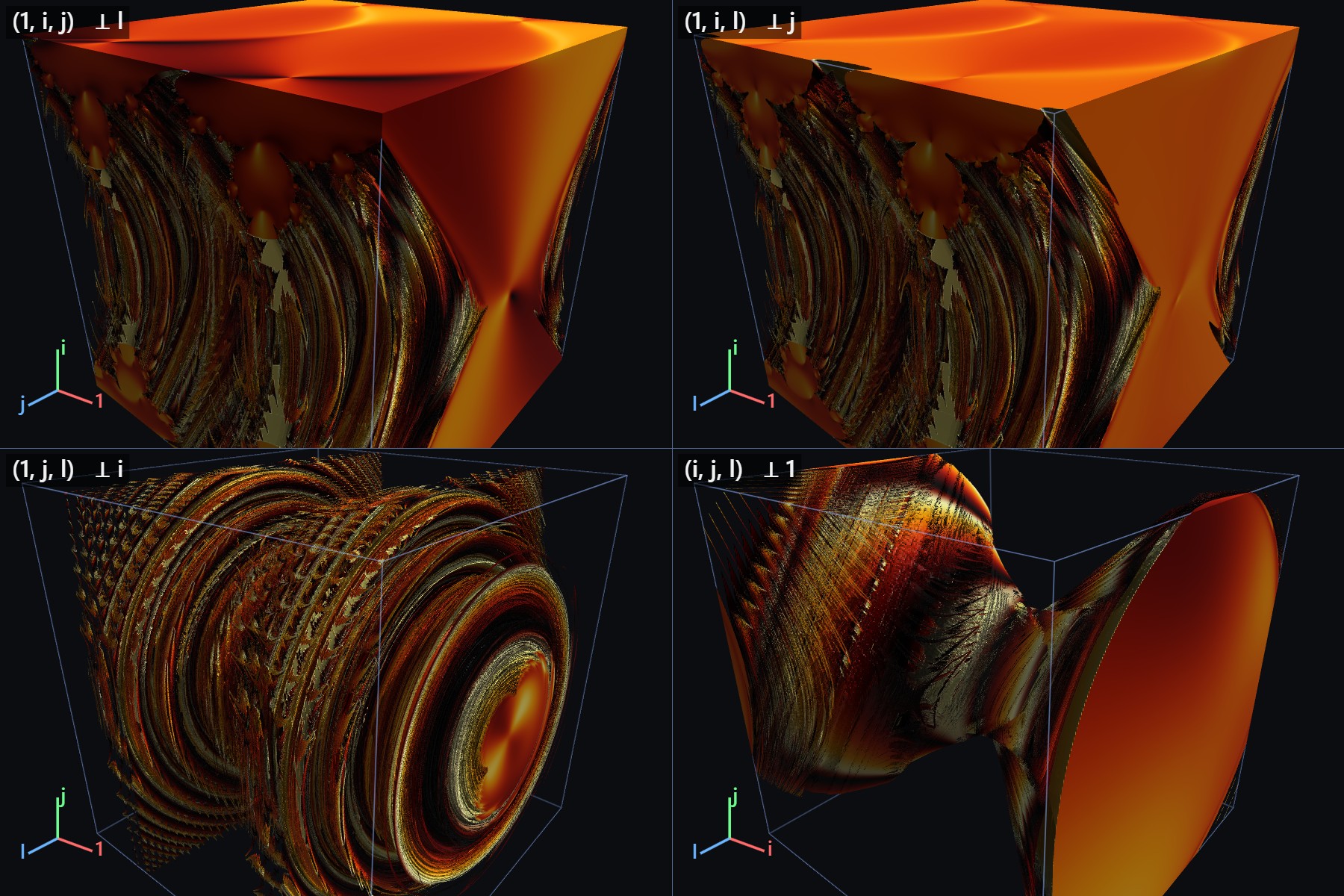}\par\vspace{1mm}
\noindent{\footnotesize \textbf{Map} $z\mapsto c\,\sin z$ $\cdot$ \textbf{Algebra} Split-octonions ($i^2=+1$) (Cayley--Dickson, $\gamma=(1,-1,-1)$) $\cdot$ $c=(0.976,\,\allowbreak -0.091,\,\allowbreak -0.052,\,\allowbreak -0.054,\,\allowbreak -0.081,\,\allowbreak -0.157,\,\allowbreak -0.129,\,\allowbreak 0.021)$ $\cdot$ \textbf{Frame} $X,Y,Z,N=1,i,j,l$ $\cdot$ half-extent 3 $\cdot$ 12 iterations $\cdot$ \textbf{Coloring} 8D orbit trap (point), gradient Fire}\par\vspace{1mm}
\noindent{\footnotesize\begin{tabular}{@{}lllll@{}}\toprule Slice (X, Y, Z) & Body [95\%] & Exact ($J\geq0.9995$) & Approximate ($0.97\leq J<0.9995$) & Dim\\\midrule
($1, i, j$) $\perp l$ & 88.0\% {\scriptsize[87.8, 88.2]} & $P$ & -- & 4 \\
($1, i, l$) $\perp j$ & 87.9\% {\scriptsize[87.7, 88.1]} & $P$ & -- & 4 \\
($1, j, l$) $\perp i$ & 62.0\% {\scriptsize[61.7, 62.3]} & $P$ & -- & 4 \\
($i, j, l$) $\perp 1$ & 43.0\% {\scriptsize[42.6, 43.3]} & $P$ & -- & 4 \\
\bottomrule\end{tabular}}\par\vspace{2mm}
With $i^2=+1$ the plane of $1$ and $i$, which carries the dominant part of the constant, becomes split-complex, and along it the sine behaves like a hyperbolic function. The quarters $(1,i,j)$ and $(1,i,l)$ are filled cubes (body $88\%$), whose faces show the bulbs of the sine family on one side and a fibrous, wavy structure on the other; the quarter $(1,j,l)$, where $j$ and $l$ are elliptic, shows the nested shells and rings of the sine family ($62\%$), and the imaginary quarter an hourglass, with a large flat face at one end ($43\%$).

The point reflection $P$ is exact in all quarters, as for every odd map; the mirror $M_y$ and the half-turn about the $i$-axis come close to the threshold in the two upper quarters ($J=0.967$--$0.969$) but are not symmetries. The orbits lie in a subalgebra of dimension $4$. This plate and Plate~F08 show the same exchange of roles between the complex and the hyperbolic planes.
\clearpage

\subsection{Part M. Four-dimensional mappings}\label{atl:part:M}
In this part a four-dimensional mapping of HyperFractal Studio is applied to the sample point before the iteration. As explained in Remark~\ref{atl:rem:mappings}, these transformations act only on the coordinates $1,i,j,k$, so that the $l$-directions are left untouched and the mapping may be invisible in some quarters; for the same reason they break, artificially, the symmetries of the root planes of the Okubo algebra. These plates do not respect the algebra, and they are not measured in this version: the tables show the render-time measurements of the first version ($2400$ points, threshold $0.97$), which are only indicative. The base bodies are the octonionic Julia set of Plate~A01 and the Okubonic Julia set with the generic constant of Plate~B07.

\subsubsection{M01. Tesseract kaleidoscope on the Okubo body}\label{atl:pl:M01}\mbox{}\par\vspace{0.5mm}
\noindent\includegraphics[width=\textwidth]{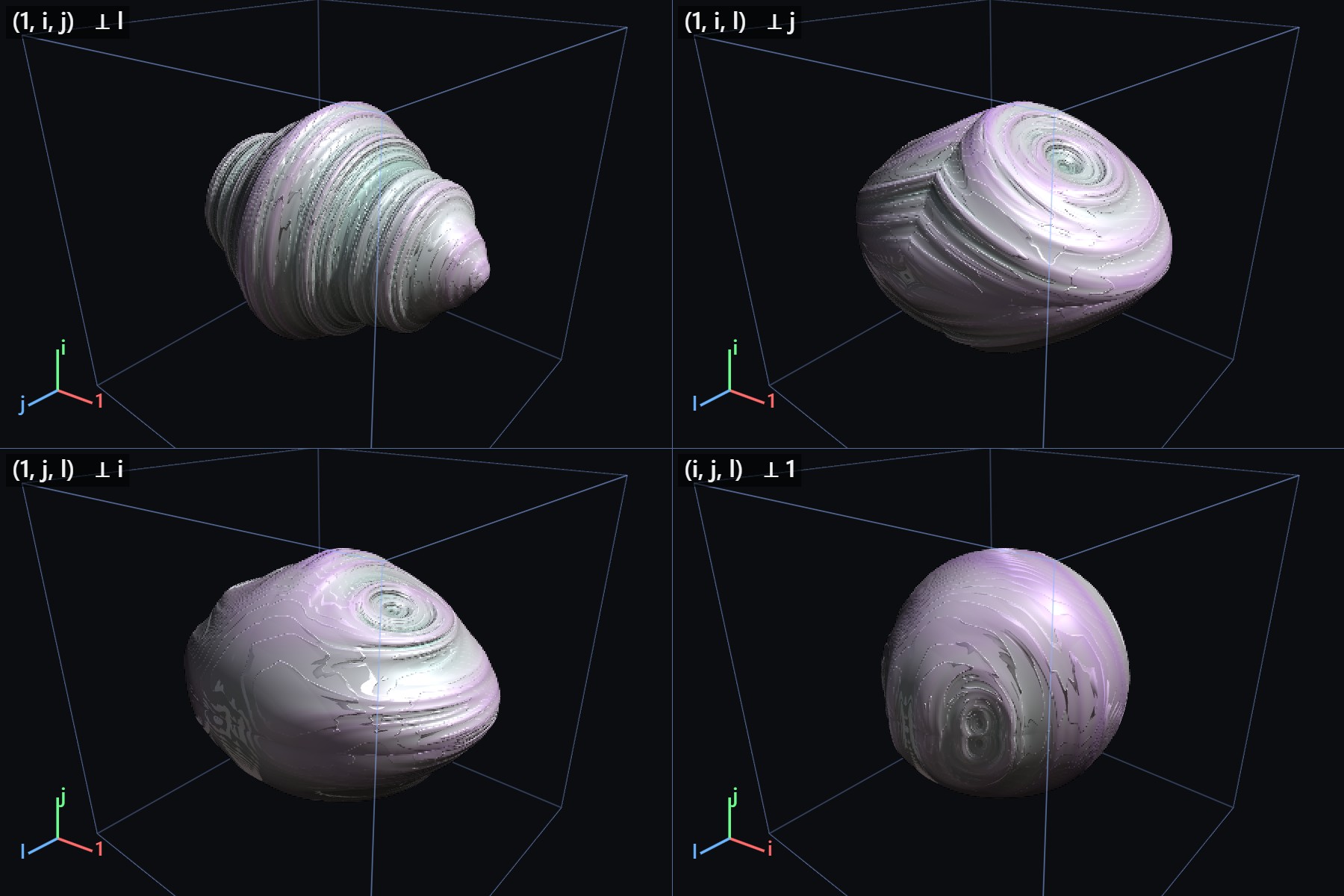}\par\vspace{1mm}
\noindent{\footnotesize \textbf{Map} $z\mapsto z\diamond z+c$ $\cdot$ \textbf{Algebra} Okubo algebra (Petersson, $\gamma=(-1,-1,-1)$, $\theta=120^\circ$) $\cdot$ $c=(-0.397,\,\allowbreak -0.021,\,\allowbreak 0.049,\,\allowbreak 0.016,\,\allowbreak 0.004,\,\allowbreak 0.16,\,\allowbreak -0.113,\,\allowbreak 0.133)$ ($=0.95\,s^*$, $s^*=0.4900$) $\cdot$ \textbf{Frame} $X,Y,Z,N=1,i,j,l$ $\cdot$ half-extent 1.4 $\cdot$ 20 iterations $\cdot$ \textbf{Coloring} 8D orbit trap (point), gradient Aurora $\cdot$ \textbf{Mapping} \texttt{Kaleidoscope4D}}\par\vspace{1mm}
\noindent{\footnotesize\begin{tabular}{@{}lrlll@{}}\toprule Slice (X, Y, Z) & Body & Mirrors & Rotations ($J\geq0.97$) & Orbit span\\\midrule
($1, i, j$) $\perp l$ & 29\% & xyz & $R_x$: 180; $R_y$: 180; $R_z$: 180 & 8 \\
($1, i, l$) $\perp j$ & 41\% & xyz & $R_x$: 180; $R_y$: 180; $R_z$: 180 & 8 \\
($1, j, l$) $\perp i$ & 41\% & xy & $R_z$: 180 & 8 \\
($i, j, l$) $\perp 1$ & 33\% & xy & $R_z$: 180 & 8 \\
\bottomrule\end{tabular}}\par{\scriptsize Render-time measurements of version~1 ($2400$ points, threshold $0.97$): indicative only (Remark~\ref{atl:rem:mappings}).}\par\vspace{2mm}
The tesseract kaleidoscope folds the coordinates $1,i,j,k$ of the sample onto their positive side, $x_m\mapsto|x_m|$, before the iteration: each quarter shows the points whose folded image lies in the Okubonic body of Plate~M06 (generic constant, frame $F_A$), and it acquires by construction the mirrors in the folded coordinates. The pictures differ little from Plate~M06, because the base body is already close to these reflections; the fold is visible mainly as a step on the left side of the quarter $(1,i,l)$.

The plates of Part~M are not measured (Remark~\ref{atl:rem:mappings}); the render-time check, on $2400$ points with threshold $0.97$ and therefore only indicative, finds the expected mirrors: all three in the quarter $(1,i,j)$, the first two in the quarters $(1,i,l)$, $(1,j,l)$ and $(i,j,l)$, whose third axis $l$ is not touched by the mapping. In $(1,i,l)$ it also reports the mirror in $l$ at $J=0.97$, inherited from a near-symmetry of the base body ($0.98$ in Plate~M06), not an exact one.

The mirrors are imposed by the mapping, not by the algebra: for a generic constant the stabilizer is a torus whose root planes are in general position (Plate~R09), and no coordinate symmetry is to be expected from the dynamics.
\clearpage

\subsubsection{M02. Tesseract kaleidoscope on the octonion body}\label{atl:pl:M02}\mbox{}\par\vspace{0.5mm}
\noindent\includegraphics[width=\textwidth]{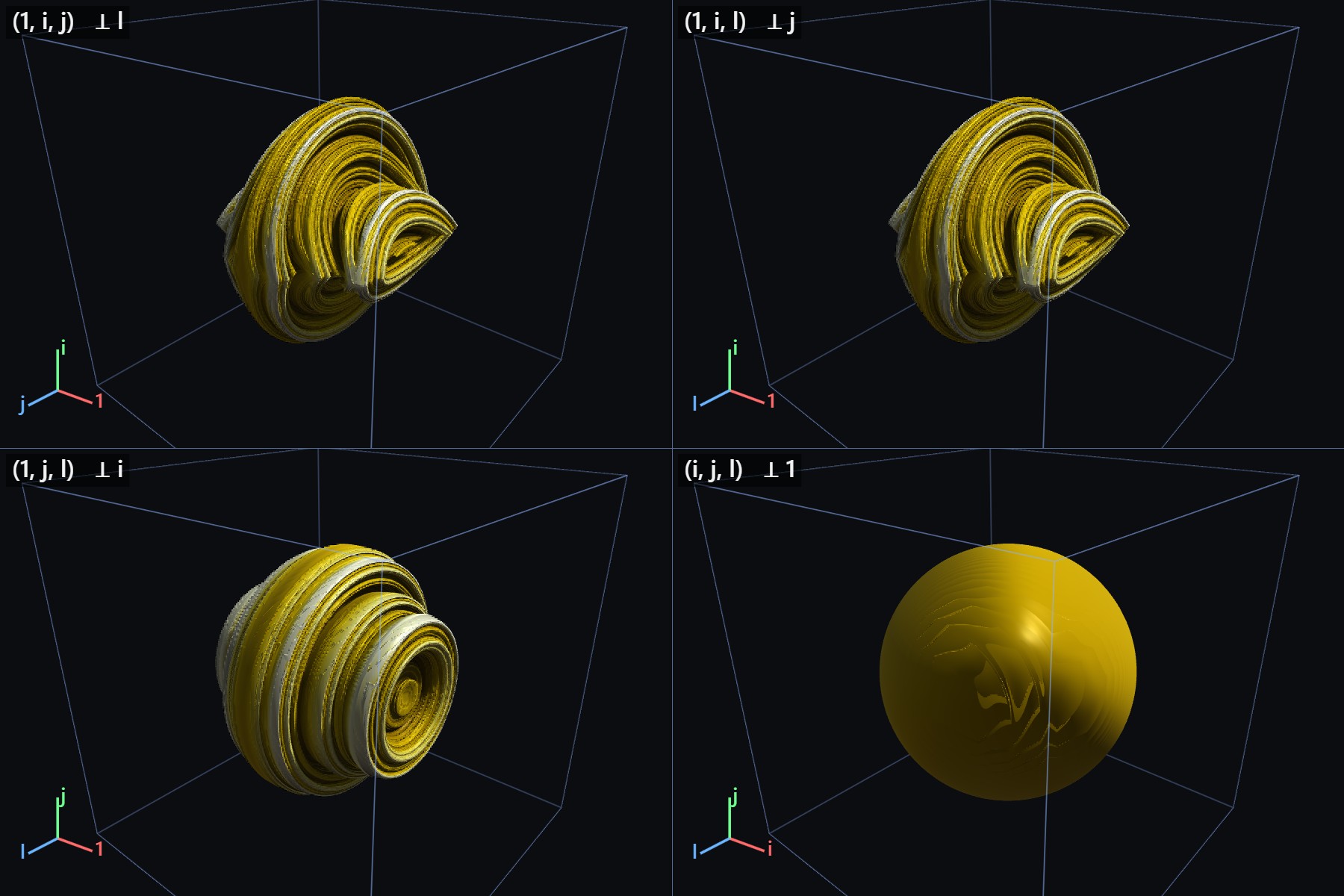}\par\vspace{1mm}
\noindent{\footnotesize \textbf{Map} $z\mapsto z\diamond z+c$ $\cdot$ \textbf{Algebra} Octonions (Cayley--Dickson, $\gamma=(-1,-1,-1)$) $\cdot$ $c=(-0.231,\,\allowbreak 0.62,\,\allowbreak 0,\,\allowbreak 0,\,\allowbreak 0,\,\allowbreak 0,\,\allowbreak 0,\,\allowbreak 0)$ ($=0.97\,s^*$, $s^*=0.6824$) $\cdot$ \textbf{Frame} $X,Y,Z,N=1,i,j,l$ $\cdot$ half-extent 1.4 $\cdot$ 20 iterations $\cdot$ \textbf{Coloring} 8D orbit trap (point), gradient Gold $\cdot$ \textbf{Mapping} \texttt{Kaleidoscope4D}}\par\vspace{1mm}
\noindent{\footnotesize\begin{tabular}{@{}lrlll@{}}\toprule Slice (X, Y, Z) & Body & Mirrors & Rotations ($J\geq0.97$) & Orbit span\\\midrule
($1, i, j$) $\perp l$ & 13\% & xyz & $R_x$: 180; $R_y$: 180; $R_z$: 180 & 3 \\
($1, i, l$) $\perp j$ & 13\% & xyz & $R_x$: 180; $R_y$: 180; $R_z$: 180 & 3 \\
($1, j, l$) $\perp i$ & 23\% & xyz & $R_x$: all; $R_y$: 180; $R_z$: 180 & 3 \\
($i, j, l$) $\perp 1$ & 35\% & xyz & $R_x$: all; $R_y$: all; $R_z$: all & 3 \\
\bottomrule\end{tabular}}\par{\scriptsize Render-time measurements of version~1 ($2400$ points, threshold $0.97$): indicative only (Remark~\ref{atl:rem:mappings}).}\par\vspace{2mm}
The same kaleidoscope on the octonionic body of Plate~A01 ($c\in\CC$). The quarters $(1,i,j)$ and $(1,i,l)$ are identical: a compact shell with spiral ridges. The quarter $(1,j,l)$ shows a ringed body and the imaginary quarter a smooth ball, as in Plate~A01.

The fold changes nothing in the two lower quarters. By Proposition~\ref{atl:prop:equivariance} the base body is invariant under $O(6)$, which contains the reflection of $j$; and by Remark~\ref{atl:rem:minus} $z\mapsto-z$ is a symmetry, whose composition with a half-turn of $O(6)$ reflects the coordinate $1$ (or $i$) alone within these sections. So the folded point always lies in the body exactly when the original one does. In the quarters $(1,i,j)$ and $(1,i,l)$ the base body has the mirror in the third axis and the half-turn about it, and the fold adds the mirrors in $1$ and $i$. The render-time check (indicative only; Part~M is not measured) agrees: three mirrors in the upper quarters, a body of revolution about $1$ in $(1,j,l)$, a ball in $(i,j,l)$.

The octonionic plate shows the kaleidoscope acting on an already symmetric object, the Okubonic Plate~M01 shows it imposing mirrors that the dynamics of a generic constant does not have.
\clearpage

\subsubsection{M03. Inversion in a 3-sphere, Okubo}\label{atl:pl:M03}\mbox{}\par\vspace{0.5mm}
\noindent\includegraphics[width=\textwidth]{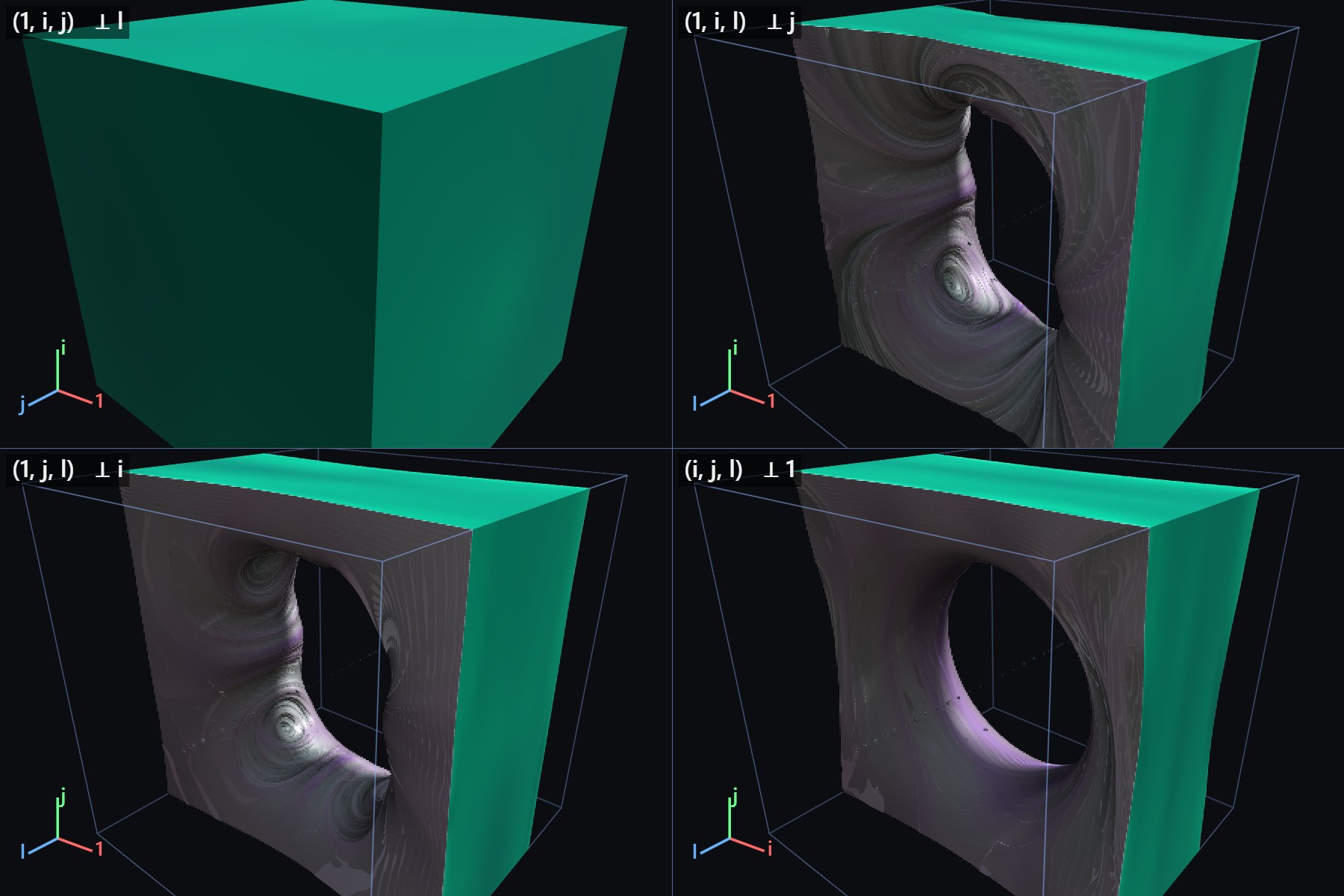}\par\vspace{1mm}
\noindent{\footnotesize \textbf{Map} $z\mapsto z\diamond z+c$ $\cdot$ \textbf{Algebra} Okubo algebra (Petersson, $\gamma=(-1,-1,-1)$, $\theta=120^\circ$) $\cdot$ $c=(-0.397,\,\allowbreak -0.021,\,\allowbreak 0.049,\,\allowbreak 0.016,\,\allowbreak 0.004,\,\allowbreak 0.16,\,\allowbreak -0.113,\,\allowbreak 0.133)$ ($=0.95\,s^*$, $s^*=0.4900$) $\cdot$ \textbf{Frame} $X,Y,Z,N=1,i,j,l$ $\cdot$ half-extent 1.8 $\cdot$ 20 iterations $\cdot$ \textbf{Coloring} 8D orbit trap (point), gradient Aurora $\cdot$ \textbf{Mapping} \texttt{Inversion4D} (radius$=1$)}\par\vspace{1mm}
\noindent{\footnotesize\begin{tabular}{@{}lrlll@{}}\toprule Slice (X, Y, Z) & Body & Mirrors & Rotations ($J\geq0.97$) & Orbit span\\\midrule
($1, i, j$) $\perp l$ & 67\% & -- & $R_x$: 37 & 8 \\
($1, i, l$) $\perp j$ & 28\% & z & $R_z$: 180 & 8 \\
($1, j, l$) $\perp i$ & 29\% & -- & -- & 8 \\
($i, j, l$) $\perp 1$ & 25\% & -- & -- & 8 \\
\bottomrule\end{tabular}}\par{\scriptsize Render-time measurements of version~1 ($2400$ points, threshold $0.97$): indicative only (Remark~\ref{atl:rem:mappings}).}\par\vspace{2mm}
The inversion in the unit 3-sphere of the coordinates $1,i,j,k$ exchanges the inside and the outside of the sphere. In the quaternionic quarter the visible faces of the cube (half-extent $1.8$) are entirely filled: the far points are sent close to the origin, inside the Okubonic body, while the points near the centre are sent far away and escape, so that the empty part of the quarter ($33\%$ at render time) is hidden inside the cube.

In the quarters containing $l$ only the first four coordinates are inverted (Remark~\ref{atl:rem:mappings}), and the pictures show thick slabs bounded by the faces of the cube and pierced by a tunnel around the $l$-axis, where the inverted coordinates are small and their images large. The walls of the tunnel carry vortices in the quarters $(1,i,l)$ and $(1,j,l)$; in the imaginary quarter they are smooth, with an oval mouth. The render-time check (indicative only) finds no symmetry except an approximate mirror in $l$ in the quarter $(1,i,l)$.

The mapping is not an automorphism of the algebra, and the pictures are no longer sections of a Julia set but of its inverse image under the mapping.
\clearpage

\subsubsection{M04. Inversion in a 3-sphere, octonions}\label{atl:pl:M04}\mbox{}\par\vspace{0.5mm}
\noindent\includegraphics[width=\textwidth]{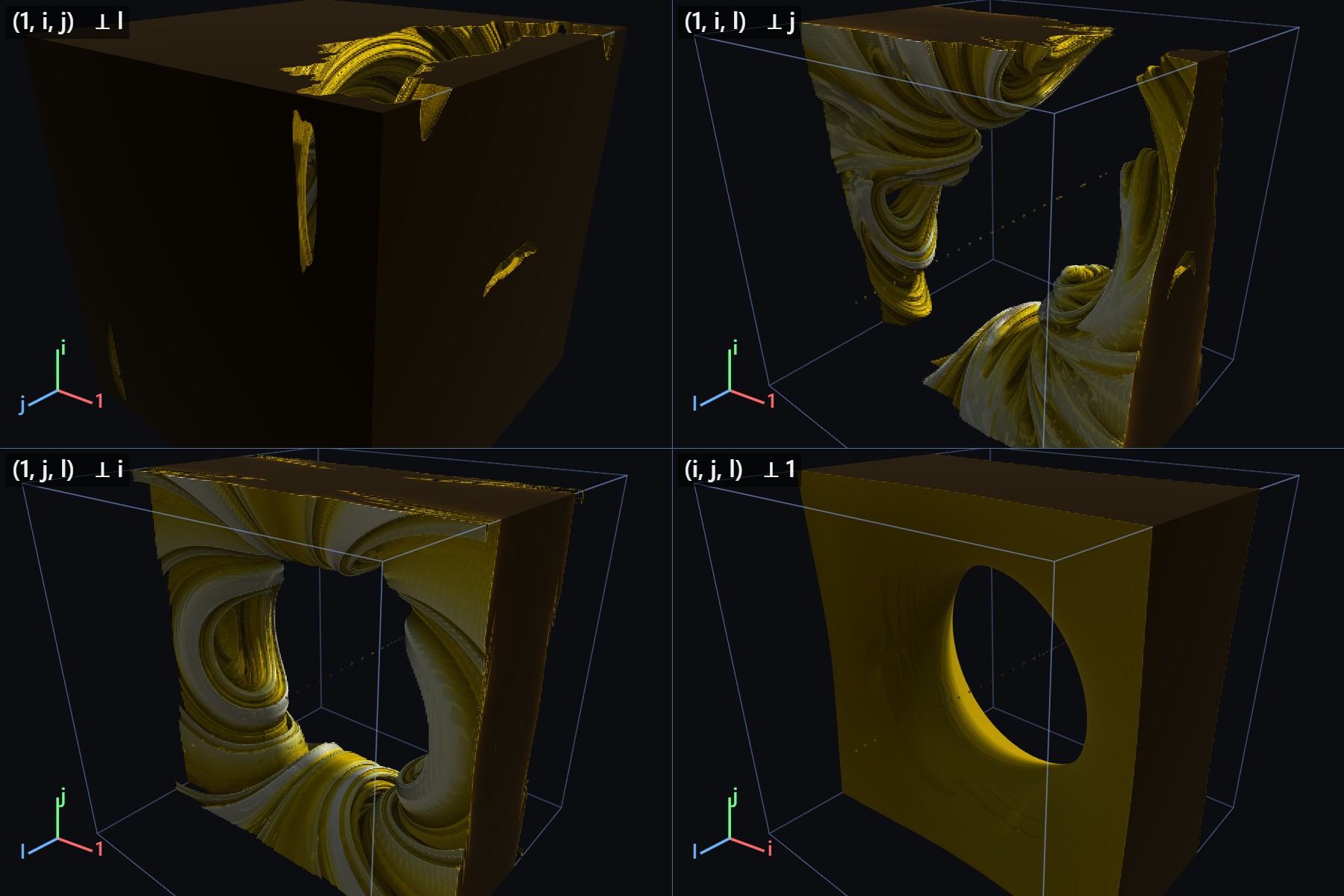}\par\vspace{1mm}
\noindent{\footnotesize \textbf{Map} $z\mapsto z\diamond z+c$ $\cdot$ \textbf{Algebra} Octonions (Cayley--Dickson, $\gamma=(-1,-1,-1)$) $\cdot$ $c=(-0.231,\,\allowbreak 0.62,\,\allowbreak 0,\,\allowbreak 0,\,\allowbreak 0,\,\allowbreak 0,\,\allowbreak 0,\,\allowbreak 0)$ ($=0.97\,s^*$, $s^*=0.6824$) $\cdot$ \textbf{Frame} $X,Y,Z,N=1,i,j,l$ $\cdot$ half-extent 1.8 $\cdot$ 20 iterations $\cdot$ \textbf{Coloring} 8D orbit trap (point), gradient Gold $\cdot$ \textbf{Mapping} \texttt{Inversion4D} (radius$=1$)}\par\vspace{1mm}
\noindent{\footnotesize\begin{tabular}{@{}lrlll@{}}\toprule Slice (X, Y, Z) & Body & Mirrors & Rotations ($J\geq0.97$) & Orbit span\\\midrule
($1, i, j$) $\perp l$ & 48\% & z & $R_z$: 180 & 3 \\
($1, i, l$) $\perp j$ & 15\% & z & $R_z$: 180 & 3 \\
($1, j, l$) $\perp i$ & 14\% & xyz & $R_x$: 180; $R_y$: 180; $R_z$: 180 & 3 \\
($i, j, l$) $\perp 1$ & 27\% & xyz & $R_x$: 180; $R_y$: 180; $R_z$: all & 3 \\
\bottomrule\end{tabular}}\par{\scriptsize Render-time measurements of version~1 ($2400$ points, threshold $0.97$): indicative only (Remark~\ref{atl:rem:mappings}).}\par\vspace{2mm}
The same inversion on the octonionic body ($c\in\CC$). The quaternionic quarter is an almost full cube, with a few openings on its faces where the swirling surface of the complement shows through. The quarter $(1,j,l)$ is a square frame with swirling walls; the quarter $(1,i,l)$ is more fragmented, two twisted walls separated by wide gaps; the imaginary quarter is a slab with a smooth circular hole around the $l$-axis.

The hole is round because the octonionic body is a ball in the imaginary quarter (Proposition~\ref{atl:prop:imaginary}) and the inversion of the coordinates $i,j$ commutes with the rotations of their plane, while $l$ is untouched: the quarter is a body of revolution about $l$, and the render-time check (indicative only) finds all the rotations about that axis. In the quarter $(1,j,l)$ it finds the three mirrors, in the two upper quarters the mirror in the third axis. In the Okubonic Plate~M03, with a generic constant, the mouth of the tunnel is oval and no rotation is found.
\clearpage

\subsubsection{M05. Sphere fold, Okubo}\label{atl:pl:M05}\mbox{}\par\vspace{0.5mm}
\noindent\includegraphics[width=\textwidth]{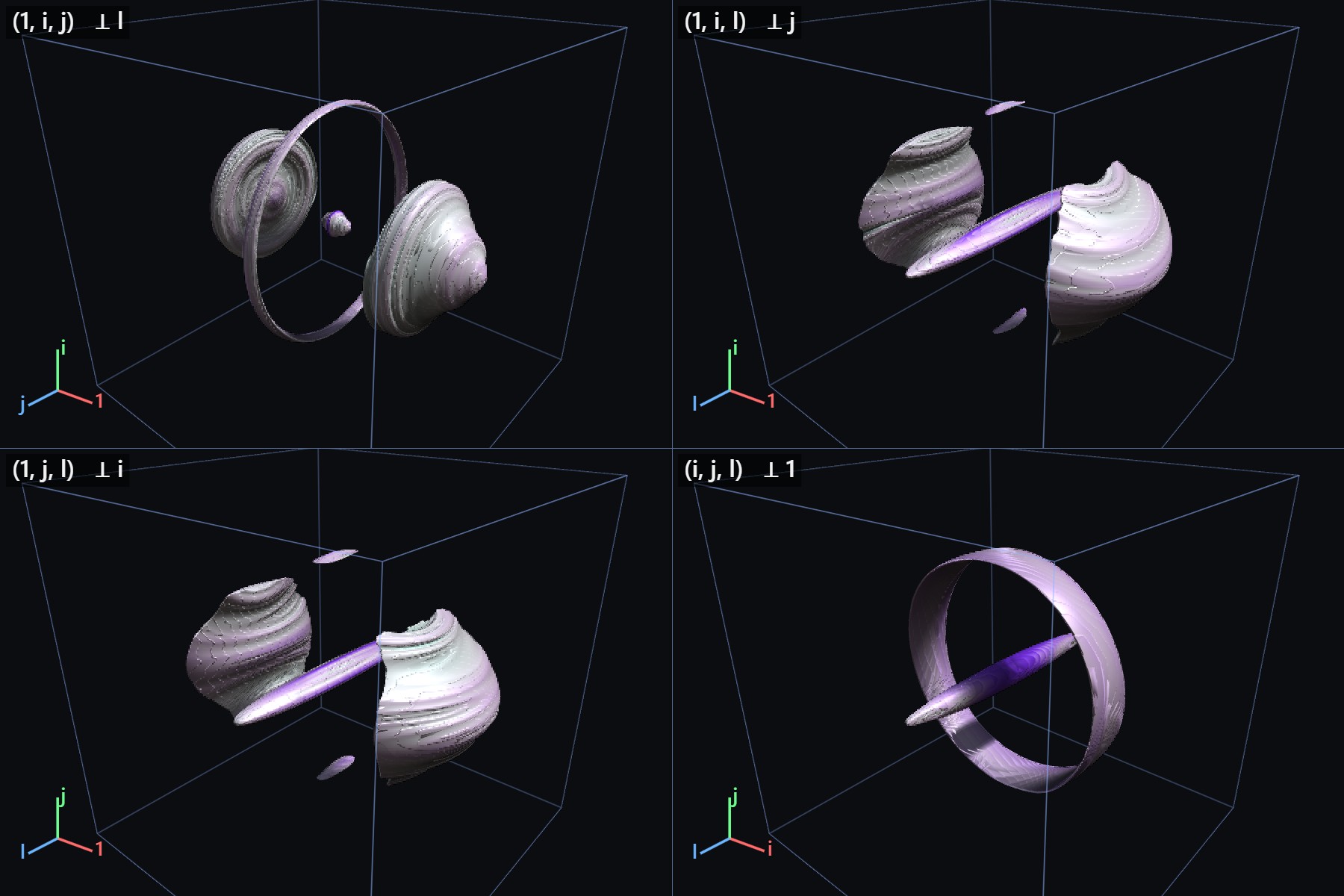}\par\vspace{1mm}
\noindent{\footnotesize \textbf{Map} $z\mapsto z\diamond z+c$ $\cdot$ \textbf{Algebra} Okubo algebra (Petersson, $\gamma=(-1,-1,-1)$, $\theta=120^\circ$) $\cdot$ $c=(-0.397,\,\allowbreak -0.021,\,\allowbreak 0.049,\,\allowbreak 0.016,\,\allowbreak 0.004,\,\allowbreak 0.16,\,\allowbreak -0.113,\,\allowbreak 0.133)$ ($=0.95\,s^*$, $s^*=0.4900$) $\cdot$ \textbf{Frame} $X,Y,Z,N=1,i,j,l$ $\cdot$ half-extent 1.4 $\cdot$ 20 iterations $\cdot$ \textbf{Coloring} 8D orbit trap (point), gradient Aurora $\cdot$ \textbf{Mapping} \texttt{SphereFold4D} (rmin$=0.3$, rfix$=0.9$)}\par\vspace{1mm}
\noindent{\footnotesize\begin{tabular}{@{}lrlll@{}}\toprule Slice (X, Y, Z) & Body & Mirrors & Rotations ($J\geq0.97$) & Orbit span\\\midrule
($1, i, j$) $\perp l$ & 5\% & -- & -- & 8 \\
($1, i, l$) $\perp j$ & 7\% & -- & -- & 8 \\
($1, j, l$) $\perp i$ & 7\% & -- & -- & 8 \\
($i, j, l$) $\perp 1$ & 1\% & -- & -- & 8 \\
\bottomrule\end{tabular}}\par{\scriptsize Render-time measurements of version~1 ($2400$ points, threshold $0.97$): indicative only (Remark~\ref{atl:rem:mappings}).}\par\vspace{2mm}
The sphere fold of the Mandelbox, in the four coordinates $1,i,j,k$, inverts the shell between the minimum radius $0.3$ and the fixed radius $0.9$ in the sphere of radius $0.9$, and scales up the inside of the minimum radius by the factor $9$. The Okubonic body is cut into pieces: in the quaternionic quarter two shells face each other across a thin ring, with a small body at the centre; in the quarters containing $l$ two shells are separated by a thin disc, with a few small fragments; in the imaginary quarter a ring surrounds a disc.

The body shares drop to a few percent ($1\%$ to $7\%$ at render time, indicative only), and no symmetry is found. This plate shows how a mapping of the sample space can split the sections of a filled Julia set into separate pieces.
\clearpage

\subsubsection{M06. Twist through \texorpdfstring{$W$}{W}, Okubo: an invisible mapping}\label{atl:pl:M06}\mbox{}\par\vspace{0.5mm}
\noindent\includegraphics[width=\textwidth]{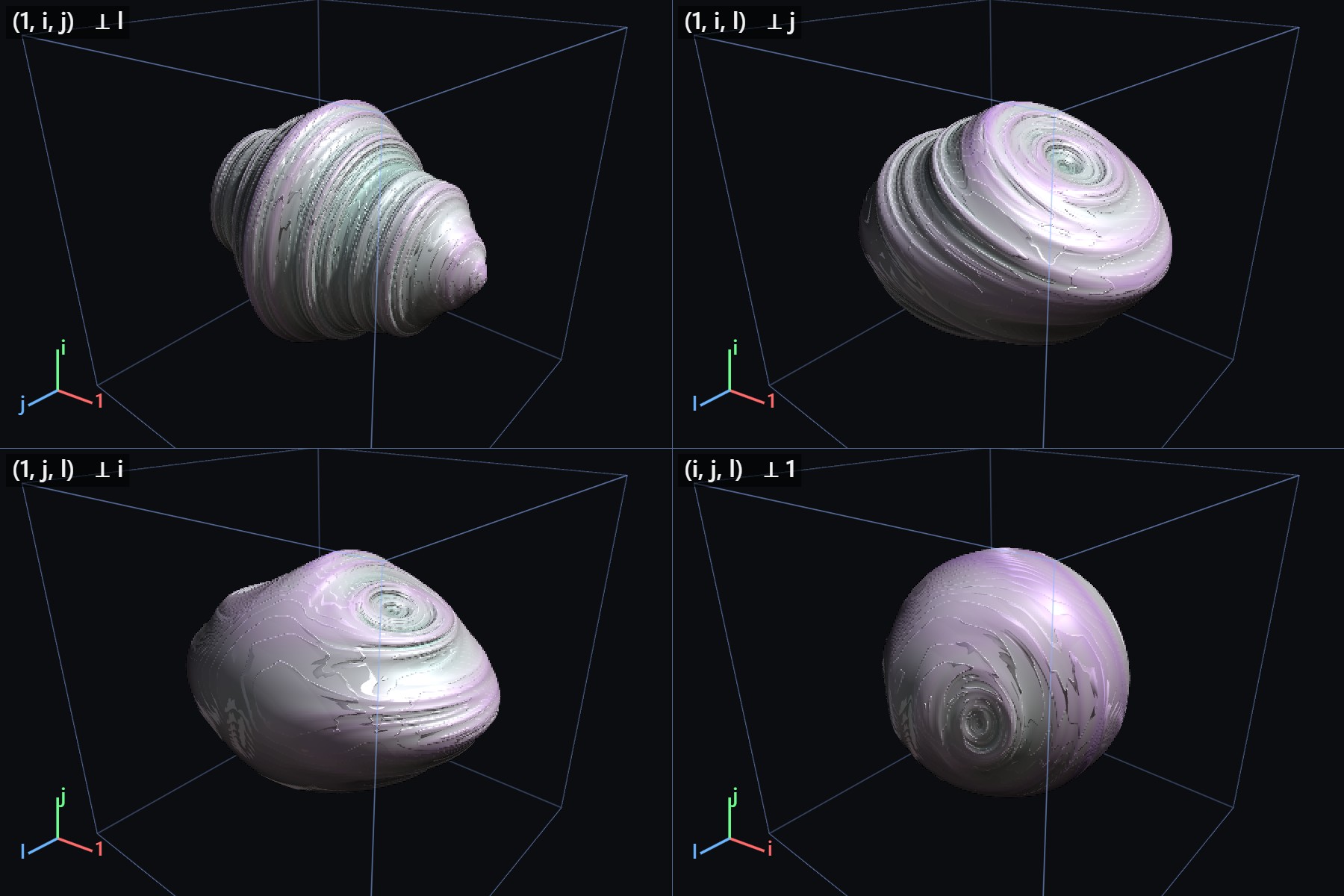}\par\vspace{1mm}
\noindent{\footnotesize \textbf{Map} $z\mapsto z\diamond z+c$ $\cdot$ \textbf{Algebra} Okubo algebra (Petersson, $\gamma=(-1,-1,-1)$, $\theta=120^\circ$) $\cdot$ $c=(-0.397,\,\allowbreak -0.021,\,\allowbreak 0.049,\,\allowbreak 0.016,\,\allowbreak 0.004,\,\allowbreak 0.16,\,\allowbreak -0.113,\,\allowbreak 0.133)$ ($=0.95\,s^*$, $s^*=0.4900$) $\cdot$ \textbf{Frame} $X,Y,Z,N=1,i,j,l$ $\cdot$ half-extent 1.4 $\cdot$ 20 iterations $\cdot$ \textbf{Coloring} 8D orbit trap (point), gradient Aurora $\cdot$ \textbf{Mapping} \texttt{Twist4D} (plane$=0$, by$=3$, amount$=120$)}\par\vspace{1mm}
\noindent{\footnotesize\begin{tabular}{@{}lrlll@{}}\toprule Slice (X, Y, Z) & Body & Mirrors & Rotations ($J\geq0.97$) & Orbit span\\\midrule
($1, i, j$) $\perp l$ & 28\% & -- & -- & 8 \\
($1, i, l$) $\perp j$ & 39\% & z & $R_z$: 180 & 8 \\
($1, j, l$) $\perp i$ & 40\% & -- & -- & 8 \\
($i, j, l$) $\perp 1$ & 34\% & x & $R_x$: 180 & 8 \\
\bottomrule\end{tabular}}\par{\scriptsize Render-time measurements of version~1 ($2400$ points, threshold $0.97$): indicative only (Remark~\ref{atl:rem:mappings}).}\par\vspace{2mm}
The twist through the fourth dimension rotates the plane $XY$ of the first two coordinates by an angle proportional to the coordinate $W$, which is $k$ in the eight-dimensional sample. All quarters of the frame $F_A$ have $k=0$, so the twist is the identity on them (Remark~\ref{atl:rem:mappings}), and this plate shows the untransformed Okubonic body with the generic constant in the frame $F_A$.

It serves as the reference for the other Okubonic plates of this part: a ribbed shell in the quaternionic quarter, an oblate body with a spiral vortex in $(1,i,l)$, a rounded body with a vortex on top in $(1,j,l)$ and a roundish body with a vortex in the imaginary quarter. The render-time check (indicative only, since Part~M is not measured) finds only approximate mirrors, in $l$ for $(1,i,l)$ and in $i$ for $(i,j,l)$ ($J=0.97$--$0.98$), and orbits spanning all eight dimensions. No exact symmetry is expected in this frame: for a generic constant the stabilizer is a torus whose root planes are in general position (Plate~R09).
\clearpage

\subsubsection{M07. Waves, para-octonions}\label{atl:pl:M07}\mbox{}\par\vspace{0.5mm}
\noindent\includegraphics[width=\textwidth]{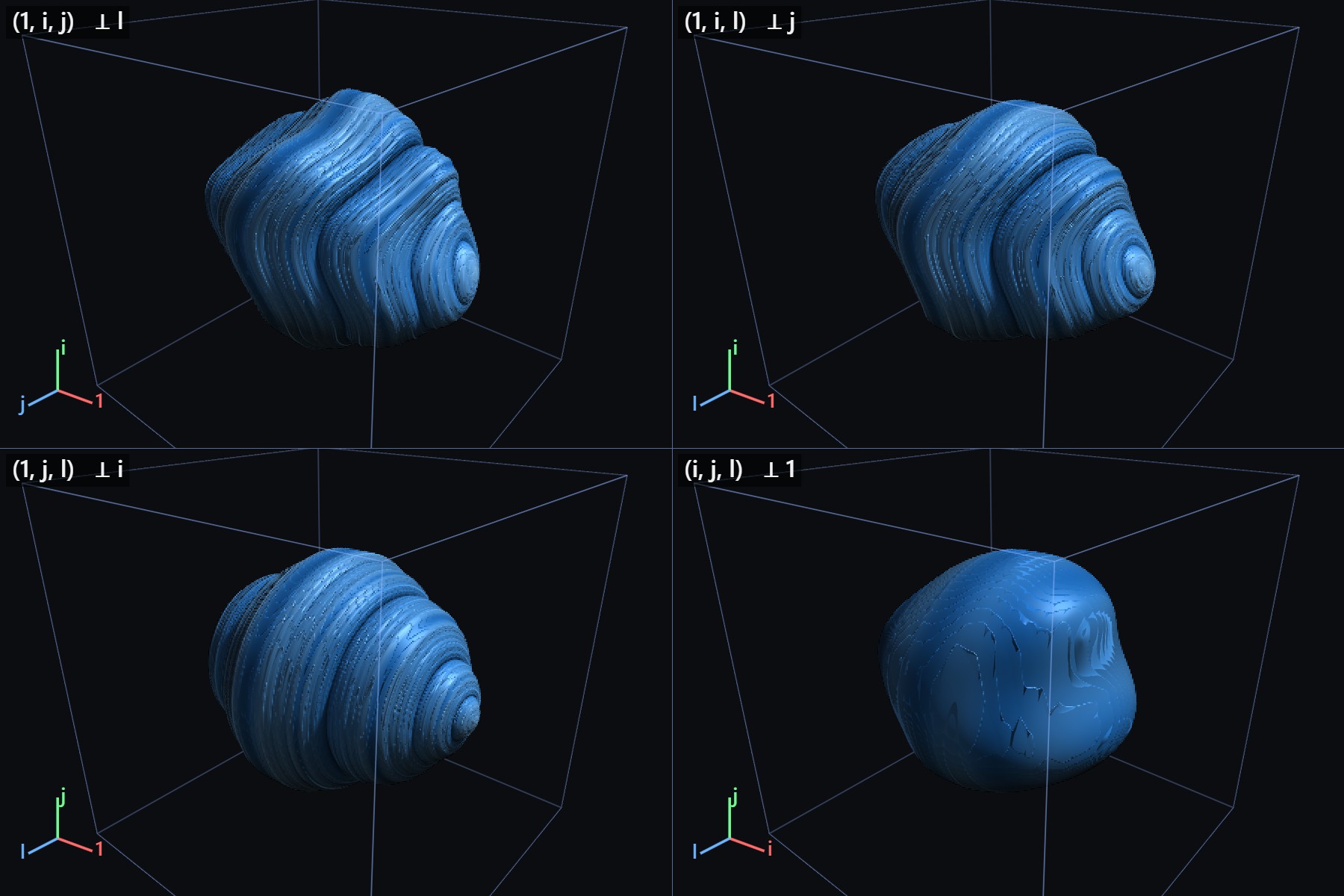}\par\vspace{1mm}
\noindent{\footnotesize \textbf{Map} $z\mapsto z\diamond z+c$ $\cdot$ \textbf{Algebra} Para-octonions (para-Hurwitz $\bar x\bar y$, $\gamma=(-1,-1,-1)$) $\cdot$ $c=(-0.235,\,\allowbreak -0.013,\,\allowbreak 0.029,\,\allowbreak 0.01,\,\allowbreak 0.002,\,\allowbreak 0.095,\,\allowbreak -0.067,\,\allowbreak 0.079)$ ($=0.95\,s^*$, $s^*=0.2901$) $\cdot$ \textbf{Frame} $X,Y,Z,N=1,i,j,l$ $\cdot$ half-extent 1.4 $\cdot$ 20 iterations $\cdot$ \textbf{Coloring} 8D orbit trap (point), gradient Deep Blue $\cdot$ \textbf{Mapping} \texttt{Waves4D} (amplitude$=0.08$, frequency$=6$)}\par\vspace{1mm}
\noindent{\footnotesize\begin{tabular}{@{}lrlll@{}}\toprule Slice (X, Y, Z) & Body & Mirrors & Rotations ($J\geq0.97$) & Orbit span\\\midrule
($1, i, j$) $\perp l$ & 33\% & -- & -- & 3 \\
($1, i, l$) $\perp j$ & 33\% & z & $R_z$: 180 & 3 \\
($1, j, l$) $\perp i$ & 32\% & z & $R_x$: 60,37; $R_z$: 180 & 3 \\
($i, j, l$) $\perp 1$ & 32\% & z & $R_z$: 180 & 3 \\
\bottomrule\end{tabular}}\par{\scriptsize Render-time measurements of version~1 ($2400$ points, threshold $0.97$): indicative only (Remark~\ref{atl:rem:mappings}).}\par\vspace{2mm}
The waves displace each of the coordinates $1,i,j,k$ by a sine of the next one ($x_0$ by $x_1$, $x_1$ by $x_2$, $x_2$ by $x_3$, $x_3$ by $x_0$), with amplitude $0.08$ and frequency $6$. The base body is the para-octonionic Julia set of the generic constant, in the frame $F_A$. The effect is a gentle rippling: the three quarters through $1$ show ribbed shells whose ridges are wavy, and the imaginary quarter a lumpy, rounded body. The mapping preserves the overall shape.

The render-time check (on $2400$ points with threshold $0.97$; Part~M is not measured) reports in the quarter $(1,j,l)$ the rotations about the $1$-axis by $37^\circ$ ($J=0.98$) and $60^\circ$ ($0.97$), but not those by $72^\circ$, $90^\circ$ and $120^\circ$. An exact continuous symmetry would pass every angle, and this sample cannot distinguish an exact symmetry from an approximate one: we record a near-symmetry of unknown size, not a symmetry. A plausible origin is the following. The para-octonionic map is equivariant under the rotations fixing $1$ and $\im c$ (Proposition~\ref{atl:prop:equivariance}); the projection of $\im c$ onto the plane $\langle j,l\rangle$ is only about a fifth of $\im c$, so the base body is nearly invariant under the rotations of that plane; and in this quarter the waves only add the small coordinates $x_1=0.08\sin 6x_2$ and $x_3=0.08\sin 6x_0$. The same check finds the mirror in $l$ in the three quarters containing $l$ ($J=0.99$--$1.00$): $l$ is not touched by the waves, and the $l$-component of the constant is very small.

The orbits of the displaced samples lie in the three-dimensional space spanned by $1$, $\im z_0$ and $\im c$ (Proposition~\ref{atl:prop:span}).
\clearpage

\subsubsection{M08. Tesseract fold, Okubo}\label{atl:pl:M08}\mbox{}\par\vspace{0.5mm}
\noindent\includegraphics[width=\textwidth]{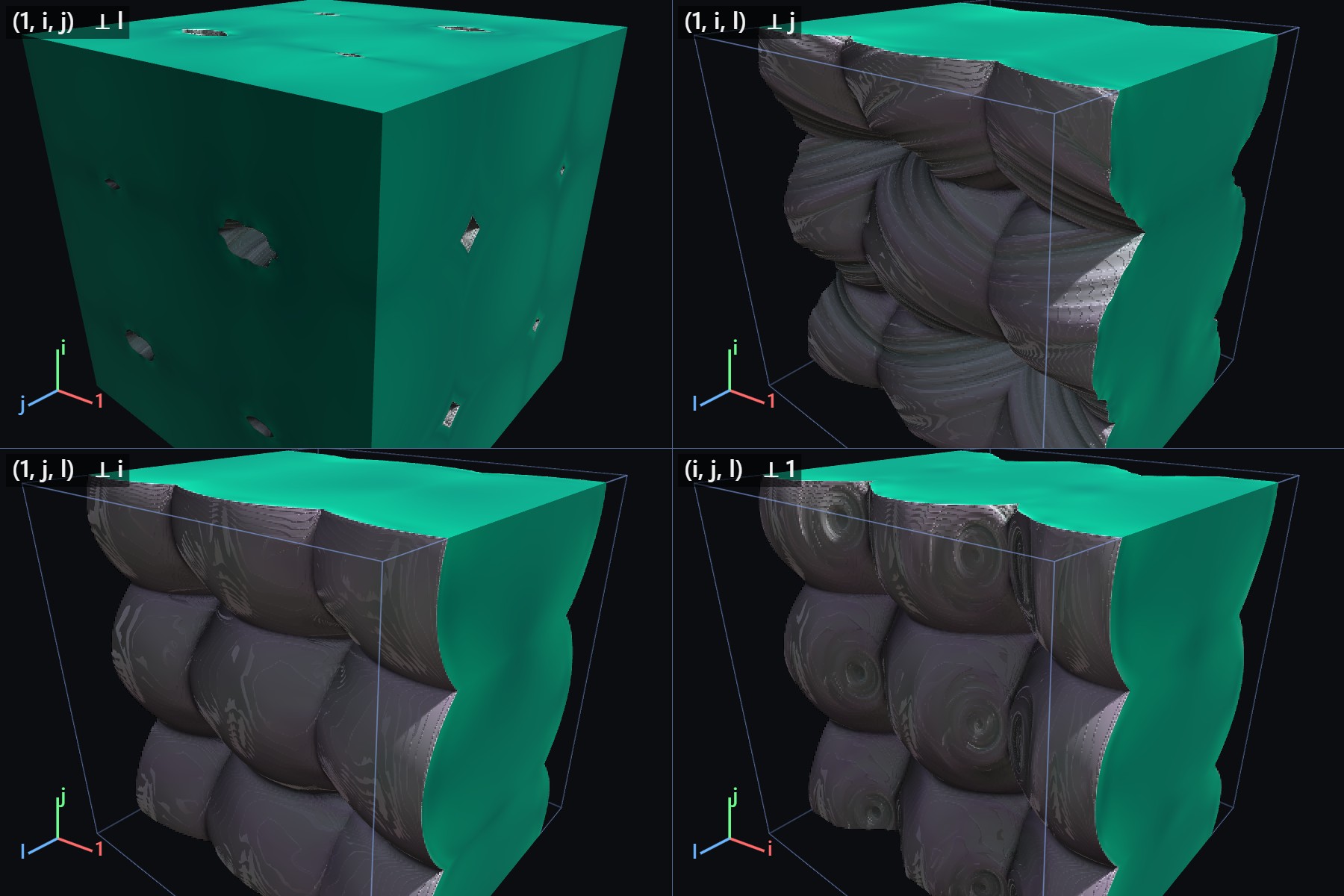}\par\vspace{1mm}
\noindent{\footnotesize \textbf{Map} $z\mapsto z\diamond z+c$ $\cdot$ \textbf{Algebra} Okubo algebra (Petersson, $\gamma=(-1,-1,-1)$, $\theta=120^\circ$) $\cdot$ $c=(-0.397,\,\allowbreak -0.021,\,\allowbreak 0.049,\,\allowbreak 0.016,\,\allowbreak 0.004,\,\allowbreak 0.16,\,\allowbreak -0.113,\,\allowbreak 0.133)$ ($=0.95\,s^*$, $s^*=0.4900$) $\cdot$ \textbf{Frame} $X,Y,Z,N=1,i,j,l$ $\cdot$ half-extent 1.4 $\cdot$ 20 iterations $\cdot$ \textbf{Coloring} 8D orbit trap (point), gradient Aurora $\cdot$ \textbf{Mapping} \texttt{TesseractFold4D} (limit$=0.6$, repeat$=1$, scale$=1$)}\par\vspace{1mm}
\noindent{\footnotesize\begin{tabular}{@{}lrlll@{}}\toprule Slice (X, Y, Z) & Body & Mirrors & Rotations ($J\geq0.97$) & Orbit span\\\midrule
($1, i, j$) $\perp l$ & 88\% & -- & -- & 8 \\
($1, i, l$) $\perp j$ & 79\% & z & $R_z$: 180 & 8 \\
($1, j, l$) $\perp i$ & 81\% & x & $R_x$: 180 & 8 \\
($i, j, l$) $\perp 1$ & 76\% & x & $R_x$: 180 & 8 \\
\bottomrule\end{tabular}}\par{\scriptsize Render-time measurements of version~1 ($2400$ points, threshold $0.97$): indicative only (Remark~\ref{atl:rem:mappings}).}\par\vspace{2mm}
The tesseract fold reflects each of the coordinates $1,i,j,k$ that exceeds the limit $0.6$ back inside the tesseract. The Okubonic body is replicated by reflection: in the quaternionic quarter the cube is almost entirely filled, with small holes on its faces; in the quarters containing $l$, where the fold acts only on two coordinates, the pictures show a regular grid of rounded, pillow-like cells, the mirrored copies of the body, bounded on the right by a smooth wavy surface.

The body shares are large ($76\%$ to $88\%$ at render time, indicative only), and the check finds only approximate mirrors ($J=0.98$) in the quarters containing $l$.
\clearpage

\subsubsection{M09. Rotation in the plane \texorpdfstring{$XW$}{XW}, octonions}\label{atl:pl:M09}\mbox{}\par\vspace{0.5mm}
\noindent\includegraphics[width=\textwidth]{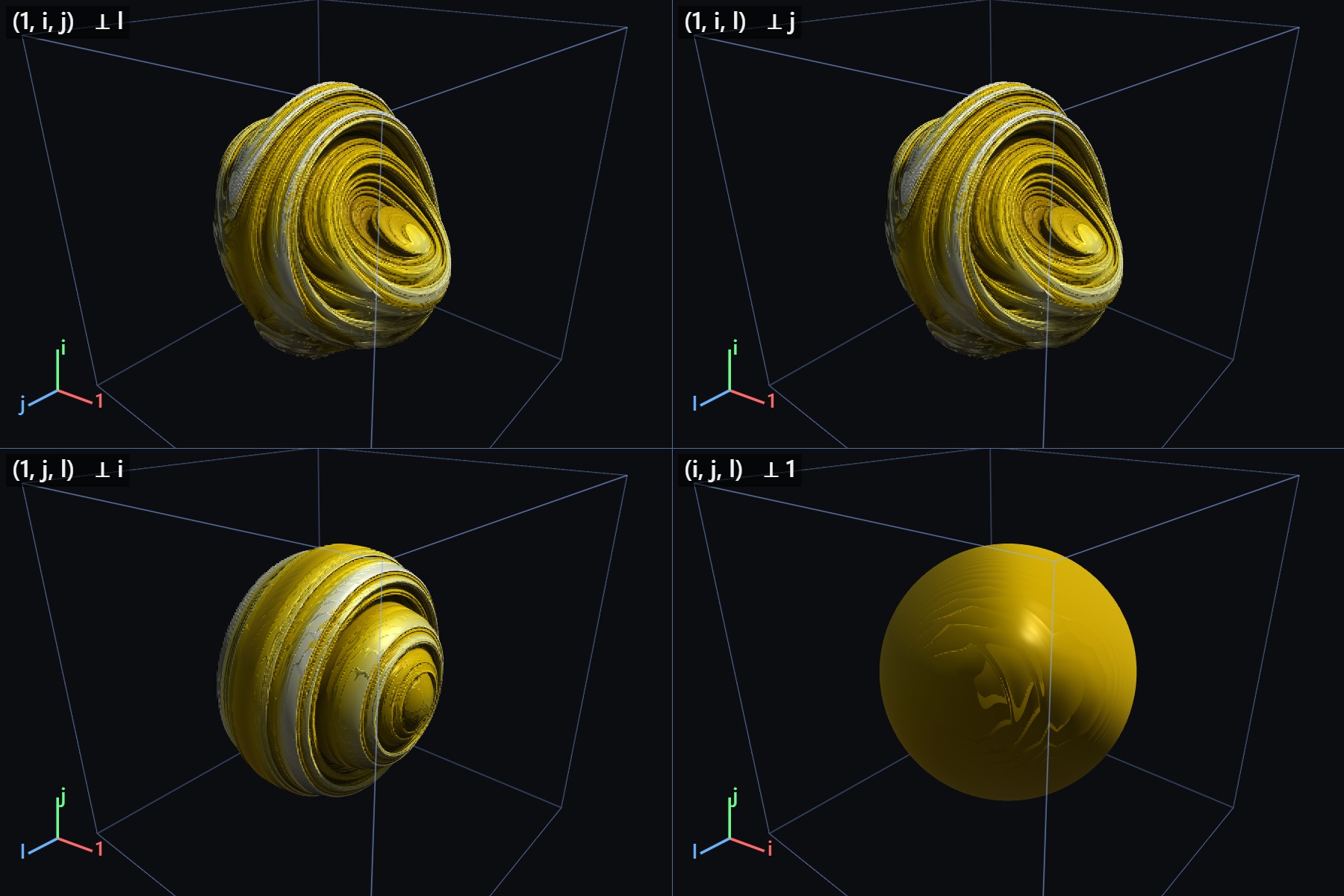}\par\vspace{1mm}
\noindent{\footnotesize \textbf{Map} $z\mapsto z\diamond z+c$ $\cdot$ \textbf{Algebra} Octonions (Cayley--Dickson, $\gamma=(-1,-1,-1)$) $\cdot$ $c=(-0.231,\,\allowbreak 0.62,\,\allowbreak 0,\,\allowbreak 0,\,\allowbreak 0,\,\allowbreak 0,\,\allowbreak 0,\,\allowbreak 0)$ ($=0.97\,s^*$, $s^*=0.6824$) $\cdot$ \textbf{Frame} $X,Y,Z,N=1,i,j,l$ $\cdot$ half-extent 1.4 $\cdot$ 20 iterations $\cdot$ \textbf{Coloring} 8D orbit trap (point), gradient Gold $\cdot$ \textbf{Mapping} \texttt{Rotation4D} (xw$=45$)}\par\vspace{1mm}
\noindent{\footnotesize\begin{tabular}{@{}lrlll@{}}\toprule Slice (X, Y, Z) & Body & Mirrors & Rotations ($J\geq0.97$) & Orbit span\\\midrule
($1, i, j$) $\perp l$ & 23\% & z & $R_z$: 180 & 3 \\
($1, i, l$) $\perp j$ & 23\% & z & $R_z$: 180 & 3 \\
($1, j, l$) $\perp i$ & 24\% & xyz & $R_x$: all; $R_y$: 180; $R_z$: 180 & 3 \\
($i, j, l$) $\perp 1$ & 35\% & xyz & $R_x$: all; $R_y$: all; $R_z$: all & 3 \\
\bottomrule\end{tabular}}\par{\scriptsize Render-time measurements of version~1 ($2400$ points, threshold $0.97$): indicative only (Remark~\ref{atl:rem:mappings}).}\par\vspace{2mm}
The rotation by $45^\circ$ in the plane $XW$, i.e., of the coordinates $1$ and $k$, is not a symmetry of the octonionic map, because it moves the unit. The shells of the two upper quarters are twisted into a spiral, but they remain identical, and the quarter $(1,j,l)$ keeps its ringed appearance; the imaginary quarter, where $x_0=k=0$, is not moved at all and is the ball of Plate~A01.

The explanation is Proposition~\ref{atl:prop:equivariance}: the rotation mixes $1$ and $k$, but the map remains equivariant under the orthogonal transformations of the span of $j,l,il,jl,kl$, which fix $1$, $i$ and $k$ and commute with the mapping. The exchange of $j$ and $l$ is one of them, and it maps the rotated section $(1,i,j)$ onto the rotated section $(1,i,l)$; the rotations of the plane $\langle j,l\rangle$ are others, and the quarter $(1,j,l)$ remains a body of revolution about $1$. The render-time check (indicative only; Part~M is not measured) confirms both, and finds in the upper quarters the mirror in the third axis.
\clearpage

\subsubsection{M10. Rotation in the plane \texorpdfstring{$XW$}{XW}, Okubo}\label{atl:pl:M10}\mbox{}\par\vspace{0.5mm}
\noindent\includegraphics[width=\textwidth]{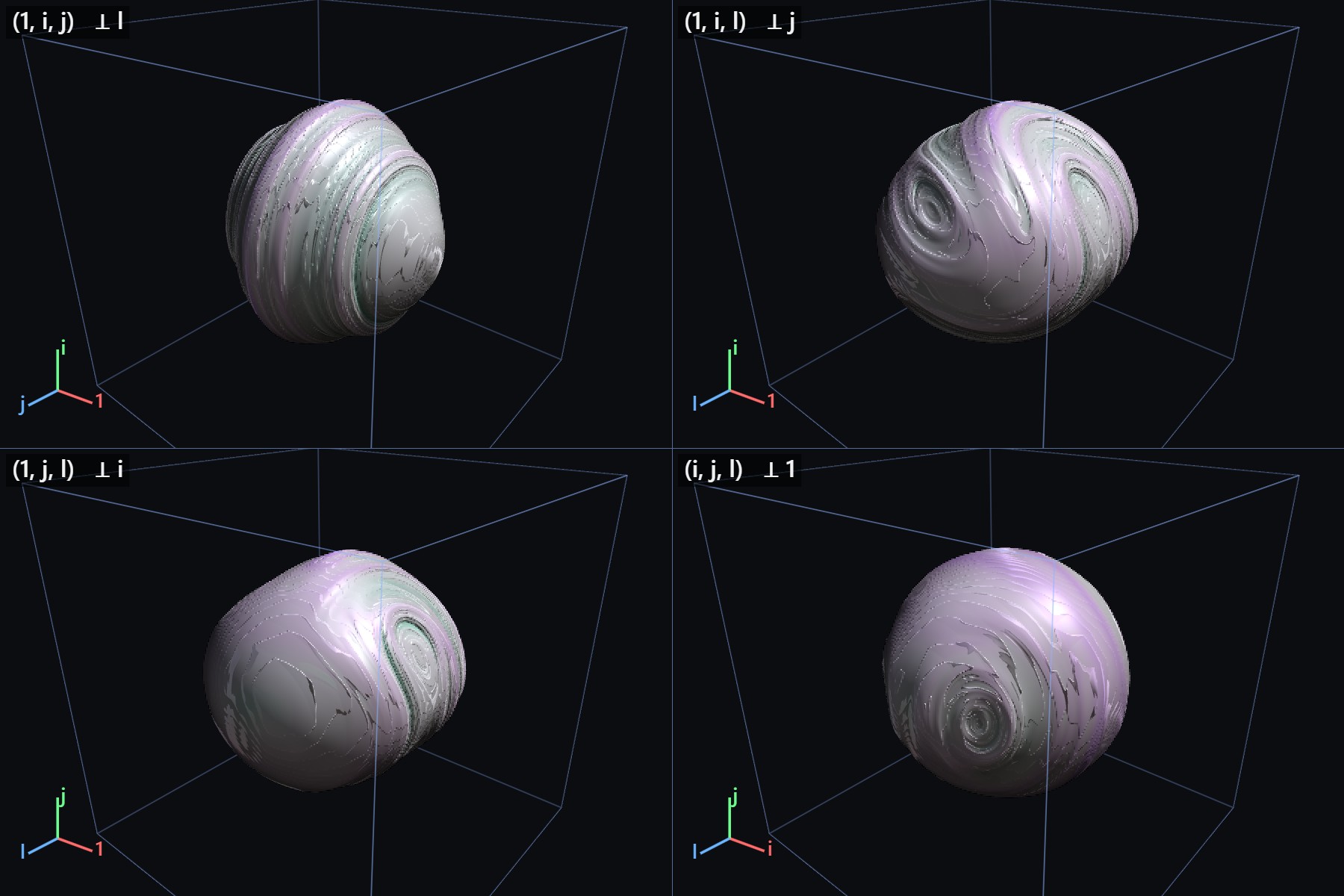}\par\vspace{1mm}
\noindent{\footnotesize \textbf{Map} $z\mapsto z\diamond z+c$ $\cdot$ \textbf{Algebra} Okubo algebra (Petersson, $\gamma=(-1,-1,-1)$, $\theta=120^\circ$) $\cdot$ $c=(-0.397,\,\allowbreak -0.021,\,\allowbreak 0.049,\,\allowbreak 0.016,\,\allowbreak 0.004,\,\allowbreak 0.16,\,\allowbreak -0.113,\,\allowbreak 0.133)$ ($=0.95\,s^*$, $s^*=0.4900$) $\cdot$ \textbf{Frame} $X,Y,Z,N=1,i,j,l$ $\cdot$ half-extent 1.4 $\cdot$ 20 iterations $\cdot$ \textbf{Coloring} 8D orbit trap (point), gradient Aurora $\cdot$ \textbf{Mapping} \texttt{Rotation4D} (xw$=45$)}\par\vspace{1mm}
\noindent{\footnotesize\begin{tabular}{@{}lrlll@{}}\toprule Slice (X, Y, Z) & Body & Mirrors & Rotations ($J\geq0.97$) & Orbit span\\\midrule
($1, i, j$) $\perp l$ & 25\% & -- & -- & 8 \\
($1, i, l$) $\perp j$ & 32\% & -- & -- & 8 \\
($1, j, l$) $\perp i$ & 33\% & -- & -- & 8 \\
($i, j, l$) $\perp 1$ & 34\% & x & $R_x$: 180 & 8 \\
\bottomrule\end{tabular}}\par{\scriptsize Render-time measurements of version~1 ($2400$ points, threshold $0.97$): indicative only (Remark~\ref{atl:rem:mappings}).}\par\vspace{2mm}
The same rotation on the Okubonic body of Plate~M06. The quaternionic quarter is a ribbed shell, the quarters $(1,i,l)$ and $(1,j,l)$ are rounded bodies with spiral vortices, and the imaginary quarter, where $x_0=k=0$, is not moved: it is the same body as in Plate~M06. The render-time check (indicative only) finds no symmetry except the approximate mirror of the imaginary quarter already present in Plate~M06.

The comparison with Plate~M09 must be read with care. There the constant is complex, and the large group of Proposition~\ref{atl:prop:equivariance} contains maps that commute with the rotation; here the constant is generic, its stabilizer is a torus whose root planes are in general position (Plate~R09), and there is no coordinate symmetry for the rotation to preserve. The Okubonic body is simply seen along a different section.
\clearpage

\subsubsection{M11. Planar kaleidoscope, Okubo}\label{atl:pl:M11}\mbox{}\par\vspace{0.5mm}
\noindent\includegraphics[width=\textwidth]{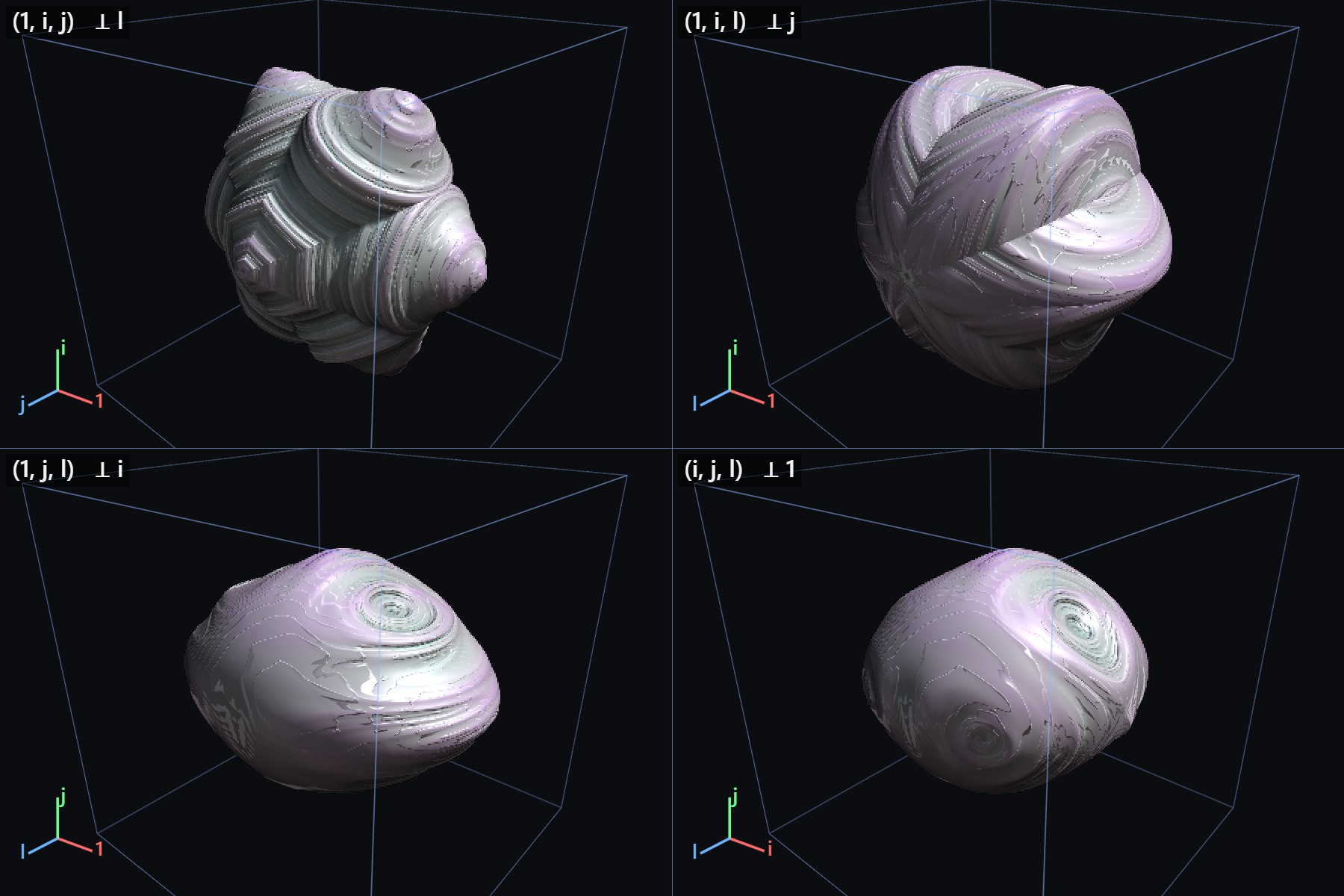}\par\vspace{1mm}
\noindent{\footnotesize \textbf{Map} $z\mapsto z\diamond z+c$ $\cdot$ \textbf{Algebra} Okubo algebra (Petersson, $\gamma=(-1,-1,-1)$, $\theta=120^\circ$) $\cdot$ $c=(-0.397,\,\allowbreak -0.021,\,\allowbreak 0.049,\,\allowbreak 0.016,\,\allowbreak 0.004,\,\allowbreak 0.16,\,\allowbreak -0.113,\,\allowbreak 0.133)$ ($=0.95\,s^*$, $s^*=0.4900$) $\cdot$ \textbf{Frame} $X,Y,Z,N=1,i,j,l$ $\cdot$ half-extent 1.4 $\cdot$ 20 iterations $\cdot$ \textbf{Coloring} 8D orbit trap (point), gradient Aurora $\cdot$ \textbf{Mapping} \texttt{Kaleidoscope}}\par\vspace{1mm}
\noindent{\footnotesize\begin{tabular}{@{}lrlll@{}}\toprule Slice (X, Y, Z) & Body & Mirrors & Rotations ($J\geq0.97$) & Orbit span\\\midrule
($1, i, j$) $\perp l$ & 36\% & xy & $R_z$: 180,120,60 & 8 \\
($1, i, l$) $\perp j$ & 55\% & xy & $R_z$: 180,120,60 & 8 \\
($1, j, l$) $\perp i$ & 40\% & x & -- & 8 \\
($i, j, l$) $\perp 1$ & 35\% & x & -- & 8 \\
\bottomrule\end{tabular}}\par{\scriptsize Render-time measurements of version~1 ($2400$ points, threshold $0.97$): indicative only (Remark~\ref{atl:rem:mappings}).}\par\vspace{2mm}
The planar kaleidoscope of Ultra Fractal acts on the plane of the first two coordinates, $1$ and $i$, and folds it into six mirrored wedges. The quarters $(1,i,j)$ and $(1,i,l)$, which contain this plane, acquire by construction a six-fold dihedral symmetry about the axis orthogonal to it: the render-time check (indicative only; Part~M is not measured) reports the mirrors in $1$ and $i$ and the rotations by $60^\circ$, $120^\circ$ and $180^\circ$ about the third axis. The quaternionic quarter shows a lobed body with conical bulbs arranged around the $j$-axis, the quarter $(1,i,l)$ a body with creased faces meeting at star-shaped points, a crystalline structure foreign to the Okubo algebra.

The two lower quarters contain only one of the two coordinates. In $(1,j,l)$ the fold reduces to $x_0\mapsto|x_0|$, and the quarter is the base body of Plate~M06 mirrored in $1$; in $(i,j,l)$ the points of the $i$-axis are sent to the line at $30^\circ$ in the plane of $1$ and $i$, so the quarter shows a tilted section of the base body, mirrored in $i$. The check finds only these mirrors.
\clearpage

\subsubsection{M12. Planar twirl, octonions}\label{atl:pl:M12}\mbox{}\par\vspace{0.5mm}
\noindent\includegraphics[width=\textwidth]{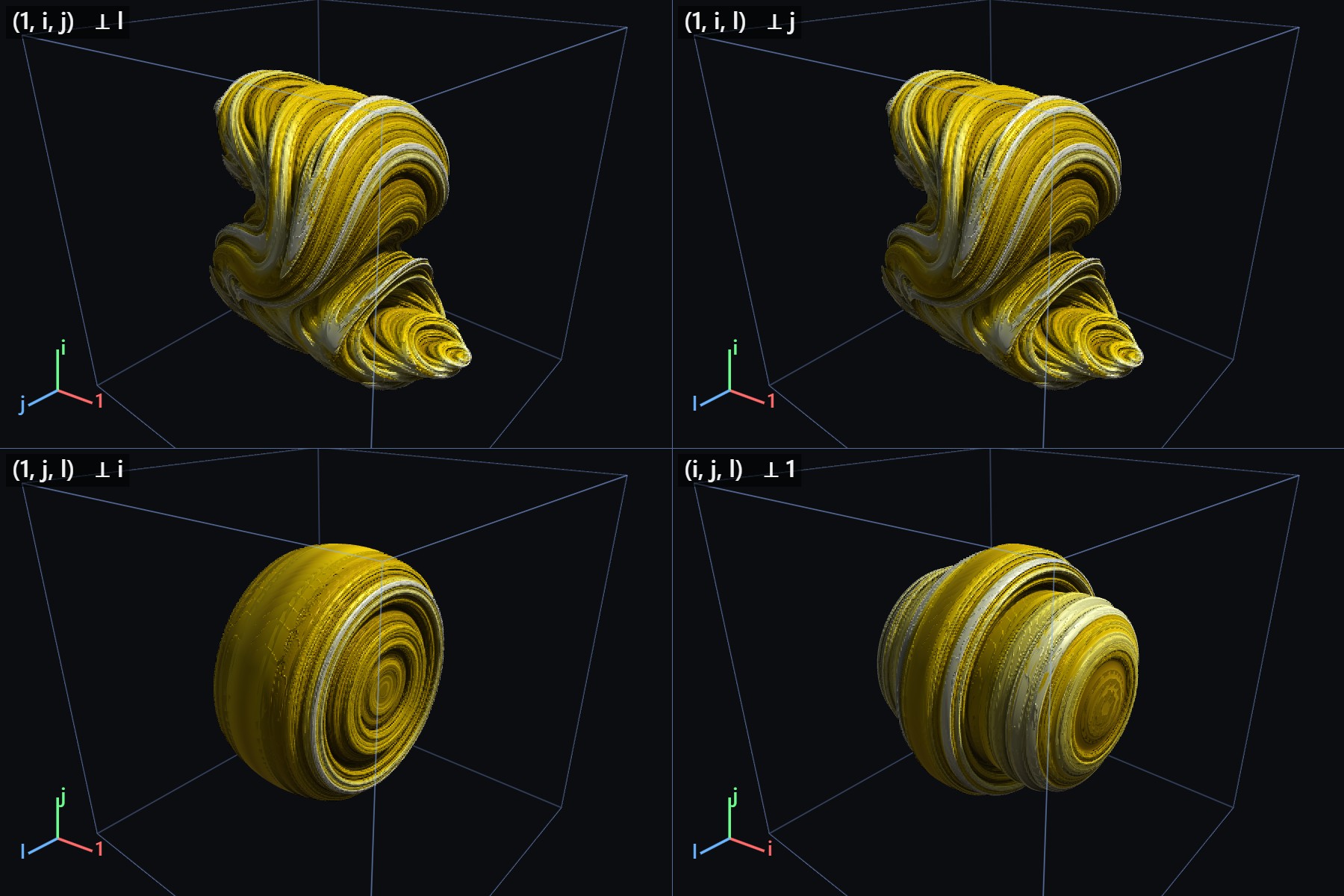}\par\vspace{1mm}
\noindent{\footnotesize \textbf{Map} $z\mapsto z\diamond z+c$ $\cdot$ \textbf{Algebra} Octonions (Cayley--Dickson, $\gamma=(-1,-1,-1)$) $\cdot$ $c=(-0.231,\,\allowbreak 0.62,\,\allowbreak 0,\,\allowbreak 0,\,\allowbreak 0,\,\allowbreak 0,\,\allowbreak 0,\,\allowbreak 0)$ ($=0.97\,s^*$, $s^*=0.6824$) $\cdot$ \textbf{Frame} $X,Y,Z,N=1,i,j,l$ $\cdot$ half-extent 1.4 $\cdot$ 20 iterations $\cdot$ \textbf{Coloring} 8D orbit trap (point), gradient Gold $\cdot$ \textbf{Mapping} \texttt{Twirl}}\par\vspace{1mm}
\noindent{\footnotesize\begin{tabular}{@{}lrlll@{}}\toprule Slice (X, Y, Z) & Body & Mirrors & Rotations ($J\geq0.97$) & Orbit span\\\midrule
($1, i, j$) $\perp l$ & 24\% & z & $R_z$: 180 & 3 \\
($1, i, l$) $\perp j$ & 24\% & z & $R_z$: 180 & 3 \\
($1, j, l$) $\perp i$ & 22\% & xyz & $R_x$: all; $R_y$: 180; $R_z$: 180 & 3 \\
($i, j, l$) $\perp 1$ & 28\% & xyz & $R_x$: all; $R_y$: 180; $R_z$: 180 & 3 \\
\bottomrule\end{tabular}}\par{\scriptsize Render-time measurements of version~1 ($2400$ points, threshold $0.97$): indicative only (Remark~\ref{atl:rem:mappings}).}\par\vspace{2mm}
The planar twirl rotates the plane of $1$ and $i$ by an angle that decreases with the distance from the centre. On the octonionic body ($c\in\CC$) it twists the shells of the two upper quarters into swirling shapes, identical in both quarters, since the twirl commutes with the exchange of $j$ and $l$ (Proposition~\ref{atl:prop:equivariance}).

The two lower quarters are ringed bodies. The twirl commutes with the rotations of the plane $\langle j,l\rangle$, which are symmetries of the octonionic map, so the quarter $(1,j,l)$ remains a body of revolution about $1$ and the imaginary quarter becomes a body of revolution about $i$; it is no longer a ball, because the twirl moves its samples out of the imaginary hyperplane. The render-time check (indicative only; Part~M is not measured) finds both revolutions.
\clearpage

\subsubsection{M13. Mirror in a hyperplane, Okubo}\label{atl:pl:M13}\mbox{}\par\vspace{0.5mm}
\noindent\includegraphics[width=\textwidth]{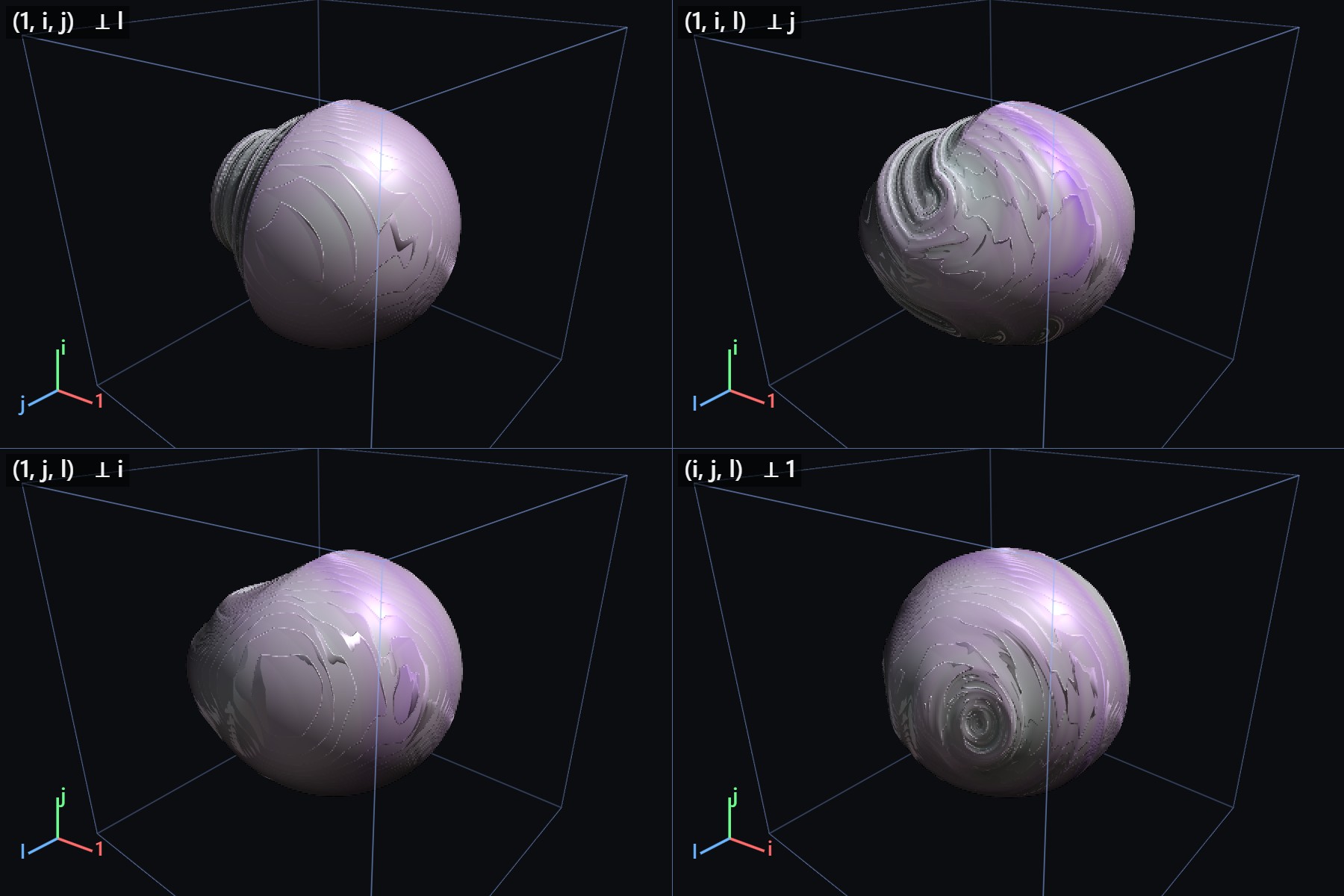}\par\vspace{1mm}
\noindent{\footnotesize \textbf{Map} $z\mapsto z\diamond z+c$ $\cdot$ \textbf{Algebra} Okubo algebra (Petersson, $\gamma=(-1,-1,-1)$, $\theta=120^\circ$) $\cdot$ $c=(-0.397,\,\allowbreak -0.021,\,\allowbreak 0.049,\,\allowbreak 0.016,\,\allowbreak 0.004,\,\allowbreak 0.16,\,\allowbreak -0.113,\,\allowbreak 0.133)$ ($=0.95\,s^*$, $s^*=0.4900$) $\cdot$ \textbf{Frame} $X,Y,Z,N=1,i,j,l$ $\cdot$ half-extent 1.4 $\cdot$ 20 iterations $\cdot$ \textbf{Coloring} 8D orbit trap (point), gradient Aurora $\cdot$ \textbf{Mapping} \texttt{Mirror4D}}\par\vspace{1mm}
\noindent{\footnotesize\begin{tabular}{@{}lrlll@{}}\toprule Slice (X, Y, Z) & Body & Mirrors & Rotations ($J\geq0.97$) & Orbit span\\\midrule
($1, i, j$) $\perp l$ & 30\% & y & $R_x$: 60,37 & 8 \\
($1, i, l$) $\perp j$ & 36\% & -- & -- & 8 \\
($1, j, l$) $\perp i$ & 38\% & -- & -- & 8 \\
($i, j, l$) $\perp 1$ & 34\% & x & $R_x$: 180 & 8 \\
\bottomrule\end{tabular}}\par{\scriptsize Render-time measurements of version~1 ($2400$ points, threshold $0.97$): indicative only (Remark~\ref{atl:rem:mappings}).}\par\vspace{2mm}
The mirror in the hyperplane orthogonal to $(1,0,0,1)$ replaces the half-space in front of it by the mirror image of the other half. For the samples of the frame $F_A$, whose coordinate $k$ vanishes, a point with $x_0>0$ is sent to the point with $x_0=0$ and $k=-x_0$; the points with $x_0\le0$ are not moved. In the quaternionic quarter one half is the ribbed shell of Plate~M06 and the other a smooth dome, the body at the points with $k=-x_0$; the quarters $(1,i,l)$ and $(1,j,l)$ show similar joins of a smooth half and a swirled half. The imaginary quarter, where $x_0=0$, is not moved and coincides with that of Plate~M06.

The render-time check (indicative only; Part~M is not measured) finds no clear symmetry; in the quaternionic quarter it reports two rotations about the $1$-axis at $J=0.98$, but not the others, which on $2400$ points is not evidence of a symmetry. The plate illustrates the difficulty of reading a four-dimensional mapping in an eight-dimensional section.
\clearpage

\subsection{Part N. Coloring algorithms on one body}\label{atl:part:N}
The same Okubonic body, the Julia set of Plate~B07, is drawn with nine coloring algorithms. The geometry and the measurements are identical in all plates of this part; only the coloring changes, and with it the information that reaches the eye.

\subsubsection{N01. Smooth iteration count (8D norm)}\label{atl:pl:N01}\mbox{}\par\vspace{0.5mm}
\noindent\includegraphics[width=\textwidth]{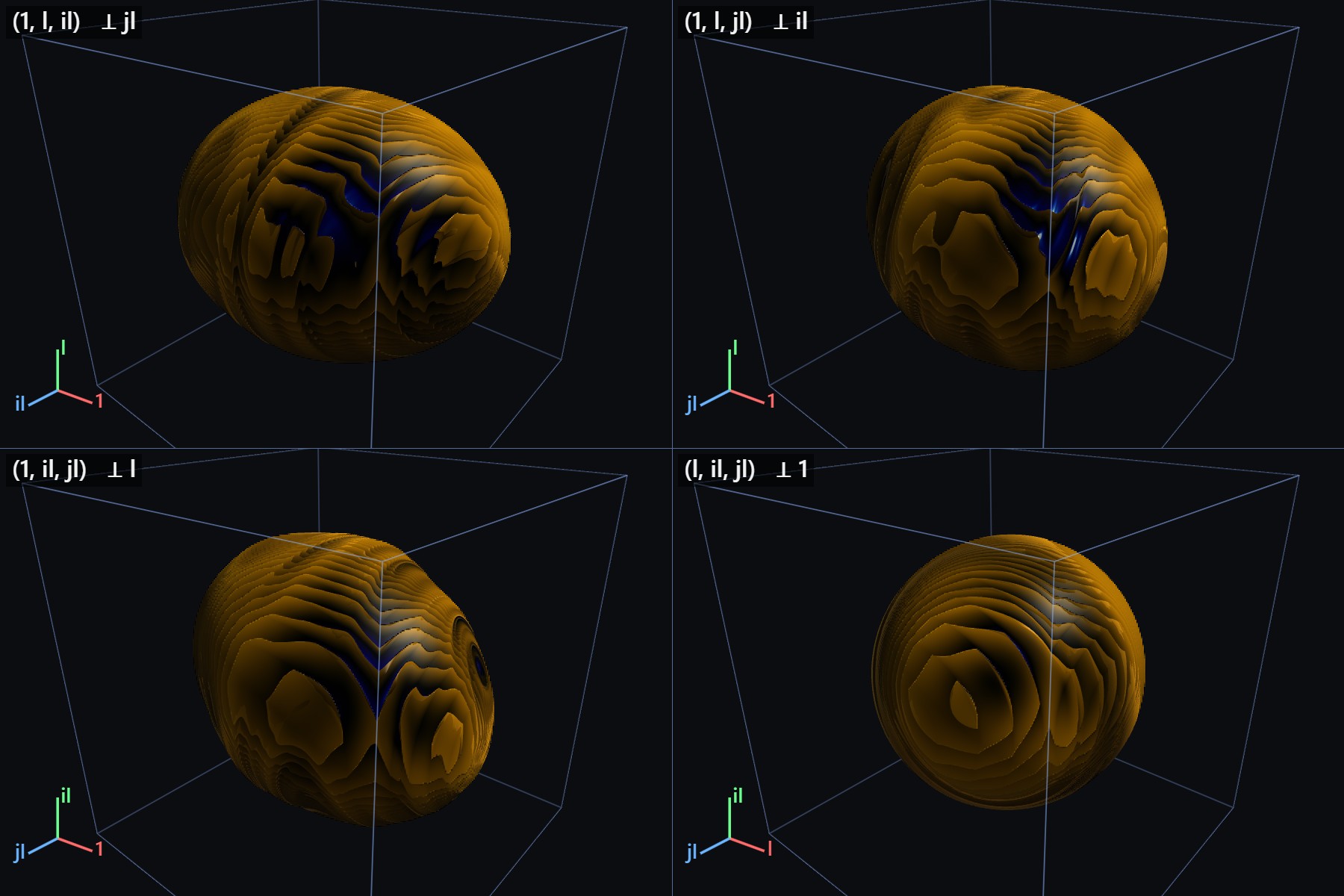}\par\vspace{1mm}
\noindent{\footnotesize \textbf{Map} $z\mapsto z\diamond z+c$ $\cdot$ \textbf{Algebra} Okubo algebra (Petersson, $\gamma=(-1,-1,-1)$, $\theta=120^\circ$) $\cdot$ $c=(-0.397,\,\allowbreak -0.021,\,\allowbreak 0.049,\,\allowbreak 0.016,\,\allowbreak 0.004,\,\allowbreak 0.16,\,\allowbreak -0.113,\,\allowbreak 0.133)$ ($=0.95\,s^*$, $s^*=0.4900$) $\cdot$ \textbf{Frame} $X,Y,Z,N=1,l,il,jl$ $\cdot$ half-extent 1.4 $\cdot$ 20 iterations $\cdot$ \textbf{Coloring} smooth iteration count (8D norm), gradient Classic}\par\vspace{1mm}
\noindent{\footnotesize\begin{tabular}{@{}lllll@{}}\toprule Slice (X, Y, Z) & Body [95\%] & Exact ($J\geq0.9995$) & Approximate ($0.97\leq J<0.9995$) & Dim\\\midrule
($1, l, il$) $\perp jl$ & 52.6\% {\scriptsize[52.3, 52.9]} & $P$ & $M$: y; $R_y$: 180 (min.\ 0.987) & 8 \\
($1, l, jl$) $\perp il$ & 51.9\% {\scriptsize[51.6, 52.2]} & $P$ & -- & 8 \\
($1, il, jl$) $\perp l$ & 51.5\% {\scriptsize[51.2, 51.8]} & $P$ & -- & 8 \\
($l, il, jl$) $\perp 1$ & 42.3\% {\scriptsize[42.0, 42.6]} & $P$ & $M$: x; $R_x$: 180,37; $R_y$: 90,72,60,37; $R_z$: 90,72,60,37 (min.\ 0.972) & 8 \\
\bottomrule\end{tabular}}\par\vspace{2mm}
Part~N keeps one body fixed, the Okubonic Julia set of Plate~B07 with the generic constant in the frame $F_B$, and changes only the coloring. The measurements are those of Plate~B07 up to the Monte Carlo error, as they must be: bodies of $53\%$, $52\%$, $52\%$ and $42\%$ of the ball, the point reflection $P$ as the only exact symmetry of every quarter (Remark~\ref{atl:rem:minus}), and orbits in a subalgebra of dimension $8$. The shapes are the same in Plates~N01--N09: three fat eggs in the quarters through $1$ and a nearly round body in the imaginary quarter $(l,il,jl)$, whose approximate rotations ($J$ between $0.97$ and $0.99$) do not reach exactness. It is worth noticing that the thin, irregular terrace lines visible on the surfaces recur identically in every plate of the part, whatever the coloring: they belong to the relief of the rendered surface, not to the coloring.

Here the surface is colored by the smooth iteration count computed from the eight-dimensional Euclidean norm of the orbit value, the continuous escape time familiar from the complex Mandelbrot set \cite{Milnor06}. With the Classic gradient the bodies are amber, crossed by dark bands that follow the terraces of the surface; in the three quarters through $1$ a deep blue appears in the V-shaped crease at the centre of the visible face. The imaginary quarter shows amber and dark concentric bands without blue.
\clearpage

\subsubsection{N02. 8D orbit trap: a point}\label{atl:pl:N02}\mbox{}\par\vspace{0.5mm}
\noindent\includegraphics[width=\textwidth]{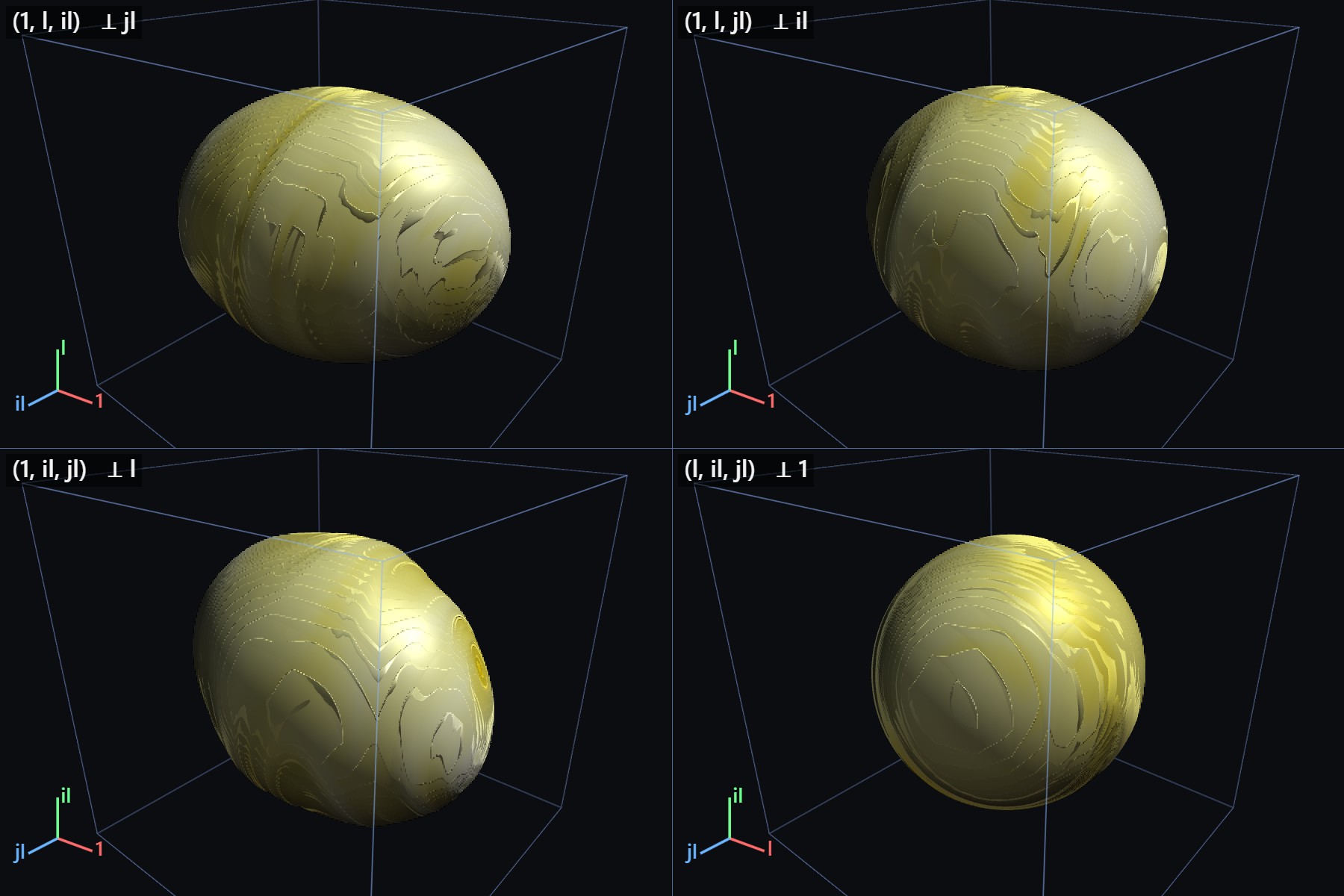}\par\vspace{1mm}
\noindent{\footnotesize \textbf{Map} $z\mapsto z\diamond z+c$ $\cdot$ \textbf{Algebra} Okubo algebra (Petersson, $\gamma=(-1,-1,-1)$, $\theta=120^\circ$) $\cdot$ $c=(-0.397,\,\allowbreak -0.021,\,\allowbreak 0.049,\,\allowbreak 0.016,\,\allowbreak 0.004,\,\allowbreak 0.16,\,\allowbreak -0.113,\,\allowbreak 0.133)$ ($=0.95\,s^*$, $s^*=0.4900$) $\cdot$ \textbf{Frame} $X,Y,Z,N=1,l,il,jl$ $\cdot$ half-extent 1.4 $\cdot$ 20 iterations $\cdot$ \textbf{Coloring} 8D orbit trap (point), gradient Gold}\par\vspace{1mm}
\noindent{\footnotesize\begin{tabular}{@{}lllll@{}}\toprule Slice (X, Y, Z) & Body [95\%] & Exact ($J\geq0.9995$) & Approximate ($0.97\leq J<0.9995$) & Dim\\\midrule
($1, l, il$) $\perp jl$ & 52.7\% {\scriptsize[52.4, 53.0]} & $P$ & $M$: y; $R_y$: 180 (min.\ 0.986) & 8 \\
($1, l, jl$) $\perp il$ & 51.8\% {\scriptsize[51.5, 52.1]} & $P$ & -- & 8 \\
($1, il, jl$) $\perp l$ & 51.3\% {\scriptsize[51.0, 51.6]} & $P$ & -- & 8 \\
($l, il, jl$) $\perp 1$ & 42.2\% {\scriptsize[41.9, 42.5]} & $P$ & $M$: x; $R_x$: 180,37; $R_y$: 90,72,60,37; $R_z$: 90,72,60,37 (min.\ 0.972) & 8 \\
\bottomrule\end{tabular}}\par\vspace{2mm}
The eight-dimensional orbit trap with a point at the origin colors each point by the minimum distance of its orbit from the origin. This is the coloring of Plate~B07, here with the Gold gradient instead of Aurora, and the plate shows the same body (bodies of $53\%$, $52\%$, $51\%$ and $42\%$, only $P$ exact, dimension $8$).

With these settings the trap varies little over the surface: the bodies are an almost uniform pale gold, modulated mainly by the lighting, with a small brighter whorl at the right end of the quarter $(1,il,jl)$. The thin light lines that cover the surfaces are the terrace lines of the relief, common to every plate of Part~N (Plate~N01), and not contours of the trap. The comparison with Plates~N05, N06 and N08 shows that the colorings based on the direction of the orbit vary much more over the surface than the distance to a point.
\clearpage

\subsubsection{N03. 8D orbit trap: a hyperplane}\label{atl:pl:N03}\mbox{}\par\vspace{0.5mm}
\noindent\includegraphics[width=\textwidth]{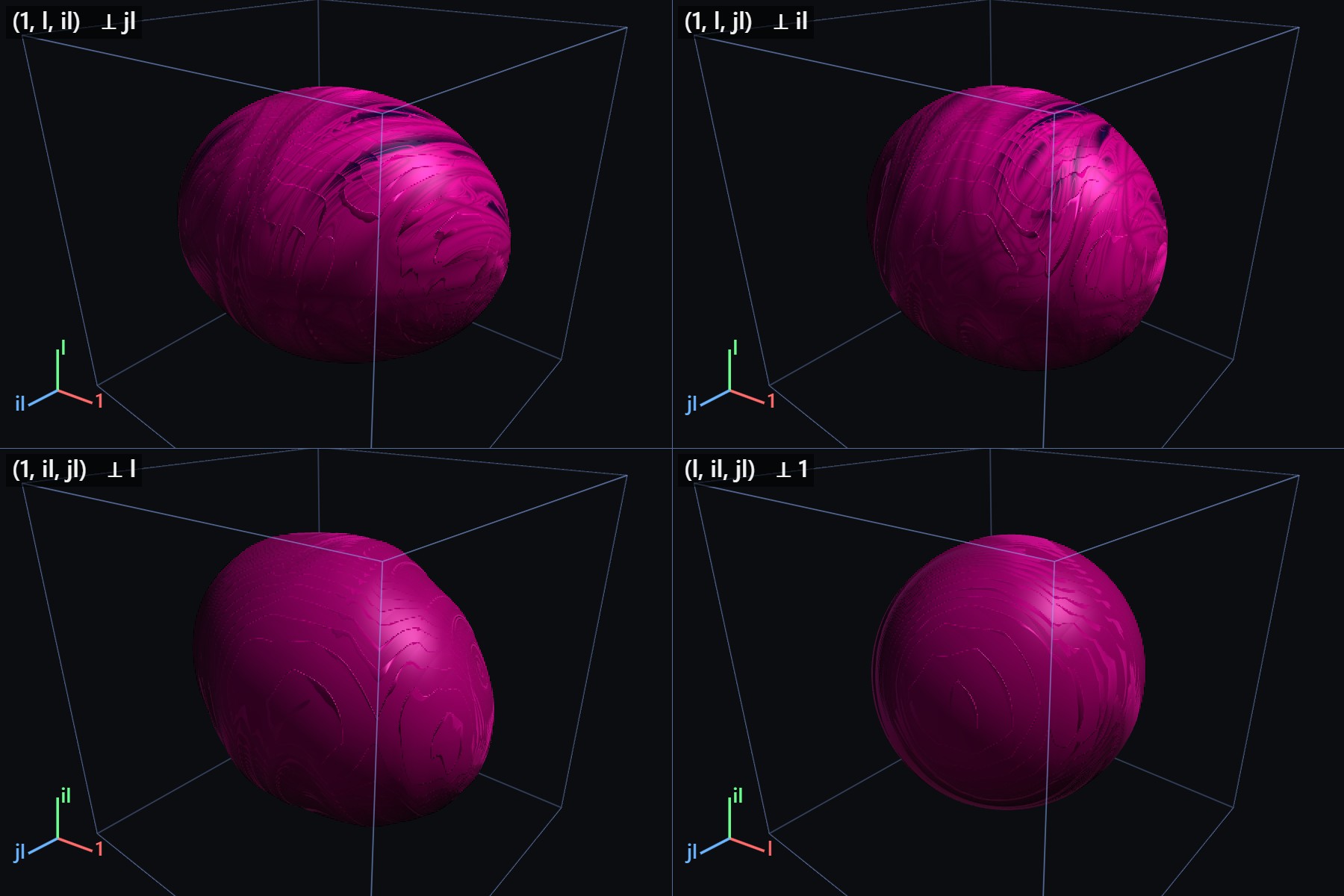}\par\vspace{1mm}
\noindent{\footnotesize \textbf{Map} $z\mapsto z\diamond z+c$ $\cdot$ \textbf{Algebra} Okubo algebra (Petersson, $\gamma=(-1,-1,-1)$, $\theta=120^\circ$) $\cdot$ $c=(-0.397,\,\allowbreak -0.021,\,\allowbreak 0.049,\,\allowbreak 0.016,\,\allowbreak 0.004,\,\allowbreak 0.16,\,\allowbreak -0.113,\,\allowbreak 0.133)$ ($=0.95\,s^*$, $s^*=0.4900$) $\cdot$ \textbf{Frame} $X,Y,Z,N=1,l,il,jl$ $\cdot$ half-extent 1.4 $\cdot$ 20 iterations $\cdot$ \textbf{Coloring} 8D orbit trap (hyperplane $\perp l$), gradient Neon}\par\vspace{1mm}
\noindent{\footnotesize\begin{tabular}{@{}lllll@{}}\toprule Slice (X, Y, Z) & Body [95\%] & Exact ($J\geq0.9995$) & Approximate ($0.97\leq J<0.9995$) & Dim\\\midrule
($1, l, il$) $\perp jl$ & 52.5\% {\scriptsize[52.2, 52.9]} & $P$ & $M$: y; $R_y$: 180 (min.\ 0.986) & 8 \\
($1, l, jl$) $\perp il$ & 51.7\% {\scriptsize[51.4, 52.0]} & $P$ & -- & 8 \\
($1, il, jl$) $\perp l$ & 51.2\% {\scriptsize[50.9, 51.5]} & $P$ & -- & 8 \\
($l, il, jl$) $\perp 1$ & 42.0\% {\scriptsize[41.7, 42.3]} & $P$ & $M$: x; $R_x$: 180,37; $R_y$: 90,72,60,37; $R_z$: 90,72,60,37 (min.\ 0.973) & 8 \\
\bottomrule\end{tabular}}\par\vspace{2mm}
The trap on the hyperplane $x_4=0$, orthogonal to $l$, colors each point by the minimum distance of its orbit from that hyperplane. The body is again that of Plate~B07 (bodies of $53\%$, $52\%$, $51\%$ and $42\%$, only $P$ exact, dimension $8$).

With the Neon gradient the surfaces are magenta. In the two upper quarters, $(1,l,il)$ and $(1,l,jl)$, whose vertical axis is $l$, the upper part of the body is marked by darker swirling streaks and brighter pink patches; the two lower quarters are an almost uniform magenta, with only the terrace lines of the relief (Plate~N01). It is natural that the trap reacts most where the initial point itself has a large component along $l$, but the plate does not allow a finer reading: this coloring singles out one hidden direction and is otherwise rather uniform.
\clearpage

\subsubsection{N04. 8D orbit trap: a 7-sphere}\label{atl:pl:N04}\mbox{}\par\vspace{0.5mm}
\noindent\includegraphics[width=\textwidth]{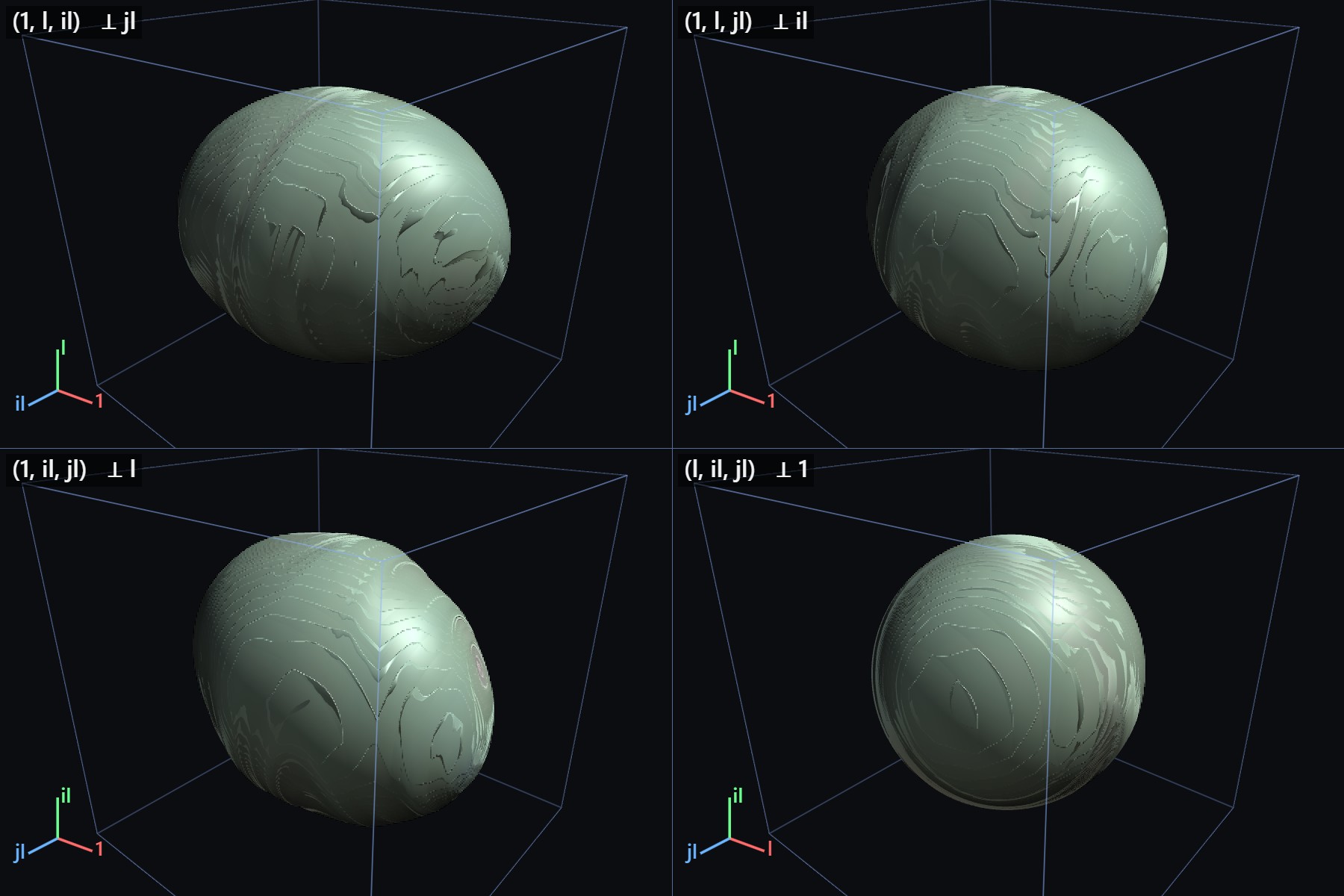}\par\vspace{1mm}
\noindent{\footnotesize \textbf{Map} $z\mapsto z\diamond z+c$ $\cdot$ \textbf{Algebra} Okubo algebra (Petersson, $\gamma=(-1,-1,-1)$, $\theta=120^\circ$) $\cdot$ $c=(-0.397,\,\allowbreak -0.021,\,\allowbreak 0.049,\,\allowbreak 0.016,\,\allowbreak 0.004,\,\allowbreak 0.16,\,\allowbreak -0.113,\,\allowbreak 0.133)$ ($=0.95\,s^*$, $s^*=0.4900$) $\cdot$ \textbf{Frame} $X,Y,Z,N=1,l,il,jl$ $\cdot$ half-extent 1.4 $\cdot$ 20 iterations $\cdot$ \textbf{Coloring} 8D orbit trap (7-sphere $r=0.6$), gradient Pastel}\par\vspace{1mm}
\noindent{\footnotesize\begin{tabular}{@{}lllll@{}}\toprule Slice (X, Y, Z) & Body [95\%] & Exact ($J\geq0.9995$) & Approximate ($0.97\leq J<0.9995$) & Dim\\\midrule
($1, l, il$) $\perp jl$ & 52.7\% {\scriptsize[52.4, 53.1]} & $P$ & $M$: y; $R_y$: 180 (min.\ 0.986) & 8 \\
($1, l, jl$) $\perp il$ & 52.0\% {\scriptsize[51.7, 52.3]} & $P$ & -- & 8 \\
($1, il, jl$) $\perp l$ & 51.6\% {\scriptsize[51.3, 51.9]} & $P$ & -- & 8 \\
($l, il, jl$) $\perp 1$ & 42.2\% {\scriptsize[41.9, 42.6]} & $P$ & $M$: x; $R_x$: 180,37; $R_y$: 90,72,60,37; $R_z$: 90,72,60,37 (min.\ 0.972) & 8 \\
\bottomrule\end{tabular}}\par\vspace{2mm}
The trap on the 7-sphere of radius $0.6$ colors each point by the minimum distance of its orbit from that sphere. The body is that of Plate~B07 (bodies of $53\%$, $52\%$, $52\%$ and $42\%$, only $P$ exact, dimension $8$).

With the Pastel gradient the bodies are an almost uniform grey-green; the visible lines are those of the relief (Plate~N01), and the only trace of the coloring is a faint pink tint in the small whorl at the right end of the quarter $(1,il,jl)$. With these settings the trap is nearly constant on the surface, and for this body it is the least informative coloring of the part.
\clearpage

\subsubsection{N05. Hidden part: orbit average}\label{atl:pl:N05}\mbox{}\par\vspace{0.5mm}
\noindent\includegraphics[width=\textwidth]{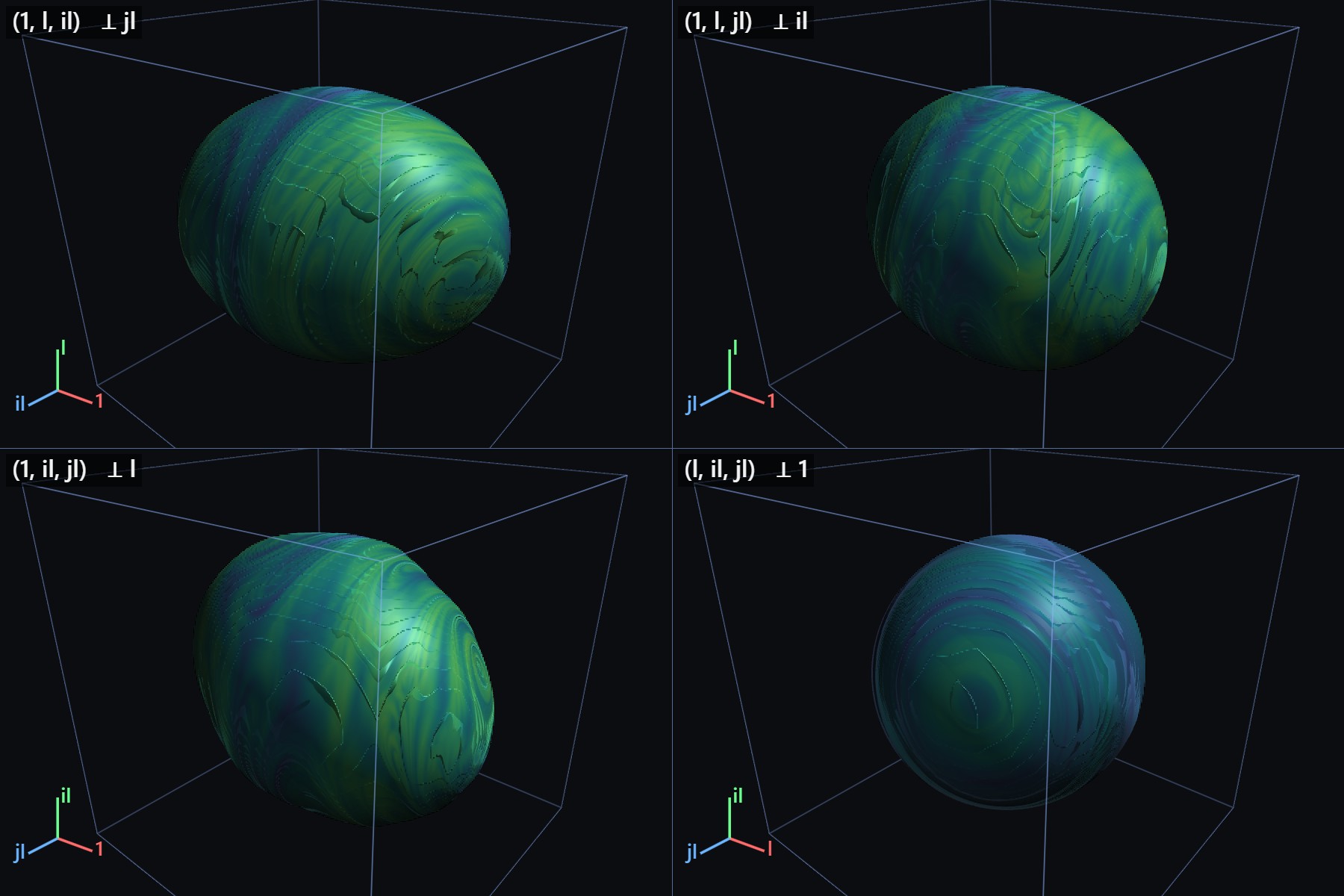}\par\vspace{1mm}
\noindent{\footnotesize \textbf{Map} $z\mapsto z\diamond z+c$ $\cdot$ \textbf{Algebra} Okubo algebra (Petersson, $\gamma=(-1,-1,-1)$, $\theta=120^\circ$) $\cdot$ $c=(-0.397,\,\allowbreak -0.021,\,\allowbreak 0.049,\,\allowbreak 0.016,\,\allowbreak 0.004,\,\allowbreak 0.16,\,\allowbreak -0.113,\,\allowbreak 0.133)$ ($=0.95\,s^*$, $s^*=0.4900$) $\cdot$ \textbf{Frame} $X,Y,Z,N=1,l,il,jl$ $\cdot$ half-extent 1.4 $\cdot$ 20 iterations $\cdot$ \textbf{Coloring} hidden part (orbit average), gradient Viridis-like}\par\vspace{1mm}
\noindent{\footnotesize\begin{tabular}{@{}lllll@{}}\toprule Slice (X, Y, Z) & Body [95\%] & Exact ($J\geq0.9995$) & Approximate ($0.97\leq J<0.9995$) & Dim\\\midrule
($1, l, il$) $\perp jl$ & 52.5\% {\scriptsize[52.2, 52.9]} & $P$ & $M$: y; $R_y$: 180 (min.\ 0.986) & 8 \\
($1, l, jl$) $\perp il$ & 51.8\% {\scriptsize[51.4, 52.1]} & $P$ & -- & 8 \\
($1, il, jl$) $\perp l$ & 51.3\% {\scriptsize[51.0, 51.6]} & $P$ & -- & 8 \\
($l, il, jl$) $\perp 1$ & 42.0\% {\scriptsize[41.7, 42.4]} & $P$ & $M$: x; $R_x$: 180,37; $R_y$: 90,72,60,37; $R_z$: 90,72,60,37 (min.\ 0.973) & 8 \\
\bottomrule\end{tabular}}\par\vspace{2mm}
The ``hidden part'' coloring measures, at each iterate, the angle between the orbit value and the quaternionic subspace $\langle1,i,j,k\rangle$, i.e.\ how much of the orbit lies in the hidden half $\langle l,il,jl,kl\rangle$, and averages it over the orbit. It was designed for this atlas. The body is that of Plate~B07 (bodies of $53\%$, $52\%$, $51\%$ and $42\%$, only $P$ exact, orbits in a subalgebra of dimension $8$).

With the Viridis-like gradient the three quarters through $1$ show a marbled texture of greens, with blue-violet streaks and yellow-green highlights, organized in bands that wrap around the bodies; the imaginary quarter is darker, in blues and teals, with concentric bands. The variety of the texture is the visible counterpart of the dimension $8$: the orbits move in and out of the quaternionic half. The comparison with Plate~N10, where the same coloring is applied to an octonionic body, is instructive.
\clearpage

\subsubsection{N06. Hidden part: final value}\label{atl:pl:N06}\mbox{}\par\vspace{0.5mm}
\noindent\includegraphics[width=\textwidth]{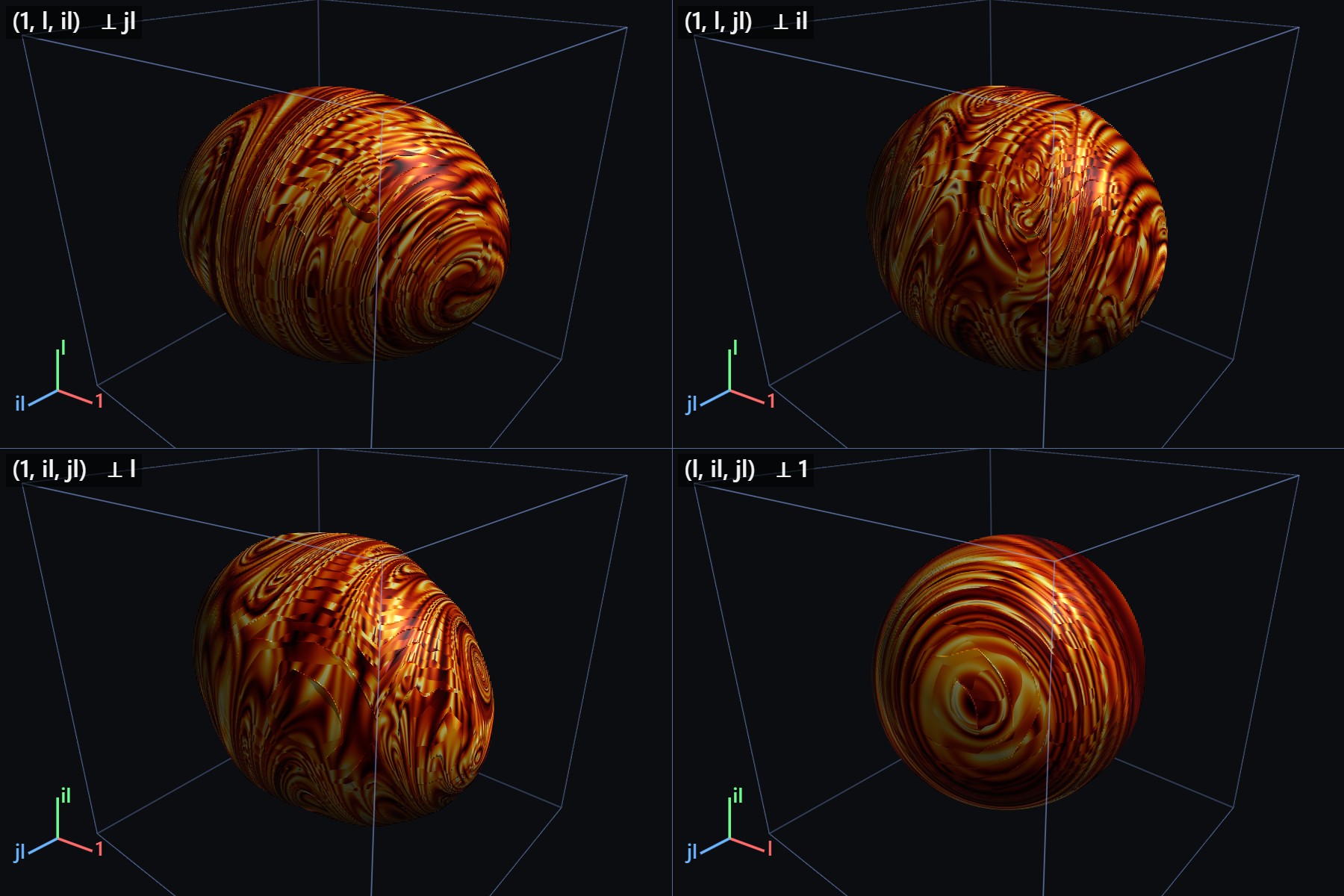}\par\vspace{1mm}
\noindent{\footnotesize \textbf{Map} $z\mapsto z\diamond z+c$ $\cdot$ \textbf{Algebra} Okubo algebra (Petersson, $\gamma=(-1,-1,-1)$, $\theta=120^\circ$) $\cdot$ $c=(-0.397,\,\allowbreak -0.021,\,\allowbreak 0.049,\,\allowbreak 0.016,\,\allowbreak 0.004,\,\allowbreak 0.16,\,\allowbreak -0.113,\,\allowbreak 0.133)$ ($=0.95\,s^*$, $s^*=0.4900$) $\cdot$ \textbf{Frame} $X,Y,Z,N=1,l,il,jl$ $\cdot$ half-extent 1.4 $\cdot$ 20 iterations $\cdot$ \textbf{Coloring} hidden part (final value), gradient Fire}\par\vspace{1mm}
\noindent{\footnotesize\begin{tabular}{@{}lllll@{}}\toprule Slice (X, Y, Z) & Body [95\%] & Exact ($J\geq0.9995$) & Approximate ($0.97\leq J<0.9995$) & Dim\\\midrule
($1, l, il$) $\perp jl$ & 52.7\% {\scriptsize[52.4, 53.0]} & $P$ & $M$: y; $R_y$: 180 (min.\ 0.986) & 8 \\
($1, l, jl$) $\perp il$ & 52.0\% {\scriptsize[51.7, 52.3]} & $P$ & -- & 8 \\
($1, il, jl$) $\perp l$ & 51.4\% {\scriptsize[51.1, 51.8]} & $P$ & -- & 8 \\
($l, il, jl$) $\perp 1$ & 42.2\% {\scriptsize[41.9, 42.5]} & $P$ & $M$: x; $R_x$: 180,37; $R_y$: 90,72,60,37; $R_z$: 90,72,60,37 (min.\ 0.971) & 8 \\
\bottomrule\end{tabular}}\par\vspace{2mm}
The hidden part of the final orbit value alone, instead of the average of Plate~N05, gives a much more contrasted texture on the same body (bodies of $53\%$, $52\%$, $51\%$ and $42\%$, only $P$ exact, dimension $8$). With the Fire gradient the surfaces are covered by bands of orange, yellow and dark brown, with swirls and small vortices in the quarters through $1$; the imaginary quarter shows concentric rings around a dark vortex on its visible face.

The comparison with Plate~N05 shows how the average over the orbit smooths the sensitive dependence of the final value on the initial point: the value at a single iterate changes quickly from point to point, and the bands become narrow and numerous.
\clearpage

\subsubsection{N07. Ultra Fractal smooth coloring}\label{atl:pl:N07}\mbox{}\par\vspace{0.5mm}
\noindent\includegraphics[width=\textwidth]{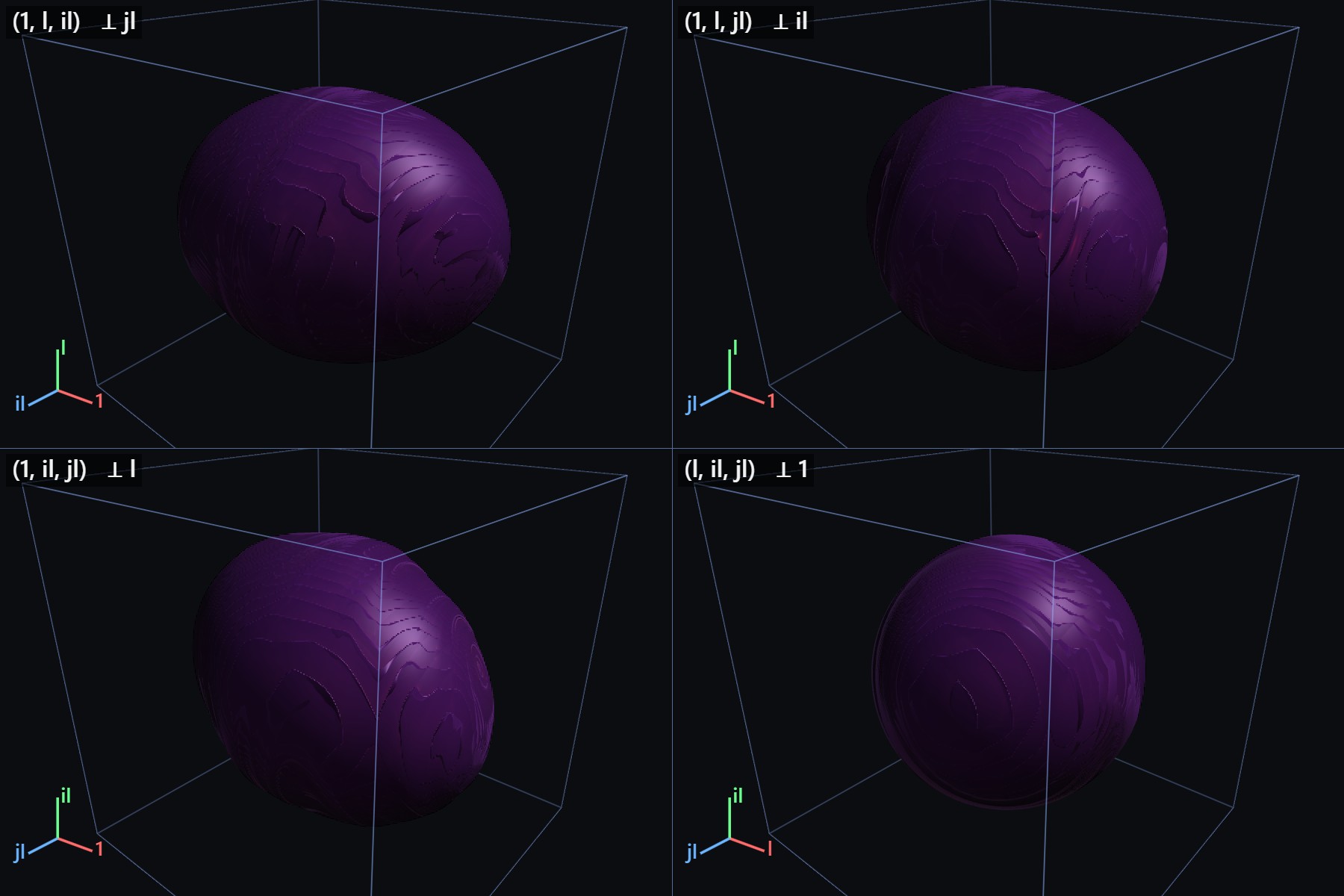}\par\vspace{1mm}
\noindent{\footnotesize \textbf{Map} $z\mapsto z\diamond z+c$ $\cdot$ \textbf{Algebra} Okubo algebra (Petersson, $\gamma=(-1,-1,-1)$, $\theta=120^\circ$) $\cdot$ $c=(-0.397,\,\allowbreak -0.021,\,\allowbreak 0.049,\,\allowbreak 0.016,\,\allowbreak 0.004,\,\allowbreak 0.16,\,\allowbreak -0.113,\,\allowbreak 0.133)$ ($=0.95\,s^*$, $s^*=0.4900$) $\cdot$ \textbf{Frame} $X,Y,Z,N=1,l,il,jl$ $\cdot$ half-extent 1.4 $\cdot$ 20 iterations $\cdot$ \textbf{Coloring} Ultra Fractal smooth (2D projection of $z$), gradient Sunset}\par\vspace{1mm}
\noindent{\footnotesize\begin{tabular}{@{}lllll@{}}\toprule Slice (X, Y, Z) & Body [95\%] & Exact ($J\geq0.9995$) & Approximate ($0.97\leq J<0.9995$) & Dim\\\midrule
($1, l, il$) $\perp jl$ & 52.7\% {\scriptsize[52.4, 53.0]} & $P$ & $M$: y; $R_y$: 180 (min.\ 0.986) & 8 \\
($1, l, jl$) $\perp il$ & 52.0\% {\scriptsize[51.7, 52.3]} & $P$ & -- & 8 \\
($1, il, jl$) $\perp l$ & 51.4\% {\scriptsize[51.1, 51.7]} & $P$ & -- & 8 \\
($l, il, jl$) $\perp 1$ & 42.1\% {\scriptsize[41.8, 42.4]} & $P$ & $M$: x; $R_x$: 180,37; $R_y$: 90,72,60,37; $R_z$: 90,72,60,37 (min.\ 0.972) & 8 \\
\bottomrule\end{tabular}}\par\vspace{2mm}
The Smooth coloring of the standard library of Ultra Fractal is a formula written for the complex plane, and in HyperFractal~8D it sees only a complex projection of the eight-dimensional orbit value. On the body of Plate~B07 (bodies of $53\%$, $52\%$, $51\%$ and $42\%$, only $P$ exact, dimension $8$) the result is an almost uniform dark violet, on which only the lighting and the terrace lines of the relief (Plate~N01) are visible. The plate shows the limits of a two-dimensional coloring applied to an eight-dimensional dynamics; the comparison with Plate~N01, which uses the full eight-dimensional norm in the same smooth formula, makes the point.
\clearpage

\subsubsection{N08. Quaternion direction}\label{atl:pl:N08}\mbox{}\par\vspace{0.5mm}
\noindent\includegraphics[width=\textwidth]{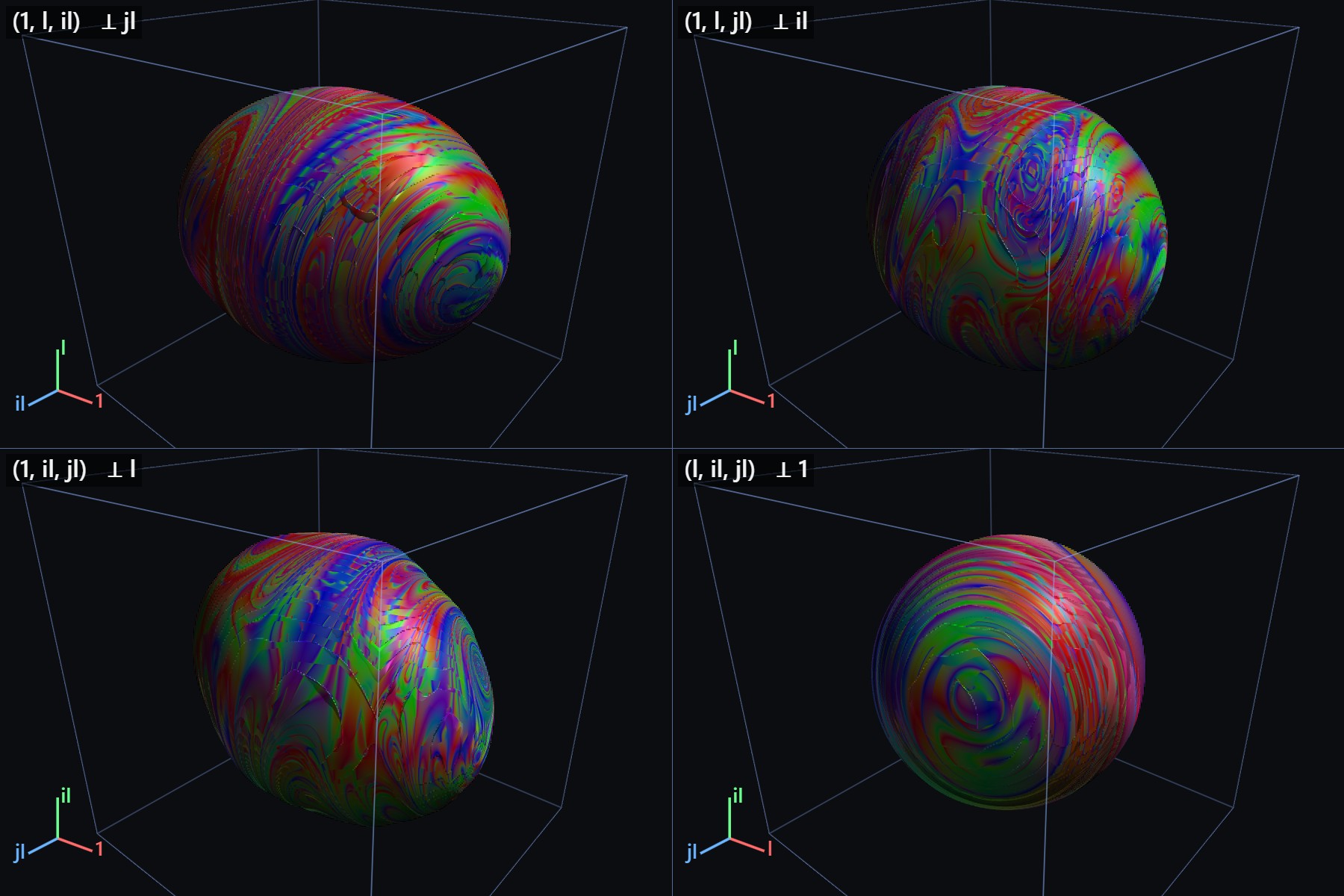}\par\vspace{1mm}
\noindent{\footnotesize \textbf{Map} $z\mapsto z\diamond z+c$ $\cdot$ \textbf{Algebra} Okubo algebra (Petersson, $\gamma=(-1,-1,-1)$, $\theta=120^\circ$) $\cdot$ $c=(-0.397,\,\allowbreak -0.021,\,\allowbreak 0.049,\,\allowbreak 0.016,\,\allowbreak 0.004,\,\allowbreak 0.16,\,\allowbreak -0.113,\,\allowbreak 0.133)$ ($=0.95\,s^*$, $s^*=0.4900$) $\cdot$ \textbf{Frame} $X,Y,Z,N=1,l,il,jl$ $\cdot$ half-extent 1.4 $\cdot$ 20 iterations $\cdot$ \textbf{Coloring} quaternion direction (4D, first four components), gradient Classic}\par\vspace{1mm}
\noindent{\footnotesize\begin{tabular}{@{}lllll@{}}\toprule Slice (X, Y, Z) & Body [95\%] & Exact ($J\geq0.9995$) & Approximate ($0.97\leq J<0.9995$) & Dim\\\midrule
($1, l, il$) $\perp jl$ & 52.8\% {\scriptsize[52.5, 53.2]} & $P$ & $M$: y; $R_y$: 180 (min.\ 0.986) & 8 \\
($1, l, jl$) $\perp il$ & 52.1\% {\scriptsize[51.8, 52.4]} & $P$ & -- & 8 \\
($1, il, jl$) $\perp l$ & 51.6\% {\scriptsize[51.3, 51.9]} & $P$ & -- & 8 \\
($l, il, jl$) $\perp 1$ & 42.2\% {\scriptsize[41.9, 42.5]} & $P$ & $M$: x; $R_x$: 180,37; $R_y$: 90,72,60,37; $R_z$: 90,72,60,37 (min.\ 0.972) & 8 \\
\bottomrule\end{tabular}}\par\vspace{2mm}
The quaternion direction coloring does not use the gradient: it assigns to each point a color whose red, green and blue components are the absolute values of the $i$, $j$ and $k$ components of the final orbit value, normalized by the length of that imaginary part. The body is that of Plate~B07 (bodies of $53\%$, $52\%$, $52\%$ and $42\%$, only $P$ exact, dimension $8$).

The bodies become iridescent: interleaved red, green, blue and violet bands, with swirls and small vortices in the quarters through $1$ and concentric rings in the imaginary quarter. The coloring ignores the real part and the four hidden components, yet it shows vividly how the direction of the orbit in the quaternionic half changes from point to point.
\clearpage

\subsubsection{N09. 4D orbit trap}\label{atl:pl:N09}\mbox{}\par\vspace{0.5mm}
\noindent\includegraphics[width=\textwidth]{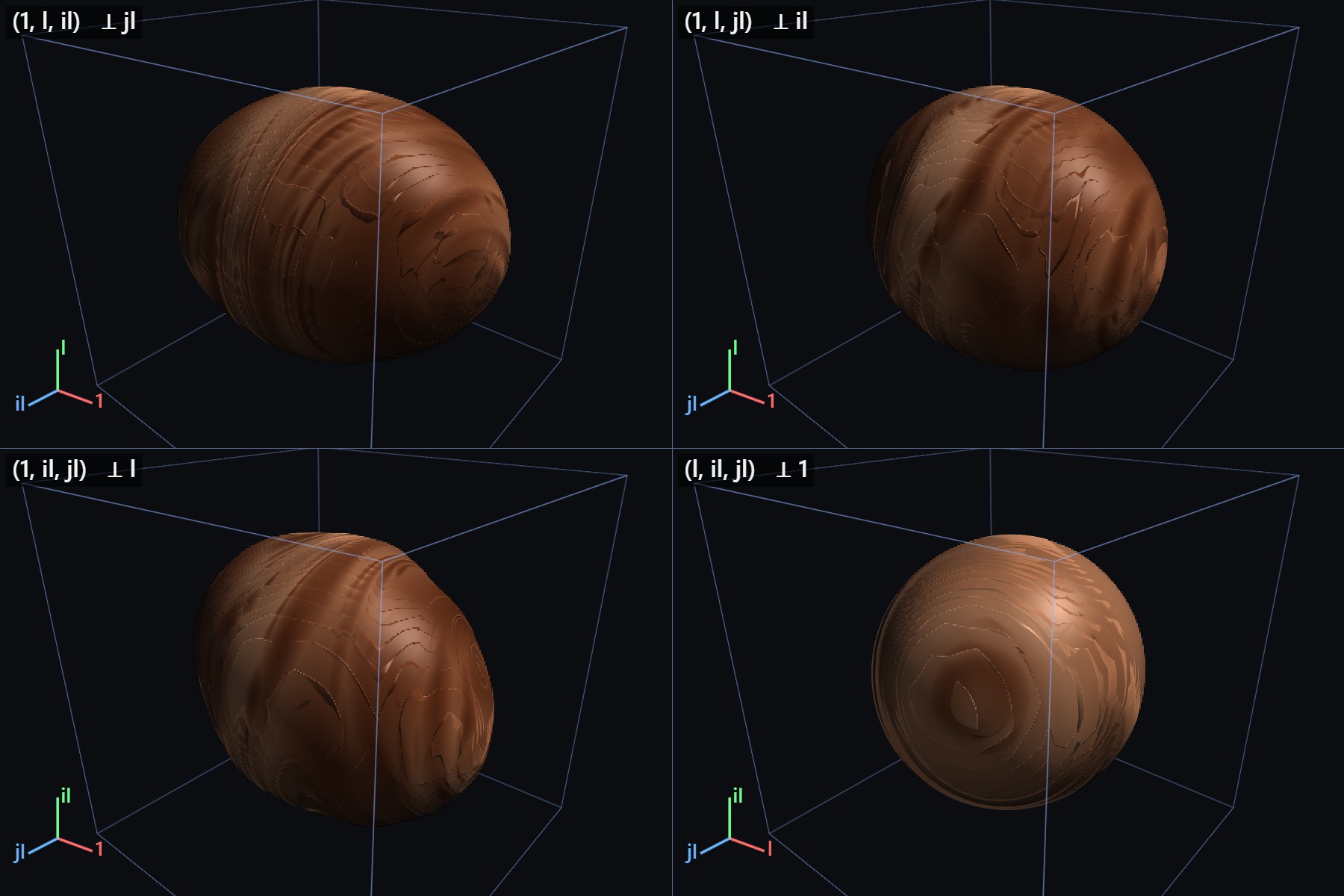}\par\vspace{1mm}
\noindent{\footnotesize \textbf{Map} $z\mapsto z\diamond z+c$ $\cdot$ \textbf{Algebra} Okubo algebra (Petersson, $\gamma=(-1,-1,-1)$, $\theta=120^\circ$) $\cdot$ $c=(-0.397,\,\allowbreak -0.021,\,\allowbreak 0.049,\,\allowbreak 0.016,\,\allowbreak 0.004,\,\allowbreak 0.16,\,\allowbreak -0.113,\,\allowbreak 0.133)$ ($=0.95\,s^*$, $s^*=0.4900$) $\cdot$ \textbf{Frame} $X,Y,Z,N=1,l,il,jl$ $\cdot$ half-extent 1.4 $\cdot$ 20 iterations $\cdot$ \textbf{Coloring} 4D orbit trap (first four components), gradient Copper}\par\vspace{1mm}
\noindent{\footnotesize\begin{tabular}{@{}lllll@{}}\toprule Slice (X, Y, Z) & Body [95\%] & Exact ($J\geq0.9995$) & Approximate ($0.97\leq J<0.9995$) & Dim\\\midrule
($1, l, il$) $\perp jl$ & 52.7\% {\scriptsize[52.4, 53.0]} & $P$ & $M$: y; $R_y$: 180 (min.\ 0.985) & 8 \\
($1, l, jl$) $\perp il$ & 52.0\% {\scriptsize[51.6, 52.3]} & $P$ & -- & 8 \\
($1, il, jl$) $\perp l$ & 51.3\% {\scriptsize[51.0, 51.7]} & $P$ & -- & 8 \\
($l, il, jl$) $\perp 1$ & 42.3\% {\scriptsize[42.0, 42.6]} & $P$ & $M$: x; $R_x$: 180,37; $R_y$: 90,72,60,37; $R_z$: 90,72,60,37 (min.\ 0.973) & 8 \\
\bottomrule\end{tabular}}\par\vspace{2mm}
The four-dimensional orbit trap of HyperFractal Studio uses only the first four components $1,i,j,k$ of the orbit, and colors by their minimum distance from the origin. The body is that of Plate~B07 (bodies of $53\%$, $52\%$, $51\%$ and $42\%$, only $P$ exact, dimension $8$).

With the Copper gradient the surfaces are copper-brown, with broad darker streaks and bright highlights that follow the curvature of the bodies, and concentric bands in the imaginary quarter. The coloring varies more than the eight-dimensional trap of Plate~N02, which is almost uniform: the distance of the quaternionic part of the orbit from the origin can become small even when the full orbit stays away from it, and the difference between the two plates is a rough indication of how much of the Okubonic dynamics takes place in the hidden half.
\clearpage

\subsubsection{N10. Hidden part on an octonionic body}\label{atl:pl:N10}\mbox{}\par\vspace{0.5mm}
\noindent\includegraphics[width=\textwidth]{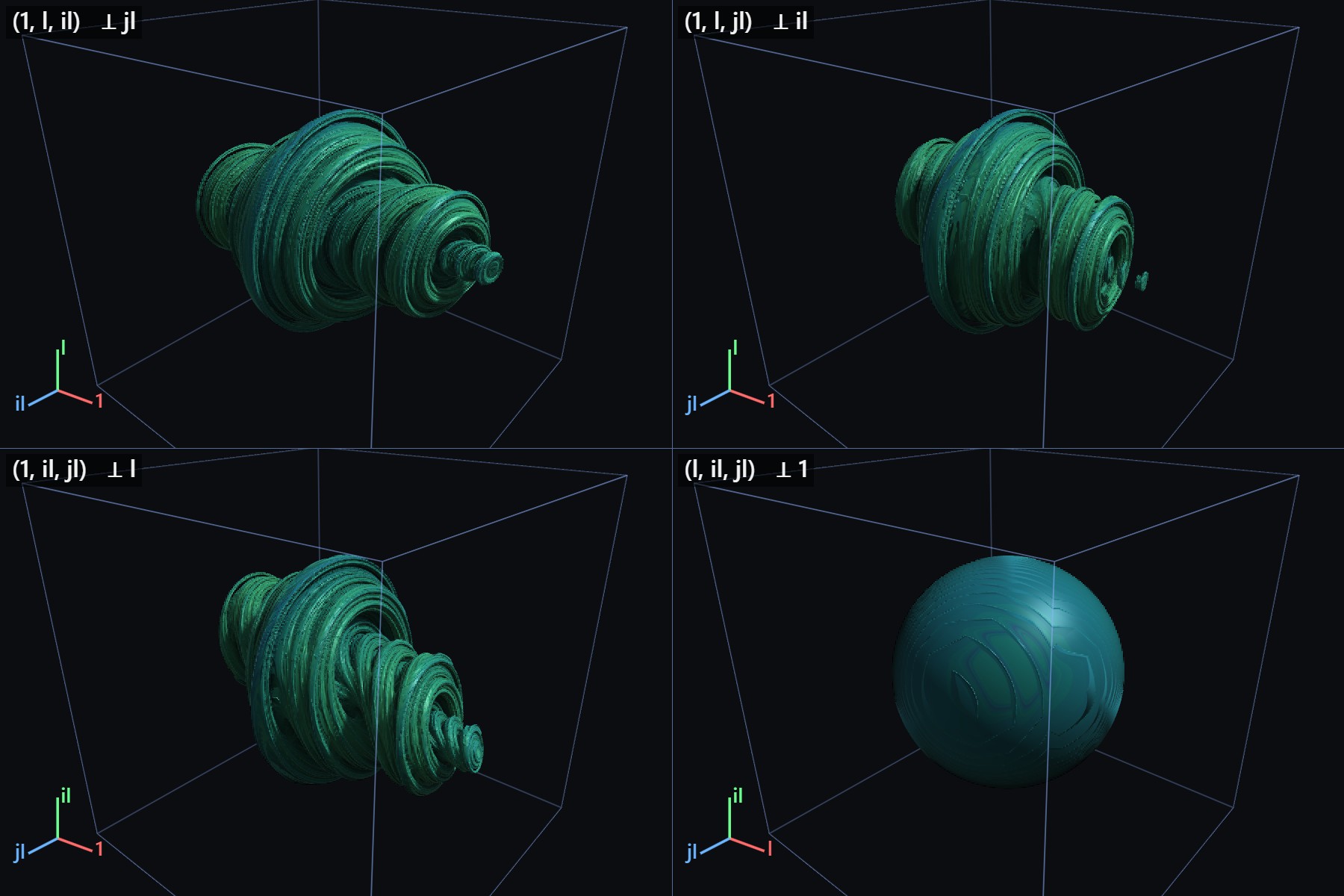}\par\vspace{1mm}
\noindent{\footnotesize \textbf{Map} $z\mapsto z\diamond z+c$ $\cdot$ \textbf{Algebra} Octonions (Cayley--Dickson, $\gamma=(-1,-1,-1)$) $\cdot$ $c=(-0.598,\,\allowbreak -0.032,\,\allowbreak 0.075,\,\allowbreak 0.025,\,\allowbreak 0.006,\,\allowbreak 0.241,\,\allowbreak -0.17,\,\allowbreak 0.2)$ ($=0.95\,s^*$, $s^*=0.7383$) $\cdot$ \textbf{Frame} $X,Y,Z,N=1,l,il,jl$ $\cdot$ half-extent 1.4 $\cdot$ 20 iterations $\cdot$ \textbf{Coloring} hidden part (orbit average), gradient Viridis-like}\par\vspace{1mm}
\noindent{\footnotesize\begin{tabular}{@{}lllll@{}}\toprule Slice (X, Y, Z) & Body [95\%] & Exact ($J\geq0.9995$) & Approximate ($0.97\leq J<0.9995$) & Dim\\\midrule
($1, l, il$) $\perp jl$ & 16.7\% {\scriptsize[16.5, 16.9]} & $P$ & $M$: y; $R_y$: 180 (min.\ 0.973) & 3 \\
($1, l, jl$) $\perp il$ & 16.5\% {\scriptsize[16.2, 16.7]} & $P$ & $M$: y; $R_y$: 180 (min.\ 0.976) & 3 \\
($1, il, jl$) $\perp l$ & 17.3\% {\scriptsize[17.0, 17.5]} & $P$ & -- & 3 \\
($l, il, jl$) $\perp 1$ & 22.1\% {\scriptsize[21.8, 22.3]} & $M$: xyz; $R_x$: all; $R_y$: all; $R_z$: all; $P$ & -- & 3 \\
\bottomrule\end{tabular}}\par\vspace{2mm}
The hidden part coloring of Plate~N05 is here applied to an octonionic body: the Julia set of Plate~B01, with the generic constant in the frame $F_B$. The three quarters through $1$ show spindle-shaped bodies built of stacked, ribbed rings that taper towards one end ($17\%$ of the ball each); the upper-right quarter shows also a small detached fragment near the tip. The imaginary quarter is a ball ($22\%$), exact under every tested rotation, as predicted by Proposition~\ref{atl:prop:imaginary}. In the quarters through $1$ only $P$ is exact, with approximate mirrors $M_y$ in $(1,l,il)$ and $(1,l,jl)$ ($J=0.973$ and $0.976$).

The coloring is a nearly uniform green-teal, with little of the marbling of Plate~N05. For the octonions the orbits lie in a subalgebra of dimension $3$, the span of $1$, $\im z_0$ and $\im c$ (Proposition~\ref{atl:prop:span}), while for the Okubo algebra of Plate~N05 the dimension is $8$: the octonionic orbit has far fewer ways to move in and out of the quaternionic half, and the hidden part varies less over the surface.
\clearpage

\subsection{Part O. Moving along the hidden axis}\label{atl:part:O}
The four slices are sections through the center of the box. Here the center is moved along the hidden axis $l$, so that the three quarters containing $l$ are simply translated, while the quarter orthogonal to $l$ becomes a new section.

\subsubsection{O01. Okubo, section at \texorpdfstring{$l=0.3$}{l=0.3}}\label{atl:pl:O01}\mbox{}\par\vspace{0.5mm}
\noindent\includegraphics[width=\textwidth]{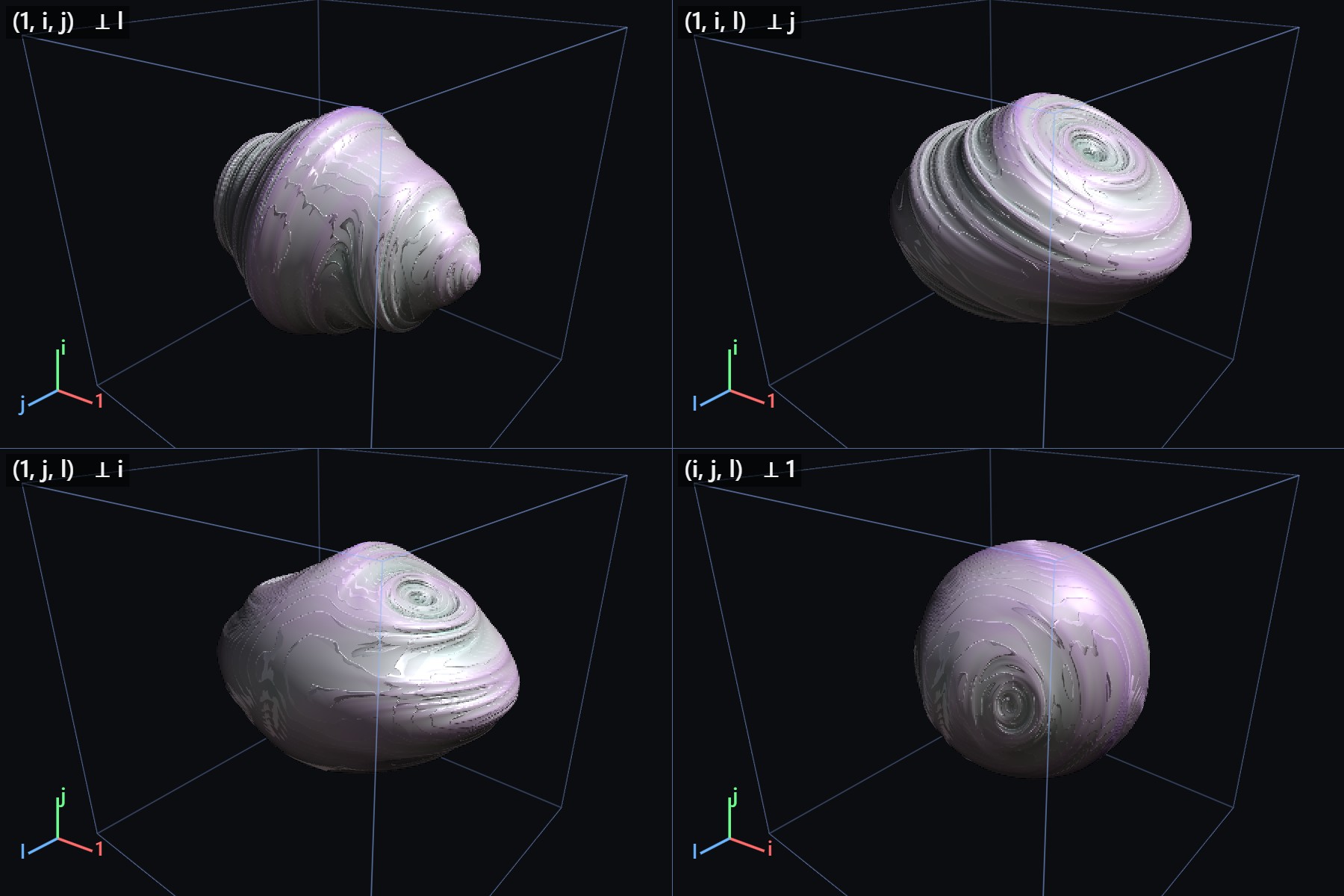}\par\vspace{1mm}
\noindent{\footnotesize \textbf{Map} $z\mapsto z\diamond z+c$ $\cdot$ \textbf{Algebra} Okubo algebra (Petersson, $\gamma=(-1,-1,-1)$, $\theta=120^\circ$) $\cdot$ $c=(-0.397,\,\allowbreak -0.021,\,\allowbreak 0.049,\,\allowbreak 0.016,\,\allowbreak 0.004,\,\allowbreak 0.16,\,\allowbreak -0.113,\,\allowbreak 0.133)$ ($=0.95\,s^*$, $s^*=0.4900$) $\cdot$ \textbf{Frame} $X,Y,Z,N=1,i,j,l$ $\cdot$ center $(0,\,\allowbreak 0,\,\allowbreak 0,\,\allowbreak 0,\,\allowbreak 0.3,\,\allowbreak 0,\,\allowbreak 0,\,\allowbreak 0)$ $\cdot$ half-extent 1.4 $\cdot$ 20 iterations $\cdot$ \textbf{Coloring} 8D orbit trap (point), gradient Aurora}\par\vspace{1mm}
\noindent{\footnotesize\begin{tabular}{@{}lllll@{}}\toprule Slice (X, Y, Z) & Body [95\%] & Exact ($J\geq0.9995$) & Approximate ($0.97\leq J<0.9995$) & Dim\\\midrule
($1, i, j$) $\perp l$ & 26.9\% {\scriptsize[26.6, 27.2]} & -- & -- & 8 \\
($1, i, l$) $\perp j$ & 39.1\% {\scriptsize[38.8, 39.4]} & -- & $R_z$: 180 (min.\ 0.974) & 8 \\
($1, j, l$) $\perp i$ & 40.0\% {\scriptsize[39.7, 40.3]} & -- & -- & 8 \\
($i, j, l$) $\perp 1$ & 33.5\% {\scriptsize[33.3, 33.8]} & -- & -- & 8 \\
\bottomrule\end{tabular}}\par\vspace{2mm}
The three plates of Part~O move the centre of the box away from the origin along $l$. Moving it to $l=0.3$, the upper-left quarter $(1,i,j)$ becomes the section of the Okubonic body of Plate~M06 (the generic constant in the frame $F_A$) at the depth $l=0.3$: a ribbed, shell-like body, somewhat smoother on its right side ($27\%$ of the ball). The three quarters containing $l$ show the same sections as in Plate~M06, since the section is the same three-dimensional subspace, but seen in a window moved by $0.3$ along the local $l$-axis: a flattened body with a vortex in $(1,i,l)$, a rounded cone with a vortex in $(1,j,l)$ and a round body with a vortex in the imaginary quarter ($39\%$, $40\%$ and $34\%$).

These sections do not pass through the origin, and the universal symmetry $z\mapsto-z$ of Remark~\ref{atl:rem:minus} does not map them to themselves: the point reflection $P$ through the centre of the box is not a symmetry, and indeed the measurements find no exact symmetry in any quarter. The only approximate one is the half-turn $R_z(180^\circ)$ about the $l$-axis of the quarter $(1,i,l)$ ($J=0.974$): a rotation about an axis parallel to $l$ commutes with the translation along $l$, and it is the composition of $P$ with the mirror in $l$, which is itself only approximate for the central section (it scored $0.98$ in the render-time measurements of version~1 of Plate~M06, $2400$ points, indicative only; Part~M is not measured in version~2). The orbits lie in a subalgebra of dimension $8$.
\clearpage

\subsubsection{O02. Okubo, section at \texorpdfstring{$l=0.6$}{l=0.6}}\label{atl:pl:O02}\mbox{}\par\vspace{0.5mm}
\noindent\includegraphics[width=\textwidth]{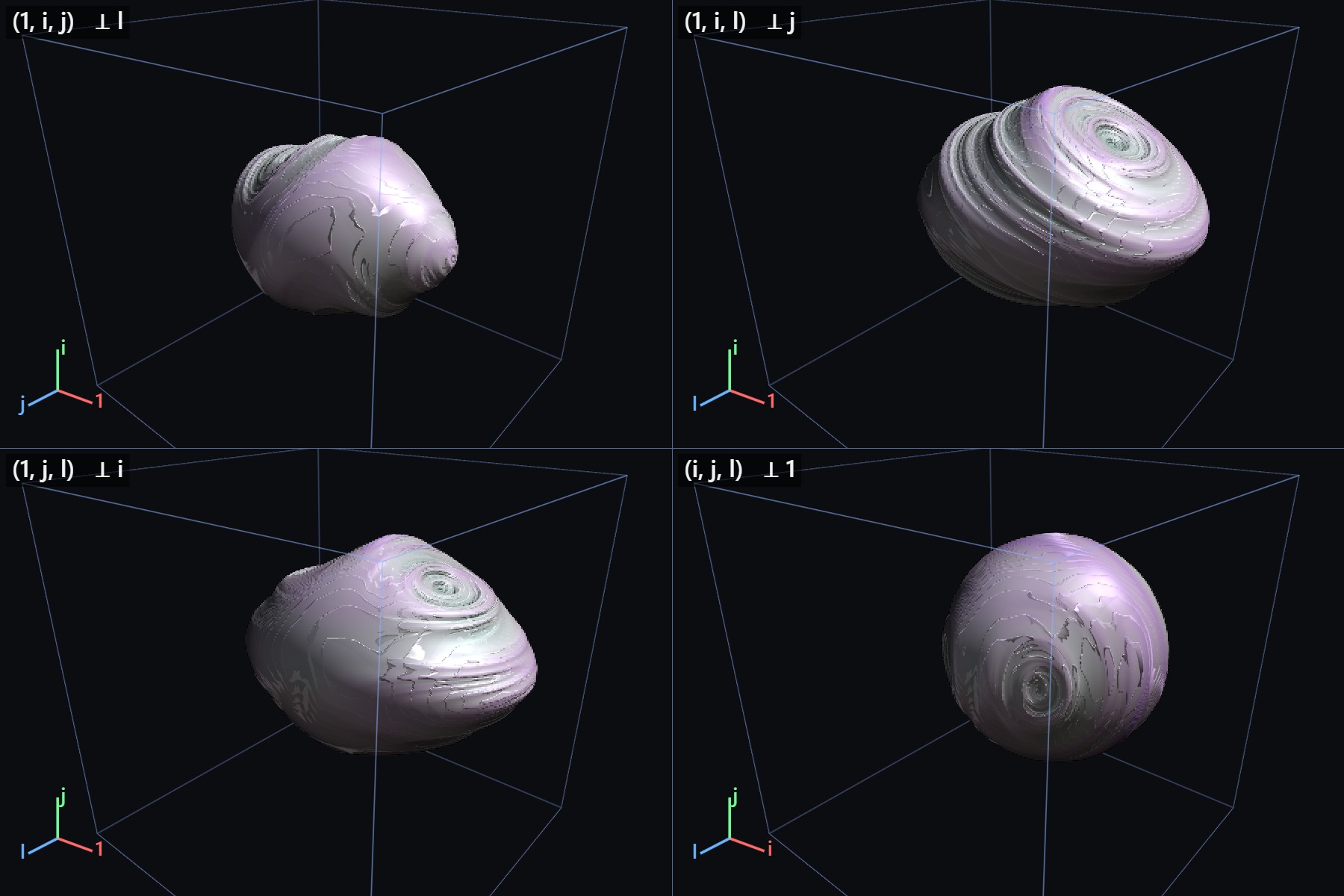}\par\vspace{1mm}
\noindent{\footnotesize \textbf{Map} $z\mapsto z\diamond z+c$ $\cdot$ \textbf{Algebra} Okubo algebra (Petersson, $\gamma=(-1,-1,-1)$, $\theta=120^\circ$) $\cdot$ $c=(-0.397,\,\allowbreak -0.021,\,\allowbreak 0.049,\,\allowbreak 0.016,\,\allowbreak 0.004,\,\allowbreak 0.16,\,\allowbreak -0.113,\,\allowbreak 0.133)$ ($=0.95\,s^*$, $s^*=0.4900$) $\cdot$ \textbf{Frame} $X,Y,Z,N=1,i,j,l$ $\cdot$ center $(0,\,\allowbreak 0,\,\allowbreak 0,\,\allowbreak 0,\,\allowbreak 0.6,\,\allowbreak 0,\,\allowbreak 0,\,\allowbreak 0)$ $\cdot$ half-extent 1.4 $\cdot$ 20 iterations $\cdot$ \textbf{Coloring} 8D orbit trap (point), gradient Aurora}\par\vspace{1mm}
\noindent{\footnotesize\begin{tabular}{@{}lllll@{}}\toprule Slice (X, Y, Z) & Body [95\%] & Exact ($J\geq0.9995$) & Approximate ($0.97\leq J<0.9995$) & Dim\\\midrule
($1, i, j$) $\perp l$ & 17.1\% {\scriptsize[16.9, 17.3]} & -- & -- & 8 \\
($1, i, l$) $\perp j$ & 33.8\% {\scriptsize[33.6, 34.1]} & -- & $R_z$: 180 (min.\ 0.982) & 8 \\
($1, j, l$) $\perp i$ & 34.4\% {\scriptsize[34.1, 34.7]} & -- & -- & 8 \\
($i, j, l$) $\perp 1$ & 30.8\% {\scriptsize[30.5, 31.1]} & -- & -- & 8 \\
\bottomrule\end{tabular}}\par\vspace{2mm}
At the depth $l=0.6$ the quaternionic quarter shows a smaller and rounder body ($17\%$ of the ball, against $27\%$ at $l=0.3$), still ribbed on one flank and ending in a small tip: the section shrinks as it moves away from the centre of the body. The other three quarters are again the sections of Plate~M06 seen in a translated window; their shares ($34\%$, $34\%$ and $31\%$) are smaller than in Plate~O01 because the ball of the measurement, centred at $l=0.6$, now misses part of the body.

As in Plate~O01 there is no exact symmetry, since the sections do not pass through the origin; the half-turn $R_z(180^\circ)$ of the quarter $(1,i,l)$ is approximate ($J=0.982$), and the orbits lie in a subalgebra of dimension $8$. Plates~M06, O01 and O02 form a tomography of the body along $l$.
\clearpage

\subsubsection{O03. Octonions, section at \texorpdfstring{$l=0.4$}{l=0.4}}\label{atl:pl:O03}\mbox{}\par\vspace{0.5mm}
\noindent\includegraphics[width=\textwidth]{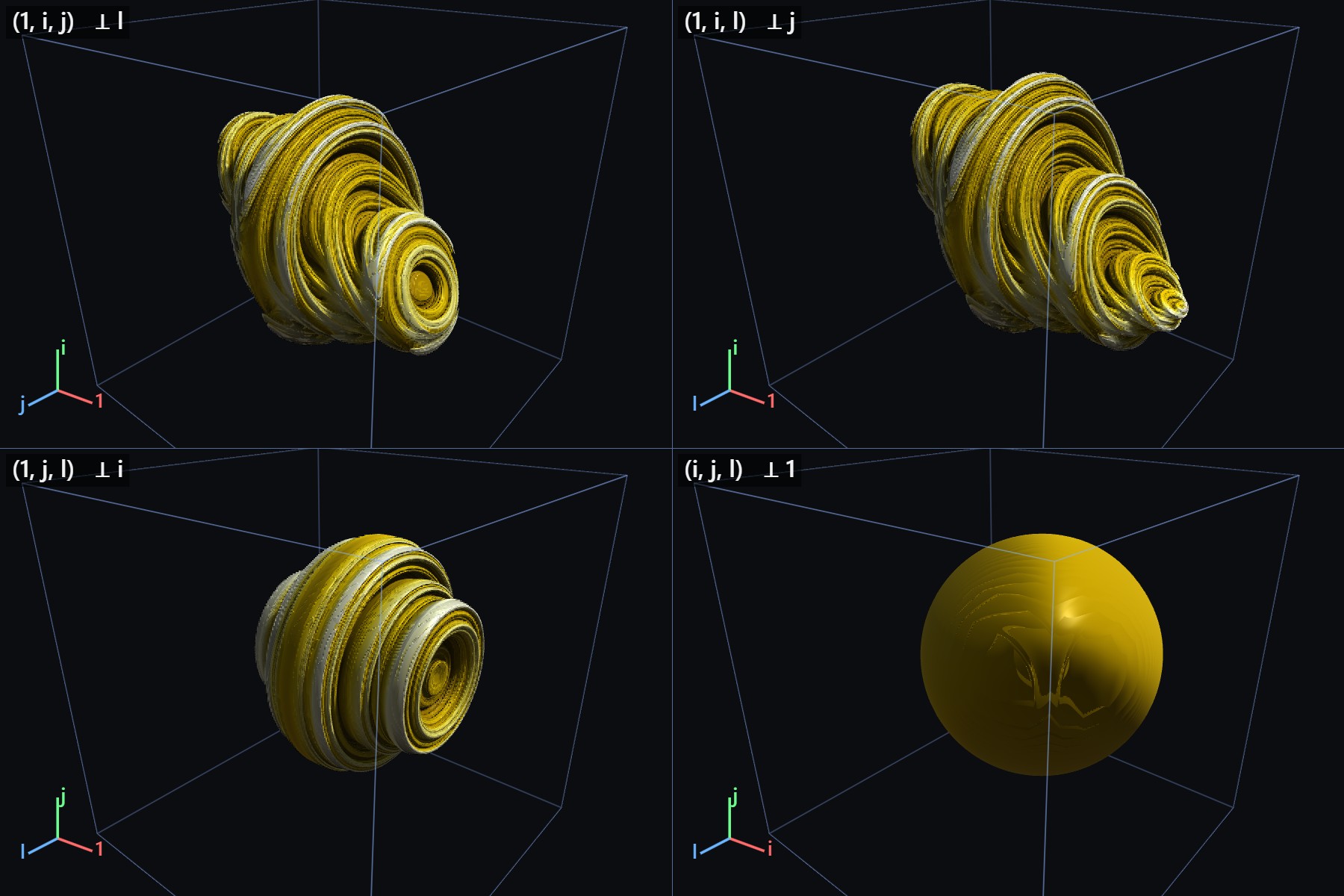}\par\vspace{1mm}
\noindent{\footnotesize \textbf{Map} $z\mapsto z\diamond z+c$ $\cdot$ \textbf{Algebra} Octonions (Cayley--Dickson, $\gamma=(-1,-1,-1)$) $\cdot$ $c=(-0.231,\,\allowbreak 0.62,\,\allowbreak 0,\,\allowbreak 0,\,\allowbreak 0,\,\allowbreak 0,\,\allowbreak 0,\,\allowbreak 0)$ ($=0.97\,s^*$, $s^*=0.6824$) $\cdot$ \textbf{Frame} $X,Y,Z,N=1,i,j,l$ $\cdot$ center $(0,\,\allowbreak 0,\,\allowbreak 0,\,\allowbreak 0,\,\allowbreak 0.4,\,\allowbreak 0,\,\allowbreak 0,\,\allowbreak 0)$ $\cdot$ half-extent 1.4 $\cdot$ 20 iterations $\cdot$ \textbf{Coloring} 8D orbit trap (point), gradient Gold}\par\vspace{1mm}
\noindent{\footnotesize\begin{tabular}{@{}lllll@{}}\toprule Slice (X, Y, Z) & Body [95\%] & Exact ($J\geq0.9995$) & Approximate ($0.97\leq J<0.9995$) & Dim\\\midrule
($1, i, j$) $\perp l$ & 17.5\% {\scriptsize[17.3, 17.7]} & $M$: z; $R_z$: 180; $P$ & -- & 3 \\
($1, i, l$) $\perp j$ & 24.3\% {\scriptsize[24.1, 24.6]} & $R_z$: 180 & -- & 3 \\
($1, j, l$) $\perp i$ & 22.8\% {\scriptsize[22.5, 23.1]} & $M$: xy; $R_z$: 180 & -- & 3 \\
($i, j, l$) $\perp 1$ & 34.9\% {\scriptsize[34.6, 35.2]} & $M$: xy; $R_z$: all & -- & 3 \\
\bottomrule\end{tabular}}\par\vspace{2mm}
For the octonions the section at the depth $l=0.4$ can be predicted exactly. By Proposition~\ref{atl:prop:equivariance} a point $x_0+x_1i+x_2j+0.4\,l$ has the same orbit, up to a rotation fixing $1$ and $i$, as $x_0+x_1i+\sqrt{x_2^2+0.16}\;j$. The section at depth $0.4$ is therefore obtained from the central section of Plate~A01 by keeping the part with $|x_2|\ge0.4$ and pushing it inwards along $j$. The quaternionic quarter confirms it: a shorter spiral shell ($18\%$ of the ball, against $25\%$ in Plate~A01). The other quarters show the bodies of Plate~A01 in a window moved by $0.4$ along $l$ (shares $24\%$, $23\%$ and $35\%$, as in Plate~A01).

The symmetries are instructive, because they show which part of the octonionic symmetry survives the translation. The point reflection $z\mapsto-z$ moves the depth $0.4$ to $-0.4$; composed with the reflection $l\mapsto-l$, which belongs to the group of Proposition~\ref{atl:prop:equivariance}, it becomes a symmetry of the translated section. Accordingly, the quaternionic quarter keeps $M_z$, $R_z(180^\circ)$ and $P$ exactly. In the quarters containing $l$ exactly those symmetries survive that fix the $l$-coordinate: in $(1,i,l)$ the half-turn about $l$; in $(1,j,l)$ the mirrors in $1$ and $j$ and the half-turn about $l$, while the revolution about $1$ of Plate~A01 is lost, since its axis no longer passes through the centre of the box; in the imaginary quarter the mirrors in $i$ and $j$ and the rotations about $l$ by every tested angle. The ball of the imaginary quarter, seen off-centre, keeps its revolution about the $l$-axis but loses its mirror in $l$ and its point symmetry. The orbits lie in a subalgebra of dimension $3$. The contrast with Plates~O01 and O02, where nothing survives, reflects the much larger symmetry group of the octonionic map.
\clearpage

\subsection{Part P. Further algebras and frames}\label{atl:part:P}
We close the plates of the first version with algebras and frames not covered by the previous parts: degenerate para-Hurwitz and Petersson algebras, an anisotropic octonion algebra, a split Petersson algebra, and the purely imaginary frame $F_D=(i,j,l\,|\,il)$ in the compact and split Okubo algebras.

\subsubsection{P01. Semi-para-octonions}\label{atl:pl:P01}\mbox{}\par\vspace{0.5mm}
\noindent\includegraphics[width=\textwidth]{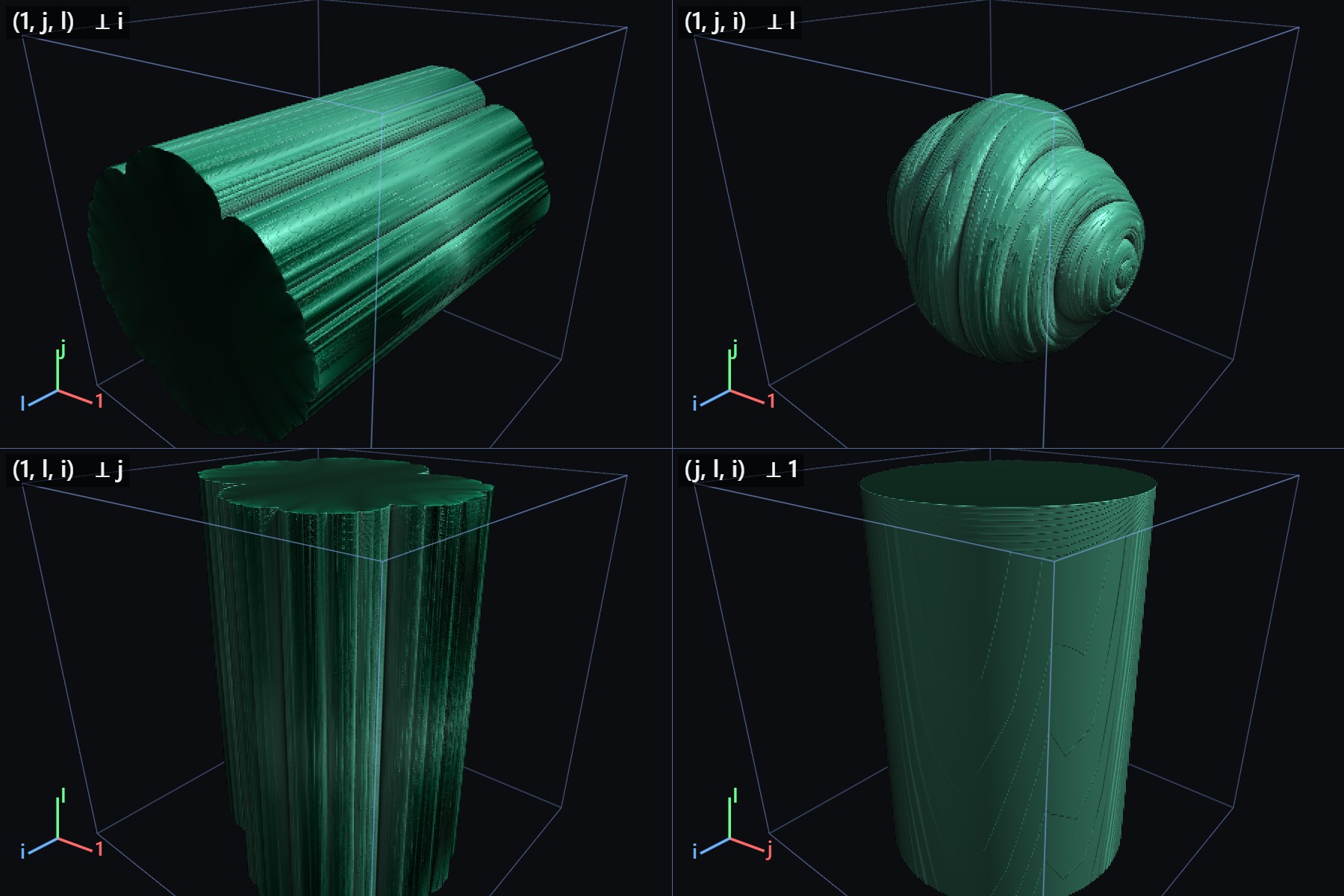}\par\vspace{1mm}
\noindent{\footnotesize \textbf{Map} $z\mapsto z\diamond z+c$ $\cdot$ \textbf{Algebra} custom: para-Hurwitz $\bar x\bar y$, $\gamma=(-1,-1,0)$ $\cdot$ $c=(-0.086,\,\allowbreak 0.231,\,\allowbreak 0,\,\allowbreak 0,\,\allowbreak 0,\,\allowbreak 0,\,\allowbreak 0,\,\allowbreak 0)$ ($=0.97\,s^*$, $s^*=0.2540$) $\cdot$ \textbf{Frame} $X,Y,Z,N=1,j,l,i$ $\cdot$ half-extent 1.4 $\cdot$ 20 iterations $\cdot$ \textbf{Coloring} 8D orbit trap (point), gradient Emerald}\par\vspace{1mm}
\noindent{\footnotesize\begin{tabular}{@{}lllll@{}}\toprule Slice (X, Y, Z) & Body [95\%] & Exact ($J\geq0.9995$) & Approximate ($0.97\leq J<0.9995$) & Dim\\\midrule
($1, j, l$) $\perp i$ & 66.7\% {\scriptsize[66.4, 67.0]} & $M$: xyz; $R_x$: 180; $R_y$: 180; $R_z$: 180; $P$ & -- & 3 \\
($1, j, i$) $\perp l$ & 35.8\% {\scriptsize[35.5, 36.1]} & $M$: y; $R_y$: 180; $P$ & -- & 3 \\
($1, l, i$) $\perp j$ & 66.7\% {\scriptsize[66.4, 67.0]} & $M$: y; $R_y$: 180; $P$ & -- & 3 \\
($j, l, i$) $\perp 1$ & 64.9\% {\scriptsize[64.6, 65.2]} & $M$: xyz; $R_x$: 180; $R_y$: all; $R_z$: 180; $P$ & -- & 3 \\
\bottomrule\end{tabular}}\par\vspace{2mm}
The degenerate para-octonions ($\gamma_3=0$ with the para-Hurwitz product) are shown in the frame $F_C$, with the complex constant. The quarter $(1,j,l)$ is a horizontal fluted cylinder along $l$ ($67\%$ of the ball), with a scalloped cross-section; the quarter $(1,j,i)$ is the para-quaternionic shell ($36\%$), the same body as the quaternionic quarter of the para-octonions in Plate~A05 (also $36\%$); the quarter $(1,l,i)$ is a vertical column along $l$ with a scalloped top face ($67\%$), whose section at $l=0$ is the filled Julia set of the Mandelbar map $\bar z^2+c$ in the plane $\CC$; the imaginary quarter $(j,l,i)$ is a smooth circular cylinder along $l$ ($65\%$).

The measurements confirm the cylinders: the imaginary quarter is exactly invariant under the rotations about $l$ by every tested angle, and the mirror in $l$ is exact wherever $l$ is an axis. The quarter $(1,j,l)$ has the three mirrors and the three half-turns, the quarters $(1,j,i)$ and $(1,l,i)$ the mirror in their second axis, the corresponding half-turn and $P$ (Remark~\ref{atl:rem:minus}). The orbits lie in a subalgebra of dimension $3$. The cylinders are those of the semi-octonions (Plate~A04): the para-Hurwitz conjugations do not change the mechanism of Proposition~\ref{atl:prop:rescaling}, by which the $l$-part of the orbit is driven linearly by the quaternionic part.
\clearpage

\subsubsection{P02. Semi-Okubo algebra}\label{atl:pl:P02}\mbox{}\par\vspace{0.5mm}
\noindent\includegraphics[width=\textwidth]{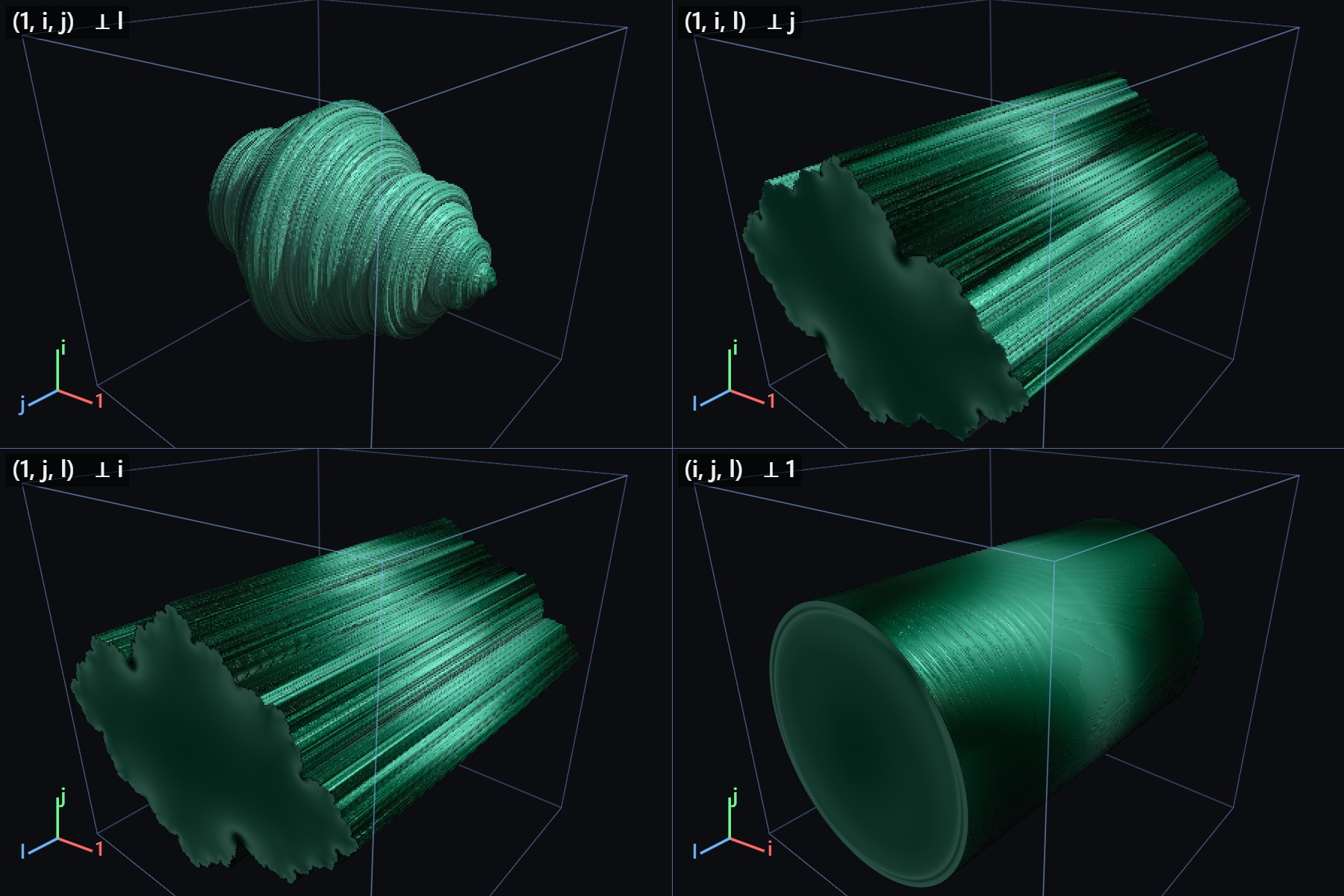}\par\vspace{1mm}
\noindent{\footnotesize \textbf{Map} $z\mapsto z\diamond z+c$ $\cdot$ \textbf{Algebra} custom: Petersson, $\gamma=(-1,-1,0)$, $\theta=120^\circ$ $\cdot$ $c=(-0.433,\,\allowbreak -0.023,\,\allowbreak 0.054,\,\allowbreak 0.018,\,\allowbreak 0.004,\,\allowbreak 0.174,\,\allowbreak -0.123,\,\allowbreak 0.145)$ ($=0.95\,s^*$, $s^*=0.5350$) $\cdot$ \textbf{Frame} $X,Y,Z,N=1,i,j,l$ $\cdot$ half-extent 1.4 $\cdot$ 20 iterations $\cdot$ \textbf{Coloring} 8D orbit trap (point), gradient Emerald}\par\vspace{1mm}
\noindent{\footnotesize\begin{tabular}{@{}lllll@{}}\toprule Slice (X, Y, Z) & Body [95\%] & Exact ($J\geq0.9995$) & Approximate ($0.97\leq J<0.9995$) & Dim\\\midrule
($1, i, j$) $\perp l$ & 28.0\% {\scriptsize[27.7, 28.3]} & $P$ & -- & 7 \\
($1, i, l$) $\perp j$ & 62.4\% {\scriptsize[62.1, 62.7]} & $M$: z; $R_z$: 180; $P$ & -- & 7 \\
($1, j, l$) $\perp i$ & 62.4\% {\scriptsize[62.1, 62.7]} & $M$: z; $R_z$: 180; $P$ & -- & 7 \\
($i, j, l$) $\perp 1$ & 61.6\% {\scriptsize[61.3, 61.9]} & $M$: xyz; $R_x$: 180; $R_y$: 180; $R_z$: all; $P$ & -- & 7 \\
\bottomrule\end{tabular}}\par\vspace{2mm}
The degenerate Okubo algebra ($\gamma_3=0$, $\theta=120^\circ$) with the generic constant is shown here in the frame $F_A$; the plate was rendered in a second run (see the reproducibility appendix), because in the frame $F_B$ every quarter contains at least two degenerate directions and the cube was entirely filled, and the measurements refer to this second run. The quaternionic quarter $(1,i,j)$ is a ribbed conical shell ($28\%$ of the ball). The quarters $(1,i,l)$ and $(1,j,l)$ are columns along $l$ with irregular, jagged cross-sections ($62\%$ each), and the imaginary quarter is a smooth circular cylinder along $l$ with a flat, rimmed base ($62\%$).

The measurements confirm the columns. The imaginary quarter is exactly invariant under the rotations about $l$ by every tested angle, and under the three mirrors; in $(1,i,l)$ and $(1,j,l)$ the mirror in $l$ and the half-turn about $l$ are exact, together with $P$, which is the only exact symmetry of the quaternionic quarter. The degenerate directions play the same role as in the semi-octonions (Plate~A04, Proposition~\ref{atl:prop:rescaling}), and the Petersson twist does not destroy the cylinders. The orbits lie in a subalgebra of dimension $7$, one less than for the Okubo algebra.
\clearpage

\subsubsection{P03. An anisotropic octonion algebra}\label{atl:pl:P03}\mbox{}\par\vspace{0.5mm}
\noindent\includegraphics[width=\textwidth]{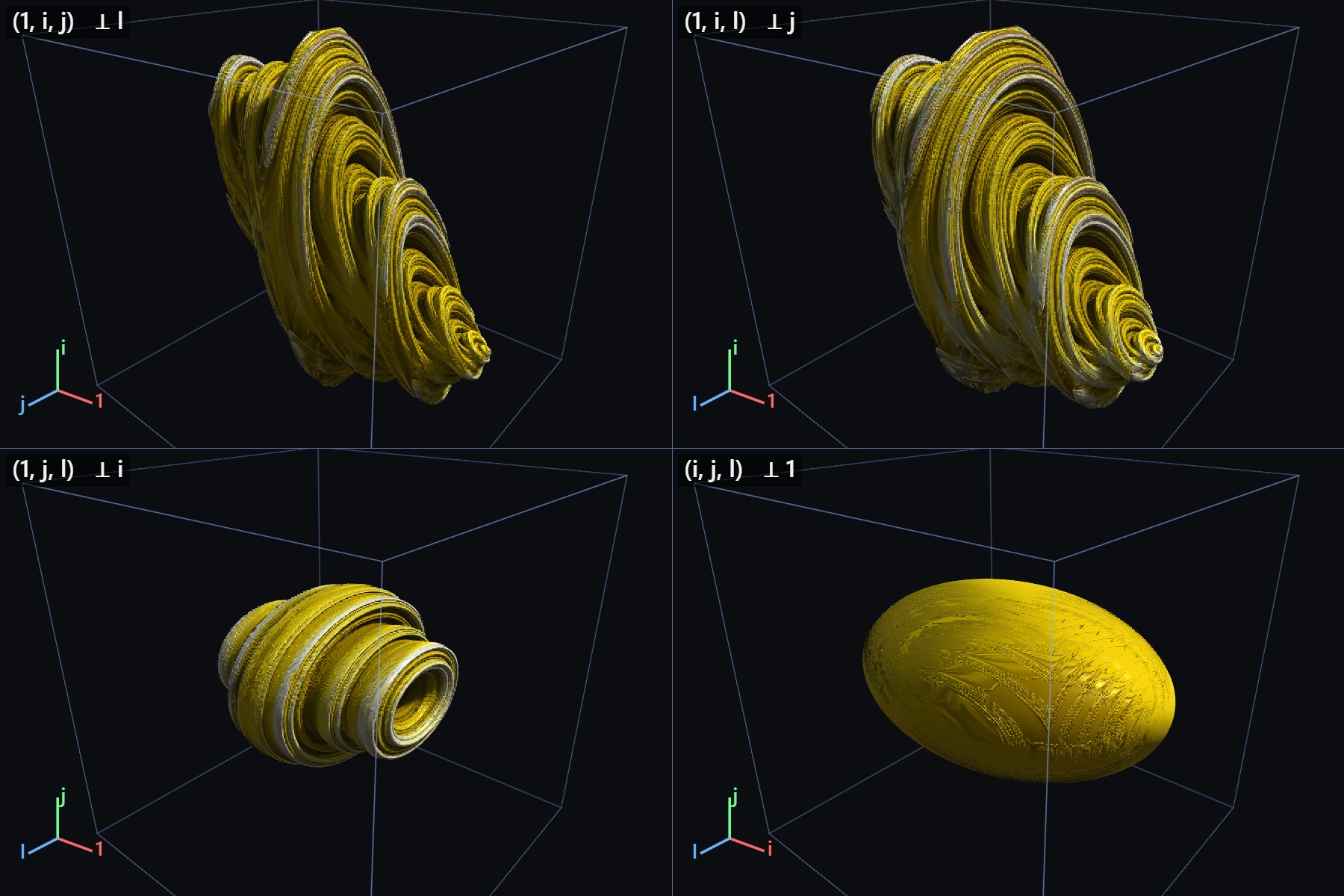}\par\vspace{1mm}
\noindent{\footnotesize \textbf{Map} $z\mapsto z\diamond z+c$ $\cdot$ \textbf{Algebra} custom: Cayley--Dickson, $\gamma=(-0.5,-2,-1)$ $\cdot$ $c=(-0.315,\,\allowbreak 0.845,\,\allowbreak 0,\,\allowbreak 0,\,\allowbreak 0,\,\allowbreak 0,\,\allowbreak 0,\,\allowbreak 0)$ ($=0.97\,s^*$, $s^*=0.9299$) $\cdot$ \textbf{Frame} $X,Y,Z,N=1,i,j,l$ $\cdot$ half-extent 1.4 $\cdot$ 20 iterations $\cdot$ \textbf{Coloring} 8D orbit trap (point), gradient Gold}\par\vspace{1mm}
\noindent{\footnotesize\begin{tabular}{@{}lllll@{}}\toprule Slice (X, Y, Z) & Body [95\%] & Exact ($J\geq0.9995$) & Approximate ($0.97\leq J<0.9995$) & Dim\\\midrule
($1, i, j$) $\perp l$ & 21.7\% {\scriptsize[21.5, 22.0]} & $M$: z; $R_z$: 180; $P$ & -- & 3 \\
($1, i, l$) $\perp j$ & 30.8\% {\scriptsize[30.5, 31.1]} & $M$: z; $R_z$: 180; $P$ & -- & 3 \\
($1, j, l$) $\perp i$ & 13.9\% {\scriptsize[13.7, 14.1]} & $M$: xyz; $R_x$: 180; $R_y$: 180; $R_z$: 180; $P$ & -- & 3 \\
($i, j, l$) $\perp 1$ & 30.5\% {\scriptsize[30.2, 30.7]} & $M$: xyz; $R_x$: 180; $R_y$: 180; $R_z$: 180; $P$ & -- & 3 \\
\bottomrule\end{tabular}}\par\vspace{2mm}
The anisotropic octonion algebra $\gamma=(-0.5,-2,-1)$ is isomorphic to the octonions: by a substitution analogous to that of Proposition~\ref{atl:prop:rescaling}, rescaling $i$ and $j$ turns it into $\OO(-1,-1,-1)$. Its Julia sets are therefore linear images of octonionic Julia sets, stretched along $i$ and compressed along $j$. The imaginary quarter shows it most clearly: an ellipsoid ($31\%$) instead of a ball, since the norm $0.5x_1^2+2x_2^2+x_4^2$ has semi-axes in the ratio $\sqrt2:1/\sqrt2:1$ along $i$, $j$, $l$; the picture shows a body elongated along $i$ and flattened along $j$. The quarters through $1$ are spiral shells stretched accordingly.

The two upper quarters are no longer identical, as they were in Plate~A01 ($22\%$ and $31\%$), because $j$ and $l$ are now scaled differently. The measurements report mirrors, half-turns and $P$ but no revolution: the circles of the octonionic bodies of revolution have become ellipses. In the quarter $(1,j,l)$, a body of revolution about $1$ in the octonions, and in the imaginary quarter, a ball there, only the three mirrors, the three half-turns and $P$ remain exact. The orbits lie in a subalgebra of dimension $3$, as in the octonions.
\clearpage

\subsubsection{P04. A split Petersson algebra at \texorpdfstring{$60^\circ$}{60\textdegree{}}}\label{atl:pl:P04}\mbox{}\par\vspace{0.5mm}
\noindent\includegraphics[width=\textwidth]{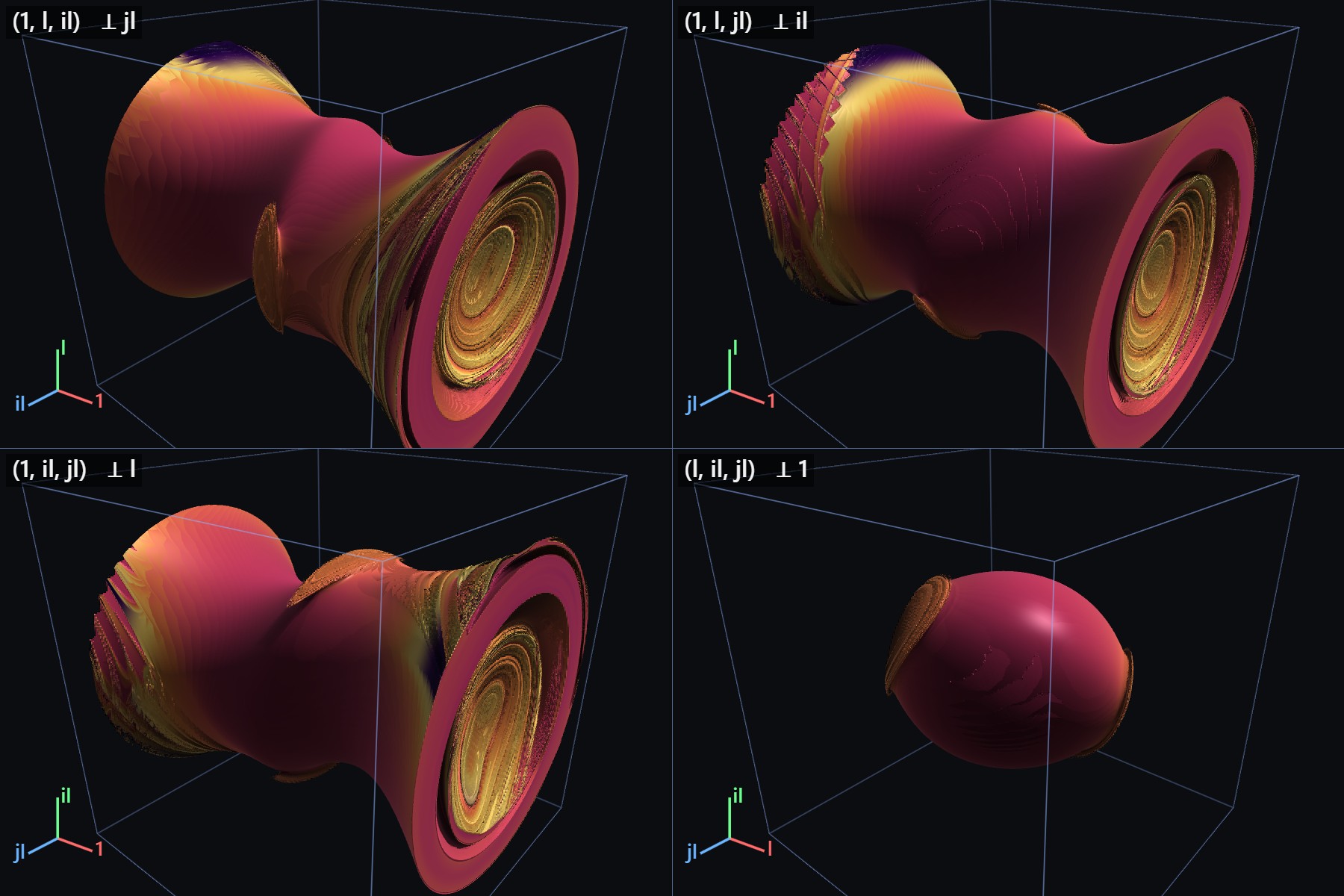}\par\vspace{1mm}
\noindent{\footnotesize \textbf{Map} $z\mapsto z\diamond z+c$ $\cdot$ \textbf{Algebra} custom: Petersson, $\gamma=(-1,-1,1)$, $\theta=60^\circ$ $\cdot$ $c=(-0.235,\,\allowbreak -0.013,\,\allowbreak 0.029,\,\allowbreak 0.01,\,\allowbreak 0.002,\,\allowbreak 0.094,\,\allowbreak -0.067,\,\allowbreak 0.079)$ ($=0.95\,s^*$, $s^*=0.2899$) $\cdot$ \textbf{Frame} $X,Y,Z,N=1,l,il,jl$ $\cdot$ half-extent 1.4 $\cdot$ 20 iterations $\cdot$ \textbf{Coloring} 8D orbit trap (point), gradient Sunset}\par\vspace{1mm}
\noindent{\footnotesize\begin{tabular}{@{}lllll@{}}\toprule Slice (X, Y, Z) & Body [95\%] & Exact ($J\geq0.9995$) & Approximate ($0.97\leq J<0.9995$) & Dim\\\midrule
($1, l, il$) $\perp jl$ & 42.7\% {\scriptsize[42.4, 43.0]} & $P$ & -- & 8 \\
($1, l, jl$) $\perp il$ & 39.9\% {\scriptsize[39.5, 40.2]} & $P$ & -- & 8 \\
($1, il, jl$) $\perp l$ & 41.0\% {\scriptsize[40.7, 41.3]} & $P$ & -- & 8 \\
($l, il, jl$) $\perp 1$ & 17.9\% {\scriptsize[17.7, 18.1]} & $P$ & -- & 8 \\
\bottomrule\end{tabular}}\par\vspace{2mm}
The split Petersson algebra at $\theta=60^\circ$ with the generic constant, in the frame $F_B$. The three quarters through $1$ show large hourglass-shaped bodies ($43\%$, $40\%$ and $41\%$ of the ball): a bulb at one end, a narrow waist, and a wide flared bell at the other end, cut by the faces of the cube, whose flat cross-sections show concentric rings. The imaginary quarter $(l,il,jl)$ holds a small egg ($18\%$) with thin, ragged caps at its two ends. The flared bodies are the signature of the split norm, which has signature $(1,2)$ on the sections through $1$, while on $(l,il,jl)$ it is definite and the body is compact.

The only exact symmetry is $P$ (Remark~\ref{atl:rem:minus}), and the orbits lie in a subalgebra of dimension $8$. The comparison with Plate~B08 ($\theta=120^\circ$, the split Okubo algebra) shows bodies of the same kind, differing in detail (the hourglasses of Plate~B08 have a second collar near the waist): nothing in these pictures singles out the Okubo angle.
\clearpage

\subsubsection{P05. Okubo in the imaginary frame \texorpdfstring{$F_D$}{FD}}\label{atl:pl:P05}\mbox{}\par\vspace{0.5mm}
\noindent\includegraphics[width=\textwidth]{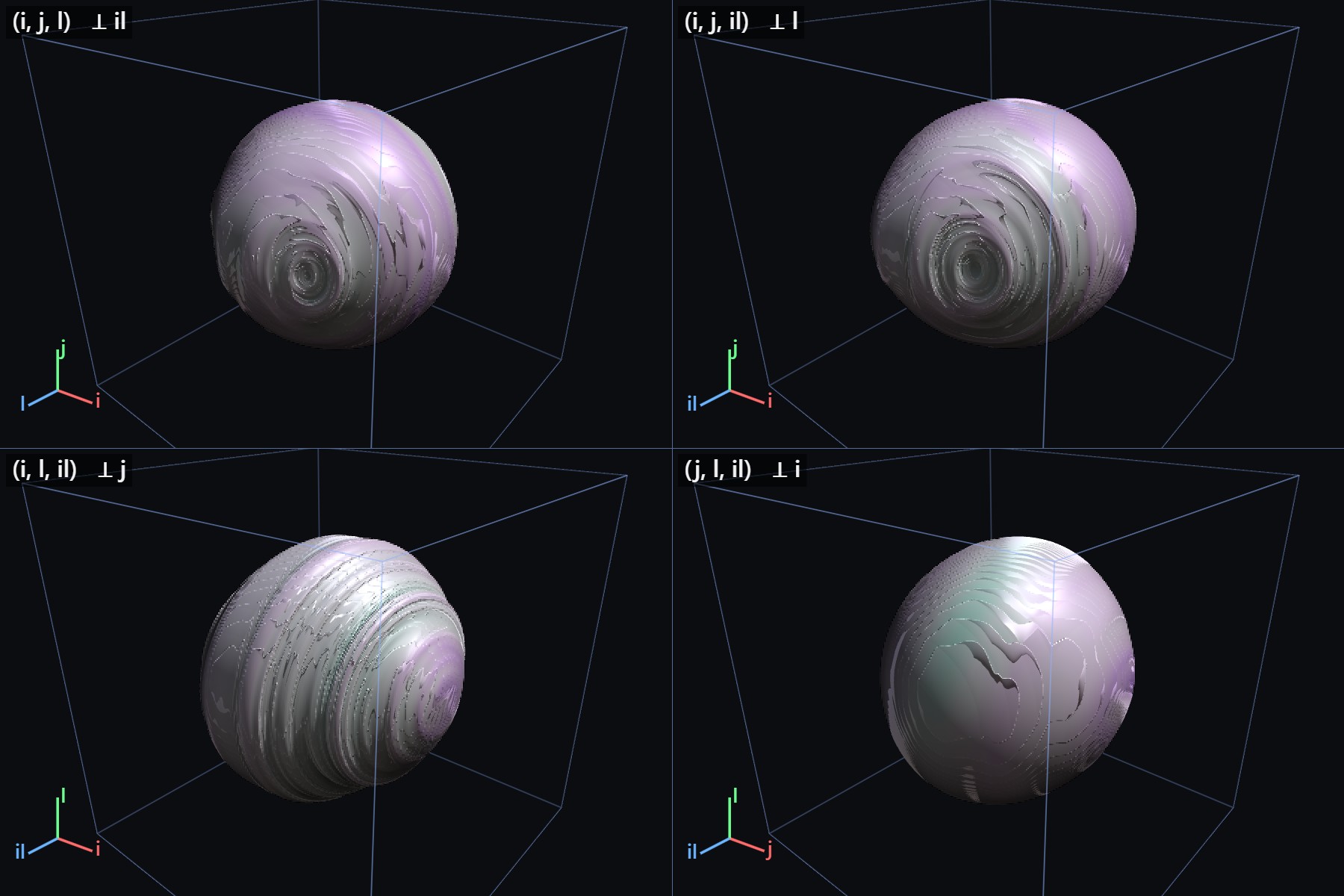}\par\vspace{1mm}
\noindent{\footnotesize \textbf{Map} $z\mapsto z\diamond z+c$ $\cdot$ \textbf{Algebra} Okubo algebra (Petersson, $\gamma=(-1,-1,-1)$, $\theta=120^\circ$) $\cdot$ $c=(-0.397,\,\allowbreak -0.021,\,\allowbreak 0.049,\,\allowbreak 0.016,\,\allowbreak 0.004,\,\allowbreak 0.16,\,\allowbreak -0.113,\,\allowbreak 0.133)$ ($=0.95\,s^*$, $s^*=0.4900$) $\cdot$ \textbf{Frame} $X,Y,Z,N=i,j,l,il$ $\cdot$ half-extent 1.4 $\cdot$ 20 iterations $\cdot$ \textbf{Coloring} 8D orbit trap (point), gradient Aurora}\par\vspace{1mm}
\noindent{\footnotesize\begin{tabular}{@{}lllll@{}}\toprule Slice (X, Y, Z) & Body [95\%] & Exact ($J\geq0.9995$) & Approximate ($0.97\leq J<0.9995$) & Dim\\\midrule
($i, j, l$) $\perp il$ & 33.4\% {\scriptsize[33.1, 33.7]} & $P$ & -- & 8 \\
($i, j, il$) $\perp l$ & 33.5\% {\scriptsize[33.2, 33.8]} & $P$ & -- & 8 \\
($i, l, il$) $\perp j$ & 34.9\% {\scriptsize[34.6, 35.2]} & $P$ & $M$: y; $R_y$: 180 (min.\ 0.972) & 8 \\
($j, l, il$) $\perp i$ & 37.8\% {\scriptsize[37.5, 38.1]} & $P$ & $R_x$: 37 (min.\ 0.977) & 8 \\
\bottomrule\end{tabular}}\par\vspace{2mm}
The imaginary frame $F_D=(i,j,l\,|\,il)$ contains no direction of the unit. With the generic constant, the Okubo algebra gives four round bodies of nearly the same size ($33\%$ to $38\%$ of the ball): the quarters $(i,j,l)$ and $(i,j,il)$ each show a single vortex on the visible face, the quarter $(i,l,il)$ is covered by ribs ending in a whorl, and $(j,l,il)$ is smoother, with a few folds. In the octonions all four would be balls (Proposition~\ref{atl:prop:imaginary}); here the identity $z^2=-\n(z)$ for imaginary $z$ has no counterpart, and the bodies are only nearly round.

The only exact symmetry is $P$ (Remark~\ref{atl:rem:minus}); the mirror $M_y$ of $(i,l,il)$ ($J=0.972$) and the rotation $R_x(37^\circ)$ of $(j,l,il)$ ($J=0.977$) are approximate, and the orbits lie in a subalgebra of dimension $8$. It is worth noticing that the quarters $(i,l,il)$ and $(j,l,il)$ contain the root plane $V_2=\langle l,il\rangle$: for a constant in $\CC$ they would be exact bodies of revolution about their first axis (Proposition~\ref{atl:prop:torus}). For the generic constant the stabilizer is a conjugate torus whose root planes are in general position, and the revolutions are not visible in the coordinate frame (Plate~R09).
\clearpage

\subsubsection{P06. Split Okubo in the imaginary frame \texorpdfstring{$F_D$}{FD}}\label{atl:pl:P06}\mbox{}\par\vspace{0.5mm}
\noindent\includegraphics[width=\textwidth]{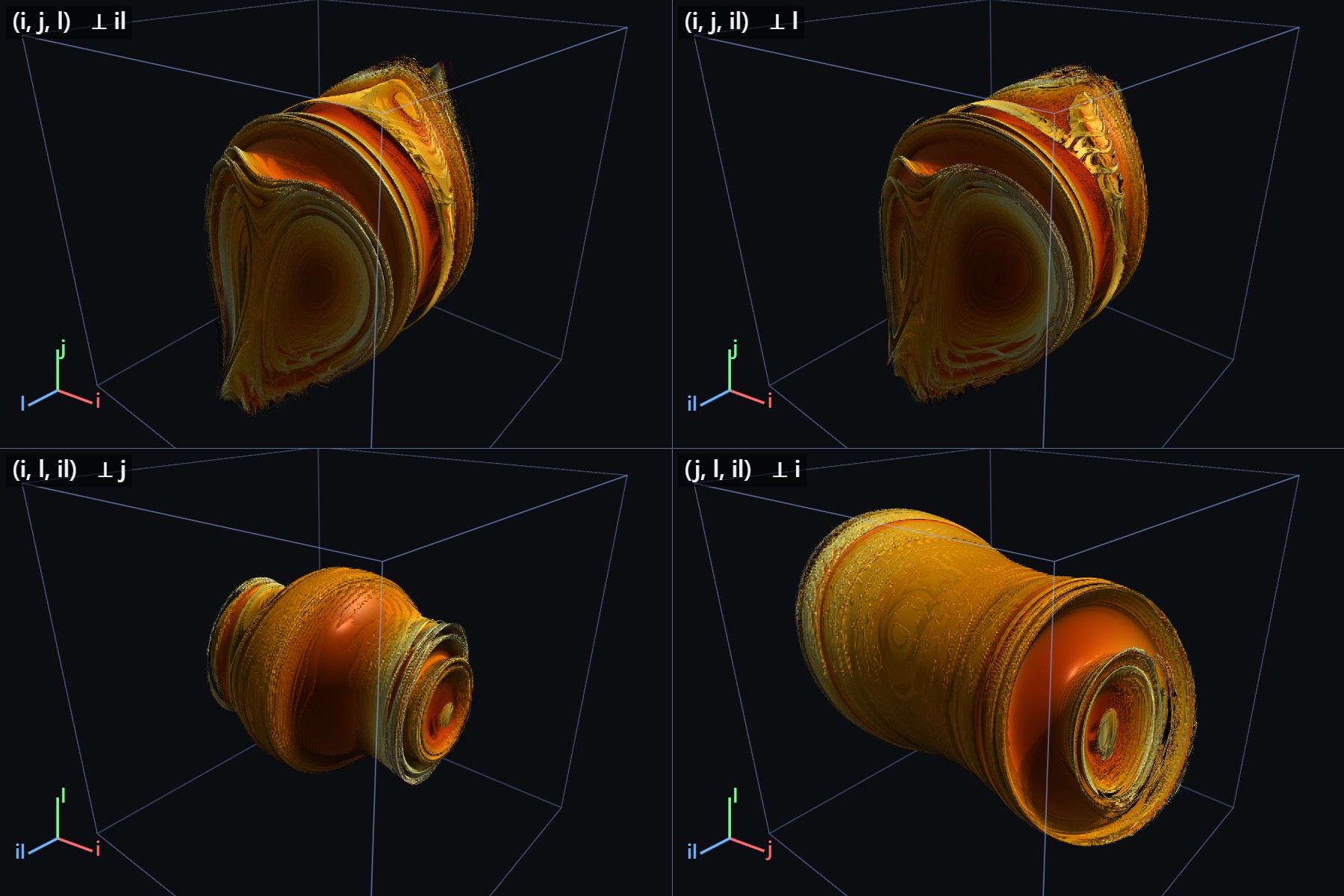}\par\vspace{1mm}
\noindent{\footnotesize \textbf{Map} $z\mapsto z\diamond z+c$ $\cdot$ \textbf{Algebra} Split Okubo algebra (Petersson, $\gamma=(-1,-1,1)$, $\theta=120^\circ$) $\cdot$ $c=(-0.433,\,\allowbreak -0.023,\,\allowbreak 0.054,\,\allowbreak 0.018,\,\allowbreak 0.004,\,\allowbreak 0.174,\,\allowbreak -0.123,\,\allowbreak 0.145)$ ($=0.95\,s^*$, $s^*=0.5339$) $\cdot$ \textbf{Frame} $X,Y,Z,N=i,j,l,il$ $\cdot$ half-extent 1.4 $\cdot$ 20 iterations $\cdot$ \textbf{Coloring} 8D orbit trap (point), gradient Fire}\par\vspace{1mm}
\noindent{\footnotesize\begin{tabular}{@{}lllll@{}}\toprule Slice (X, Y, Z) & Body [95\%] & Exact ($J\geq0.9995$) & Approximate ($0.97\leq J<0.9995$) & Dim\\\midrule
($i, j, l$) $\perp il$ & 29.3\% {\scriptsize[29.0, 29.5]} & $P$ & -- & 8 \\
($i, j, il$) $\perp l$ & 28.5\% {\scriptsize[28.3, 28.8]} & $P$ & -- & 8 \\
($i, l, il$) $\perp j$ & 18.9\% {\scriptsize[18.7, 19.2]} & $P$ & -- & 8 \\
($j, l, il$) $\perp i$ & 34.2\% {\scriptsize[33.9, 34.5]} & $P$ & -- & 8 \\
\bottomrule\end{tabular}}\par\vspace{2mm}
In the split Okubo algebra the imaginary frame shows four quite different bodies. The quarters $(i,j,l)$ and $(i,j,il)$ are tilted shells ($29\%$ and $28\%$ of the ball) whose broad face carries a large spiral, with ragged fringes at the edges; the quarter $(i,l,il)$ is a small jar ($19\%$) with ringed, open ends; the quarter $(j,l,il)$ is a larger barrel-shaped jar ($34\%$) whose mouth shows concentric rings.

The split norm has signature $(2,1)$ on the first two sections and $(1,2)$ on the last two, and the bodies are accordingly open along some directions, although none of them is a simple hyperboloid. The only exact symmetry is $P$ (Remark~\ref{atl:rem:minus}), with no approximate one, and the orbits lie in a subalgebra of dimension $8$ in every quarter. As for all split algebras, the Euclidean bailout of the renderer is not an escape criterion, and the rendered bodies are approximations.
\clearpage

\subsection{Part R. Tests of the predictions}\label{atl:part:R}
The plates of this part were rendered \emph{after} the statements of the main text were formulated, in order to test them; each one carries a prediction, which the image and the measurements confirm or refute. Unless stated otherwise the constant is $c_A=-0.086+0.231\,i$, the constant of Plates~A07 and C03, whose argument $110.4^\circ$ lies at $9.6^\circ$ from the Weyl wall at $120^\circ$. Plates~R01--R04 test the root planes of the Okubo algebra, with the octonionic term of comparison; R05--R08 follow the constant towards the Weyl wall at fixed modulus; R09 is the generic constant of Plate~B07 in the frame adapted to its stabilizer; R10--R12 test the subspace $W$ at $\theta=90^\circ$ and its controls; R13--R16 complete the Petersson path in the frame $F_C$; R17 repeats Plate~J01 with $20$ iterations, as a null test of Kaplansky's trick.

\subsubsection{R01. The displaced revolution: the root plane \texorpdfstring{$\langle l,il\rangle$}{<l,il>} of the Okubo algebra}\label{atl:pl:R01}\mbox{}\par\vspace{0.5mm}
\noindent\includegraphics[width=\textwidth]{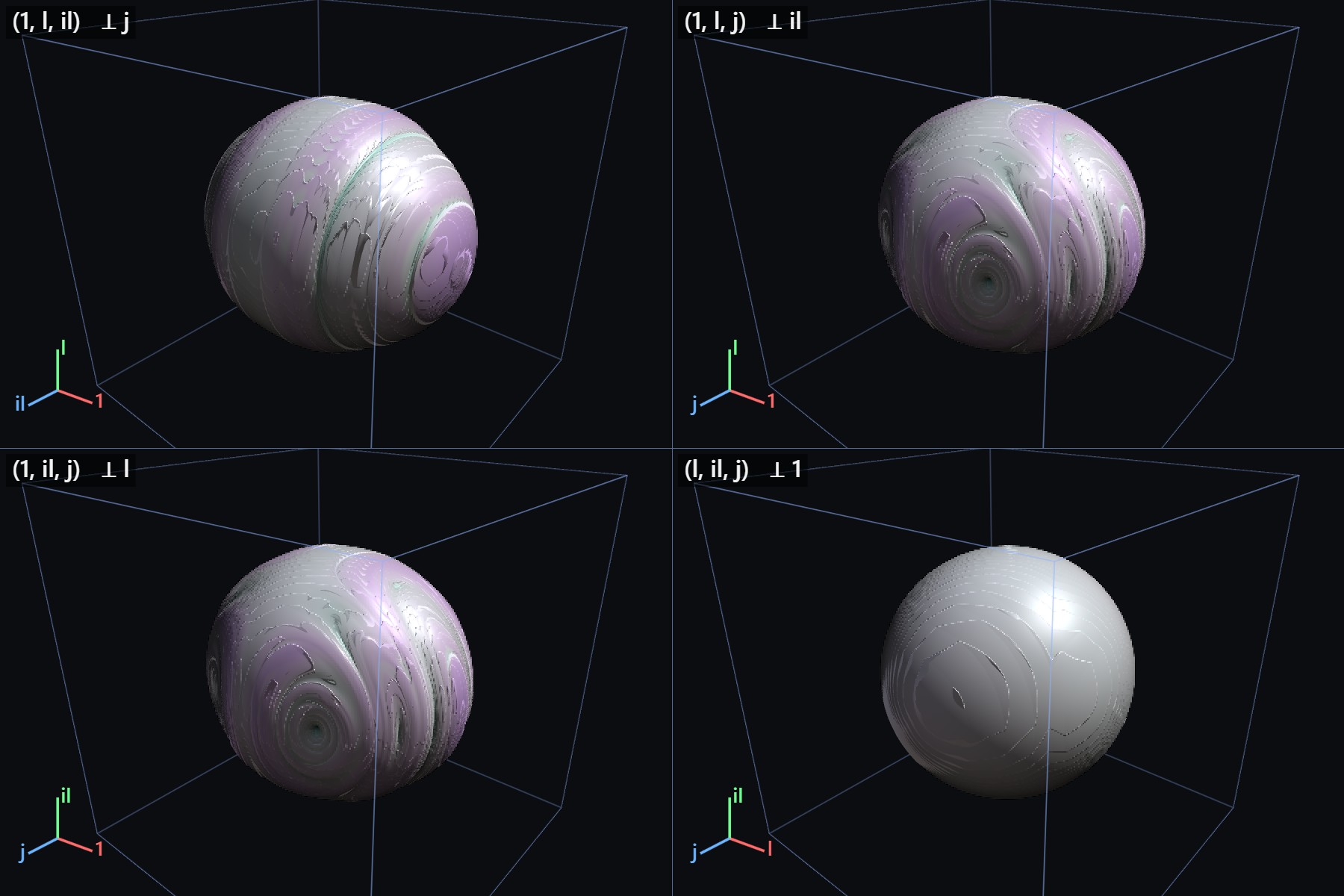}\par\vspace{1mm}
\noindent{\footnotesize \textbf{Map} $z\mapsto z\diamond z+c$ $\cdot$ \textbf{Algebra} Okubo algebra (Petersson, $\gamma=(-1,-1,-1)$, $\theta=120^\circ$) $\cdot$ $c=(-0.086,\,\allowbreak 0.231,\,\allowbreak 0,\,\allowbreak 0,\,\allowbreak 0,\,\allowbreak 0,\,\allowbreak 0,\,\allowbreak 0)$ $\cdot$ \textbf{Frame} $X,Y,Z,N=1,l,il,j$ $\cdot$ half-extent 1.4 $\cdot$ 20 iterations $\cdot$ \textbf{Coloring} 8D orbit trap (point), gradient Aurora}\par\vspace{1mm}
\noindent{\footnotesize\setlength{\tabcolsep}{2.7pt}\begin{tabular}{@{}lllll@{}}\toprule Slice (X, Y, Z) & Body [95\%] & Exact ($J\geq0.9995$) & Approximate ($0.97\leq J<0.9995$) & Dim\\\midrule
($1, l, il$) $\perp j$ & 37.9\% {\scriptsize[37.6, 38.2]} & $M$: xyz; $R_x$: all; $R_y$: 180; $R_z$: 180; $P$ & -- & 3 \\
($1, l, j$) $\perp il$ & 37.3\% {\scriptsize[37.0, 37.6]} & $M$: xyz; $R_x$: 180; $R_y$: 180; $R_z$: 180; $P$ & -- & 5 \\
($1, il, j$) $\perp l$ & 37.3\% {\scriptsize[37.0, 37.6]} & $M$: xyz; $R_x$: 180; $R_y$: 180; $R_z$: 180; $P$ & -- & 5 \\
($l, il, j$) $\perp 1$ & 34.6\% {\scriptsize[34.3, 34.9]} & $M$: xyz; $R_x$: 180; $R_y$: 180; $R_z$: all; $P$ & $R_x$: 120,90,72,60,37; $R_y$: 120,90,72,60,37 (min.\ 0.981) & 5 \\
\bottomrule\end{tabular}}\par\vspace{2mm}
This plate was rendered to test Theorem~A. The constant is $c_A$, the constant of Plates~A07 and C03, and the frame $(1,l,il\,|\,j)$ is chosen so that two of its quarters contain the root plane $V_2=\langle l,il\rangle$, which the torus of Proposition~\ref{atl:prop:torus} rotates by the angle $\beta-\alpha$ \eqref{atl:eq:rootangles}. The prediction is precise: the upper-left quarter $(1,l,il)$ is a body of revolution about $1$; the lower-right quarter $(l,il,j)$ is a body of revolution about $j$, since the circle $\alpha=0$ of the torus rotates $V_2$ and fixes $j$; the two mixed quarters $(1,l,j)$ and $(1,il,j)$ are not bodies of revolution, but they are congruent, being exchanged by the element of that circle that turns $V_2$ by a right angle.

The image agrees at a glance. The upper-left body carries ridges that run as circles about the axis $1$, closing into a set of concentric rings where that axis leaves the body on the right. The two mixed quarters show the same rounded body, with the same vortex near the bottom, and the lower-right quarter is an almost perfect sphere with a smooth surface.

The measurements confirm the prediction without exception. In $(1,l,il)$ the rotations about $1$ are exact by every tested angle, including $37^\circ$, and nothing else is; in $(l,il,j)$ the rotations about $j$ are exact by every angle, while those about $l$ and $il$ are only approximate ($J\approx0.98$--$0.99$): the lower-right quarter is nearly, but not exactly, a ball. The mixed quarters have the same body share, $37.3\%$, and only the finite symmetries: the three mirrors, the three half-turns and the point reflection, which by Remark~\ref{atl:rem:minus} come in pairs. The rotations about $1$ in the mixed quarters stop at $J=0.959$. The orbits lie in a subalgebra of dimension $3$ in the upper-left quarter and of dimension $5$ in the other three.

It is worth noticing what this plate changes in the reading of Plates~A07 and C03: their frames contain no root plane, and this is the only reason why they show no revolution. The revolutions of the Okubo algebra are not lost; they are displaced to the root planes. The octonionic term of comparison, with the same constant and the same frame, is Plate~R04.
\clearpage

\subsubsection{R02. The root plane \texorpdfstring{$\langle j,k\rangle$}{<j,k>} of the Okubo algebra}\label{atl:pl:R02}\mbox{}\par\vspace{0.5mm}
\noindent\includegraphics[width=\textwidth]{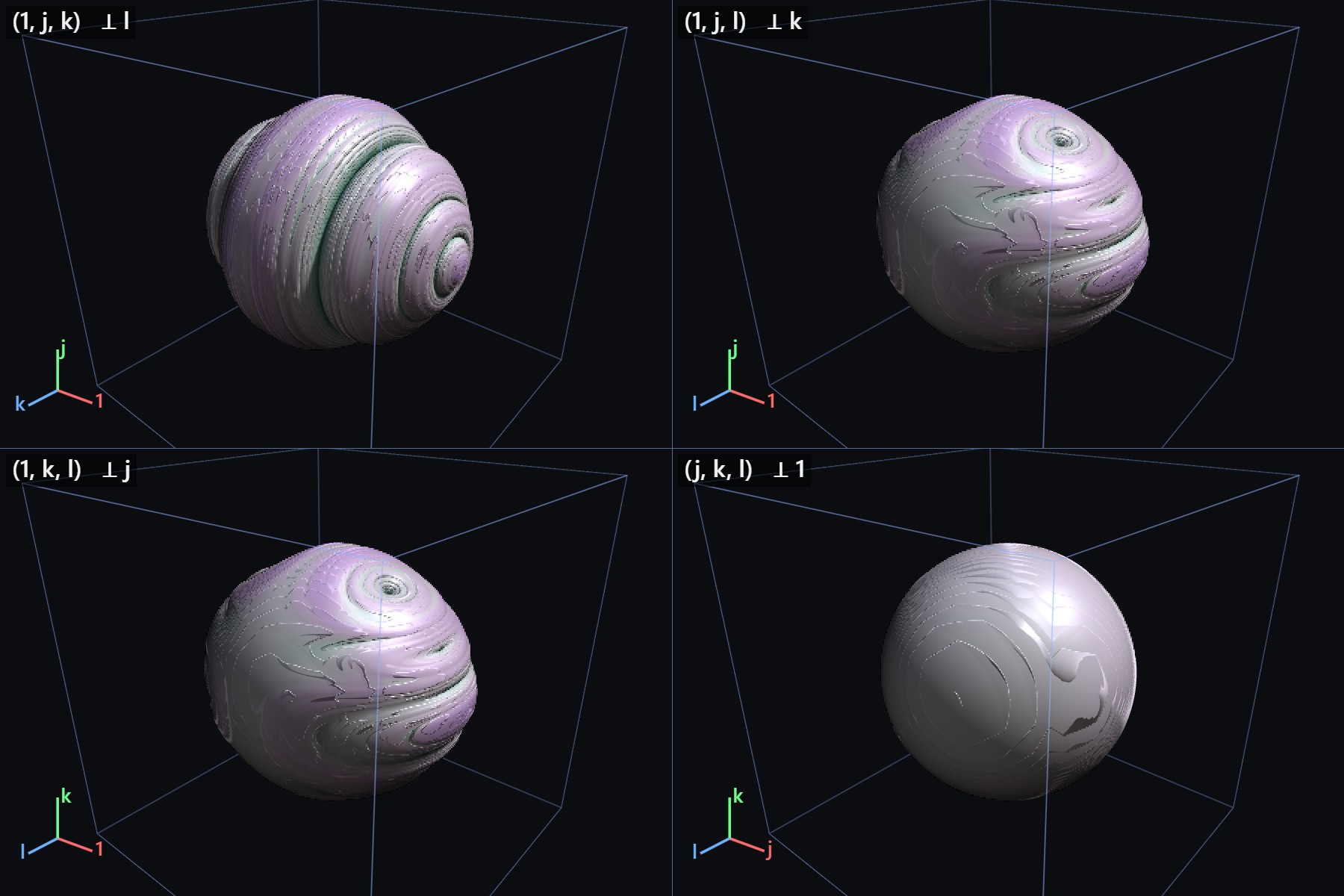}\par\vspace{1mm}
\noindent{\footnotesize \textbf{Map} $z\mapsto z\diamond z+c$ $\cdot$ \textbf{Algebra} Okubo algebra (Petersson, $\gamma=(-1,-1,-1)$, $\theta=120^\circ$) $\cdot$ $c=(-0.086,\,\allowbreak 0.231,\,\allowbreak 0,\,\allowbreak 0,\,\allowbreak 0,\,\allowbreak 0,\,\allowbreak 0,\,\allowbreak 0)$ $\cdot$ \textbf{Frame} $X,Y,Z,N=1,j,k,l$ $\cdot$ half-extent 1.4 $\cdot$ 20 iterations $\cdot$ \textbf{Coloring} 8D orbit trap (point), gradient Aurora}\par\vspace{1mm}
\noindent{\footnotesize\setlength{\tabcolsep}{2.7pt}\begin{tabular}{@{}lllll@{}}\toprule Slice (X, Y, Z) & Body [95\%] & Exact ($J\geq0.9995$) & Approximate ($0.97\leq J<0.9995$) & Dim\\\midrule
($1, j, k$) $\perp l$ & 34.7\% {\scriptsize[34.4, 35.0]} & $M$: xyz; $R_x$: all; $R_y$: 180; $R_z$: 180; $P$ & -- & 3 \\
($1, j, l$) $\perp k$ & 37.3\% {\scriptsize[37.0, 37.6]} & $M$: xyz; $R_x$: 180; $R_y$: 180; $R_z$: 180; $P$ & -- & 5 \\
($1, k, l$) $\perp j$ & 37.3\% {\scriptsize[37.0, 37.6]} & $M$: xyz; $R_x$: 180; $R_y$: 180; $R_z$: 180; $P$ & -- & 5 \\
($j, k, l$) $\perp 1$ & 35.1\% {\scriptsize[34.8, 35.4]} & $M$: xyz; $R_x$: 180; $R_y$: 180; $R_z$: all; $P$ & $R_x$: 120,90,72,60,37; $R_y$: 120,90,72,60,37 (min.\ 0.982) & 5 \\
\bottomrule\end{tabular}}\par\vspace{2mm}
The second root plane, $V_1=\langle j,k\rangle$, rotated by the angle $2\alpha$ \eqref{atl:eq:rootangles}, in the frame $(1,j,k\,|\,l)$ with $c=c_A$. The prediction is the same as for Plate~R01: the upper-left quarter $(1,j,k)$ is a body of revolution about $1$; the lower-right quarter $(j,k,l)$ is a body of revolution about $l$, since the circle $\beta=\alpha$ of the torus rotates $V_1$ and fixes $V_2$ pointwise; the mixed quarters $(1,j,l)$ and $(1,k,l)$ are congruent, and not of revolution.

The upper-left body is the most visibly ``turned'' of the plate: its ridges are circles about the axis $1$, with a small concentric cap on the right. The two mixed quarters show the same body, with a vortex near the top, and the lower-right quarter is again a near-sphere, with a small notch on its surface.

The measurements confirm the prediction. The rotations about $1$ in $(1,j,k)$ and about $l$ in $(j,k,l)$ are exact by every tested angle; in $(j,k,l)$ the rotations about $j$ and $k$ are approximate ($J\approx0.98$--$0.99$), so that this quarter too is nearly, but not exactly, a ball. The mixed quarters have the same body share, $37.3\%$, the same values of $J$ and only the finite symmetries (mirrors, half-turns, point reflection). The body shares are $34.7\%$ in $(1,j,k)$ and $35.1\%$ in $(j,k,l)$, the dimensions of the subalgebras $3$ and $5$.

The profile of the upper-left body is not that of Plate~R01, and its ridges are more pronounced. The three root planes see the constant $c_A$ under different angles, and the torus rotates them at different speeds; nothing forces the three bodies of revolution to coincide. The quarter $(1,j,k)$ lies in the quaternions, where the Okubo product is the para-quaternionic one (Proposition~\ref{atl:prop:torus}): it is also a section of the para-octonionic Julia set with the same constant.
\clearpage

\subsubsection{R03. The root plane \texorpdfstring{$\langle jl,kl\rangle$}{<jl,kl>} of the Okubo algebra}\label{atl:pl:R03}\mbox{}\par\vspace{0.5mm}
\noindent\includegraphics[width=\textwidth]{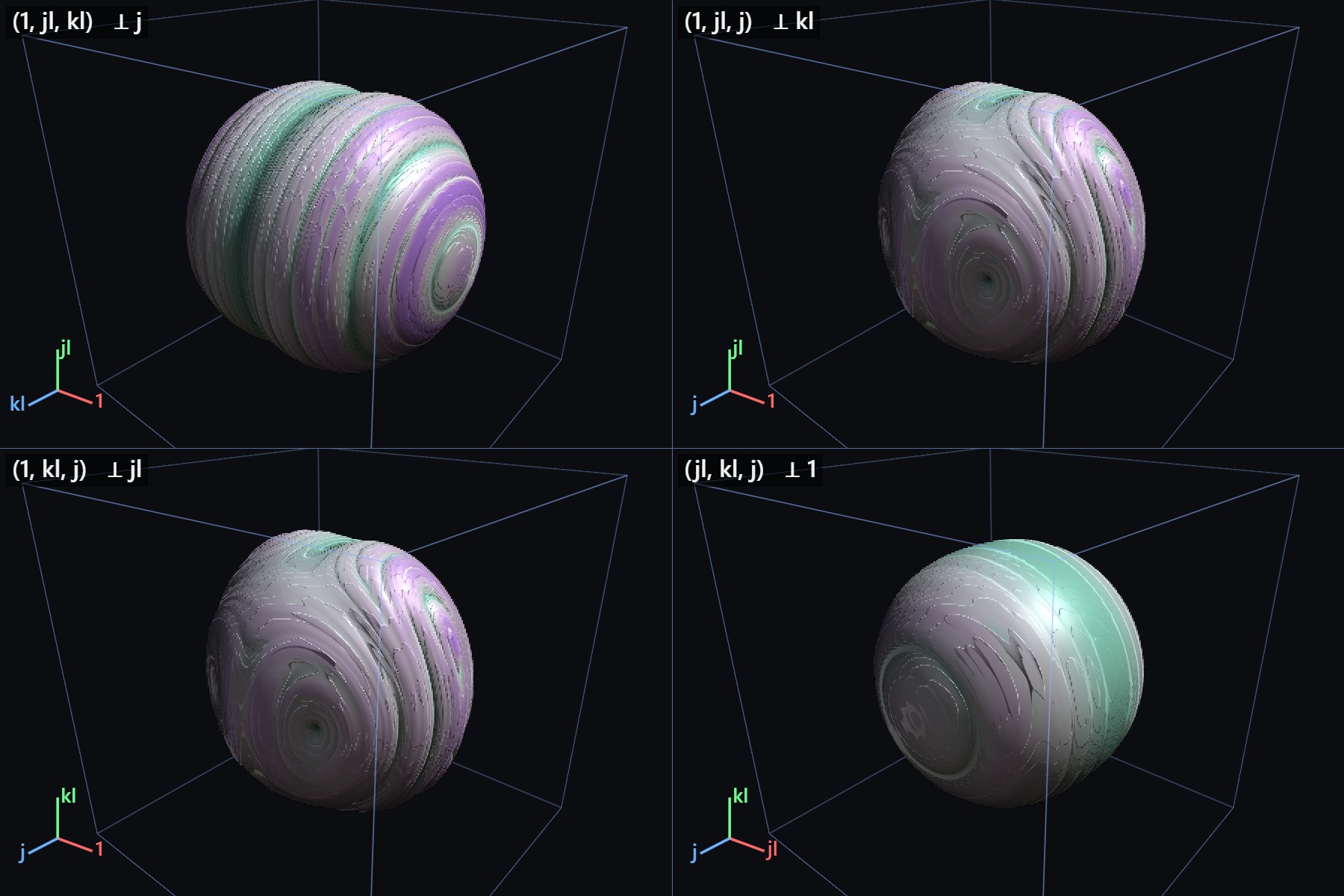}\par\vspace{1mm}
\noindent{\footnotesize \textbf{Map} $z\mapsto z\diamond z+c$ $\cdot$ \textbf{Algebra} Okubo algebra (Petersson, $\gamma=(-1,-1,-1)$, $\theta=120^\circ$) $\cdot$ $c=(-0.086,\,\allowbreak 0.231,\,\allowbreak 0,\,\allowbreak 0,\,\allowbreak 0,\,\allowbreak 0,\,\allowbreak 0,\,\allowbreak 0)$ $\cdot$ \textbf{Frame} $X,Y,Z,N=1,jl,kl,j$ $\cdot$ half-extent 1.4 $\cdot$ 20 iterations $\cdot$ \textbf{Coloring} 8D orbit trap (point), gradient Aurora}\par\vspace{1mm}
\noindent{\footnotesize\begin{tabular}{@{}lllll@{}}\toprule Slice (X, Y, Z) & Body [95\%] & Exact ($J\geq0.9995$) & Approximate ($0.97\leq J<0.9995$) & Dim\\\midrule
($1, jl, kl$) $\perp j$ & 50.6\% {\scriptsize[50.3, 50.9]} & $M$: xyz; $R_x$: all; $R_y$: 180; $R_z$: 180; $P$ & -- & 3 \\
($1, jl, j$) $\perp kl$ & 42.2\% {\scriptsize[41.9, 42.5]} & $M$: xyz; $R_x$: 180; $R_y$: 180; $R_z$: 180; $P$ & -- & 5 \\
($1, kl, j$) $\perp jl$ & 42.2\% {\scriptsize[41.9, 42.5]} & $M$: xyz; $R_x$: 180; $R_y$: 180; $R_z$: 180; $P$ & -- & 5 \\
($jl, kl, j$) $\perp 1$ & 39.3\% {\scriptsize[39.0, 39.6]} & $M$: xyz; $R_x$: 180; $R_y$: 180; $R_z$: all; $P$ & -- & 5 \\
\bottomrule\end{tabular}}\par\vspace{2mm}
The third root plane, $V_3=\langle jl,kl\rangle$, rotated by the angle $\alpha+\beta$ \eqref{atl:eq:rootangles}, in the frame $(1,jl,kl\,|\,j)$ with $c=c_A$. The prediction: the upper-left quarter $(1,jl,kl)$ is a body of revolution about $1$; the lower-right quarter $(jl,kl,j)$ is a body of revolution about $j$, fixed by the circle $\alpha=0$, which rotates $V_3$; the mixed quarters $(1,jl,j)$ and $(1,kl,j)$ are congruent and not of revolution.

The upper-left body is the largest of the three root-plane bodies, a ball-like body girdled by bands that run as circles about the axis $1$ and close into concentric rings on the right. The mixed quarters show the same body, with a vortex at the bottom, and the lower-right quarter is a round body whose surface pattern is centred on the front-left, where the axis $j$ points.

The measurements confirm the prediction: the rotations about $1$ in $(1,jl,kl)$ and about $j$ in $(jl,kl,j)$ are exact by every tested angle. Unlike Plates~R01 and R02, the lower-right quarter is not even approximately a ball: the rotations about $jl$ and $kl$ stay at $J\le0.95$. The mixed quarters have the same body share, $42.2\%$, and only the finite symmetries; the upper-left quarter occupies $50.6\%$ of the ball, against $37.9\%$ and $34.7\%$ for the corresponding quarters of Plates~R01 and R02.

The three root planes thus give three different bodies of revolution. It is worth noticing that the difference can be read off the dynamics as well, as the main text shows: along a bounded orbit the radii in $V_1$, $V_2$, $V_3$, which the torus leaves invariant, decay at three different rates.
\clearpage

\subsubsection{R04. Octonions in the frame \texorpdfstring{$(1,l,il\,|\,j)$}{(1,l,il | j)}: the term of comparison}\label{atl:pl:R04}\mbox{}\par\vspace{0.5mm}
\noindent\includegraphics[width=\textwidth]{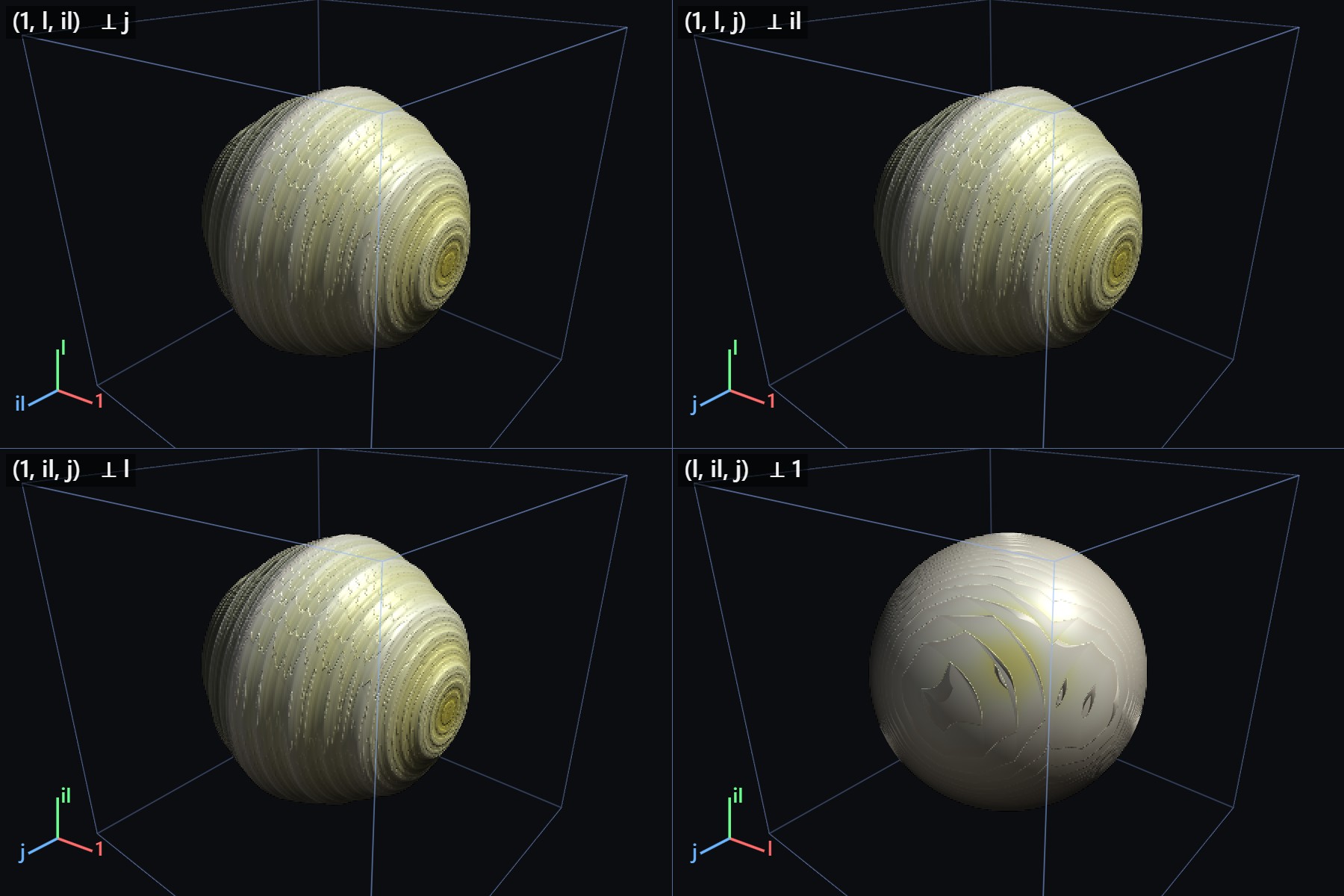}\par\vspace{1mm}
\noindent{\footnotesize \textbf{Map} $z\mapsto z\diamond z+c$ $\cdot$ \textbf{Algebra} Octonions (Cayley--Dickson, $\gamma=(-1,-1,-1)$) $\cdot$ $c=(-0.086,\,\allowbreak 0.231,\,\allowbreak 0,\,\allowbreak 0,\,\allowbreak 0,\,\allowbreak 0,\,\allowbreak 0,\,\allowbreak 0)$ $\cdot$ \textbf{Frame} $X,Y,Z,N=1,l,il,j$ $\cdot$ half-extent 1.4 $\cdot$ 20 iterations $\cdot$ \textbf{Coloring} 8D orbit trap (point), gradient Gold}\par\vspace{1mm}
\noindent{\footnotesize\begin{tabular}{@{}lllll@{}}\toprule Slice (X, Y, Z) & Body [95\%] & Exact ($J\geq0.9995$) & Approximate ($0.97\leq J<0.9995$) & Dim\\\midrule
($1, l, il$) $\perp j$ & 39.5\% {\scriptsize[39.2, 39.8]} & $M$: xyz; $R_x$: all; $R_y$: 180; $R_z$: 180; $P$ & -- & 3 \\
($1, l, j$) $\perp il$ & 39.5\% {\scriptsize[39.2, 39.8]} & $M$: xyz; $R_x$: all; $R_y$: 180; $R_z$: 180; $P$ & -- & 3 \\
($1, il, j$) $\perp l$ & 39.5\% {\scriptsize[39.2, 39.8]} & $M$: xyz; $R_x$: all; $R_y$: 180; $R_z$: 180; $P$ & -- & 3 \\
($l, il, j$) $\perp 1$ & 44.7\% {\scriptsize[44.4, 45.0]} & $M$: xyz; $R_x$: all; $R_y$: all; $R_z$: all; $P$ & -- & 3 \\
\bottomrule\end{tabular}}\par\vspace{2mm}
The octonions, with the same constant $c_A$ and the same frame $(1,l,il\,|\,j)$ as Plate~R01. Here the prediction comes from Proposition~\ref{atl:prop:equivariance}: the quadratic map is equivariant under the group $O(6)$ of the imaginary directions orthogonal to $c$, which acts transitively on their orthonormal pairs. Hence the three quarters through $1$ are the same body, a body of revolution about $1$, and the purely imaginary quarter $(l,il,j)$ is a ball (Proposition~\ref{atl:prop:imaginary}).

The image shows exactly this. The three quarters through $1$ are indistinguishable: the same ribbed body, with ridges running as circles about the axis $1$ and closing into a concentric pattern on the right. The lower-right quarter is a sphere; the markings on its surface come from the colouring.

The measurements confirm it: the three quarters through $1$ have the same body share, $39.5\%$, and their rotations about $1$ are exact by every tested angle; in $(l,il,j)$ every rotation about every axis is exact. The orbits lie in a subalgebra of dimension $3$ in all four quarters (Proposition~\ref{atl:prop:span}).

The comparison with Plate~R01 is the point of this plate. The upper-left quarters of the two plates are both bodies of revolution about $1$, but the octonionic one is shared by every section through $1$, whereas the Okubonic one needs the root plane: in the Okubo algebra the mixed quarters lose their revolution, and the imaginary quarter keeps only the revolution about $j$ and is no longer a ball.
\clearpage

\subsubsection{R05. Towards the Weyl wall: \texorpdfstring{$\arg c=100^\circ$}{arg c=100\textdegree{}}}\label{atl:pl:R05}\mbox{}\par\vspace{0.5mm}
\noindent\includegraphics[width=\textwidth]{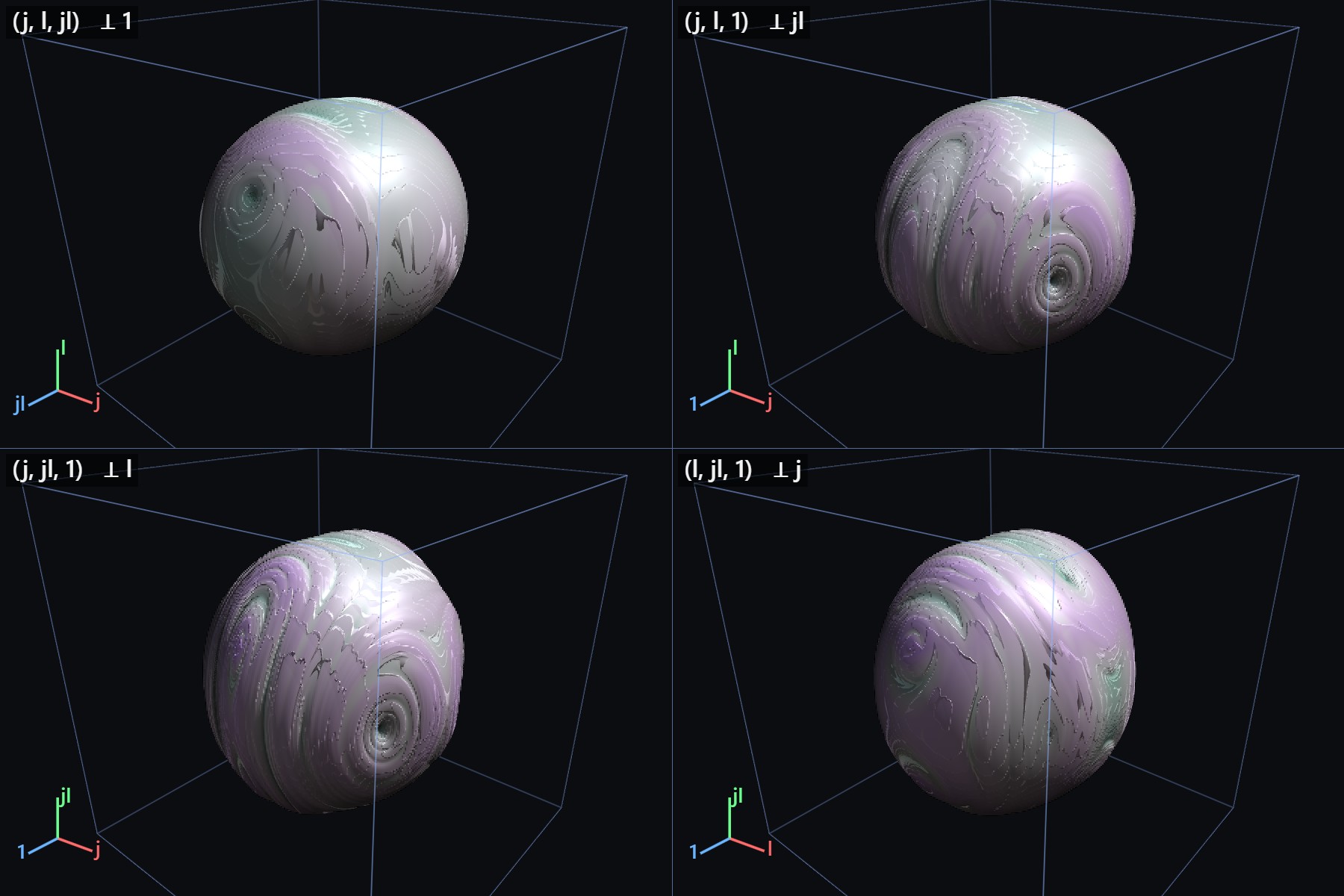}\par\vspace{1mm}
\noindent{\footnotesize \textbf{Map} $z\mapsto z\diamond z+c$ $\cdot$ \textbf{Algebra} Okubo algebra (Petersson, $\gamma=(-1,-1,-1)$, $\theta=120^\circ$) $\cdot$ $c=(-0.043,\,\allowbreak 0.243,\,\allowbreak 0,\,\allowbreak 0,\,\allowbreak 0,\,\allowbreak 0,\,\allowbreak 0,\,\allowbreak 0)$ $\cdot$ \textbf{Frame} $X,Y,Z,N=j,l,jl,1$ $\cdot$ half-extent 1.4 $\cdot$ 20 iterations $\cdot$ \textbf{Coloring} 8D orbit trap (point), gradient Aurora}\par\vspace{1mm}
\noindent{\footnotesize\begin{tabular}{@{}lllll@{}}\toprule Slice (X, Y, Z) & Body [95\%] & Exact ($J\geq0.9995$) & Approximate ($0.97\leq J<0.9995$) & Dim\\\midrule
($j, l, jl$) $\perp 1$ & 38.2\% {\scriptsize[37.9, 38.5]} & $M$: xyz; $R_x$: 180; $R_y$: 180; $R_z$: 180; $P$ & $R_z$: 37 (min.\ 0.974) & 8 \\
($j, l, 1$) $\perp jl$ & 37.2\% {\scriptsize[36.9, 37.5]} & $M$: xyz; $R_x$: 180; $R_y$: 180; $R_z$: 180; $P$ & -- & 5 \\
($j, jl, 1$) $\perp l$ & 42.1\% {\scriptsize[41.8, 42.4]} & $M$: xyz; $R_x$: 180; $R_y$: 180; $R_z$: 180; $P$ & -- & 5 \\
($l, jl, 1$) $\perp j$ & 42.4\% {\scriptsize[42.1, 42.7]} & $M$: xyz; $R_x$: 180; $R_y$: 180; $R_z$: 180; $P$ & -- & 5 \\
\bottomrule\end{tabular}}\par\vspace{2mm}
The first of four plates that follow the constant towards a Weyl wall at fixed modulus $|c|=|c_A|$, with $\arg c=100^\circ$, $110.4^\circ$, $115^\circ$ and $120^\circ$ (Plates~R05--R08), in the frame $(j,l,jl\,|\,1)$. The quarter of interest is the upper-left one, $(j,l,jl)$, which mixes the three root planes and contains none of them. The prediction is that it is not a body of revolution about $jl$ off the wall, and that it becomes one on the wall $\arg c=120^\circ$, where the stabilizer of $c$ grows from the torus to $U(2)$ and a circle of $U(2)$ turns the plane $\langle j,l\rangle$ as a complex line of $V_1\oplus V_2$ while fixing $jl$.

At $100^\circ$, $20^\circ$ from the wall, the four quarters are rounded bodies with swirling surface patterns. In the upper-left quarter a single vortex appears at the upper left of the body, away from the point where the axis $jl$ comes towards the viewer; nothing in the picture suggests a revolution about $jl$.

The measurements agree. The quarter $(j,l,jl)$ has the finite symmetries only (mirrors, half-turns, point reflection); the rotation by $90^\circ$ about $jl$ scores $J=0.961$, and among the rotations about $jl$ only the one by $37^\circ$ reaches the approximate range, $J=0.974$. The other quarters have the finite symmetries as well, with body shares between $37\%$ and $42\%$. The orbits of the upper-left quarter span the whole algebra (dimension $8$), those of the other quarters a subalgebra of dimension $5$.
\clearpage

\subsubsection{R06. Towards the Weyl wall: \texorpdfstring{$c=c_A$}{c=cA}, \texorpdfstring{$\arg c=110.4^\circ$}{arg c=110.4\textdegree{}}}\label{atl:pl:R06}\mbox{}\par\vspace{0.5mm}
\noindent\includegraphics[width=\textwidth]{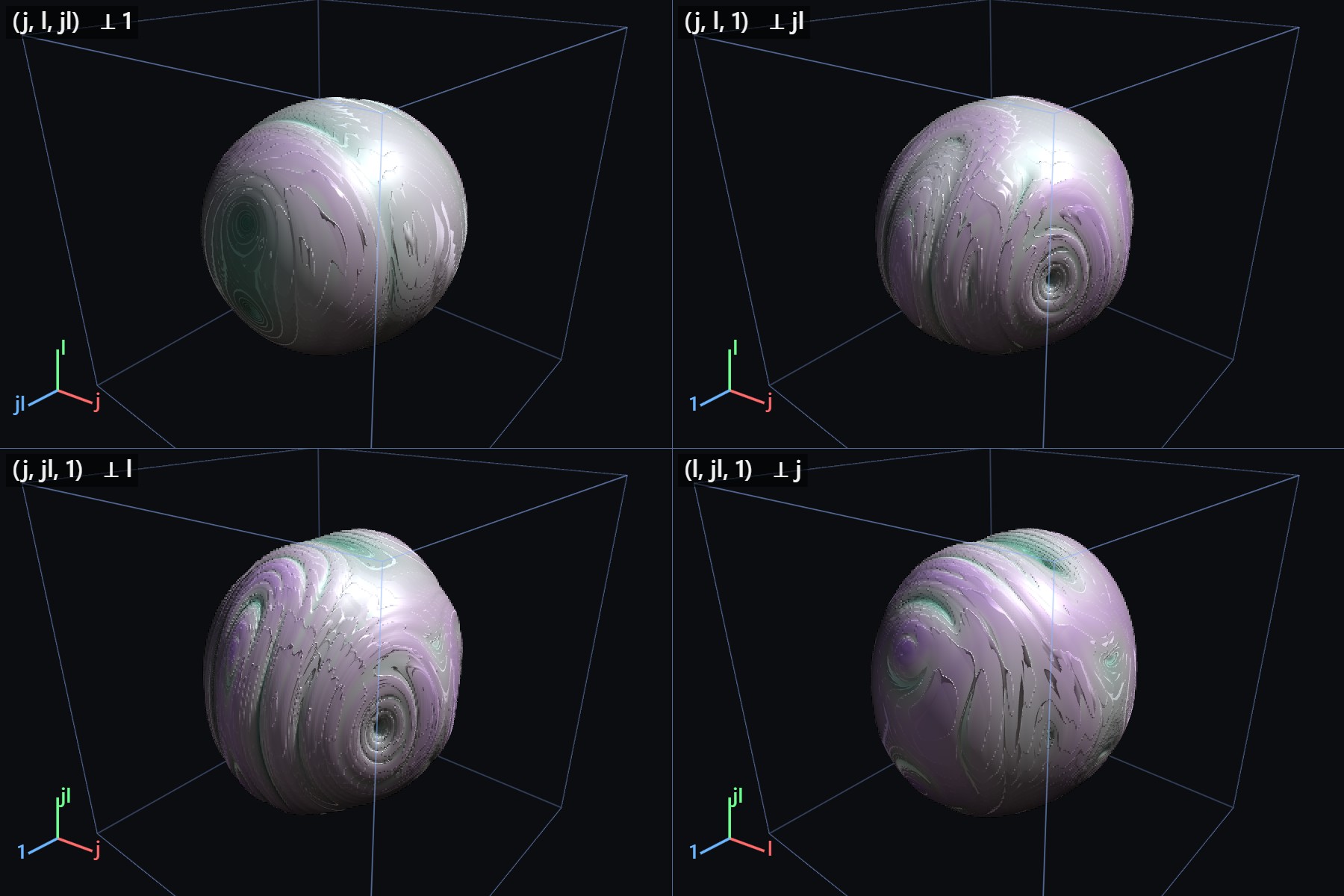}\par\vspace{1mm}
\noindent{\footnotesize \textbf{Map} $z\mapsto z\diamond z+c$ $\cdot$ \textbf{Algebra} Okubo algebra (Petersson, $\gamma=(-1,-1,-1)$, $\theta=120^\circ$) $\cdot$ $c=(-0.086,\,\allowbreak 0.231,\,\allowbreak 0,\,\allowbreak 0,\,\allowbreak 0,\,\allowbreak 0,\,\allowbreak 0,\,\allowbreak 0)$ $\cdot$ \textbf{Frame} $X,Y,Z,N=j,l,jl,1$ $\cdot$ half-extent 1.4 $\cdot$ 20 iterations $\cdot$ \textbf{Coloring} 8D orbit trap (point), gradient Aurora}\par\vspace{1mm}
\noindent{\footnotesize\begin{tabular}{@{}lllll@{}}\toprule Slice (X, Y, Z) & Body [95\%] & Exact ($J\geq0.9995$) & Approximate ($0.97\leq J<0.9995$) & Dim\\\midrule
($j, l, jl$) $\perp 1$ & 37.3\% {\scriptsize[37.0, 37.6]} & $M$: xyz; $R_x$: 180; $R_y$: 180; $R_z$: 180; $P$ & $R_z$: 120,90,72,60,37 (min.\ 0.972) & 8 \\
($j, l, 1$) $\perp jl$ & 37.3\% {\scriptsize[37.0, 37.6]} & $M$: xyz; $R_x$: 180; $R_y$: 180; $R_z$: 180; $P$ & -- & 5 \\
($j, jl, 1$) $\perp l$ & 42.1\% {\scriptsize[41.8, 42.4]} & $M$: xyz; $R_x$: 180; $R_y$: 180; $R_z$: 180; $P$ & -- & 5 \\
($l, jl, 1$) $\perp j$ & 44.0\% {\scriptsize[43.7, 44.3]} & $M$: xyz; $R_x$: 180; $R_y$: 180; $R_z$: 180; $P$ & -- & 5 \\
\bottomrule\end{tabular}}\par\vspace{2mm}
The constant $c_A$ of Plates~A07 and C03, whose argument $110.4^\circ$ lies $9.6^\circ$ from the Weyl wall at $120^\circ$, in the frame $(j,l,jl\,|\,1)$ of Plates~R05--R08. The first version of the atlas measured, at render time and on $2400$ points, a rotation about $jl$ in the quarter $(j,l,jl)$ with $J\approx0.98$, and the question was whether it was a symmetry. The prediction of Theorem~A is that it is not: $c_A$ is not on a wall, its stabilizer is the torus, and the quarter contains no root plane.

In the image the upper-left body shows two vortices, one above the other, on its left side; the pattern is more regular than in Plate~R05, but it is not centred on the axis $jl$.

The measurements settle the question. The rotations about $jl$ are all approximate, not exact: $J=0.972$ for $90^\circ$, between $0.973$ and $0.983$ for the other angles. The symmetry is therefore a near-symmetry, explained by the proximity of $c_A$ to the wall, and it grows as the constant approaches it (Plates~R07 and R08). The other quarters have only the finite symmetries; the upper-left quarter generates the whole algebra (dimension $8$), the others a subalgebra of dimension $5$.
\clearpage

\subsubsection{R07. Towards the Weyl wall: \texorpdfstring{$\arg c=115^\circ$}{arg c=115\textdegree{}}}\label{atl:pl:R07}\mbox{}\par\vspace{0.5mm}
\noindent\includegraphics[width=\textwidth]{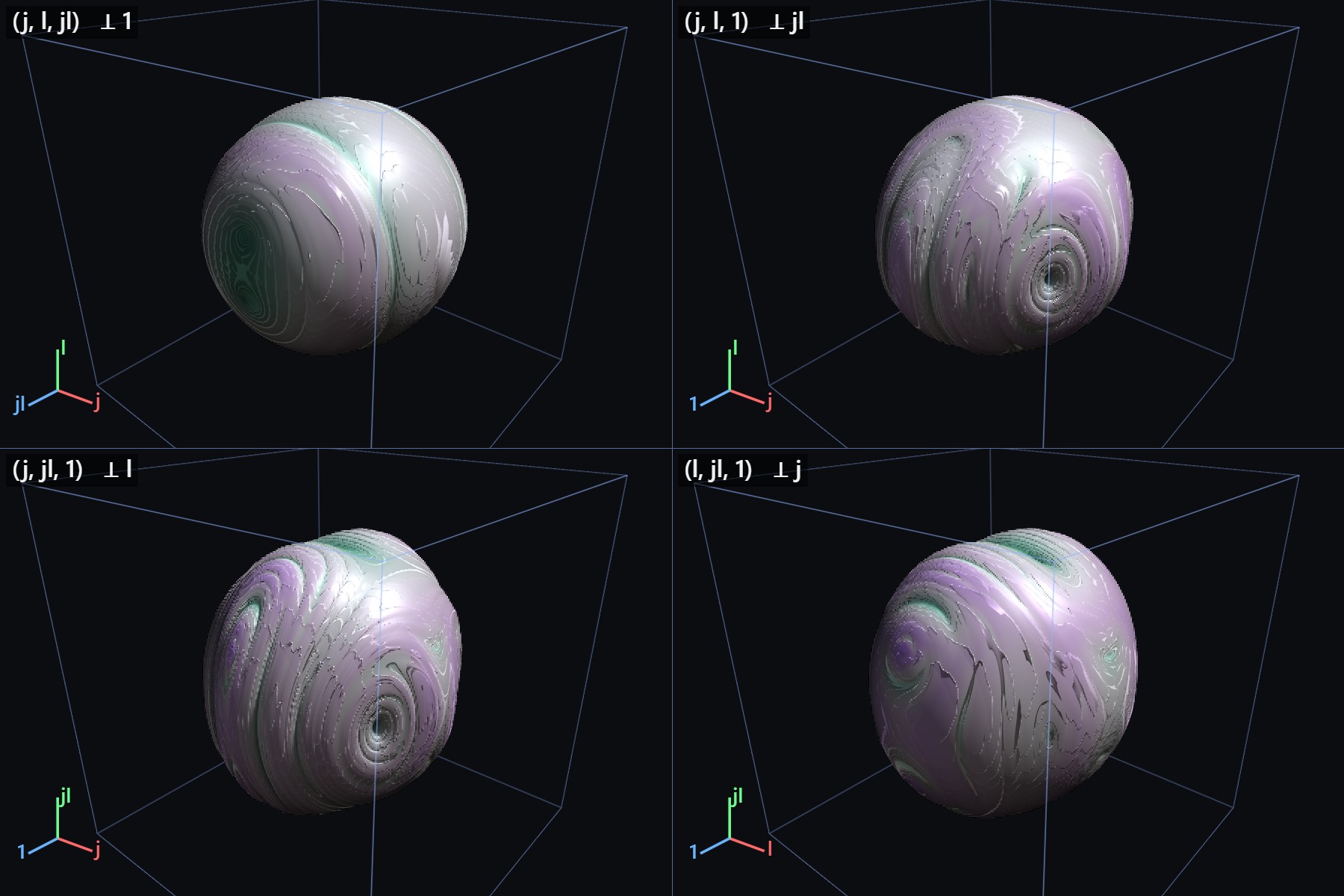}\par\vspace{1mm}
\noindent{\footnotesize \textbf{Map} $z\mapsto z\diamond z+c$ $\cdot$ \textbf{Algebra} Okubo algebra (Petersson, $\gamma=(-1,-1,-1)$, $\theta=120^\circ$) $\cdot$ $c=(-0.104,\,\allowbreak 0.223,\,\allowbreak 0,\,\allowbreak 0,\,\allowbreak 0,\,\allowbreak 0,\,\allowbreak 0,\,\allowbreak 0)$ $\cdot$ \textbf{Frame} $X,Y,Z,N=j,l,jl,1$ $\cdot$ half-extent 1.4 $\cdot$ 20 iterations $\cdot$ \textbf{Coloring} 8D orbit trap (point), gradient Aurora}\par\vspace{1mm}
\noindent{\footnotesize\begin{tabular}{@{}lllll@{}}\toprule Slice (X, Y, Z) & Body [95\%] & Exact ($J\geq0.9995$) & Approximate ($0.97\leq J<0.9995$) & Dim\\\midrule
($j, l, jl$) $\perp 1$ & 37.2\% {\scriptsize[36.9, 37.5]} & $M$: xyz; $R_x$: 180; $R_y$: 180; $R_z$: 180; $P$ & $R_z$: 120,90,72,60,37 (min.\ 0.982) & 8 \\
($j, l, 1$) $\perp jl$ & 37.6\% {\scriptsize[37.3, 37.9]} & $M$: xyz; $R_x$: 180; $R_y$: 180; $R_z$: 180; $P$ & -- & 5 \\
($j, jl, 1$) $\perp l$ & 42.2\% {\scriptsize[41.9, 42.5]} & $M$: xyz; $R_x$: 180; $R_y$: 180; $R_z$: 180; $P$ & -- & 5 \\
($l, jl, 1$) $\perp j$ & 44.7\% {\scriptsize[44.4, 45.0]} & $M$: xyz; $R_x$: 180; $R_y$: 180; $R_z$: 180; $P$ & -- & 5 \\
\bottomrule\end{tabular}}\par\vspace{2mm}
The constant at $\arg c=115^\circ$, $5^\circ$ from the Weyl wall, in the frame $(j,l,jl\,|\,1)$. The prediction is the same as for Plate~R06: no exact revolution about $jl$ in the quarter $(j,l,jl)$, but a near-symmetry that improves as the constant approaches the wall.

In the upper-left quarter the two vortices of Plate~R06 have drawn together into a single elongated system of closed curves on the left of the body, closer to the point where the axis $jl$ comes towards the viewer. The other three quarters are practically unchanged along the series.

The measurements confirm the trend: the rotations about $jl$ are approximate, with $J=0.982$ for $90^\circ$ and between $0.983$ and $0.988$ for the other angles, and none is exact. The index of the rotation by $90^\circ$ grows monotonically along the series, $0.961$, $0.972$, $0.982$ at $100^\circ$, $110.4^\circ$, $115^\circ$. The other quarters keep the finite symmetries only, with body shares between $37.6\%$ and $44.7\%$.
\clearpage

\subsubsection{R08. On the Weyl wall: \texorpdfstring{$\arg c=120^\circ$}{arg c=120\textdegree{}} and the birth of a revolution}\label{atl:pl:R08}\mbox{}\par\vspace{0.5mm}
\noindent\includegraphics[width=\textwidth]{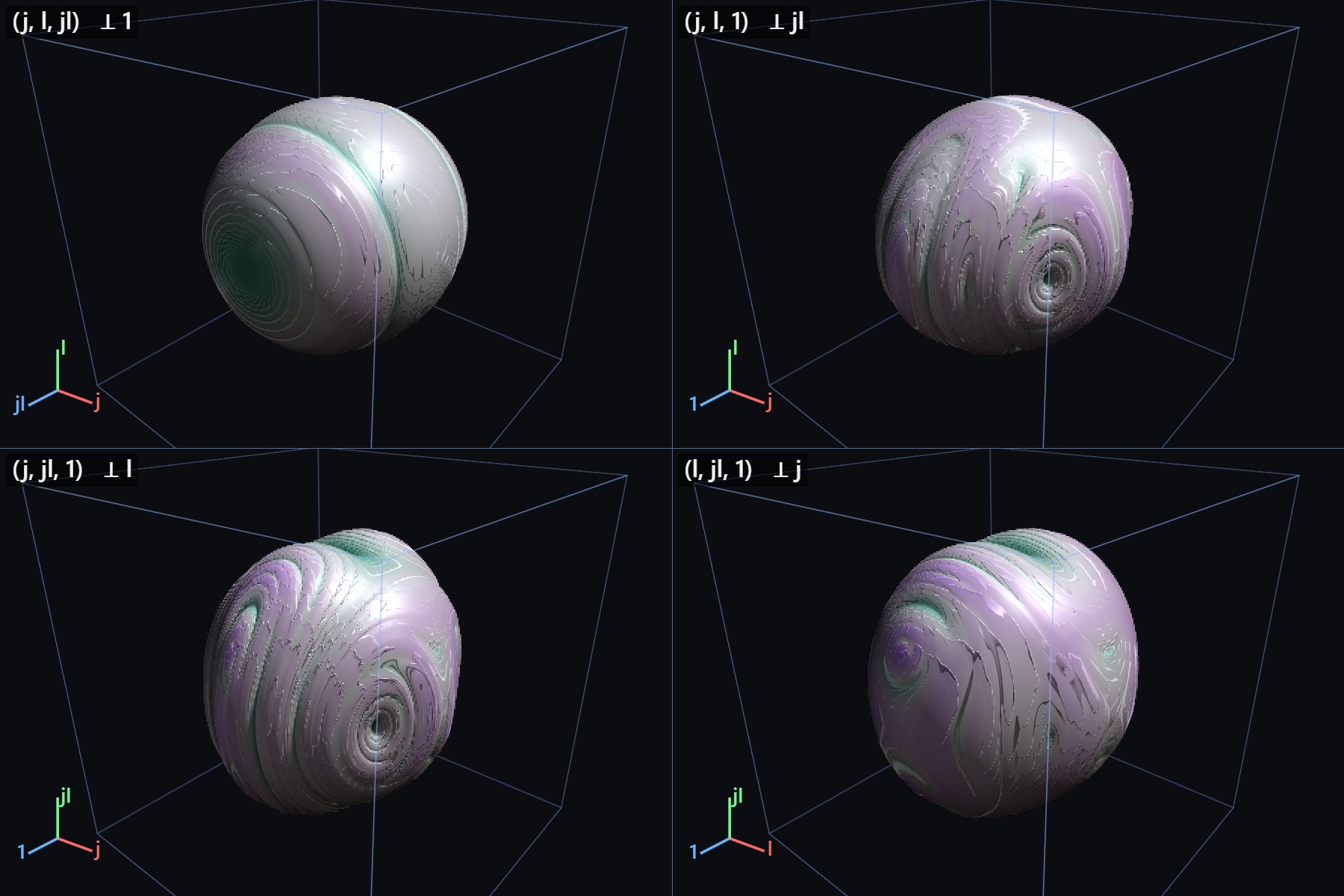}\par\vspace{1mm}
\noindent{\footnotesize \textbf{Map} $z\mapsto z\diamond z+c$ $\cdot$ \textbf{Algebra} Okubo algebra (Petersson, $\gamma=(-1,-1,-1)$, $\theta=120^\circ$) $\cdot$ $c=(-0.123,\,\allowbreak 0.213,\,\allowbreak 0,\,\allowbreak 0,\,\allowbreak 0,\,\allowbreak 0,\,\allowbreak 0,\,\allowbreak 0)$ $\cdot$ \textbf{Frame} $X,Y,Z,N=j,l,jl,1$ $\cdot$ half-extent 1.4 $\cdot$ 20 iterations $\cdot$ \textbf{Coloring} 8D orbit trap (point), gradient Aurora}\par\vspace{1mm}
\noindent{\footnotesize\begin{tabular}{@{}lllll@{}}\toprule Slice (X, Y, Z) & Body [95\%] & Exact ($J\geq0.9995$) & Approximate ($0.97\leq J<0.9995$) & Dim\\\midrule
($j, l, jl$) $\perp 1$ & 36.9\% {\scriptsize[36.6, 37.2]} & $M$: xyz; $R_x$: 180; $R_y$: 180; $R_z$: all; $P$ & -- & 8 \\
($j, l, 1$) $\perp jl$ & 37.7\% {\scriptsize[37.4, 38.0]} & $M$: xyz; $R_x$: 180; $R_y$: 180; $R_z$: 180; $P$ & -- & 5 \\
($j, jl, 1$) $\perp l$ & 42.0\% {\scriptsize[41.7, 42.3]} & $M$: xyz; $R_x$: 180; $R_y$: 180; $R_z$: 180; $P$ & -- & 5 \\
($l, jl, 1$) $\perp j$ & 45.5\% {\scriptsize[45.2, 45.8]} & $M$: xyz; $R_x$: 180; $R_y$: 180; $R_z$: 180; $P$ & -- & 5 \\
\bottomrule\end{tabular}}\par\vspace{2mm}
The constant on the Weyl wall, $\arg c=120^\circ$, with $|c|=|c_A|$, in the frame $(j,l,jl\,|\,1)$. Here the stabilizer of $c$ is $U(2)$, and the prediction is that the quarter $(j,l,jl)$, which contains no root plane, becomes an exact body of revolution about $jl$: the new circle of $U(2)$ turns $\langle j,l\rangle$ as a complex line of $V_1\oplus V_2$ and fixes $jl$.

The image shows it clearly. In the upper-left quarter the swirls of Plates~R05--R07 have become a single system of concentric circles, centred on the front-left of the body, where the axis $jl$ comes towards the viewer. The three other quarters look as in the previous plates.

The measurements confirm the prediction: in $(j,l,jl)$ the rotations about $jl$ are exact by every tested angle, including $37^\circ$, whereas the rotations about $j$ and $l$ stay below $0.94$. The series of Plates~R05--R08 is thus complete: the index of the rotation by $90^\circ$ about $jl$ is $0.961$, $0.972$, $0.982$ and finally $1$ on the wall. The other three quarters keep only the finite symmetries.

It is worth noticing that the growth of the index is gradual, while the symmetry itself appears only on the wall. The near-symmetries of $c_A$ observed in the first version were the shadow of a revolution that exists only $9.6^\circ$ away.
\clearpage

\subsubsection{R09. The generic constant of Plate~B07 in the frame of its stabilizer}\label{atl:pl:R09}\mbox{}\par\vspace{0.5mm}
\noindent\includegraphics[width=\textwidth]{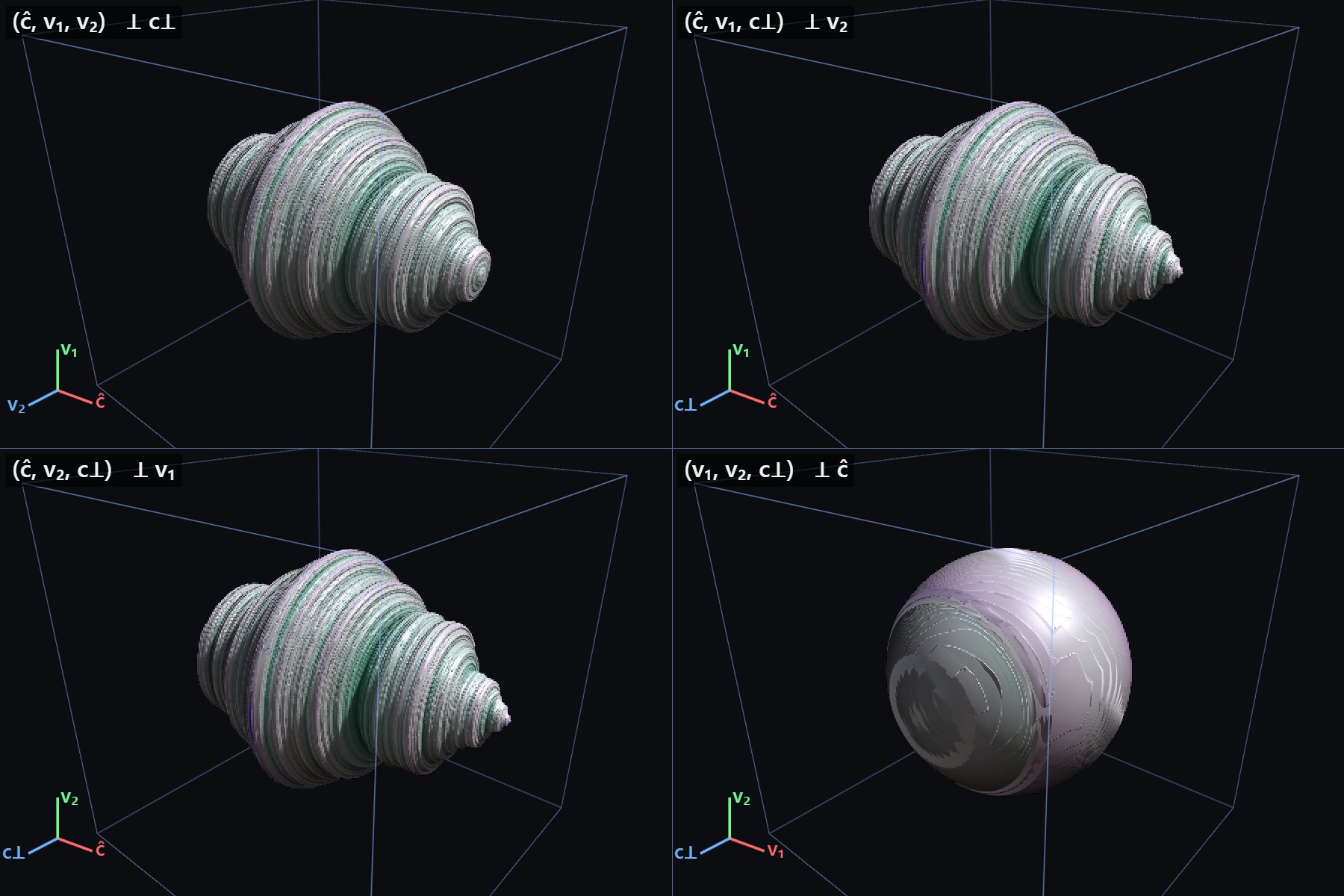}\par\vspace{1mm}
\noindent{\footnotesize \textbf{Map} $z\mapsto z\diamond z+c$ $\cdot$ \textbf{Algebra} Okubo algebra (Petersson, $\gamma=(-1,-1,-1)$, $\theta=120^\circ$) $\cdot$ $c=(-0.397,\,\allowbreak -0.021,\,\allowbreak 0.049,\,\allowbreak 0.016,\,\allowbreak 0.004,\,\allowbreak 0.16,\,\allowbreak -0.113,\,\allowbreak 0.133)$ $\cdot$ \textbf{Frame} $X,Y,Z,N=\hat c,v_1,v_2,c^{\perp}$ (explicit orthonormal frame, see the commentary) $\cdot$ half-extent 1.4 $\cdot$ 20 iterations $\cdot$ \textbf{Coloring} 8D orbit trap (point), gradient Aurora}\par\vspace{1mm}
\noindent{\footnotesize\setlength{\tabcolsep}{2.7pt}\begin{tabular}{@{}lllll@{}}\toprule Slice (X, Y, Z) & Body [95\%] & Exact ($J\geq0.9995$) & Approximate ($0.97\leq J<0.9995$) & Dim\\\midrule
($\hat c, v_1, v_2$) $\perp c^{\perp}$ & 26.2\% {\scriptsize[25.9, 26.5]} & $M$: xyz; $R_x$: all; $R_y$: 180; $R_z$: 180; $P$ & -- & 8 \\
($\hat c, v_1, c^{\perp}$) $\perp v_2$ & 27.1\% {\scriptsize[26.8, 27.4]} & $M$: y; $R_y$: 180; $P$ & -- & 8 \\
($\hat c, v_2, c^{\perp}$) $\perp v_1$ & 27.1\% {\scriptsize[26.8, 27.4]} & $M$: y; $R_y$: 180; $P$ & -- & 8 \\
($v_1, v_2, c^{\perp}$) $\perp \hat c$ & 31.4\% {\scriptsize[31.2, 31.7]} & $M$: xyz; $R_x$: 180; $R_y$: 180; $R_z$: all; $P$ & $R_x$: 120,90,72,60,37; $R_y$: 120,90,72,60,37 (min.\ 0.973) & 3 \\
\bottomrule\end{tabular}}\par\vspace{2mm}
The generic constant of Plate~B07 over the Okubo algebra, in an orthonormal frame computed from its stabilizer. For a generic $c\in\RR^8$ the stabilizer of $c$ is a two-dimensional torus $T_c$, conjugate to the torus of Proposition~\ref{atl:prop:torus}; it fixes a plane that contains $c$, and it rotates three root planes in general position with respect to the coordinate axes. The frame of this plate is $(\hat c,v_1,v_2\,|\,c^\perp)$, where $\hat c=c/|c|$, $\langle v_1,v_2\rangle$ is the root plane of largest weight of $T_c$, and $c^\perp$ spans the fixed plane together with $\hat c$. Theorem~A predicts two exact revolutions: $(\hat c,v_1,v_2)$ about $\hat c$ and $(v_1,v_2,c^\perp)$ about $c^\perp$; the mixed quarters $(\hat c,v_1,c^\perp)$ and $(\hat c,v_2,c^\perp)$ should be congruent and not of revolution.

The upper-left quarter is a ribbed body made of stacked rings along the axis $\hat c$, tapering towards the right. The two mixed quarters show the same body as each other, a ribbed body with a sharper tip; at first sight they could be taken for bodies of revolution as well. The lower-right quarter is a near-sphere with a shallow depression on its front-left side, where the axis $c^\perp$ points.

The measurements confirm the prediction. In $(\hat c,v_1,v_2)$ the rotations about $\hat c$ are exact by every tested angle, and those about $v_1$ and $v_2$ do not exceed $J=0.67$. In $(v_1,v_2,c^\perp)$ the rotations about $c^\perp$ are exact by every angle, and those about $v_1$ and $v_2$ are approximate ($J\approx0.97$--$0.99$): again a nearly, but not exactly, round quarter. The mixed quarters have the same body share, $27.1\%$, and only the mirror orthogonal to their second axis, the corresponding half-turn and the point reflection; their rotations about $\hat c$ stop at $J=0.906$, so the impression of revolution given by the image is not correct. The orbits generate the whole algebra (dimension $8$) in the three quarters through $\hat c$, and a subalgebra of dimension $3$ in the lower-right one.

The absence of revolutions in Plate~B07 is therefore an effect of the frame: the continuous symmetry of the generic constant is there, but it turns planes that no coordinate frame contains.
\clearpage

\subsubsection{R10. The twisted revolutions of \texorpdfstring{$W$}{W} at \texorpdfstring{$\theta=90^\circ$}{theta=90\textdegree{}}}\label{atl:pl:R10}\mbox{}\par\vspace{0.5mm}
\noindent\includegraphics[width=\textwidth]{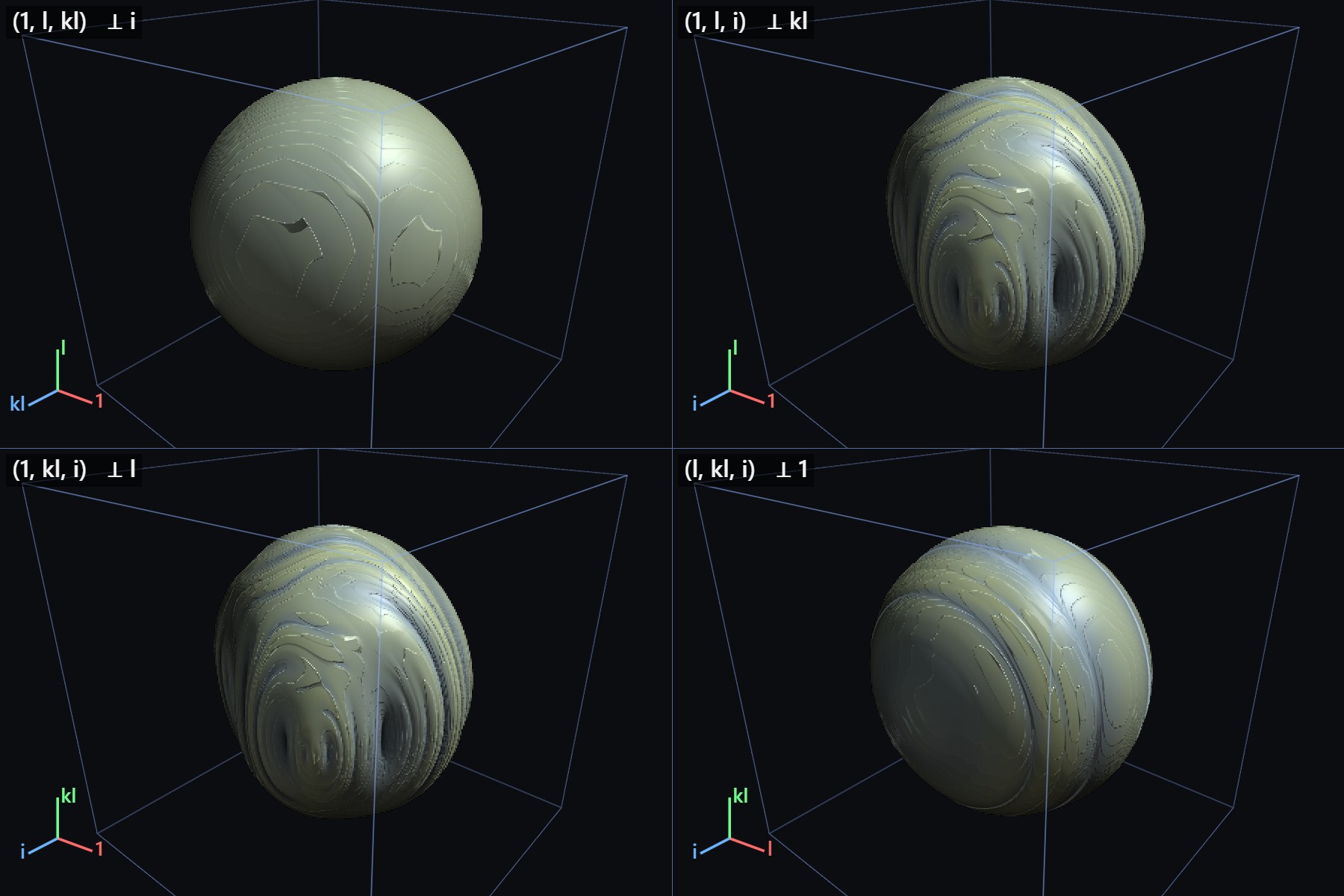}\par\vspace{1mm}
\noindent{\footnotesize \textbf{Map} $z\mapsto z\diamond z+c$ $\cdot$ \textbf{Algebra} custom: Petersson, $\gamma=(-1,-1,-1)$, $\theta=90^\circ$ $\cdot$ $c=(-0.086,\,\allowbreak 0.231,\,\allowbreak 0,\,\allowbreak 0,\,\allowbreak 0,\,\allowbreak 0,\,\allowbreak 0,\,\allowbreak 0)$ $\cdot$ \textbf{Frame} $X,Y,Z,N=1,l,kl,i$ $\cdot$ half-extent 1.4 $\cdot$ 20 iterations $\cdot$ \textbf{Coloring} 8D orbit trap (point), gradient Pastel}\par\vspace{1mm}
\noindent{\footnotesize\begin{tabular}{@{}lllll@{}}\toprule Slice (X, Y, Z) & Body [95\%] & Exact ($J\geq0.9995$) & Approximate ($0.97\leq J<0.9995$) & Dim\\\midrule
($1, l, kl$) $\perp i$ & 51.8\% {\scriptsize[51.5, 52.1]} & $M$: xyz; $R_x$: all; $R_y$: all; $R_z$: all; $P$ & -- & 4 \\
($1, l, i$) $\perp kl$ & 44.6\% {\scriptsize[44.3, 44.9]} & $M$: y; $R_y$: 180; $P$ & -- & 3 \\
($1, kl, i$) $\perp l$ & 44.6\% {\scriptsize[44.3, 44.9]} & $M$: y; $R_y$: 180; $P$ & -- & 3 \\
($l, kl, i$) $\perp 1$ & 46.4\% {\scriptsize[46.1, 46.7]} & $M$: xyz; $R_x$: 180; $R_y$: 180; $R_z$: all; $P$ & -- & 4 \\
\bottomrule\end{tabular}}\par\vspace{2mm}
The Petersson algebra at $\theta=90^\circ$, $c=c_A$, frame $(1,l,kl\,|\,i)$. The frame is chosen to test Proposition~\ref{atl:prop:theta90}. The upper-left quarter $(1,l,kl)$ lies in $\RR\oplus\HH l$, where $z*z=x_0^2+|q|^2$: the first iterate depends only on the norm, and the quarter must be a ball. The lower-right quarter $(l,kl,i)$ lies in the subspace $W=\langle1,i\rangle\oplus\langle l,kl\rangle$, closed under squaring, where the rotations $R_\varphi$ of the plane $\langle l,kl\rangle$ satisfy the twisted equivariance $f(R_\varphi z)=R_{-\varphi}f(z)$ \eqref{atl:eq:W}: it must be a body of revolution about $i$. The mixed quarters $(1,l,i)$ and $(1,kl,i)$, also contained in $W$, are exchanged by $R_{90^\circ}$ and must have the same shape.

The image agrees. The upper-left quarter is a smooth ball, with only a few colouring lines on its surface. The two mixed quarters show the same egg-shaped body, with two vortices at its lower end. The lower-right quarter is a rounded, slightly elongated body with fine ridges on its surface.

The measurements confirm all three predictions. In $(1,l,kl)$ every rotation about every axis is exact: the quarter is a ball, occupying $51.8\%$ of the sample. In $(l,kl,i)$ the rotations about $i$ are exact by every tested angle, and those about $l$ and $kl$ stay below $0.92$. The mixed quarters have the same body share, $44.6\%$, the same values of $J$, and only the mirror orthogonal to their second axis, the corresponding half-turn and the point reflection.

Quite surprisingly, none of these symmetries comes from an automorphism: at $\theta=90^\circ$ the stabilizer of $c$ in the automorphism group is still the torus of Proposition~\ref{atl:prop:torus}, which does not preserve the plane $\langle l,kl\rangle$, since $l$ and $kl$ belong to different root planes. The revolutions of this plate come from the twisted equivariance on $W$. Compare with the controls at $\theta=60^\circ$ and $\theta=120^\circ$, Plates~R11 and R12.
\clearpage

\subsubsection{R11. Control at \texorpdfstring{$\theta=60^\circ$}{theta=60\textdegree{}}: the subspace \texorpdfstring{$W$}{W} without revolutions}\label{atl:pl:R11}\mbox{}\par\vspace{0.5mm}
\noindent\includegraphics[width=\textwidth]{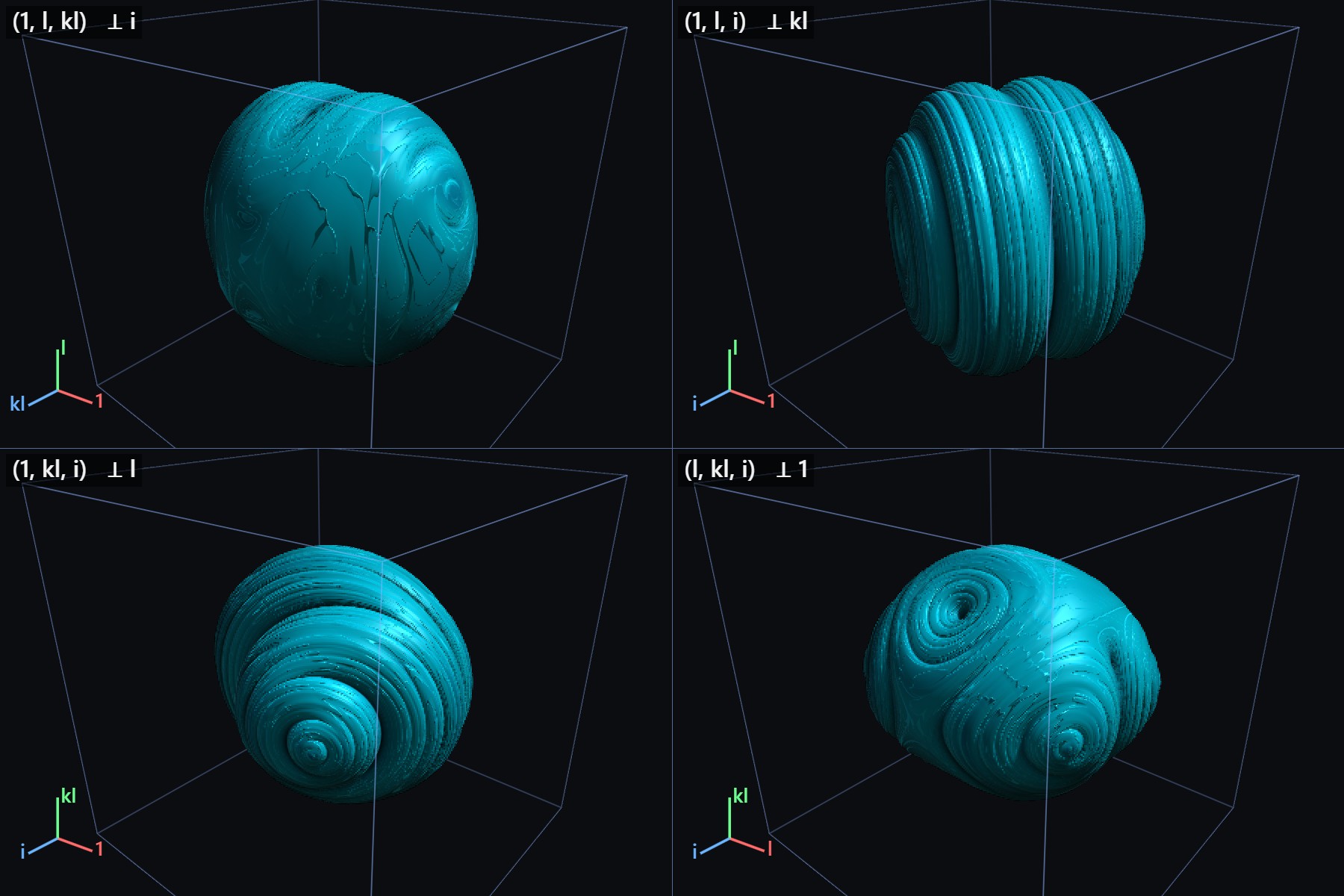}\par\vspace{1mm}
\noindent{\footnotesize \textbf{Map} $z\mapsto z\diamond z+c$ $\cdot$ \textbf{Algebra} custom: Petersson, $\gamma=(-1,-1,-1)$, $\theta=60^\circ$ $\cdot$ $c=(-0.086,\,\allowbreak 0.231,\,\allowbreak 0,\,\allowbreak 0,\,\allowbreak 0,\,\allowbreak 0,\,\allowbreak 0,\,\allowbreak 0)$ $\cdot$ \textbf{Frame} $X,Y,Z,N=1,l,kl,i$ $\cdot$ half-extent 1.4 $\cdot$ 20 iterations $\cdot$ \textbf{Coloring} 8D orbit trap (point), gradient Neon}\par\vspace{1mm}
\noindent{\footnotesize\begin{tabular}{@{}lllll@{}}\toprule Slice (X, Y, Z) & Body [95\%] & Exact ($J\geq0.9995$) & Approximate ($0.97\leq J<0.9995$) & Dim\\\midrule
($1, l, kl$) $\perp i$ & 44.1\% {\scriptsize[43.8, 44.4]} & $M$: xyz; $R_x$: 180; $R_y$: 180; $R_z$: 180; $P$ & -- & 5 \\
($1, l, i$) $\perp kl$ & 46.2\% {\scriptsize[45.9, 46.5]} & $M$: y; $R_y$: 180; $P$ & -- & 3 \\
($1, kl, i$) $\perp l$ & 34.8\% {\scriptsize[34.5, 35.1]} & $M$: y; $R_y$: 180; $P$ & -- & 3 \\
($l, kl, i$) $\perp 1$ & 38.8\% {\scriptsize[38.5, 39.1]} & $M$: xyz; $R_x$: 180; $R_y$: 180; $R_z$: 180; $P$ & -- & 5 \\
\bottomrule\end{tabular}}\par\vspace{2mm}
The same frame $(1,l,kl\,|\,i)$ and the same constant $c_A$ as Plate~R10, at $\theta=60^\circ$. This is a control: the identities of Proposition~\ref{atl:prop:theta90} hold only at $\theta=90^\circ$, and the prediction is that the ball of the upper-left quarter and the revolution of the lower-right one disappear, and that the mixed quarters are no longer congruent.

The image differs from Plate~R10 in every quarter. The upper-left body is roughly round but irregular, with cracks and swirls on its surface; the upper-right quarter $(1,l,i)$ is a ribbed body with two lobes separated by a groove; the lower-left quarter $(1,kl,i)$ is a body with a spiral of rings; the lower-right quarter shows two separate vortices.

The measurements confirm the prediction. No quarter has a rotation other than the half-turns: $(1,l,kl)$ and $(l,kl,i)$ have the three mirrors, the three half-turns and the point reflection, the mixed quarters only one mirror, one half-turn and the point reflection. The upper-left quarter is not a ball (the best rotation, by $90^\circ$ about $kl$, has $J=0.947$), and the mixed quarters have different body shares, $46.2\%$ and $34.8\%$.

It is worth noticing that the measurements of this plate and of Plate~R12 ($\theta=120^\circ$) coincide up to the exchange of $l$ and $kl$: the body shares $46.2\%$ and $34.8\%$ of the mixed quarters are swapped, and so are the indices of the rotations in the other two quarters; in the images, the two-lobed body and the spiral body change places. This suggests that the two algebras, with this constant, are related by an isometry exchanging the roles of $l$ and $kl$; we have not verified it algebraically.
\clearpage

\subsubsection{R12. Control in the Okubo algebra: the subspace \texorpdfstring{$W$}{W} without revolutions}\label{atl:pl:R12}\mbox{}\par\vspace{0.5mm}
\noindent\includegraphics[width=\textwidth]{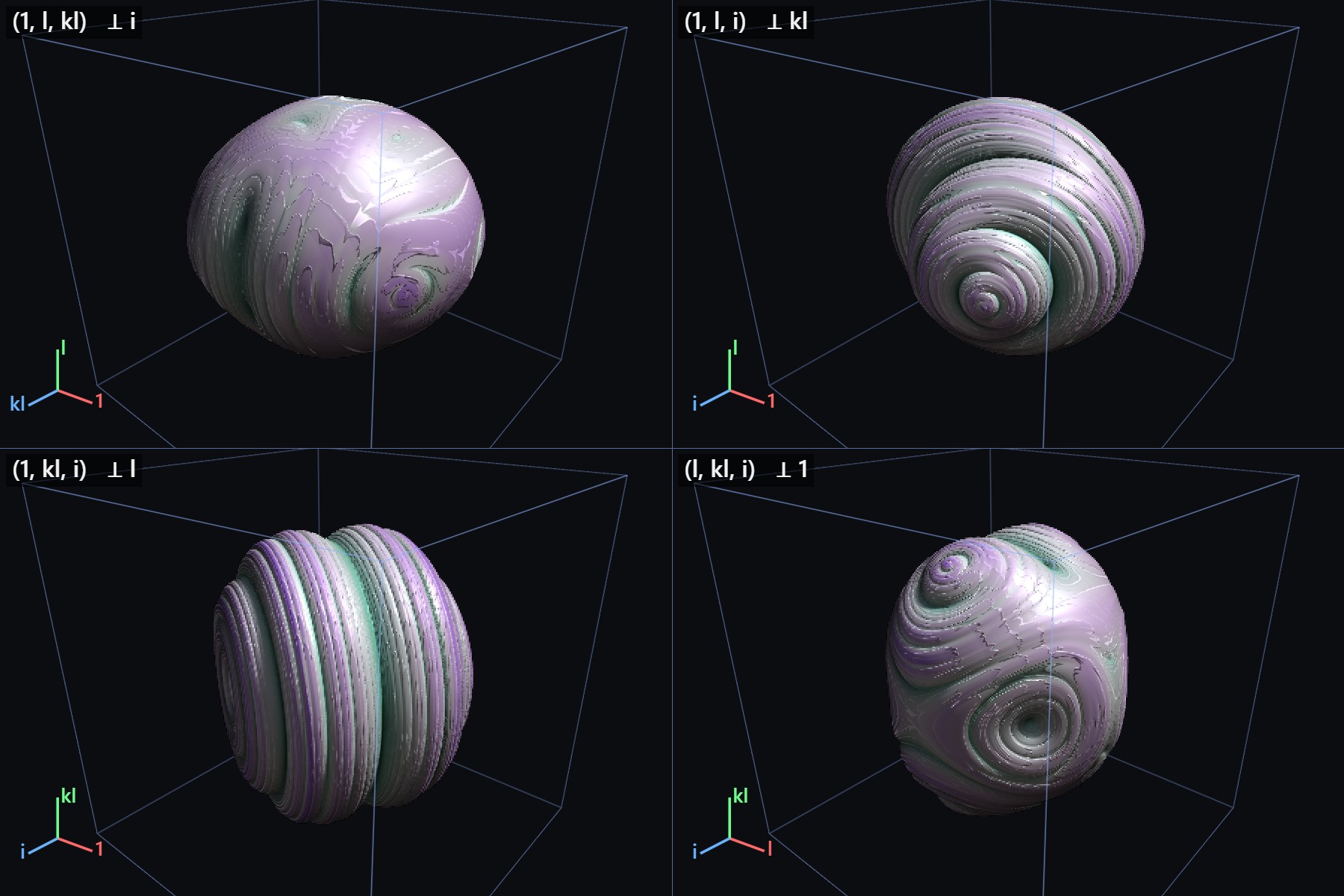}\par\vspace{1mm}
\noindent{\footnotesize \textbf{Map} $z\mapsto z\diamond z+c$ $\cdot$ \textbf{Algebra} Okubo algebra (Petersson, $\gamma=(-1,-1,-1)$, $\theta=120^\circ$) $\cdot$ $c=(-0.086,\,\allowbreak 0.231,\,\allowbreak 0,\,\allowbreak 0,\,\allowbreak 0,\,\allowbreak 0,\,\allowbreak 0,\,\allowbreak 0)$ $\cdot$ \textbf{Frame} $X,Y,Z,N=1,l,kl,i$ $\cdot$ half-extent 1.4 $\cdot$ 20 iterations $\cdot$ \textbf{Coloring} 8D orbit trap (point), gradient Aurora}\par\vspace{1mm}
\noindent{\footnotesize\begin{tabular}{@{}lllll@{}}\toprule Slice (X, Y, Z) & Body [95\%] & Exact ($J\geq0.9995$) & Approximate ($0.97\leq J<0.9995$) & Dim\\\midrule
($1, l, kl$) $\perp i$ & 44.0\% {\scriptsize[43.7, 44.3]} & $M$: xyz; $R_x$: 180; $R_y$: 180; $R_z$: 180; $P$ & -- & 5 \\
($1, l, i$) $\perp kl$ & 34.8\% {\scriptsize[34.5, 35.1]} & $M$: y; $R_y$: 180; $P$ & -- & 3 \\
($1, kl, i$) $\perp l$ & 46.2\% {\scriptsize[45.9, 46.5]} & $M$: y; $R_y$: 180; $P$ & -- & 3 \\
($l, kl, i$) $\perp 1$ & 38.8\% {\scriptsize[38.5, 39.1]} & $M$: xyz; $R_x$: 180; $R_y$: 180; $R_z$: 180; $P$ & -- & 5 \\
\bottomrule\end{tabular}}\par\vspace{2mm}
The Okubo algebra ($\theta=120^\circ$) with $c=c_A$ in the frame $(1,l,kl\,|\,i)$ of Plate~R10: the second control for the subspace $W$. The frame contains no root plane of the torus (it contains $l$ from $V_2$ and $kl$ from $V_3$), so Theorem~A predicts no revolution, and the identities of Proposition~\ref{atl:prop:theta90} do not hold at this angle.

The quarters are compact rounded bodies. The upper-left quarter $(1,l,kl)$ is irregular, with swirls and cracks; the upper-right quarter $(1,l,i)$ is a body with a spiral of rings; the lower-left quarter $(1,kl,i)$ is a ribbed body with two lobes separated by a groove; the lower-right quarter is an ovoid with two vortices.

The measurements confirm the prediction: only the finite symmetries remain (in $(1,l,kl)$ and $(l,kl,i)$ the three mirrors, the three half-turns and the point reflection; in the mixed quarters one mirror, one half-turn and the point reflection), and no rotation by a smaller angle reaches $0.95$. The mixed quarters are not congruent: $34.8\%$ and $46.2\%$. The quarter $(1,l,i)$ is the same section as the lower-left quarter of Plate~C03, with the same constant, and it has the same body share, $34.8\%$.

As noted for Plate~R11, the two controls correspond to each other: the numbers of this plate are those of Plate~R11 with $l$ and $kl$ exchanged, and the two-lobed and the spiral bodies change places.
\clearpage

\subsubsection{R13. The Petersson path at \texorpdfstring{$c_A$}{cA}: \texorpdfstring{$\theta=30^\circ$}{theta=30\textdegree{}}}\label{atl:pl:R13}\mbox{}\par\vspace{0.5mm}
\noindent\includegraphics[width=\textwidth]{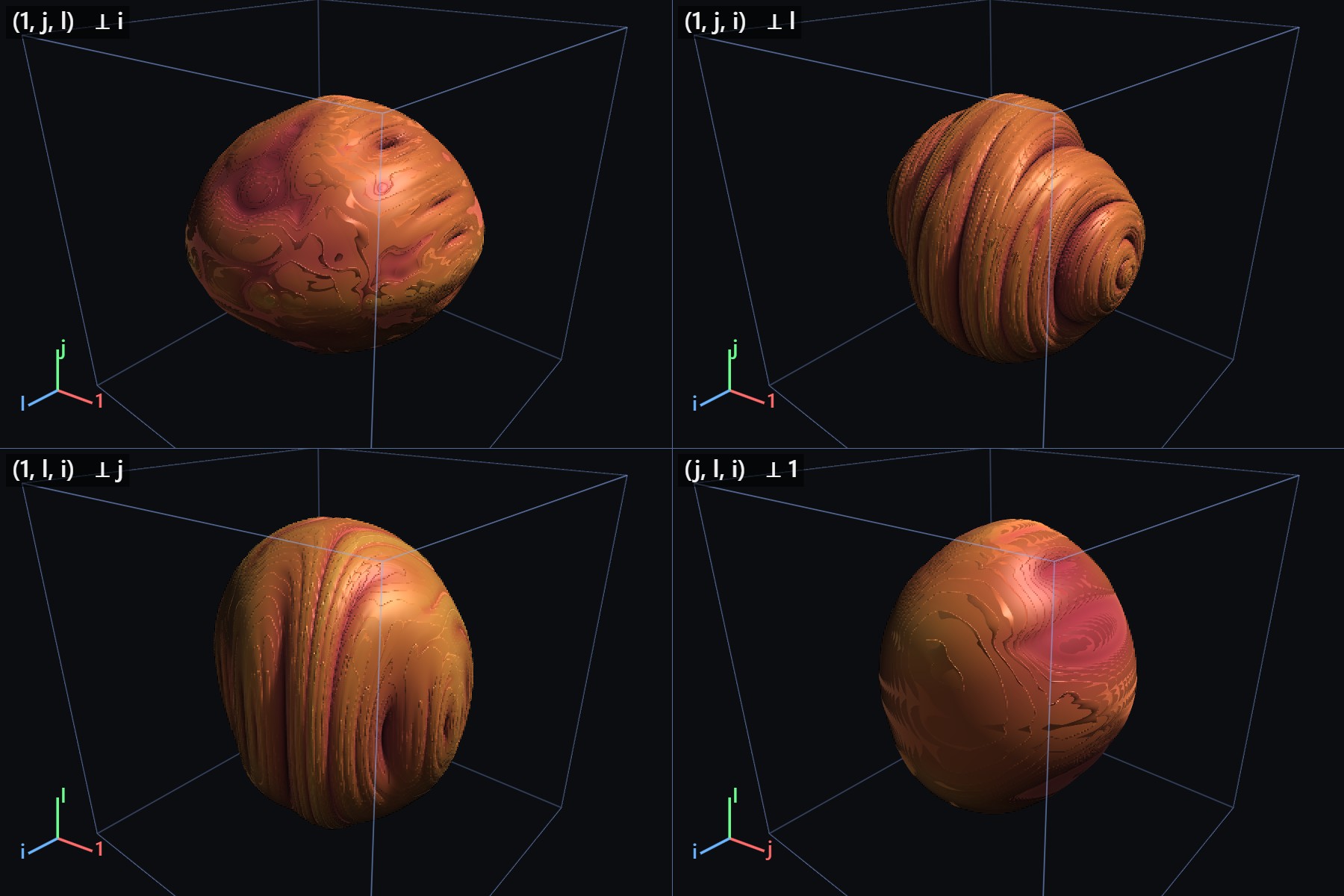}\par\vspace{1mm}
\noindent{\footnotesize \textbf{Map} $z\mapsto z\diamond z+c$ $\cdot$ \textbf{Algebra} custom: Petersson, $\gamma=(-1,-1,-1)$, $\theta=30^\circ$ $\cdot$ $c=(-0.086,\,\allowbreak 0.231,\,\allowbreak 0,\,\allowbreak 0,\,\allowbreak 0,\,\allowbreak 0,\,\allowbreak 0,\,\allowbreak 0)$ $\cdot$ \textbf{Frame} $X,Y,Z,N=1,j,l,i$ $\cdot$ half-extent 1.4 $\cdot$ 20 iterations $\cdot$ \textbf{Coloring} 8D orbit trap (point), gradient Sunset}\par\vspace{1mm}
\noindent{\footnotesize\begin{tabular}{@{}lllll@{}}\toprule Slice (X, Y, Z) & Body [95\%] & Exact ($J\geq0.9995$) & Approximate ($0.97\leq J<0.9995$) & Dim\\\midrule
($1, j, l$) $\perp i$ & 42.1\% {\scriptsize[41.8, 42.4]} & $M$: xyz; $R_x$: 180; $R_y$: 180; $R_z$: 180; $P$ & -- & 5 \\
($1, j, i$) $\perp l$ & 35.9\% {\scriptsize[35.6, 36.2]} & $M$: y; $R_y$: 180; $P$ & -- & 3 \\
($1, l, i$) $\perp j$ & 48.0\% {\scriptsize[47.7, 48.3]} & $M$: y; $R_y$: 180; $P$ & -- & 3 \\
($j, l, i$) $\perp 1$ & 40.1\% {\scriptsize[39.8, 40.4]} & $M$: xyz; $R_x$: 180; $R_y$: 180; $R_z$: 180; $P$ & -- & 5 \\
\bottomrule\end{tabular}}\par\vspace{2mm}
The first of the new plates of the Petersson path at fixed $c=c_A$, in the frame $F_C=(1,j,l\,|\,i)$ of Plates~C02 ($\theta=0^\circ$), C03 ($\theta=120^\circ$) and C05 ($\theta=90^\circ$); Plates~R13--R16 add $\theta=30^\circ$, $60^\circ$, $90^\circ$ and $150^\circ$. For every angle the torus of Proposition~\ref{atl:prop:torus} and the quaternionic subalgebra are there, so the frame $F_C$, which contains no root plane, should show the same finite symmetry pattern at every angle; the section $(j,l,i)$ is predicted to be a body of revolution about $l$ only at $\theta=0^\circ$, where the group $O(6)$ of the para-octonions acts, and at $\theta=90^\circ$, by the twisted equivariance of Proposition~\ref{atl:prop:theta90}.

At $\theta=30^\circ$ the upper-left and lower-right quarters are irregular ovoids with swirls and small vortices; the lower-left quarter $(1,l,i)$ is an elongated body with vertical ridges and vortices at its lower end. The upper-right quarter $(1,j,i)$, a ribbed body with a concentric cap, lies in the quaternions: the product there is $\bar x\,\bar y$ for every $\theta$, and this quarter is the same in all the plates of the series and in Plates~C02, C03 and C05.

The measurements confirm the prediction. The quarter $(j,l,i)$ has only the finite symmetries (mirrors, half-turns, point reflection); its best rotation about $l$, by $37^\circ$, has $J=0.918$. The quarters $(1,j,l)$ and $(j,l,i)$ have the three mirrors, the three half-turns and the point reflection, the quarters $(1,j,i)$ and $(1,l,i)$ one mirror, one half-turn and the point reflection: the same pattern as Plate~C03. The quaternionic quarter occupies $35.9\%$ of the sample, as in all the plates of the series ($35.6$--$36.1\%$).
\clearpage

\subsubsection{R14. The Petersson path at \texorpdfstring{$c_A$}{cA}: \texorpdfstring{$\theta=60^\circ$}{theta=60\textdegree{}}}\label{atl:pl:R14}\mbox{}\par\vspace{0.5mm}
\noindent\includegraphics[width=\textwidth]{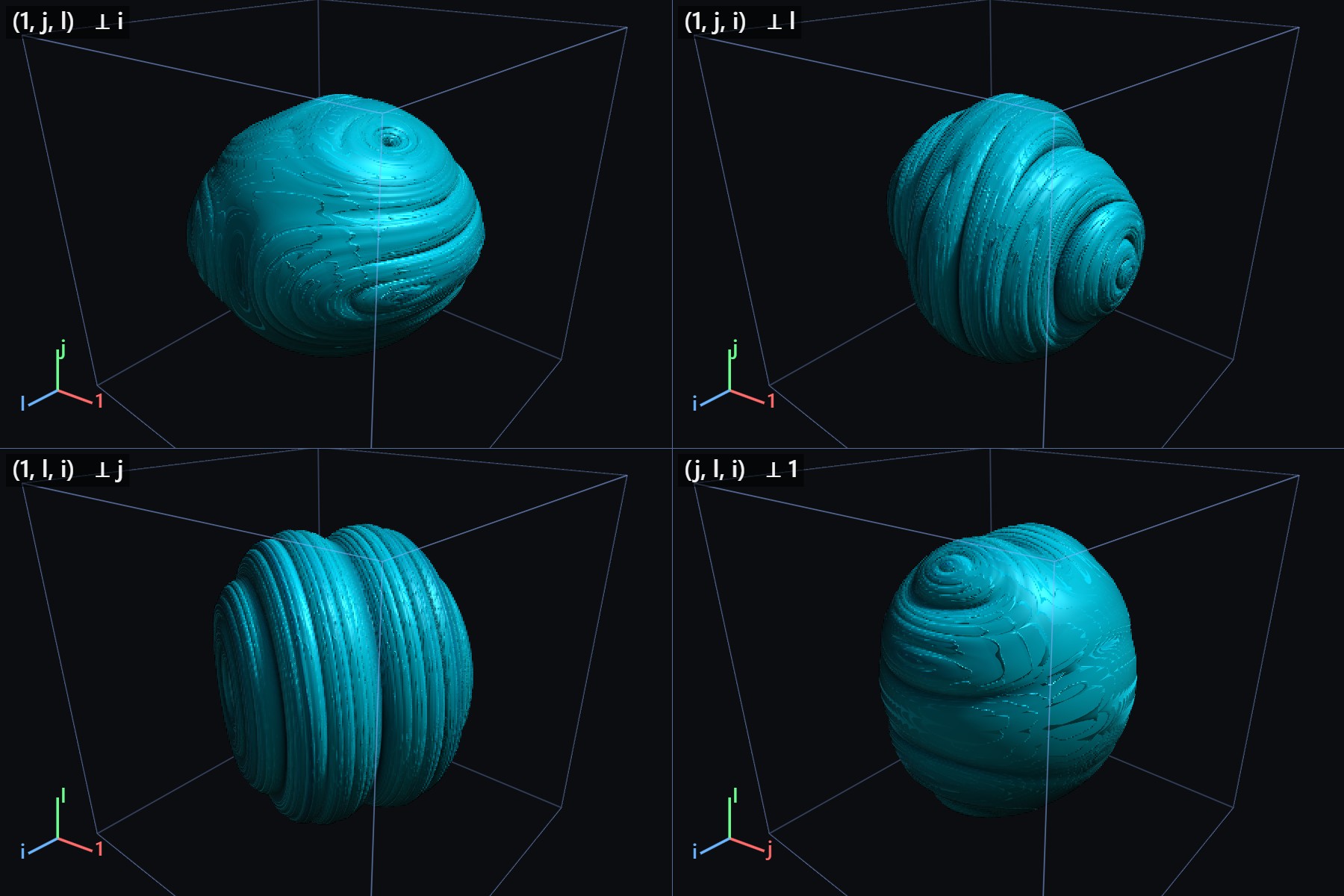}\par\vspace{1mm}
\noindent{\footnotesize \textbf{Map} $z\mapsto z\diamond z+c$ $\cdot$ \textbf{Algebra} custom: Petersson, $\gamma=(-1,-1,-1)$, $\theta=60^\circ$ $\cdot$ $c=(-0.086,\,\allowbreak 0.231,\,\allowbreak 0,\,\allowbreak 0,\,\allowbreak 0,\,\allowbreak 0,\,\allowbreak 0,\,\allowbreak 0)$ $\cdot$ \textbf{Frame} $X,Y,Z,N=1,j,l,i$ $\cdot$ half-extent 1.4 $\cdot$ 20 iterations $\cdot$ \textbf{Coloring} 8D orbit trap (point), gradient Neon}\par\vspace{1mm}
\noindent{\footnotesize\begin{tabular}{@{}lllll@{}}\toprule Slice (X, Y, Z) & Body [95\%] & Exact ($J\geq0.9995$) & Approximate ($0.97\leq J<0.9995$) & Dim\\\midrule
($1, j, l$) $\perp i$ & 42.4\% {\scriptsize[42.1, 42.7]} & $M$: xyz; $R_x$: 180; $R_y$: 180; $R_z$: 180; $P$ & -- & 5 \\
($1, j, i$) $\perp l$ & 36.1\% {\scriptsize[35.8, 36.4]} & $M$: y; $R_y$: 180; $P$ & -- & 3 \\
($1, l, i$) $\perp j$ & 46.3\% {\scriptsize[46.0, 46.6]} & $M$: y; $R_y$: 180; $P$ & -- & 3 \\
($j, l, i$) $\perp 1$ & 41.7\% {\scriptsize[41.4, 42.0]} & $M$: xyz; $R_x$: 180; $R_y$: 180; $R_z$: 180; $P$ & -- & 5 \\
\bottomrule\end{tabular}}\par\vspace{2mm}
The Petersson path at $c=c_A$, $\theta=60^\circ$, frame $F_C=(1,j,l\,|\,i)$. The prediction is the same as for Plate~R13: the finite symmetry pattern of the frame $F_C$, and no revolution about $l$ in the quarter $(j,l,i)$.

The quaternionic quarter $(1,j,i)$ is the same as in Plate~R13. The lower-left quarter $(1,l,i)$ is the two-lobed ribbed body already seen in the upper-right quarter of Plate~R11, which shows the same section at the same angle. The lower-right quarter $(j,l,i)$ is an ovoid with a spiral near its top; the upper-left quarter a roundish body with swirls.

The measurements confirm the prediction: every quarter has only the finite symmetries, in the same pattern as Plates~R13 and C03, and the best rotation about $l$ in $(j,l,i)$, by $37^\circ$, has $J=0.937$. The body share of $(1,l,i)$, $46.3\%$, agrees with that of the same section in Plate~R11, $46.2\%$, within the sampling error; the quaternionic quarter has $36.1\%$.
\clearpage

\subsubsection{R15. The Petersson path at \texorpdfstring{$c_A$}{cA}: \texorpdfstring{$\theta=90^\circ$}{theta=90\textdegree{}}}\label{atl:pl:R15}\mbox{}\par\vspace{0.5mm}
\noindent\includegraphics[width=\textwidth]{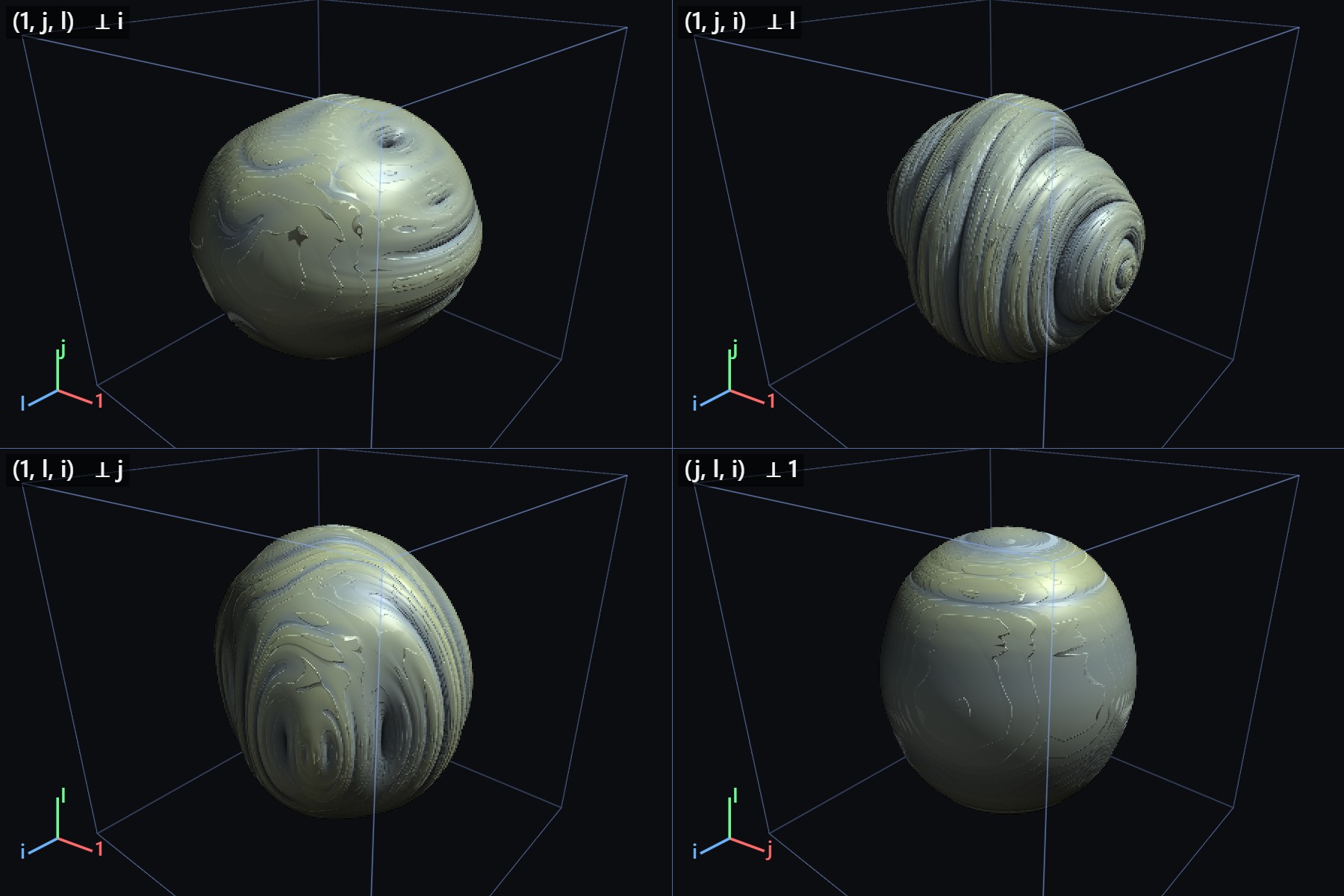}\par\vspace{1mm}
\noindent{\footnotesize \textbf{Map} $z\mapsto z\diamond z+c$ $\cdot$ \textbf{Algebra} custom: Petersson, $\gamma=(-1,-1,-1)$, $\theta=90^\circ$ $\cdot$ $c=(-0.086,\,\allowbreak 0.231,\,\allowbreak 0,\,\allowbreak 0,\,\allowbreak 0,\,\allowbreak 0,\,\allowbreak 0,\,\allowbreak 0)$ $\cdot$ \textbf{Frame} $X,Y,Z,N=1,j,l,i$ $\cdot$ half-extent 1.4 $\cdot$ 20 iterations $\cdot$ \textbf{Coloring} 8D orbit trap (point), gradient Pastel}\par\vspace{1mm}
\noindent{\footnotesize\begin{tabular}{@{}lllll@{}}\toprule Slice (X, Y, Z) & Body [95\%] & Exact ($J\geq0.9995$) & Approximate ($0.97\leq J<0.9995$) & Dim\\\midrule
($1, j, l$) $\perp i$ & 43.5\% {\scriptsize[43.2, 43.8]} & $M$: xyz; $R_x$: 180; $R_y$: 180; $R_z$: 180; $P$ & -- & 5 \\
($1, j, i$) $\perp l$ & 35.6\% {\scriptsize[35.3, 35.9]} & $M$: y; $R_y$: 180; $P$ & -- & 3 \\
($1, l, i$) $\perp j$ & 44.4\% {\scriptsize[44.1, 44.7]} & $M$: y; $R_y$: 180; $P$ & -- & 3 \\
($j, l, i$) $\perp 1$ & 40.8\% {\scriptsize[40.5, 41.1]} & $M$: xyz; $R_x$: 180; $R_y$: all; $R_z$: 180; $P$ & -- & 5 \\
\bottomrule\end{tabular}}\par\vspace{2mm}
The Petersson path at $c=c_A$, $\theta=90^\circ$, frame $F_C=(1,j,l\,|\,i)$. The algebra, the constant and the frame are those of Plate~C05: this plate was rendered again to complete the series, and its measurements, on an independent sample, agree with those of Plate~C05 within the sampling error. The prediction is that the quarter $(j,l,i)$ is an exact body of revolution about $l$: since $z*z\in W$ for $z\in\langle i,j,l\rangle$, the twisted equivariance of Proposition~\ref{atl:prop:theta90} applies.

The lower-right quarter is an egg-shaped body with a concentric system of rings around its top, where the vertical axis $l$ leaves it. The quaternionic quarter $(1,j,i)$ is the same as in the other plates of the series, and the lower-left quarter $(1,l,i)$ is the egg with two vortices of the mixed quarters of Plate~R10, which contain the same section at the same angle.

The measurements confirm the prediction: in $(j,l,i)$ the rotations about $l$ are exact by every tested angle, including $37^\circ$, and those about $j$ and $i$ stay below $0.91$. The other quarters have only the finite symmetries of the frame $F_C$. The body share of $(j,l,i)$ is $40.8\%$ ($41.1\%$ in Plate~C05), that of $(1,l,i)$ $44.4\%$ ($44.6\%$ in Plate~R10).

Among the angles of the series, $0^\circ$ (Plate~C02, where $(j,l,i)$ is even a ball), $30^\circ$, $60^\circ$, $90^\circ$, $120^\circ$ (Plate~C03) and $150^\circ$, only $0^\circ$ and $90^\circ$ produce a revolution about $l$. The angle $90^\circ$ is an exception, not a stage of a continuous transition, and the symmetry is not due to an automorphism.
\clearpage

\subsubsection{R16. The Petersson path at \texorpdfstring{$c_A$}{cA}: \texorpdfstring{$\theta=150^\circ$}{theta=150\textdegree{}}}\label{atl:pl:R16}\mbox{}\par\vspace{0.5mm}
\noindent\includegraphics[width=\textwidth]{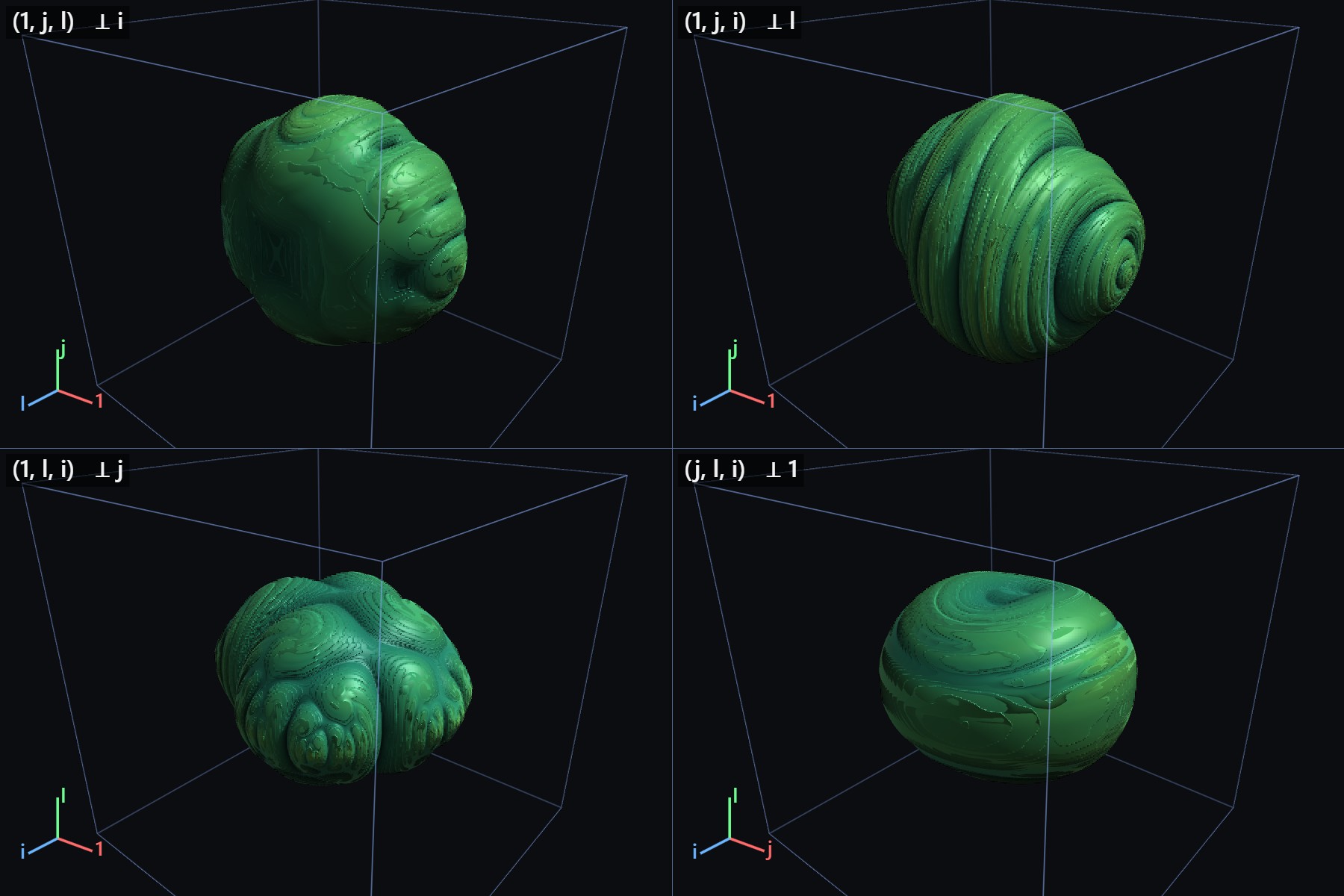}\par\vspace{1mm}
\noindent{\footnotesize \textbf{Map} $z\mapsto z\diamond z+c$ $\cdot$ \textbf{Algebra} custom: Petersson, $\gamma=(-1,-1,-1)$, $\theta=150^\circ$ $\cdot$ $c=(-0.086,\,\allowbreak 0.231,\,\allowbreak 0,\,\allowbreak 0,\,\allowbreak 0,\,\allowbreak 0,\,\allowbreak 0,\,\allowbreak 0)$ $\cdot$ \textbf{Frame} $X,Y,Z,N=1,j,l,i$ $\cdot$ half-extent 1.4 $\cdot$ 20 iterations $\cdot$ \textbf{Coloring} 8D orbit trap (point), gradient Viridis-like}\par\vspace{1mm}
\noindent{\footnotesize\begin{tabular}{@{}lllll@{}}\toprule Slice (X, Y, Z) & Body [95\%] & Exact ($J\geq0.9995$) & Approximate ($0.97\leq J<0.9995$) & Dim\\\midrule
($1, j, l$) $\perp i$ & 28.6\% {\scriptsize[28.3, 28.9]} & $M$: xyz; $R_x$: 180; $R_y$: 180; $R_z$: 180; $P$ & -- & 5 \\
($1, j, i$) $\perp l$ & 36.0\% {\scriptsize[35.7, 36.3]} & $M$: y; $R_y$: 180; $P$ & -- & 3 \\
($1, l, i$) $\perp j$ & 25.8\% {\scriptsize[25.5, 26.0]} & $M$: y; $R_y$: 180; $P$ & -- & 3 \\
($j, l, i$) $\perp 1$ & 27.4\% {\scriptsize[27.1, 27.6]} & $M$: xyz; $R_x$: 180; $R_y$: 180; $R_z$: 180; $P$ & -- & 5 \\
\bottomrule\end{tabular}}\par\vspace{2mm}
The Petersson path at $c=c_A$, $\theta=150^\circ$, beyond the Okubo angle, frame $F_C=(1,j,l\,|\,i)$. The prediction: the finite symmetry pattern of the frame $F_C$, and no revolution about $l$ in $(j,l,i)$.

The bodies of this plate are visibly smaller than in the rest of the series, except the quaternionic quarter $(1,j,i)$, which is again the same. The upper-left quarter is a lumpy body with folds on its right side; the lower-left quarter $(1,l,i)$ is a compact body with several rounded protrusions at its lower end; the lower-right quarter $(j,l,i)$ is a flattened ovoid with a whorl on top.

The measurements confirm the prediction: every quarter has only the finite symmetries, in the same pattern as the other angles, and the best rotation about $l$ in $(j,l,i)$, by $37^\circ$, has $J=0.942$. The body shares are $28.6\%$, $25.8\%$ and $27.4\%$ in the three non-quaternionic quarters, against $40$--$48\%$ at $\theta=30^\circ$, $60^\circ$ and $90^\circ$; the quaternionic quarter keeps its $36.0\%$. For this constant the filled set shrinks considerably past the Okubo angle.
\clearpage

\subsubsection{R17. Kaplansky's trick at equal iterations: a null test}\label{atl:pl:R17}\mbox{}\par\vspace{0.5mm}
\noindent\includegraphics[width=\textwidth]{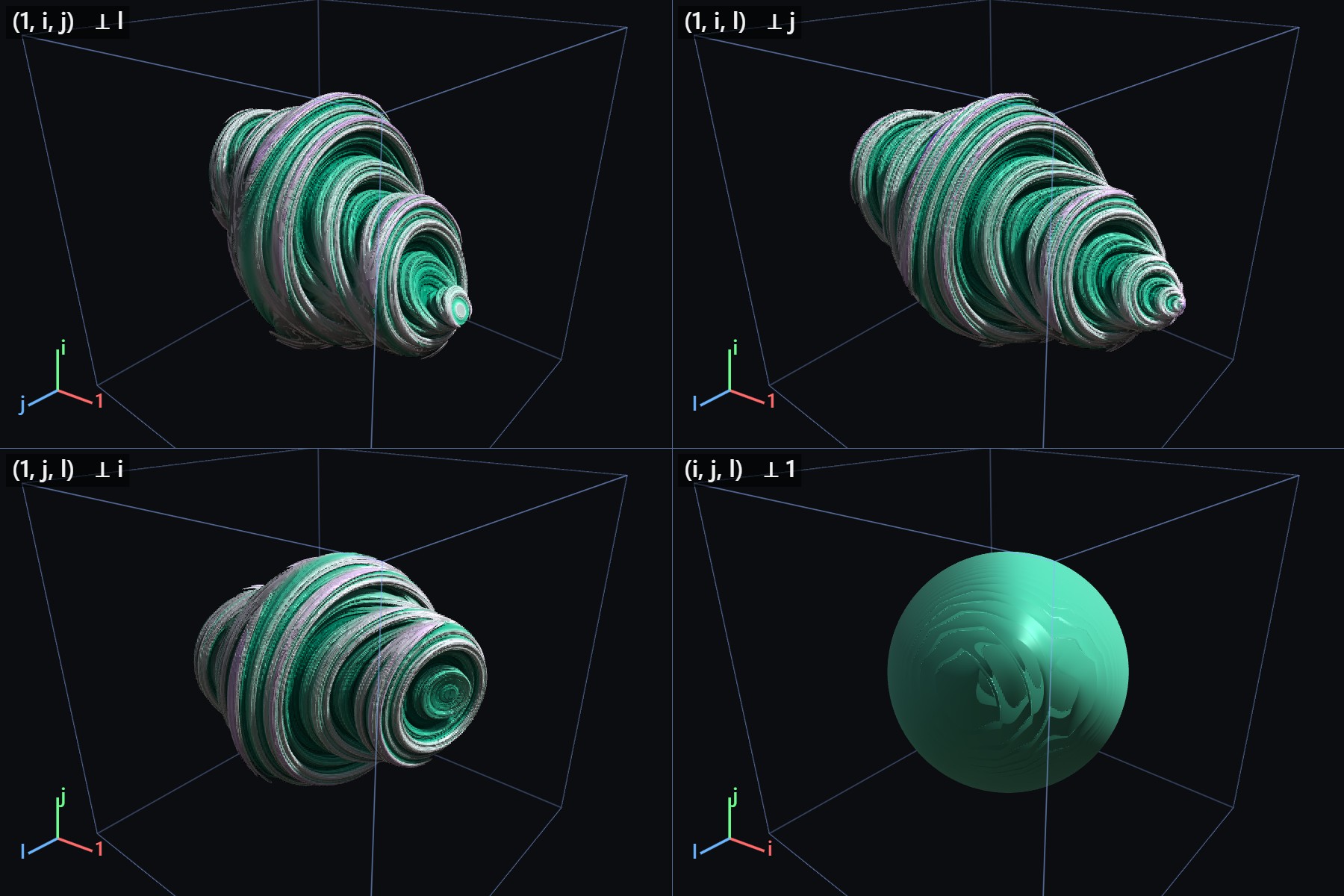}\par\vspace{1mm}
\noindent{\footnotesize \textbf{Map} $z\mapsto(a\diamond z)\diamond(z\diamond b)+c$ $\cdot$ \textbf{Algebra} Okubo algebra (Petersson, $\gamma=(-1,-1,-1)$, $\theta=120^\circ$) $\cdot$ $c=(-0.35,\,\allowbreak 0.45,\,\allowbreak 0,\,\allowbreak 0,\,\allowbreak 0.35,\,\allowbreak 0,\,\allowbreak 0,\,\allowbreak 0)$ $\cdot$ $a=(1,\,\allowbreak 0,\,\allowbreak 0,\,\allowbreak 0,\,\allowbreak 0,\,\allowbreak 0,\,\allowbreak 0,\,\allowbreak 0)$ $\cdot$ $b=(1,\,\allowbreak 0,\,\allowbreak 0,\,\allowbreak 0,\,\allowbreak 0,\,\allowbreak 0,\,\allowbreak 0,\,\allowbreak 0)$ $\cdot$ \textbf{Frame} $X,Y,Z,N=1,i,j,l$ $\cdot$ half-extent 1.4 $\cdot$ 20 iterations $\cdot$ \textbf{Coloring} 8D orbit trap (point), gradient Aurora}\par\vspace{1mm}
\noindent{\footnotesize\begin{tabular}{@{}lllll@{}}\toprule Slice (X, Y, Z) & Body [95\%] & Exact ($J\geq0.9995$) & Approximate ($0.97\leq J<0.9995$) & Dim\\\midrule
($1, i, j$) $\perp l$ & 21.4\% {\scriptsize[21.1, 21.6]} & $M$: z; $R_z$: 180; $P$ & -- & 8 \\
($1, i, l$) $\perp j$ & 22.4\% {\scriptsize[22.2, 22.7]} & $P$ & -- & 4 \\
($1, j, l$) $\perp i$ & 21.0\% {\scriptsize[20.7, 21.2]} & $M$: y; $R_y$: 180; $P$ & -- & 8 \\
($i, j, l$) $\perp 1$ & 29.1\% {\scriptsize[28.8, 29.4]} & $M$: xyz; $R_x$: all; $R_y$: all; $R_z$: all; $P$ & -- & 8 \\
\bottomrule\end{tabular}}\par\vspace{2mm}
The bilateral map $(a*z)*(z*b)+c$ over the Okubo algebra with $a=b=1$ and the constant of Plate~J01, rendered with $20$ iterations, as Plate~J02, instead of the $16$ of Plate~J01. By Kaplansky's trick (Proposition~\ref{atl:prop:kaplansky}), $(1*x)*(y*1)=xy$ is the octonionic product: the map does not merely resemble the octonionic quadratic map, it coincides with it. The prediction is that, at equal iterations, this plate and Plate~J02 are the same picture, and that the difference between Plates~J01 and J02 was due only to the iteration count.

The image cannot be told apart from Plate~J02: the three quarters through $1$ are ribbed bodies whose ridges wind about tilted axes, and the imaginary quarter $(i,j,l)$ is a ball, as Proposition~\ref{atl:prop:imaginary} requires for the octonions. According to the comparison made for the main text, the two images differ by JPEG noise only.

The measurements confirm the prediction exactly. On the same $10^5$ points per quarter, the bodies of this plate and of Plate~J02 differ in no point at all, whereas Plate~J01, with $16$ iterations, differs from Plate~J02 in $318$, $334$ and $323$ points in the quarters through $1$ and in $73$ points in the imaginary quarter, i.e., by about $0.3$ percentage points of body share. The tables agree as well: the body shares are $21.4\%$, $22.4\%$, $21.0\%$, $29.1\%$ here and $21.4\%$, $22.3\%$, $20.9\%$, $29.1\%$ in Plate~J02 (on independent samples); the symmetries are the same, one mirror, one half-turn and the point reflection in $(1,i,j)$ and $(1,j,l)$, only the point reflection in $(1,i,l)$, and every rotation about every axis in the ball $(i,j,l)$.

One entry of the table must be read with care: the dimension of the subalgebra generated by the point and by $a$, $b$, $c$ under the Okubo product is $8$ ($4$ in the quarter $(1,i,l)$), while Plate~J02 reports $3$. The orbits of the two maps are the same, and they lie in the three-dimensional space spanned by $1$, $\im z_0$ and $\im c$ (Proposition~\ref{atl:prop:span}); the Okubonic subalgebra that contains them is simply larger than the space they actually visit.
\clearpage

\subsection{Synthesis}\label{atl:sec:synthesis}

The atlas is, in essence, a long comparison with the octonionic case, and it can be summarized by a small number of visual signatures. Table~\ref{atl:tab:signatures} collects them.

\begin{table}[h]
\centering\small
\begin{tabular}{@{}p{0.29\textwidth}p{0.33\textwidth}p{0.3\textwidth}@{}}
\toprule
Visual signature & Algebraic cause & Plates\\
\midrule
Identical quarters through $1,i$ & $O(6)$-equivariance of $z^2+c$ (Prop.~\ref{atl:prop:equivariance}) & A01, A05, B01, B05, H07\\
Bodies of revolution about $1$ & rotations fixing $1$ and $c$ & A01, C01, F01, F02, G01\\
Balls in the imaginary quarter & $z^2=-\n(z)$ for $\re z=0$ (Prop.~\ref{atl:prop:imaginary}) & A01, B05, E01, F01, H07\\
Hyperboloids, horns, saddles & indefinite norm, zero divisors & A02, A06, B02, D05--D07, F05, I05\\
Flat faces, square Mandelbrot sections & split-complex plane ($i^2=+1$), null lines $1\pm i$ & A03 (Julia), F08, L04\\
Exact cylinders & $\gamma_3=0$: $q_{n+1}=2\re(p_n)q_n+c_l$ (Prop.~\ref{atl:prop:rescaling}) & A04, B04, D04, F07, L03, P01, P02\\
Stretched bodies & $\gamma_3<0$ is a rescaled octonion algebra & D02, D03, P03\\
Centrally symmetric quarters; mirrors paired with half-turns & $z\mapsto-z$ for even maps (Rem.~\ref{atl:rem:minus}) & every plate centred at $0$\\
No revolution in the frames $F_A$, $F_C$ & the frames contain no root plane of $\mathfrak{su}(3)$ (Prop.~\ref{atl:prop:torus}) & A07, C03, C08\\
Revolutions in the root planes & the torus rotates $V_1,V_2,V_3$ (Prop.~\ref{atl:prop:torus}) & R01--R03\\
Revolutions born on the Weyl walls & stabilizer $U(2)$ for $\arg c\in60^\circ\mathbb Z$ & R05--R08\\
Revolutions of a generic constant & conjugate torus, adapted frame & B07, R09\\
Subalgebra of dimension 5--8 & $z_0,c$ generate a large $\circ$-subalgebra & A07, B07, C08, E02--E06, F11\\
Three-armed Mandelbrot sections & the tricorn moved by $PSU(3)$ (Prop.~\ref{atl:prop:okubocubic}, Thm.~B) & F03, F06, F11\\
Perfect balls, identical in para and Okubo & $(z*a)*z=\n(z)a$: a recursion on the norm (Prop.~\ref{atl:prop:symmetric}) & I02, I03\\
Revolutions about the direction of $c$, identical in para and Okubo & $(z*z)*z+c=\n(z)z+c$ commutes with the isometries fixing $c$ & H02, H03, H06\\
Okubo map equal to octonion map & Kaplansky's trick (Prop.~\ref{atl:prop:kaplansky}) & J01, J02, R17\\
Bilateral bodies with tilted ribs & bilateral map $=$ octonionic map moved by a similarity (split cases: a pseudo-isometry and a dilation) & J04--J07\\
Four identical balls at $\theta=90^\circ$ & $x*x=x_0^2+|q|^2$ (Prop.~\ref{atl:prop:theta90}) & E04\\
Revolutions at $\theta=90^\circ$ not due to automorphisms & twisted equivariance on $W$ (Prop.~\ref{atl:prop:theta90}) & C05, R10, R15\\
Tricorn of revolution & $x\bullet x=\overline{x^2}$; $G_2\cdot\mathcal T$ (Thm.~B) & F02, G01\\
\bottomrule
\end{tabular}
\caption{\emph{Synoptic table of the visual signatures of the atlas} and of the algebraic properties that account for them.}
\label{atl:tab:signatures}
\end{table}

\subsubsection*{Conclusions and Future Developments}
In this work we have presented an exploratory atlas of the dynamics of very simple maps over the real eight-dimensional composition algebras, drawn in a four-slices mode that displays four orthogonal sections at once and measured quarter by quarter. It is mesmerizing and surprising how faithfully the pictures report the algebra: nearly every trait of the images --- the coincidence of two quarters, a body of revolution, a sphere, a cylinder, a hyperbolic horn, a three-armed body --- is the visible trace of an identity that can be stated in one line. The octonions, with their quadratic structure, appear as the algebra of revolutions, whose dynamics is equivariant under a group much larger than $G_2$; the para-octonions are their mirror image, with tricorns in place of Mandelbrot sets; the Okubo algebra, which has no unit and no para-unit, appears as the algebra of the torus, whose revolutions are displaced to the root planes of $\mathfrak{su}(3)$, whose orbits genuinely explore up to eight dimensions and whose Mandelbrot sets are the tricorn moved by $PSU(3)$.

The scope of the atlas must be stated with care. The measurements are Monte Carlo estimates on finite-iteration bodies; the propositions and the theorems of the main text explain the symmetries of algebraic origin, but we have not proved that no further symmetry exists; the four-dimensional mappings act only on the first four coordinates; for the split algebras the Euclidean bailout is not an escape criterion, and the rendered bodies are approximations. The questions that the first version of the atlas left open --- a description of the Okubonic Mandelbrot sets through the invariants $\n(c)$ and $\n(c,c*c)$, and a classification of the symmetries of the Okubonic sections through the stabilizers of $c$ --- are answered in the main text. Other questions remain open: the Mandelbrot sets with a seed $a\neq0$, which depend on the joint invariants of $(a,c)$; the Petersson algebras of generic angle, whose Mandelbrot orbits leave the plane of $c$ and $c*c$; the split Okubo algebra, whose Cartan subgroups are not all conjugate. The next steps of our investigation will clearly include an eight-dimensional version of the mappings and the extension of the atlas to the Albert algebras built over the Okubo algebra, that we plan to cover in a forthcoming work.

\subsection{Reproducibility of the plates}\label{atl:app:repro}
The plates were produced by the following scripts, which are among the ancillary files described in Appendix~\ref{app:repro}:
\begin{itemize}
\item \path{plates_v1.mjs} and \path{plates_v2.mjs}: the specification of the plates (v2 moves Plate P02 to the frame $F_A$);
\item \path{atlas_render_v1.mjs}: rendering and measurements. The rendering uses WebGL~2 on an NVIDIA RTX 4070 SUPER through ANGLE/Direct3D~11, at $1800\times1200$ pixels with $4$ samples per pixel; the measurements run on the CPU in double precision;
\item \path{verify_propositions_v3.py}: numerical verification of the identities of Section~\ref{atl:sec:algebras} on random samples (numpy, double precision), with the certified algebra of the HF8D~0.2 prototype;
\item \path{summarize_v1.py} and \path{build_atlas_v1.py}: digest of the measurements, assembly and compilation of the first version;
\item \path{plates_v3.mjs} and \path{atlas_render_v2.mjs}: the plates R01--R17 of version~2, with arbitrary orthonormal frames;
\item \path{hf8np_v1.py} and \path{measure_v2.py}: the independent implementation of the algebras and the measurements of version~2 (the tables of this atlas);
\item \path{verify_theory_v2.py} and \path{verify_extra_v1.py}: the verification of the statements of the main text;
\item \path{build_supplement_v2.py}: assembly and compilation of this atlas.
\end{itemize}
The recorded runs are in the folder \path{results} of the ancillary files. Each one has a \path{manifest.json} with the command, the environment, the commit of HyperFractal~8D, the input and output hashes, the conclusions and the limitations.
\begin{itemize}
\item \path{run_atlas_render_2026-09-29T17-56-41}: the 116 plates of \path{plates_v1.mjs} (Plate P02 of this run is superseded);
\item \path{run_atlas_render_2026-09-29T18-06-48}: Plate P02 in the frame $F_A$;
\item \path{run_verify_propositions_2026-09-29T18-30-40}: 21 checks, all passed (the runs \path{...T18-29-48} and \path{...T18-31-18} give the same \path{results.json}; the manifest of the first one quotes the command of version~2 by mistake);
\item \path{run_atlas_render_2026-09-30T07-23-18}: the plates R01--R17 (their render-time measurements are those of version~1, $2400$ points; the tables of this atlas use the next run);
\item \path{run_measure_v2_2026-09-30T13-09-48}: the measurements of all plates except Part~M ($10^5$ points per quarter);
\item \path{run_verify_theory_2026-09-30T07-21-50} and \path{run_verify_extra_2026-09-30T07-47-42}: the checks of the main text.
\end{itemize}
Two earlier verification runs document a correction. Run \path{...T17-59-45} (version~1) recorded one failure, the conjecture that in the dual plane the critical orbit escapes exactly at $\re c=-3/4$; the corrected statement is Remark~\ref{atl:rem:dual}, checked by version~2 in \path{...T18-01-45}. The renderer is the application HyperFractal~8D~1.0 \cite{HF8D}, which is not among the ancillary files; its eight-dimensional shaders were verified against its double-precision CPU kernel on 270 scenes (run \path{2026-09-29T17-05-08_hf8d-v1-gpu}).

\endgroup

\end{document}